\documentclass[12pt]{book}  
  
  \usepackage{makeidx}
  
  \newtheorem{lem}{Lemma}[section]  
 \newtheorem{thm}{Theorem}[section]  
 \newtheorem{rmk}{Remark}[section]  
 \newtheorem{cor}{Corollary}[section]  
 \newtheorem{dfn}{Definition}[section]  
 \newtheorem{exm}{Example}[section]

 \newtheorem{cnj}{Conjecture}[section]  

\newcommand{\mod}{\hbox{\bf mod}}
\newcommand{\tr}{\hbox{\bf Tr}}
\newcommand{\Sing}{\hbox{\bf Sing}}
\newcommand{\dete}{\hbox{\bf det}}

\newcommand{\Mod}{\hbox{\bf Mod}}
\newcommand{\Tors}{\hbox{\bf Tors}}
\newcommand{\Coh}{\hbox{\bf Coh}}
\newcommand{\Spec}{\hbox{\bf Spec}}
\newcommand{\Irred}{\hbox{\bf Irred}}

\newcommand{\ka}{\hbox{\bf k}}
\newcommand{\n}{\hbox{\bf n}}

\input{amssym.tex}  

\makeindex
\begin{document}


\title{NOTES ON NONCOMMUTATIVE GEOMETRY}
\author{Igor ~Nikolaev}

\date{{\sf E-mail: igor.v.nikolaev@gmail.com}}

\maketitle


\pagenumbering{roman}




\chapter*{Foreword}
These  notes  are neither  an introduction nor  a survey (if only a brief) of  noncommutative 
geometry -- later NCG;  rather,  they  strive  to answer  some naive  but  vital  
questions:

\bigskip
{\it   What is the purpose of NCG  and  what is it good for? 
Why a number theorist or an algebraic geometer should
care about the  NCG?  Can NCG solve open problems of classical geometry 
inaccessible  otherwise?    In other words,  why does NCG matter?  
What is it anyway?}

\bigskip
Good answer means good examples.  A sweetheart  of  NCG   called  
noncommutative torus  captures classical  geometry of elliptic curves 
 because such a  torus is a coordinate ring  for elliptic  curves. 
  In other words,  one deals with  a  functor from algebraic  geometry to the NCG;   
  such functors  are at the heart  of  our  book.

What is NCG anyway?    It is  a  calculus   of  functors  on
the classical spaces  (e.g. algebraic, geometric, topological,  etc) with the values 
in NCG.   Such an approach   departs  from  the  tradition 
 of recasting   geometry of the classical space $X$  in terms of  the   $C^*$-algebra $C(X)$ 
 of  continuous   complex-valued functions on $X$,  see  the monograph by   [Connes 1994]  \cite{C}.

 \index{NCG}
 \index{elliptic curve}
 \index{noncommutative torus}
 \index{functor}

\medskip
\bigskip\noindent
\centerline{{\it Boston,  March  2014 \hskip8cm  Igor Nikolaev}}






\chapter*{Introduction}
It is  not easy  to write an elementary introduction to the NCG because
the simplest  non-trivial  examples are mind-blowing and involve the $KK$-groups,   subfactors,  
Sklyanin algebras, etc.   Such a  material  cannot be shrink-wrapped into a single graduate course
being   the result of a slow roasting of ideas from  different (and distant)  mathematical areas.  
But this universality of NCG is a thrill  giving the reader a long lost sense 
of unity of  mathematics.
In writing these notes the author had in mind a graduate student in (say) 
number theory  eager to learn  something new,   e.g.   a  noncommutative
torus with real multiplication;  it will soon transpire that  such an object 
is linked to the $K$-rational points of  elliptic curves and the Langlands program. 

 \index{Langlands program}
 \index{real multiplication}
 \index{elliptic curve}
 \index{$K$-rational points}

The book has three parts.   Part I  is preparatory:  Chapter 1 deals with  the
simplest examples of functors arising in algebraic geometry, number theory
and topology;  the functors take value in a category of the $C^*$-algebras known as
noncommutative tori.   Using these functors one gets a set of  noncommutative invariants          
for elliptic curves and Anosov's  automorphisms of the two-dimensional  torus.   
Chapter 2 is a brief introduction to the categories, functors and natural transformations;
they will be used throughout the book.  Chapter 3 covers an essential information 
about the category of $C^*$-algebras and their $K$-theory;  we introduce 
certain important classes of the $C^*$-algebras:  the AF-algebras, the UHF-algebras
and the Cuntz-Krieger algebras.  Our choice of the $C^*$-algebras is motivated 
by their applications in Part II.

Part II deals with the noncommutative invariants obtained from the functors acting
on various classical spaces.  Chapter 4 is devoted to such functors on the topological
spaces with values in the category of the so-called stationary AF-algebras;  
the noncommutative invariants  are the Handelman triples $(\Lambda, [I], K)$,
where $\Lambda$ is an order in a real algebraic number field $K$ and $[I]$ an
equivalence class of the ideals of $\Lambda$.  Chapter 5 deals with the examples of
 functors arising in projective algebraic geometry and their noncommutative invariants. 
 Finally,  Chapter 6 covers functors in   number theory and the corresponding  
invariants.

Part III is a brief  survey of the NCG;   the survey is cursory  yet  an extensive guide to the
  literature has been compiled at the end of each chapter.  
We hope that the reader can instruct himself by looking at the original
publications;  we  owe an apology to the authors whose works are not on the 
list.

There exist  several  excellent textbooks  on the NCG.  The  first and foremost 
is  the  monograph  by A.~Connes {\it ``G\'eom\'etrie Non Commutative''},  Paris, 1990 and 
its English edition {\it ``Noncommutative Geometry''},  Academic Press, 1994. 
The books  by J.~Madore {\it ``An Introduction to Noncommutative  Differential Geometry \& its Applications''}, 
Cambridge Univ. Press, 1995,      by J.~M.~Gracia-Bondia, J.~C.~Varilly and H.~Figueroa {\it ``Elements of Noncommutative Geometry''},
  Birkh\"auser, 2000 and  by M.~Khalkhali   {\it ``Basic Noncommutative  Geometry''},  EMS Series of Lectures
  in Mathematics, 2007   treat  particular  aspects of   Connes'   monograph.   A different approach to  the NCG is covered in a small
  but  instructive  book  by  Yu.~I.~Manin   {\it ``Topics  in Noncommutative Geometry''},  Princeton Univ. Press, 1991.  
Finally,  a more specialized {\it ``Noncommutative Geometry, Quantum Fields and Motives''},  AMS Colloquium Publications, 2008 
by A.~Connes and M.~Marcolli  is devoted to the links to physics and  number theory.  None of these books  treat  the NCG
 as a functor \cite{Nik0}.

I thank  the organizers, participants and sponsors of the  Spring Institute on {\it  Noncommutative Geometry and 
Operator Algebras (NCGOA)}  held annually at   the  Vanderbilt University in  Nashville,  Tennessee;  
these notes  grew from efforts  to find  out  what  is going on  there.   (I still  don't  have an answer.)   
  I  am grateful  to folks  who helped me with the project;   among them are  D.~Anosov,  P.~Baum, 
D.~Bisch,  B.~Blackadar,  O.~Bratteli,   A.~Connes,  J.~Cuntz,  G.~Elliott,  K.~Goodearl,  D.~Handelman, N.~Higson,
B.~Hughes, V.~F.~R.~Jones,  M.~Kapranov,  M.~Khalkhali,  W.~Krieger,  Yu.~Manin,  V.~Manuilov,  M.~Marcolli,  
V.~Mathai,  A.~Mishchenko,  S.~Novikov, N. ~C. ~Phillips,  M.~Rieffel,  W.~Thurston,  V.~Troitsky,   G.~Yu and others.



\tableofcontents  
\pagenumbering{arabic}                         




\part{BASICS}


\chapter{Model  Examples}
We shall start with  the simplest  functors arising in algebraic geometry,
number theory and topology;   all  these  functors range  in  a  category
of the $C^*$-algebras called  noncommutative tori.    We strongly believe that a handful of
simple  examples  tell  more than lengthy  theories  based on them;  we encourage
the reader   to keep  these model  functors  in mind for the rest    of the book.     
No special knowledge of the $C^*$-algebras,  elliptic curves or Anosov automorphisms
(beyond  an intuitive level)   is required at this point;   the interested reader can look up the missing
definitions in the standard literature  indicated  at the end of each section.

 \index{noncommutative torus}

\section{Noncommutative torus}
The noncommutative torus is an associative  algebra over ${\Bbb C}$
of particular simplicity and beauty;  such an algebra can be defined in several equivalent
ways,  e.g. as the universal algebra ${\Bbb C}\langle u, v\rangle$ on two unitary generators
$u$ and $v$ satisfying the unique commutation relation $vu=e^{2\pi i\theta} uv$, 
where $\theta$ is a real number.  There is a more geometric introduction as  a deformation
of the commutative algebra $C^{\infty} (T^2)$ of smooth  complex-valued functions on the 
two-dimensional torus $T^2$;  we shall pick up the latter because it clarifies the origin
and notation for such algebras.    
Roughly speaking,  one starts with the commutative algebra $C^{\infty}(T^2)$ of infinitely
differentiable complex-valued functions on  $T^2$ endowed with the usual pointwise
sum  and product  of two functions.  The idea is to replace the 
commutative product $f(x)g(x)$ of functions $f,g\in C^{\infty}(T^2)$  by a non-commutative product 
$f(x)\ast_{\hbar} g(x)$  depending on a continuous deformation parameter $\hbar$,  so that $\hbar=0$
corresponds to the usual product $f(x)g(x)$;  the product   $f(x)\ast_{\hbar} g(x)$ must be associative 
for each value of $\hbar$.  To  achieve the goal,  it is sufficient to construct the Poisson bracket 
$\{f, g\}$   on $C^{\infty}(T^2)$,  i.e. a binary operation satisfying the identities $\{f,f\}=0$ and  
$\{f, \{g,h\}\} +\{h, \{f, g\}\}+\{g, \{h,f\}\}=0$;   the claim  is a special case of Kontsevich's Theorem for 
the Poisson manifolds,  see Section 14.3.  The algebra $C_{\hbar}^{\infty} (T^2)$ equipped with the usual 
sum  $f(x)+g(x)$ and a non-commutative associative  product  $f(x)\ast_{\hbar} g(x)$ is called a 
{\it deformation quantization} of algebra $C^{\infty}(T^2)$. 

 \index{deformation quantization}

 \index{Poisson bracket}

The required Poisson bracket can be constructed as follows.   For a real number  $\theta$ define a 
 bracket on   $C^{\infty}(T^2)$ by the formula
\displaymath
\{f,g\}_{\theta}:=\theta\left({\partial f\over\partial x}  {\partial g\over\partial y}-
{\partial f\over\partial y} {\partial g\over\partial x}\right).
\enddisplaymath
The reader is encouraged to verify that the bracket satisfies the identities $\{f,f\}_{\theta}=0$ and  
$\{f, \{g,h\}_{\theta}\}_{\theta} +\{h, \{f, g\}_{\theta}\}_{\theta}+\{g, \{h,f\}_{\theta}\}_{\theta}=0$,  
i.e. is the Poisson bracket.  The Kontsevich Theorem
says that there exists an associative product   $f\ast_{\hbar} g$ on $C^{\infty}(T^2)$
obtained  from the bracket $\{f, g\}_{\theta}$.   Namely,   let $\varphi$ and $\psi$ denote  the Fourier 
transform of functions $f$ and $g$ respectively;   one can define an $\hbar$-family of  products
between  the Fourier transforms according to the formula
\displaymath
(\varphi\ast_{\hbar}\psi)(p)=\sum_{q\in {\Bbb Z}^2}
\varphi(q)\psi(p-q) e^{-\pi i\hbar ~k(p,q)},
\enddisplaymath
where  $k(p,q)=\theta (pq-qp)$ is the kernel of the Fourier transform of the Poisson bracket  
$\{f,g\}_{\theta}$,  i.e.  an expression defined by  the formula
\displaymath
\{\varphi, \psi\}_{\theta}=-4\pi^2 \sum_{q\in {\Bbb Z}^2}
\varphi(q) \psi(p-q) k(q, p-q).
\enddisplaymath
The product   $f\ast_{\hbar} g $  is defined as a pull back of the product 
$(\varphi\ast_{\hbar}\psi)$;   the resulting associative algebra   $C^{\infty}_{\hbar, \theta}(T^2)$ 
is called the deformation quantization of  $C^{\infty}(T^2)$ in the 
{\it direction $\theta$}  defined by the Poisson bracket $\{f,g\}_{\theta}$,
see Fig. 1.1.  
\begin{rmk}
\textnormal{
The algebra $C^{\infty}_{\hbar, \theta}(T^2)$ is endowed with the natural involution
coming from the complex conjugation on $C^{\infty}(T^2)$,  so that 
$\varphi^*(p):=\bar\varphi(-p)$ for all $p\in {\Bbb Z}^2$.  The natural norm 
on  $C^{\infty}_{\hbar, \theta}(T^2)$ comes from the operator norm of the Schwartz functions
$\varphi,\psi\in {\cal S}({\Bbb Z}^2)$  acting  on the  Hilbert space $\ell^2({\Bbb Z}^2)$.  
 }
\end{rmk}
\begin{dfn}
By a noncommutative torus ${\cal A}_{\theta}^{geometric}$ one understands the
$C^*$-algebra obtained from the norm closure of the  $^*$-algebra   $C^{\infty}_{1, ~\theta}(T^2)$.  
\end{dfn}
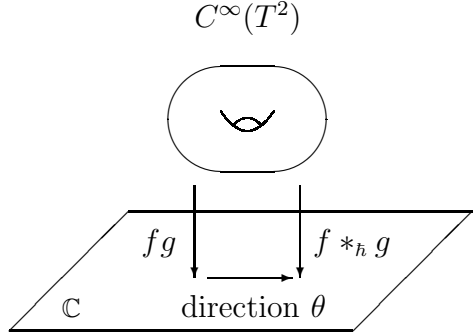
\begin{figure}[here]
\begin{picture}(400,180)(0,40)

\put(190,185){$C^{\infty}(T^2)$}
\put(170,100){$fg$}
\put(235,100){$f\ast_{\hbar} g$}
\put(140,75){${\Bbb C}$}
\put(185,75){direction $\theta$}


\put(210,150){\oval(60,40)}


\put(120,70){\line(1,0){130}}
\put(120,70){\line(1,1){45}}
\put(250,70){\line(1,1){45}}
\put(165,115){\line(1,0){130}}

\put(190,125){\vector(0,-1){35}}
\put(230,125){\vector(0,-1){35}}
\put(195,90){\vector(1,0){32}}


\qbezier(200,153)(210,138)(220,153)
\qbezier(205,148)(210,153)(215,148)

\end{picture}
\caption{Deformation of algebra $C^{\infty}(T^2)$.}
\end{figure}

\bigskip\noindent
The less visual {\it analytic} definition of the noncommutative torus involves 
bounded linear operators acting on the Hilbert space ${\cal H}$;  the reader can
think of the operators as the infinite-dimensional matrices over ${\Bbb C}$.
Namely,    let $S^1$ be the unit circle;   denote by $L^2(S^1)$ the Hilbert space of the 
square integrable complex valued functions on $S^1$.  Fix a real number $\theta\in [0,1)$;
for every $f(e^{2\pi it})\in L^2(S^1)$ we shall consider two bounded
linear operators $U$ and $V$  acting  by the formula  
\displaymath
\left\{
\begin{array}{cc}
U[f(e^{2\pi it})] &= f(e^{2\pi i(t-\theta)})\\
V[f(e^{2\pi it})]  &= e^{2\pi it}f(e^{2\pi it}).
\end{array}
\right.
\enddisplaymath
It is verified directly that 
\displaymath
\left\{
\begin{array}{cc}
VU  &= e^{2\pi i\theta}UV,\\
UU^* &= U^*U = E,\\
VV^* &= V^*V = E,
\end{array}
\right.
\enddisplaymath
where $U^*$ and $V^*$ are the adjoint operators of $U$ and $V$, respectively,
and $E$ is the identity operator. 
\begin{dfn}
By a noncommutative torus ${\cal A}_{\theta}^{analytic}$ one understands the
$C^*$-algebra  generated by the operators $U$ and $V$ acting on the 
Hilbert space $L^2(S^1)$.
 \end{dfn}

\bigskip\noindent
The {\it algebraic}  definition of the noncommutative torus is the shortest;
it involves the universal  algebras,  i.e. the associative  algebras given by 
the generators and relations.  Namely,    let ${\Bbb C}\langle x_1, x_2, x_3, x_4\rangle$ 
be the polynomial ring  in four non-commuting variables $x_1, x_2, x_3$ and $x_4$.  Consider 
a two-sided ideal, $I_{\theta}$, generated by the relations
\displaymath
\left\{
\begin{array}{cc}
x_3x_1  &= e^{2\pi i\theta}x_1x_3,\\
x_1x_2 &= x_2x_1 = e,\\
x_3x_4  &= x_4x_3 = e.
\end{array}
\right.
\enddisplaymath
\begin{dfn}
By a noncommutative torus ${\cal A}_{\theta}^{algebraic}$ one understands the
$C^*$-algebra  given  by the norm  closure of the $\ast$-algebra 
\linebreak
${\Bbb C}\langle x_1, x_2, x_3, x_4\rangle/I_{\theta}$, where the involution
acts on  the generators according to the formula $x_1^*=x_2$ and $x_3^*=x_4$.    
  \end{dfn}

\bigskip
\begin{thm}
${\cal A}_{\theta}^{geometric}\cong{\cal A}_{\theta}^{analytic}\cong{\cal A}_{\theta}^{algebraic}$
\end{thm}
{\it Proof.}  The isomorphism ${\cal A}_{\theta}^{analytic}\cong{\cal A}_{\theta}^{algebraic}$
is obvious,  because one can write $x_1=U, x_2=U^*, x_3=V$ and $x_4=V^*$.  
The isomorphism  ${\cal A}_{\theta}^{geometric}\cong{\cal A}_{\theta}^{analytic}$
is  established  by the identification of functions  $t\mapsto e^{2\pi i t p}$ 
of ${\cal A}_{\theta}^{geometric}$ with the unitary operators $U_p$ 
for each $p\in {\Bbb Z}^2$;   then the generators of ${\Bbb Z}^2$ 
will correspond to the operators $U$ and $V$.
$\square$

\begin{rmk}
\textnormal{
We shall write ${\cal A}_{\theta}$ to denote an abstract noncommutative torus
independent of its geometric, analytic or algebraic realization.  
 }
\end{rmk}

\bigskip\noindent
The noncommutative torus ${\cal A}_{\theta}$  has a plethora of remarkable properties;   for the moment we 
shall dwell on the most fundamental:  {\it Morita equivalence} and {\it real multiplication}.  
Roughly speaking,  the first property  presents  a basic equivalence relation in the category of noncommutative tori;
such a relation  indicates that   ${\cal A}_{\theta}$ and ${\cal A}_{\theta'}$ are identical   from the standpoint 
of noncommutative geometry.   The second property is rare;  only a countable family of non-equivalent
${\cal A}_{\theta}$ can have real multiplication. The property means that the ring of endomorphisms of
${\cal A}_{\theta}$ is non-trivial, i.e.  it  exceeds the ring ${\Bbb Z}$.   To give an exact definition,     
denote by ${\cal K}$ the $C^*$-algebra of all compact operators.
\begin{dfn}
The noncommutative torus ${\cal A}_{\theta}$ is said to be  stably isomorphic (Morita equivalent) to a noncommutative 
torus ${\cal A}_{\theta'}$ whenever   ${\cal A}_{\theta}\otimes {\cal K}\cong {\cal A}_{\theta'}\otimes {\cal K}$.
\end{dfn}
 \index{Morita equivalence}
 \index{stable isomorphism}
Recall that the Morita equivalence means that  the associative algebras $A$ and $A'$ have the
same category of projective modules, i.e. {\bf Mod} $(A)\cong$ {\bf Mod} $(A')$.  It is
notoriously  hard to tell (in intrinsic terms)  when  two non-isomorphic algebras are Morita equivalent;  of course,
if $A\cong A'$ then $A$ is Morita equivalent to $A'$.      The following remarkable result provides
a clear and definitive  solution  to the Morita equivalence problem for the noncommutative tori; 
it would be futile to talk about any links to the elliptic curves (complex tori) if a weaker or fuzzier 
result were true.        
\begin{thm}
{\bf (Rieffel)}
The noncommutative tori  ${\cal A}_{\theta}$ and  ${\cal A}_{\theta'}$ are stably isomorphic
(Morita equivalent)  if and only if 
\displaymath
\theta'={a\theta+b\over c\theta+d} \quad \hbox{for some matrix} 
 \quad\left(\matrix{a & b\cr c & d}\right)  \in SL_2({\Bbb Z}).
 \enddisplaymath
 \end{thm}

 \index{real multiplication}

\bigskip\noindent
The second fundamental property of the algebra ${\cal A}_{\theta}$ is the so-called 
{\it real multiplication};  such a multiplication signals  exceptional symmetry of ${\cal A}_{\theta}$. 
Recall that the Weierstrass  uniformization of elliptic curves by the lattices $L_{\tau}:={\Bbb Z}+{\Bbb Z}\tau$
gives rise to the {\it complex multiplication},  i.e.  phenomenon of an unusual behavior of
 the endomorphism ring of $L_{\tau}$;  the noncommutative torus ${\cal A}_{\theta}$ demonstrates  the same behavior
 with  (almost)  the same name. To introduce real multiplication, denote by ${\cal Q}$ the set of all quadratic irrational numbers, 
 i.e. the  irrational roots of  all quadratic polynomials with integer coefficients.  
\begin{thm}
{\bf (Manin)}
The  endomorphism  ring of a noncommutative torus ${\cal A}_{\theta}$ is given by the 
 formula
\displaymath
End~({\cal A}_{\theta})\cong
\cases{{\Bbb Z},  & \hbox{if} ~$\theta\in {\Bbb R}-({\cal Q}\cup {\Bbb Q})$\cr
           {\Bbb Z}+f O_k,  & \hbox{if}  ~$\theta\in {\cal Q}$,
           }
\enddisplaymath
where integer $f\ge 1$ is  conductor of an order in the ring of integers $O_k$
of  the real quadratic  field $k={\Bbb Q}(\sqrt{D})$.
\end{thm}
\begin{dfn}
The noncommutative torus ${\cal A}_{\theta}$ is said to have  real multiplication if  
$End~({\cal A}_{\theta})$ is bigger than ${\Bbb Z}$,  i.e.  $\theta$ is a quadratic 
irrationality;  we shall write ${\cal A}_{RM}^{(D, f)}$ to denote noncommutative tori with real
multiplication by an order of conductor $f$ in the quadratic field ${\Bbb Q}(\sqrt{D})$,
\end{dfn}
\begin{rmk}
\textnormal{
It is easy to see that real multiplication is an invariant of the stable isomorphism 
(Morita equivalence) class of noncommutative torus ${\cal A}_{\theta}$.   
 }
\end{rmk}

 \index{quadratic irrationality}
 \index{irrational rotation algebra}

\vskip1cm\noindent
{\bf Guide to the literature.}
For an authentic introduction to the noncommutative tori  we encourage the reader to  start  with the survey paper 
by  [Rieffel 1990]   \cite{Rie1}.    The noncommutative torus is also known as the irrational rotation 
algebra,   see  [Pimsner \& Voiculescu  1980]   \cite{PiVo1} and  [Rieffel 1981]  \cite{Rie2}.  
The real multiplication has been introduced  and studied by  [Manin 2003] \cite{Man1}.

 \index{elliptic curve}

\section{Elliptic curves}
Elliptic curves are so fundamental that they hardly need any introduction;
many basic facts and open problems of complex analysis, algebraic 
geometry and number theory can be reformulated in terms of such 
curves.  Perhaps it is the single most ancient mathematical object so
well explored yet hiding the deepest unsolved problems, e.g. the Birch and
Swinnerton-Dyer Conjecture.   Unless otherwise stated, we deal with the elliptic 
curves over the   field  ${\Bbb C}$  of complex numbers;  recall that   an {\it elliptic curve}
is the subset of the complex projective plane of the form
\displaymath
{\cal E}({\Bbb C})=\{(x,y,z)\in {\Bbb C}P^2 ~|~ y^2z=4x^3+axz^2+bz^3\},
\enddisplaymath
where $a$ and  $b$  are some constant complex numbers. The  real points of 
${\cal E}({\Bbb C})$  are depicted in Figure 1.2.
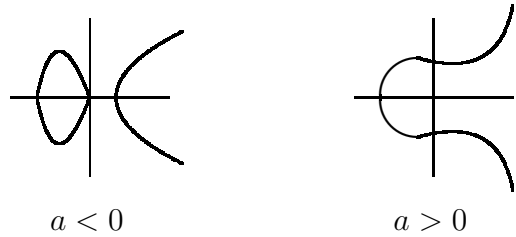
\begin{figure}[here]
\begin{picture}(300,100)(0,-10)

\put(100,50){\line(1,0){60}}
\put(230,50){\line(1,0){60}}
\put(130,20){\line(0,1){60}}
\put(260,20){\line(0,1){60}}

\thicklines
\qbezier(165,25)(115,50)(165,75)
\qbezier(110,50)(117,15)(130,50)
\qbezier(110,50)(117,85)(130,50)

\put(255,50){\oval(30,30)[l]}
\qbezier(254,65)(285,56)(290,85)
\qbezier(254,35)(285,44)(290,15)

\put(115,0){$a<0$}
\put(245,0){$a>0$}
\end{picture}
\caption{The real points of an affine elliptic curve $y^2=4x^3+ax$.}
\end{figure}
 \index{quadric surface}
\begin{rmk}
\textnormal{
It is known  that each elliptic curve ${\cal E}({\Bbb C})$  is isomorphic to the set
of points of intersection of two {\it quadric surfaces} in the complex projective space ${\Bbb C}P^3$
given by the  system of homogeneous equations 
\displaymath
\left\{
\begin{array}{ccc}
u^2+v^2+w^2+z^2 &=&  0,\\
Av^2+Bw^2+z^2  &=&  0,   
\end{array}
\right.
\enddisplaymath
where $A$ and $B$ are some constant  complex numbers and  
$(u,v,w,z)\in {\Bbb C}P^3$;  the system  is called
the {\it Jacobi form} of elliptic curve ${\cal E}({\Bbb C})$.   
}
\end{rmk}
 \index{Jacobi elliptic curve}
 \index{complex torus}
 \index{complex modulus}
\begin{dfn}
 By a  complex torus  one understands  the space  ${\Bbb C}/({\Bbb Z}\omega_1+{\Bbb Z}\omega_2)$,  
 where  $\omega_1$ and $\omega_2$ are   linearly independent
vectors in the complex plane ${\Bbb C}$,  see Fig. 1.3;  the ratio $\tau=\omega_2/\omega_1$ is called a  
complex modulus.
 \end{dfn}
\begin{figure}[here]
\begin{picture}(300,100)(0,0)

\put(90,90){\line(1,0){80}}
\put(80,70){\line(1,0){80}}
\put(60,30){\line(1,0){80}}

\put(85,20){\line(1,2){40}}
\put(65,20){\line(1,2){40}}
\put(105,20){\line(1,2){40}}
\put(125,20){\line(1,2){40}}

\put(180,50){\vector(1,0){40}}

\thicklines
\put(70,50){\line(1,0){80}}
\put(100,20){\line(0,1){80}}

\put(100,50){\vector(1,0){20}}
\put(100,50){\vector(1,2){10}}
\put(110,70){\line(1,0){20}}
\put(120,50){\line(1,2){10}}


\put(290,50){\oval(60,40)}
\qbezier(280,53)(290,38)(300,53)
\qbezier(285,48)(290,53)(295,48)


\put(260,90){${\Bbb C}/({\Bbb Z}+{\Bbb Z}\tau)$}
\put(70,90){${\Bbb C}$}
\put(120,35){$1$}
\put(105,75){${\tau}$}

\put(185,55){factor}
\put(190,40){map}

\end{picture}
\caption{Complex torus  ${\Bbb C}/({\Bbb Z}+{\Bbb Z}\tau)$.}
\end{figure}
\begin{rmk}
\textnormal{
Two complex  tori   ${\Bbb C}/({\Bbb Z}+{\Bbb Z}\tau)$ and   ${\Bbb C}/({\Bbb Z}+{\Bbb Z}\tau')$ 
are isomorphic   if and only if 
\displaymath
\tau'={a\tau+b\over c\tau+d} \quad \hbox{for some matrix} 
 \quad\left(\matrix{a & b\cr c & d}\right)  \in SL_2({\Bbb Z}).
 \enddisplaymath
 (We leave the proof  to the reader.  Hint:  notice that $z\mapsto\alpha z$
is an invertible holomorphic map for each $\alpha\in {\Bbb C}-\{0\}$.)
}
\end{rmk}

\bigskip\noindent
One may wonder if the complex analytic manifold ${\Bbb C}/({\Bbb Z}+{\Bbb Z}\tau)$
can be embedded into an  $n$-dimensional complex  projective space as an algebraic variety;
it turns out that the answer is emphatically yes even for the case $n=2$.  The following 
classical result relates complex torus    ${\Bbb C}/({\Bbb Z}+{\Bbb Z}\tau)$ with an elliptic
curve ${\cal E}({\Bbb C})$ in the projective plane ${\Bbb C}P^2$.  
\begin{thm}
{\bf (Weierstrass)}
There exists a holomorphic embedding  $${\Bbb C}/({\Bbb Z}+{\Bbb Z}\tau) \hookrightarrow  {\Bbb C}P^2$$ given by the formula
\displaymath
z\mapsto
\cases{(\wp(z), \wp'(z), 1) & \hbox{for}  $z\not\in L_{\tau}:={\Bbb Z}+{\Bbb Z}\tau$,\cr
          (0,1,0) & \hbox{for} $z\in L_{\tau}$},
\enddisplaymath
which is an isomorphism  between  complex torus ${\Bbb C}/({\Bbb Z}+{\Bbb Z}\tau)$ and elliptic curve 
\displaymath
{\cal E}({\Bbb C})=\{(x,y,z)\in {\Bbb C}P^2 ~|~ y^2z=4x^3+axz^2+bz^3\},
\enddisplaymath
where $\wp(z)$ is  the  Weierstrass  function defined by the convergent series
\displaymath
\wp(z)={1\over z^2}+\sum_{\omega\in L_{\tau}-\{0\}} \left({1\over (z-\omega)^2}-{1\over\omega^2}\right).
\enddisplaymath
and
\displaymath
\left\{
\begin{array}{ccc}
a &=&  -60\sum_{\omega\in  L_{\tau}-\{0\}} {1\over\omega^4},\\
b &=&  -140\sum_{\omega\in L_{\tau}-\{0\}} {1\over\omega^6}.   
\end{array}
\right.
\enddisplaymath
\end{thm}
 \index{Weierstrass $\wp$ function}
\begin{rmk}
\textnormal{
Roughly speaking,  the Weierstrass Theorem   identifies elliptic curves ${\cal E}({\Bbb C})$ 
and  complex tori  ${\Bbb C}/({\Bbb Z}\omega_1+{\Bbb Z}\omega_2)$;
we shall write ${\cal E}_{\tau}$ to denote  elliptic curve corresponding to the
complex torus of modulus $\tau=\omega_2/\omega_1$. 
}
\end{rmk}
\begin{rmk}
{\bf (First encounter with functors)}
\textnormal{
The declared purpose of our notes were functors with the range in category of 
associative algebras;  as a model example we picked the category of noncommutative
tori.  Based on what is known about the algebra ${\cal A}_{\theta}$,  one  cannot avoid
 the following fundamental question: 
 {\it  Why the isomorphisms of  elliptic curves look exactly 
  the same as   the stable isomorphisms (Morita equivalences) of noncommutative tori?
 In other words,  what makes the diagram in Fig. 1.4 
 commute?}
}
\end{rmk}
\begin{figure}[here]
\begin{picture}(300,110)(-120,-5)
\put(20,70){\vector(0,-1){35}}
\put(130,70){\vector(0,-1){35}}
\put(45,23){\vector(1,0){60}}
\put(45,83){\vector(1,0){60}}
\put(15,20){${\cal A}_{\theta}$}
\put(0,50){$F$}
\put(145,50){$F$}
\put(123,20){${\cal A}_{\theta'={a\theta+b\over c\theta+d}}$}
\put(17,80){${\cal E}_{\tau}$}
\put(122,80){${\cal E}_{\tau'={a\tau+b\over c\tau+d}}$}
\put(60,30){\sf stably}
\put(50,10){\sf  isomorphic}
\put(50,90){\sf isomorphic}
\end{picture}
\caption{Fundamental phenomenon.}
\end{figure}
 \index{fundamental phenomenon}

\bigskip\noindent
To settle the problem,  we presume that the observed phenomenon is the part
of a  categorical correspondence between elliptic curves and noncommutative tori;  
the following theorem says that it is indeed so.
\begin{thm}
There exists a covariant functor $F$ from  the category of all elliptic curves ${\cal E}_{\tau}$
to the category of  noncommutative tori ${\cal A}_{\theta}$,  such that if  
${\cal E}_{\tau}$ is isomorphic to ${\cal E}_{\tau'}$ then 
${\cal A}_{\theta}=F({\cal E}_{\tau})$ is stably isomorphic (Morita equivalent) 
to ${\cal A}_{\theta'}=F({\cal E}_{\tau'})$.  
 \end{thm}
{\it Proof.}
We shall give an algebraic proof of this fact based on the notion of  a {\it Sklyanin algebra};
there exists a geometric proof using the notion of measured foliations and the Teichm\"uller 
theory, see Section 5.1.2.   
Recall that  the Sklyanin algebra   $S(\alpha,\beta,\gamma)$  is  a free   ${\Bbb C}$-algebra   on   four  generators  
$x_1,\dots,x_4$   and  six  quadratic relations: 
\displaymath
\left\{
\begin{array}{ccc}
x_1x_2-x_2x_1 &=& \alpha(x_3x_4+x_4x_3),\\
x_1x_2+x_2x_1 &=& x_3x_4-x_4x_3,\\
x_1x_3-x_3x_1 &=& \beta(x_4x_2+x_2x_4),\\
x_1x_3+x_3x_1 &=& x_4x_2-x_2x_4,\\
x_1x_4-x_4x_1 &=& \gamma(x_2x_3+x_3x_2),\\ 
x_1x_4+x_4x_1 &=& x_2x_3-x_3x_2,
\end{array}
\right.
\enddisplaymath
where $\alpha+\beta+\gamma+\alpha\beta\gamma=0$,   see  e.g.  [Smith \& Stafford  1992]  
\cite{SmiSta1},  p. 260.  The algebra $S(\alpha,\beta,\gamma)$ is isomorphic to a 
{\it (twisted homogeneous) coordinate ring}   of elliptic curve ${\cal E}_{\tau}\subset {\Bbb C}P^3$ 
 given in its  Jacobi form
\displaymath
\left\{
\begin{array}{ccc}
u^2+v^2+w^2+z^2 &=&  0,\\
{1-\alpha\over 1+\beta}v^2+{1+\alpha\over 1-\gamma}w^2+z^2  &=&  0;  
\end{array}
\right.
\enddisplaymath
 \index{Sklyanin algebra}
the latter means that   $S(\alpha,\beta,\gamma)$ satisfies the fundamental isomorphism
\displaymath
\hbox{{\bf Mod}}~(S(\alpha,\beta,\gamma))/\hbox{{\bf Tors}}\cong \hbox{{\bf Coh}}~({\cal E}_{\tau}),
\enddisplaymath
 where {\bf Coh} is  the category of quasi-coherent sheaves on ${\cal E}_{\tau}$, 
  {\bf Mod}  the category of graded left modules over the graded ring $S(\alpha,\beta,\gamma)$
 and  {\bf Tors}  the full sub-category of {\bf Mod} consisting of the
torsion modules,  see [Serre 1955] \cite{Ser1}.  
The algebra $S(\alpha,\beta,\gamma)$ defines a natural {\it automorphism} 
$\sigma: {\cal E}_{\tau}\to {\cal E}_{\tau}$ of the elliptic curve ${\cal E}_{\tau}$,  
see e.g.   [Stafford \& van ~den ~Bergh  2001]  \cite{StaVdb1}, p. 173. 
Fix  an automorphism $\sigma$ of order 4,  i.e. $\sigma^4=1$;    in  this   case  $\beta=1, ~\gamma=-1$
and  it is known that system of quadratic relations for the Sklyanin algebra $S(\alpha,\beta,\gamma)$ 
can be brought to a skew symmetric form
\displaymath
\left\{
\begin{array}{cc}
x_3x_1 &= \mu e^{2\pi i\theta}x_1x_3,\\
x_4x_2 &= {1\over \mu} e^{2\pi i\theta}x_2x_4,\\
x_4x_1 &= \mu e^{-2\pi i\theta}x_1x_4,\\
x_3x_2 &= {1\over \mu} e^{-2\pi i\theta}x_2x_3,\\
x_2x_1 &= x_1x_2,\\
x_4x_3 &= x_3x_4,
\end{array}
\right.
\enddisplaymath
where $\theta=Arg~(q)$ and $\mu=|q|$ for some complex number $q\in {\Bbb C} - \{0\}$,
see   [Feigin \& Odesskii  1989]  \cite{FeOd1},   Remark 1.

On the other hand,  the system of  relations involved in the algebraic 
definition of noncommutative torus ${\cal A}_{\theta}$ is equivalent to the following
system of quadratic relations 
\displaymath
\left\{
\begin{array}{cc}
x_3x_1 &=  e^{2\pi i\theta}x_1x_3,\\
x_4x_2 &=  e^{2\pi i\theta}x_2x_4,\\
x_4x_1 &=  e^{-2\pi i\theta}x_1x_4,\\
x_3x_2 &=   e^{-2\pi i\theta}x_2x_3,\\
x_2x_1 &= x_1x_2=e,\\
x_4x_3 &= x_3x_4=e.
\end{array}
\right.
\enddisplaymath
(We leave the proof to the reader as an exercise in  non-commutative algebra.)
Comparing these relations with the skew symmetric relations   for the Sklyanin algebra 
$S(\alpha, 1, -1)$,  one  concludes  that they are almost   identical;  to pin down the difference
we shall  add  two extra relations $$x_1x_3=x_3x_4={1\over\mu}e$$ to
relations of the Sklyanin algebra and bring it (by multiplication and cancellations) 
to the following equivalent form
\displaymath
\left\{
\begin{array}{cc}
x_3x_1x_4 &= e^{2\pi i\theta}x_1,\\
x_4 &= e^{2\pi i\theta}x_2x_4x_1,\\
x_4x_1x_3 &= e^{-2\pi i\theta}x_1,\\
x_2 &= e^{-2\pi i\theta}x_4x_2x_3,\\
x_2x_1 &= x_1x_2={1\over\mu}e,\\
x_4x_3 &= x_3x_4={1\over\mu}e.
\end{array}
\right.
\enddisplaymath
Doing the same type of equivalent transformations to the system
of relations for the noncommutative torus, one  brings the system to the form
\displaymath
\left\{
\begin{array}{cc}
x_3x_1x_4 &= e^{2\pi i\theta}x_1,\\
x_4 &= e^{2\pi i\theta}x_2x_4x_1,\\
x_4x_1x_3 &= e^{-2\pi i\theta}x_1,\\
x_2 &= e^{-2\pi i\theta}x_4x_2x_3,\\
x_2x_1 &= x_1x_2=e,\\
x_4x_3 &= x_3x_4=e.
\end{array}
\right.
\enddisplaymath
Thus the only difference between relations for the Sklyanin algebra
(modulo the ideal $I_{\mu}$ generated by relations  
 $x_1x_3=x_3x_4={1\over\mu}e$)  and such for the noncommutative
 torus ${\cal A}_{\theta}$ is a  {\it scaling of the unit} ${1\over\mu} e$.   
Thus one obtains the following remarkable isomorphism 
\displaymath
{\cal A}_{\theta}\cong S(\alpha,1,-1) ~/~ I_{\mu}.
\enddisplaymath

\smallskip
\begin{rmk}
{\bf (Noncommutative torus as coordinate ring of ${\cal E}_{\tau}$)}
\textnormal{
Roughly speaking,  the above formula  says  that modulo the ideal
$I_{\mu}$ the noncommutative torus ${\cal A}_{\theta}$  is a coordinate ring of elliptic curve 
${\cal E}_{\tau}$. 
}
\end{rmk}
The required functor $F$ can be obtained as a quotient map  of the fundamental (Serre)  isomorphism
\displaymath
I_{\mu}\backslash\hbox{{\bf Coh}}~({\cal E}_{\tau})\cong
\hbox{{\bf Mod}}~(I_{\mu}\backslash S(\alpha,1,-1))/\hbox{{\bf Tors}}\cong 
\hbox{{\bf Mod}}~({\cal A}_{\theta})/\hbox{{\bf Tors}}
\enddisplaymath
and the fact that the isomorphisms in category {\bf Mod} $({\cal A}_{\theta})$ correspond 
to the stable isomorphisms (Morita equivalences) of category ${\cal A}_{\theta}$.  
 $\square$

\vskip1cm\noindent
{\bf Guide to the literature.}
The reader can  enjoy a  plenty  of excellent literature 
 introducing  elliptic curves;   see e.g.   [Husem\"oller 1986]  \cite{H2}, 
  [Knapp 1992] \cite{K1},  [Koblitz 1984]  \cite{K2},   [Silverman 1985]  \cite{S1},  
   [Silverman 1994]   \cite{S2},  [Silverman \& Tate  1992]  \cite{ST}  and others.  
More advanced  topics are covered in the survey papers [Cassels 1966]  \cite{Cas1},
[Mazur  1986]  \cite{Maz1} and [Tate 1974]  \cite{Tat1}.  
Noncommutative tori as coordinate rings of elliptic curves were studied
in \cite{Nik1} and \cite{Nik2};  the higher genus curves were considered in
\cite{Nik3}.

 \index{complex multiplication}
 \index{Kronecker-Weber Theorem}
\section{Complex multiplication}
A central problem of algebraic number theory is to give an
explicit construction for the abelian extensions of a given field $k$.  
For instance,  if $k\cong {\Bbb Q}$ is the field of rational numbers,
the Kronecker-Weber Theorem says that the maximal abelian extension
of ${\Bbb Q}$ is the union of all cyclotomic extensions;  thus we have an
explicit {\it class field theory} over ${\Bbb Q}$.  If $k={\Bbb Q}(\sqrt{-D})$ 
is an imaginary quadratic field,  then  {\it complex multiplication} 
  realizes the class field theory for $k$.  Let us recall some
useful definition.       
\begin{dfn}
Elliptic curve ${\cal E}_{\tau}$ is said  to have complex multiplication if the endomorphism 
ring of ${\cal E}_{\tau}$  is bigger than ${\Bbb Z}$,  i.e. $\tau$ is a quadratic irrationality (see below);  
we shall write ${\cal E}_{CM}^{(-D, f)}$ to denote elliptic  curves with complex
multiplication by an order of conductor $f$ in the quadratic field ${\Bbb Q}(\sqrt{-D})$.  
\end{dfn}
\begin{rmk}
\textnormal{
The  endomorphism  ring 
$$End~({\cal E}_{\tau}):=\{\alpha\in {\Bbb C} ~:~ \alpha L_{\tau}\subseteq L_{\tau}\}$$
of elliptic curve  ${\cal E}_{\tau}={\Bbb C}/L_{\tau}$  is given by the  formula
\displaymath
End~({\cal E}_{\tau})\cong
\cases{{\Bbb Z},  & \hbox{if} ~$\tau\in {\Bbb C}-{\cal Q}$\cr
           {\Bbb Z}+f O_k,  & \hbox{if}  ~$\tau\in {\cal Q}$,
           }
\enddisplaymath
where ${\cal Q}$ is the set of all imaginary quadratic numbers and  integer $f\ge 1$ is  
conductor of an order in the ring of integers $O_k$  of  the imaginary quadratic  number field 
$k={\Bbb Q}(\sqrt{-D})$.
}
\end{rmk}

\medskip\noindent
In previous section we constructed a functor $F$ on elliptic curves ${\cal E}_{\tau}$ with the range in the
category of noncommutative tori  ${\cal A}_{\theta}$;  roughly speaking the following 
theorem characterizes the restriction of $F$ to elliptic curves with complex multiplication. 
\begin{thm}
{\bf (\cite{Nik4}, \cite{Nik5})}
$F({\cal E}_{CM}^{(-D,f)})={\cal A}_{RM}^{(D,f)}$.
 \end{thm}
\begin{rmk}
\textnormal{
It follows from the theorem that the associative algebra ${\cal A}_{RM}^{(D,f)}$
is a coordinate ring of elliptic curve  ${\cal E}_{CM}^{(-D,f)}$;   in other words, 
each geometric property of curve ${\cal E}_{CM}^{(-D,f)}$ can be expressed 
in terms of the noncommutative torus  ${\cal A}_{RM}^{(D,f)}$.  
}
\end{rmk}
 \index{rank of elliptic curve}
The noncommutative invariants of algebra ${\cal A}_{RM}^{(D,f)}$ are linked to
 {\it ranks} of the $K$-rational elliptic curves,  i.e.  elliptic curves  ${\cal E}(K)$
 over an algebraic field $K$;  to illustrate
the claim,  let us recall some basic definitions and facts.  Let $k={\Bbb Q}(\sqrt{-D})$
be an imaginary quadratic number field and let $j({\cal E}_{CM}^{(-D,f)})$
be the $j$-invariant of elliptic curve ${\cal E}_{CM}^{(-D,f)}$;   it is well
known from  complex multiplication, that 
\displaymath
{\cal E}_{CM}^{(-D,f)}\cong
{\cal E}(K),
\enddisplaymath
where $K=k (j({\cal E}_{CM}^{(-D,f)}))$ is the Hilbert class field of $k$, see
e.g.   [Silverman 1994]   \cite{S2}, p. 95.   The Mordell-Weil theorem  says that the set of 
the $K$-rational  points of ${\cal E}_{CM}^{(-D,f)}$ is a finitely generated abelian group,
see  [Tate 1974]  \cite{Tat1}, p. 192;  the rank of such a group will be denoted by   
 $rk~({\cal E}_{CM}^{(-D,f)})$.    
For the sake of simplicity,  we further  restrict to the following class  of curves.  
If     $({\cal E}_{CM}^{(-D,f)})^{\sigma},  ~\sigma\in Gal~(K|{\Bbb Q})$ is 
the  Galois conjugate of the curve  ${\cal E}_{CM}^{(-D,f)}$,  then
by a {\it ${\Bbb Q}$-curve}  one understands elliptic curve   ${\cal E}_{CM}^{(-D,f)}$,
such that  there exists  an  isogeny between $({\cal E}_{CM}^{(-D,f)})^{\sigma}$ and   ${\cal E}_{CM}^{(-D,f)}$
for  each    $\sigma\in Gal~(K|{\Bbb Q})$,  see e.g.  [Gross 1980]   \cite{G}.
Let ${\goth P}_{3 ~\mod ~4}$ be the set of all primes  
$p=3 ~\mod ~4$;    it is known that  ${\cal E}_{CM}^{(-p,1)}$ is a ${\Bbb Q}$-curve 
whenever $p\in  {\goth P}_{3 ~\mod ~4}$,  see [Gross 1980]   \cite{G}, p. 33. 
The rank of  ${\cal E}_{CM}^{(-p,1)}$  is always divisible by $2h_k$,  where $h_k$ is the 
 class number of field $k={\Bbb Q}(\sqrt{-p})$,  see [Gross 1980]   \cite{G},  p. 49;  
  by a {\it ${\Bbb Q}$-rank}  of  ${\cal E}_{CM}^{(-p,1)}$
 one understands  the integer 
 \index{${\Bbb Q}$-curve}
 \index{isogeny}
 \displaymath
 rk_{\Bbb Q}({\cal E}_{CM}^{(-p,1)}):={1\over 2h_k}~rk~({\cal E}_{CM}^{(-p,1)}).
 \enddisplaymath
\begin{dfn}
Suppose that  $[a_0, \overline{a_1,a_2,\dots,a_2,a_1, 2a_0}]$
 is the periodic continued fraction of  $\sqrt{D}$, see e.g.  [Perron 1954]  \cite{P},  p.83;    
 then   by an  arithmetic complexity $c({\cal A}_{RM}^{(D,f)})$  of  torus   ${\cal A}_{RM}^{(D,f)}$
one understands  the  total  number of independent  $a_i$  in  its  period 
$(a_1, a_2, \dots, a_2, a_1, 2a_0)$,  see Section 6.3.2 for the details.   
 \end{dfn}
 \index{arithmetic complexity}
\begin{rmk}
\textnormal{
It is easy to  see that  arithmetic complexity $c({\cal A}_{RM}^{(D,f)})$ is an invariant of the stable isomorphism 
(Morita equivalence) class of the noncommutative torus    ${\cal A}_{RM}^{(D,f)}$;
in other words,  $c({\cal A}_{RM}^{(D,f)})$ is a {\it noncommutative invariant}
of  torus  ${\cal A}_{RM}^{(D,f)}$.
}
\end{rmk}

\bigskip\noindent
The declared purpose of our notes were noncommutative invariants related
to the classical geometry of elliptic curves;   theorem below is one of
such statements for the  ${\Bbb Q}$-curves. 
\begin{thm}
{\bf (\cite{Nik6})}
$rk_{\Bbb Q}~({\cal E}_{CM}^{(-p,1)}) +1=c({\cal A}_{RM}^{(p,1)})$,
whenever $p=3 ~\mod ~4$.  
\end{thm}
\begin{rmk}
\textnormal{
It is known that there are infinitely many pairwise  non-isomorphic ${\Bbb Q}$-curves,  
see e.g.  [Gross 1980]   \cite{G};   all pairwise non-isomorphic ${\Bbb Q}$-curves 
${\cal E}_{CM}^{(-p,1)}$ with  $p<100$ and their noncommutative invariant
$c({\cal A}_{RM}^{(p,1)})$  are calculated  in Fig.1.5.
}
\end{rmk}
\begin{figure}[here]
\begin{tabular}{c|c|c|c}
\hline
&&&\\
$p\equiv 3~\mod~4$ & $rk_{\Bbb Q}({\cal E}_{CM}^{(-p,1)})$ & $\sqrt{p}$ & $c({\cal A}_{RM}^{(p,1)})$\\
&&&\\
\hline
$3$ & $1$ & $[1,\overline{1,2}]$ & $2$\\
\hline
$7$ & $0$ & $[2,\overline{1,1,1,4}]$ & $1$\\
\hline
$11$ & $1$ & $[3,\overline{3,6}]$ & $2$\\
\hline
$19$ & $1$ & $[4,\overline{2,1,3,1,2,8}]$ & $2$\\
\hline
$23$ & $0$ & $[4,\overline{1,3,1,8}]$ & $1$\\
\hline
$31$ & $0$ & $[5,\overline{1,1,3,5,3,1,1,10}]$ & $1$\\
\hline
$43$ & $1$ & $[6,\overline{1,1,3,1,5,1,3,1,1,12}]$ & $2$\\
\hline
$47$ & $0$ & $[6,\overline{1,5,1,12}]$ & $1$\\
\hline
$59$ & $1$ & $[7,\overline{1,2,7,2,1,14}]$ & $2$\\
\hline
$67$ & $1$ & $[8,\overline{5,2,1,1,7,1,1,2,5,16}]$ & $2$\\
\hline
$71$ & $0$ & $[8,\overline{2,2,1,7,1,2,2,16}]$ & $1$\\
\hline
$79$ & $0$ & $[8,\overline{1,7,1,16}]$ & $1$\\
\hline
$83$ & $1$ & $[9,\overline{9,18}]$ & $2$\\
\hline
\end{tabular}
\caption{The ${\Bbb Q}$-curves ${\cal E}_{CM}^{(-p,1)}$  with  $p<100$.}
\end{figure}

\vskip1cm\noindent
{\bf Guide to the literature.}
D.~Hilbert counted complex multiplication as not only the most beautiful
part of mathematics but also of entire science;  it surely  does as it links
complex analysis and number theory.  One cannot beat [Serre 1967] \cite{Ser2}
for an introduction,  but more comprehensive   [Silverman 1994]   \cite{S2}, Ch. 2
is the must.  Real multiplication has been introduced in [Manin 2004]  \cite{Man1}.  
The link between the two was  the subject of \cite{Nik4}.

 \index{Anosov automorphism}
\section{Anosov automorphisms}
Roughly speaking topology studies invariants of continuous maps $f: X\to Y$ between
the topological spaces $X$ and $Y$;  if $f$ is invertible and $X=Y$ we shall call it
an {\it automorphism}.  The automorphisms $f,f': X\to X$ are said to be {\it conjugate}
if there exists an automorphism $h: X\to X$  such that  $f'=h\circ f\circ h^{-1}$,  where
$f\circ f'$ means the composition of $f$ and $f'$;   the conjugation means a 
``change of coordinate system'' for the topological space $X$ and each property 
of $f$ invariant under the conjugation is {\it intrinsic},  i.e. a topological invariant 
of $f$.  The conjugation problem is unsolved even when $X$ is a topological surface
(compact two-dimensional manifold);  such a solution would imply topological 
classification of the three-dimensional manifolds,  see e.g. [Hemion 1979]  \cite{Hem1}.  
The automorphism $f: X\to X$ is said to have an {\it infinite order} if $f^n\ne Id$ 
for each $n\in {\Bbb Z}$.  Further we shall focus on the topological invariants 
of automorphisms $f$ when $X=T^2$ is the two-dimensional torus.  Because 
$T^2\cong {\Bbb R}^2/{\Bbb Z}^2$,  each automorphism of $T^2$ can be given by
an invertible map (isomorphism) of lattice ${\Bbb Z}^2\subset {\Bbb R}^2$, 
see Fig. 1.6;   in other words,  the automorphism $f: T^2\to T^2$ can be written 
in the matrix form
\displaymath
A_f=\left(\matrix{a_{11} & a_{12}\cr a_{21} & a_{22}}\right)\in GL(2, {\Bbb Z}).  
\enddisplaymath
\begin{figure}[here]
\begin{picture}(300,100)(0,0)
\put(180,50){\vector(1,0){40}}

\put(70,80){\line(1,0){80}}
\put(70,65){\line(1,0){80}}
\put(70,35){\line(1,0){80}}
\put(70,20){\line(1,0){80}}

\put(140,10){\line(0,1){80}}
\put(125,10){\line(0,1){80}}
\put(95,10){\line(0,1){80}}
\put(80,10){\line(0,1){80}}

\thicklines
\put(70,50){\line(1,0){80}}
\put(110,10){\line(0,1){80}}


\put(290,50){\oval(60,40)}
\qbezier(280,53)(290,38)(300,53)
\qbezier(285,48)(290,53)(295,48)


\put(270,90){${\Bbb R}^2/{\Bbb Z}^2$}
\put(150,90){${\Bbb Z}^2\subset {\Bbb R}^2$}

\put(185,55){factor}
\put(190,40){map}

\end{picture}
\caption{Topological  torus  $T^2\cong {\Bbb R}^2/{\Bbb Z}^2$.}
\end{figure}
\begin{dfn}
An infinite order automorphism $f: T^2\to T^2$ is called {\it Anosov} if 
 its  matrix form $A_f$ satisfies the inequality $|a_{11}+a_{22}|>2$. 
 \end{dfn}
\begin{rmk}
\textnormal{
The definition of Anosov's automorphism does not depend on the conjugation,
because the trace $a_{11}+a_{22}$ is an invariant of the latter. 
Moreover,  it is easy to see that  ``almost all''  automorphisms of $T^2$ are  Anosov's; 
they constitute the most interesting part among all automorphisms of the torus. 
  }
\end{rmk}
 \index{noncommutative torus}
To study topological invariants   we shall construct 
a functor $F$ on the set of all Anosov automorphisms with the values in a category
of the noncommutative tori such that the diagram in Fig. 1.7 is commutative;
in other  words, if $f$ and $f'$ are conjugate Anosov automorphisms,   then the corresponding
noncommutative tori   ${\cal  A}_{\theta}$ and  ${\cal  A}_{\theta'}$ are stably isomorphic 
(Morita equivalent).  
\begin{figure}[here]
\begin{picture}(300,110)(-110,-5)
\put(20,70){\vector(0,-1){35}}
\put(130,70){\vector(0,-1){35}}
\put(45,23){\vector(1,0){53}}
\put(45,83){\vector(1,0){53}}
\put(10,20){${\cal A}_{\theta}$}
\put(0,50){$F$}
\put(140,50){$F$}
\put(120,20){${\cal A}_{\theta'}$}
\put(17,80){$f$}
\put(117,80){$f'=h\circ f\circ h^{-1}$}
\put(58,30){\sf stable}
\put(45,10){\sf isomorphism}
\put(44,90){\sf conjugation}
\end{picture}
\caption{Functor  $F$.}
\end{figure}
 \index{Perron-Frobenius eigenvalue}
 \index{Perron-Frobenius eigenvector}
The required map $F: A_f\mapsto {\cal A}_{\theta}$ can be constructed  as follows.
For simplicity,  we shall assume that  $a_{11}+a_{22}>2$;   the case $a_{11}+a_{22}<-2$ 
is treated similarly.  Moreover,  we can assume that $A_f$ is a positive matrix since
each class of conjugation of  the Anosov automorphism $f$ contains such a representative;
denote by $\lambda_{A_f}$ the Perron-Frobenius eigenvalue of  positive matrix $A_f$. 
The noncommutative torus ${\cal A}_{\theta}=F(A_f)$ is defined by the normalized 
Perron-Frobenius eigenvector $(1,\theta)$ of  the matrix $A_f$,  i.e. 
 \displaymath
 A_f\left(\matrix{1\cr\theta}\right)=\lambda_{A_f}\left(\matrix{1\cr\theta}\right).
 \enddisplaymath
 \index{positive matrix}
\begin{rmk}
\textnormal{
We leave it to the reader  to prove that  if  $f$ is Anosov's,   then   $\theta$ is a quadratic 
irrationality given  by the formula
 \displaymath
 \theta={a_{22}-a_{11}+\sqrt{(a_{11}+a_{22})^2-4}\over 2a_{12}}. 
 \enddisplaymath
 }
\end{rmk}

 \index{real multiplication}

\bigskip\noindent
It follows from the above formula   that   map $F$  takes values in the  noncommutative tori 
with real multiplication;    the following theorem says that our  map  $F: A_f\mapsto {\cal A}_{\theta}$ 
is  actually  a functor.  (The proof of this fact is an easy exercise for anyone familiar with 
the notion of the $AF$-algebra of stationary type, the Bratteli diagram, etc.;   we refer the interested 
reader to \cite{Nik0}, p.153  for a short proof.) 
\begin{thm}
{\bf (\cite{Nik7})}
If $f$ and $f'$ are conjugate Anosov automorphisms,   then  the noncommutative torus  
${\cal A}_{\theta}=F(A_f)$ is stably  isomorphic  (Morita equivalent)  to ${\cal A}_{\theta'}=F(A_{f'})$.  
 \end{thm}
Thus  the problem of conjugation for the Anosov automorphisms can  be recast   in  terms of the 
noncommutative tori;   namely, one  needs   to find  invariants of the stable isomorphism class of a noncommutative
torus with real multilplication.   Such a noncommutative  invariant has been calculated in 
[Handelman 1981]  \cite{Han2};  namely,  consider the eigenvalue problem for  a  matrix $A_f\in GL(2, {\Bbb Z})$,
i.e.  $A_f v_A=\lambda_{A_f} v_A$,  where $\lambda_{A_f}>1$ is the Perron-Frobenius eigenvalue  and 
$v_A=(v_A^{(1)},v_A^{(2)})$ the corresponding eigenvector with the positive entries  normalized so that 
$v_A^{(i)}\in K={\Bbb Q}(\lambda_{A_f})$.   Denote by ${\goth m}={\Bbb Z}v_A^{(1)}+{\Bbb Z}v_A^{(2)}$
 a ${\Bbb Z}$-module in the number field $K$. Recall that the coefficient
ring,  $\Lambda$,  of module ${\goth m}$ consists of the elements $\alpha\in K$
such that $\alpha {\goth m}\subseteq {\goth m}$.  It is known that 
$\Lambda$ is an order in $K$ (i.e. a subring of $K$  containing $1$) and, with no restriction, 
one can assume that  ${\goth m}\subseteq\Lambda$. It follows  from the definition, that ${\goth m}$
coincides with an ideal, $I$, whose equivalence class in $\Lambda$ we shall denote
by $[I]$. 
\begin{thm}
{\bf (Handelman's noncommutative invariant)}
The triple $(\Lambda, [I], K)$  is an arithmetic invariant of the  stable isomorphism class of 
the noncommutative torus ${\cal  A}_{\theta}$ with real multiplication,  i.e.  the tori  ${\cal  A}_{\theta}$ and  
${\cal A}_{\theta'}$ are stably isomorphic (Morita equivalent)  if and only if 
$\Lambda=\Lambda', [I]=[I']$ and $K=K'$. 
 \end{thm}
 \index{Handelman invariant}
\begin{rmk}
\textnormal{
The Handelman Theorem was proved for the so-called $AF$-algebras of a stationary type;
such algebras and the noncommutative tori with real multiplication are known to have the 
same $K_0^+$ semigroup  and therefore the same classes of stable isomorphisms.    
Similar problem for matrices  was solved in  [Latimer  \&  MacDuffee 1933]   \cite{LaMa1}
and  [Wallace 1984]   \cite{Wal1}. 
 }
\end{rmk}
Handelman's Invariant  $(\Lambda, [I], K)$  gives rise to  a series of numerical invariants 
of the conjugation class of Anosov's  automorphisms;  we shall consider one such invariant
called {\it module determinant} $\Delta ({\goth m})$.  Let   ${\goth m}={\Bbb Z}v_A^{(1)}+{\Bbb Z}v_A^{(2)}$
be  the  module attached to $\Lambda$;    consider the symmetric bilinear form
 \displaymath
 q(x,y)=\sum_{i=1}^2 \sum_{j=1}^2  Tr~(v_A^{(i)} v_A^{(j)})  x_i x_j,
 \enddisplaymath
where $Tr~(v_A^{(i)} v_A^{(j)})$ is the trace of the algebraic number  $v_A^{(i)} v_A^{(j)}$.  
\begin{dfn}
By a determinant of module ${\goth m}$ one understands the determinant of the bilinear
form $q(x,y)$, i.e. the rational integer
 \displaymath
 \Delta ({\goth m}):=  
 Tr~(v_A^{(1)} v_A^{(1)})   
 ~Tr~(v_A^{(2)} v_A^{(2)})   
-  Tr^2~(v_A^{(1)} v_A^{(2)}).
  \enddisplaymath
\end{dfn}
\begin{rmk}
\textnormal{
The rational integer  $\Delta ({\goth m})$ is a numerical    invariant 
of Anosov's  automorphisms,   because it does not depend on the basis of module 
 ${\goth m}={\Bbb Z}v_A^{(1)}+{\Bbb Z}v_A^{(2)}$;   we leave the proof to the reader.
 }
\end{rmk}

 \index{Alexander polynomial}

\medskip\noindent
In conclusion,  we  calculate the noncommutative invariant   $\Delta ({\goth m})$
for the concrete automorphisms  $f$ of $T^2$;  the reader can see that in both  cases our
invariant  $\Delta ({\goth m})$ is {\it stronger}  than  the classical Alexander polynomial $\Delta(t)$,
i.e.   $\Delta ({\goth m})$ detects the topological classes of $f$   which invariant  $\Delta(t)$ cannot see. 
\begin{exm}
\textnormal{
Consider Anosov's automorphisms $f_A, f_B: T^2\to T^2$ given by matrices  
\displaymath
A=\left(\matrix{5 & 2\cr 2 & 1}\right)\qquad  \hbox{and}
\qquad B=\left(\matrix{5 & 1\cr 4 & 1}\right), 
\enddisplaymath
respectively. The Alexander polynomials of $f_A$ and $f_B$ are identical   
$\Delta_A(t)=\Delta_B(t)= t^2-6t+1$;  yet  the automorphisms $f_A$ and $f_B$ are 
{\it not}  conjugate. 
Indeed, the Perron-Frobenius eigenvector of matrix $A$ is  $v_A=(1, \sqrt{2}-1)$
while of the matrix $B$ is $v_B=(1, 2\sqrt{2}-2)$. The  bilinear forms for the modules 
${\goth m}_A={\Bbb Z}+(\sqrt{2}-1){\Bbb Z}$ and 
${\goth m}_B={\Bbb Z}+(2\sqrt{2}-2){\Bbb Z}$ can be written as
\displaymath
q_A(x,y)= 2x^2-4xy+6y^2,\qquad q_B(x,y)=2x^2-8xy+24y^2,
\enddisplaymath
respectively.  The  modules ${\goth m}_A, {\goth m}_B$ are not similar in the number field 
$K={\Bbb Q}(\sqrt{2})$,  since   their  determinants $\Delta({\goth m}_A)=8$ and 
$\Delta({\goth m}_B)=32$ are not equal.   Therefore,    matrices $A$ and $B$ are not similar
in the group  $GL(2,{\Bbb Z})$.        
}
\end{exm}
\begin{exm}
\textnormal{
Consider Anosov's automorphisms $f_A, f_B: T^2\to T^2$ given by matrices 
\displaymath
A=\left(\matrix{4 & 3\cr 5 & 4}\right)\qquad  \hbox{and}
\qquad B=\left(\matrix{4 & 15\cr 1 & 4}\right), 
\enddisplaymath
respectively.  The Alexander polynomials of $f_A$ and $f_B$ are identical   
$\Delta_A(t)=\Delta_B(t)= t^2-8t+1$;    yet  the automorphisms $f_A$ and $f_B$ are 
not conjugate.  Indeed, the Perron-Frobenius eigenvector of matrix $A$ is  $v_A=(1,  {1\over 3}\sqrt{15})$
while of the matrix $B$ is $v_B=(1, {1\over 15}\sqrt{15})$.  
The corresponding  modules are  ${\goth m}_A={\Bbb Z}+( {1\over 3}\sqrt{15}){\Bbb Z}$ and 
${\goth m}_B={\Bbb Z}+( {1\over 15}\sqrt{15}){\Bbb Z}$;  therefore 
\displaymath
q_A(x,y)= 2x^2+18y^2,\qquad q_B(x,y)=2x^2 +450y^2,
\enddisplaymath
respectively.  The  modules ${\goth m}_A, {\goth m}_B$ are not similar in the number field 
$K={\Bbb Q}(\sqrt{15})$,  since  the module  determinants $\Delta({\goth m}_A)=36$ and 
$\Delta({\goth m}_B)=900$ are not equal. Therefore,  matrices $A$ and $B$ are not similar
in the group  $GL(2,{\Bbb Z})$.     
}
\end{exm}

 \index{Geometrization Conjecture}

\vskip1cm\noindent
{\bf Guide to the literature.}
The topology of  surface  automorphisms   is the fundamental and the oldest part
of geometric topology;  it dates back to the works of  J.~Nielsen  [Nielsen 1927; 1929; 1932]  \cite{Nie1} 
and M.~Dehn  [Dehn 1938]  \cite{Deh1}.     W.~Thurston proved that that there are only
three types of such automorphisms:  they are either of finite order,  or of the Anosov
type (called {\it pseudo-Anosov}) or else a mixture of the two,  see e.g. [Thurston 1988]
\cite{Thu1};    the topological classification of pseudo-Anosov automorphisms 
is the next  problem after the {\it Geometrization Conjecture}  proved by G.~Perelman,
see [Thurston 1982]   \cite{Thu2}.   An excellent introduction to the subject are the books
[Fathi, Laudenbach \& Po\'enaru 1979]   \cite{FLP}  and  [Casson \& Bleiler  1988]   \cite{CaB}.  
The noncommutative invariants of pseudo-Anosov automorphisms were constructed
in \cite{Nik7}.

\section*{Exercises}

\begin{enumerate}

\item
Show that the  bracket on   $C^{\infty}(T^2)$ defined by the formula
\displaymath
\{f,g\}_{\theta}:=\theta\left({\partial f\over\partial x}  {\partial g\over\partial y}-
{\partial f\over\partial y} {\partial g\over\partial x}\right).
\enddisplaymath
is the Poisson bracket, i.e.  satisfies the identities $\{f,f\}_{\theta}=0$ and  
$\{f, \{g,h\}_{\theta}\}_{\theta} +\{h, \{f, g\}_{\theta}\}_{\theta}+\{g, \{h,f\}_{\theta}\}_{\theta}=0$.

\item
Prove that real multiplication is an invariant of the stable isomorphism 
(Morita equivalence) class of noncommutative torus ${\cal A}_{\theta}$.

\item
Prove that  complex  tori   ${\Bbb C}/({\Bbb Z}+{\Bbb Z}\tau)$ a  ${\Bbb C}/({\Bbb Z}+{\Bbb Z}\tau')$ 
are isomorphic   if and only if 
$\tau'={a\tau+b\over c\tau+d}$ for some matrix  $\left(\small\matrix{a & b\cr c & d}\right)  \in SL_2({\Bbb Z})$.
(Hint:  notice that $z\mapsto\alpha z$
is an invertible holomorphic map for each $\alpha\in {\Bbb C}-\{0\}$.)

\item
Prove that   the system of  relations
\displaymath
\left\{
\begin{array}{cc}
x_3x_1  &= e^{2\pi i\theta}x_1x_3,\\
x_1x_2 &= x_2x_1 = e,\\
x_3x_4  &= x_4x_3 = e.
\end{array}
\right.
\enddisplaymath
involved in the algebraic 
definition of noncommutative torus ${\cal A}_{\theta}$ is equivalent to the following
system of quadratic relations 
\displaymath
\left\{
\begin{array}{cc}
x_3x_1 &=  e^{2\pi i\theta}x_1x_3,\\
x_4x_2 &=  e^{2\pi i\theta}x_2x_4,\\
x_4x_1 &=  e^{-2\pi i\theta}x_1x_4,\\
x_3x_2 &=   e^{-2\pi i\theta}x_2x_3,\\
x_2x_1 &= x_1x_2=e,\\
x_4x_3 &= x_3x_4=e.
\end{array}
\right.
\enddisplaymath

 \index{complex modulus}
 \index{complex multiplication}
 \index{imaginary quadratic number}

\item
Prove that 
elliptic curve ${\cal E}_{\tau}$ has complex multiplication if and only if
the complex modulus $\tau\in {\Bbb Q}(\sqrt{-D})$ is an imaginary 
quadratic number.
(Hint:  Let $\alpha\in {\Bbb C}$ be such that $\alpha({\Bbb Z}+{\Bbb Z}\tau)\subseteq {\Bbb Z}+{\Bbb Z}\tau$.
That is there exist $m,n,r,s\in {\Bbb Z}$,  such that
\displaymath
\left\{
\begin{array}{cc}
\alpha &= m+n\tau,\\
\alpha\tau &= r+s\tau.
\end{array}
\right.
\enddisplaymath
One can divide the second equation by the first, so that one gets $\tau={r+s\tau\over m+n\tau}$.
Thus $n\tau^2+(m-s)\tau-r=0$;  in other words,  $\tau$ is an imaginary quadratic number. 
Conversely, if $\tau$ is the imaginary quadratic number, then it is easy to see that 
$End~({\Bbb C}/({\Bbb Z}+{\Bbb Z}\tau))$ is non-trivial.)

\item
Prove that the period $\overline{(a_1,a_2,\dots,a_2,a_1, 2a_0)}$  (and the arithmetic complexity)  is an invariant of the stable 
isomorphism  (Morita equivalence) class of the noncommutative torus    ${\cal A}_{RM}^{(D,f)}$.

\item
Prove that   the rational integer  $\Delta ({\goth m})$ is a numerical   invariant 
of Anosov's  automorphisms.  (Hint:  $\Delta ({\goth m})$ does not depend on the basis of module 
 ${\goth m}={\Bbb Z}v_A^{(1)}+{\Bbb Z}v_A^{(2)}$.)

\item
 Prove that  the matrices 
\displaymath
A=\left(\matrix{5 & 2\cr 2 & 1}\right)\qquad  \hbox{and}
\qquad B=\left(\matrix{5 & 1\cr 4 & 1}\right)
\enddisplaymath
are not similar by  using the {\it Gauss method},  i.e. the method  of continued fractions. 
(Hint:  Find the fixed points  $Ax=x$ and $Bx=x$, which gives us $x_A=1+\sqrt{2}$ and
$x_B={1+\sqrt{2}\over 2}$, respectively. Then one unfolds the fixed points into a periodic continued fraction,
which gives us $x_A=[2,2,2,\dots]$ and $x_B=[1,4,1,4,\dots]$. Since the period $(\overline{2})$ of $x_A$
differs from the period $(\overline{1,4})$ of $B$,  one concludes that matrices $A$ and $B$ belong to different
similarity  classes in $GL(2, {\Bbb Z})$.) 

 \index{Gauss method}

\item
 Repeat the exercise  for  matrices
\displaymath
A=\left(\matrix{4 & 3\cr 5 & 4}\right)\qquad  \hbox{and}
\qquad B=\left(\matrix{4 & 15\cr 1 & 4}\right).  
\enddisplaymath
 
 \end{enumerate}




\chapter{Categories and Functors}
Categories were  designed to extend  the ideas and methods of
algebraic topology and homological algebra to the rest of mathematics;  roughly,  
 categories  can be described as  a formal but amazingly helpful   ``calculus of arrows''  between certain 
``objects'' . 
Categories, functors and natural transformations  were introduced and studied  by 
[Eilenberg \& MacLane 1942] \cite{EiMa1}.  
 We refer the interested reader to the classic monograph [Cartan \& Eilenberg 1956] 
\cite{CE} for the basics of homological algebra;  for modern exposition,   see
[Gelfand \&  Manin 1994]  \cite{GM}.   
We briefly review the categories in Section 2.1,  functors in Section 2.2 and natural 
transformations in Section 2.3;   our exposition  follows the classical  book by  
[MacLane 1971]   \cite{ML}.  A set of exercises can be found at the end of the chapter.

 \index{category}

\section{Categories}
Category theory deals with certain diagrams representing sets 
and maps between the sets;  the maps can be composed with
each other, i.e. they  constitute a monoid with the operation 
of composition.  The sets are represented by points of the diagram;
the maps between the sets correspond to the arrows (i.e. directed 
edges) of the diagram.  

 \index{commutative diagram}

The diagram  is called {\it commutative} if for any choice of sets $A_i$
and maps $f_i: A_{i-1}\to A_i$ the resulting composition map from $A_0$
to $A_n$ is the same.  The examples of commutative diagrams are given
in Fig. 2.1.

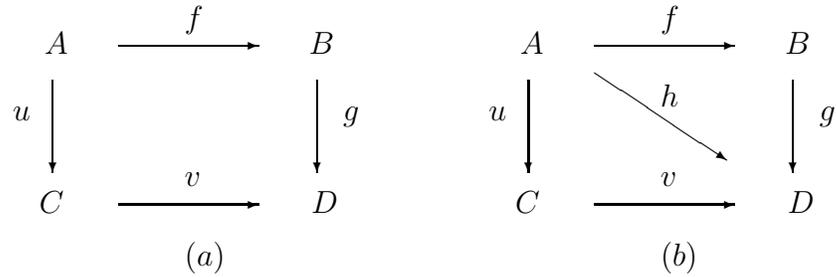
\begin{figure}[here]
\begin{picture}(500,110)(-40,-5)

\put(20,70){\vector(0,-1){35}}
\put(120,70){\vector(0,-1){35}}
\put(45,23){\vector(1,0){53}}
\put(45,83){\vector(1,0){53}}
\put(15,20){$C$}
\put(5,55){$u$}
\put(130,55){$g$}
\put(118,20){$D$}
\put(17,80){$A$}
\put(117,80){$B$}
\put(70,90){$f$}
\put(70,30){$v$}

\put(200,70){\vector(0,-1){35}}
\put(300,70){\vector(0,-1){35}}
\put(225,23){\vector(1,0){53}}
\put(225,83){\vector(1,0){53}}
\put(195,20){$C$}
\put(185,55){$u$}
\put(310,55){$g$}
\put(298,20){$D$}
\put(197,80){$A$}
\put(297,80){$B$}
\put(250,90){$f$}
\put(250,30){$v$}
\put(250,60){$h$}
\put(225,73){\vector(3,-2){50}}

\put(70,0){$(a)$}

\put(250,0){$(b)$}

\end{picture}

\caption{Two examples: (a) $vu=gf$ and (b) $h=vu=gf$.}
\end{figure}

\begin{exm}
\textnormal{
The commutative diagrams allow to prove that functions on the 
cartesian product of two sets $X$ and $Y$ are uniquely determined 
by such on the sets $X$ and $Y$.  Indeed,  consider the cartesian
product $X\times Y$ of two sets, consisting as usual of all ordered 
pairs $(x,y)$ of elements $x\in X$ and $y\in Y$.  The projections
$(x,y)\to x$ and $(x,y)\to y$ of the product define the functions
$p: X\times Y\to X$ and $q: X\times Y\to Y$.  Any function
$h: W\to X\times Y$  from a third set $W$ is uniquely determined by 
the composites $p\circ h$ and $q\circ h$.  Conversely,  given $W$ and
two functions $f$ and $g$,  there is a unique function $h$ which 
fits in the commutative diagram of Fig. 2.2. 
}
\end{exm}
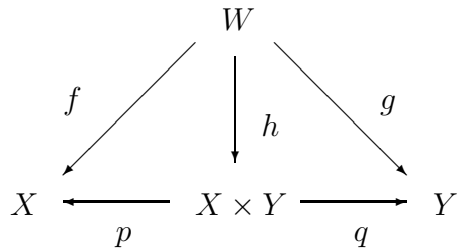
\begin{figure}[here]
\begin{picture}(500,110)(0,-5)

\put(205,90){\vector(0,-1){40}}
\put(230,35){\vector(1,0){40}}
\put(180,35){\vector(-1,0){40}}
\put(220,95){\vector(1,-1){50}}
\put(190,95){\vector(-1,-1){50}}

\put(200,100){$W$}
\put(190,30){$X\times Y$}
\put(120,30){$X$}
\put(280,30){$Y$}

\put(160,20){$p$}
\put(250,20){$q$}
\put(140,70){$f$}
\put(215,60){$h$}
\put(260,70){$g$}

\end{picture}
\caption{The universal property of  cartesian product.}
\end{figure}

\begin{rmk}
\textnormal{
The construction of cartesian product is called a {\it functor} because
it applies not only to sets but also to the functions between them. 
Indeed, two functions $f: X\to X'$ and $g: Y\to Y'$ define a function
$f\times g$ as their cartesian product according to the formula
\displaymath
f\times g: X\times Y\to X'\times Y', \quad (x,y)\mapsto (gx, fy).
\enddisplaymath
}
\end{rmk}
Below are examples illustrating the universal property of cartesian
product.  
\begin{exm}
\textnormal{
Let $X, Y$ and $W$ be the topological spaces and $p,q, f$ and $g$
be continuous maps between them.  Then the space $X\times Y$
is the product of topological spaces $X$ and $Y$.  Such a product
satisfies the universal property of cartesian product.  
}
\end{exm}
\begin{exm}
\textnormal{
Let $X, Y$ and $W$ be the abelian groups and $p,q, f$ and $g$
be the homomorphisms between the groups. 
Consider the direct sum (product)  $X\oplus Y$ of the abelian
groups $X$ and $Y$.  In view of the canonical embeddings 
$X\to X\oplus Y$ and $Y\to X\oplus Y$,  the direct sum satisfies
the universal property of cartesian product.  
}
\end{exm}
\begin{rmk}
\textnormal{
The direct sum of non-abelain groups does not satisfy the
universal property of cartesian products.  Indeed, suppose that
$W$ has two subgroups isomorphic to $X$ and $Y$, respectively, 
but whose elements do not commute with one another;  let $f$ and
$g$ be isomorphisms of $X$ and $Y$ with these subgroups.  
Since elements $x\in X$ and $y\in Y$ commute in the direct sum
$X\oplus Y$,  the diagram for cartesian product will not be commutative
for any homomorphism $h$.  
}
\end{rmk}
\begin{exm}
\textnormal{
Let $X, Y$ and $W$ be any groups and $X\ast Y$ be the
free product of groups $X$ and $Y$. It is left as an exercise 
to the reader to prove that the product   $X\ast Y$ satisfies
the universal property of cartesian product. 
}
\end{exm}
The above examples involved sets and certain maps between the sets.
One never needed to consider what kind of elements our sets were made of,
or how these elements transformed under the maps.  The only thing we needed
was that maps can be composed with one another and that such maps can be
arranged into a commutative diagram.  Such a standpoint can be axiomatized  
in the notion of a category.  
\begin{dfn}
A category ${\cal C}$ consists of the following data: 

\medskip
(i)  a set $Ob ~{\cal C}$ whose elements are called the objects of ${\cal C}$;

\smallskip
(ii)  for any $A, B\in Ob~{\cal C}$ a set $H(A, B)$ whose elements are called 
the morphisms of ${\cal C}$ from $A$ to $B$;  

\smallskip
(iii) for any $A,B,C\in Ob~{\cal C}$ and any $f\in H(A,B)$ and $g\in H(B,C)$
a morphism $h\in H(A,C)$ is defined and which is called the composite $g\circ f$
of $g$ and $f$;

\smallskip
(iv)  for any $A\in Ob~{\cal C}$ a morphism $1_A\in H(A,A)$ is defined 
and called the identity morphism,  so that $f\circ 1_A=1_B\circ f=f$ for
each $f\in H(A,B)$;

\smallskip
(v) morphisms are associative,  i.e. $h\circ (g\circ f)=(h\circ g)\circ f$
for all $f\in H(A,B), g\in H(B,C)$ and $h\in H(C,D)$.         
\end{dfn}
 \index{arrows}
 \index{objects}
 \index{morphism of category}
The morphism of in category ${\cal C}$ are often called {\it arrows},
since they are represented by such in  the commutative diagrams.
Notice that since the objects in a category correspond exactly to its 
identity arrows, it is technically possible to dispense altogether with
the objects and deal only with the arrows. 
The data for an  {\it arrows-only} category ${\cal C}$ consists of arrows,
certain ordered pairs $(g,f)$ called the composable pairs of arrows, and
an operation assigning to each composable pair $(g,f)$ an arrow $g\circ f$
called the composite. With these data one defines an identity of ${\cal C}$
to be an arrow $u$ such that $f\circ u=f$ whenever the composite $f\circ u$
is defined and $u\circ g=g$ whenever $u\circ g$ is defined.  The data is required 
to satisfy the following three axioms:  (i)  the composite $(k\circ g)\circ f$ is defined
if and only if $k\circ (g\circ f)$ is defined and we write it as $k\circ g\circ f$, 
(ii)  the triple composite $k\circ  g\circ f$ is defined whenever both composites
$k\circ g$ and $g\circ f$ are defined and (iii) for each arrow $g$ of ${\cal C}$
there exists identity arrows $u$ and $u'$ of ${\cal C}$ such that $u' \circ g$ 
and $g\circ u$ are defined.    
\begin{exm}
\textnormal{
An {\it empty} category {\bf 0} contains no objects and no arrows.
}
\end{exm}
\begin{exm}
\textnormal{
A category is {\it discrete} when every arrow is an identity arrow. 
Every set $X$ is the set of objects of a discrete category and every
discrete category is determined by its set of objects. Thus, discrete
categories are sets.  
}
\end{exm}
\begin{exm}
\textnormal{
Consider the category {\bf Ind} whose objects are arbitrary subsets of 
a given set $X$ and whose arrows are the inclusions maps between them
so that $H(A,B)$ is either empty or consists of a single element. 
}
\end{exm}
\begin{exm}
\textnormal{
The category {\bf Top} whose objects are topological spaces and whose
arrows are continuous maps between them. 
}
\end{exm}
\begin{exm}
\textnormal{
The category {\bf Grp} whose objects are all groups and whose 
arrows are homomorphisms between the groups.
}
\end{exm}
\begin{exm}
\textnormal{
The category {\bf Ab} whose objects are additive abelian groups and
arrows are homomorphisms of the abelian groups. 
}
\end{exm}
\begin{exm}
\textnormal{
The category {\bf Rng} whose objects are rings and arrows are homomorphisms
between the rings preserving the unit. 
}
\end{exm}
\begin{exm}
\textnormal{
The category {\bf CRng} whose objects are commutative rings and whose
arrows are homomorphisms between the rings.  
}
\end{exm}
\begin{exm}
\textnormal{
The category {\bf Mod}$_{R}$ whose objects are modules over a given
ring $R$,  and whose arrows are homomorphisms between them.
The category    {\bf Mod}$_{\Bbb Z}$  of abelian  groups was denoted 
by {\bf Ab}. 
}
\end{exm}

\vskip1cm\noindent
{\bf Guide to the literature.}
Categories were introduced and studied  by [Eilenberg \& MacLane 1942] \cite{EiMa1}.  
We refer the interested reader to the classic monograph [Cartan \& Eilenberg 1956] 
\cite{CE} for the basics of homological algebra;  for modern exposition,   see
[Gelfand \&  Manin 1994]  \cite{GM}.   Our exposition  of categories follows the 
classical  book by  [MacLane 1971]   \cite{ML}.

 \index{functor}

\section{Functors}
A powerful tool in mathematics is furnished by the so-called invariant
or natural constructions;  a formalization of such constructions leads
to the notion of a {\it functor}.  Roughly speaking,  the functor is a morphism
of categories.  
\begin{dfn}
A covariant functor from a category ${\cal C}$ to a category ${\cal D}$ 
consists of two maps (denoted by the same letter $F$):

\medskip
(i)  $F: Ob~{\cal C}\to Ob~{\cal D}$ and 

\smallskip
(ii)  $F:  H(A,B)\to H(F(A), F(B))$ for all $A,B\in Ob~ {\cal C}$,

\medskip\noindent
which satisfy the following conditions:

\medskip
(i) $F(1_A)=1_{F(A)}$ for all $A\in  Ob~{\cal C}$;

\smallskip
(ii)  $F(f\circ g)=F(f)\circ F(g)$, whenever $f\circ g$ is defined in ${\cal C}$.   
\end{dfn}
\begin{rmk}
\textnormal{
A {\it contravariant functor}  is also given by a map $F: Ob~{\cal C}\to Ob~{\cal D}$ 
but it defines a reverse order map:
\displaymath
F: H(A,B)\to H(F(B), F(A))
\enddisplaymath
for all $A,B\in Ob~{\cal C}$;  the latter  must satisfy the reverse order conditions:
\displaymath
F(1_A)=1_{F(A)}   \quad\hbox{and} \quad F(f\circ g)=F(g)\circ F(f).
\enddisplaymath
 }
\end{rmk}
 \index{contravariant funcor}
 \index{covariant functor}

\bigskip\noindent
A functor,  like a category,  can be described in the ``arrows-only'' mode;  
roughly speaking,  $F$ is a function from arrows $f$ of ${\cal C}$ to arrows 
$Ff$ of ${\cal D}$ carrying each identity of ${\cal C}$ to an identity of ${\cal D}$
and each composable pair $(g,f)$ in ${\cal C}$ to a composable pair 
$(Fg, Ff)$ in ${\cal D}$ with $Fg\circ Ff=F(g\circ f)$. 

\begin{exm}
\textnormal{
The singular $n$-dimensional homology assigns to each topological space $X$
an abelian group $H_n(X)$ called the $n$-th homology group. 
Each continuous map $f: X\to Y$ between the topological spaces $X$ and $Y$
gives rise to a homomorphism $f_*: H_n(X)\to H_n(Y)$ between the corresponding
$n$-th homology groups. Therefore, the $n$-th singular homology is a covariant 
functor $H_n:$ {\bf Top} $\to$ {\bf Ab}  between the category of topological spaces
and the category of abelian groups. 
}
\end{exm}
\begin{exm}
\textnormal{
The $n$-th homotopy groups $\pi_n(X)$ of the topological space $X$
can be regarded as functors.  Since $\pi_n(X)$ depend on the choice 
of a base point in $X$, they are functors {\bf Top}$_*\to$ {\bf Grp} 
between the category of topological spaces with a distinguished point
and the category of groups. Note that the groups are abelian unless
$n=1$.  
}
\end{exm}
\begin{exm}
\textnormal{
Let $X$ be a set, e.g. a Hausdorff topological space or the homogeneous 
space of a Lie group $G$.  Denote by ${\cal F}(X, {\Bbb C})$ the space
of continuous complex-valued functions on $X$.  Since any map $f: X\to Y$
 takes the function $\varphi\in {\cal F}(X, {\Bbb C})$ into a function
$\varphi'\in {\cal F}(X, {\Bbb C})$, it follows that ${\cal F}(X, {\Bbb C})$
is a contravariant functor  from category of sets to the category of vector
spaces (usually infinite-dimensional).  The functor maps any transformation
of $X$ into an invertible linear operator on the space ${\cal F}(X, {\Bbb C})$.
In particular, if $X=G$ is a Lie group, one can consider the action of $G$ on itself
by the left translations;  thus one gets the regular representation of $G$ by the
linear operators on ${\cal F}(X, {\Bbb C})$.  
}
\end{exm}
\begin{exm}\label{exm2.3.4}
\textnormal{
Let $R$ be a commutative ring.  Consider the multiplicative
group $GL_n(R)$ of all non-singular $n\times n$ matrices with
entries in $R$.  Because each homomorphism $f: R\to R'$
produces in the evident way a group homomorphism 
$GL_n(R)\to GL_n(R')$, one gets a functor 
$GL_n:$ {\bf CRng} $\to$ {\bf Grp} from the category of 
commutative rings to the category of groups. 
}
\end{exm}
\begin{rmk}
\textnormal{
Functors can be composed.  Namely,  given functors
\displaymath
{\cal C}
\buildrel\rm
T
\over\longrightarrow 
{\cal B}
\buildrel\rm
S
\over\longrightarrow
{\cal A}
\enddisplaymath
between categories ${\cal A}, {\cal B}$ and ${\cal C}$,
the composite functions 
\displaymath
c\mapsto S(Tc), \qquad f\mapsto S(Tf)
\enddisplaymath
on the objects $c\in Ob~{\cal C}$ and arrows of ${\cal C}$
define a functor $S\circ T: {\cal C}\to {\cal A}$ called the 
{\it composite} of $S$ with $T$.  The composition is associative.
For each category ${\cal B}$ there is an identity functor 
$I_{\cal B}: {\cal B}\to {\cal B}$ which acts as an identity for this
composition. Thus one gets a category of all categories endowed 
with the composition arrow. 
}    
\end{rmk}
\begin{dfn}
A functor $F$ is called forgetful if $F$ forgets some or all of the
structure of an algebraic object. 
\end{dfn}
\begin{exm}
\textnormal{
The functor $F:$ {\bf Rng} $\to$ {\bf Ab} assigns to each ring $R$
the additive abelian group of $R$; such a functor is a forgetful
functor, because $F$ forgets the multiplicative structure of the ring $R$.
}
\end{exm}
\begin{dfn}
A functor $F: {\cal C}\to {\cal B}$ is called full when to every pair of 
objects $C,C'\in Ob~{\cal C}$ and every arrow $g: FC\to FC'$ of ${\cal B}$,
there is an arrow $f: C\to C'$ of ${\cal C}$ with $g=Ff$. 
\end{dfn}
 \index{faithful functor}
\begin{dfn}
A functor $F: {\cal C}\to {\cal B}$ is called faithful (or an embedding)
when to every pair $C, C'\in Ob~{\cal C}$ and to every pair 
$f_1,f_2\in H(C, C')$ of parallel arrows of ${\cal C}$ the equality
$Ff_1=Ff_2: F(C)\to F(C')$ implies $f_1=f_2$. 
\end{dfn}
\begin{dfn}
An isomorphism $F: {\cal C}\to {\cal B}$ of categories is a functor $F$ which 
is a bijection both on the objects and arrows of the respective categories. 
\end{dfn}
 \index{isomorphic categories}
\begin{rmk}
\textnormal{
Every full and faithful functor $F: {\cal C}\to {\cal B}$ is an
isomorphism of categories ${\cal C}$ and ${\cal B}$ provided
there are no objects of ${\cal B}$ not in the image of $F$.  
}    
\end{rmk}

\vskip1cm\noindent
{\bf Guide to the literature.}
Functors were introduced and studied  by [Eilenberg \& MacLane 1942] \cite{EiMa1}.  
We refer the interested reader to the classic monograph [Cartan \& Eilenberg 1956] 
\cite{CE} for the basics of homological algebra;  for modern exposition,   see
[Gelfand \&  Manin 1994]  \cite{GM}.   Our exposition  of categories follows the 
classical  book by  [MacLane 1971]   \cite{ML}.

\section{Natural transformations}
Let $S, T: {\cal C}\to {\cal B}$ be two functors between the same 
categories ${\cal C}$ and ${\cal B}$.  Roughly speaking, a natural
transformation $\tau$ is a set of arrows of ${\cal B}$,  which translates
the ``picture $S$'' of category ${\cal C}$ to the ``picture $T$'' of ${\cal C}$.  In 
other words, a natural transformation is a morphism of functors.  
\begin{dfn}
Given two functors $S, T: {\cal C}\to {\cal B}$, a natural transformation 
$\tau: S \buildrel\rm \bullet \over\rightarrow T$ is a function which assigns
to each object $c\in {\cal C}$ an arrow $\tau_c: Sc\to Tc$ of ${\cal B}$ 
in such a way that every arrow $f: c\to c'$ in ${\cal C}$ yields the commutative
diagram shown in Fig. 2.3.  
\end{dfn}
\begin{figure}[here]
\begin{picture}(300,100)(-130,0)

\put(25,70){\vector(0,-1){35}}
\put(120,70){\vector(0,-1){35}}
\put(45,23){\vector(1,0){53}}
\put(45,83){\vector(1,0){53}}
\put(15,20){$Sc'$}
\put(0,55){$Sf$}
\put(130,55){$Tf$}
\put(118,20){$Tc'$}
\put(17,80){$Sc$}
\put(117,80){$Tc$}
\put(70,90){$\tau_c$}
\put(70,30){$\tau_{c'}$}

\end{picture}
\caption{The natural transformation.}
\end{figure}
\begin{rmk}
\textnormal{
A natural transformation $\tau$ with every component $\tau_c$ invertible
in category ${\cal B}$ is called a {\it natural isomorphism} of functors 
$S$ and $T$;  such an isomorphism is denoted by $\tau: S\cong T$.  
 }    
\end{rmk}
\begin{exm}
\textnormal{
Let $R$ be a commutative ring. Denote by $R^*$ the group of units,
i.e. the invertible elements of $R$.  Let $S: R\to R^*$ be a functor 
{\bf CRng} $\to$ {\bf Grp} which is the embedding of $R^*$ into $R$. 
Recall that there exists another functor $GL_n:$ {\bf CRng} $\to$ {\bf Grp}, 
which assigns to $R$ the group $GL_n(R)$ consisting of the $n\times n$
matrices with the entries in the ring $R$;  we shall put $T=GL_n$.  
A natural transformation $\tau: S \buildrel\rm \bullet \over\rightarrow T$
is defined by the function $det_R: M_n(R)\to R^*$ in category {\bf Grp}, 
which assigns to each matrix $M_n(R)$ its determinant;  since $M_n(R)$
is a non-singular matrix, the $det_R$ takes valued in $R^*$.  Thus one gets 
the commutative diagram shown in Fig. 2.4.  
}
\end{exm}
\begin{figure}[here]
\begin{picture}(300,100)(-30,0)

\put(120,93){\vector(1,0){60}}
\put(110,80){\vector(1,-1){40}}
\put(160,40){\vector(1,1){40}}

\put(100,90){$R$}
\put(200,90){$R^*$}
\put(140,20){$GL_n(R)$}

\put(150,99){$S$}
\put(110,50){$T$}
\put(190,50){$det_R$}

\end{picture}
\caption{Construction of the natural transformation $det_R$.}
\end{figure}

\vskip1cm\noindent
{\bf Guide to the literature.}
Natural transformations  were introduced and studied  by [Eilenberg \& MacLane 1942] \cite{EiMa1}.  
We refer the interested reader to the classic monograph [Cartan \& Eilenberg 1956] 
\cite{CE} for the basics of homological algebra;  for modern exposition,   see
[Gelfand \&  Manin 1994]  \cite{GM}.   Our exposition  of categories follows the 
classical  book by  [MacLane 1971]   \cite{ML}.

\section*{Exercises}

\begin{enumerate}

\item
Prove that the map from  the set of all integral domains to their 
quotient fields is a functor.

\item Prove that correspondence between the Lie groups and their
Lie algebras is a functor.

\item
Prove that the correspondence between affine algebraic varieties and their 
coordinate rings is a functor.  (Hint:  For missing definitions, facts and a proof,  
see e.g.  [Hartshorne 1977], \cite{H1}, Chapter 1.)

\item
Show that the category of all finite-dimensional vector spaces over 
a field $F$ is equivalent to the category of matrices over $F$.

 \end{enumerate}




\chapter{$C^*$-Algebras}
One can think of the $C^*$-algebras  as a deep  generalization of the
field  of  complex numbers ${\Bbb C}$;   the fundamental role of  ${\Bbb C}$ 
in   algebra, analysis and geometry matches the one of the $C^*$-algebras in 
noncommutative geometry.   Operator algebras grew from the problems of theory 
of unitary group representations and quantum mechanics,   see the 
seminal paper by [Murray \& von Neumann 1936]   \cite{MuNeu1}. 
The $C^*$-algebras  {\it per se}  were introduced and studied in [Gelfand \& Naimark 1943]
\cite{GelNai1}.  There exists plenty of excellent textbooks  on the subject, 
see [Arveson 1976]  \cite{A},   [Blackadar 1986]  \cite{B},   [Davidson ]  \cite{D},   [Fillmore 1996] \cite{F},   [Murphy 1990]  \cite{M},     
[R\o rdam,  Larsen \& Laustsen  2000]  \cite{RLL},  [Wegge-Olsen 1993]  \cite{W1} and others.
Section 3.1 contains basic definitions related to the $C^*$-algebras. The important class 
of $C^*$-algebras called {\it crossed products} is reviewed in Section 3.2.  
The {\it $K$-theory}  of $C^*$-algebras revolutionized the subject and its brief review
can be found in Section 3.3.  The {\it $n$-dimensional noncommutative tori}
are introduced in Section 3.4;  along with the functors,  such tori  
are used throughout the book.   The {\it $AF$-algebras}    are
an  important class of the non-trivial $C^*$-algebras,  which can be classified 
and studied;  their brief review is given in Section 3.5.  The {\it UHF-algebras}
and the {\it Cuntz-Krieger algebras}  are introduced in Sections 3.6 and 3.7,
respectively.  Our choice of the $C^*$-algebras is prompted  by their applications
in Part II.

 \index{$C^*$-algebra}

\section{Basic definitions}
\begin{dfn}
A $C^*$-algebra $A$ is an algebra over ${\Bbb C}$ with a norm $a\mapsto ||a||$
and an involution $a\mapsto a^*$ on $a\in A$,  such that  $A$ is complete with
respect to the norm and such that for every $a,b\in A$:

\medskip
(i)  $||ab||\le ||a|| ||b||$;

\smallskip
(ii) $||a^*a||=||a||^2$.

\bigskip\noindent
The $C^*$-algebra $A$ is called unital if it has a multiplicative identity.  
\end{dfn}
\begin{rmk}
\textnormal{
The $C^*$-algebras is an abstraction of the algebras generated by  bounded linear 
operators acting on a Hilbert space ${\cal H}$.  The space ${\cal H}$  comes 
with a scalar product and,  therefore,  one gets   the norm $||\bullet||$ and 
involution $a\mapsto a^*$  on the operators  $a\in A$.   Notice that algebra $A$
is essentially non-commutative.  
}
\end{rmk}
 \index{projection}
 \index{unitary}
 \index{partial isometry}
\begin{dfn}
An element $p$ in a $C^*$-algebra $A$ is called projection, if
$p=p^*=p^2$;  two projections $p$ and $q$ are said to be orthogonal if $pq=0$
and we write $p\perp q$ in this case.  The element $u$ in a unital $C^*$-algebra $A$ 
is said to be unitary if $uu^*=u^*u=1$.  An element $s$ in $A$ is called a partial
isometry when $s^*s$ is a projection.  
\end{dfn}
\begin{exm}
\textnormal{
The field of complex numbers ${\Bbb C}$ is a $C^*$-algebra with involution 
given by the complex conjugation $z\mapsto \bar z$ of complex numbers $z\in {\Bbb C}$.  
}
\end{exm}
\begin{exm}
\textnormal{
The matrix algebra $M_n({\Bbb C})$ of all $n\times n$ matrices $(z_{ij})$ 
with complex entries $z_{ij}$  endowed with the usual matrix norm and 
involution $(z_{ij})^*=(\bar z_{ij})$  is a $C^*$-algebra. It is a finite-dimensional 
$C^*$-algebra,  since it is isomorphic to the algebra of bounded linear operators
on the finite-dimensional Hilbert space ${\Bbb C}^n$.  The $C^*$-algebra
$M_n({\Bbb C})$ is non-commutative unless $n=1$. 
}
\end{exm}
\begin{exm}
\textnormal{
The matrix algebra $M_n({\Bbb C})$ of all $n\times n$ matrices $(z_{ij})$ 
with complex entries $z_{ij}$  endowed with the usual matrix norm and 
involution $(z_{ij})^*=(\bar z_{ij})$  is a $C^*$-algebra. It is a finite-dimensional 
$C^*$-algebra,  since it is isomorphic to the algebra of bounded linear operators
on the finite-dimensional Hilbert space ${\Bbb C}^n$.  The $C^*$-algebra
$M_n({\Bbb C})$ is non-commutative unless $n=1$. 
}
\end{exm}
\begin{exm}
\textnormal{
If ${\cal H}$ is a Hilbert space,  then the algebra of all compact linear operators on 
${\cal H}$ is a $C^*$-algebra denoted by ${\cal K}$;  the $C^*$-algebra ${\cal K}$
can be viewed as the limit of finite-dimensional $C^*$-algebras $M_n({\Bbb C})$
when $n$ tends to infinity. 
}
\end{exm}
\begin{exm}
\textnormal{
If ${\cal H}$ is a Hilbert space, then the algebra of all bounded linear
operators on ${\cal H}$ is a $C^*$-algebra denoted by $B({\cal H})$; 
it is known that every $C^*$-algebra can be thought of as a $C^*$-subalgebra
of some $B({\cal H})$ (Gelfand-Naimark Theorem). 
}
\end{exm}
\begin{exm}
\textnormal{
Let $A$ be a $C^*$-algebra.  Consider an algebra $M_n(A)$
consisting of all $n\times n$ matrices $(a_{ij})$ with entries 
$a_{ij}\in A$;  the involution on $M_n(A)$ is defined  by the 
formula $(a_{ij})^*=(a_{ji}^*)$.  There exists a unique norm on
$M_n(A)$,  such that $M_n(A)$ is a $C^*$-algebra. 
}
\end{exm}
\begin{rmk}
\textnormal{
The $C^*$-algebra $M_n(A)$ is isomorphic to  a tensor product
$C^*$-algebra $M_n({\Bbb C})\otimes A$;   since the algebra
${\cal K}=\lim_{n\to\infty} M_n({\Bbb C})$,  one can talk about the
tensor product $C^*$-algebra $A\otimes {\cal K}$.  
}
\end{rmk}
 \index{stable isomorphism}
\begin{dfn}
The $C^*$-algebras $A$ and $A'$ are said to be stably isomorphic,
if 
\displaymath
A\otimes {\cal K}\cong A'\otimes {\cal K}. 
\enddisplaymath
\end{dfn}
We proceed with further examples of the $C^*$-algebras in the
next sections;  below we shall focus on the special case of 
commutative $C^*$-algebras.

Let $A$ be a Banach algebra,  i.e. a complete norm algebra over ${\Bbb C}$. 
When $A$ is commutative, a non-zero homomorphism $\tau: A\to {\Bbb C}$
is called a {\it character} of $A$;  we shall denote by $\Omega(A)$ the space
of all characters on $A$ endowed with the weak topology.  Define  a function
$\hat a: \Omega(A)\to {\Bbb C}$ by the formula $\tau\mapsto \tau(a)$.  
\begin{thm}
{\bf (Gelfand)}
 If $A$ is a commutative Banach algebra, then:
 
 \medskip
 (i)  $\Omega(A)$ is a locally compact Hausdorff topological space,
 which is compact whenever $A$ is unital;
 
 \smallskip
 (ii)  the map $A\mapsto C_0(\Omega(A))$ given by the formula 
 $a\mapsto \hat a$ is a norm-decreasing homomorphism.
 \end{thm}
\begin{cor}
Each commutative $C^*$-algebra is isomorphic to the algebra  $C_0(X)$
of continuous complex-valued functions vanishing at the infinity of a locally
compact Hausdorff topological space $X$.   
\end{cor}
For non-commutative $C^*$-algebras no general classification is known;
however, the Gelfand-Naimark-Segal (GNS-) construction implies that
each $C^*$-algebra has a concrete  realization  as a $C^*$-subalgebra of the algebra
$B({\cal H})$ for some Hilbert space ${\cal H}$.  Below is a brief account 
of the GNS-construction. 

 \index{GNS-construction}

A {\it representation} of a $C^*$-algebra $A$ is a pair $({\cal H}, \varphi)$,
where ${\cal H}$ is a Hilbert space and $\varphi: A\to B({\cal H})$ is 
a $\ast$-homomorphism.   The representation $({\cal H}, \varphi)$ is  said 
to be {\it faithful} if $\varphi$ is injective.  One can sum up the representations
as follows.  If  $({\cal H}_{\lambda}, \varphi_{\lambda})_{\lambda\in\Lambda}$
is a family of representation of the $C^*$-algebra $A$,  a {\it direct sum} 
is  the representation  $({\cal H}, \varphi)$ got by setting 
${\cal H}=\oplus_{\lambda\in\Lambda} {\cal H}_{\lambda}$ and 
$\varphi(a)((x_{\lambda})_{\lambda})=(\varphi_{\lambda}(a)(x_{\lambda}))_{\lambda}$
for all $a\in A$ and $(x_{\lambda})_{\lambda}\in {\cal H}$.

For each positive linear functional $\tau: A\to {\cal C}$ on the $C^*$-algebra $A$,  
there is  associated a representation of $A$.  Indeed,  let 
$N_{\tau}=\{a\in A ~|~ \tau(a^*a)=0\}$ be a subset of $A$;
it is easy to verify that $N_{\tau}$ is a closed left ideal of $A$.   
We shall define a map $(A/N_{\tau}, A/N_{\tau})\to {\Bbb C}$
by the formula:
\displaymath
(a+N_{\tau},  b+N_{\tau})\mapsto \tau(b^*a);
\enddisplaymath
such a map is a well-defined inner product on the inner product space $A/N_{\tau}$.
We denote by ${\cal H}_{\tau}$ the Hilbert space completion of $A/N_{\tau}$. 
 If $a\in A$, then one can define an operator $\varphi(a)\in B({\cal H}_{\tau})$
 by setting
\displaymath
\varphi(a)(b+N_{\tau})=ab+N_{\tau}.
\enddisplaymath
Thus one gets a $\ast$-homomorphism
\displaymath
\varphi_{\tau}: A\to B({\cal H}_{\tau}),
\enddisplaymath
given by the formula $a\mapsto \varphi_{\tau}(a)$;  
the  representation $({\cal H}_{\tau}, \varphi_{\tau})$ of $A$
is called the {\it Gelfand-Naimark-Segal represetation}
associated to the positive linear functional $\tau$.  A {\it universal}
representation of $A$ is defined as the direct sum of 
representations   $({\cal H}_{\tau}, \varphi_{\tau})$ as 
$\tau$ ranges over the space $S(A)$ of all positive linear functionals
on  $A$.  
\begin{thm}
{\bf (Gelfand-Naimark)}
If $A$ is a $C^*$-algebra, then it has a faithful representation
in the space $B({\cal H})$ of bounded linear operators on a Hilbert   
 space ${\cal H}$;  such a  representation coincides with the universal 
 representation of $A$. 
 \end{thm}

 \index{Gelfand-Naimark Theorem}

\vskip1cm\noindent
{\bf Guide to the literature.}
The definition of an abstract $C^*$-algebra is due to   [Gelfand \& Naimark 1943]  \cite{GelNai1}.  
For other  facts  regarding  the subject,  the reader is encouraged to consult  the textbooks by 
[Arveson 1976]  \cite{A},  [Davidson ]  \cite{D},   [Fillmore 1996] \cite{F},   [Murphy 1990]  \cite{M},     
[R\o rdam,  Larsen \& Laustsen  2000]  \cite{RLL} and  [Wegge-Olsen 1993]  \cite{W1}.

 \index{crossed product}

\section{Crossed products}
The crossed product is a popular construction of new operator algebras 
from the given one;   the construction dates back to the works of Murray and
von Neimann.  
To give an idea, let $A$ be a $C^*$-algebra and $G$ a locally compact group.  
We shall consider a continuous homomorphism $\alpha$ from $G$ to the 
group $Aut~A$ of $\ast$-automorphisms of $A$ endowed with the topology
of pointwise norm-convergence.  
Roughly speaking,  the idea of the crossed
product construction  is to embed $A$ into a larger $C^*$-algebra 
in which the automorphism becomes the inner automorphism.  
We shall pass to a detailed description of the crossed product 
construction.  

 \index{covariant representation}

A {\it covariant representation} of the 
triple $(A,G,\alpha)$ is a pair of representations $(\pi,\rho)$ of $A$ and $G$ on the
same Hilbert space ${\cal H}$, such that
\displaymath
\rho(g)\pi(a)\rho(g)^*=\pi(\alpha_g(a))
\enddisplaymath
for all $a\in A$ and $g\in G$. Each covariant representation of 
$(A,G,\alpha)$ gives rise to a convolution algebra $C(G,A)$
of continuous functions from $G$ to $A$;  the completion of
$C(G,A)$  in the norm topology is a $C^*$-algebra $A\rtimes_{\alpha} G$
called a {\it crossed product} of $A$ by $G$. If $\alpha$ is a single
automorphism of $A$, one gets an action of ${\Bbb Z}$ on $A$;  
the crossed product in this case is called simply the crossed product
of $A$ by ${\alpha}$.  
\begin{exm}
\textnormal{
Let $A\cong {\Bbb C}$ the field of complex numbers.  The non-degenerate 
representations  $\rho$ of $C(G)$ correspond  to unitary representations 
$U$ of group $G$ via the formula
\displaymath
\rho(f)=\int_G f(s) U_s ds,
\enddisplaymath
 for all $f\in C(G)$.  The completion of $C(G)$ in the operator norm is 
 called the {\it group $C^*$-algebra} and denoted by $C^*(G)$.
 Thus ${\Bbb C}\rtimes_{Id} G\cong C^*(G)$. 
}
\end{exm}
\begin{exm}
\textnormal{
Let $A\cong C(S^1)$ be the commutative $C^*$-algebra of continuous 
complex-valued functions on the unit circle $S^1$.  Let ${\Bbb Z}/n {\Bbb Z}$
be the cyclic group of order $n$ acting on $S^1$ by rotation through the 
rational angle ${2\pi\over n}$.  In this case it is known that
\displaymath
C(S^1)\rtimes {\Bbb Z}/n{\Bbb Z}\cong M_n(C(S^1)),
\enddisplaymath
where  $M_n(C(S^1))$ is the $C^*$-algebra of all $n\times n$ matrices
with the entries in the $C^*$-algebra $C(S^1)$. 
}
\end{exm}
\begin{exm}
\textnormal{
Let $A\cong C(S^1)$ be the commutative $C^*$-algebra of continuous 
complex-valued functions on the unit circle $S^1$. 
Let $\theta\in {\Bbb R}-{\Bbb Q}$ be an irrational number and 
consider an automorphism of $C(S^1)$ generated by the rotation
of $S^1$ by the irrational angle $2\pi\theta$.  In this case $G\cong {\Bbb Z}$
and it is known that 
\displaymath
C(S^1)\rtimes {\Bbb Z}\cong {\cal A}_{\theta},
\enddisplaymath
where ${\cal A}_{\theta}$ is the {\it irrational rotation $C^*$-algebra}
also called the {\it noncommutative torus}.  Such algebras are crucial
examples of the noncommutative spaces and  extremely important
for the rest of this book. 
}
\end{exm}
One can recover the $C^*$-algebra $A$ from its crossed product by the action
of a locally compact abelian  group $G$;  the corresponding construction 
is known as the {\it Landstadt-Takai duality}.  Let us  give a brief description of 
this important construction. 

 \index{Landstadt-Takai duality}

Let $(A,G,\alpha)$ be a $C^*$-dynamical system with $G$ locally compact
abelian group; let $\hat G$ be the dual of $G$. For each $\gamma\in \hat G$,
one can define a map $\hat a_{\gamma}: C(G,A)\to C(G,A)$
given by the formula:
\displaymath
\hat a_{\gamma}(f)(s)=\bar\gamma(s)f(s), \qquad\forall s\in G.
 \enddisplaymath
In fact, $\hat a_{\gamma}$ is a $\ast$-homomorphism, since it
respects the convolution product and involution on $C_c(G,A)$
[Williams  2007]  \cite{W}.  Because the crossed product $A\rtimes_{\alpha}G$
is the closure of $C(G,A)$, one gets an extension of $\hat a_{\gamma}$
to an element of $Aut~(A\rtimes_{\alpha}G)$ and, therefore, a 
homomorphism:
\displaymath
\hat\alpha: \hat G\to Aut~(A\rtimes_{\alpha}G).
 \enddisplaymath
\begin{thm}
{\bf (Landstadt-Takai duality)}
If $G$ is a locally compact abelian group, then
\displaymath
A\rtimes_{\alpha} G)\rtimes_{\hat\alpha}\hat G\cong A\otimes {\cal K}(L^2(G)),
\enddisplaymath
where ${\cal K}(L^2(G))$ is the algebra of compact operators on the
Hilbert space $L^2(G)$. 
 \end{thm}

\vskip1cm\noindent
{\bf Guide to the literature.}
The crossed product construction  (for von Neumann algebras)  was    
already popular in the classical paper by [Murray \& von Neumann 1936]   \cite{MuNeu1}.
The detailed account of the crossed products of $C^*$-algebras 
can be found in the monograph by [Williams 2007]  \cite{W}.

 \index{K-theory}

\section{K-theory of $C^*$-algebras }
Like the $n$-dimensional manifolds or discrete groups, 
the non-commutative $C^*$-algebras are hardly classifiable;
yet there exist certain covariant functors  on the category {\bf C*-Alg} 
with values in the category of abelian groups.   Typically, such functors
are not faithful but capture a good deal of information about the $C^*$-algebras;
in special cases,  e.g. for the so-called AF-algebras,  the functors are full and 
faithful,  i.e. they establish  an isomorphism between the categories {\bf C*-Alg}
and {\bf Ab}.  In this section we shall introduce the following covariant functors:
\displaymath
\left\{
\begin{array}{ccc}
K_0^+ : \hbox{{\bf C*-Alg}} &\to& \hbox{{\bf Ab-Semi}},\\
K_0 : \hbox{{\bf C*-Alg}}   &\to&  \hbox{{\bf Ab}},\\
K_1 :  \hbox{{\bf C*-Alg}}   &\to& \hbox{{\bf Ab}},
\end{array}
\right.
\enddisplaymath
where {\bf C*-Alg} is the category of $C^*$-algebras and homomorphisms
between them, {\bf Ab-Semi} the category of abelian semigroups and homomorphisms
between them and {\bf Ab} the category of abelian groups and the respective 
homomorphisms.

Let $A$ be a unital $C^*$-algebra;   the definition of $K_0$-group of $A$
requires simultaneous consideration of all matrix algebras with the entries
in $A$.  
\begin{dfn}
By $M_{\infty}(A)$ one understands the algebraic direct limit of the $C^*$-algebras
$M_n(A)$ under the embeddings $a\mapsto diag ~(a, 0)$.  
\end{dfn}
\begin{rmk}
\textnormal{
The  direct limit $M_{\infty}(A)$ can be thought of as the $C^*$-algebra of
infinite-dimensional matrices whose entries are all zero except for a finite
number of the non-zero entries taken from the $C^*$-algebra $A$. 
}
\end{rmk}
\begin{dfn}
They say that projections $p.q\in M_{\infty}(A)$ are equivalent, if
there exists an element $v\in M_{\infty}(A)$ such that
\displaymath
p=v^*v \quad\hbox{and} \quad q=vv^*.
\enddisplaymath
The equivalence relation is denoted by $\sim$ and the corresponding
equivalence class of projection $p$ is denoted by $[p]$.   
\end{dfn}
\begin{dfn}
We shall write $V(A)$ to denote all equivalence classes of projections 
in the $C^*$-algebra $M_{\infty}(A)$;  in other words,
\displaymath
V(A):=\{ [p] ~: ~ p=p^*=p^2\in M_{\infty}(A)\}.  
\enddisplaymath
\end{dfn}
\begin{rmk}
\textnormal{
The set $V(A)$ has the natural structure of an abelian semigroup 
with the addition operation defined by the formula
\displaymath
[p]+[q]:= [diag~(p,q)]=[p'\oplus q'],
\enddisplaymath
where $p'\sim p, q'\sim q$ and $p'\perp q'$. The identity of the semigroup
$V(A)$ is given by $[0]$, where $0$ is the zero projection. 
}
\end{rmk}
\begin{dfn}
The semigroup $V(A)$ is said to have the cancellation property
if the equality $[p]+[r]=[q]+[r]$ implies $[p]=[q]$ for any elements 
$[p],[q], [r]\in V(A)$
\end{dfn}
\begin{exm}
\textnormal{
Let $A\cong {\Bbb C}$ be the field of complex numbers.  It is not hard to
see that $V({\Bbb C})$ has the cancellation property and  
$V({\Bbb C})\cong {\Bbb N}\cup \{0\}$ is the additive semigroup of natural
numbers with the zero. 
}
\end{exm}
 \index{functor $K_0^+$}
\begin{thm}
{\bf (Functor $K_0^+$)}
If $h: A\to A'$ is a homomorphism of the $C^*$-algebras $A$ and $A'$,
then the induced map $h_*: V(A)\to V(A')$ given by the formula
$[p]\mapsto [h(p)]$ is a well defined homomorphism of the abelian 
semigroups. In other words,  the  correspondence
\displaymath
K_0^+:  A\rightarrow V(A)
\enddisplaymath
is a covariant functor from the category {\bf C*-Alg} to the 
category {\bf Ab-Semi}.  
 \end{thm}
 \index{abelian semigroup}
Let $(S, +)$ be an abelian semigroup. One can associate to every
$(S, +)$ an abelian group as follows.  Define an equivalence relation
$\sim$ on $S\times S$ by $(x_1,y_1)\sim (x_2,y_2)$  if there exists
$z\in S$ such that $x_1+y_2+z=x_2+y_1+z$.  We shall write 
\displaymath
G(S):=S\times S/\sim
\enddisplaymath
and let $[x,y]$ be the equivalence class of $(x,y)$.  The operation
\displaymath
[x_1,y_1]+[x_2,y_2]=[x_1+x_2, y_1+y_2]
\enddisplaymath
is well defined and turns the semigroup $(S,+)$ into an abelian 
group $G(S)$.  
\begin{dfn}
The abelian group $G(S)$ is called the Grothendieck group of 
the abelian semigroup $(S,+)$.  The map  $\gamma_S: S\to G(S)$
given by the formula
\displaymath
x\mapsto [x+y, y]
\enddisplaymath
is independent of the choice of $y\in S$ and called the Grothendieck map. 
\end{dfn}
 \index{Grothendieck map}
\begin{rmk}
\textnormal{
The Grothendieck map $\gamma_S:  S\to G(S)$  is injective if and only
if the semigroup $(S, +)$ has the cancellation property. Thus although 
the Grothendieck group is always well defined,  it  is hard to recover 
the initial semigroup  $(S, +)$ from it,  if  the cancellation property
was absent.  
}
\end{rmk}
\begin{dfn}
By the $K_0$-group $K_0(A)$ of the unital $C^*$-algebra $A$ 
one understands the Grothendieck group of the abelian semigroup
$V(A)$.    
\end{dfn}
\begin{exm}
\textnormal{
Let $A\cong {\Bbb C}$ be the field of complex numbers.  It is not hard to
see that   $K_0({\Bbb C})\cong {\Bbb Z}$ is the infinite cyclic group. 
}
\end{exm}
\begin{exm}
\textnormal{
If $A\cong M_n({\Bbb C})$, then  $K_0(M_n({\Bbb C}))\cong {\Bbb Z}$. 
}
\end{exm}
\begin{exm}
\textnormal{
If $A\cong C(X)$, then  $K_0(C(X))\cong {\Bbb Z}$,   whenever $X$
is a contractible  topological space. 
}
\end{exm}
\begin{exm}
\textnormal{
If $A\cong {\cal K}$, then  $K_0({\cal K})\cong {\Bbb Z}$. 
}
\end{exm}

\begin{exm}
\textnormal{
If $A\cong B({\cal H})$, then  $K_0(B({\cal H}))\cong 0$. 
}
\end{exm}
\begin{exm}
\textnormal{
If $A$ is a $C^*$-algebra, then  $K_0(M_n(A))\cong K_0(A)$. 
}
\end{exm}
 \index{functor $K_0$}
\begin{thm}
{\bf (Functor $K_0$)}
If $h: A\to A'$ is a homomorphism of the $C^*$-algebras $A$ and $A'$,
then the induced map $h_*: K_0(A)\to K_0(A')$ given by the formula
$[p]\mapsto [h(p)]$ is a well defined homomorphism of the abelian 
groups. In other words,  the  correspondence
\displaymath
K_0:  A\rightarrow K_0(A)
\enddisplaymath
is a covariant functor from the category {\bf C*-Alg} to the 
category {\bf Ab}.  
 \end{thm}
Let $A$ be a unital $C^*$-algebra.  Consider the multiplicative group $GL_n(A)$
consisting of all invertible $n\times n$ matrices with the entries in $A$;  the 
$GL_n(A)$ is a topological group, since $A$ is endowed with a norm.  By 
$GL_{\infty}(A)$ one understands the direct limit
\displaymath
GL_{\infty}(A):= \lim_{n\to\infty} GL_n(A),
\enddisplaymath
defined by the embeddings $GL_n(A) \hookrightarrow GL_{n+1}$ given by
the formula $x\to diag (x,1)$.  (Note that this embedding is the ``exponential''
of the embedding of $M_n(A)$ into $M_{n+1}(A)$ used in the construction
of the $K_0$-group of $A$.)   Let $GL_{\infty}^0(A)$ be the connected 
component of $GL_{\infty}(A)$, which contains the unit of the group $GL_{\infty}(A)$.  
\begin{dfn}
By the $K_1$-group $K_1(A)$ of the unital $C^*$-algebra $A$ 
one understands the multiplicative abelian group consisting of the 
connected components of the group $GL_{\infty}(A)$,   i.e.
\displaymath
K_1(A):= GL_{\infty}(A)/GL_{\infty}^0(A). 
\enddisplaymath
\end{dfn}
\begin{exm}
\textnormal{
Let $A\cong {\Bbb C}$ be the field of complex numbers.  It is not hard to
see that   $K_1({\Bbb C})\cong 0$ is the infinite cyclic group. 
}
\end{exm}
\begin{exm}
\textnormal{
If $A\cong M_n({\Bbb C})$, then  $K_1(M_n({\Bbb C}))\cong 0$. 
}
\end{exm}
\begin{exm}
\textnormal{
If $A\cong C(X)$, then  $K_1(C(X))\cong 0$,   whenever $X$
is a contractible  topological space. 
}
\end{exm}
\begin{exm}
\textnormal{
If $A\cong {\cal K}$, then  $K_1({\cal K})\cong 0$. 
}
\end{exm}
\begin{exm}
\textnormal{
If $A\cong B({\cal H})$, then  $K_1(B({\cal H}))\cong 0$. 
}
\end{exm}
\begin{exm}
\textnormal{
If $A$ is a $C^*$-algebra, then  $K_1(M_n(A))\cong K_1(A)$. 
}
\end{exm}
 \index{functor $K_1$}
\begin{thm}
{\bf (Functor $K_1$)}
If $h: A\to A'$ is a homomorphism of the $C^*$-algebras $A$ and $A'$,
then the induced map $h_*: K_0(A)\to K_0(A')$ 
is a homomorphism of the multiplicative abelian  groups. In other words,  the  correspondence
\displaymath
K_1:  A\rightarrow K_0(A)
\enddisplaymath
is a covariant functor from the category {\bf C*-Alg} to the 
category {\bf Ab}.  
 \end{thm}
\begin{rmk}
\textnormal{
The functors $K_0$ and $K_1$ are related to each other by the formula
\displaymath
K_1(A)\cong K_0(C_0({\Bbb R})\otimes A),
\enddisplaymath
where $C_0({\Bbb R})$ is the commutative $C^*$-algebra of continuous 
complex-valued functions with the compact support on ${\Bbb R}$; 
the tensor product $S(A)=C_0({\Bbb R})\otimes A$ is called a 
{\it suspension} $C^*$-algebra of $A$.  The suspension can be regarded
as a natural transformation implementing the morphism between two functors
$K_0$ and $K_1$,  because the diagram shown in Fig. 3.1 is commutative. 
}
\end{rmk}
 \index{suspension}
\begin{figure}[here]
\begin{picture}(300,100)(-150,0)

\put(25,70){\vector(0,-1){35}}
\put(120,70){\vector(0,-1){35}}
\put(45,23){\vector(1,0){53}}
\put(45,83){\vector(1,0){53}}
\put(-10,20){$K_0(S(A))$}
\put(0,50){$S$}
\put(130,50){$S$}
\put(110,20){$K_0(S(A'))$}
\put(7,80){$K_1(A)$}
\put(110,80){$K_1(A')$}
\put(70,90){$h_*$}
\put(60,30){$S\circ h_*$}

\end{picture}
\caption{Suspension as a natural transformation.}
\end{figure}
For  practical calculations of the $K_0$ and $K_1$-groups 
of the crossed product $C^*$-algebras, the following result is extremely useful.
\begin{thm}
{\bf (Pimsner and Voiculescu)}
Let $A$ be a $C^*$-algebra and $\alpha\in Aut~(A)$;  consider the crossed
product $C^*$-algebra $A\rtimes_{\alpha} {\Bbb Z}$ and let 
$i: A\to A\rtimes_{\alpha} {\Bbb Z}$ be the canonical embedding. 
Then there exists  a cyclic six-term exact sequence of the abelian groups shown in Fig. 3.2. 
 \end{thm}
 \index{Pimsner-Voiculescu Theorem}
\begin{figure}[here]
\begin{picture}(300,100)(-80,0)

\put(25,35){\vector(0,1){35}}
\put(220,70){\vector(0,-1){35}}
\put(100,23){\vector(-1,0){53}}
\put(45,83){\vector(1,0){53}}
\put(150,83){\vector(1,0){53}}
\put(205,23){\vector(-1,0){53}}

\put(-15,20){$K_1(A\rtimes_{\alpha} {\Bbb Z})$}
\put(210,80){$K_0(A\rtimes_{\alpha} {\Bbb Z})$}
\put(215,20){$K_1(A)$}
\put(110,20){$K_1(A)$}
\put(7,80){$K_0(A)$}
\put(110,80){$K_0(A)$}
\put(60,90){$1-\alpha_*$}
\put(65,30){$i_*$}
\put(165,30){$1-\alpha_*$}
\put(170,90){$i_*$}

\end{picture}
\caption{The Pimsner-Voiculescu exact sequence.}
\end{figure}

\vskip1cm\noindent
{\bf Guide to the literature.}
The $K$-theory  is  the single most important tool in the $C^*$-algebra theory.
(Albeit such a theory is not very useful for  the von Neumann algebras.) 
The $K$-theory alone classifies certain types of the $C^*$-algebras and is
linked  with the geometry and topology  of manifolds,  see a review of 
the Index Theory in Chapter 10.  An  excellent account of the $K$-theory for
$C^*$-algebras  is given by   [Blackadar 1986]  \cite{B}.   For a friendly 
approach, we refer the reader to  [Wegge-Olsen 1993]  \cite{W1}.  
The textbooks by  [Fillmore 1996] \cite{F} and  [R\o rdam,  Larsen \& Laustsen  2000]  \cite{RLL}
give   an encyclopedic and detailed accounts,  respectively.

 \index{higher-dimensional\linebreak 
 noncommutative torus}

\section{Noncommutative tori}
The noncommutative tori  (and  related functors)    is a centerpiece  of these notes;  
such $C^*$-algebras appear naturally in many areas ranging from the dynamical systems
and  foliations  to the unitary  representations of groups, see  the survey  by  [Rieffel 1990]   \cite{Rie1}.   
The noncommutative tori are {\it universal $C^*$-algebras},  i.e. can be given by a finite number
of generators and relations;  such a property is extremely useful for applications in   
non-commutative  algebraic geometry,  see Chapter 5. 
The higher-dimensional noncommutative tori are  treated in Section 3.4.1; 
the case of two-dimensional torus is the subject of Section 3.4.2.

\subsection{General case}
This type of $C^*$-algebras is  an elegant example of the non-commutative  yet 
accessible $C^*$-algebra  with the well-understood structure.   Let 
\displaymath
\Theta=\left(
\matrix{0              & \theta_{12}  & \dots & \theta_{1n}\cr
             -\theta_{12} & 0  & \dots & \theta_{2n}\cr
              \vdots         & \vdots         & \ddots   &\vdots\cr
             -\theta_{1n} & -\theta_{2n} & \dots & 0 }
              \right)
\enddisplaymath
be a skew-symmetric matrix with the real entries $\theta_{ij}\in {\Bbb R}$.  
\begin{dfn}
An $n$-dimensional   noncommutative  torus  ${\cal A}_{\Theta}$  is the universal 
$C^*$-algebra  generated by $n$ unitary operators $u_1,\dots, u_n$
satisfying the commutation relations:  
\displaymath
u_ju_i=e^{2\pi i \theta_{ij}} u_iu_j, \qquad 1\le i,j\le n.
\enddisplaymath
\end{dfn}
\begin{dfn}
Let $A,B,C,D$ be the $n\times n$ matrices
with integer entries;  by $SO(n,n~|~{\Bbb Z})$ we shall understand 
a subgroup of the matrix group $GL_{2n}({\Bbb Z})$ consisting of all matrices 
of the form
\displaymath
\left(\matrix{A & B\cr C & D}\right),
\enddisplaymath
such that the matrices $A,B,C,D\in GL_k({\Bbb Z})$ satisfy the conditions:
\displaymath
A^TD+C^TB=I,\quad A^TC+C^TA=0=B^TD+D^TB,
\enddisplaymath
where $A^T, B^T, C^T$ and $D^T$ are the transpose of the 
matrices $A,B,C$ and $D$,  respectively,  and $I$ is the identity
matrix.   
\end{dfn}
\begin{rmk}
\textnormal{
The group $SO(n, n ~| ~{\Bbb Z})$ can be equivalently defined as a
subgroup of the group  $SO(n, n ~| ~{\Bbb R})$ consisting of linear transformations 
of the space ${\Bbb R}^{2n}$,  which  preserve the quadratic form 
\displaymath
x_1x_{n+1}+x_2x_{n+2}+\dots+x_nx_{2n}.
\enddisplaymath 
}
\end{rmk}
 \index{Rieffel-Schwarz Theorem}
\begin{thm}
{\bf (Rieffel and  Schwarz)}
The $n$-dimensional noncommutative   tori  ${\cal A}_{\Theta}$ and ${\cal A}_{\Theta'}$ 
are stably  isomorphic,  if  the matrices $\Theta$ and $\Theta'$
belong to the same orbit of the group $SO(n,n~|~{\Bbb Z})$  acting  on 
$\Theta$ by the formula
\displaymath
\Theta'={A\Theta+B\over  C\Theta+D}, 
\enddisplaymath
where  $A, B,  C$ and  $D$ are integer matrices. 
\end{thm}
The $K$-theory of ${\cal A}_{\Theta}$ is amazingly simple;   it basically coincides
with such for the $n$-dimensional topological tori.  
\begin{thm}
$K_0({\cal A}_{\Theta})\cong K_1({\cal A}_{\Theta})\cong {\Bbb Z}^{2^{n-1}}$. 
\end{thm}
{\it Proof.}  Let $T^{n-1}$ be an $(n-1)$-dimensional  topological torus. 
 Because ${\cal A}_{\Theta}\cong C(T^{n-1})\rtimes_{\Theta} {\Bbb Z}$
is the crossed product $C^*$-algebra,  one can apply the Pimsner-Voiculescu 
exact sequence and  get the diagram of Fig. 3.3

\begin{figure}[here]
\begin{picture}(300,100)(-50,0)

\put(25,35){\vector(0,1){35}}
\put(240,70){\vector(0,-1){35}}
\put(100,23){\vector(-1,0){53}}
\put(60,83){\vector(1,0){43}}
\put(185,83){\vector(1,0){33}}
\put(225,23){\vector(-1,0){33}}

\put(-5,20){$K_1({\cal A}_{\Theta})$}
\put(230,80){$K_0({\cal A}_{\Theta})$}
\put(235,20){$K_1(C(T^{n-1}))$}
\put(110,20){$K_1(C(T^{n-1}))$}
\put(-10,80){$K_0(C(T^{n-1}))$}
\put(110,80){$K_0(C(T^{n-1}))$}
\put(60,90){$1-\Theta_*$}
\put(65,30){$i_*$}
\put(195,30){$1-\Theta_*$}
\put(200,90){$i_*$}

\end{picture}
\caption{The Pimsner-Voiculescu exact sequence for ${\cal A}_{\Theta}$.}
\end{figure}

It is known, that for $i\in \{0,1\}$ it holds $K_i(C(T^{n-1}))\cong K^{i-1}(T^{n-1})\cong {\Bbb Z}^{2^{n-2}}$,
where $K^{\bullet}(T^{n-1})$ are the topological $K$-groups of torus $T^{n-1}$. 
Thus one gets two split exact  sequences
\displaymath
0\rightarrow  {\Bbb Z}^{2^{n-2}}\rightarrow K_i({\cal A}_{\Theta})\rightarrow  {\Bbb Z}^{2^{n-2}}\rightarrow 0.
\enddisplaymath
Hence $K_i({\cal A}_{\Theta})\cong  {\Bbb Z}^{2^{n-2}}\oplus  {\Bbb Z}^{2^{n-2}}\cong  {\Bbb Z}^{2^{n-1}}$.
$\square$

\begin{rmk}
\textnormal{
Unless $n<5$,  the abelian semigroup $V({\cal A}_{\Theta})$ does not have cancellation and,
therefore,  $V({\cal A}_{\Theta})$ cannot be embedded into the group $K_0({\cal A}_{\Theta})\cong {\Bbb Z}^{2^{n-1}}$.
However, for  $n=2$  the semigroup  $V({\cal A}_{\Theta})$ does embed into $K_0({\cal A}_{\Theta})\cong {\Bbb Z}^2$
turning $K_0({\cal A}_{\Theta})$ into a totally ordered abelian group;  in this case the functor 
\displaymath
{\cal A}_{\Theta}\longrightarrow V({\cal A}_{\Theta})
\enddisplaymath
yields an equivalence of the categories of 2-dimensional noncommutative tori  and
totally ordered abelian groups of rank 2. 
}
\end{rmk}

 \index{noncommutative torus}

\subsection{Two-dimensional  torus}
In case $n=2$,  we have
\displaymath
\Theta=\left(\matrix{0 & \theta_{12}\cr -\theta_{12} & 0}\right)
\enddisplaymath
and,  therefore, the noncommutative torus ${\cal A}_{\Theta}$ depends on a 
single real parameter $\theta_{12}:=\theta$;  we shall simply write ${\cal A}_{\theta}$
in this case.  Thus one gets the following  
\begin{dfn}
The noncommutative torus ${\cal A}_{\theta}$ is the universal $C^*$-algebra
on the generators $u, u^*, v, v^*$ and relations
\displaymath
\left\{
\begin{array}{ccc}
vu &=& e^{2\pi i\theta} uv\\ 
u^*u &=&  uu^*=1\\
v^*v &=& vv^*=1.  
\end{array}
\right.
\enddisplaymath
\end{dfn}
\begin{thm}
The noncommutative tori ${\cal A}_{\theta}$ and ${\cal A}_{\theta'}$ are stably isomorphic
if and only if
\displaymath
\theta'={a\theta+b\over c\theta+d},
\enddisplaymath
where $a,b.c,d\in {\Bbb Z}$ are such that $ad-bc=\pm 1$.   
\end{thm}
{\it Proof.}  We shall use the general formula involving  group 
$SO(n, n ~| ~{\Bbb Z})$.  The skew-symmerix matrix $\Theta$ has the form
\displaymath
\Theta=\left(\matrix{0 & \theta\cr -\theta & 0}\right).
\enddisplaymath
The $2\times 2$ integer matrices $A,B,C,D$  take the form
\displaymath
\left\{
\begin{array}{ccc}
A &=& \left(\matrix{a & 0\cr 0 & a}\right) \\ 
B &=&  \left(\matrix{0 & b \cr -b & 0}\right)\\
C &=&  \left(\matrix{0 & -c\cr  c & 0}\right)\\
D &=& \left(\matrix{d & 0 \cr  0 & d}\right).  
\end{array}
\right.
\enddisplaymath
The reader can verify, that $A,B,C$ and $D$ satisfy the conditions
\displaymath
A^TD+C^TB=I,\quad A^TC+C^TA=0=B^TD+D^TB.
\enddisplaymath
Therefore,  the matrix 
\displaymath
\left(\matrix{A & B\cr C & D}\right) \in SO(2, 2 ~| ~{\Bbb Z}).
\enddisplaymath
The direct computation gives the following equations
\displaymath
\Theta'={A\Theta+B\over  C\Theta+D}=
\left(\matrix{0 & {a\theta+b\over c\theta+d}\cr
 -{a\theta+b\over c\theta+d} & 0}\right).
 \enddisplaymath
In other words,  one gets $\theta'={a\theta+b\over c\theta+d}$.
$\square$

\begin{thm}
$K_0({\cal A}_{\theta})\cong K_1({\cal A}_{\theta})\cong {\Bbb Z}^{2}$. 
\end{thm}
{\it Proof.}  Follows from the general formula when $n=2$. 
$\square$

\begin{thm}
{\bf (Pimsner-Voiculescu,  Rieffel)}
The abelian semigroup $V({\cal A}_{\theta})$ embeds into the group 
$K_0({\cal A}_{\theta})\cong {\Bbb Z}^{2}$;  the image of the injective map
\displaymath
V({\cal A}_{\theta})\longrightarrow K_0({\cal A}_{\theta})
\enddisplaymath
is given by the formula 
\displaymath
\{(m,n)\in {\Bbb Z}^2 ~|~ m+n\theta>0\}.
\enddisplaymath
\end{thm}
\begin{cor}
The abelian group $K_0({\cal A}_{\theta})$ gets an order structure $\ge$ coming from
the semigroup $V({\cal A}_{\theta})$ and given by the formula
\displaymath
 z_1\ge z_2 ~\hbox{iff} ~ z_1-z_2\in V({\cal A}_{\theta})\subset K_0({\cal A}_{\theta}).
\enddisplaymath
If {\bf NTor}  is the category of 2-dimensional noncommutative tori whose arrows are 
stable homomorphisms of the  tori and {\bf Ab-Ord} is the category 
of ordered abelian groups of rank 2 whose arrows are order-homomorphisms of the 
groups,  then the functor
\displaymath
K_0^+: \hbox{{\bf NTor}}\to \hbox{{\bf Ab-Ord}}
\enddisplaymath
is full and faithful, i.e. an isomorphism of the categories.
\end{cor}

\vskip1cm\noindent
{\bf Guide to the literature.}
We encourage the reader to  start  with the survey paper 
by  [Rieffel 1990]   \cite{Rie1}.    The noncommutative torus is also known as the irrational rotation 
algebra,   see  [Pimsner \& Voiculescu  1980]   \cite{PiVo1} and  [Rieffel 1981]  \cite{Rie2}.  
The higher-dimensional noncommutative tori and their $K$-theory  were considered
by [Rieffel \& Schwarz 1999]  \cite{RiSch1}.
The real multiplication has been introduced  and studied by  [Manin 2003] \cite{Man1};
a higher-dimensional analog of real multiplication was considered in \cite{Nik14}.

 \index{AF-algebra}

\section{AF-algebras }
This vast class of the $C^*$-algebras  was introduced by [Bratteli 1972]  \cite{Bra1};   
it is a higher-rank generalization of the $UHF$-algebras studied by [Glimm 1960]  \cite{Gli1}.
The $AF$-algebras are classified by their
$K$-theory;   the latter is known as a {\it dimension group} introduced by  [Elliott 1976]  \cite{Ell1}. 
Such a  property is  highly used in the construction
of covariant functors from different geometric and topological spaces  to the  category
of $AF$-algebras.  Unlike the nonommutative tori,  the $AF$-algebras are {\it not}  universal,
i.e. cannot be given by the generators and relations.  Section 3.5.1 covers       
general $AF$-algebras;   Section 3.5.2  is devoted to an important  special 
of {\it stationary} $AF$-algebras.

\subsection{General case}
The field of rational numbers ${\Bbb Q}$ admits a completion consisting
of all real numbers ${\Bbb R}$;  similarly,  the set of all finite-dimensional 
$C^*$-algebras $M_n({\Bbb C})$ has a completion by the so-called 
{\it approximately finite-dimensional $C^*$-algebras},  or AF-algebras,  for brevity.  
Unlike the noncommutative tori,  the AF-algebras are not universal (i.e. cannot
be given by the generators and relations)  yet have  remarkable $K$-theory
and classification.  We begin with the following important
\begin{thm}
{\bf (Bratteli,  Elliott)}
Any finite-dimensional $C^*$-algebra $A$ is isomorphic to
\displaymath
M_{n_1}({\Bbb C})\oplus  M_{n_2}({\Bbb C})\oplus\dots\oplus   M_{n_k}({\Bbb C})
\enddisplaymath
for some positive integers $n_1, n_2, \dots, n_k$.  Moreover,
\displaymath
K_0(A)\cong {\Bbb Z}^k ~\hbox{and} ~K_1(A)\cong 0.
\enddisplaymath
\end{thm}
\begin{dfn}
By   an AF-algebra  one understands a $C^*$-algebra ${\Bbb A}$ which is 
the  norm closure of an ascending sequence of the embeddings $\varphi_i$ of  the finite-dimensional
$C^*$-algebras
\displaymath
A_1\buildrel\rm\varphi_1\over\hookrightarrow A_2
   \buildrel\rm\varphi_2\over\hookrightarrow\dots,
\enddisplaymath
The set-theoretic limit  ${\Bbb A}=\lim A_i$ is endowed with a natural algebraic structure given by the formula
$a_m+b_k\to a+b$; here $a_m\to a,b_k\to b$ for the  sequences $a_m\in A_m,b_k\in A_k$.  
\end{dfn}
\begin{rmk}
\textnormal{
To keep track on all  embeddings $\varphi_i: A_i\hookrightarrow A_{i+1}$
it is convenient to arrange them into a graph called the  Bratteli diagram.  
 Let  $A_i=M_{i_1}\oplus\dots\oplus M_{i_k}$ and 
$A_{i'}=M_{i_1'}\oplus\dots\oplus M_{i_k'}$ be 
the semi-simple $C^*$-algebras and $\varphi_i: A_i\to A_{i'}$ the  homomorphism. 
One has the two sets of vertices $V_{i_1},\dots, V_{i_k}$ and $V_{i_1'},\dots, V_{i_k'}$
joined by the $a_{rs}$ edges, whenever the summand $M_{i_r}$ contains $a_{rs}$
copies of the summand $M_{i_s'}$ under the embedding $\varphi_i$. 
As $i$ varies, one obtains an infinite graph;  the graph is called a {\it Bratteli diagram} of the
$AF$-algebra ${\Bbb A}$.    The Bratteli diagram defines a unique  $AF$-algebra,  yet the converse is false.
}
\end{rmk}
 \index{Bratteli diagram}
\begin{exm}
\textnormal{
The embeddings
\displaymath
M_1({\Bbb C})\buildrel\rm Id\over\hookrightarrow M_2({\Bbb C})
   \buildrel\rm Id\over\hookrightarrow\dots,
\enddisplaymath
define the $C^*$-algebra ${\cal K}$ of all compact operators on a Hilbert space ${\cal H}$.
The corresponding Bratteli diagram is shown in Fig. 3.4.     
}
\end{exm}
\begin{figure}[here]
\begin{picture}(300,60)(0,0)

\put(120,38){$\bullet$}
\put(140,38){$\bullet$}
\put(160,38){$\bullet$}

\put(120,40){\line(1,0){20}}
\put(140,40){\line(1,0){20}}

\put(180,40){$\dots$}

\put(125,52){$1$}
\put(145,52){$1$}

\end{picture}

\caption{Bratteli diagram of the AF-algebra ${\cal K}$.}
\end{figure}
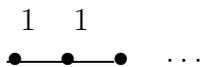

\begin{exm}\label{exm3.5.2} 
\textnormal{
Let $\theta\in {\Bbb R}$ be an irrational number and consider the regular continued fraction:
\displaymath
\theta=a_0+{1\over\displaystyle a_1+
{\strut 1\over\displaystyle a_2+\dots}}
:=[a_0, a_1, a_2, \dots].
\enddisplaymath
One can define an AF-algebra ${\Bbb A}_{\theta}$ by the Bratteli diagram shown in Fig. 3.5;
the AF-algebra is called an {\it Effros-Shen algebra}.  
}
\end{exm}
 \index{Effros-Shen algebra}
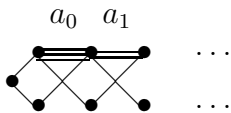
\begin{figure}[here]
\begin{picture}(300,60)(0,0)
\put(108,27){$\bullet$}
\put(118,18){$\bullet$}
\put(138,18){$\bullet$}
\put(158,18){$\bullet$}
\put(118,38){$\bullet$}
\put(138,38){$\bullet$}
\put(158,38){$\bullet$}

\put(110,30){\line(1,1){10}}
\put(110,30){\line(1,-1){10}}
\put(120,42){\line(1,0){20}}
\put(120,40){\line(1,0){20}}
\put(120,38){\line(1,0){20}}
\put(120,40){\line(1,-1){20}}
\put(120,20){\line(1,1){20}}
\put(140,41){\line(1,0){20}}
\put(140,39){\line(1,0){20}}
\put(140,40){\line(1,-1){20}}
\put(140,20){\line(1,1){20}}

\put(180,20){$\dots$}
\put(180,40){$\dots$}

\put(125,52){$a_0$}
\put(145,52){$a_1$}

\end{picture}

\caption{The Effros-Shen algebra  ${\Bbb A}_{\theta}$.}
\end{figure}
\begin{rmk}
\textnormal{
Because $K_1(A_i)\cong 0$ and the $K_1$-functor is continuous on the inductive limits,
one gets
\displaymath
K_1({\Bbb A})\cong 0,
\enddisplaymath
for any AF-algebra ${\Bbb A}$.  The  group $K_0({\Bbb A})$ is more interesting,
because each embedding $\varphi_i: A_i\hookrightarrow  A_{i+1}$ induces a homomorphism
$(\varphi_i)_*: K_0(A_i)\to K_0(A_{i+1})$ of the abelian groups  $K_0(A_i)$; 
we need  the following
}
\end{rmk}
 \index{dimension group}
\begin{dfn}\label{dfn3.5.2}
By a dimension group one understands an ordered abelian group $(G, G^+)$  which is the 
limit of the sequence of ordered abelian groups ${\Bbb Z}^{n_i}$  and positive
homomorphisms   $(\varphi_i)_*: {\Bbb Z}^{n_i}\to {\Bbb Z}^{n_{i+1}}$: 
\displaymath
{\Bbb Z}^{n_1}\buildrel\rm(\varphi_1)_*\over\rightarrow 
{\Bbb Z}^{n_2} \buildrel\rm(\varphi_2)_*\over\rightarrow\dots,
\enddisplaymath
where the positive cone of ${\Bbb Z}^{n_i}$ is defined by the formula
\displaymath
({\Bbb Z}^{n_i})^+=\{(x_1,\dots,x_{n_i})\in {\Bbb Z}^{n_i} : x_j\ge 0\}. 
\enddisplaymath
\end{dfn}
 \index{Elliott Theorem}
\begin{thm}\label{thm3.5.2}
{\bf (Elliott)}  If ${\Bbb A}$ is an AF-algebra,  then
\displaymath
K_0({\Bbb A})\cong G   ~\hbox{and} ~V({\Bbb A})\cong G^+.
\enddisplaymath
Moreover,   if {\bf AF-Alg} is the category of the AF-algebras and 
homomorphisms between them and {\bf Ab-Ord} is the category 
of ordered abelian groups (with the scaled units)  and order-preserving homomorphisms
between them, then the functor
\displaymath
K_0^+: \hbox{{\bf AF-Alg}}\to \hbox{{\bf Ab-Ord}}
\enddisplaymath
is full and faithful, i.e. an isomorphism of the categories.
\end{thm}
\begin{exm}
\textnormal{
If ${\Bbb A}\cong {\cal K}$,  then $K_0({\cal K})\cong {\Bbb Z}$
and $V({\cal K})\cong {\Bbb Z}^+$,  where ${\Bbb Z}^+$ is an additive
semigroup of the non-negative integers. 
}
\end{exm}
\begin{exm}
\textnormal{
For the Effros-Shen algebra ${\Bbb A}_{\theta}$,  it has been proved that
\displaymath
K_0({\Bbb A}_{\theta})\cong {\Bbb Z}^2 ~\hbox{and}  
~V({\Bbb A}_{\theta})\cong\{(m,n)\in {\Bbb Z}^2 ~|~m+n\theta>0\}. 
\enddisplaymath
}
\end{exm}
\begin{thm}\label{thm3.5.3}
{\bf (Pimsner-Voiculescu)}  
Let ${\cal A}_{\theta}$ be a two-dimensional noncommutative torus
and ${\Bbb A}_{\theta}$ the corresponding Effros-Shen algebra;
there exists a natural embedding
\displaymath
{\cal A}_{\theta}\hookrightarrow {\Bbb A}_{\theta},
\enddisplaymath
which induces an order-isomorphism of the corresponding
semigroups  $V({\cal A}_{\theta})\cong V({\Bbb A}_{\theta})$.  
\end{thm}

 \index{stationary AF-algebra}

\subsection{Stationary AF-algebras}
The following type of the AF-algebras will be  critical for applications in topology and 
number theory;  roughly speaking,   such AF-algebras are depicted by the periodic Bratteli 
diagrams.  
\begin{dfn}
By a stationary AF-algebra ${\Bbb A}_{\varphi}$ one understands an  AF-algebra for which the embedding 
homomorphisms  $\varphi_1 =\varphi_2=\dots=Const$;  in other words,   the stationary
AF-algebra has the form:   
\displaymath
A_1\buildrel\rm\varphi\over\hookrightarrow A_2
   \buildrel\rm\varphi\over\hookrightarrow\dots,
\enddisplaymath
where $A_i$ are  some finite-dimensional $C^*$-algebras.  
\end{dfn}
\begin{rmk}
\textnormal{
The corresponding dimension group $(G, G^+_{\varphi_*})$ can be written as
\displaymath
{\Bbb Z}^{k}\buildrel\rm
\varphi_*
\over\longrightarrow {\Bbb Z}^{k}
\buildrel\rm
\varphi_*
\over\longrightarrow
{\Bbb Z}^{k}\buildrel\rm
\varphi_*
\over\longrightarrow \dots
\enddisplaymath
where $\varphi_*$ is a matrix with non-negative integer entries;  
one can take a minimal power of  $\varphi_*$ to obtain a strictly positive integer 
matrix  $B$ -- which we always assume to be the case. 
}
\end{rmk}
 \index{golden mean}
\begin{exm}
\textnormal{
Let $\theta\in {\Bbb R}$ be a quadratic irrationality  called the {\it golden mean}:
\displaymath
{1+\sqrt{5}\over 2}=1+{1\over\displaystyle 1+
{\strut 1\over\displaystyle 1+\dots}}
=[1, 1, 1, \dots].
\enddisplaymath
The corresponding Effros-Shen algebra has a periodic Bratteli diagram  shown in Fig. 3.6;
the dimension group $(G, G^+_B)$ is given by the following sequence of positive homomorphisms:
\displaymath
{\Bbb Z}^{2}\buildrel\rm
\left(\small\matrix{1 & 1\cr 1 & 0}\right)
\over\longrightarrow {\Bbb Z}^{2}
\buildrel\rm
\left(\small\matrix{1 & 1\cr 1 & 0}\right)
\over\longrightarrow
{\Bbb Z}^{2}\buildrel\rm
\left(\small\matrix{1 & 1\cr 1 & 0}\right)
\over\longrightarrow \dots
\enddisplaymath
}
\end{exm}
\begin{figure}[here]
\begin{picture}(300,60)(0,0)
\put(108,27){$\bullet$}
\put(118,18){$\bullet$}
\put(138,18){$\bullet$}
\put(158,18){$\bullet$}
\put(118,38){$\bullet$}
\put(138,38){$\bullet$}
\put(158,38){$\bullet$}

\put(110,30){\line(1,1){10}}
\put(110,30){\line(1,-1){10}}
\put(120,40){\line(1,0){20}}
\put(120,40){\line(1,-1){20}}
\put(120,20){\line(1,1){20}}
\put(140,40){\line(1,0){20}}
\put(140,40){\line(1,-1){20}}
\put(140,20){\line(1,1){20}}

\put(180,20){$\dots$}
\put(180,40){$\dots$}


\end{picture}

\caption{Stationary Effros-Shen algebra  ${\Bbb A}_{{1+\sqrt{5}\over 2}}$.}
\end{figure}
 \index{shift automorphism}
\begin{dfn}\label{dfn3.5.4}
If ${\Bbb A}_{\varphi}$ is a stationary AF-algebra,  consider an order-automorphism of its
dimension group $(G,G_{\varphi_*}^+)$ generated by the 1-shift of the diagram:   
\begin{picture}(300,60)(-50,0)

\put(100,18){${\Bbb Z}^k$}
\put(138,18){${\Bbb Z}^k$}
\put(170,18){${\Bbb Z}^k$}
\put(100,38){${\Bbb Z}^k$}
\put(138,38){${\Bbb Z}^k$}
\put(170,38){${\Bbb Z}^k$}

\put(110,40){\vector(1,0){20}}
\put(110,20){\vector(1,0){20}}
\put(110,38){\vector(2,-1){25}}
\put(150,40){\vector(1,0){20}}
\put(150,20){\vector(1,0){20}}
\put(150,38){\vector(2,-1){25}}

\put(190,20){$\dots$}
\put(190,40){$\dots$}

\put(115,45){$\varphi_*$}
\put(155,45){$\varphi_*$}

\put(115,10){$\varphi_*$}
\put(155,10){$\varphi_*$}

\end{picture}

\noindent
and let $\sigma_{\varphi}: {\Bbb A}_{\varphi}\to {\Bbb A}_{\varphi}$ be the corresponding automorphism
of ${\Bbb A}_{\varphi}$;  the $\sigma_{\varphi}$ is called the shift automorphism. 
\end{dfn}
 \index{stationary AF-algebra}
 \index{stationary dimension group}
To classify all stationary AF-algebras,  one  needs  the following set of invariants. 
Let $B\in M_k({\Bbb Z})$  be a matrix with the strictly positive entries corresponding to
a stationary dimension group $(G, G^+_B)$:
\displaymath
{\Bbb Z}^{k}\buildrel\rm
B
\over\longrightarrow {\Bbb Z}^{k}
\buildrel\rm
B
\over\longrightarrow
{\Bbb Z}^{k}\buildrel\rm
B
\over\longrightarrow \dots
\enddisplaymath
By the Perron-Frobenius theory,  matrix $B$ has a real eigenvalue $\lambda_B>1$, which exceeds
the absolute values of other roots of the characteristic polynomial of $B$.
Note that $\lambda_B$ is an  algebraic integer.  Consider
the  real algebraic number field $K={\Bbb Q}(\lambda_B)$ obtained as 
an extension of the field of the rational numbers by the algebraic 
number $\lambda_B$. Let $(v^{(1)}_B,\dots,v^{(k)}_B)$ be the eigenvector
corresponding to the eigenvalue $\lambda_B$. One can  normalize the eigenvector so 
that $v^{(i)}_B\in K$.  Consider  the  ${\Bbb Z}$-module
${\goth m}={\Bbb Z}v^{(1)}_B+\dots+{\Bbb Z}v^{(k)}_B$. The module ${\goth m}$
brings in two new arithmetic objects: (i) the ring $\Lambda$ of the endomorphisms of ${\goth m}$  
and (ii) an ideal $I$ in the ring $\Lambda$,  such that $I={\goth m}$  after a scaling. 
The ring $\Lambda$ is an order in the algebraic number field $K$ and therefore one can talk about the ideal
classes in $\Lambda$. The ideal class of $I$ is denoted by $[I]$. 
\begin{thm}
{\bf (Handelman)}
The triple $(\Lambda, [I], K)$ is an invariant of the  stable isomorphism class of the  stationary AF-algebra 
${\Bbb A}_{\varphi}$. 
\end{thm}
 \index{Handelman invariant}

\vskip1cm\noindent
{\bf Guide to the literature.}
The $AF$-algebras  were  introduced by [Bratteli 1972]  \cite{Bra1};   
it is a higher-rank generalization of the $UHF$-algebras studied by [Glimm 1960]  \cite{Gli1}.
The   dimension groups of the $AF$-algebras were  introduced by  [Elliott 1976]  \cite{Ell1};
he showed that such groups classify the $AF$-algebras {\it ibid}. 
The stationary $AF$-algebras are covered in [Effros 1981] \cite{E},  Chapter 6.  
The arithmetic invariant $(\Lambda, [I], K)$ of stationary $AF$-algebras
was introduced by [Handelman 1981]  \cite{Han2}.

 \index{UHF-algebra}

\section{UHF-algebras}
The uniformly hyper-finite $C^*$-algebras (the {\it UHF-algebras}, for brevity)  is a special
type of the AF-algebras;  they were the first AF-algebras studied and classified by an 
invariant called a {\it supernatural number}.  
\begin{dfn}
The  UHF-algebra  is an AF-algebra which is isomorphic to the inductive limit of the
sequence of finite-dimensional $C^*$-algebras of the form:
\displaymath
M_{k_1}({\Bbb C})\to M_{k_1}({\Bbb C})\otimes M_{k_2}({\Bbb C})\to 
M_{k_1}({\Bbb C})\otimes M_{k_2}({\Bbb C})\otimes M_{k_3}({\Bbb C})\to\dots,
\enddisplaymath
where $M_{k_i}({\Bbb C})$ is a matrix $C^*$-algebra and $k_i\in\{1, 2, 3,\dots\}$;
we shall denote the UHF-algebra by $M_{\ka}$, where $\ka=(k_1,k_2,k_3,\dots)$.
\end{dfn}
\begin{exm}
\textnormal{
Let $p$ be a prime number  and consider the UHF-algebra
\displaymath
M_{p^{\infty}}:= M_p({\Bbb C})\otimes M_p({\Bbb C})\otimes\dots
\enddisplaymath
For $p=2$ the algebra $M_{2^{\infty}}$  is known as   a  Canonical Anticommutation Relations $C^*$-algebra 
(the {\it CAR} or {\it Fermion algebra});  its Bratteli diagram is shown in Fig. 3.7.  
}
\end{exm}
 \index{Fermion algebra}
\begin{figure}[here]
\begin{picture}(300,60)(0,0)

\put(120,38){$\bullet$}
\put(140,38){$\bullet$}
\put(160,38){$\bullet$}

\put(120,39){\line(1,0){20}}
\put(120,41){\line(1,0){20}}
\put(140,39){\line(1,0){20}}
\put(140,41){\line(1,0){20}}

\put(180,40){$\dots$}

\put(125,52){$2$}
\put(145,52){$2$}

\end{picture}

\caption{Bratteli diagram of the CAR algebra $M_{2^{\infty}}$.}
\end{figure}
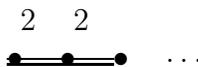

\bigskip\noindent
To classify the UHF-algebras up to the stable isomorphism, one needs
the following construction. Let $p$ be a prime number and 
$n=\sup~\{0\le j\le \infty : p^j~|~\prod_{i=1}^{\infty}k_i\}$;
denote by $\n=(n_1,n_2,\dots)$ an infinite sequence of $n_i$ as $p_i$
runs through the ordered set of all primes.  
\begin{dfn}
By ${\Bbb Q}(\n)$ we understand an additive subgroup of ${\Bbb Q}$
consisting of rational numbers  whose denominators divide the ``supernatural
number'' $p_1^{n_1}p_2^{n_2}\dots$,   where each $n_j$ belongs to the set
$\{0,1,2,\dots,\infty\}$.
\end{dfn}
 \index{supernatural number}
\begin{rmk}
\textnormal{
The ${\Bbb Q}(\n)$ is a dense subgroup  of ${\Bbb Q}$ and every dense subgroup of ${\Bbb Q}$ 
  containing ${\Bbb Z}$ is given by ${\Bbb Q}(\n)$  for some $\n$.
}
\end{rmk}
\begin{thm}
{\bf (Glimm)}  
 $K_0(M_{\ka})\cong {\Bbb Q}(\n)$.  
\end{thm}
\begin{exm}
\textnormal{
For the CAR algebra $M_{2^{\infty}}$,  one gets 
\displaymath
K_0(M_{2^{\infty}})\cong {\Bbb Z}\left[{1\over 2}\right],
\enddisplaymath
where ${\Bbb Z}[{1\over 2}]$ are the so-called {\it dyadic}  rationals. 
}
\end{exm}
 \index{dyadic number}
 \index{Glimm Theorem}
\begin{thm}
{\bf (Glimm)}  
 The UHF-algebras  $M_{\ka}$ and $M_{\ka'}$ are stably isomorphic
if and only if  $r{\Bbb Q}(\n)=s{\Bbb Q}(\n')$ for  some positive integers $r$ and $s$. 
\end{thm}

\vskip1cm\noindent
{\bf Guide to the literature.}
The UHF-algebras were introduced and classified by [Glimm 1960]  \cite{Gli1}.
The supernatural numbers are covered in  [Effros 1981]  \cite{E},  p.28  
and  [R\o rdam,  Larsen \& Laustsen  2000]  \cite{RLL},   Section 7.4.

 \index{Cuntz-Krieger algebra}

\section{Cuntz-Krieger algebras}
Roughly speaking,  the Cuntz-Krieger algebras are  universal $C^*$-algebras
with  a (typically) finite $K$-theory;  they  have an interesting crossed product
structure linking them to the stationary AF-algebras.  Let
\displaymath
B=\left(
\matrix{b_{11}             & b_{12}  & \dots & b_{1n}\cr
             b_{21} & b_{22}  & \dots & b_{2n}\cr
              \vdots         & \vdots         & \ddots   &\vdots\cr
             b_{n1} & b_{n2} & \dots & b_{nn} }
              \right)
\enddisplaymath
be a square matrix, such that $b_{ij}\in \{0, 1, 2, \dots\}$.   
\begin{dfn}
By a   Cuntz-Krieger algebra, ${\cal O}_B$ one understands the $C^*$-algebra
generated by the  partial isometries $s_1,\dots, s_n$ which satisfy  the relations
\displaymath
\left\{
\begin{array}{ccc}
s_1^*s_1 &=& b_{11} s_1s_1^*+b_{12} s_2s_2^*+\dots+b_{1n}s_ns_n^*\\ 
s_2^*s_2 &=& b_{21} s_1s_1^*+b_{22} s_2s_2^*+\dots+b_{2n}s_ns_n^*\\ 
                  &\vdots&\\
s_n^*s_n &=& b_{n1} s_1s_1^*+b_{n2} s_2s_2^*+\dots+b_{nn}s_ns_n^*.                   
\end{array}
\right.
\enddisplaymath
\end{dfn}
\begin{exm}
\textnormal{
Let 
\displaymath
B=\left(
\matrix{1             & 1 & \dots & 1\cr
             1 & 1 & \dots & 1\cr
              \vdots         & \vdots         & \ddots   &\vdots\cr
             1 & 1 & \dots & 1}
              \right).
\enddisplaymath
Then ${\cal O}_B$ is called a {\it Cuntz algebra} and denoted by ${\cal O}_n$.
}
\end{exm}
 \index{Cuntz algebra}
\begin{rmk}
\textnormal{
 It is known,  that the $C^*$-algebra ${\cal O}_B$ is simple, 
whenever matrix $B$ is irreducible, i.e. a certain power
of $B$ is a strictly positive integer matrix. 
}
\end{rmk}
 \index{Cuntz-Krieger Theorem}
\begin{thm}
{\bf (Cuntz \& Krieger)}
If ${\cal O}_B$ is a Cuntz-Krieger algebra,  then
\displaymath
\left\{
\begin{array}{ccc}
K_0({\cal O}_B) &\cong& {{\Bbb Z}^n \over (I-B^T){\Bbb Z}^n}\\ 
&&\\
K_1({\cal O}_B)  &\cong & Ker~(I-B^T),                   
\end{array}
\right.
\enddisplaymath
where $B^T$ is a transpose of the matrix $B$.
\end{thm}
\begin{rmk}
\textnormal{
It is not difficult to see,  that whenever $det~(I-B^T)\ne 0$, the $K_0({\cal O}_B)$
is a finite abelian group and $K_1({\cal O}_B)=0$.  The both groups are invariants of the
stable isomorphism class of the Cuntz-Krieger algebra.
 }
\end{rmk}
 \index{Cuntz-Krieger Crossed Product Theorem}
\begin{thm}
{\bf (Cuntz \& Krieger)}
Let ${\Bbb A}_{\varphi}$ be a stationary  AF-algebra,  such that $\varphi_*\cong B$;
let $\sigma_{\varphi}: {\Bbb A}_{\varphi}\to {\Bbb A}_{\varphi}$ be the corresponding
shift automorphism of ${\Bbb A}_{\varphi}$.  Then
\displaymath
{\cal O}_B\otimes {\cal K}\cong {\Bbb A}_{\varphi}\rtimes_{\sigma_{\varphi}}{\Bbb Z}, 
\enddisplaymath
where the crossed product is taken by the shift automorphism $\sigma_{\varphi}$.  
\end{thm}

\vskip1cm\noindent
{\bf Guide to the literature.}
The Cuntz-Krieger algebras were introduced and studied in [Cuntz \& Krieger 1980]
\cite{CuKr1};   they generalize the {\it Cuntz algebras}  introduced by  [Cuntz 1977]  
\cite{Cun1}.  The crossed product structure and remarkable $K$-theory of ${\cal O}_B$ 
make such algebras important in the non-commutative localization theory,  see
Section 6.3.

\section*{Exercises}

\begin{enumerate}

 \index{spectrum of $C^*$-algebra}

\item
Recall that {\it spectrum} {\bf Sp} $(a)$ of an element $a$ of 
the $C^*$-algebra $A$ is the set of complex numbers $\lambda$,
such that $a-\lambda e$ is not invertible;   show that if $p\in  A$ is 
a projection, then  {\bf Sp} $(p)\subseteq\{0, 1\}$.

\item
Show that if $u\in  A$ is a unitary element,  then {\bf Sp} 
$(u) \subseteq \{z\in {\Bbb C} ~:~  |z|=1\}$.

\item
Show that $||p-q||\le 1$ for every pair of projections  in the $C^*$-algebra $A$.

\item
Show that $||u-v||\le 2$ for every pair of unitary elements in the $C^*$-algebra $A$.

\item
Recall that projections $p,q\in A$ are {\it orthogonal} $p\perp q$,  if $pq=0$;
show that the following three conditions are equivalent:

\subitem
(i) $p\perp q$;

\subitem
(ii) $p+q$ is a projection;

\subitem
(iii)  $p+q\le 1$.

\item 
 Recall that $v$ is a {\it partial isometry} if $v^*v$ is a projection;
 show that $v=vv^*v$ and conclude that $vv^*$ is a projection.

\item
Prove that if $A\cong {\Bbb C}$,  then  $K_0({\Bbb C})\cong {\Bbb Z}$ is the 
infinite cyclic group.

\item
Prove that if $A\cong M_n({\Bbb C})$,  then  $K_0(M_n({\Bbb C}))\cong {\Bbb Z}$. 

\item
Prove that if $A\cong {\Bbb C}$,  then  $K_1({\Bbb C})\cong 0$.

\item
Prove that if $A\cong M_n({\Bbb C})$,  then  $K_1(M_n({\Bbb C}))\cong 0$.

\item
Calculate the $K$-theory of noncommutative torus ${\cal A}_{\theta}$ using
its crossed product structure. (Hint:  use the Pimsner-Voiculscu exact sequence
for crossed products.)

 \item
 Calculate the Handelman invariant  $(\Lambda, [I], K)$ of the stationary
 Effros-Shen algebra ${\Bbb A}_{1+\sqrt{5}\over 2}$ (golden mean algebra).

\item
Calculate the $K$-theory of the Cuntz-Krieger algebra ${\cal O}_B$, 
where 
\displaymath
B=\left(\matrix{5 & 1\cr 4 & 1}\right).
\enddisplaymath
(Hint: use the formulas 
$K_0({\cal O}_B) \cong {{\Bbb Z}^n \over (I-B^T){\Bbb Z}^n}$
and  $K_1({\cal O}_B)  \cong  Ker~(I-B^T)$;  bring the matrix $I-B^T$
to the {\it Smith normal form},  see e.g.  [Lind \& Marcus 1995]  \cite{LM}.)

\item
Repeat  the  exercise for matrix
\displaymath
B=\left(\matrix{5 & 2\cr 2 & 1}\right).
\enddisplaymath

 \end{enumerate}





\part{NONCOMMUTATIVE  INVARIANTS}
\chapter{Topology}
In this chapter we shall construct several functors on the topological spaces
with values in the category of stationary AF-algebras or the Cuntz-Krieger
algebras.  The functors give rise to a set of noncommutative invariants 
some of which can be explicitly calculated; all the invariants are  homotopy
invariants of the corresponding topological space.  
The chapter is written for a topologist and we assume that all topological 
 facts  are  known to the reader;  for otherwise, a reference 
list is compiled at the end of each section.

 \index{surface map}

\section{Classification of  the  surface  maps}
We assume that $X$ is a compact oriented surface of genus $g\ge 1$;  
we shall be interested in the continuous invertible self-maps (automorphisms)
of $X$,  i.e.
\displaymath
\phi: X\to X.
\enddisplaymath
As it was shown in  the model example for $X\cong T^2$,    there exists a
functor on the set of all Anosov's maps $\phi$ with values in the category
of noncommutative tori  with real multiplication;   the functor sends the 
conjugate Anosov's maps to the stably isomorphic (Morita equivalent)
noncommutative tori ${\cal A}_{\theta}$.  Roughly speaking,  in this section we extend this
result to the higher genus surfaces,  i.e  for $g\ge 2$.   However, instead of 
using ${\cal A}_{\theta}$'s  as the target category,  we shall use the category
of stationary AF-algebras introduced in Section 3.5.2;  in the case $g=1$ 
the two categories are order-isomorphic because their $K_0^+$ semigroups
are,  see the end of Section 3.5.1. 
  
 \index{pseudo-Anosov map}

\subsection{Pseudo-Anosov maps of a surface}
Let $Mod~(X)$ be the mapping class group of a compact  surface $X$, i.e. the 
group of orientation preserving automorphisms of $X$ modulo the trivial ones.
Recall that $\phi,\phi'\in Mod~(X)$ are  conjugate automorphisms,  whenever $\phi'=h\circ\phi\circ h^{-1}$
for an  $h\in Mod~(X)$. It is not hard to see  that  conjugation is an equivalence relation 
which splits the mapping class group  into  disjoint classes of conjugate automorphisms.  
The construction of invariants of the conjugacy classes in $Mod~(X)$ is an important and difficult 
problem studied by [Hemion 1979] \cite{Hem1},  [Mosher 1986]  \cite{Mos1},  and others;
it is important to understand that any knowledge of  such invariants leads to  a topological classification 
of  three-dimensional  manifolds  [Thurston 1982]  \cite{Thu2}.   
It is known that any $\phi\in Mod ~(X)$ is isotopic to an automorphism
$\phi'$, such that either (i) $\phi'$ has a finite order, or
(ii) $\phi'$ is a {\it pseudo-Anosov}  (aperiodic) automorphism, or else
(iii) $\phi'$ is reducible by a system of curves $\Gamma$ surrounded
 by the small tubular neighborhoods $N(\Gamma)$, such that on
 $X \setminus  N(\Gamma)$ $\phi'$ satisfies either (i) or (ii).
Let $\phi$ be a representative of the equivalence class 
of a pseudo-Anosov automorphism. 
Then there exist a pair consisting of the stable ${\cal F}_s$
and unstable ${\cal F}_u$ mutually orthogonal measured foliations on the surface $X$,
such that $\phi({\cal F}_s)={1\over\lambda_{\phi}}{\cal F}_s$ 
and $\phi({\cal F}_u)=\lambda_{\phi}{\cal F}_u$, where $\lambda_{\phi}>1$
is called a dilatation of $\phi$. The foliations ${\cal F}_s,{\cal F}_u$ are minimal,
uniquely ergodic and describe the automorphism $\phi$ up to a power.  
In the sequel,  we shall focus on the conjugacy problem for the pseudo-Anosov 
automorphisms of a surface $X$;   we shall try to solve  the problem using functors with
values in the NCG.
Namely,  we shall    assign to each pseudo-Anosov map  $\phi$ an AF-algebra,  ${\Bbb A}_{\phi}$,
so that for  every $h\in Mod~(X)$  the  diagram in Fig. 4.1 is  commutative.
In  words, if $\phi$ and  $\phi'$ are  conjugate pseudo-Anosov automorphisms,  then
the AF-algebras ${\Bbb A}_{\phi}$ and  ${\Bbb A}_{\phi'}$ are stably isomorphic.  For the sake of 
clarity,  we  shall consider  an example illustrating  the idea in the  case  $X\cong T^2$.    
\begin{figure}[here]
\begin{picture}(300,110)(-120,-5)
\put(20,70){\vector(0,-1){35}}
\put(130,70){\vector(0,-1){35}}
\put(45,23){\vector(1,0){53}}
\put(45,83){\vector(1,0){53}}
\put(15,20){${\Bbb A}_{\phi}$}
\put(128,20){${\Bbb A}_{\phi'}$}
\put(17,80){$\phi$}
\put(117,80){$\phi'=h\circ\phi\circ h^{-1}$}
\put(60,30){\sf stable}
\put(50,10){\sf isomorphism}
\put(54,90){\sf conjugacy}
\end{picture}
\caption{Conjugation of the pseudo-Anosov maps.}
\end{figure}
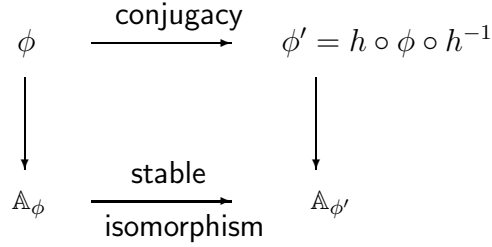

\begin{exm}
{\bf (case $X\cong T^2$)}
\textnormal{
Let $\phi\in Mod~(T^2)$ be the Anosov automorphism given by a non-negative matrix 
$A_{\phi}\in SL(2, {\Bbb Z})$. 
 Consider a stationary AF-algebra,  ${\Bbb A}_{\phi}$,  given by the periodic Bratteli diagram
 shown in Fig. 4.2,  where $a_{ij}$ indicate the multiplicity of the respective edges of the graph. 
(We encourage  the reader to verify that $F: \phi\mapsto {\Bbb A}_{\phi}$
is a well-defined function on the set of Anosov automorphisms given by
the hyperbolic matrices with the non-negative entries.)
Let us show that if $\phi,\phi'\in Mod~(T^2)$ are  conjugate Anosov automorphisms,
then ${\Bbb A}_{\phi},{\Bbb A}_{\phi'}$ are  stably isomorphic
AF-algebras.  Indeed, let $\phi'=h\circ\phi\circ h^{-1}$ for
an $h\in Mod~(X)$. Then $A_{\phi'}=TA_{\phi}T^{-1}$ for
a matrix $T\in SL_2({\Bbb Z})$. Note that $(A_{\phi}')^n=(TA_{\phi}T^{-1})^n=
TA_{\phi}^nT^{-1}$, where $n\in {\Bbb N}$. We shall use the following
criterion:  the AF-algebras ${\Bbb A},{\Bbb A}'$
are stably isomorphic if and only if their Bratteli diagrams contain a 
common block of an arbitrary length, see  [Effros   1981]  \cite{E}, Theorem 2.3
and   recall  that an order-isomorphism mentioned in the theorem is equivalent
to the condition that the corresponding Bratteli diagrams have the same
infinite tails -- i.e. a common block of infinite length.  Consider the following sequences  of 
matrices
\displaymath
\left\{
\begin{array}{c}
\underbrace{A_{\phi}A_{\phi}\dots A_{\phi}}_n\\
T\underbrace{A_{\phi}A_{\phi}\dots A_{\phi}}_nT^{-1},
\end{array}
\right.
\enddisplaymath
which mimic the Bratteli diagrams of ${\Bbb A}_{\phi}$ and ${\Bbb A}_{\phi'}$.
Letting  $n\to\infty$, we conclude that  
${\Bbb A}_{\phi}\otimes {\cal K}\cong {\Bbb A}_{\phi'}\otimes {\cal K}$.
}
\end{exm}
\begin{figure}[here]
\begin{picture}(350,100)(60,0)
\put(97,47){$\bullet$}

\put(117,27){$\bullet$}
\put(117,67){$\bullet$}

\put(157,27){$\bullet$}
\put(157,67){$\bullet$}

\put(197,27){$\bullet$}
\put(197,67){$\bullet$}

\put(237,27){$\bullet$}
\put(237,67){$\bullet$}


\put(100,50){\line(1,1){20}}
\put(100,50){\line(1,-1){20}}

\put(120,30){\line(1,0){40}}

\put(120,32){\line(1,1){40}}
\put(120,28){\line(1,1){40}}

\put(120,72){\line(1,0){40}}
\put(120,70){\line(1,0){40}}
\put(120,68){\line(1,0){40}}

\put(120,70){\line(1,-1){40}}


\put(160,30){\line(1,0){40}}

\put(160,32){\line(1,1){40}}
\put(160,28){\line(1,1){40}}

\put(160,72){\line(1,0){40}}
\put(160,70){\line(1,0){40}}
\put(160,68){\line(1,0){40}}

\put(160,70){\line(1,-1){40}}


\put(200,30){\line(1,0){40}}

\put(200,32){\line(1,1){40}}
\put(200,28){\line(1,1){40}}

\put(200,72){\line(1,0){40}}
\put(200,70){\line(1,0){40}}
\put(200,68){\line(1,0){40}}

\put(200,70){\line(1,-1){40}}

\put(250,30){$\dots$}
\put(250,70){$\dots$}

\put(137,78){$a_{11}$}
\put(177,78){$a_{11}$}
\put(217,78){$a_{11}$}


\put(116,55){$a_{12}$}
\put(156,55){$a_{12}$}
\put(196,55){$a_{12}$}


\put(113,44){$a_{21}$}
\put(153,44){$a_{21}$}
\put(193,44){$a_{21}$}


\put(137,22){$a_{22}$}
\put(177,22){$a_{22}$}
\put(217,22){$a_{22}$}

\put(290,50){
$A_{\phi}=\left(\matrix{a_{11} & a_{12}\cr a_{21} & a_{22}}\right)$,}

\end{picture}

\caption{The AF-algebra  ${\Bbb A}_{\phi}$.}
\end{figure}
\begin{rmk}
{\bf (Handelman's  invariant of the AF-algebra  ${\Bbb A}_{\phi}$)}
\textnormal{
One can reformulate the conjugacy  problem for the  automorphisms $\phi: T^2\to T^2$
in terms of  the AF-algebras  ${\Bbb A}_{\phi}$;   namely,   one needs to  
find  invariants of  the  stable isomorphism (Morita equivalence) classes  of the stationary
AF-algebras  ${\Bbb A}_{\phi}$.  
One such  invariant was  introduced in Section 1.4;  let us recall its definition and properties.   
Consider an eigenvalue problem for the matrix $A_{\phi}\in SL(2, {\Bbb Z})$,  i.e. 
$A_{\phi}v_A=\lambda_Av_A$,  where $\lambda_A>1$ is the Perron-Frobenius eigenvalue  and 
$v_A=(v_A^{(1)},v_A^{(2)})$ the corresponding eigenvector with the positive entries 
normalized so that $v_A^{(i)}\in K={\Bbb Q}(\lambda_A)$. 
Denote by ${\goth m}={\Bbb Z}v_A^{(1)}+{\Bbb Z}v_A^{(2)}$
 the  ${\Bbb Z}$-module in the number field $K$.   The coefficient
ring, $\Lambda$, of module ${\goth m}$ consists of the elements $\alpha\in K$
such that $\alpha {\goth m}\subseteq {\goth m}$. It is known that 
$\Lambda$ is an order in $K$ (i.e. a subring of $K$
containing $1$) and, with no restriction, one can assume that 
${\goth m}\subseteq\Lambda$. It follows  from the definition, that ${\goth m}$
coincides with an ideal, $I$, whose equivalence class in $\Lambda$ we shall denote
by $[I]$.  The triple $(\Lambda, [I], K)$ is an arithmetic invariant of the 
stable isomorphism class of ${\Bbb A}_{\phi}$: the ${\Bbb A}_{\phi},{\Bbb A}_{\phi'}$
are stably isomorphic AF-algebras if and only if $\Lambda=\Lambda', [I]=[I']$ and $K=K'$,
see [Handelman 1981]  \cite{Han2}. 
}
\end{rmk}

 \index{Handelman invariant}

\subsection{Functors and invariants}
Denote by ${\cal F}_{\phi}$ the stable foliation of a pseudo-Anosov automorphism
$\phi\in Mod~(X)$. For brevity, we assume that ${\cal F}_{\phi}$ is an oriented 
foliation given by the trajectories of a closed $1$-form $\omega\in H^1(X; {\Bbb R})$.
Let $v^{(i)}=\int_{\gamma_i}\omega$, where $\{\gamma_1,\dots,\gamma_n\}$ is a basis
in the relative homology $H_1(X, \Sing~{\cal F}_{\phi}; {\Bbb Z})$, such that
$\theta=(\theta_1,\dots,\theta_{n-1})$ is a vector with  positive coordinates 
$\theta_i=v^{(i+1)} / v^{(1)}$. 
\begin{rmk}
\textnormal{
The constants $\theta_i$ depend on a basis in the  homology group,  but the 
${\Bbb Z}$-module generated by the  $\theta_i$ does not. 
}
\end{rmk}
Consider the infinite Jacobi-Perron continued  fraction of $\theta$:
\displaymath
\left(\matrix{1\cr \theta}\right)=
\lim_{k\to\infty} \left(\matrix{0 & 1\cr I & b_1}\right)\dots
\left(\matrix{0 & 1\cr I & b_k}\right)
\left(\matrix{0\cr {\Bbb I}}\right),
\enddisplaymath
where $b_i=(b^{(i)}_1,\dots, b^{(i)}_{n-1})^T$ is a vector of the nonnegative integers,  
$I$ the unit matrix and ${\Bbb I}=(0,\dots, 0, 1)^T$;   we refer the reader to  [Bernstein 1971] 
 \cite{BE}  for the definition of  the Jacobi-Perron algorithm and related  
 continued  fractions.  
 \index{Jacobi-Perron continued fraction}
\begin{dfn}
By  ${\Bbb A}_{\phi}$ one understands  the   AF-algebra given by the 
Bratteli diagram  defined by  the incidence matrices 
$B_k=\left(\small\matrix{0 & 1\cr I & b_k}\right)$ for  $k=1,\dots, \infty$.
\end{dfn}
\begin{rmk}
\textnormal{
We encourage the reader to verify,  that   ${\Bbb A}_{\phi}$  coincides with the one
for the Anosov maps.  (Hint:  the Jacobi-Perron fractions  of  dimension $n=2$ coincide 
with the regular continued fractions.)   
}
\end{rmk}
 \index{Anosov map}
\begin{dfn}
For a matrix $A\in GL_n({\Bbb Z})$ with positive
entries,  we shall  denote by $\lambda_A$  the Perron-Frobenius eigenvalue
and let $(v^{(1)}_A,\dots, v^{(n)}_A)$ be  the corresponding normalized eigenvector
such that    $v^{(i)}_A\in K={\Bbb Q}(\lambda_A)$.
The coefficient (endomorphism) ring of the module   ${\goth m}={\Bbb Z}v^{(1)}_A+\dots+{\Bbb Z}v^{(n)}_A$ 
will shall write as  $\Lambda$;   the equivalence class of ideal $I$ in $\Lambda$
 will be written as   $[I]$.  We shall denote by $\Delta=\dete~(a_{ij})$ and $\Sigma$ the determinant
and signature of the symmetric bilinear form $q(x,y)=\sum_{i,j}^na_{ij}x_ix_j$,
where $a_{ij}=\tr~(v^{(i)}_Av^{(j)}_A)$ and $\tr~ (\bullet)$ 
the trace function.  
\end{dfn}
\begin{thm}\label{thm4.1.1}
${\Bbb A}_{\phi}$ is a stationary AF-algebra. 
\end{thm}
Let $\Phi$ be a category of all pseudo-Anosov (Anosov, resp.) automorphisms
of a surface of the genus $g\ge 2$ ($g=1$,  resp.);  the arrows (morphisms)
are conjugations between the automorphisms. Likewise,
let ${\cal A}$ be the category of all  stationary  AF-algebras ${\Bbb A}_{\phi}$,
where  $\phi$ runs over the set $\Phi$; the arrows of ${\cal A}$ are  stable 
isomorphisms among the algebras ${\Bbb A}_{\phi}$. 
\begin{thm}\label{thm4.1.2}
{\bf (Functor on pseudo-Anosov maps)}
Let  $F:\Phi\to {\cal A}$  be a map given by the formula $\phi\mapsto {\Bbb A}_{\phi}$.   Then:

\medskip
(i) $F$ is a functor  which  maps  conjugate pseudo-Anosov automorphisms to  stably
isomorphic AF-algebras;

\smallskip
(ii) $Ker~F=[\phi]$, where $[\phi]=\{\phi'\in \Phi~|~(\phi')^m=\phi^n, ~m,n\in {\Bbb N}\}$
is the commensurability class of the pseudo-Anosov automorphism $\phi$.  
\end{thm}
\begin{cor}\label{cor4.1.1}
{\bf (Noncommutative invariants)}
The following are invariants of the conjugacy classes of the pseudo-Anosov 
automorphisms:

\medskip
(i) triples  $(\Lambda, [I], K)$;

\smallskip
(ii) integers $\Delta$ and $\Sigma$.  
\end{cor}
\begin{rmk}
{\bf 
(Effectiveness of invariants  $(\Lambda, [I], K)$, $\Delta$ and $\Sigma$)}
\textnormal{
How to calculate invariants  $(\Lambda, [I], K)$, $\Delta$ and $\Sigma$?
There is no obvious  way;  the problem is similar  to that of numerical invariants 
of the fundamental group of a knot. A step in this direction would be 
computation of the matrix $A$; the latter is similar  to the matrix $\rho(\phi)$,
where $\rho: Mod~(X)\to PIL$ is a faithful representation of the mapping class
group as a group of the piecewise-integral-linear (PIL) transformations
[Penner  1984]   \cite{Pen1},  p.45. The entries of $\rho(\phi)$ are the linear combinations 
of the Dehn twists along the $(3g-1)$ (Lickorish) curves on the surface $X$. 
Then one can effectively determine whether the $\rho(\phi)$ and $A$ are
similar matrices (over ${\Bbb Z}$) by bringing the polynomial matrices
$\rho(\phi)-xI$ and $A-xI$ to the Smith normal form; when the similarity
is established, the numerical invariants $\Delta$ and $\Sigma$ become the
polynomials in the Dehn twists. A tabulation of the simplest elements 
of $Mod~(X)$ is possible in terms of $\Delta$ and $\Sigma$.  
}
\end{rmk}

\bigskip\noindent
Theorems \ref{thm4.1.1},  \ref{thm4.1.2} and Corollary \ref{cor4.1.1} 
will be proved in Section 4.1.5;  the necessary background is developed 
in the sections  below.

 \index{Jacobian of measured foliation}
 \index{measured foliation}
\subsection{Jacobian of  measured foliations}
Let ${\cal F}$ be a measured foliation on a compact surface
$X$ \cite{Thu1}. For the sake of brevity, we shall always assume
that ${\cal F}$ is an oriented foliation, i.e. given by the 
trajectories of a closed $1$-form $\omega$ on $X$. (The assumption
is not a  restriction -- each measured foliation is oriented on a 
surface $\widetilde X$, which is a double cover of $X$ ramified 
at the singular points of the half-integer index of the non-oriented
foliation [Hubbard \& Masur  1979]   \cite{HuMa1}.)  Let $\{\gamma_1,\dots,\gamma_n\}$ be a
basis in the relative homology group $H_1(X, \Sing~{\cal F}; {\Bbb Z})$,
where $\Sing~{\cal F}$ is the set of singular points of the foliation ${\cal F}$.
It is well known that $n=2g+m-1$, where $g$ is the genus of $X$ and $m= |\Sing~({\cal F})|$. 
The periods of $\omega$ in the above basis will be written
\displaymath
\lambda_i=\int_{\gamma_i}\omega.
\enddisplaymath
The real numbers $\lambda_i$ are coordinates of ${\cal F}$ in the space of all 
measured foliations on $X$ with the fixed set of  singular points,  see e.g.
[Douady \& Hubbard  1975]  \cite{DoHu1}.  
\begin{dfn}
By a jacobian $Jac~({\cal F})$ of the measured foliation ${\cal F}$, we
understand a ${\Bbb Z}$-module ${\goth m}={\Bbb Z}\lambda_1+\dots+{\Bbb Z}\lambda_n$
regarded as a subset of the real line ${\Bbb R}$. 
\end{dfn}
An importance of the jacobians stems from an observation that
although the periods, $\lambda_i$,  depend on the choice of basis in 
$H_1(X, \Sing~{\cal F}; {\Bbb Z})$,  the jacobian does not.
Moreover,  up to a scalar multiple,   the jacobian is an invariant
of the equivalence class of the foliation ${\cal F}$. 
We formalize these observations in the following two lemmas. 
\begin{lem}\label{lm4.1.1}
The ${\Bbb Z}$-module ${\goth m}$  is independent of choice of 
basis
\linebreak
in $H_1(X, \Sing~{\cal F}; {\Bbb Z})$ and depends solely on the foliation ${\cal F}$.
\end{lem}
{\it Proof.} Indeed, let $A=(a_{ij})\in GL(n, {\Bbb Z})$ and let
$$
\gamma_i'=\sum_{j=1}^na_{ij}\gamma_j
$$
be a new basis in $H_1(X, \Sing~{\cal F}; {\Bbb Z})$. 
Then using the integration rules:
\begin{eqnarray}
\lambda_i'  &= \int_{\gamma_i'}\omega &= \int_{\sum_{j=1}^na_{ij}\gamma_j}\omega=\nonumber \\
 &= \sum_{j=1}^n\int_{\gamma_j}\omega  &= \sum_{j=1}^na_{ij}\lambda_j.\nonumber 
\end{eqnarray}

\bigskip
To prove that ${\goth m}={\goth m}'$, consider the following equations:
\begin{eqnarray}
{\goth m}'  &= \sum_{i=1}^n{\Bbb Z}\lambda_i' &= \sum_{i=1}^n {\Bbb Z} \sum_{j=1}^n a_{ij}\lambda_j=\nonumber \\
 &= \sum_{j=1}^n \left(\sum_{i=1}^n a_{ij}{\Bbb Z}\right)\lambda_j  &\subseteq  {\goth m}. \nonumber
\end{eqnarray}
Let $A^{-1}=(b_{ij})\in GL(n, {\Bbb Z})$ be an inverse to the matrix $A$.
Then $\lambda_i=\sum_{j=1}^nb_{ij}\lambda_j'$ and 
\begin{eqnarray}
{\goth m}  &= \sum_{i=1}^n{\Bbb Z}\lambda_i &= \sum_{i=1}^n {\Bbb Z} \sum_{j=1}^n b_{ij}\lambda_j'=\nonumber \\
 &= \sum_{j=1}^n \left(\sum_{i=1}^n b_{ij}{\Bbb Z}\right)\lambda_j'  &\subseteq  {\goth m}'.\nonumber 
\end{eqnarray}
Since both ${\goth m}'\subseteq {\goth m}$ and ${\goth m}\subseteq {\goth m}'$, we conclude
that ${\goth m}' = {\goth m}$. Lemma \ref{lm4.1.1} follows.
$\square$

\begin{dfn}
Two  measured foliations ${\cal F}$ and ${\cal F}'$ are said to 
{\it equivalent},  if there exists an automorphism $h\in Mod~(X)$,
which sends the leaves of the foliation ${\cal F}$ to the leaves of the 
foliation ${\cal F}'$.
\end{dfn}
\begin{rmk}
\textnormal{
The  equivalence relation involves  the topological foliations,  i.e.  projective
classes of the measured foliations,  see  [Thurston 1988]  \cite{Thu1} for 
the details. 
}
\end{rmk}
\begin{lem}\label{lm4.1.2}
Let ${\cal F}, {\cal F}'$  be the equivalent measured foliations 
on a surface $X$. Then 
\displaymath
Jac~({\cal F}')=\mu ~Jac~({\cal F}),
\enddisplaymath
where $\mu>0$ is a real number.  
\end{lem}
{\it Proof.} 
Let $h: X\to X$ be an automorphism of the surface $X$. Denote
by $h_*$ its action on $H_1(X, \Sing~({\cal F}); {\Bbb Z})$
and by $h^*$ on $H^1(X; {\Bbb R})$ connected  by the formula: 
$$
\int_{h_*(\gamma)}\omega=\int_{\gamma}h^*(\omega), \qquad\forall\gamma\in H_1(X, \Sing~({\cal F}); {\Bbb Z}), 
\qquad\forall\omega\in H^1(X; {\Bbb R}).
$$
Let $\omega,\omega'\in H^1(X; {\Bbb R})$ be the closed $1$-forms whose
trajectories define the foliations ${\cal F}$ and ${\cal F}'$, respectively.
Since ${\cal F}, {\cal F}'$ are equivalent measured foliations,
$$
\omega'= \mu ~h^*(\omega)
$$
for a $\mu>0$.

Let $Jac~({\cal F})={\Bbb Z}\lambda_1+\dots+{\Bbb Z}\lambda_n$ and 
$Jac~({\cal F}')={\Bbb Z}\lambda_1'+\dots+{\Bbb Z}\lambda_n'$. Then:
$$
\lambda_i'=\int_{\gamma_i}\omega'=\mu~\int_{\gamma_i}h^*(\omega)=
\mu~\int_{h_*(\gamma_i)}\omega, \qquad 1\le i\le n.
$$
By lemma \ref{lm4.1.1},   it holds:
$$
Jac~({\cal F})=\sum_{i=1}^n{\Bbb Z}\int_{\gamma_i}\omega=
\sum_{i=1}^n{\Bbb Z}\int_{h_*(\gamma_i)}\omega.
$$
Therefore:
$$
Jac~({\cal F}')=\sum_{i=1}^n{\Bbb Z}\int_{\gamma_i}\omega'=
\mu~\sum_{i=1}^n{\Bbb Z}\int_{h_*(\gamma_i)}\omega=\mu~Jac~({\cal F}).
$$
Lemma \ref{lm4.1.2} follows.
$\square$

 \index{equivalent foliations}
\subsection{Equivalent foliations} 
Recall that for a measured foliation  ${\cal F}$,   we constructed 
 an AF-algebra, ${\Bbb A}_{\cal F}$.    Our goal is to prove   
commutativity of the  diagram in Fig.  4.3.;    in   other words, 
two equivalent measured foliations map  to the stably isomorphic 
(Morita equivalent) AF-algebras   ${\Bbb A}_{\cal F}$.  
\begin{figure}[here]
\begin{picture}(300,110)(-120,-5)
\put(20,70){\vector(0,-1){35}}
\put(130,70){\vector(0,-1){35}}
\put(45,23){\vector(1,0){53}}
\put(45,83){\vector(1,0){53}}
\put(15,20){${\Bbb A}_{\cal F}$}
\put(128,20){${\Bbb A}_{{\cal F}'}$}
\put(17,80){${\cal F}$}
\put(125,80){${\cal F}'$}
\put(60,30){\sf stably}
\put(50,10){\sf isomorphic}

\put(50,90){\sf equivalent}
\end{picture}
\caption{Functor on measured foliations.}
\end{figure}
\begin{lem}\label{lm4.1.3}
{\bf (Perron)}
Let ${\goth m}={\Bbb Z}\lambda_1+\dots+{\Bbb Z}\lambda_n$
and   ${\goth m}'={\Bbb Z}\lambda_1'+\dots+{\Bbb Z}\lambda_n'$
be two ${\Bbb Z}$-modules, such that ${\goth m}'=\mu {\goth m}$ for a $\mu>0$. 
Then the Jacobi-Perron continued fractions of the vectors $\lambda$ and $\lambda'$
coincide except,   possibly,  at a finite number of terms. 
\end{lem}
{\it Proof.}
Let ${\goth m}={\Bbb Z}\lambda_1+\dots+{\Bbb Z}\lambda_n$ and 
${\goth m}'={\Bbb Z}\lambda_1'+\dots+{\Bbb Z}\lambda_n'$. Since
${\goth m}'=\mu {\goth m}$, where $\mu$ is a positive real,
one gets the following identity of the ${\Bbb Z}$-modules:
\displaymath
{\Bbb Z}\lambda_1'+\dots+{\Bbb Z}\lambda_n'={\Bbb Z}(\mu\lambda_1)+\dots+{\Bbb Z}(\mu\lambda_n).
\enddisplaymath
One can always assume that $\lambda_i$ and $\lambda_i'$ are positive reals.
For obvious reasons, there exists a basis $\{\lambda_1^{''},\dots,\lambda_n^{''}\}$
of the module ${\goth m}'$, such that:
\displaymath
\left\{
\begin{array}{cc}
\lambda'' &= A(\mu\lambda) \nonumber\\
\lambda'' &= A'\lambda',
\end{array}
\right.
\enddisplaymath
where $A,A'\in GL^+(n, {\Bbb Z})$ are the matrices, whose entries 
are non-negative integers.  In view of [Bauer  1996]  \cite{Bau1}, 
 Proposition 3,  we have 
\displaymath
\left\{
\begin{array}{cc}
A &=  \left(\matrix{0 & 1\cr I & b_1}\right)\dots
\left(\matrix{0 & 1\cr I & b_k}\right)\nonumber\\
A' &= \left(\matrix{0 & 1\cr I & b_1'}\right)\dots
\left(\matrix{0 & 1\cr I & b_l'}\right),
\end{array}
\right.
\enddisplaymath
where $b_i, b_i'$ are non-negative integer vectors.
Since the  Jacobi-Perron continued fraction for the vectors
$\lambda$ and $\mu\lambda$ coincide for any $\mu>0$ 
(see e.g.  [Bernstein  1971]  \cite{BE}),   we conclude that: 
\displaymath
\left\{
\begin{array}{cc}
\left(\matrix{1\cr \theta}\right)
 &=  \left(\matrix{0 & 1\cr I & b_1}\right)\dots
\left(\matrix{0 & 1\cr I & b_k}\right)
\left(\matrix{0 & 1\cr I & a_1}\right)
\left(\matrix{0 & 1\cr I & a_2}\right)\dots
\left(\matrix{0\cr {\Bbb I}}\right)
\nonumber\\
\left(\matrix{1\cr \theta'}\right)
 &= \left(\matrix{0 & 1\cr I & b_1'}\right)\dots
\left(\matrix{0 & 1\cr I & b_l'}\right)
\left(\matrix{0 & 1\cr I & a_1}\right)
\left(\matrix{0 & 1\cr I & a_2}\right)\dots
\left(\matrix{0\cr {\Bbb I}}\right),
\end{array}
\right.
\enddisplaymath
where 
\displaymath
\left(\matrix{1\cr \theta''}\right)=
\lim_{i\to\infty} \left(\matrix{0 & 1\cr I & a_1}\right)\dots
\left(\matrix{0 & 1\cr I & a_i}\right)
\left(\matrix{0\cr {\Bbb I}}\right). 
\enddisplaymath
In other words, the continued fractions of the vectors $\lambda$ and $\lambda'$
coincide  but at a finite number of terms.    Lemma \ref{lm4.1.3}  follows. 
$\square$

\begin{lem}\label{lm4.1.4}
{\bf (Basic lemma)}
Let ${\cal F}$ and  $ {\cal F}'$ be  equivalent measured foliations on a surface $X$.
Then the AF-algebras ${\Bbb A}_{\cal F}$ and  ${\Bbb A}_{{\cal F}'}$
are stably isomorphic.  
\end{lem}
{\it Proof.}
Notice that lemma \ref{lm4.1.2} implies that
 equivalent measured foliations ${\cal F}, {\cal F}'$ have
 proportional jacobians, i.e. ${\goth m}'=\mu {\goth m}$
for a $\mu>0$.  On the other hand, by  lemma \ref{lm4.1.3}
the continued fraction expansion of the basis vectors
of the proportional jacobians must coincide, except
a finite number of terms. Thus, the AF-algebras
${\Bbb A}_{\cal F}$ and ${\Bbb A}_{{\cal F}'}$ 
are given by  the Bratteli diagrams, which are identical,
except a finite part of the diagram.  
It is well  known (see e.g. [Effros 1981]  \cite{E}, Theorem 2.3),   that the AF-algebras,  
which have such a  property,   are stably isomorphic.    Lemma \ref{lm4.1.4} follows. 
$\square$

\subsection{Proofs} 
\subsubsection{Proof of theorem \ref{thm4.1.1}} 
Let $\phi\in Mod~(X)$ be a pseudo-Anosov automorphism of the surface $X$.
Denote by ${\cal F}_{\phi}$ the invariant foliation of $\phi$.
By definition of such a foliation, $\phi({\cal F}_{\phi})=\lambda_{\phi}{\cal F}_{\phi}$,
where $\lambda_{\phi}>1$ is the dilatation of $\phi$.
Consider the jacobian $Jac~({\cal F}_{\phi})={\goth m}_{\phi}$
of  foliation ${\cal F}_{\phi}$. 
Since ${\cal F}_{\phi}$ is an invariant foliation of the pseudo-Anosov automorphism $\phi$, 
one gets the following equality  of the ${\Bbb Z}$-modules:
\displaymath
{\goth m}_{\phi}=\lambda_{\phi}{\goth m}_{\phi}, \qquad \lambda_{\phi}\ne\pm 1.
\enddisplaymath 
Let $\{v^{(1)},\dots,v^{(n)}\}$ be a basis in  module  ${\goth m}_{\phi}$,
such that $v^{(i)}>0$;  from the above equation,  one obtains the following
system of linear equations:
\displaymath
\left\{
\begin{array}{ccc}
\lambda_{\phi}v^{(1)}  &=& a_{11}v^{(1)}+a_{12}v^{(2)}+\dots+a_{1n}v^{(n)}\\
\lambda_{\phi}v^{(2)}  &=& a_{21}v^{(1)}+a_{22}v^{(2)}+\dots+a_{2n}v^{(n)}\\
\vdots && \\
\lambda_{\phi}v^{(n)}  &=& a_{n1}v^{(1)}+a_{n2}v^{(2)}+\dots+a_{nn}v^{(n)},
\end{array}
\right.
\enddisplaymath
where $a_{ij}\in {\Bbb Z}$. The matrix $A=(a_{ij})$ is invertible. Indeed,
since  foliation ${\cal F}_{\phi}$ is minimal, real numbers  $v^{(1)},\dots,v^{(n)}$ 
are linearly independent over ${\Bbb Q}$. So do numbers $\lambda_{\phi}v^{(1)},\dots,\lambda_{\phi}v^{(n)}$,
which therefore can be taken for a basis of the module ${\goth m}_{\phi}$. 
Thus,  there exists an integer matrix $B=(b_{ij})$, such that $v^{(j)}=\sum_{i,j}w^{(i)}$,
where $w^{(i)}=\lambda_{\phi}v^{(i)}$. Clearly, $B$ is an inverse to  matrix $A$.
Therefore, $A\in GL(n, {\Bbb Z})$.

Moreover, without loss of the generality one can assume that $a_{ij}\ge0$. 
Indeed, if it is not yet the case, consider the conjugacy class $[A]$
of the matrix $A$. Since $v^{(i)}>0$,  there exists a matrix $A^+\in [A]$
whose entries are non-negative integers. One has to replace 
 basis $v=(v^{(1)},\dots,v^{(n)})$ in the module ${\goth m}_{\phi}$ 
by a  basis  $Tv$, where $A^+=TAT^{-1}$. It will be further assumed that
 $A=A^+$. 
\begin{lem}\label{lm4.1.5}
Vector $(v^{(1)}, \dots, v^{(n)})$ is the  limit of a periodic  
Jacobi-Perron continued fraction. 
\end{lem}
{\it Proof.} It follows from the discussion above,  that there exists
a non-negative integer matrix $A$,  such that $Av=\lambda_{\phi}v$. 
In view of [Bauer 1996]  \cite{Bau1},  Proposition 3,  matrix $A$ admits 
the unique factorization
\displaymath
A=
\left(\matrix{0 & 1\cr I & b_1}\right)\dots
\left(\matrix{0 & 1\cr I & b_k}\right),
\enddisplaymath
where $b_i=(b^{(i)}_1,\dots, b^{(i)}_{n})^T$ are vectors of the  non-negative integers.
Let us consider the  periodic Jacobi-Perron continued fraction
\displaymath
Per
~\overline{
 \left(\matrix{0 & 1\cr I & b_1}\right)\dots
\left(\matrix{0 & 1\cr I & b_k}\right)
}
\left(\matrix{0\cr {\Bbb I}}\right).
\enddisplaymath
According to [Perron 1907]  \cite{Per1}, {\bf Satz XII},  the above  fraction converges to vector $w=(w^{(1)},\dots, w^{(n)})$,
such that $w$ satisfies equation $(B_1B_2\dots B_k)w=Aw=\lambda_{\phi}w$.
In view of  equation $Av=\lambda_{\phi}v$, we conclude that  vectors $v$ and $w$
are collinear.  Therefore, the Jacobi-Perron continued fractions of $v$ and $w$
must coincide.  Lemma \ref{lm4.1.5} follows.  
$\square$

\bigskip\noindent
It is easy to see,  that the AF-algebra attached to  foliation
${\cal F}_{\phi}$ is stationary.   Indeed, by lemma \ref{lm4.1.5}, 
the vector of periods $v^{(i)}=\int_{\gamma_i}\omega$ unfolds  into a periodic Jacobi-Perron 
continued fraction. By  definition, the Bratteli diagram of the AF-algebra  ${\Bbb A}_{\phi}$
is periodic as well.  In other words, the AF-algebra  ${\Bbb A}_{\phi}$ is stationary. Theorem \ref{thm4.1.1}
is proved.  
$\square$

\subsubsection{Proof of theorem \ref{thm4.1.2}} 
(i) Let us prove the first statement. For the sake of completeness,
let us give a proof of the following well-known lemma. 
\begin{lem}\label{lm4.1.6}
If  $\phi$ and $\phi'$ are  conjugate pseudo-Anosov automorphisms
of a surface $X$,   then their invariant measured  foliations ${\cal F}_{\phi}$
and ${\cal F}_{\phi'}$ are equivalent.  
\end{lem}
{\it Proof.} Let $\phi,\phi'\in Mod~(X)$ be conjugate, i.e 
$\phi'=h\circ\phi\circ h^{-1}$ for an automorphism $h\in Mod~(X)$.
Since $\phi$ is the pseudo-Anosov automorphism, there exists  a measured foliation
${\cal F}_{\phi}$,  such that $\phi({\cal F}_{\phi})=\lambda_{\phi}{\cal F}_{\phi}$.
Let us evaluate the automorphism $\phi'$ on the foliation $h({\cal F}_{\phi})$:
\begin{eqnarray}
\phi'(h({\cal F}_{\phi}))  &= h\phi h^{-1}(h({\cal F}_{\phi})) &= 
h\phi({\cal F}_{\phi})=\nonumber \\
 &= h \lambda_{\phi} {\cal F}_{\phi}  &= \lambda_{\phi} (h({\cal F}_{\phi}))\nonumber. 
\end{eqnarray}
Thus, ${\cal F}_{\phi'}=h({\cal F}_{\phi})$ is the invariant foliation for the 
pseudo-Anosov automorphism $\phi'$ and ${\cal F}_{\phi}, {\cal F}_{\phi'}$
are equivalent foliations. Note also that the pseudo-Anosov automorphism $\phi'$ has 
the same dilatation as the automorphism $\phi$.  Lemma \ref{lm4.1.6} follows.  
$\square$

\medskip
To finish the proof of item (i), 
suppose that $\phi$ and $\phi'$ are  conjugate pseudo-Anosov 
automorphisms.  Functor $F$ acts by the formulas $\phi\mapsto {\Bbb A}_{\phi}$
and $\phi'\mapsto {\Bbb A}_{\phi'}$, where ${\Bbb A}_{\phi}, {\Bbb A}_{\phi'}$
are the AF-algebras corresponding to   invariant foliations ${\cal F}_{\phi}, {\cal F}_{\phi'}$. 
In view of lemma \ref{lm4.1.6}, ${\cal F}_{\phi}$ and ${\cal F}_{\phi'}$ are
 equivalent measured foliations.  Then,  by lemma \ref{lm4.1.4}, the AF-algebras ${\Bbb A}_{\phi}$
and ${\Bbb A}_{\phi'}$ are stably isomorphic AF-algebras.  Item (i)
follows.

\bigskip
(ii) Let us prove the second statement. We start with an elementary observation.
Let $\phi\in Mod~(X)$ be a pseudo-Anosov automorphism. Then there exists a unique
measured foliation, ${\cal F}_{\phi}$, such that $\phi({\cal F}_{\phi})=\lambda_{\phi}{\cal F}_{\phi}$,
where $\lambda_{\phi}>1$ is an algebraic integer. Let us evaluate  automorphism
$\phi^2\in Mod~(X)$ on the foliation ${\cal F}_{\phi}$: 
\begin{eqnarray}
\phi^2({\cal F}_{\phi}) &= \phi (\phi({\cal F}_{\phi})) &= 
\phi(\lambda_{\phi} {\cal F}_{\phi})=\nonumber \\
= \lambda_{\phi} \phi({\cal F}_{\phi}) &= \lambda_{\phi}^2{\cal F}_{\phi}  &= 
\lambda_{\phi^2}{\cal F}_{\phi}, \nonumber
\end{eqnarray}
where $\lambda_{\phi^2}:= \lambda_{\phi}^2$. Thus,  foliation ${\cal F}_{\phi}$
is an invariant foliation for the automorphism $\phi^2$ as well. By induction,
one concludes that ${\cal F}_{\phi}$ is an invariant foliation of  the 
automorphism $\phi^n$ for any $n\ge 1$.

Even more is true. Suppose that $\psi\in Mod~(X)$ is a pseudo-Anosov
automorphism, such that $\psi^m=\phi^n$ for some $m\ge 1$ and $\psi\ne\phi$.
Then ${\cal F}_{\phi}$ is an invariant foliation for the automorphism 
$\psi$. Indeed, ${\cal F}_{\phi}$ is  invariant foliation of  the 
automorphism $\psi^m$. If there exists ${\cal F}'\ne {\cal F}_{\phi}$
such that the foliation  ${\cal F}'$ is an invariant foliation  of $\psi$, then 
the foliation ${\cal F}'$ is
also an invariant foliation  of the pseudo-Anosov automorphism $\psi^m$. 
Thus, by the uniqueness, 
${\cal F}'={\cal F}_{\phi}$.  We have just proved the following lemma. 
\begin{lem}\label{lm4.1.7}
If  $[\phi]$ is   the set of all  pseudo-Anosov automorphisms $\psi$
of $X$,  such that $\psi^m=\phi^n$ for some positive integers
$m$ and $n$,  then the pseudo-Anosov foliation ${\cal F}_{\phi}$
is an invariant foliation for every pseudo-Anosov automorphism $\psi\in [\phi]$. 
\end{lem}
In view of lemma \ref{lm4.1.7},  one gets  the following
identities for  the AF-algebras
\displaymath
{\Bbb A}_{\phi}={\Bbb A}_{\phi^2}=\dots={\Bbb A}_{\phi^n}=
{\Bbb A}_{\psi^m}=\dots={\Bbb A}_{\psi^2}={\Bbb A}_{\psi}.
\enddisplaymath
Thus,   functor $F$ is not an 
injective functor: the preimage, $Ker~F$, of  algbera
${\Bbb A}_{\phi}$ consists of a countable set of the pseudo-Anosov
automorphisms $\psi\in [\phi]$, commensurable with the automorphism
$\phi$.  Theorem \ref{thm4.1.2}  is proved.
$\square$

\subsubsection{Proof of corollary \ref{cor4.1.1}}
(i) Theorem \ref{thm4.1.1}  says   that ${\Bbb A}_{\phi}$
is a stationary AF-algebra.   An arithmetic invariant of the stable
isomorphism classes of the stationary AF-algebras has been found 
by D.~Handelman in [Handelman 1981]  \cite{Han1}. Summing up his results, the invariant
is as follows. 
Let $A\in GL(n, {\Bbb Z})$ be a matrix with the strictly positive
entries,  such that $A$ is  equal to the minimal period of the Bratteli diagram
of the stationary AF-algebra. (In case the matrix $A$ has  zero entries,
it is necessary to take a proper minimal power of the matrix $A$.) By the Perron-Frobenius
theory,  matrix $A$ has a real eigenvalue $\lambda_A>1$, which exceeds
the absolute values of other roots of the characteristic polynomial of $A$.
Note that $\lambda_A$ is an invertible algebraic integer (the unit).  Consider
the  real algebraic number field $K={\Bbb Q}(\lambda_A)$ obtained as 
an extension of the field of the rational numbers by the algebraic 
number $\lambda_A$. Let $(v^{(1)}_A,\dots,v^{(n)}_A)$ be the eigenvector
corresponding to the eigenvalue $\lambda_A$. One can  normalize the eigenvector so 
that $v^{(i)}_A\in K$. 
The   departure point  of Handelman's invariant is the  ${\Bbb Z}$-module
${\goth m}={\Bbb Z}v^{(1)}_A+\dots+{\Bbb Z}v^{(n)}_A$. The module ${\goth m}$
brings in two new arithmetic objects: (i) the ring $\Lambda$ of the endomorphisms of ${\goth m}$  
and (ii) an ideal $I$ in the ring $\Lambda$,  such that $I={\goth m}$  after a scaling,
see e.g. [Borevich \& Shafarevich 1966]  \cite{BS},  Lemma 1,  p. 88. 
The ring $\Lambda$ is an order in the algebraic number field $K$ and therefore one can talk about the ideal
classes in $\Lambda$. The ideal class of $I$ is denoted by $[I]$. 
Omitting the embedding question for the field $K$,  the triple $(\Lambda, [I], K)$ is an invariant of
the  stable isomorphism class of the  stationary AF-algebra ${\Bbb A}_{\phi}$,   see  
[Handelman 1981]  \cite{Han1}, \S 5.  
Item (i) follows.

\bigskip
(ii)  Numerical invariants of the stable isomorphism
classes of the stationary AF-algebras can be  derived from the  triple $(\Lambda, [I], K)$.
These invariants are the rational integers -- called the determinant and signature -- can be  obtained as 
follows.  
Let ${\goth m}, {\goth m}'$ be the full ${\Bbb Z}$-modules in an algebraic
number field $K$. It follows from (i), that if ${\goth m}\ne {\goth m}'$
are distinct as the ${\Bbb Z}$-modules,
then the corresponding AF-algebras cannot be stably isomorphic. 
We wish to find the numerical invariants, which discern the case ${\goth m}\ne {\goth m}'$. 
It is assumed that a ${\Bbb Z}$-module is given by the  set of generators
$\{\lambda_1,\dots,\lambda_n\}$. Therefore, the  problem can be formulated as follows: find a number
attached to the set of generators $\{\lambda_1,\dots,\lambda_n\}$,  which does not change   
on the set of generators $\{\lambda_1',\dots,\lambda_n'\}$ of  the same ${\Bbb Z}$-module. 
One such invariant is associated with the trace function on the algebraic number field $K$. 
Recall that 
\displaymath
Tr: K\to {\Bbb Q}
\enddisplaymath
 is a linear function on  $K$ such  that $Tr ~(\alpha+\beta)=\tr~(\alpha)+ \tr~(\beta)$ and
$\tr~(a\alpha)=a ~\tr~(\alpha)$ for $\forall\alpha,\beta\in K$ and
$\forall a\in {\Bbb Q}$.   Let ${\goth m}$ be a full ${\Bbb Z}$-module in the field $K$.   
The trace function defines a symmetric bilinear form  $q(x,y): {\goth m}\times {\goth m}\to {\Bbb Q}$ by 
the formula
\displaymath
(x,y)\longmapsto \tr~(xy), \qquad \forall x,y\in {\goth m}.
\enddisplaymath
The form $q(x,y)$ depends on the basis $\{\lambda_1,\dots,\lambda_n\}$
in the module ${\goth m}$
\displaymath
q(x,y)=\sum_{j=1}^n\sum_{i=1}^na_{ij}x_iy_j, \qquad\hbox{where} \quad a_{ij}=\tr~(\lambda_i\lambda_j).
\enddisplaymath
However, the general theory of the bilinear forms (over the fields ${\Bbb Q}, {\Bbb R}, {\Bbb C}$ or the ring
of rational integers ${\Bbb Z}$)  tells us that certain numerical quantities will not depend on the 
choice of such a basis. 
\begin{dfn}
By a determinant  of the bilinear form  $q(x,y)$ one understands
 the rational integer number 
\displaymath
\Delta= \dete~(\tr~(\lambda_i\lambda_j)).
\enddisplaymath
\end{dfn}
\begin{lem}\label{lm4.1.8}
The determinant $\Delta({\goth m})$  is independent of  the choice of the basis 
$\{\lambda_1,\dots,\lambda_n\}$ in the  module ${\goth m}$. 
\end{lem}
{\it Proof.}   Consider a symmetric matrix $A$  corresponding to the bilinear form $q(x,y)$,
i.e. 
\displaymath
A=\left(
\matrix{a_{11} & a_{12} & \dots & a_{1n}\cr
        a_{12} & a_{22} & \dots & a_{2n}\cr
        \vdots  &        &       & \vdots\cr
        a_{1n} & a_{2n} & \dots & a_{nn}}
\right).
\enddisplaymath
It is known that the matrix $A$, written in a new basis, will take the
form $A'=U^TAU$,   where ~$U\in GL(n, {\Bbb Z})$. 
Then $\dete~(A')=\dete~(U^TAU)=\dete~(U^T) \dete~(A) \dete~(U) =\dete~(A)$.
Therefore, the rational integer number 
\displaymath
\Delta= \dete~(\tr~(\lambda_i\lambda_j)),
\enddisplaymath
does not depend on the  choice of the basis $\{\lambda_1,\dots,\lambda_n\}$ in the
module ${\goth m}$.   Lemma \ref{lm4.1.8} follows.
$\square$

 \index{$p$-adic invariant}
\begin{rmk}
{\bf ($p$-adic invariants)}
\textnormal{
Roughly speaking, Lemma \ref{lm4.1.8} says that determinant  
$\Delta({\goth m})$ discerns  two distinct modules,  i.e.  ${\goth m}\ne {\goth m}'$.
Note that if $\Delta({\goth m})=\Delta({\goth m}')$ for the modules ${\goth m}$ and  
${\goth m}'$,
one cannot conclude that ${\goth m}={\goth m}'$.  The problem of equivalence
of the  symmetric bilinear forms over ${\Bbb Q}$ (i.e. the existence of a linear
substitution over ${\Bbb Q}$,  which transforms one form to the other),
is a fundamental question of number theory. The Minkowski-Hasse theorem
says that two such forms are equivalent if and only if they are
equivalent over the $p$-adic field {\bf Q}$_p$ for every prime number $p$
and over the field ${\Bbb R}$. Clearly, the resulting $p$-adic quantities
will give new invariants of the stable isomorphism classes of the AF-algebras.
The question is much similar to the Minkowski units attached to knots, see e.g.
[Reidemeister  1932]  \cite{R}. 
}
\end{rmk}
\begin{dfn}
By a signature of the bilinear form $q(x,y)$  one understands the rational integer
$\Sigma = (\# a_i^+) - (\# a_i^-)$,    where  $a_i^+$ are  the positive and $a_i^-$ the 
negative entries in the diagonal form
\displaymath
a_1x_1^2+a_2x_2^2+\dots+a_nx_n^2
\enddisplaymath
of  $q(x,y)$;   recall that each  $q(x,y)$ can be brought by an integer 
linear  transformation  to the diagonal form.
\end{dfn}
\begin{lem}\label{lm4.1.9}
The   signature $\Sigma({\goth m})$  is  independent   of the choice of
 basis in the module ${\goth m}$ and, therefore,  $\Sigma({\goth m})\ne \Sigma({\goth m}')$ 
 implies  ${\goth m}\ne {\goth m}'$.  
\end{lem}
{\it Proof.} 
 The claim follows from the Law of Inertia for the signature of the bilinear form $q(x,y)$.
 $\square$

\bigskip
Corollary \ref{cor4.1.1} follows from Lemmas \ref{lm4.1.8} and \ref{lm4.1.9}.
$\square$

 \index{Anosov map}

\subsection{Anosov maps of the torus}
We shall calculate the noncommutative  invariants $\Delta({\goth m})$ and $\Sigma({\goth m})$
for the Anosov automorphisms of the two-dimensional torus;  we construct 
concrete examples    of  Anosov automorphisms with
the same Alexander polynomial $\Delta(t)$   but different invariant  $\Delta({\goth m})$,
i.e. showing that $\Sigma({\goth m})$ is {\it finer} than $\Delta(t)$.  
Recall that isotopy classes of the orientation-preserving diffeomorphisms
of the torus $T^2$ are bijective with the $2\times 2$ matrices with integer entries
and determinant $+1$, i.e. $Mod~ (T^2)\cong SL(2,{\Bbb Z})$. Under the identification,
the non-periodic automorphisms correspond to the matrices $A\in SL(2,{\Bbb Z})$ 
with $|\tr~A|>2$. 
Let $K={\Bbb Q}(\sqrt{d})$ be a quadratic extension of the field of rational numbers ${\Bbb Q}$.
Further we suppose that $d$ is a positive square free integer. Let
\displaymath
\omega=\cases{{1+\sqrt{d}\over 2} & if $d\equiv 1 ~mod~4$,\cr
               \sqrt{d} & if $d\equiv 2,3 ~mod~4$.}
\enddisplaymath
\begin{rmk}
\textnormal{
Recall that  if  $f$ is a positive integer  then every order  in $K$  has the form  
$\Lambda_f={\Bbb Z} +(f\omega){\Bbb Z}$,   where  $f$ is the conductor of $\Lambda_f$,  
see e.g.  [Borevich \& Shafarevich 1966]  \cite{BS},   pp. 130-132. 
This formula   allows to classify the similarity classes of the full modules 
in the field $K$.  Indeed,  there exists a finite number of ${\goth m}_f^{(1)},\dots,
{\goth m}_f^{(s)}$ of the non-similar full modules in the field  $K$  whose coefficient
ring is the order $\Lambda_f$,   see  [Borevich \& Shafarevich 1966]  \cite{BS},  Chapter 2.7, Theorem 3. 
Thus one gets  a finite-to-one classification of the similarity classes of full modules in the field $K$.  
}
\end{rmk}

\subsubsection{Numerical invariants of the Anosov maps}
Let $\Lambda_f$ be an order  in $K$ with
the  conductor $f$. Under the addition operation, the order  $\Lambda_f$ is a full module, 
which we denote by ${\goth m}_f$.  
Let us evaluate  the invariants $q(x,y)$, $\Delta$  and $\Sigma$ 
on the module  ${\goth m}_f$. To calculate  $(a_{ij})=\tr~(\lambda_i\lambda_j)$,  we let 
$\lambda_1=1,\lambda_2=f\omega$. Then:
\begin{eqnarray}
a_{11} &=& 2, \quad a_{12}=a_{21}= f, \quad a_{22}= {1\over 2} f^2(d+1)\quad \hbox{if} \quad d\equiv 1 ~mod~4\nonumber\\
a_{11} &=&  2, \quad a_{12}=a_{21}= 0, \quad a_{22}= 2f^2d \quad \hbox{if} \quad d\equiv 2,3 ~mod~4, \nonumber
\end{eqnarray}
and 
\begin{eqnarray}
q(x,y) &=& 2x^2 +2f xy +{1\over 2}f^2(d+1)y^2\quad \hbox{if} \quad d\equiv 1 ~mod~4\nonumber\\
q(x,y) &=&  2x^2+2f^2dy^2\quad  \hbox{if} \quad d\equiv 2,3 ~mod~4.\nonumber
\end{eqnarray}
Therefore
\displaymath
\Delta({\goth m}_f)=\cases{f^2d  & if $d\equiv 1 ~mod~4$,\cr
               4f^2d & if $d\equiv 2,3 ~mod~4$}
\enddisplaymath
and
\displaymath
 \Sigma({\goth m}_f)=+2.
\enddisplaymath
\begin{exm}
\textnormal{
Consider the Anosov maps $\phi_A,\phi_B: T^2\to T^2$ given  by matrices
\displaymath
A=\left(\matrix{5 & 2\cr 2 & 1}\right)\qquad  \hbox{and}
\qquad B=\left(\matrix{5 & 1\cr 4 & 1}\right), 
\enddisplaymath
respectively.  The reader can verify that the Alexander polynomials of 
$\phi_A$ and $\phi_B$ are identical  and  equal to $\Delta_A(t)=\Delta_B(t)= t^2-6t+1$;
yet $\phi_A$ and $\phi_B$ are {\it not} conjugate. 
Indeed, the Perron-Frobenius eigenvector of matrix $A$ is  $v_A=(1, \sqrt{2}-1)$
while of the matrix $B$ is $v_B=(1, 2\sqrt{2}-2)$. The  bilinear forms for the modules 
${\goth m}_A={\Bbb Z}+(\sqrt{2}-1){\Bbb Z}$ and 
${\goth m}_B={\Bbb Z}+(2\sqrt{2}-2){\Bbb Z}$ can be written as
\displaymath
q_A(x,y)= 2x^2-4xy+6y^2,\qquad q_B(x,y)=2x^2-8xy+24y^2,
\enddisplaymath
respectively.   The  modules ${\goth m}_A, {\goth m}_B$ are not similar in the number field 
$K={\Bbb Q}(\sqrt{2})$,  since   their  determinants $\Delta({\goth m}_A)=8$ and 
$\Delta({\goth m}_B)=32$ are not equal.  Therefore,    matrices $A$ and $B$ are not similar
in the group  $SL(2,{\Bbb Z})$.   Note that the class number $h_K=1$ for the field $K$.     
}
\end{exm}
\begin{rmk}
{\bf (Gauss method)}
\textnormal{
The reader can verify that $A$ and $B$ are non-similar by  using the {\it method of periods}, 
which dates back  to C. -F.~Gauss. According to the algorithm,   we have to find the fixed points  
$Ax=x$ and $Bx=x$,   which gives us $x_A=1+\sqrt{2}$ and
$x_B={1+\sqrt{2}\over 2}$, respectively. Then one unfolds the fixed points into a periodic continued fraction,
which gives us $x_A=[2,2,2,\dots]$ and $x_B=[1,4,1,4,\dots]$. Since the period $(\overline{2})$ of $x_A$
differs from the period $(\overline{1,4})$ of $B$,  the matrices $A$ and $B$ belong to different
similarity classes in $SL(2, {\Bbb Z})$.
}
\end{rmk}
\begin{exm}
\textnormal{
Consider the Anosov maps $\phi_A, \phi_B:  T^2\to T^2$ given by matrices 
\displaymath
A=\left(\matrix{4 & 3\cr 5 & 4}\right)\qquad  \hbox{and}
\qquad B=\left(\matrix{4 & 15\cr 1 & 4}\right), 
\enddisplaymath
respectively.  The Alexander polynomials of $\phi_A$ and $\phi_B$ are identical   
$\Delta_A(t)=\Delta_B(t)= t^2-8t+1$;    yet  the automorphisms $\phi_A$ and $\phi_B$ are 
not conjugate.   Indeed, the Perron-Frobenius eigenvector of matrix $A$ is  $v_A=(1,  {1\over 3}\sqrt{15})$
while of the matrix $B$ is $v_B=(1, {1\over 15}\sqrt{15})$.  
The corresponding  modules are  ${\goth m}_A={\Bbb Z}+( {1\over 3}\sqrt{15}){\Bbb Z}$ and 
${\goth m}_B={\Bbb Z}+( {1\over 15}\sqrt{15}){\Bbb Z}$;  therefore 
\displaymath
q_A(x,y)= 2x^2+18y^2,\qquad q_B(x,y)=2x^2 +450y^2,
\enddisplaymath
respectively.  The  modules ${\goth m}_A, {\goth m}_B$ are not similar in the number field 
$K={\Bbb Q}(\sqrt{15})$,  since  the module  determinants $\Delta({\goth m}_A)=36$ and 
$\Delta({\goth m}_B)=900$ are not equal. Therefore,  matrices $A$ and $B$ are not similar
in the group  $SL(2,{\Bbb Z})$.     
}
\end{exm}
 \index{Alexander polynomial}
\begin{exm}
{\bf ([Handelman 2009]  \cite{Han3},  p.12})
\textnormal{
Let $a, b$ be a pair of positive integers satisfying the Pell equation
$a^2-8b^2=1$;  the latter has infinitely many  solutions, e.g. $a=3, b=1$,
{\it etc.}   Denote by   $\phi_A, \phi_B: T^2\to T^2$ the Anosov maps given by 
matrices
\displaymath
A=\left(\matrix{a & 4b\cr 2b & a}\right)\qquad  \hbox{and}
\qquad B=\left(\matrix{a & 8b\cr b & a}\right), 
\enddisplaymath
respectively. 
The Alexander polynomials of $\phi_A$ and $\phi_B$ are identical   
$\Delta_A(t)=\Delta_B(t)= t^2-2a t+1$;   yet  maps  $\phi_A$ and $\phi_B$ are 
{\it not}  conjugate.  Indeed, the Perron-Frobenius eigenvector of matrix $A$ is  $v_A=(1,  {1\over 4b}\sqrt{a^2-1})$
while of the matrix $B$ is $v_B=(1, {1\over 8b}\sqrt{a^2-1})$.  
The corresponding  modules are  ${\goth m}_A={\Bbb Z}+( {1\over 4b}\sqrt{a^2-1}){\Bbb Z}$ and 
${\goth m}_B={\Bbb Z}+( {1\over 8b}\sqrt{a^2-1}){\Bbb Z}$.  It is easy to see  that  the discriminant $d=a^2-1\equiv 3~mod~4$
for all  $a\ge 2$.  Indeed, $d=(a-1)(a+1)$ and,   therefore,  integer   $a\not\equiv 1; 3 ~mod~4$;  
hence $a\equiv 2~mod~4$ so that  $a-1\equiv 1~mod ~4$
and $a+1\equiv 3 ~mod  ~4$ and, thus,  $d=a^2-1\equiv 3~mod~4$.    
 Therefore the  corresponding conductors are  $f_A=4b$ and $f_B=8b$, 
and
\displaymath
q_A(x,y)= 2x^2+32b^2(a^2-1)y^2,\quad q_B(x,y)=2x^2 +128b^2(a^2-1)y^2,
\enddisplaymath
respectively.  The  modules ${\goth m}_A, {\goth m}_B$ are not similar in the number field 
$K={\Bbb Q}(\sqrt{a^2-1})$, since   their  determinants $\Delta({\goth m}_A)=64~b^2(a^2-1)$ and 
$\Delta({\goth m}_B)=256~b^2(a^2-1)$ are not equal.  Therefore,  matrices $A$ and $B$ are not similar
in the group  $SL(2,{\Bbb Z})$. 
}
\end{exm}

\vskip1cm\noindent
{\bf Guide to the literature.}
The topology of  surface  automorphisms   is  the oldest part
of geometric topology;  it dates back to the works of  J.~Nielsen  [Nielsen 1927; 1929; 1932]  \cite{Nie1} 
and M.~Dehn  [Dehn 1938]  \cite{Deh1}.     W.~Thurston proved that that there are only
three types of such automorphisms:  they are either of finite order,  or pseudo-Anosov
or else a mixture of the two, see [Thurston 1988]
\cite{Thu1};    the topological classification of pseudo-Anosov automorphisms 
is the next biggest problem after the {\it Geometrization Conjecture}  proved by G.~Perelman,
see [Thurston 1982]   \cite{Thu2}.   An excellent introduction to the subject are the books
[Fathi, Laudenbach \& Po\'enaru 1979]   \cite{FLP}  and  [Casson \& Bleiler  1988]   \cite{CaB}.
The measured foliations on compact surfaces were introduced in 1970's by W.~Thurston  
[Thurston 1988]  \cite{Thu1}   and covered in [Hubbard \& Masur  1979]   \cite{HuMa1}.  
The Jacobi-Perron algorithm can be found in [Perron 1907]  \cite{Per1} and 
[Bernstein 1971]  \cite{BE}. 
The noncommutative invariants of pseudo-Anosov automorphisms were constructed
in \cite{Nik7}.

\section{Torsion in the torus bundles}
We assume that $M_{\alpha}$ is a {\it torus bundle},   i.e.  an $(n+1)$-dimensional
manifold fibering over the circle with monodromy $\alpha: T^n\to T^n$,  where 
$T^n$ is the  $n$-dimensional torus.  Roughly speaking,   we want to construct
a covariant functor on such bundles with the values in a category of the Cuntz-Krieger
algebras,  see Section 3.7;    such a  functor must map homeomorphic
bundles $M_{\alpha}$ to the stably isomorphic (Morita equivalent) Cuntz-Krieger
algebras.    The functor (called a {\it Cuntz-Krieger functor})  is constructed below 
and it is proved that the $K$-theory of the Cuntz-Krieger algebra is linked to the 
{\it torsion subgroup}  of the first homology group  of $M_{\alpha}$. 
The Cuntz-Krieger functor can be regarded as an ``abelianized''  version 
of functor {\bf F} $: \alpha\mapsto {\Bbb A}_{\alpha}$ constructed in 
Section 4.1,  see Remark  \ref{rmk4.2.2}.

\subsection{Cuntz-Krieger functor}
\begin{dfn}
If $T^n$ is a  torus of dimension $n\ge 1$,   then by a  torus bundle one understands   
an  $(n+1)$-dimensional manifold  
\displaymath
M_{\alpha}=\{T^n\times [0,1] ~|~ (T^n,0)=(\alpha(T^n),1)\},
\enddisplaymath
where $\alpha: T^n\to T^n$ is an automorphism of $T^n$. 
\end{dfn}
 \index{torus bundle}
\begin{rmk}
\textnormal{
The torus  bundles $M_{\alpha}$ and $M_{\alpha'}$ are homeomorphic, 
if and only if   the automorphisms $\alpha$ and $\alpha'$ are conjugate, 
i.e. $\alpha'=\beta\circ\alpha\circ\beta^{-1}$ for an automorphism 
$\beta: T^n\to T^n$. 
}
\end{rmk}
 Let $H_1(T^n;  {\Bbb Z})\cong {\Bbb Z}^n$ be the first homology of
torus;  consider the group $Aut~(T^n)$ of (homotopy classes of) automorphisms
of $T^n$.   Any $\alpha\in Aut~(T^n)$  induces a linear transformation 
of  $H_1(T^n; {\Bbb Z})$,  given by an invertible $n\times n$  matrix $A$
with the integer entries;  conversely,  each $A\in GL(n, {\Bbb Z})$  
defines an automorphism $\alpha: T^n\to T^n$.  In this  matrix 
representation, the conjugate  automorphisms $\alpha$ and  $\alpha'$  define  similar 
matrices  $A,A'\in GL(n, {\Bbb Z})$, i.e. such that $A'=BAB^{-1}$ for a matrix 
$B\in GL(n, {\Bbb Z})$.  Each class of matrices, similar to  
 a matrix $A\in GL(n, {\Bbb Z})$ and such that $tr~(A)\ge 0$ ($tr~(A)\le 0$), 
contains a matrix with only  the non-negative (non-positive) entries.
We always assume, that our bundle $M_{\alpha}$ is given by a non-negative
matrix $A$;  the matrices with  $tr~(A)\le 0$ can be  reduced to this case 
by switching the sign (from negative to positive) in the respective non-positive
representative. 
\begin{dfn}
Denote by ${\cal M}$ a category of torus bundles (of fixed dimension)  endowed with
homeomorphisms between the bundles;  denote by ${\cal A}$ 
a category of the Cuntz-Krieger algebras  ${\cal O}_A$ with  $det~(A)=\pm 1$,
endowed with stable isomorphisms between the algebras.  
By a Cuntz-Krieger map  $F:  {\cal M}\to {\cal A}$ one understands the map given 
by the formula
\displaymath
M_{\alpha}\mapsto {\cal O}_A.
\enddisplaymath
\end{dfn}
\begin{thm}\label{thm4.2.1}
{\bf (Functor on torus bundles)}
The map $F$ is a covariant functor,  which induces an isomorphism between  the 
abelian groups
\displaymath
H_1(M_{\alpha}; {\Bbb Z})\cong {\Bbb Z}\oplus K_0(F(M_{\alpha})).
\enddisplaymath
\end{thm}
\begin{rmk}\label{rmk4.2.2}
\textnormal{
The functor $F:  M_{\alpha}\mapsto {\cal O}_A$  can be obtained 
from  functor  {\bf F} $:\alpha\mapsto {\Bbb A}_{\alpha}$  on automorphisms
$\alpha: T^n\to T^n$ with values in the stationary AF-algebras ${\Bbb A}_{\alpha}$,
see Section 4.1;   the correspondence between $F$ and {\bf F} comes from the
canonical isomorphism   
\displaymath
{\cal O}_A\otimes {\cal K}\cong {\Bbb A}_{\alpha}\rtimes_{\sigma} {\Bbb Z},
\enddisplaymath
where $\sigma$ is the shift automorphism of ${\Bbb A}_{\alpha}$,
see Section 3.7.   Thus, one can interpret the invariant  ${\Bbb Z}\oplus K_0({\cal O}_A)$
as ``abelianized'' Handelman's  invariant $(\Lambda, [I], K)$ of algebra ${\Bbb A}_{\alpha}$;
here we assume that  $(\Lambda, [I], K)$ is an analog of the fundamental group $\pi_1(M_{\alpha})$. 
}
\end{rmk}

 \index{abelianized Handelman invariant}

\subsection{Proof of theorem \ref{thm4.2.1}}
The idea of proof consists in a reduction of the conjugacy problem for the
automorphisms of $T^n$ to the Cuntz-Krieger theorem on the flow equivalence 
of the subshifts of finite type, see e.g. [Lind \& Marcus 1995]  \cite{LM}.   
For  a different proof of  Theorem \ref{thm4.2.1},  see  
[Rodrigues \& Ramos  2005]   \cite{RoRa1}.

 \index{subshift of finite type}

\medskip
(i) The main reference to the subshifts of finite type (SFT) is  [Lind \& Marcus 1995]  \cite{LM}.
Recall that 
a  {\it full} Bernoulli $n$-shift is the set $X_n$ of bi-infinite sequences
$x=\{x_k\}$, where $x_k$ is a symbol taken from a set $S$ of cardinality $n$.  
The set $X_n$ is endowed with the product topology, making $X_n$ a Cantor set.
The shift homeomorphism $\sigma_n:X_n\to X_n$ is given by the formula
$\sigma_n(\dots x_{k-1}x_k x_{k+1}\dots)=(\dots x_{k}x_{k+1}x_{k+2}\dots)$
The homeomorphism defines a (discrete) dynamical system $\{X_n,\sigma_n\}$
given by the iterations of $\sigma_n$. 
Let $A$ be an $n\times n$ matrix, whose entries $a_{ij}:= a(i,j)$ are $0$ or $1$.
Consider a subset $X_A$ of $X_n$ consisting of the bi-infinite sequences,
which satisfy the restriction $a(x_k, x_{k+1})=1$ for all $-\infty<k<\infty$.
(It takes a moment to verify that $X_A$ is indeed a subset of $X_n$ and 
$X_A=X_n$, if and only if,  all the entries of $A$ are $1$'s.) By definition,
$\sigma_A=\sigma_n~|~X_A$ and the pair $\{X_A,\sigma_A\}$ is called a {\it SFT}.
A standard edge shift construction described in [Lind \& Marcus 1995]  \cite{LM}
allows to extend the notion of SFT to any matrix $A$ with the non-negative entries. 
It is well known that the SFT's $\{X_A,\sigma_A\}$ and $\{X_B,\sigma_B\}$
are topologically conjugate (as the dynamical systems), if and only if, 
the matrices $A$ and $B$ are {\it strong shift equivalent} (SSE), see 
[Lind \& Marcus 1995]  \cite{LM} for 
the corresponding definition. The SSE of two matrices is a difficult
algorithmic problem, which motivates the consideration of a weaker 
equivalence between the matrices called a {\it shift equivalence} (SE). 
Recall,  that the matrices $A$ and $B$ are said to be shift equivalent
(over ${\Bbb Z}^+$),  when there exist non-negative matrices $R$ and $S$
and a positive integer $k$ (a lag), satisfying the equations 
$AR=RB, BS=SA, A^k=RS$ and $SR=B^k$. Finally, the SFT's    
 $\{X_A,\sigma_A\}$ and $\{X_B,\sigma_B\}$ (and the matrices $A$ and $B$)
are said to be {\it flow equivalent} (FE),  if the suspension flows of the SFT's
act on the topological spaces, which are homeomorphic under a homeomorphism
that respects the orientation of the orbits of the suspension flow. 
We shall use the following implications
\displaymath
SSE\Rightarrow  SE\Rightarrow  FE.
\enddisplaymath
\begin{rmk}
\textnormal{
The first implication is rather classical, while for the second
we refer the reader to  [Lind \& Marcus 1995]  \cite{LM}, p. 456.
}
\end{rmk}
We further restrict to the SFT's given by the matrices with determinant
$\pm 1$. In view of  [Wagoner 1999]  \cite{Wag1},  Corollary 2.13,  the matrices $A$ and
$B$ with $det~(A)=\pm 1$ and $det~(B)=\pm 1$ are SE (over ${\Bbb Z}^+$),
if and only if,  matrices $A$ and $B$ are similar  in $GL_n({\Bbb Z})$.
Let now $\alpha$ and $\alpha'$ be a pair of conjugate automorphisms of $T^n$.
Since the corresponding matrices $A$ and $A'$ are similar in $GL_n({\Bbb Z})$,
one concludes that the SFT's $\{X_{A},\sigma_{A}\}$ and 
$\{X_{A'},\sigma_{A'}\}$ are SE. In particular,  the SFT's $\{X_{A},\sigma_{A}\}$ and 
$\{X_{A'},\sigma_{A'}\}$ are FE.
One can now apply the known result due to Cuntz and Krieger; it says, 
that the $C^*$-algebra ${\cal O}_A\otimes {\cal K}$ is an invariant
of the flow equivalence of the irreducible SFT's,   see   
[Cuntz \& Krieger 1980]  \cite{CuKr1},   p. 252
and its proof in Section  4 of the same work.  Thus, the map  $F$ sends the 
conjugate automorphisms of $T^n$ into the stably
isomorphic Cuntz-Krieger algebras,   i.e. $F$ is a functor.

\begin{figure}[here]
\begin{picture}(300,110)(-120,-5)
\put(23,70){\vector(0,-1){35}}
\put(122,70){\vector(0,-1){35}}
\put(45,23){\vector(1,0){53}}
\put(45,83){\vector(1,0){53}}
\put(17,20){${\cal O}_{A}$}
\put(0,55){$F$}
\put(137,55){$F$}
\put(118,20){${\cal O}_{BAB^{-1}},$}
\put(17,80){$A$}
\put(117,80){$A'=BAB^{-1}$}
\put(57,30){\sf stable}
\put(50,12){\sf isomorphism}
\put(54,90){\sf similarity}
\end{picture}
\caption{Cuntz-Krieger functor.}
\end{figure}

\medskip\noindent
Let us show that $F$ is a covariant functor.   Consider the  commutative
diagram in  Fig. 4.4,    where $A,B\in GL_n({\Bbb Z})$ and ${\cal O}_A, {\cal O}_{BAB^{-1}}\in {\cal A}$.  
Let $g_1,g_2$ be the arrows (similarity
of matrices) in the upper category and $F(g_1), F(g_2)$ the corresponding
arrows (stable  isomorphisms) in the lower category. In view of the
diagram, we have the following identities:
\begin{eqnarray}
F(g_1g_2) &=& {\cal O}_{B_2B_1AB_1^{-1}B_2^{-1}}= {\cal O}_{B_2(B_1AB_1^{-1})B_2^{-1}}\nonumber\\
          &=& {\cal O}_{B_2A'B_2^{-1}}=F(g_1)F(g_2),\nonumber
\end{eqnarray}
where $F(g_1)({\cal O}_A)={\cal O}_{A'}$ and $F(g_2)({\cal O}_{A'})={\cal O}_{A''}$.
Thus, $F$ does not reverse the arrows and is,  therefore,  a covariant functor.  
The first statement of Theorem \ref{thm4.2.1} is proved.

\bigskip
(ii) Let $M_{\alpha}$ be a torus bundle with a monodromy,  given by the matrix $A\in GL(n, {\Bbb Z})$.
It can be calculated, e.g. using the Leray spectral sequence for the fiber bundles, 
that $H_1(M_{\alpha}; {\Bbb Z})\cong {\Bbb Z}\oplus  {\Bbb Z}^n / (A-I){\Bbb Z}^n$. 
Comparing this  calculation with the $K$-theory of the Cuntz-Krieger algebra, one concludes
that $H_1(M_{\alpha}; {\Bbb Z})\cong {\Bbb Z}\oplus K_0({\cal O}_{A})$, where ${\cal O}_{A}=F(M_{\alpha})$. 
The second statement   of  Theorem \ref{thm4.2.1} follows.
$\square$   

 \index{Cuntz-Krieger algebra}

\subsection{Noncommutative invariants of torus bundles}
To illustrate Theorem  \ref{thm4.2.1},  we shall consider concrete examples of 
the torus bundles and calculate the noncommutative invariant  $K_0({\cal O}_{A})$
for them.   The reader can see,  that in some cases $K_0({\cal O}_{A})$ is  complete
invariant of a family of the torus bundles.  We compare $K_0({\cal O}_{A})$ with the
corresponding Alexander polynomial $\Delta(t)$ of the bundle and prove that  
$K_0({\cal O}_{A})$ is {\it finer} than $\Delta(t)$ (in some cases). 
\begin{exm}
\textnormal{
Consider a three-dimensional torus bundle 
\displaymath
A_1^n = \left(\matrix{1 & n\cr 0 & 1}\right), \quad n\in {\Bbb Z}.
\enddisplaymath
Using the reduction of  matrix  to its Smith normal form (see e.g. [Lind \& Marcus 1995]  \cite{LM}),
one can easily calculate
\displaymath
 K_0({\cal O}_{A_1^n})\cong {\Bbb Z}\oplus {\Bbb Z}_n.
\enddisplaymath
}
\end{exm}
\begin{rmk}
\textnormal{
The Cuntz-Krieger invariant  $K_0({\cal O}_{A_1^n})\cong {\Bbb Z}\oplus {\Bbb Z}_n$
is a complete topological invariant  of the family of bundles  $M_{\alpha_1^n}$; 
thus,  such an invariant solves the classification problem for such bundles.
}
\end{rmk}
 \index{Cuntz-Krieger invariant}
\begin{exm}
\textnormal{
Consider a three-dimensional torus bundle 
\displaymath
A_2 = \left(\matrix{5 & 2\cr 2 & 1}\right).
\enddisplaymath
Using the reduction of  matrix  to its Smith normal form,   one gets
\displaymath
 K_0({\cal O}_{A_2})\cong {\Bbb Z}_2\oplus {\Bbb Z}_2.
\enddisplaymath
}
\end{exm}
\begin{exm}
\textnormal{
Consider a three-dimensional torus bundle 
\displaymath
A_3 = \left(\matrix{5 & 1\cr 4 & 1}\right).
\enddisplaymath
Using the reduction of  matrix  to its Smith normal form,
one obtains
\displaymath
K_0({\cal O}_{A_3})\cong {\Bbb Z}_4.
\enddisplaymath
}
\end{exm}
\begin{rmk}
{\bf ($K_0({\cal O}_A)$ versus the Alexander polynomial)}
\textnormal{
Note that for the bundles $M_{\alpha_2}$ and $M_{\alpha_3}$ the Alexander polynomial:
\begin{equation}
\Delta_{A_2}(t)=\Delta_{A_3}(t)=t^2-6t+1.
\end{equation}
Therefore,  the Alexander polynomial   cannot  distinguish between the bundles
$M_{\alpha_2}$ and $M_{\alpha_3}$;  however, since   
$K_0({\cal O}_{A_2})\not\cong  K_0({\cal O}_{A_3})$,  Theorem \ref{thm4.2.1} says 
that the torus bundles $M_{\alpha_2}$ and $M_{\alpha_3}$ are topologically distinct.
Thus the noncommutative invariant $K_0({\cal O}_A)$ is {\it finer} than the Alexander
polynomial.  
}
\end{rmk}
\begin{rmk}
\textnormal{
According to  the  Thurston  Geometrization Theorem,  the torus bundle $M_{\alpha_1^n}$ is  a 
{\it  nilmanifold}  for any  $n$,  while torus bundles $M_{\alpha_2}$ and  $M_{\alpha_3}$ are
{\it solvmanifolds},   see [Thurston 1982]   \cite{Thu2}.
}
\end{rmk}

\vskip1cm\noindent
{\bf Guide to the literature.}
For an excellent introduction to the subshifts of finite type we refer the reader to 
the book by  [Lind \& Marcus 1995]  \cite{LM} and survey  by [Wagoner 1999]  \cite{Wag1}.
The Cuntz-Krieger algebras ${\cal O}_A$, the abelian group $K_0({\cal O}_A)$ and their
connection to the subshifts of finite type were introduced in  [Cuntz \& Krieger 1980]  \cite{CuKr1}.
Note that Theorem \ref{thm4.2.1} follows from the results by [Rodrigues \& Ramos  2005]   \cite{RoRa1};
however, our argument  is  different and the proof is more direct and shorter  than in the above cited work.   
The Cuntz-Krieger functor was constructed in \cite{Nik8}.

 \index{Anosov bundle}

\section{Obstruction theory for Anosov's  bundles}
We shall  use functors ranging in  the category of AF-algebras  to  study    the   {\it Anosov  bundles}
$M_{\varphi}$,  i.e. mapping tori of the Anosov diffeomorphisms
\displaymath
\varphi: M\to M
\enddisplaymath
of a smooth manifold $M$,   see e.g. [Smale 1967]  \cite{Sma1}, p. 757.   
Namely,  we construct a covariant functor $F$ form the category of Anosov's bundles to a category
of stationary AF-algebras;   the functor sends each continuous map between
the bundles to  a stable homomorphism between  the corresponding AF-algebras. 
We develop an {\it obstruction theory}  for  continuous maps between
Anosov's bundles;  such a theory  exploits  noncommutative invariants   derived from 
the triple $(\Lambda, [I], K)$ attached to  stationary AF-algebras. 
We illustrate the obstruction theory by concrete examples of dimension 2, 3 and 4.

 \index{fundamental AF-algebra}

\subsection{Fundamental AF-algebra}
By a $q$-dimensional,  class  $C^r$ foliation of an $m$-dimensional manifold $M$ 
one understands a decomposition of $M$ into a union of
disjoint connected subsets $\{ {\cal L}_{\alpha}\}_{\alpha\in A}$, called
the {\it leaves}, see e.g.  [Lawson 1974]   \cite{Law1}.  They must satisfy the
following property: each point in $M$ has a neighborhood $U$
and a system of local class $C^r$ coordinates 
$x=(x^1,\dots, x^m): U\to {\Bbb R}^m$ such that for each leaf 
${\cal L}_{\alpha}$, the components of $U\cap {\cal L}_{\alpha}$
are described by the equations $x^{q+1}=Const, \dots, x^m=Const$.
Such a foliation is denoted by ${\cal F}=\{ {\cal L_{\alpha}}\}_{\alpha\in A}$.
The number $k=m-q$ is called a  codimension of the foliation.
An  example of a codimension $k$ foliation ${\cal F}$ 
is given by a closed $k$-form $\omega$ on $M$:  the leaves of ${\cal F}$
are tangent to a plane defined by the normal vector  $\omega(p)=0$ at each point $p$ of $M$. 
The $C^r$-foliations ${\cal F}_0$ and ${\cal F}_1$ 
of codimension $k$
are said to be $C^s$-conjugate ($0\le s\le r$), if there exists an (orientation-preserving) 
diffeomorphism of $M$, of class $C^s$, which maps the leaves of ${\cal F}_0$
onto the leaves of ${\cal F}_1$;  when $s=0$, ${\cal F}_0$ and ${\cal F}_1$
are {\it topologically conjugate}. 
Denote by  $f: N\to M$  a map of class $C^s$ ($1\le s\le r$) of a manifold
$N$ into $M$;  the map $f$ is said to be {\it transverse} to ${\cal F}$,
if for all $x\in N$ it holds $T_y(M)=\tau_y({\cal F})+f_* T_x(N)$,
where $\tau_y({\cal F})$ are the vectors of $T_y(M)$ tangent
to ${\cal F}$ and  $f_*: T_x(N)\to T_y(M)$ is  the linear map
on tangent vectors induced by $f$, where $y=f(x)$. 
If map  $f: N\to M$ is transverse to a foliation ${\cal F}'=\{{\cal L}\}_{\alpha\in A}$
on $M$, then $f$  induces a class $C^s$ foliation ${\cal F}$ on $N$,
where the leaves are defined as $f^{-1}({\cal L}_{\alpha})$ for all 
$\alpha\in A$; it is immediate, that $codim~({\cal F})=codim~({\cal F}')$. 
We shall call ${\cal F}$ an {\it induced foliation}.  
When $f$ is a submersion, it is transverse  to any foliation of $M$;
in this case, the induced foliation ${\cal F}$ is correctly defined for all ${\cal F}'$
on $M$,   see   [Lawson 1974]   \cite{Law1}, p.373. 
Notice,  that for $M=N$  the above definition corresponds to  
topologically conjugate foliations ${\cal F}$ and ${\cal F}'$.  
To introduce measured foliations, denote by $P$ and $Q$  two $k$-dimensional
submanifolds of $M$,  which are everywhere transverse to a foliation ${\cal F}$ of codimension $k$. 
Consider a collection of $C^r$ homeomorphisms between subsets of $P$ and $Q$
induced by  a return map along the leaves of ${\cal F}$.
The collection of all such homeomorphisms between subsets of all possible pairs of 
transverse manifolds generates a {\it holonomy pseudogroup} of ${\cal F}$
under composition of the homeomorphisms, see  [Plante 1975]  \cite{Pla1}, p.329.  
A foliation ${\cal F}$  is said to have measure preserving holonomy,
if its holonomy pseudogroup has a non-trivial invariant measure,  which is
finite on compact sets;  for brevity, we call ${\cal F}$ a  {\it measured foliation}.
An example of measured foliation is a foliation, determined by 
closed $k$-form $\omega$; the restriction of $\omega$ to a transverse $k$-dimensional
manifold determines a volume element, which gives a positive invariant measure on open sets.     
Each measured foliation ${\cal F}$ defines an element of the cohomology
group $H^k(M; {\Bbb R})$,  see  [Plante 1975]  \cite{Pla1};  in the case of ${\cal F}$ given by a closed
$k$-form $\omega$, such an element coincides with the de Rham cohomology class of $\omega$,  
{\it ibid}. In view of the isomorphism $H^k(M; {\Bbb R})\cong Hom~(H_k(M), {\Bbb R})$,
foliation ${\cal F}$ defines a linear map $h$ from the $k$-th homology group $H_k(M)$ 
to ${\Bbb R}$. 
\begin{dfn}
By a Plante group $P({\cal F})$ of measured foliation ${\cal F}$ 
one  understand the finitely generated abelian subgroup $h(H_k(M)/Tors)\subset  {\Bbb R}$. 
\end{dfn}
\begin{rmk}
\textnormal{
 If $\{\gamma_i\}$ is a basis of the homology group $H_k(M)$,
then the periods $\lambda_i=\int_{\gamma_i}\omega$ are  generators of the group  
$P({\cal F})$,   see  [Plante 1975]  \cite{Pla1}.
}
\end{rmk}

\smallskip\noindent
Let  $\lambda=(\lambda_1,\dots,\lambda_n)$ be a basis of the Plante group $P({\cal F})$
of  a  measured foliation ${\cal F}$, such that $\lambda_i>0$.  
Take a vector $\theta=(\theta_1,\dots,\theta_{n-1})$
with $\theta_i=\lambda_{i+1} / \lambda_1$;   the Jacobi-Perron continued fraction of
vector $(1, \theta)$ (or, projective class of vector $\lambda$)  is given by the formula
\displaymath
\left(\matrix{1\cr \theta}\right)=
\lim_{i\to\infty} \left(\matrix{0 & 1\cr I & b_1}\right)\dots
\left(\matrix{0 & 1\cr I & b_i}\right)
\left(\matrix{0\cr {\Bbb I}}\right)=
\lim_{i\to\infty} B_i\left(\matrix{0\cr {\Bbb I}}\right),
\enddisplaymath
where $b_i=(b^{(i)}_1,\dots, b^{(i)}_{n-1})^T$ is a vector of the non-negative integers,  
$I$ the unit matrix and ${\Bbb I}=(0,\dots, 0, 1)^T$,  see [Bernstein 1971]  \cite{BE}, p.13;  
the $b_i$ are obtained from $\theta$  by the Euclidean algorithm, {\it ibid.},   pp.2-3.
\begin{dfn}
An AF-algebra given by the Bratteli diagram with the incidence matrices $B_i$
will be called associated to the measured  foliation ${\cal F}$;  we shall denote such 
an algebra by ${\Bbb A}_{\cal F}$.  
 \end{dfn}
\begin{rmk}
\textnormal{
Taking another basis of the Plante group $P({\cal F})$ gives 
an AF-algebra  which is stably isomorphic (Morita equivalent)  to  ${\Bbb A}_{\cal F}$;  
this is an algebraic recast of the main  property  of the Jacobi-Perron fractions. 
}
\end{rmk}
If ${\cal F}'$ is a measured foliation on a manifold $M$ and 
$f:N\to M$ is a submersion,  then  induced foliation 
${\cal F}$ on $N$ is a measured foliation.
We shall denote by {\bf MFol}  the category of all manifolds
with measured foliations (of fixed codimension),  whose arrows are submersions of
the manifolds;  by {\bf MFol}$_0$  we understand  a subcategory of {\bf MFol}, 
consisting  of manifolds,   whose foliations have a  unique transverse measure. 
Let   {\bf AF-alg}  be a category 
of the (isomorphism classes of)  AF-algebras given by {\it convergent} Jacobi-Perron fractions,
so that the arrows of {\bf AF-alg}  are  stable homomorphisms of the AF-algebras.
By $F:$  {\bf MFol}$_0\to$ {\bf AF-alg}  we denote a map  given by the
formula ${\cal F}\mapsto {\Bbb A}_{\cal F}$.  Notice, that $F$ is correctly defined,
since foliations with the unique measure have the convergent Jacobi-Perron
fractions; this assertion follows from [Bauer  1996]  \cite{Bau1}.    The following 
result will be proved in Section 4.3.2. 
\begin{thm}\label{thm4.3.1}
The map   $F:$  {\bf MFol}$_0\to$ {\bf AF-alg}   is a functor,
which sends any pair of induced foliations to a pair of stably homomorphic
AF-algebras.
\end{thm}
Let $M$ be an $m$-dimensional manifold and $\varphi: M\to M$ a diffeomorphisms  of $M$;
recall,  that  an orbit of point $x\in M$ is the subset $\{\varphi^n(x) ~|~n\in {\Bbb Z}\}$
of $M$. The finite orbits $\varphi^m(x)=x$ are called periodic; 
when $m=1$, $x$ is a {\it fixed point} of diffeomorphism $\varphi$.
The fixed point $p$ is {\it hyperbolic} if the eigenvalues $\lambda_i$ of 
the linear map $D\varphi(p): T_p(M)\to T_p(M)$ do not lie at the unit circle. 
If $p\in M$ is a hyperbolic fixed point of a diffeomorphism $\varphi: M\to M$,
denote by  $T_p(M)=V^s+V^u$ the corresponding decomposition of the tangent space 
under the linear map $D\varphi(p)$, where $V^s$ ($V^u$) is the eigenspace of $D\varphi(p)$
corresponding to $|\lambda_i|>1$ ($|\lambda_i|<1$). 
For a sub-manifold $W^s(p)$
there exists a contraction $g: W^s(p)\to W^s(p)$ with fixed point $p_0$
and an injective equivariant immersion $J: W^s(p)\to M$, 
such that $J(p_0)=p$ and $DJ(p_0): T_{p_0}(W^s(p))\to T_p(M)$
is an isomorphism; the image of $J$ defines an immersed submanifold 
$W^s(p)\subset M$ called a {\it stable manifold} of $\varphi$ at $p$. 
Clearly, $dim~(W^s(p))=dim~(V^s)$. 
\begin{dfn}
{\bf ([Anosov 1967]  \cite{Ano1})}
A diffeomorphism  $\varphi:M\to M$ is called  Anosov  if  there exists a splitting of the tangent 
bundle $T(M)$ into a continuous Whitney sum $T(M)=E^s+E^u$
invariant under $D\varphi: T(M)\to T(M)$,  so that $D\varphi: E^s\to E^s$
is contracting and $D\varphi: E^u\to E^u$ is expanding map.
\end{dfn}
 \index{Anosov diffeomorphism}
\begin{rmk}
\textnormal{
The Anosov diffeomorphism
imposes a restriction on topology of manifold $M$, in the sense
that not each manifold can support such a diffeomorphism; 
however, if one Anosov diffeomorphism exists on $M$, there 
are infinitely many (conjugacy classes of) such diffeomorphisms
on $M$. It is an open problem of S.~Smale, which $M$ can carry
an Anosov diffeomorphism;  so far,  it is proved that    
 the hyperbolic  diffeomorphisms of $m$-dimensional tori  and 
certain automorphisms of the nilmanifolds are Anosov's, see e.g.  
[Smale 1967]  \cite{Sma1}. 
}
\end{rmk}
Let  $p$ be a fixed point
of the Anosov diffeomorphism $\varphi:M\to M$ and  $W^s(p)$
its stable manifold. Since $W^s(p)$ cannot have self-intersections 
 or limit compacta,   $W^s(p)\to M$ is a dense immersion, i.e. 
 the closure of $W^s(p)$ is the entire $M$.  Moreover, if $q$
is a periodic point of $\varphi$ of period $n$, then  $W^s(q)$ is
a translate of $W^s(p)$, i.e. locally they look like two parallel lines. 
Consider a foliation ${\cal F}$ of $M$, whose leaves are the
translates of $W^s(p)$; the ${\cal F}$ is a continuous foliation, 
which is invariant under the action of diffeomorphism $\varphi$ on its
leaves, i.e. $\varphi$ moves leaves of ${\cal F}$ to the leaves of ${\cal F}$,  see 
[Smale 1967]  \cite{Sma1}, p.760.  
  The holonomy of ${\cal F}$ preserves the Lebesgue measure and, therefore,
${\cal F}$ is a measured foliation; we shall call it an {\it invariant measured
foliation} and  denote by ${\cal F}_{\varphi}$.
\begin{dfn}
By a fundamental AF-algebra we shall understand  the AF-algebra of foliation 
${\cal F}_{\varphi}$,  where  $\varphi: M\to M$  is  an Anosov diffeomorphism of a 
manifold $M$;  the fundamental AF-algebra will be denoted by ${\Bbb A}_{\varphi}$. 
\end{dfn}
\begin{thm}\label{thm4.3.2}
The ${\Bbb A}_{\varphi}$ is a stationary AF-algebra. 
\end{thm}
Consider the mapping torus of the Anosov  diffeomorphism $\varphi$,  i.e. a manifold  
\displaymath
M_{\varphi}:= M\times [0,1]~/\sim,  ~\hbox{where}  ~(x,0)\sim (\varphi(x), 1), ~\forall x\in M.  
\enddisplaymath
Let {\bf AnoBnd} be a category of the mapping tori
of all Anosov's  diffeomorphisms; 
the arrows  of {\bf AnoBnd}  are continuous maps between the mapping tori. Likewise,  
let {\bf Fund-AF} 
be a category of all fundamental AF-algebras; the arrows  of {\bf Fund-AF} 
are stable homomorphisms between the fundamental AF-algebras. By $F:$ {\bf AnoBnd} $\to$
{\bf Fund-AF}  we understand a map given by the formula $M_{\varphi}\mapsto {\Bbb A}_{\varphi}$,  where $M_{\varphi}\in$
{\bf AnoBnd}  and ${\Bbb A}_{\varphi}\in$  {\bf Fund-AF}.   The following theorem says that $F$
is a functor. 
\begin{thm}\label{thm4.3.3}
{\bf (Functor on Anosov's bundles)}
The map $F$ is a covariant functor, which sends each continuous map $N_{\psi}\to M_{\varphi}$ to a stable
homomorphism ${\Bbb A}_{\psi}\to {\Bbb A}_{\varphi}$ of the corresponding fundamental AF-algebras.
\end{thm}
\begin{rmk}
{\bf (Obstruction theory)}
\textnormal{
Theorem \ref{thm4.3.3} can be used  e.g. in the obstruction theory,
because  stable homomorphisms of the fundamental AF-algebras are easier 
to detect,  than  continuous maps between  manifolds $N_{\psi}$ and $M_{\varphi}$;
such homomorphisms are bijective with the inclusions of certain ${\Bbb Z}$-modules
  belonging to  a  real  algebraic number field.  Often it is possible to prove, that no inclusion
is possible and, thus, draw a  topological conclusion about the  maps,  see Section  4.3.3. 
}
\end{rmk}

 \index{obstruction theory}

\subsection{Proofs}
\subsubsection{Proof of Theorem  \ref{thm4.3.1}}
Let ${\cal F}'$ be measured foliation on $M$, given by a closed form
$\omega'\in H^k(M; {\Bbb R})$; let ${\cal F}$ be measured foliation on $N$,
induced by a submersion $f: N\to M$.  Roughly speaking, we have to prove,
that diagram in Fig. 4.5  is commutative;  the proof amounts to the fact, that 
the periods of form $\omega'$ are contained among the periods of form
$\omega\in H^k(N; {\Bbb R})$ corresponding to the foliation ${\cal F}$.  
\begin{figure}[here]
\begin{picture}(300,110)(-110,-5)
\put(20,70){\vector(0,-1){35}}
\put(130,70){\vector(0,-1){35}}
\put(45,23){\vector(1,0){60}}
\put(45,83){\vector(1,0){60}}
\put(15,20){${\Bbb A}_{\cal F}$}
\put(128,20){${\Bbb A}_{\cal F'}$}
\put(17,80){${\cal F}$}
\put(125,80){${\cal F}'$}
\put(60,30){\sf stable}
\put(45,10){\sf homomorphism}
\put(54,90){\sf induction}
\end{picture}
\caption{$F:$  {\bf MFol}$_0\to$ {\bf AF-alg} .}
\end{figure}
The map $f$ defines a homomorphism $f_*: H_k(N)\to H_k(M)$ of the $k$-th
homology groups; let $\{e_i\}$ and $\{e_i'\}$ be a basis in $H_k(N)$ and
$H_k(M)$, respectively. Since $H_k(M)=H_k(N)~/~ker~(f_*)$, we shall denote 
by $[e_i]:= e_i+ker~(f_*)$ a coset representative of $e_i$;  these can be
identified with the elements $e_i\not\in ker~(f_*)$.  The integral
$\int_{e_i}\omega$ defines a scalar product $H_k(N)\times H^k(N; {\Bbb R})\to {\Bbb R}$,
so that $f_*$ is a linear self-adjoint operator; thus, we can write:
\displaymath
\lambda_i'=\int_{e_i'}\omega'=\int_{e_i'}f^*(\omega)=\int_{f_*^{-1}(e_i')}\omega=
\int_{[e_i]}\omega\in P({\cal F}),
\enddisplaymath
where $P({\cal F})$ is the Plante group (the group of periods) of foliation ${\cal F}$.
Since $\lambda_i'$ are generators of $P({\cal F}')$, we conclude that 
$P({\cal F}')\subseteq P({\cal F})$.  Note,  that $P({\cal F}')=P({\cal F})$
if and only if  $f_*$ is an isomorphism. 

 \index{Plante group}

One can apply a criterion of the stable  homomorphism of AF-algebras;  namely,
${\Bbb A}_{\cal F}$ and ${\Bbb A}_{\cal F'}$ are stably homomorphic, if and only if,
there exists a positive homomorphism $h: G\to H$ between their dimension
groups $G$ and $H$,  see  [Effros 1981]  \cite{E},  p.15.  But $G\cong P({\cal F})$ and $H\cong P({\cal F}')$,
while $h=f_*$.  Thus, ${\Bbb A}_{\cal F}$ and ${\Bbb A}_{\cal F'}$ are stably homomorphic.

The functor $F$ is compatible with the composition;  indeed,
let $f:N\to M$ and $f':L\to N$ be submersions.  If ${\cal F}$
is a measured foliation of $M$,  one gets the induced 
foliations ${\cal F}'$ and ${\cal F}''$ on $N$ and $L$, respectively;
these foliations fit the diagram
$(L, {\cal F}'')\buildrel f'\over\longrightarrow (N, {\cal F}')\buildrel f\over\longrightarrow (M, {\cal F})$
and the corresponding Plante groups are included: $P({\cal F}'')\supseteq P({\cal F}')\supseteq P({\cal F})$.
Thus, $F(f'\circ f)=F(f')\circ F(f)$, since the inclusion of the Plante groups corresponds to the
composition of homomorphisms;  Theorem \ref{thm4.3.1} is proved.
$\square$

\subsubsection{Proof of Theorem \ref{thm4.3.2}}
Let $\varphi:M\to M$ be an Anosov diffeomorphism; we proceed by showing, that
invariant foliation ${\cal F}_{\varphi}$ is given by form $\omega\in H^k(M; {\Bbb R})$,
which is an eigenvector of the linear map $[\varphi]:  H^k(M; {\Bbb R})\to  H^k(M; {\Bbb R})$
induced by $\varphi$.
Indeed, let $0<c<1$ be contracting constant of the stable sub-bundle $E^s$ of
diffeomorphism $\varphi$ and $\Omega$ the corresponding volume element; by definition,
$\varphi(\Omega)=c\Omega$. Note, that
$\Omega$ is given by restriction of form $\omega$ to a $k$-dimensional 
manifold, transverse to the leaves of ${\cal F}_{\varphi}$. The leaves of ${\cal F}_{\varphi}$
are fixed by $\varphi$ and, therefore, $\varphi(\Omega)$ is given by a multiple $c\omega$
of form $\omega$. Since $\omega\in H^k(M; {\Bbb R})$ is a vector, whose
coordinates define ${\cal F}_{\varphi}$ up to a scalar, we conclude, that $[\varphi](\omega)=c\omega$,
i.e. $\omega$ is an eigenvector of the linear map $[\varphi]$.  
Let  $(\lambda_1,\dots,\lambda_n)$ be a basis of the Plante group $P({\cal F}_{\varphi})$,
such that $\lambda_i>0$. Notice, that $\varphi$ acts on $\lambda_i$ as multiplication
by constant $c$; indeed, since $\lambda_i=\int_{\gamma_i}\omega$, we have:
\displaymath
\lambda_i'=\int_{\gamma_i}[\varphi](\omega)=\int_{\gamma_i}c\omega=c\int_{\gamma_i}\omega=c\lambda_i,
\enddisplaymath
where $\{\gamma_i\}$ is a basis in $H_k(M)$.  Since $\varphi$ preserves the leaves
of ${\cal F}_{\varphi}$, one concludes that $\lambda_i'\in P({\cal F}_{\varphi})$;
therefore, $\lambda_j'=\sum b_{ij}\lambda_i$ for a non-negative 
integer matrix $B=(b_{ij})$. According to  [Bauer 1996]  \cite{Bau1},  matrix $B$ can be 
written as a finite product:
\displaymath
B=
\left(\matrix{0 & 1\cr I & b_1}\right)\dots
\left(\matrix{0 & 1\cr I & b_p}\right):=B_1\dots B_p,
\enddisplaymath
where $b_i=(b^{(i)}_1,\dots, b^{(i)}_{n-1})^T$ is a vector of non-negative 
integers and   $I$ the unit matrix.   Let $\lambda=(\lambda_1,\dots,\lambda_n)$.  
Consider a purely periodic Jacobi-Perron continued fraction: 
\displaymath
\lim_{i\to\infty} 
\overline{B_1\dots B_p}
\left(\matrix{0\cr {\Bbb I}}\right),
\enddisplaymath
where  ${\Bbb I}=(0,\dots, 0, 1)^T$;  by a basic property of such fractions, 
it converges to an eigenvector  $\lambda'=(\lambda_1',\dots,\lambda_n')$ of matrix $B_1\dots B_p$, 
see [Bernstein 1971]  \cite{BE},  Chapter 3.   But $B_1\dots B_p=B$ and $\lambda$ is an eigenvector of matrix $B$;
therefore, vectors $\lambda$ and $\lambda'$ are collinear.  The collinear vectors  
are known to have the same continued fractions;  thus, we have     
\displaymath
\left(\matrix{1\cr \theta}\right)=
\lim_{i\to\infty} 
\overline{B_1\dots B_p}
\left(\matrix{0\cr {\Bbb I}}\right),
\enddisplaymath
where  $\theta=(\theta_1,\dots,\theta_{n-1})$ and  $\theta_i=\lambda_{i+1}/\lambda_1$. 
Since vector $(1,\theta)$ unfolds into a periodic Jacobi-Perron  continued fraction,
we conclude, that the AF-algebra ${\Bbb A}_{\varphi}$ is stationary.
Theorem  \ref{thm4.3.2} is proved.   
$\square$

\subsubsection{Proof of Theorem \ref{thm4.3.3}}
Let $\psi: N\to N$ and $\varphi: M\to M$ be a pair of Anosov diffeomorphisms;
denote by $(N, {\cal F}_{\psi})$ and $(M, {\cal F}_{\varphi})$ the corresponding
invariant  foliations of manifolds $N$ and $M$,  respectively. 
In view of Theorem \ref{thm4.3.1},  it is sufficient to prove, that the diagram in
Fig.4.6  is commutative. We shall split  the proof in a series of lemmas.
\begin{figure}[here]
\begin{picture}(300,110)(-100,-5)
\put(20,70){\vector(0,-1){35}}
\put(130,70){\vector(0,-1){35}}
\put(45,23){\vector(1,0){60}}
\put(45,83){\vector(1,0){60}}
\put(0,20){$(N, {\cal F}_{\psi})$}
\put(115,20){$(M, {\cal F}_{\varphi})$}
\put(10,80){$N_{\psi}$}
\put(120,80){$M_{\varphi}$}
\put(60,30){\sf induced}
\put(55,10){\sf foliations}
\put(54,90){\sf continuous}
\put(74,70){\sf map}
\end{picture}
\caption{Mapping tori and invariant foliations.}
\end{figure}
\begin{lem}\label{lem4.3.1}
There exists a continuous map $N_{\psi}\to M_{\varphi}$,   whenever
$f\circ\varphi=\psi\circ f$ for a submersion $f: N\to M$.  
\end{lem}
{\it Proof.} (i) Suppose, that $h: N_{\psi}\to M_{\varphi}$ is a continuous
map;  let us show, that there exists a submersion $f: N\to M$, such that
$f\circ\varphi=\psi\circ f$.  Both $N_{\psi}$ and $M_{\varphi}$ fiber over
the circle $S^1$ with the projection map  $p_{\psi}$ and $p_{\varphi}$,  respectively; 
therefore, the diagram in Fig. 4.7  is commutative.  Let $x\in S^1$; 
since $p_{\psi}^{-1}=N$ and $p_{\varphi}^{-1}=M$, the restriction of $h$
to $x$ defines a submersion $f: N\to M$, i.e. $f=h_x$. Moreover, 
since $\psi$ and $\varphi$ are monodromy maps of  the bundle, it holds:
\displaymath
\left\{
\begin{array}{cc}
p_{\psi}^{-1}(x+2\pi)  &= \psi(N),\\
p_{\varphi}^{-1}(x+2\pi)  &= \varphi(M).
\end{array}
\right.
\enddisplaymath
From the diagram in Fig. 4.7,  we get 
$\psi(N)=p_{\psi}^{-1}(x+2\pi)=f^{-1}(p_{\varphi}^{-1}(x+2\pi))=
f^{-1}(\varphi(M))=f^{-1}(\varphi(f(N)))$;  thus, $f\circ\psi=\varphi\circ f$. 
The necessary condition of Lemma  \ref{lem4.3.1}  follows.

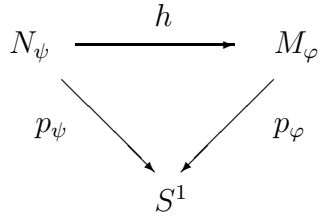
\begin{figure}[here]
\begin{picture}(300,110)(-100,-5)
\put(40,70){\vector(1,-1){35}}
\put(120,70){\vector(-1,-1){35}}
\put(45,83){\vector(1,0){60}}
\put(20,80){$N_{\psi}$}
\put(120,80){$M_{\varphi}$}
\put(75,20){$S^1$}
\put(75,90){$h$}
\put(30,50){$p_{\psi}$}
\put(120,50){$p_{\varphi}$}
\end{picture}
\caption{The fiber bundles $N_{\psi}$ and $M_{\varphi}$ over $S^1$.}
\end{figure}

\bigskip
(ii) Suppose, that $f: N\to M$ is a submersion, such that $f\circ\varphi=\psi\circ f$;
we have to construct a continuous map $h: N_{\psi}\to M_{\varphi}$. Recall,  that
\displaymath
\left\{
\begin{array}{cc}
N_{\psi}  &= \{N\times [0,1]~|~(x,0)\sim(\psi(x),1)\}, \\  
M_{\varphi} &= \{M\times [0,1]~|~(y,0)\sim(\varphi(y),1)\}.
\end{array}
\right.
\enddisplaymath
We shall identify the points of $N_{\psi}$ and $M_{\varphi}$ using
the substitution $y=f(x)$; it remains to verify, that such an identification
will satisfy the gluing condition $y\sim\varphi(y)$. In view of condition
$f\circ\varphi=\psi\circ f$, we have:
\displaymath
y=f(x)\sim f(\psi(x))=\varphi(f(x))=\varphi(y).
\enddisplaymath
Thus, $y\sim\varphi(y)$ and, therefore, the map $h: N_{\psi}\to M_{\varphi}$
is continuous. The sufficient condition of lemma \ref{lem4.3.1}  is proved.
$\square$

\begin{figure}[here]
\begin{picture}(300,110)(-100,-5)
\put(20,70){\vector(0,-1){35}}
\put(130,70){\vector(0,-1){35}}
\put(45,23){\vector(1,0){60}}
\put(45,83){\vector(1,0){60}}
\put(-10,20){$H^k(M, {\Bbb R})$}
\put(120,20){$H^k(M, {\Bbb R})$}
\put(-10,80){$H^k(N; {\Bbb R})$}
\put(120,80){$H^k(N, {\Bbb R})$}
\put(75,30){$[\varphi]$}
\put(75,90){$[\psi]$}
\put(0,50){$[f]$}
\put(138,50){$[f]$}
\end{picture}
\caption{The linear maps $[\psi], [\varphi]$ and $[f]$.}
\end{figure}
\begin{lem}\label{lem4.3.2}
If a submersion $f: N\to M$ satisfies condition $f\circ\varphi=\psi\circ f$
for the Anosov diffeomorphisms $\psi: N\to N$ and $\varphi: M\to M$,
then the invariant foliations $(N, {\cal F}_{\psi})$ and $(M, {\cal F}_{\varphi})$
are induced by $f$. 
 \end{lem}
{\it Proof.} 
The invariant foliations ${\cal F}_{\psi}$ and ${\cal F}_{\varphi}$
are measured;  we shall denote by $\omega_{\psi}\in H^k(N; {\Bbb R})$
and $\omega_{\varphi}\in H^k(M; {\Bbb R})$ the corresponding cohomology
class,  respectively.  The linear maps on $H^k(N; {\Bbb R})$ and $H^k(M; {\Bbb R})$
induced by $\psi$ and $\varphi$, we shall denote by $[\psi]$ and $[\varphi]$; 
the linear map between $H^k(N; {\Bbb R})$ and $H^k(M; {\Bbb R})$
induced by $f$, we write as $[f]$.  Notice, that $[\psi]$ and $[\varphi]$
are isomorphisms,  while $[f]$ is generally a homomorphism. It was shown
earlier, that $\omega_{\psi}$ and $\omega_{\varphi}$ are eigenvectors of
linear maps $[\psi]$ and $[\varphi]$, respectively; in other words, we have:
\displaymath
\left\{
\begin{array}{cc}
\hbox{$[\psi]$}\omega_{\psi}  &=  c_1\omega_{\psi},\\  
\hbox{$[\varphi]$}\omega_{\varphi} &= c_2\omega_{\varphi},
\end{array}
\right.
\enddisplaymath
where $0<c_1<1$ and $0<c_2<1$.  Consider a diagram in Fig. 4.8 ,  which involves
the linear maps $[\psi], [\varphi]$ and $[f]$;  the diagram is commutative,
since  condition $f\circ\varphi=\psi\circ f$ implies,  that 
$[\varphi]\circ [f]=[f]\circ [\psi]$. Take the eigenvector $\omega_{\psi}$
and consider its image under the linear map $[\varphi]\circ [f]$:
\displaymath
[\varphi]\circ [f](\omega_{\psi})=[f]\circ [\psi](\omega_{\psi})=
[f](c_1\omega_{\psi})=c_1\left([f](\omega_{\psi})\right).
\enddisplaymath
Therefore, vector $[f](\omega_{\psi})$ is an eigenvector of the
linear map $[\varphi]$; let  compare it with the eigenvector $\omega_{\varphi}$:
\displaymath
\left\{
\begin{array}{cc}
\hbox{$[\varphi]$}\left(\hbox{$[f]$}(\omega_{\psi})\right)
 &=
 c_1\left(\hbox{$[f]$}(\omega_{\psi})\right),\\ 
\hbox{$[\varphi]$}\omega_{\varphi}  &=  c_2\omega_{\varphi}.
\end{array}
\right.
\enddisplaymath
We conclude, therefore, that $\omega_{\varphi}$ and $[f](\omega_{\psi})$
are collinear vectors,  such that $c_1^m=c_2^n$ for some integers $m,n>0$;
a scaling gives us $[f](\omega_{\psi})=\omega_{\varphi}$. 
The latter is an analytic formula, which says that the submersion
$f:N\to M$ induces foliation $(N, {\cal F}_{\psi})$  from   the
foliation $(M, {\cal F}_{\varphi})$.  Lemma \ref{lem4.3.2} is proved.
$\square$

\bigskip
To finish the proof of Theorem \ref{thm4.3.3},  let $N_{\psi}\to M_{\varphi}$
be a continuous map;  by Lemma \ref{lem4.3.1}, there exists a submersion
$f: N\to M$, such that $f\circ\varphi=\psi\circ f$. Lemma \ref{lem4.3.2}
says, that in this case the invariant measured foliations $(N, {\cal F}_{\psi})$
and $(M, {\cal F}_{\varphi})$ are induced. On the other hand, from Theorem  \ref{thm4.3.2}
we know, that the Jacobi-Perron continued fraction connected to foliations ${\cal F}_{\psi}$
and ${\cal F}_{\varphi}$ are periodic and, hence,  convergent,  see e.g [Bernstein 1971]  \cite{BE};
therefore, one can apply Theorem  \ref{thm4.3.1}   which says  that the AF-algebra
${\Bbb A}_{\psi}$ is stably homomorphic to the AF-algebra ${\Bbb A}_{\varphi}$.  The latter are,
 by definition, the fundamental AF-algebras of the Anosov diffeomorphisms $\psi$ and 
$\varphi$,  respectively.   Theorem \ref{thm4.3.3} is proved.
$\square$

\subsection{Obstruction theory}
Let ${\Bbb A}_{\psi}$ be  a  fundamental AF-algebra and $B$ its primitive incidence matrix,
i.e.  $B$ is not a power of some positive integer matrix.  Suppose that the characteristic polynomial of $B$ is irreducible
and let  $K_{\psi}$ be its  splitting field; then  $K_{\psi}$ is a Galois extension of ${\Bbb Q}$. 
\begin{dfn}
We call $Gal~({\Bbb A}_{\psi}):=Gal~(K_{\psi}|{\Bbb Q})$ the Galois group of the fundamental
AF-algebra  ${\Bbb A}_{\psi}$;  such a group  is determined by the AF-algebra  ${\Bbb A}_{\psi}$.
\end{dfn}
The second algebraic field is connected to the Perron-Frobenius eigenvalue $\lambda_B$
of the matrix $B$; we shall denote this field ${\Bbb Q}(\lambda_B)$. Note, that
${\Bbb Q}(\lambda_B)\subseteq K_{\psi}$ and ${\Bbb Q}(\lambda_B)$  is not, in general, 
a Galois extension of ${\Bbb Q}$; the reason being complex roots the polynomial 
$char~(B)$ may have and if there are no such roots ${\Bbb Q}(\lambda_B)=K_{\psi}$.
There is still a group $Aut~({\Bbb Q}(\lambda_B))$ of automorphisms of ${\Bbb Q}(\lambda_B)$
fixing the field ${\Bbb Q}$ and $Aut~({\Bbb Q}(\lambda_B))\subseteq Gal~(K_{\psi})$
is a subgroup inclusion.  
\begin{lem}\label{lem4.3.3}
If $h: {\Bbb A}_{\psi}\to {\Bbb A}_{\varphi}$ is a stable homomorphism, 
then ${\Bbb Q}(\lambda_{B'})\subseteq K_{\psi}$  is a field inclusion.    
\end{lem}
{\it Proof.} Notice  that the non-negative matrix $B$ becomes strictly positive, 
when a proper power of it is taken; we always assume $B$ positive.  Let 
$\lambda=(\lambda_1,\dots,\lambda_n)$ be a basis of the Plante group $P({\cal F}_{\psi})$.
Following the proof of Theorem \ref{thm4.3.2},  one concludes that $\lambda_i\in K_{\psi}$;
indeed, $\lambda_B\in K_{\psi}$ is the Perron-Frobenius eigenvalue of $B$ , while $\lambda$   
the corresponding eigenvector. The latter can be scaled so, that $\lambda_i\in K_{\psi}$. 
Any stable homomorphism $h: {\Bbb A}_{\psi}\to {\Bbb A}_{\varphi}$ induces a positive
homomorphism of their dimension groups $[h]: G\to H$;  but $G\cong P({\cal F}_{\psi})$
and $H\cong P({\cal F}_{\varphi})$. From inclusion $P({\cal F}_{\varphi})\subseteq P({\cal F}_{\psi})$,
one gets ${\Bbb Q}(\lambda_{B'})\cong P({\cal F}_{\varphi})\otimes {\Bbb Q}\subseteq 
P({\cal F}_{\psi})\otimes {\Bbb Q}\cong {\Bbb Q}(\lambda_B)\subseteq K_f$ and, therefore, 
${\Bbb Q}(\lambda_{B'})\subseteq K_{\psi}$. Lemma \ref{lem4.3.3} follows.
$\square$

\begin{cor}\label{cor4.3.1}
If $h: {\Bbb A}_{\psi}\to {\Bbb A}_{\varphi}$ is a stable homomorphism, then
$Aut~({\Bbb Q}(\lambda_{B'}))$ (or,  $Gal~({\Bbb A}_{\varphi})$)  is a subgroup 
(or, a  normal subgroup)  of $Gal~({\Bbb A}_{\psi})$.
\end{cor}
{\it Proof.} 
The (Galois) subfields of the Galois field $K_{\psi}$ are bijective with the (normal) subgroups 
of the group $Gal~(K_{\psi})$, see e.g. [Morandi 1996]  \cite{MO}.
$\square$

 \index{tight hyperbolic matrix}

\bigskip\noindent
Let $T^m\cong {\Bbb R}^m/{\Bbb Z}^m$ be an $m$-dimensional torus; let
$\psi_0$ be  a  $m\times m$ integer matrix with $det~(\psi_0)=1$, such that it
is similar to a positive matrix.
The  matrix $\psi_0$ defines  a linear transformation of ${\Bbb R}^m$,
which preserves the lattice $L\cong {\Bbb Z}^m$ of points with integer 
coordinates. There is an induced diffeomorphism $\psi$ of the quotient $T^m\cong {\Bbb R}^m/{\Bbb Z}^m$
onto itself; this diffeomorphism $\psi: T^m\to T^m$ has a fixed point $p$ corresponding 
to the origin of ${\Bbb R}^m$. 
Suppose that $\psi_0$ is hyperbolic, i.e. there are no eigenvalues of $\psi_0$ at the
unit circle; then $p$ is a hyperbolic fixed point of $\psi$ and the stable manifold 
$W^s(p)$ is the image of the corresponding eigenspace of $\psi_0$ under the 
projection ${\Bbb R}^m\to T^m$. If $codim~W^s(p)=1$,  the hyperbolic linear transformation
 $\psi_0$  (and the diffeomorphism $\psi$)  will be  called {\it tight}. 
\begin{lem}\label{lem4.3.4}
The tight hyperbolic matrix  $\psi_0$ is similar to the matrix $B$ of the fundamental AF-algebra 
${\Bbb A}_{\psi}$.
\end{lem}
{\it Proof.} Since $H_k(T^m; {\Bbb R})\cong {\Bbb R}^{{m!\over k!(m-k)!}}$,
one gets $H_{m-1}(T^m; {\Bbb R})\cong {\Bbb R}^m$; in view of the Poincar\'e 
duality, $H^1(T^m; {\Bbb R})=H_{m-1}(T^m; {\Bbb R})\cong {\Bbb R}^m$. 
Since $codim~W^s(p)=1$, measured foliation ${\cal F}_{\psi}$ is given by a closed
form $\omega_{\psi}\in H^1(T^m; {\Bbb R})$, such that $[\psi]\omega_{\psi}=\lambda_{\psi}\omega_{\psi}$,
where $\lambda_{\psi}$ is the eigenvalue of the linear transformation 
$[\psi]: H^1(T^m; {\Bbb R})\to H^1(T^m; {\Bbb R})$.  It is easy to see that $[\psi]={\psi}_0$, because $H^1(T^m; {\Bbb R})\cong {\Bbb R}^m$
is the universal cover for $T^m$, where  the eigenspace $W^u(p)$ of ${\psi}_0$ is the span of the
eigenform $\omega_{\psi}$. 
On the other hand, from the proof of Theorem  \ref{thm4.3.2} we know that the Plante group $P({\cal F}_{\psi})$ 
is generated by the coordinates of vector $\omega_{\psi}$;  the matrix $B$ is nothing but the matrix
$\psi_0$ written in a new basis of $P({\cal F}_{\psi})$.  
Each change of basis in the ${\Bbb Z}$-module $P({\cal F}_{\psi})$ is given by an integer invertible matrix $S$;
therefore, $B=S^{-1}\psi_0S$. Lemma \ref{lem4.3.4} follows.
$\square$

\bigskip\noindent
Let $\psi: T^m\to T^m$ be a hyperbolic diffeomorphism;
the mapping torus $T^m_{\psi}$ will be called a (hyperbolic) {\it torus bundle}
of dimension $m$. Let $k=|Gal~({\Bbb A}_{\psi})|$;  it follows from the Galois theory,
that $1<k\le m!$. Denote $t_i$ the cardinality of a subgroup $G_i\subseteq Gal~({\Bbb A}_{\psi})$.
\begin{cor}\label{cor4.3.2}
There are no (non-trivial) continuous map $T^m_{\psi}\to T^{m'}_{\varphi}$,  whenever $t_i'\nmid k$
for all $G_i'\subseteq Gal~({\Bbb A}_{\varphi})$. 
\end{cor}
{\it Proof.} 
If $h: T^m_{\psi}\to T^{m'}_{\varphi}$ was a continuous map to a torus bundle of dimension $m'<m$, then,
by Theorem \ref{thm4.3.3} and Corollary \ref{cor4.3.1}, the $Aut~({\Bbb Q}(\lambda_{B'}))$
(or, $Gal~({\Bbb A}_{\varphi})$) were a non-trivial
subgroup (or, normal subgroup)  of the group $Gal~({\Bbb A}_{\psi})$; since $k=|Gal~({\Bbb A}_{\psi})|$, one 
concludes that one of $t_i'$ divides $k$.   This contradicts our assumption.
$\square$
\begin{dfn}
The torus bundle $T^m_{\psi}$ is called robust, if there exists $m'<m$, such that
no continuous map $T_{\psi}^m\to T_{\varphi}^{m'}$ is possible.  
\end{dfn}
\begin{rmk}
\textnormal{
Are there robust bundles? It is shown in this section, that for $m=2,3$ and $4$
there are infinitely many robust  bundles.  
}
\end{rmk}

 \index{robust torus bundle}

\subsubsection{Case $m=2$}
This case is trivial; $\psi_0$ is a hyperbolic matrix and always tight. 
The $char~(\psi_0)=char~(B)$ is an irreducible quadratic polynomial with two
real roots; $Gal~({\Bbb A}_{\psi})\cong {\Bbb Z}_2$ and, therefore, 
$|Gal~({\Bbb A}_{\psi})|=2$.  Formally, $T^2_{\psi}$ is robust, since no torus 
bundle of a smaller dimension is defined.

\subsubsection{Case $m=3$}
The $\psi_0$ is hyperbolic; it is always tight, since one root of $char~(\psi_0)$
is real and isolated inside or outside the unit circle.  
\begin{cor}\label{cor4.3.3}
Let 
\displaymath
\psi_0(b,c)=\left(\matrix{-b & 1 & 0\cr
                  -c & 0 & 1\cr
                  -1 & 0 & 0}\right)
\enddisplaymath
be such, that $char~(\psi_0(b,c))=x^3+bx^2+cx+1$ is irreducible and $-4b^3+b^2c^2+18bc-4c^3-27$ is 
the square of an integer; then $T^3_{\psi}$ admits no continuous map to any $T_{\varphi}^2$.  
\end{cor}
{\it Proof.} The $char~(\psi_0(b,c))=x^3+bx^2+cx+1$ and the discriminant
 $D=-4b^3+b^2c^2+18bc-4c^3-27$. By  [Morandi 1996]  \cite{MO}, Theorem 13.1,  we have 
 $Gal~({\Bbb A}_{\psi})\cong {\Bbb Z}_3$
and, therefore, $k=|Gal~({\Bbb A}_{\psi})|=3$. 
For $m'=2$, it was shown that $Gal~({\Bbb A}_{\varphi})\cong {\Bbb Z}_2$ and,
therefore, $t_1'=2$. Since  $2\nmid 3$, Corollary \ref{cor4.3.2} says that
no continuous map $T_{\psi}^3\to T_{\varphi}^2$ can be constructed.
Corollary \ref{cor4.3.3} follows.  
$\square$

\begin{exm}
\textnormal{
There are infinitely many matrices $\psi_0(b,c)$ satisfying the assumptions of Corollary
\ref{cor4.3.3};   below are a few numerical examples of robust bundles:
\displaymath 
\left(\matrix{0 & 1 & 0\cr
                  3 & 0 & 1\cr
                  -1 & 0 & 0}\right),
\quad
\left(\matrix{1 & 1 & 0\cr
                  2 & 0 & 1\cr
                  -1 & 0 & 0}\right),
\quad
\left(\matrix{2 & 1 & 0\cr
                  1 & 0 & 1\cr
                  -1 & 0 & 0}\right),
\quad
\left(\matrix{3 & 1 & 0\cr
                  0 & 0 & 1\cr
                  -1 & 0 & 0}\right).
\enddisplaymath
}
\end{exm}
\begin{rmk}
\textnormal{
Notice  that the above matrices are not pairwise similar; it can be gleaned from their
traces;   thus  they represent topologically distinct torus bundles.
}
\end{rmk}

\subsubsection{Case $m=4$}
Let $p(x)=x^4+ax^3+bx^2+cx+d$ be a quartic.  Consider the associated cubic polynomial
$r(x)=x^3-bx^2+(ac-4d)x+4bd-a^2d-c^2$;  denote by $L$ the splitting field of $r(x)$. 
\begin{cor}\label{cor4.3.4}
Let 
\displaymath
\psi_0(a,b,c)=\left(\matrix{-a & 1 & 0 & 0\cr
                  -b & 0 & 1 & 0\cr
                  -c & 0 & 0 & 1\cr
                  -1 & 0 & 0 & 0}\right)
\enddisplaymath
be tight and such, that $char~(\psi_0(a,b,c))=x^4+ax^3+bx^2+cx+1$ is irreducible and 
one of the following holds: (i) $L={\Bbb Q}$;
(ii) $r(x)$ has a unique root $t\in {\Bbb Q}$ and $h(x)=(x^2-tx+1)[x^2+ax+(b-t)]$
splits over $L$; (iii) $r(x)$ has a unique root $t\in {\Bbb Q}$ and $h(x)$ does not split over $L$.
Then $T_{\psi}^4$ admits no continuous map to any $T_{\varphi}^3$  with $D>0$.  
\end{cor}
{\it Proof.}
According to [Morandi 1996]  \cite{MO},  Theorem 13.4,   $Gal~({\Bbb A}_{\psi})\cong {\Bbb Z}_2\oplus {\Bbb Z}_2$
in case (i); $Gal~({\Bbb A}_{\psi})\cong {\Bbb Z}_4$ in case (ii); and $Gal~({\Bbb A}_{\psi})\cong D_4$
(the dihedral group) in case (iii). Therefore, $k=|{\Bbb Z}_2\oplus {\Bbb Z}_2|=|{\Bbb Z}_4|=4$
or $k=|D_4|=8$. On the other hand, for $m'=3$ with $D>0$ (all roots are real), we have
$t_1'=|{\Bbb Z}_3|=3$ and $t_2'=|S_3|=6$. Since $3; 6\nmid 4;8$, corollary \ref{cor4.3.2} says
that continuous map $T_{\psi}^4\to T_{\varphi}^3$ is impossible. Corollary \ref{cor4.3.4} follows. 
$\square$

\begin{exm}
\textnormal{
There are infinitely many matrices $\psi_0$,  which satisfy the assumption of 
corollary \ref{cor4.3.4};  indeed, consider a family
\displaymath
\psi_0(a,c)=\left(\matrix{-2a & 1 & 0 & 0\cr
                  -a^2-c^2 & 0 & 1 & 0\cr
                  -2c & 0 & 0 & 1\cr
                  -1 & 0 & 0 & 0}\right),
 \enddisplaymath
where $a, c\in {\Bbb Z}$. The associated cubic  becomes $r(x)=x[x^2-(a^2+c^2)x+4(ac-1)]$,
so that $t=0$ is a rational root; then $h(x)=(x^2+1)[x^2+2ax+a^2+c^2]$. 
The matrix $\psi_0(a,c)$ satisfies one of the conditions (i)-(iii) of 
corollary \ref{cor4.3.4}  for each $a,c\in {\Bbb Z}$; it remains to 
eliminate the (non-generic) matrices, which are not tight or irreducible.
Thus, $\psi_0(a,c)$ defines  a family of topologically distinct robust bundles.
}
\end{exm}

\vskip1cm\noindent
{\bf Guide to the literature.}
The Anosov diffeomorphisms were introduced and studied in [Anosov 1967]  \cite{Ano1};
for a classical  account  of the differentiable dynamical systems  see  [Smale  1967]
\cite{Sma1}.  An excellent survey of  foliations has been compiled by [Lawson  1974]  \cite{Law1}.  
The Galois theory is covered in the textbook by [Morandi 1996]  \cite{MO}.  
The original proof of Theorem \ref{thm4.3.3} and obstruction theory for Anosov's bundles   can be found 
 in \cite{Nik9}.

\section*{Exercises,  problems and conjectures}

\begin{enumerate}
\item
Verify that $F: \phi\mapsto {\Bbb A}_{\phi}$ is a well-defined function on the set of  all 
Anosov automorphisms given by  the hyperbolic matrices with the non-negative entries.

\item
Verify that the definition of the AF-algebra   ${\Bbb A}_{\phi}$  for the pseudo-Anosov maps  coincides 
with the one for the Anosov maps.  (Hint:  the Jacobi-Perron fractions  of dimension $n=2$ coincide with the
regular continued fractions.)

\item
{\bf $p$-adic invariants of  pseudo-Anosov maps. }
Let $\phi\in Mod~(X)$ be pseudo-Anosov automorphism
of a surface $X$. If $\lambda_{\phi}$ is the dilatation of
$\phi$,   then one can consider a ${\Bbb Z}$-module 
${\goth m}={\Bbb Z}v^{(1)}+\dots+{\Bbb Z}v^{(n)}$ in the
number field $K={\Bbb Q}(\lambda_{\phi})$ generated by the
normalized eigenvector $(v^{(1)},\dots,v^{(n)})$ corresponding
to the eigenvalue $\lambda_{\phi}$. The trace function on the number field $K$
gives rise to a symmetric bilinear form $q(x,y)$ on the 
module ${\goth m}$. The form is defined over the field ${\Bbb Q}$. 
It has been shown that a pseudo-Anosov automorphism $\phi'$,
conjugate to $\phi$, yields a form $q'(x,y)$, equivalent to
$q(x,y)$, i.e. $q(x,y)$ can be transformed to $q'(x,y)$ by
an invertible linear substitution with the coefficients in ${\Bbb Z}$.
It is well known that   two rational bilinear forms
$q(x,y)$ and  $q'(x,y)$ are equivalent  whenever 
the following conditions are satisfied:

\medskip
(i) $\Delta=\Delta'$,  where $\Delta$ is the determinant of the form;

\smallskip
(ii) for each prime number $p$ (including $p=\infty$) certain $p$-adic
equation between the coefficients of forms $q, q'$ must be satisfied,
see e.g. [Borevich \& Shafarevich 1966] \cite{BS}, Chapter 1, \S 7.5. 
(In fact, only a {\it finite}   number of such equations have to be verified.)

\medskip\noindent
Condition (i) has been already used to discern between the 
conjugacy classes of the  pseudo-Anosov automorphisms. One can use
condition (ii) to discern between the pseudo-Anosov automorphisms
with $\Delta=\Delta'$;  in other words,  one gets a problem:

\medskip
\centerline{{\it To define $p$-adic invariants of the pseudo-Anosov maps.}}

\item
{\bf  The signature of a pseudo-Anosov map.}
The signature is an important and well-known invariant connected to
the chirality and knotting number of knots and links, see e.g.
[Reidemeister 1932]   \cite{R}.  It will be interesting to find a geometric interpretation
of  the signature $\Sigma$ for  the pseudo-Anosov automorphisms;  
one can ask the following  question:

\medskip
\centerline{{\it To find  geometric interpretation of the invariant $\Sigma$.}} 

 \index{signature of pseudo-Anosov map}

\item
{\bf  The number of conjugacy classes of  pseudo-Anosov maps
with the  same dilatation.}
The dilatation $\lambda_{\phi}$ is an invariant of the conjugacy class of
the pseudo-Anosov automorphism $\phi\in Mod~(X)$. On the other hand, 
it is known that there exist non-conjugate pseudo-Anosov's with 
the same dilatation and the number of such classes is finite, 
see  [Thurston 1988]  \cite{Thu1},  p.428.
It is natural to expect that the invariants of operator algebras
can be used to evaluate the number;  we have the following
\begin{cnj}
Let $(\Lambda, [I], K)$ be the triple corresponding to a pseudo-Anosov
map  $\phi\in Mod~(X)$.  Then the number of the conjugacy classes
of the pseudo-Anosov automorphisms with the dilatation $\lambda_{\phi}$
is equal to the class number $h_{\Lambda}=|\Lambda/[I]|$ 
of the integral order $\Lambda$. 
\end{cnj}

 \end{enumerate}





\chapter{Algebraic Geometry}
The    NCG-valued functors arise  in algebraic geometry;
what is going on conceptually?    Remember the covariant  functor $GL_n:$  {\bf CRng} $\to$ {\bf Grp} 
from  the category of commutative rings $R$  to the  category of groups;  functor $GL_n$ 
produces a multiplicative group of  all $n\times n$ invertible matrices with entries in $R$       
and preserves homomorphisms between the objects in the respective categories,  see
Example \ref{exm2.3.4}.   The NCG-valued functors take one step further:   they deal with 
the natural embedding {\bf Grp} $\hookrightarrow$ {\bf Grp-Rng},   where {\bf Grp-Rng} 
is the category of associative group rings;  thus we have 
\displaymath
\hbox{{\bf CRng}}
\buildrel\rm GL_n\over\longrightarrow
\hbox{{\bf Grp}}\hookrightarrow \hbox{{\bf Grp-Rng}}, 
\enddisplaymath
where {\bf Grp-Rng} is an associative ring,  i.e. the NCG.   (Of course this simple observation 
would be of little use  if the objects in  {\bf Grp-Rng} was fuzzy and nothing concrete can be said 
about them;  note also that the abelianization {\bf Grp-Rng}$/[\bullet, \bullet]$  of {\bf Grp-Rng} 
is naturally isomorphic to {\bf CRng}.)   For $n=2$ and {\bf CRng} being the coordinate ring 
of elliptic curves,  the category {\bf Grp-Rng} consists of the noncommutative tori (with scaled units);  
this fact will be proved in Section 5.1  in two independent ways.   
For the higher genus algebraic curves  the category {\bf Grp-Rng} consists of
the so-called {\it toric} AF-algebras,  see Section 5.2.   The general case of complex 
projective varieties  is considered in Section 5.3 and {\bf Grp-Rng}  consists of the
{\it Serre $C^*$-algebras}.  In Section 5.4  we use  the stable isomorphism  group
of toric AF-algebras  to prove Harvey's  conjecture  on the  linearity of the mapping   
class groups.
 
 \index{toric AF-algebra}
 \index{Serre $C^*$-algebra}
 \index{elliptic curve}

\section{Elliptic curves}
Let us repeat  some  known facts.   We will be  working with the  ground field 
of complex numbers ${\Bbb C}$;  by an {\it elliptic curve} we shall understand  
the subset of the complex projective plane  of the form
\displaymath
{\cal E}({\Bbb C})=\{(x,y,z)\in {\Bbb C}P^2 ~|~ y^2z=4x^3+axz^2+bz^3\},
\enddisplaymath
where $a$ and  $b$  are some constant complex numbers.  One can visualize 
the  real points of  ${\cal E}({\Bbb C})$  as  it is shown  in Figure 5.1.
\begin{figure}[here]
\begin{picture}(300,100)(0,-10)

\put(100,50){\line(1,0){60}}
\put(230,50){\line(1,0){60}}
\put(130,20){\line(0,1){60}}
\put(260,20){\line(0,1){60}}

\thicklines
\qbezier(165,25)(115,50)(165,75)
\qbezier(110,50)(117,15)(130,50)
\qbezier(110,50)(117,85)(130,50)

\put(255,50){\oval(30,30)[l]}
\qbezier(254,65)(285,56)(290,85)
\qbezier(254,35)(285,44)(290,15)

\put(115,0){$a<0$}
\put(245,0){$a>0$}
\end{picture}
\caption{The real points of an affine elliptic curve $y^2=4x^3+ax$.}
\end{figure}
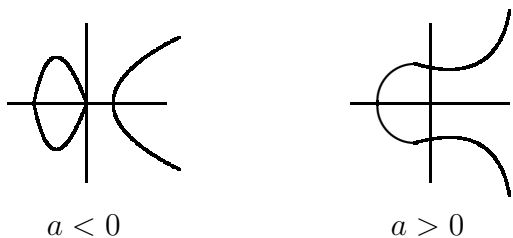
\begin{rmk}
\textnormal{
It is known  that each elliptic curve ${\cal E}({\Bbb C})$  is isomorphic to the set
of points of intersection of two {\it quadric surfaces} in the complex projective space ${\Bbb C}P^3$
given by the  system of homogeneous equations 
\displaymath
\left\{
\begin{array}{ccc}
u^2+v^2+w^2+z^2 &=&  0,\\
Av^2+Bw^2+z^2  &=&  0,   
\end{array}
\right.
\enddisplaymath
where $A$ and $B$ are some constant  complex numbers and  
$(u,v,w,z)\in {\Bbb C}P^3$;  the system  is called
the {\it Jacobi form} of elliptic curve ${\cal E}({\Bbb C})$.   
}
\end{rmk}
 \index{Jacobi elliptic curve}
\begin{dfn}
 By a  complex torus  one understands  the space  ${\Bbb C}/({\Bbb Z}\omega_1+{\Bbb Z}\omega_2)$,  
 where  $\omega_1$ and $\omega_2$ are   linearly independent
vectors in the complex plane ${\Bbb C}$,  see Fig. 5.2;   the ratio $\tau=\omega_2/\omega_1$ is called a  
complex modulus.
 \end{dfn}
 \index{complex torus}
\begin{figure}[here]
\begin{picture}(300,100)(0,0)

\put(90,90){\line(1,0){80}}
\put(80,70){\line(1,0){80}}
\put(60,30){\line(1,0){80}}

\put(85,20){\line(1,2){40}}
\put(65,20){\line(1,2){40}}
\put(105,20){\line(1,2){40}}
\put(125,20){\line(1,2){40}}

\put(180,50){\vector(1,0){40}}

\thicklines
\put(70,50){\line(1,0){80}}
\put(100,20){\line(0,1){80}}

\put(100,50){\vector(1,0){20}}
\put(100,50){\vector(1,2){10}}
\put(110,70){\line(1,0){20}}
\put(120,50){\line(1,2){10}}


\put(290,50){\oval(60,40)}
\qbezier(280,53)(290,38)(300,53)
\qbezier(285,48)(290,53)(295,48)


\put(260,90){${\Bbb C}/({\Bbb Z}+{\Bbb Z}\tau)$}
\put(70,90){${\Bbb C}$}
\put(120,35){$1$}
\put(105,75){${\tau}$}

\put(185,55){factor}
\put(190,40){map}

\end{picture}
\caption{Complex torus  ${\Bbb C}/({\Bbb Z}+{\Bbb Z}\tau)$.}
\end{figure}
\begin{rmk}
\textnormal{
Two complex  tori   ${\Bbb C}/({\Bbb Z}+{\Bbb Z}\tau)$ and   ${\Bbb C}/({\Bbb Z}+{\Bbb Z}\tau')$ 
are isomorphic   if and only if 
\displaymath
\tau'={a\tau+b\over c\tau+d} \quad \hbox{for some matrix} 
 \quad\left(\matrix{a & b\cr c & d}\right)  \in SL_2({\Bbb Z}).
 \enddisplaymath
 }
\end{rmk}

\bigskip\noindent
The complex analytic manifold ${\Bbb C}/({\Bbb Z}+{\Bbb Z}\tau)$
can be embedded into an  $n$-dimensional complex  projective space as an algebraic variety.
For $n=2$  we have the following  classical result,  which relates complex torus    
${\Bbb C}/({\Bbb Z}+{\Bbb Z}\tau)$  with an elliptic
curve ${\cal E}({\Bbb C})$ in the projective plane ${\Bbb C}P^2$.  
\begin{thm}\label{thm5.1.1}
{\bf (Weierstrass)}
There exists a holomorphic embedding  $${\Bbb C}/({\Bbb Z}+{\Bbb Z}\tau) \hookrightarrow  {\Bbb C}P^2$$ given by the formula
\displaymath
z\mapsto
\cases{(\wp(z), \wp'(z), 1) & \hbox{for}  $z\not\in L_{\tau}:={\Bbb Z}+{\Bbb Z}\tau$,\cr
          (0,1,0) & \hbox{for} $z\in L_{\tau}$},
\enddisplaymath
which is an isomorphism  between  complex torus ${\Bbb C}/({\Bbb Z}+{\Bbb Z}\tau)$ and elliptic curve 
\displaymath
{\cal E}({\Bbb C})=\{(x,y,z)\in {\Bbb C}P^2 ~|~ y^2z=4x^3+axz^2+bz^3\},
\enddisplaymath
where $\wp(z)$ is  the  Weierstrass  function defined by the convergent series
\displaymath
\wp(z)={1\over z^2}+\sum_{\omega\in L_{\tau}-\{0\}} \left({1\over (z-\omega)^2}-{1\over\omega^2}\right).
\enddisplaymath
and
\displaymath
\left\{
\begin{array}{ccc}
a &=&  -60\sum_{\omega\in  L_{\tau}-\{0\}} {1\over\omega^4},\\
b &=&  -140\sum_{\omega\in L_{\tau}-\{0\}} {1\over\omega^6}.   
\end{array}
\right.
\enddisplaymath
\end{thm}
 \index{Weierstrass $\wp$ function}
\begin{rmk}
\textnormal{
The Weierstrass Theorem   identifies elliptic curves ${\cal E}({\Bbb C})$ 
and  complex tori  ${\Bbb C}/({\Bbb Z}\omega_1+{\Bbb Z}\omega_2)$;
we shall write ${\cal E}_{\tau}$ to denote  elliptic curve corresponding to the
complex torus of modulus $\tau=\omega_2/\omega_1$. 
}
\end{rmk}
\begin{dfn}
 By {\bf Ell}  we shall mean the category of all elliptic curves   ${\cal E}_{\tau}$;
 the arrows of {\bf Ell}  are identified with the isomorphisms between 
 elliptic curves ${\cal E}_{\tau}$.  We shall write {\bf NC-Tor} to denote the
  category of all noncommutative tori  ${\cal A}_{\theta}$;
 the arrows of {\bf NC-Tor}  are identified with the stable isomorphisms
 (Morita equivalences)  between  noncommutative tori  ${\cal A}_{\theta}$.  
 \end{dfn}
\begin{figure}[here]
\begin{picture}(300,110)(-120,-5)
\put(20,70){\vector(0,-1){35}}
\put(130,70){\vector(0,-1){35}}
\put(45,23){\vector(1,0){60}}
\put(45,83){\vector(1,0){60}}
\put(15,20){${\cal A}_{\theta}$}
\put(0,50){$F$}
\put(145,50){$F$}
\put(123,20){${\cal A}_{\theta'={a\theta+b\over c\theta+d}}$}
\put(17,80){${\cal E}_{\tau}$}
\put(122,80){${\cal E}_{\tau'={a\tau+b\over c\tau+d}}$}
\put(60,30){\sf stably}
\put(50,10){\sf  isomorphic}
\put(50,90){\sf isomorphic}
\end{picture}
\caption{Functor on elliptic curves.}
\end{figure}
\begin{thm}\label{thm5.1.2}
{\bf (Functor on elliptic curves)}
There exists a covariant functor 
\displaymath
F:  \hbox{{\bf Ell}} \longrightarrow \hbox{{\bf NC-Tor}},
\enddisplaymath
which maps isomorphic elliptic curves ${\cal E}_{\tau}$ to the 
stably isomorphic (Morita equivalent) noncommutative tori ${\cal A}_{\theta}$,
see Fig. 5.3;   the functor $F$ is non-injective and $Ker~F\cong (0,\infty)$. 
\end{thm}
Theorem \ref{thm5.1.2} will be proved in Section 5.1.1 using the Sklyanin algebras
and in Section 5.1.2 using measured foliations and the Teichm\"uller theory.

 \index{Sklyanin algebra}
\subsection{Noncommutative tori  via Sklyanin algebras}
\begin{dfn}
{\bf ([Sklyanin  1982]  \cite{Skl1})}
By   the Sklyanin algebra   $S(\alpha,\beta,\gamma)$  one understands   a free   
${\Bbb C}$-algebra   on   four  generators   $x_1,\dots,x_4$   and  six  quadratic relations
\displaymath
\left\{
\begin{array}{ccc}
x_1x_2-x_2x_1 &=& \alpha(x_3x_4+x_4x_3),\\
x_1x_2+x_2x_1 &=& x_3x_4-x_4x_3,\\
x_1x_3-x_3x_1 &=& \beta(x_4x_2+x_2x_4),\\
x_1x_3+x_3x_1 &=& x_4x_2-x_2x_4,\\
x_1x_4-x_4x_1 &=& \gamma(x_2x_3+x_3x_2),\\ 
x_1x_4+x_4x_1 &=& x_2x_3-x_3x_2,
\end{array}
\right.
\enddisplaymath
where $\alpha+\beta+\gamma+\alpha\beta\gamma=0$. 
\end{dfn}
 \index{twisted homogeneous coordinate ring}
\begin{rmk}
{\bf ([Smith \& Stafford  1992]  \cite{SmiSta1},  p. 260)}
\textnormal{
The algebra $S(\alpha,\beta,\gamma)$ is isomorphic to a 
{\it  twisted homogeneous  coordinate ring}   of elliptic curve ${\cal E}_{\tau}\subset {\Bbb C}P^3$ 
 given in its  Jacobi form
\displaymath
\left\{
\begin{array}{ccc}
u^2+v^2+w^2+z^2 &=&  0,\\
{1-\alpha\over 1+\beta}v^2+{1+\alpha\over 1-\gamma}w^2+z^2  &=&  0,
\end{array}
\right.
\enddisplaymath
i.e.  $S(\alpha,\beta,\gamma)$ satisfies an  isomorphism
\displaymath
\hbox{{\bf Mod}}~(S(\alpha,\beta,\gamma))/\hbox{{\bf Tors}}\cong \hbox{{\bf Coh}}~({\cal E}_{\tau}),
\enddisplaymath
 where {\bf Coh} is  the category of quasi-coherent sheaves on ${\cal E}_{\tau}$, 
  {\bf Mod}  the category of graded left modules over the graded ring $S(\alpha,\beta,\gamma)$
 and  {\bf Tors}  the full sub-category of {\bf Mod} consisting of the
torsion modules,  see [Serre 1955]   \cite{Ser1}.  
The algebra $S(\alpha,\beta,\gamma)$ defines a natural {\it automorphism} 
$\sigma: {\cal E}_{\tau}\to {\cal E}_{\tau}$ of the elliptic curve ${\cal E}_{\tau}$,  
see e.g.   [Stafford \& van ~den ~Bergh  2001]  \cite{StaVdb1}, p. 173. 
}
\end{rmk}
\begin{lem}\label{lem5.1.1}
If $\sigma^4=Id$,  then  algebra  $S(\alpha,\beta,\gamma)$
is isomorphic to  a free algebra ${\Bbb C}\langle x_1,x_2,x_3,x_4\rangle$  modulo an ideal 
generated  by  six   skew-symmetric quadratic  relations
\displaymath
\left\{
\begin{array}{cc}
x_3x_1 &= \mu e^{2\pi i\theta}x_1x_3,\\
x_4x_2 &= {1\over \mu} e^{2\pi i\theta}x_2x_4,\\
x_4x_1 &= \mu e^{-2\pi i\theta}x_1x_4,\\
x_3x_2 &= {1\over \mu} e^{-2\pi i\theta}x_2x_3,\\
x_2x_1 &= x_1x_2,\\
x_4x_3 &= x_3x_4,
\end{array}
\right.
\enddisplaymath
where $\theta\in S^1$ and $\mu\in (0,\infty)$. 
\end{lem}
{\it Proof.}  
(i)  Since  $\sigma^4=Id$,   the Sklyanin algebra $S(\alpha,\beta,\gamma)$
is isomorphic to a free   algebra  ${\Bbb C}\langle x_1,x_2,x_3,x_4\rangle$
modulo an ideal generated by   the  skew-symmetric relations 
\displaymath
\left\{
\begin{array}{ccc}
x_3x_1 &=& q_{13} x_1x_3,\\
x_4x_2 &=&  q_{24}x_2x_4,\\
x_4x_1 &=&  q_{14}x_1x_4,\\
x_3x_2 &=&  q_{23}x_2x_3,\\
x_2x_1&=&  q_{12}x_1x_2,\\
x_4x_3&=&  q_{34}x_3x_4,
\end{array}
\right.
\enddisplaymath
where $q_{ij}\in {\Bbb C}\setminus\{0\}$,  see  [Feigin \& Odesskii  1989]  \cite{FeOd1},  Remark 1
and   [Feigin \& Odesskii  1993]  \cite{FeOd2}, \S 2 for the proof.

\bigskip
(ii) It is verified directly,  that above relations  are invariant of the involution 
$x_1^*=x_2, x_3^*=x_4$,     if and only if,  the following restrictions on the constants $q_{ij}$ hold
\displaymath
\left\{
\begin{array}{ccc}
q_{13} &=&  (\bar q_{24})^{-1},\\
q_{24} &=&  (\bar q_{13})^{-1},\\
q_{14} &= & (\bar q_{23})^{-1},\\
q_{23} &= & (\bar q_{14})^{-1},\\
q_{12} &= & \bar q_{12},\\
q_{34} &= & \bar q_{34},
\end{array}
\right.
\enddisplaymath
where $\bar q_{ij}$ means the complex conjugate of $q_{ij}\in {\Bbb C}\setminus\{0\}$.     
\begin{rmk}
\textnormal{
The skew-symmetric relations invariant of the involution $x_1^*=x_2, x_3^*=x_4$
define an involution on the Sklyanin algebra;  we shall call such an algebra
a {\it Sklyanin $\ast$-algebra}.
}
\end{rmk}

 \index{Sklyanin $\ast$-algebra}

\bigskip
(iii)
Consider a one-parameter family  $S(q_{13})$ of the Sklyanin $\ast$-algebras
defined by the following additional constraints
\displaymath
\left\{
\begin{array}{ccc}
q_{13} &=& \bar  q_{14},\\
q_{12} &=&  q_{34}=1.
\end{array}
\right.
\enddisplaymath
It is not hard to see,  that the $\ast$-algebras  $S(q_{13})$ 
are pairwise non-isomorphic for different values of complex  parameter $q_{13}$;
therefore  the family $S(q_{13})$ is   a  normal form of  the Sklyanin $\ast$-algebra 
$S(\alpha,\beta,\gamma)$ with $\sigma^4=Id$. 
It remains to notice, that one can write  complex parameter $q_{13}$
in the polar form $q_{13}=\mu e^{2\pi i\theta}$,  where $\theta=Arg~(q_{13})$
and $\mu=|q_{13}|$.   Lemma \ref{lem5.1.1} follows.
$\square$

\begin{lem}\label{lem5.1.2}
The system of  relations 
\displaymath
\left\{
\begin{array}{cc}
x_3x_1  &= e^{2\pi i\theta}x_1x_3,\\
x_1x_2 &= x_2x_1 = e,\\
x_3x_4  &= x_4x_3 = e
\end{array}
\right.
\enddisplaymath
defining the noncommutative torus ${\cal A}_{\theta}$  is equivalent to  the 
following system  of  quadratic relations
\displaymath
\left\{
\begin{array}{cc}
x_3x_1 &=  e^{2\pi i\theta}x_1x_3,\\
x_4x_2 &=  e^{2\pi i\theta}x_2x_4,\\
x_4x_1 &=  e^{-2\pi i\theta}x_1x_4,\\
x_3x_2 &=   e^{-2\pi i\theta}x_2x_3,\\
x_2x_1 &= x_1x_2=e,\\
x_4x_3 &= x_3x_4=e.
\end{array}
\right.
\enddisplaymath
\end{lem}
{\it Proof.} 
Indeed, the first  and the two last equations  of both systems coincide; 
 we shall proceed stepwise for the rest  of the equations.

\smallskip
(i) Let us prove that equations for ${\cal A}_{\theta}$  imply  $x_1x_4=e^{2\pi i\theta}x_4x_1$. 
It follows from $x_1x_2=e$ and $x_3x_4=e$ that $x_1x_2x_3x_4=e$.  Since $x_1x_2=x_2x_1$ we can bring  the last 
equation to  the form  $x_2x_1x_3x_4=e$ and multiply the  both sides by the constant $e^{2\pi i\theta}$;
thus one gets the equation $x_2(e^{2\pi i\theta}x_1x_3)x_4=e^{2\pi i\theta}$.  
But $e^{2\pi i\theta}x_1x_3=x_3x_1$ and our main equation takes the form $x_2x_3x_1x_4= e^{2\pi i\theta}$. 

We can multiply on the left the  both sides of the equation by the element $x_1$
and thus get the equation $x_1x_2x_3x_1x_4= e^{2\pi i\theta}x_1$;  since $x_1x_2=e$ 
one  arrives at the equation  $x_3x_1x_4= e^{2\pi i\theta}x_1$.

Again one can multiply on the left the both sides  by the element 
$x_4$ and thus get the equation  $x_4x_3x_1x_4= e^{2\pi i\theta}x_4x_1$; since $x_4x_3=e$
one gets the required identity  $x_1x_4= e^{2\pi i\theta}x_4x_1$.

\smallskip
(ii) Let us prove that equations  for ${\cal A}_{\theta}$  imply  $x_2x_3=e^{2\pi i\theta}x_3x_2$. 
As in the case (i),  it follows from the equations $x_1x_2=e$ and $x_3x_4=e$ that $x_3x_4x_1x_2=e$.  Since $x_3x_4=x_4x_3$ 
we can bring  the last 
equation to  the form $x_4x_3x_1x_2=e$  and multiply the  both sides by the constant $e^{-2\pi i\theta}$;
thus one gets the equation $x_4(e^{-2\pi i\theta}x_3x_1)x_2=e^{-2\pi i\theta}$.  
But $e^{-2\pi i\theta}x_3x_1=x_1x_3$ and our main equation takes the form $x_4 x_1x_3 x_2= e^{-2\pi i\theta}$.

We can multiply on the left the  both sides of the equation by the element $x_3$
and thus get the equation $x_3 x_4 x_1 x_3 x_2= e^{-2\pi i\theta} x_3$; since $x_3 x_4 =e$ 
one  arrives at the equation  $ x_1 x_3 x_2= e^{-2\pi i\theta} x_3$.

Again one can multiply on the left the both sides  by the element 
$x_2$ and thus get the equation  $x_2 x_1 x_3 x_2= e^{-2\pi i\theta} x_2 x_3$; since $x_2 x_1=e$
one gets the  equation  $x_3 x_2= e^{-2\pi i\theta} x_2 x_3$.  Multiplying both sides
by constant $e^{2\pi i\theta}$ we obtain the required identity  $x_2x_3=e^{2\pi i\theta} x_3 x_2$.

\smallskip
(iii) Let us prove that equations for ${\cal A}_{\theta}$  imply  $x_4 x_2 = e^{2\pi i\theta}x_2 x_4$. 
Indeed, it was proved in  (i) that $x_1x_4= e^{2\pi i\theta}x_4x_1$;  we shall
multiply this equation on the right   by the equation $x_2x_1=e$. Thus one arrives
at  the equation $x_1 x_4 x_2 x_1= e^{2\pi i\theta} x_4 x_1 $. 

Notice that in the last equation one can cancel $x_1$ on the right thus bringing 
it to the simpler form  $x_1 x_4 x_2 = e^{2\pi i\theta} x_4$.

We shall multiply on the left both sides of the above equation by the element $x_2$;
one gets therefore $x_2 x_1 x_4 x_2 = e^{2\pi i\theta} x_2 x_4$.  But $x_2 x_1 = e$
and the left hand side simplifies giving  the required identity
  $x_4 x_2 = e^{2\pi i\theta} x_2x_4$.
  
\smallskip
Lemma \ref{lem5.1.2} follows.    
$\square$

 \index{basic isomorphism}

\begin{lem}\label{lem5.1.3}
{\bf (Basic isomorphism)}
The system of relations for noncommutative torus ${\cal A}_{\theta}$
\displaymath
\left\{
\begin{array}{cc}
x_3x_1 &=  e^{2\pi i\theta}x_1x_3,\\
x_4x_2 &=  e^{2\pi i\theta}x_2x_4,\\
x_4x_1 &=  e^{-2\pi i\theta}x_1x_4,\\
x_3x_2 &=   e^{-2\pi i\theta}x_2x_3,\\
x_2x_1 &= x_1x_2=e,\\
x_4x_3 &= x_3x_4=e, 
\end{array}
\right.
\enddisplaymath
is equivalent to the system of relations for the Sklyanin $\ast$-algebra
\displaymath
\left\{
\begin{array}{cc}
x_3x_1 &= \mu e^{2\pi i\theta}x_1x_3,\\
x_4x_2 &= {1\over \mu} e^{2\pi i\theta}x_2x_4,\\
x_4x_1 &= \mu e^{-2\pi i\theta}x_1x_4,\\
x_3x_2 &= {1\over \mu} e^{-2\pi i\theta}x_2x_3,\\
x_2x_1 &= x_1x_2,\\
x_4x_3 &= x_3x_4,
\end{array}
\right.
\enddisplaymath
modulo the following  ``scaled unit relation'' 
 \index{scaled unit}
\displaymath
x_1x_2=x_3x_4={1\over\mu}e.
\enddisplaymath
\end{lem}
{\it Proof.} 
(i) Using the last two relations,  one can bring the noncommutative torus 
relations  to the form
\displaymath
\left\{
\begin{array}{ccc}
x_3x_1x_4 &=&  e^{2\pi i\theta}x_1,\\
x_4 &= & e^{2\pi i\theta}x_2x_4x_1,\\
x_4x_1x_3 &=&  e^{-2\pi i\theta}x_1,\\
x_2 &=&  e^{-2\pi i\theta}x_4x_2x_3,\\
x_1x_2   &=&   x_2x_1   =e,\\
 x_3x_4    &=&  x_4x_3  =e.
\end{array}
\right.
\enddisplaymath

\bigskip
(ii)  The system of relations for the Sklyanin $\ast$-algebra complemented 
by the scaled unit relation,  i.e.  
\displaymath
\left\{
\begin{array}{cc}
x_3x_1 &= \mu e^{2\pi i\theta}x_1x_3,\\
x_4x_2 &= {1\over \mu} e^{2\pi i\theta}x_2x_4,\\
x_4x_1 &= \mu e^{-2\pi i\theta}x_1x_4,\\
x_3x_2 &= {1\over \mu} e^{-2\pi i\theta}x_2x_3,\\
x_2x_1 &= x_1x_2={1\over\mu}e,\\
x_4x_3 &= x_3x_4={1\over\mu}e
\end{array}
\right.
\enddisplaymath
is equivalent to the system
\displaymath
\left\{
\begin{array}{cc}
x_3x_1x_4 &= e^{2\pi i\theta}x_1,\\
x_4 &= e^{2\pi i\theta}x_2x_4x_1,\\
x_4x_1x_3 &= e^{-2\pi i\theta}x_1,\\
x_2 &= e^{-2\pi i\theta}x_4x_2x_3,\\
x_2x_1 &= x_1x_2={1\over\mu}e,\\
x_4x_3 &= x_3x_4={1\over\mu}e
\end{array}
\right.
\enddisplaymath
by using  multiplication and cancellation involving the last two equations.

\bigskip
(iii)  For each $\mu\in (0,\infty)$ consider a {\it scaled unit}  $e':={1\over\mu} e$ of
the Sklyanin $\ast$-algebra $S(q_{13})$ and the two-sided ideal $I_{\mu}\subset S(q_{13})$
generated by the relations $x_1x_2=x_3x_4=e'$.    Comparing the defining relations for
$S(q_{13})$ with such for the noncommutative torus ${\cal A}_{\theta}$,   one gets an
isomorphism 
\displaymath
S(q_{13})~/~I_{\mu}\cong {\cal A}_{\theta},
\enddisplaymath
see items (i) and (ii).  The isomorphism maps  generators $x_1,\dots,x_4$ of 
$\ast$-algebra $S(q_{13})$ to such  of the $C^*$-algebra ${\cal A}_{\theta}$ and 
the scaled unit $e'\in S(q_{13})$ to the {\it ordinary} unit of algebra ${\cal A}_{\theta}$. 
Lemma \ref{lem5.1.3} follows.  
$\square$

\bigskip
To finish the proof of Theorem \ref{thm5.1.2},   recall that  the Sklyanin $\ast$-algebra
$S(q_{13})$ satisfies the fundamental  isomorphism  
{\bf Mod}~$(S(q_{13}))/$ {\bf Tors} ~$\cong$~ {\bf Coh}~$({\cal E}_{\tau})$. 
Using the isomorphism $S(q_{13}) / I_{\mu}\cong {\cal A}_{\theta}$ established 
in Lemma \ref{lem5.1.3},    we conclude that 
\displaymath
I_{\mu}\backslash\hbox{{\bf Coh}}~({\cal E}_{\tau})\cong
\hbox{{\bf Mod}}~(I_{\mu}\backslash S(q_{13}))/\hbox{{\bf Tors}}\cong 
\hbox{{\bf Mod}}~({\cal A}_{\theta})/\hbox{{\bf Tors}}.  
\enddisplaymath
Thus  one gets an isomorphism {\bf Coh} $({\cal E}_{\tau})/I_{\mu}\cong$
{\bf Mod} $({\cal A}_{\theta})/$ {\bf Tors},  which  defines a  map
$F:$ {\bf Ell} $\to$ {\bf NC-Tor}.   Moreover,  map $F$ is a functor
because  isomorphisms in  the category {\bf Mod} 
$({\cal A}_{\theta})$  give rise to the stable isomorphisms (Morita equivalences)
in the category {\bf NC-Tor}.   
 The second part of Theorem \ref{thm5.1.2}  is due to the fact that $F$ 
forgets scaling of the unit,  i.e. for each $\mu\in (0,\infty)$  we have a constant
map 
\displaymath
S(q_{13})\ni  e':={1\over\mu}e\longmapsto e\in {\cal A}_{\theta}. 
\enddisplaymath
Thus $Ker~F\cong (0,\infty)$.   Theorem \ref{thm5.1.2} is proved.
 $\square$

 \index{measured foliation}

\subsection{Noncommutative tori  via measured foliations}
\begin{dfn}
{\bf ([Thurston  1988]  \cite{Thu1})}
By a measured foliation ${\cal F}$  on a surface $X$ one understands
   partition of $X$ into  the singular points $x_1,\dots,x_n$ of
order $k_1,\dots, k_n$ and the regular leaves, i.e. $1$-dimensional submanifolds of $X$; 
on each  open cover $U_i$ of $X\backslash \{x_1,\dots,x_n\}$ there exists a non-vanishing
real-valued closed 1-form $\phi_i$  such that: 

\medskip
(i)  $\phi_i=\pm \phi_j$ on $U_i\cap U_j$;

\smallskip
(ii) at each $x_i$ there exists a local chart $(u,v):V\to {\Bbb R}^2$
such that for $z=u+iv$, it holds $\phi_i=Im~(z^{k_i\over 2}dz)$ on
$V\cap U_i$ for some branch of $z^{k_i\over 2}$.

\medskip\noindent
The pair $(U_i,\phi_i)$ is called an atlas for  measured foliation ${\cal F}$.
A measure $\mu$ is assigned to each segment $(t_0,t)\in U_i$;   the measure  is  
transverse to the leaves of ${\cal F}$ and is defined by  the integral $\mu(t_0,t)=\int_{t_0}^t\phi_i$. 
Such a   measure is invariant along the leaves of ${\cal F}$,  hence the name. 
\end{dfn}
\begin{rmk}
\textnormal{
In  case   $X\cong T^2$  (a  torus)  each measured foliation is given 
by a family of parallel lines of a slope $\theta>0$ as shown in Fig. 5.4. 
}
\end{rmk}
\begin{figure}[here]
\begin{picture}(300,60)(-30,0)

\put(130,10){\line(1,0){40}}
\put(130,10){\line(0,1){40}}
\put(130,50){\line(1,0){40}}
\put(170,10){\line(0,1){40}}

\put(130,40){\line(2,1){20}}
\put(130,30){\line(2,1){40}}
\put(130,20){\line(2,1){40}}
\put(130,10){\line(2,1){40}}

\put(150,10){\line(2,1){20}}

\end{picture}

\caption{A measured foliation on the torus ${\Bbb R}^2/{\Bbb Z}^2$.}
\end{figure}
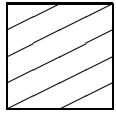

 \index{Teichm\"uller space}

\medskip\noindent
Let $T(g)$ be the Teichm\"uller space of  surface $X$ of genus $g\ge 1$,
i.e. the space of the complex structures on $X$. 
Consider the vector bundle $p: Q\to T(g)$ over $T(g)$,  whose fiber above a point 
$S\in T_g$ is the vector space $H^0(S,\Omega^{\otimes 2})$.   
Given a non-zero $q\in Q$ above $S$, we can consider the horizontal measured foliation
${\cal F}_q\in \Phi_X$ of $q$, where $\Phi_X$ denotes the space of the equivalence
classes of the measured foliations on $X$. If $\{0\}$ is the zero section of $Q$,
the above construction defines a map $Q-\{0\}\longrightarrow \Phi_X$. 
For any ${\cal F}\in\Phi_X$, let $E_{\cal F}\subset Q-\{0\}$ be the fiber
above ${\cal F}$. In other words, $E_{\cal F}$ is a subspace of the holomorphic 
quadratic forms,  whose horizontal trajectory structure coincides with the 
measured foliation ${\cal F}$. 
\begin{rmk}
\textnormal{
If ${\cal F}$ is a measured foliation with the simple zeroes (a generic case),  
then $E_{\cal F}\cong {\Bbb R}^n - 0$ and  $T(g)\cong {\Bbb R}^n$,  where $n=6g-6$ if
$g\ge 2$ and $n=2$ if $g=1$.
}
\end{rmk}
\begin{thm}\label{thm5.1.3}
{\bf ([Hubbard \& Masur  1979]  \cite{HuMa1})}
The restriction of $p$ to $E_{\cal F}$ defines a homeomorphism (an embedding)
$h_{\cal F}: E_{\cal F}\to T(g)$.
\end{thm}
\begin{cor}\label{cor5.1.1}
There exists a canonical  homeomorphism $h:\Phi_X\to T(g)-\{pt\}$,
where $pt= h_{\cal F}(0)$ and $\Phi_X\cong {\Bbb R}^n-0$ is the space of 
equivalence classes of  measured foliations ${\cal F}'$ on $X$. 
\end{cor}
{\it Proof.} 
Denote by  ${\cal F}'$ a vertical trajectory structure of  $q$. Since ${\cal F}$
and ${\cal F}'$ define $q$, and ${\cal F}=Const$ for all $q\in E_{\cal F}$, 
one gets a homeomorphism  between $T(g)-\{pt\}$ and  $\Phi_X$. 
Corollary \ref{cor5.1.1} follows. 
$\square$

\begin{rmk}
\textnormal{
The  homeomorphism $h:\Phi_X\to T(g)-\{pt\}$ depends on a foliation ${\cal F}$;
yet  there exists a  canonical homeomorphism $h=h_{\cal F}$
as follows.  Let $Sp ~(S)$ be the length spectrum of the Riemann surface
$S$ and $Sp ~({\cal F}')$ be the set positive reals $\inf \mu(\gamma_i)$,
where $\gamma_i$ runs over all simple closed curves, which are transverse 
to the foliation ${\cal F}'$. A canonical homeomorphism 
$h=h_{\cal F}: \Phi_X\to T(g)-\{pt\}$ is defined by the formula
$Sp ~({\cal F}')= Sp ~(h_{\cal F}({\cal F}'))$ for $\forall {\cal F}'\in\Phi_X$. 
}
\end{rmk}
Let $X\cong T^2$;   then  $T(1)\cong {\Bbb H}:=\{z=x+iy\in {\Bbb C} ~|~y>0\}$.
Since $q\ne 0$ there are no singular points and each $q\in H^0(S, \Omega^{\otimes 2})$
has the form $q=\omega^2$, where $\omega$ is a nowhere zero  holomorphic differential
on the complex torus $S$. 
Note that $\omega$ is just a constant times $dz$, and hence its vertical
trajectory structure is just a family of the parallel lines of a slope $\theta$,
see e.g.  [Strebel  1984]  \cite{S},  pp. 54--55.
Therefore,  $\Phi_{T^2}$ consists of the equivalence classes of the non-singular 
measured foliations on the two-dimensional torus.  It is well known (the Denjoy theory), 
that every such foliation  is measure equivalent to the foliation of a slope $\theta$ and a 
transverse measure $\mu>0$,  which is  invariant along the leaves of the foliation
Thus  one obtains  a canonical bijection 
\displaymath
h: \Phi_{T^2}\longrightarrow  {\Bbb H}-\{pt\}.
\enddisplaymath
\begin{dfn}
{\bf (Category of lattices)}
By a lattice in the complex plane ${\Bbb C}$ one understands
a triple $(\Lambda, {\Bbb C}, j)$,  where
$\Lambda\cong {\Bbb Z}^2$ and $j: \Lambda\to {\Bbb C}$ is an injective 
homomorphism with the discrete image.  A  morphism  of lattices 
$(\Lambda, {\Bbb C}, j)\to (\Lambda', {\Bbb C}, j')$  
is the identity $j\circ\psi=\varphi\circ j'$
where $\varphi$ is a group homomorphism and $\psi$ is a ${\Bbb C}$-linear
map. It is not hard to see, that any isomorphism class of a lattice contains
a representative given by $j: {\Bbb Z}^2\to  {\Bbb C}$ such that $j(1,0)=1,
j(0,1)=\tau\in {\Bbb H}$.   The category of lattices ${\cal L}$ consists of $Ob~({\cal L})$,
which are lattices $(\Lambda, {\Bbb C}, j)$ and morphisms $H(L,L')$
between $L,L'\in Ob~({\cal L})$ which coincide with the morphisms
of lattices specified above. For any  $L,L',L''\in Ob~({\cal L})$
and any morphisms $\varphi': L\to L'$, $\varphi'': L'\to L''$ a 
morphism $\phi: L\to L''$ is the  composite  of $\varphi'$ and
$\varphi''$,  which we write as $\phi=\varphi''\varphi'$. 
The  identity  morphism, $1_L$, is a morphism $H(L,L)$.
\end{dfn}
\begin{rmk}
\textnormal{
The lattices are bijective with the complex tori (and elliptic curves) via the formula
$(\Lambda, {\Bbb C}, j)\mapsto {\Bbb C}/j(\Lambda)$;  
thus   ${\cal L}\cong $ {\bf Ell}. 
}
\end{rmk}
 \index{pseudo-lattice}
\begin{dfn}
{\bf (Category of pseudo-lattices)}
By a pseudo- lattice (of rank 2) in the real line  ${\Bbb R}$ one understands
a triple $(\Lambda, {\Bbb R}, j)$, where
$\Lambda\cong {\Bbb Z}^2$ and $j: \Lambda\to {\Bbb R}$ is a homomorphism.  
A morphism of the pseudo-lattices $(\Lambda, {\Bbb R}, j)\to (\Lambda', {\Bbb R}, j')$
is  the identity $j\circ\psi=\varphi\circ j'$,
where $\varphi$ is a group homomorphism and $\psi$ is an inclusion 
map (i.e. $j'(\Lambda')\subseteq j(\Lambda)$).  
Any isomorphism class of a pseudo-lattice contains
a representative given by $j: {\Bbb Z}^2\to  {\Bbb R}$, such that $j(1,0)=\lambda_1,
j(0,1)=\lambda_2$, where $\lambda_1,\lambda_2$ are the positive reals.
The pseudo-lattices make up a category, which we denote by ${\cal PL}$.
\end{dfn}
\begin{lem}\label{lem5.1.4}
The pseudo-lattices are bijective with the measured foliations on torus
via the formula $(\Lambda, {\Bbb R}, j)\mapsto {\cal F}_{\lambda_2/\lambda_1}^{\lambda_1}$,
where  ${\cal F}_{\lambda_2/\lambda_1}^{\lambda_1}$ is a foliation of the slope $\theta=\lambda_2/\lambda_1$
and measure $\mu=\lambda_1$.  
\end{lem}
{\it Proof}. Define a pairing by the formula
$(\gamma, Re~\omega)\mapsto \int_{\gamma} Re~\omega$, where  $\gamma\in H_1(T^2, {\Bbb Z})$
and $\omega\in H^0(S; \Omega)$. The trajectories of the closed differential $\phi:=Re~\omega$ 
define a measured foliation on $T^2$. Thus, in view of the pairing, the linear spaces $\Phi_{T^2}$
and  $Hom~(H_1(T^2, {\Bbb Z}); {\Bbb R})$ are isomorphic. Notice that 
the latter space coincides with  the  space of the pseudo-lattices.
To obtain an explicit bijection formula, let us evaluate the integral:
\displaymath
\int_{{\Bbb Z}\gamma_1+{\Bbb Z}\gamma_2}\phi ={\Bbb Z}\int_{\gamma_1}\phi + 
{\Bbb Z}\int_{\gamma_2}\phi= {\Bbb Z}\int_0^1\mu dx + {\Bbb Z}\int_0^1\mu dy,
\enddisplaymath
where $\{\gamma_1,\gamma_2\}$ is a basis in $H_1(T^2, {\Bbb Z})$. 
Since ${dy\over dx}=\theta$, one gets:
\displaymath
\left\{
\begin{array}{cccc}
\int_0^1\mu dx  &= \mu  &= \lambda_1  & \nonumber\\
\int_0^1\mu dy  &= \int_0^1\mu\theta dx   &= \mu\theta  &=
\lambda_2. 
\end{array}
\right.
\enddisplaymath
Thus, $\mu=\lambda_1$ and $\theta={\lambda_2\over\lambda_1}$.
Lemma \ref{lem5.1.4} follows.  
$\square$

\begin{rmk}
\textnormal{
It follows from Lemma \ref{lem5.1.4}  and the canonical bijection 
 $h: \Phi_{T^2}\to {\Bbb H}-\{pt\}$,   that ${\cal L}\cong {\cal PL}$
are the equivalent categories.
 }
\end{rmk}
 \index{projective pseudo-lattice}
\begin{dfn}
{\bf (Category of projective pseudo-lattices)}
By a projective pseudo- lattice (of rank 2)  one understands
a triple  $(\Lambda, {\Bbb R}, j)$, where $\Lambda\cong {\Bbb Z}^2$ and $j: \Lambda\to {\Bbb R}$ 
is a  homomorphism. A morphism of the projective pseudo-lattices
$(\Lambda, {\Bbb C}, j)\to (\Lambda', {\Bbb R}, j')$
is  the identity $j\circ\psi=\varphi\circ j'$,
where $\varphi$ is a group homomorphism and $\psi$ is an ${\Bbb R}$-linear
map. (Notice, that unlike the case of the pseudo-lattices, $\psi$ is a scaling
map as opposite to an inclusion map. Thus,  the two pseudo-lattices
can  be projectively equivalent, while being distinct in the category ${\cal PL}$.) 
It is not hard to see that any isomorphism class of a projective pseudo-lattice 
contains a representative given by $j: {\Bbb Z}^2\to  {\Bbb R}$ such that $j(1,0)=1,
j(0,1)=\theta$, where $\theta$ is a positive real.
The projective pseudo-lattices make up a category, which we shall 
denote by ${\cal PPL}$.
\end{dfn}
\begin{lem}\label{lem5.1.5}
${\cal PPL}\cong$  {\bf NC-Tor},  i.e projective pseudo-lattices and noncommutative tori 
 are equivalent categories.
\end{lem}
{\it Proof.}
Notice  that  projective pseudo-lattices are bijective with the noncommutative tori,
via the formula $(\Lambda, {\Bbb R}, j)\mapsto {\cal  A}_{\theta}$.
An isomorphism $\varphi: \Lambda\to\Lambda'$ acts by the formula
$1\mapsto a+b\theta$, $\theta\mapsto c+d\theta$, where $ad-bc=1$
and $a,b,c,d\in {\Bbb Z}$. Therefore,  $\theta'={c+d\theta\over a+b\theta}$.
Thus, isomorphic projective pseudo-lattices map to the stably isomorphic (Morita equivalent)
 noncommutative tori.  Lemma \ref{lem5.1.5} follows.  
$\square$

\bigskip
To define a map $F:$ {\bf Ell} $\to$ {\bf NC-Tor},  we shall
consider a composition of the following morphisms
\displaymath
\hbox{{\bf Ell}}
\buildrel\rm\sim
\over\longrightarrow
{\cal L}
\buildrel\rm\sim
\over\longrightarrow
{\cal PL}
\buildrel\rm F
\over\longrightarrow
{\cal PPL}
\buildrel\rm\sim
\over\longrightarrow
 \hbox{{\bf NC-Tor}},
\enddisplaymath
where all the arrows, but $F$, have been defined. 
To define $F$,  let $PL\in {\cal PL}$  be a pseudo-lattice,
such that $PL=PL(\lambda_1,\lambda_2)$,  where $\lambda_1=j(1,0), \lambda_2=j(0,1)$
are positive reals.  Let $PPL\in {\cal PPL}$  be a projective pseudo-lattice,
such that $PPL=PPL(\theta)$,  where $j(1,0)=1$ and  $j(0,1)=\theta$ is a positive real.
Then $F: {\cal PL}\to {\cal PPL}$ is given by the formula 
$PL(\lambda_1, \lambda_2)\longmapsto PPL\left({\lambda_2\over\lambda_1}\right)$.  
It is easy to see, that $Ker ~F\cong (0,\infty)$ and $F$ is not an injective map. 
Since all the arrows, but $F$,   are the isomorphisms between 
the categories, one gets a map
\displaymath
F: \hbox{{\bf Ell}} \longrightarrow  \hbox{{\bf NC-Tor}}. 
\enddisplaymath
\begin{lem}\label{lem5.1.6}
{\bf (Basic lemma)} 
The map $F:$  {\bf Ell} $\to$  {\bf NC-Tor}   is a  covariant functor
which maps isomorphic complex tori to the stably isomorphic (Morita equivalent) 
noncommutative tori;  the functor is  non-injective functor and  $Ker ~F\cong (0,\infty)$.
\end{lem}
{\it Proof.} (i) Let us show that $F$ maps isomorphic
complex tori to the stably isomorphic noncommutative tori. 
Let  ${\Bbb C}/({\Bbb Z}\omega_1 +{\Bbb Z}\omega_2)$
be a complex torus.  Recall that the periods 
$\omega_1=\int_{\gamma_1}\omega_E$ and $\omega_2=\int_{\gamma_2}\omega_E$,
where $\omega_E=dz$ is an invariant (N\'eron) differential on the complex
torus and $\{\gamma_1,\gamma_2\}$ is a basis in $H_1(T^2, {\Bbb Z})$.
The map $F$ can be written as
\displaymath
{\Bbb C}/L_{(\int_{\gamma_2}\omega_E) / (\int_{\gamma_1}\omega_E)}
\buildrel\rm F\over
\longmapsto
{\cal A}_{(\int_{\gamma_2}\phi)/(\int_{\gamma_1}\phi)},
\enddisplaymath
where $L_{\omega_2/\omega_1}$ is a lattice and  $\phi=Re ~\omega$ is a closed differential defined earlier.
Note that every isomorphism in the category {\bf Ell} is induced by 
an orientation preserving automorphism,  $\varphi$, of the torus $T^2$. 
The action of $\varphi$ on the homology basis $\{\gamma_1,\gamma_2\}$
of $T^2$ is given by the formula:
\displaymath
\left\{
\begin{array}{cc}
\gamma_1' &= a\gamma_1+b\gamma_2\nonumber\\
\gamma_2' &= c\gamma_1+d\gamma_2
\end{array}
\right.
\hbox{, ~~ where}
\quad\left(\matrix{a & b\cr c & d}\small\right)\in SL_2({\Bbb Z}). 
\enddisplaymath
The functor $F$ acts   by the formula:
\displaymath
\tau={\int_{\gamma_2}\omega_E\over\int_{\gamma_1}\omega_E}
\longmapsto
\theta={\int_{\gamma_2}\phi\over\int_{\gamma_1}\phi}.
\enddisplaymath

\smallskip\noindent
(a) From the left-hand side of the above equation,  one obtains 
\displaymath
\left\{
\begin{array}{ccccc}
\omega_1' &= \int_{\gamma_1'}\omega_E &=  \int_{a\gamma_1+b\gamma_2}\omega_E  &= 
a\int_{\gamma_1}\omega_E +b\int_{\gamma_2}\omega_E &= a\omega_1+b\omega_2\nonumber\\
\omega_2' &= \int_{\gamma_2'}\omega_E &=  \int_{c\gamma_1+d\gamma_2}\omega_E  &= 
c\int_{\gamma_1}\omega_E +d\int_{\gamma_2}\omega_E &= c\omega_1+d\omega_2,
\end{array}
\right.
\enddisplaymath
and therefore $\tau'={\int_{\gamma_2'}\omega_E\over\int_{\gamma_1'}\omega_E}=
{c+d\tau\over a+b\tau}$.

\smallskip\noindent
(b) From the right-hand side, one obtains
\displaymath
\left\{
\begin{array}{ccccc}
\lambda_1' &= \int_{\gamma_1'}\phi &=  \int_{a\gamma_1+b\gamma_2}\phi  &= 
a\int_{\gamma_1}\phi +b\int_{\gamma_2}\phi &= a\lambda_1+b\lambda_2\nonumber\\
\lambda_2' &= \int_{\gamma_2'}\phi &=  \int_{c\gamma_1+d\gamma_2}\phi  &= 
c\int_{\gamma_1}\phi +d\int_{\gamma_2}\phi &= c\lambda_1+d\lambda_2,
\end{array}
\right.
\enddisplaymath
and therefore $\theta'={\int_{\gamma_2'}\phi\over\int_{\gamma_1'}\phi}=
{c+d\theta\over a+b\theta}$.  Comparing (a) and (b),  one concludes
that $F$ maps  isomorphic complex tori to the stably isomorphic 
(Morita equivalent)  noncommutative tori.

\medskip
(ii) Let us show that $F$ is a covariant functor, i.e. $F$ does not reverse 
the arrows.  Indeed, it can be verified directly using the above  formulas, that 
$F(\varphi_1\varphi_2)=\varphi_1\varphi_2=F(\varphi_1)F(\varphi_2)$
for any pair of the isomorphisms $\varphi_1,\varphi_2\in Aut~(T^2)$.

\medskip
(iii) Since  $F: {\cal PL}\to {\cal PPL}$ is given by the formula 
$PL(\lambda_1, \lambda_2)\longmapsto PPL\left({\lambda_2\over\lambda_1}\right)$, 
one gets $Ker ~F\cong (0,\infty)$ and $F$ is not an injective map. 
Lemma  \ref{lem5.1.6} is proved.
$\square$

\bigskip
Theorem \ref{thm5.1.2} follows from Lemma \ref{lem5.1.6}.
$\square$

\vskip1cm\noindent
{\bf Guide to the literature.}
 The  basics of  elliptic curves are covered by   [Husem\"oller 1986]  \cite{H2}, 
   [Knapp 1992] \cite{K1},  [Koblitz 1984]  \cite{K2},   [Silverman 1985]  \cite{S1},  
   [Silverman 1994]   \cite{S2},  [Silverman \& Tate  1992]  \cite{ST}  and others.  
More advanced  topics are discussed  in the survey papers [Cassels 1966]  \cite{Cas1},
[Mazur  1986]  \cite{Maz1} and [Tate 1974]  \cite{Tat1}.  
The Sklyanin algebras were introduced  and studied in [Sklyanin  1982]  \cite{Skl1} 
and [Sklyanin  1983]  \cite{Skl2};  for a detailed account,    see   [Feigin \& Odesskii  1989]  \cite{FeOd1}
and   [Feigin \& Odesskii  1993]  \cite{FeOd2}.  The general theory is covered by 
 [Stafford \& van ~den ~Bergh  2001]  \cite{StaVdb1}.   The basics of measured foliations
 and the Teichm\"uller theory can be found in    [Thurston  1988]  \cite{Thu1} and 
[Hubbard \& Masur  1979]  \cite{HuMa1}.  The functor from elliptic curves to
noncommutative tori was constructed in  \cite{Nik2} and \cite{Nik4} using  measured foliations   and 
in \cite{Nik10}  using  Sklyanin's  algebras.   
The  idea  of   infinite-dimensional  representations of  Sklyanin's algebras
by the linear operators on a Hilbert space  ${\cal H}$   belongs to    [Sklyanin 1982]  
 \cite{Skl1},  the end of Section 3.

 \index{complex algebraic curve}

\section{Higher genus algebraic curves}
By a  {\it complex algebraic  curve}  one understands  a subset of the complex projective 
plane  of the form
\displaymath
C=\{(x,y,z)\in {\Bbb C}P^2 ~|~ P(x,y,z)=0\},
\enddisplaymath
where $P(x,y,z)$ is a homogeneous polynomial with complex coefficients;
such curves are isomorphic to the  complex 2-dimensional manifolds,  i.e.
the {\it Riemann surfaces} $S$.   
We shall  construct a functor $F$ on a  generic set of   complex algebraic curves  with values 
in the category of {\it toric} AF-algebras;   the functor maps  isomorphic algebraic curves 
to the stably  isomorphic (Morita equivalent)  toric AF-algebras. 
For genus $g=1$ algebraic (i.e. elliptic) curves,  the toric AF-algebras are 
isomorphic to  the Effros-Shen algebras ${\Bbb A}_{\theta}$;   such  AF-algebras
are known to contain the noncommutative torus ${\cal A}_{\theta}$,  see Theorem
\ref{thm3.5.3}.    The functor $F$ will be used to prove Harvey's conjecture by 
construction  of a faithful  representation of the mapping class group of genus $g\ge 2$  
in the matrix  group $GL(6g-6, {\Bbb Z})$,   see Section 5.4.

 \index{toric AF-algebra}

\subsection{Toric AF-algebras}
We repeat some facts of the Teichm\"uller theory,  see e.g. [Hubbard \& Masur  1979]  \cite{HuMa1}.  
Denote by $T_S(g)$ the Teichm\"uller  space of genus $g\ge 1$
(i.e. the space of all complex 2-dimensional manifolds of genus $g$)  
 endowed with a distinguished point $S$. 
Let   $q\in H^0(S, \Omega^{\otimes 2})$ be a holomorphic quadratic 
differential on the Riemann surface $S$, such that all zeroes
of $q$ (if any) are simple. By $\widetilde S$ we understand a double 
cover of $S$ ramified over the zeroes of $q$ and by
$H_1^{odd}(\widetilde S)$ the odd part of the integral homology of $\widetilde S$ relatively the  zeroes.
Note that $H_1^{odd}(\widetilde S)\cong {\Bbb Z}^{n}$, where $n=6g-6$ if $g\ge 2$ and $n=2$ if $g=1$. 
It is known that 
\displaymath
T_S(g)\cong Hom~(H_1^{odd}(\widetilde S); {\Bbb R})-\{0\},
\enddisplaymath
where $0$ is the zero homomorphism [Hubbard \& Masur  1979]  \cite{HuMa1}.  
Denote by $\lambda=(\lambda_1,\dots,\lambda_{n})$ the image of a basis of 
$H_1^{odd}(\widetilde S)$ in the real line ${\Bbb R}$, such that $\lambda_1\ne0$. 
\begin{rmk}
\textnormal{
The claim $\lambda_1\ne0$ is not restrictive, 
because  the zero homomorphism is excluded.
 }
\end{rmk}
We let $\theta=(\theta_1,\dots,\theta_{n-1})$, where $\theta_i=\lambda_{i-1}/\lambda_1$. 
Recall that,  up to a scalar multiple, vector $(1,\theta)\in {\Bbb R}^{n}$ is the limit
of a generically convergent Jacobi-Perron continued fraction [Bernstein 1971]  \cite{BE}
\displaymath
\left(\matrix{1\cr \theta}\right)=
\lim_{k\to\infty} \left(\matrix{0 & 1\cr I & b_1}\right)\dots
\left(\matrix{0 & 1\cr I & b_k}\right)
\left(\matrix{0\cr {\Bbb I}}\right),
\enddisplaymath
where $b_i=(b^{(i)}_1,\dots, b^{(i)}_{n-1})^T$ is a vector of the non-negative integers,  
$I$ the unit matrix and ${\Bbb I}=(0,\dots, 0, 1)^T$.
\begin{dfn}\label{dfn5.2.1} 
By a toric AF-algebra ${\Bbb A}_{\theta}$ one understands the AF-algebra given by  
the Bratteli diagram in Fig. 5.5,   where numbers $b_j^{(i)}$  indicate the multiplicity  
of edges of   the graph. 
\end{dfn}
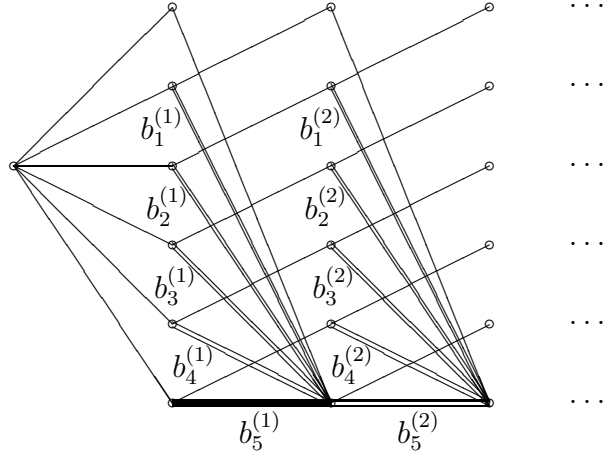
\begin{figure}[here]
\begin{picture}(300,200)(-50,0)

\put(40,110){\circle{3}}

\put(40,110){\line(1,1){60}}
\put(40,110){\line(2,1){60}}
\put(40,110){\line(1,0){60}}
\put(40,110){\line(2,-1){60}}
\put(40,110){\line(1,-1){60}}
\put(40,110){\line(2,-3){60}}

\put(100,20){\circle{3}}
\put(100,50){\circle{3}}
\put(100,80){\circle{3}}
\put(100,110){\circle{3}}
\put(100,140){\circle{3}}
\put(100,170){\circle{3}}

\put(160,20){\circle{3}}
\put(160,50){\circle{3}}
\put(160,80){\circle{3}}
\put(160,110){\circle{3}}
\put(160,140){\circle{3}}
\put(160,170){\circle{3}}

\put(220,20){\circle{3}}
\put(220,50){\circle{3}}
\put(220,80){\circle{3}}
\put(220,110){\circle{3}}
\put(220,140){\circle{3}}
\put(220,170){\circle{3}}

\put(160,20){\line(2,1){60}}
\put(160,19){\line(1,0){60}}
\put(160,21){\line(1,0){60}}
\put(160,50){\line(2,1){60}}
\put(160,49){\line(2,-1){60}}
\put(160,51){\line(2,-1){60}}
\put(160,80){\line(2,1){60}}
\put(160,79){\line(1,-1){60}}
\put(160,81){\line(1,-1){60}}
\put(160,110){\line(2,1){60}}
\put(160,109){\line(2,-3){60}}
\put(160,111){\line(2,-3){60}}
\put(160,140){\line(2,1){60}}
\put(160,139){\line(1,-2){60}}
\put(160,141){\line(1,-2){60}}
\put(160,170){\line(2,-5){60}}


\put(100,20){\line(2,1){60}}
\put(100,19){\line(1,0){60}}
\put(100,20){\line(1,0){60}}
\put(100,21){\line(1,0){60}}
\put(100,50){\line(2,1){60}}
\put(100,49){\line(2,-1){60}}
\put(100,51){\line(2,-1){60}}
\put(100,80){\line(2,1){60}}
\put(100,79){\line(1,-1){60}}
\put(100,81){\line(1,-1){60}}
\put(100,110){\line(2,1){60}}
\put(100,109){\line(2,-3){60}}
\put(100,111){\line(2,-3){60}}
\put(100,140){\line(2,1){60}}
\put(100,139){\line(1,-2){60}}
\put(100,141){\line(1,-2){60}}
\put(100,170){\line(2,-5){60}}


\put(250,20){$\dots$}
\put(250,50){$\dots$}
\put(250,80){$\dots$}
\put(250,110){$\dots$}
\put(250,140){$\dots$}
\put(250,170){$\dots$}


\put(125,5){$b_5^{(1)}$}
\put(100,30){$b_4^{(1)}$}
\put(93,60){$b_3^{(1)}$}
\put(90,90){$b_2^{(1)}$}
\put(88,120){$b_1^{(1)}$}

\put(185,5){$b_5^{(2)}$}
\put(160,30){$b_4^{(2)}$}
\put(153,60){$b_3^{(2)}$}
\put(150,90){$b_2^{(2)}$}
\put(148,120){$b_1^{(2)}$}


\end{picture}

\caption{Toric AF-algebra ${\Bbb A}_{\theta}$ (case  $g=2$).}
\end{figure}

\begin{rmk}
\textnormal{
Note that in the case $g=1$, the Jacobi-Perron 
fraction coincides with the  regular continued fraction and ${\Bbb A}_{\theta}$
becomes the Effros-Shen algebra,  see Example  \ref{exm3.5.2}.  
}
\end{rmk}
\begin{dfn}
 By {\bf Alg-Gen}  we shall mean  the maximal subset of $T_S(g)$ such that 
 for each complex algebraic  curve $C\in$ {\bf Alg-Gen}     the corresponding 
  Jacobi-Perron continued fraction is convergent;    the arrows of {\bf Alg-Gen}  
  are   isomorphisms between  complex algebraic  curves $C$.   We shall write 
  {\bf AF-Toric} to denote the  category of all toric AF-algebras  ${\Bbb  A}_{\theta}$;
 the arrows of {\bf AF-Toric}  are  stable isomorphisms (Morita equivalences)  
 between  toric AF-algebras  ${\Bbb A}_{\theta}$.  By $F$ we understand 
 a map  given by the formula  $C\mapsto {\Bbb A}_{\theta}$;
 in other words,  we have a map
 \displaymath
 F: \hbox{{\bf Alg-Gen}}\longrightarrow \hbox{{\bf AF-Toric}}. 
 \enddisplaymath
 \end{dfn}
\begin{thm}\label{thm5.2.1}
{\bf (Functor on algebraic curves)} 
The set {\bf Alg-Gen}  is a generic subset of $T_S(g)$ and the  map 
$F$ has the following properties:

\smallskip
(i) {\bf Alg-Gen}  $\cong$  {\bf AF-Toric} $\times ~(0,\infty)$ is a trivial fiber bundle, whose
projection  map $p:$ {\bf Alg-Gen}  $\to$  {\bf AF-Toric}  coincides with $F$;

\smallskip
(ii) $F:$ {\bf Alg-Gen}  $\to$  {\bf AF-Toric}  is a covariant functor, which maps isomorphic complex algebraic curves 
$C,C'\in$  {\bf Alg-Gen}  to the stably isomorphic (Morita equivalent)  toric AF-algebras 
${\Bbb A}_{\theta},{\Bbb A}_{\theta'}\in$ {\bf AF-Toric}.
\end{thm}

\subsection{Proof of Theorem \ref{thm5.2.1}}
We shall repeat some known facts and notation.  
Let $S$ be a Riemann surface, and $q\in H^0(S,\Omega^{\otimes 2})$ a holomorphic
quadratic differential on $S$. The lines $Re~q=0$ and $Im~q=0$ define a pair
of measured foliations on $R$, which are transversal to each other outside the set of 
singular points. The set of singular points is common to both foliations and coincides
with the zeroes of $q$. The above measured foliations are said to represent the  vertical and horizontal 
 trajectory structure  of $q$, respectively. 
Denote by  $T(g)$  the Teichm\"uller space of the topological surface $X$ of genus $g\ge 1$,
i.e. the space of the complex structures on $X$. 
Consider the vector bundle $p: Q\to T(g)$ over $T(g)$ whose fiber above a point 
$S\in T_g$ is the vector space $H^0(S,\Omega^{\otimes 2})$.   
Given non-zero $q\in Q$ above $S$, we can consider horizontal measured foliation
${\cal F}_q\in \Phi_X$ of $q$, where $\Phi_X$ denotes the space of equivalence
classes of measured foliations on $X$. If $\{0\}$ is the zero section of $Q$,
the above construction defines a map $Q-\{0\}\longrightarrow \Phi_X$. 
For any ${\cal F}\in\Phi_X$, let $E_{\cal F}\subset Q-\{0\}$ be the fiber
above ${\cal F}$. In other words, $E_{\cal F}$ is a subspace of the holomorphic 
quadratic forms whose horizontal trajectory structure coincides with the 
measured foliation ${\cal F}$. 
If ${\cal F}$ is a measured foliation with the simple zeroes (a generic case),  
then $E_{\cal F}\cong {\Bbb R}^n - 0$, while $T(g)\cong {\Bbb R}^n$, where $n=6g-6$ if
$g\ge 2$ and $n=2$ if $g=1$.   The restriction of $p$ to $E_{\cal F}$ defines a homeomorphism 
(an embedding)
\displaymath
h_{\cal F}: E_{\cal F}\to T(g).
\enddisplaymath
The above result implies that the measured foliations  parametrize  
the space $T(g)-\{pt\}$,  where $pt= h_{\cal F}(0)$.
Indeed, denote by  ${\cal F}'$ a vertical trajectory structure of  $q$. Since ${\cal F}$
and ${\cal F}'$ define $q$, and ${\cal F}=Const$ for all $q\in E_{\cal F}$, one gets a homeomorphism 
between $T(g)-\{pt\}$ and  $\Phi_X$, where $\Phi_X\cong {\Bbb R}^n-0$ is the space of 
equivalence classes of the measured foliations ${\cal F}'$ on $X$. 
Note that the above parametrization depends on a foliation ${\cal F}$.
However, there exists a unique canonical homeomorphism $h=h_{\cal F}$
as follows. Let $Sp ~(S)$ be the length spectrum of the Riemann surface
$S$ and $Sp ~({\cal F}')$ be the set positive reals $\inf \mu(\gamma_i)$,
where $\gamma_i$ runs over all simple closed curves, which are transverse 
to the foliation ${\cal F}'$. A canonical homeomorphism 
$h=h_{\cal F}: \Phi_X\to T(g)-\{pt\}$ is defined by the formula
$Sp ~({\cal F}')= Sp ~(h_{\cal F}({\cal F}'))$ for $\forall {\cal F}'\in\Phi_X$. 
Thus, there exists a canonical  homeomorphism 
\displaymath
h:\Phi_X\to T(g)-\{pt\}.
\enddisplaymath
Recall that $\Phi_X$ is the space of equivalence classes of measured
foliations on the topological surface $X$. Following [Douady \& Hubbard  1975] 
\cite{DoHu1},  we consider the following  coordinate system on $\Phi_X$. 
For clarity, let us make a generic assumption that $q\in H^0(S,\Omega^{\otimes 2})$
is a non-trivial holomorphic quadratic differential with only simple zeroes. 
We wish to construct a Riemann surface of $\sqrt{q}$, which is a double cover
of $S$ with ramification over the zeroes of $q$. Such a surface, denoted by
$\widetilde S$, is unique and has an advantage of carrying a holomorphic
differential $\omega$, such that $\omega^2=q$. We further denote by 
$\pi:\widetilde S\to S$ the covering projection. The vector space
$H^0(\widetilde S,\Omega)$ splits into the direct sum
$H^0_{even}(\widetilde S,\Omega)\oplus H^0_{odd}(\widetilde S,\Omega)$
in view of  the involution $\pi^{-1}$ of $\widetilde S$, and
the vector space $H^0(S,\Omega^{\otimes 2})\cong H^0_{odd}(\widetilde S,\Omega)$.
Let $H_1^{odd}(\widetilde S)$ be the odd part of the homology of $\widetilde S$
relatively  the zeroes of $q$.   Consider the pairing
$H_1^{odd}(\widetilde S)\times H^0(S, \Omega^{\otimes 2})\to {\Bbb C}$,
defined by the integration  $(\gamma, q)\mapsto \int_{\gamma}\omega$. 
We shall take the associated map
$\psi_q: H^0(S,\Omega^{\otimes 2})\to Hom~(H_1^{odd}(\widetilde S); {\Bbb C})$
and let $h_q= Re~\psi_q$. 
\begin{lem}\label{lem5.2.1}
{\bf ([Douady \& Hubbard  1975]  \cite{DoHu1})}
The map
\displaymath
h_q: H^0(S, \Omega^{\otimes 2})\longrightarrow Hom~(H_1^{odd}(\widetilde S); {\Bbb R})
\enddisplaymath
is an ${\Bbb R}$-isomorphism. 
\end{lem}
\begin{rmk}
\textnormal{
Since  each  ${\cal F}\in \Phi_X$ is the  vertical foliation 
$Re~q=0$ for a $q\in H^0(S, \Omega^{\otimes 2})$,  Lemma \ref{lem5.2.1} 
implies that $\Phi_X\cong Hom~(H_1^{odd}(\widetilde S); {\Bbb R})$.
By  formulas for the relative homology, 
one finds that $H_1^{odd}(\widetilde S)\cong {\Bbb Z}^{n}$,
where $n=6g-6$ if $g\ge 2$ and $n=2$ if $g=1$.     
Each $h\in Hom~({\Bbb Z}^{n}; {\Bbb R})$ is given
by the reals  $\lambda_1=h(e_1),\dots, \lambda_{n}=h(e_{n})$,
where $(e_1,\dots, e_{n})$ is a basis in ${\Bbb Z}^{n}$.
The numbers   $(\lambda_1,\dots,\lambda_{n})$ are the coordinates in the space $\Phi_X$
and, therefore,  in  the Teichm\"uller space $T(g)$. 
}
\end{rmk}
To prove Theorem \ref{thm5.2.1},   we shall consider the following categories:
(i) generic complex algebraic curves {\bf Alg-Gen};
(ii) pseudo-lattices ${\cal PL}$; (iii) projective pseudo-lattices
${\cal PPL}$ and (iv)  category {\bf AF-Toric}  of the toric AF-algebras.
First, we show that {\bf Alg-Gen} $\cong {\cal PL}$ are equivalent categories,  such that isomorphic
complex algebraic  curves $C,C'\in$ {\bf Alg-Gen}  map to isomorphic pseudo-lattices
$PL, PL'\in {\cal PL}$.  Next, a non-injective functor $F: {\cal PL}\to {\cal PPL}$
is constructed.  The $F$ maps isomorphic pseudo-lattices to isomorphic
projective pseudo-lattices and $Ker ~F\cong (0,\infty)$. 
Finally, it is shown that a subcategory $U\subseteq {\cal PPL}$ and {\bf AF-Toric}  
are   equivalent categories. In other words, we have the following diagram
\displaymath
\hbox{{\bf Alg-Gen}}
\buildrel\rm\alpha
\over\longrightarrow
{\cal PL}
\buildrel\rm F
\over\longrightarrow
U
\buildrel\rm \beta
\over\longrightarrow
\hbox{{\bf AF-Toric}},
\enddisplaymath
where $\alpha$ is an injective map,  $\beta$ is a bijection  and $Ker ~F\cong (0,\infty)$. 
\begin{dfn}
Let $Mod ~X$ be the mapping class group of the surface $X$. 
A {\it complex algebraic curve}  is a triple $(X, C, j)$, where
$X$ is a topological surface of genus $g\ge 1$,  $j: X\to C$ is a
complex (conformal) parametrization of $X$ and $C$ is a Riemann surface.
A {\it morphism} of complex algebraic curves  
$(X, C, j)\to (X, C', j')$   is the identity
$j\circ\psi=\varphi\circ j'$, 
where $\varphi\in Mod~X$ is a diffeomorphism of $X$  and $\psi$ is an isomorphism of Riemann surfaces. 
A  category of generic complex algebraic curves, $\hbox{{\bf Alg-Gen}}$, consists of $Ob~(\hbox{{\bf Alg-Gen}})$
which are complex algebraic curves  $C\in  T_S(g)$ and morphisms $H(C,C')$
between $C,C'\in Ob~(\hbox{{\bf Alg-Gen}})$ which coincide with the morphisms
specified above. For any  $C,C',C''\in Ob~(\hbox{{\bf Alg-Gen}})$
and any morphisms $\varphi': C\to C'$, $\varphi'': C'\to C''$ a 
morphism $\phi: C\to C''$ is the composite of $\varphi'$ and
$\varphi''$, which we write as $\phi=\varphi''\varphi'$. 
The identity  morphism, $1_C$, is a morphism $H(C,C)$.
\end{dfn}
\begin{dfn}
By a  pseudo-lattice  (of rank $n$)  one undestands the triple $(\Lambda, {\Bbb R}, j)$, where
$\Lambda\cong {\Bbb Z}^n$ and $j: \Lambda\to {\Bbb R}$ is a homomorphism.  
A morphism of pseudo-lattices $(\Lambda, {\Bbb R}, j)\to (\Lambda, {\Bbb R}, j')$
is the identity  $j\circ\psi=\varphi\circ j'$,
where $\varphi$ is a group homomorphism and $\psi$ is an inclusion 
map, i.e. $j'(\Lambda')\subseteq j(\Lambda)$.  
Any isomorphism class of a pseudo-lattice contains
a representative given by $j: {\Bbb Z}^n\to  {\Bbb R}$ such that 
$j(1,0,\dots, 0)=\lambda_1, \quad j(0,1,\dots,0)=\lambda_2,\quad \dots, \quad  j(0,0,\dots,1)=\lambda_n,$
where $\lambda_1,\lambda_2,\dots,\lambda_n$ are positive reals.
The pseudo-lattices of rank $n$  make up a category, which we denote by ${\cal PL}_n$. 
\end{dfn}
\begin{lem}\label{lem5.2.2}
{\bf (Basic lemma)}
Let $g\ge 2$ ($g=1$) and $n=6g-6$ ($n=2$). 
There exists an injective covariant functor $\alpha: \hbox{{\bf Alg-Gen}}\to {\cal PL}_{n}$,
which maps isomorphic complex algebraic curves $C,C'\in \hbox{{\bf Alg-Gen}}$ to the isomorphic
pseudo-lattices $PL,PL'\in {\cal PL}_{n}$.
\end{lem}
{\it Proof.} Recall that we have a map $\alpha: T(g)-\{pt\}\to Hom~(H_1^{odd}(\widetilde S); {\Bbb R})-\{0\}$;
  $\alpha$ is a homeomorphism and,  therefore,  $\alpha$ is injective. 
Let us find the image $\alpha(\varphi)\in Mor~({\cal PL})$ of 
$\varphi\in Mor~(\hbox{{\bf Alg-Gen}})$. Let $\varphi\in Mod~X$ be a diffeomorphism
of $X$, and let $\widetilde X\to X$ be the ramified double cover of $X$.  
We denote by $\widetilde\varphi$ the induced
mapping on $\widetilde X$. Note that $\widetilde\varphi$ is a diffeomorphism
of $\widetilde X$ modulo the covering involution ${\Bbb Z}_2$. Denote by
$\widetilde\varphi^*$ the action of $\widetilde\varphi$ on 
$H_1^{odd}(\widetilde X)\cong {\Bbb Z}^{n}$. Since $\widetilde\varphi~mod~{\Bbb Z}_2$ is
a diffeomorphism of $\widetilde X$, $\widetilde\varphi^*\in GL_{n}({\Bbb Z})$. 
Thus, $\alpha(\varphi)=\widetilde\varphi^*\in Mor~({\cal PL})$. 
Let us show that $\alpha$ is a functor.  Indeed, let $C,C'\in \hbox{{\bf Alg-Gen}}$ be  
isomorphic  complex algebraic curves,  such that $C'=\varphi(C)$ for a $\varphi\in Mod~X$. 
Let $a_{ij}$ be the elements of matrix $\widetilde\varphi^*\in GL_{n}({\Bbb Z})$. 
Recall that $\lambda_i=\int_{\gamma_i}\phi$ for a closed 1-form $\phi= Re~\omega$
and $\gamma_i\in H_1^{odd}(\widetilde X)$. Then
$\gamma_j=\sum_{i=1}^{n} a_{ij}\gamma_i, \quad j=1,\dots, n$
are the elements of a new basis in $H_1^{odd}(\widetilde X)$. 
By  integration rules we have 
\displaymath
\lambda_j'= \int_{\gamma_j}\phi=
\int_{\sum a_{ij}\gamma_i}\phi=
\sum_{i=1}^{n}a_{ij}\lambda_i.
\enddisplaymath
Let $j(\Lambda)={\Bbb Z}\lambda_1+\dots+{\Bbb Z}\lambda_{n}$
and $j'(\Lambda)={\Bbb Z}\lambda_1'+\dots+{\Bbb Z}\lambda_{n}'$.
Since $\lambda_j'= \sum_{i=1}^{n}a_{ij}\lambda_i$ and $(a_{ij})\in GL_{n}({\Bbb Z})$,
we conclude that $j(\Lambda)=j'(\Lambda)\subset {\Bbb R}$. In other words, the pseudo-lattices
$(\Lambda, {\Bbb R}, j)$ and $(\Lambda, {\Bbb R}, j')$ are isomorphic. Hence, 
$\alpha: \hbox{{\bf Alg-Gen}}\to {\cal PL}$  maps isomorphic complex algebraic curves to
the isomorphic pseudo-lattices, i.e. $\alpha$ is a functor.  
Let us show that $\alpha$ is a covariant functor. 
Indeed, let $\varphi_1,\varphi_2\in Mor (\hbox{{\bf Alg-Gen}})$. Then  
$\alpha(\varphi_1\varphi_2)= (\widetilde{\varphi_1\varphi_2})^*=\widetilde\varphi_1^*\widetilde\varphi_2^*=
\alpha(\varphi_1)\alpha(\varphi_2)$. Lemma \ref{lem5.2.2} follows.
$\square$

\begin{dfn}
By a  projective pseudo-lattice  (of rank $n$) one understands  a triple 
$(\Lambda, {\Bbb R}, j)$,  where $\Lambda\cong {\Bbb Z}^n$ and $j: \Lambda\to {\Bbb R}$ is a
homomorphism. A morphism of projective pseudo-lattices
$(\Lambda, {\Bbb C}, j)\to (\Lambda, {\Bbb R}, j')$
is the identity  $j\circ\psi=\varphi\circ j'$,
where $\varphi$ is a group homomorphism and $\psi$ is an  ${\Bbb R}$-linear
map.   It is not hard to see that any isomorphism class of a projective pseudo-lattice 
contains a representative given by $j: {\Bbb Z}^n\to  {\Bbb R}$ such that
$j(1,0,\dots,0)=1,\quad
j(0,1,\dots,0)=\theta_1,\quad  \dots, \quad  j(0,0,\dots,1)=\theta_{n-1},$
  where $\theta_i$ are  positive reals.
The projective pseudo-lattices of rank $n$  make up a category, which we denote by ${\cal PPL}_n$. 
\end{dfn}
\begin{rmk}
\textnormal{
Notice that unlike the case of pseudo-lattices,   $\psi$ is a scaling
map as opposite to an inclusion map.   This allows to the two pseudo-lattices
to be projectively equivalent, while being distinct in the category ${\cal PL}_n$. 
}
\end{rmk}
\begin{dfn}
Finally,  the toric AF-algebras ${\Bbb A}_{\theta}$,  modulo stable isomorphism
(Morita equivalences),  make up a category  which we shall denote by {\bf AF-Toric}. 
\end{dfn}
\begin{lem}\label{lem5.2.3}
Let $U_n\subseteq {\cal PPL}_n$ be a subcategory consisting of the projective
pseudo-lattices $PPL=PPL(1,\theta_1,\dots,\theta_{n-1})$ for which the Jacobi-Perron
fraction of the vector $(1,\theta_1,\dots,\theta_{n-1})$ converges to the vector.
Define a map  $\beta: U_n\to$  {\bf AF-Toric}   by the formula
$PPL(1,\theta_1,\dots,\theta_{n-1})\mapsto {\Bbb A}_{\theta},$
where $\theta=(\theta_1,\dots,\theta_{n-1})$. 
Then $\beta$ is a bijective functor, which maps isomorphic projective pseudo-lattices 
to the stably isomorphic  toric AF-algebras. 
 \end{lem}
{\it Proof.} It is evident that $\beta$ is injective and surjective. Let
us show that $\beta$ is a functor. Indeed,  according to
[Effros 1981]   \cite{E}, Corollary 4.7,   every totally 
ordered abelian group of rank $n$ has form ${\Bbb Z}+\theta_1 {\Bbb Z}+\dots+ {\Bbb Z}\theta_{n-1}$.
The latter   is a projective pseudo-lattice $PPL$  from the category $U_n$. 
On the other hand,  by  Theorem \ref{thm3.5.2}  the $PPL$ defines a stable isomorphism class of 
the toric AF-algebra ${\Bbb A}_{\theta}\in$ {\bf AF-Toric}. 
Therefore, $\beta$  maps isomorphic projective pseudo-lattices (from the set $U_n$) to the stably isomorphic toric 
AF-algebras,  and {\it vice versa}.   Lemma \ref{lem5.2.3}  follows.  
 $\square$

\begin{lem}\label{lem5.2.4}
Let $F:  {\cal PL}_n\to {\cal PPL}_n$ be a map given by formula
\displaymath
PL(\lambda_1,\lambda_2,\dots, \lambda_n)\mapsto
 PPL\left(1, {\lambda_2\over\lambda_1},\dots, {\lambda_n\over\lambda_1}\right),
\enddisplaymath
where $PL(\lambda_1,\lambda_2,\dots, \lambda_n)\in {\cal PL}_n$ and 
$PPL(1,\theta_1,\dots,\theta_{n-1})\in {\cal PPL}_n$. 
Then $Ker~F=(0,\infty)$ and  $F$ is a functor which maps  isomorphic pseudo-lattices to 
isomorphic projective  pseudo-lattices. 
\end{lem}
{\it Proof.} Indeed, $F$ can be thought as   a map from ${\Bbb R}^n$ to ${\Bbb R}P^n$. Hence 
$Ker~F= \{\lambda_1 : \lambda_1>0\}\cong (0,\infty)$.  The second part of lemma is evident.  
Lemma \ref{lem5.2.4} is proved. 
$\square$

\bigskip
Theorem \ref{thm5.2.1}  follows from Lemmas \ref{lem5.2.2} - \ref{lem5.2.4} with 
$n=6g-6$ ($n=2$) for  $g\ge 2$ ($g=1$).
$\square$

\vskip1cm\noindent
{\bf Guide to the literature.}
An excellent introduction to complex algebraic curves is the book by  [Kirwan 1992]  \cite{KI}.
For measured foliations and their relation to the Teichm\"uller theory the reader is referred 
to [Hubbard \& Masur  1979]  \cite{HuMa1}.   Functor $F:$  {\bf Alg-Gen} $\to$ {\bf AF-Toric}
was constructed in \cite{Nik3};  the term {\it toric AF-algebras} was coined by Yu.~Manin
(private communication).

 \index{complex projective variety}

\section{Complex projective varieties}
We shall generalize functors constructed in Sections 5.1 and 5.2 to arbitrary 
complex projective varieties $X$.   Namely,     for the category {\bf Proj-Alg}
of all such varieties (of fixed dimension $n$)  we construct a covariant functor
\displaymath
F:  \hbox{{\bf Proj-Alg}}\longrightarrow \hbox{{\bf C*-Serre}},
\enddisplaymath
where {\bf C*-Serre} is a category of the {\it Serre $C^*$-algebras},  ${\cal A}_X$,
attached to variety $X$.   In particular,  if $X\cong {\cal E}_{\tau}$ is an elliptic curve,
then ${\cal A}_X\cong {\cal A}_{\theta}$ is a noncommutative torus and if 
$X\cong C$ is a complex algebraic curve, then ${\cal A}_X\cong {\Bbb A}_{\theta}$ 
is a toric AF-algebra.   For $n\ge 2$ the description of ${\cal A}_X$ in terms of its 
semigroup $K_0^+({\cal A}_X)$ is less satisfactory (so far?)  but using the
Takai duality for  crossed  product $C^*$-algebras,   it is possible to prove the following general result.    
If $B$ is the commutative coordinate ring of variety $X$,   then it is well known that $X\cong$ {\bf Spec} $(B)$,
where {\bf Spec}  is the space of prime ideals of  $B$;   an analog  of this important 
formula for ${\cal A}_X$ is proved to be  $X\cong  \Irred ~({\cal A}_X\rtimes_{\hat\alpha} \hat {\Bbb Z})$,
where $\hat\alpha$ is an automorphism of ${\cal A}_X$ and $\Irred$ the space of all
irreducible representations of the crossed product $C^*$-algebra. 
We illustrate the formula in an important special case ${\cal A}_X\cong {\cal A}_{RM}$,
i.e the case of  noncommutative torus with real multiplication.

 \index{Serre $C^*$-algebra}

\subsection{Serre $C^*$-algebras}
Let $X$ be a projective scheme over a field $k$, and let ${\cal L}$ 
be the invertible sheaf ${\cal O}_X(1)$ of linear forms on $X$.  Recall
that the homogeneous coordinate ring of $X$ is a graded $k$-algebra, 
which is isomorphic to the algebra
\displaymath
B(X, {\cal L})=\bigoplus_{n\ge 0} H^0(X, ~{\cal L}^{\otimes n}). 
\enddisplaymath
Denote by $\Coh$ the category of quasi-coherent sheaves on a scheme $X$
and by $\Mod$ the category of graded left modules over a graded ring $B$.  
If $M=\oplus M_n$ and $M_n=0$ for $n>>0$, then the graded module
$M$ is called {\it right bounded}.  The  direct limit  $M=\lim M_{\alpha}$
is called a {\it torsion}, if each $M_{\alpha}$ is a right bounded graded
module. Denote by $\Tors$ the full subcategory of $\Mod$ of the torsion
modules.  The following result is the fundamental fact  about the graded 
ring $B=B(X, {\cal L})$.   
\begin{thm}\label{thm5.3.1}
{\bf ([Serre 1955]  \cite{Ser1})}
\quad $$\Mod~(B(X, {\cal L})) ~/~\Tors \cong \Coh~(X).$$
\end{thm}
\begin{dfn}\label{dfn5.3.1} 
Let $\alpha$ be an automorphism of a projective scheme $X$;   
the pullback of sheaf ${\cal L}$  along $\alpha$ will be denoted by ${\cal L}^{\alpha}$,  
i.e.   ${\cal L}^{\alpha}(U):= {\cal L}(\alpha U)$ for every $U\subset X$. 
We shall set
\displaymath
B(X, {\cal L}, \alpha)=\bigoplus_{n\ge 0} H^0(X, ~{\cal L}\otimes {\cal L}^{\alpha}\otimes\dots
\otimes  {\cal L}^{\alpha^{ n}}). 
\enddisplaymath
 The multiplication of sections is defined by the rule
 \displaymath
 ab=a\otimes b^{\alpha^m},
 \enddisplaymath
 whenever $a\in B_m$ and $b\in B_n$. 
 Given a pair $(X,\alpha)$ consisting of a Noetherian scheme $X$ and 
 an automorphism $\alpha$ of $X$,  an invertible sheaf ${\cal L}$ on $X$
 is called $\alpha$-ample, if for every coherent sheaf ${\cal F}$ on $X$,
 the cohomology group $H^q(X, ~{\cal L}\otimes {\cal L}^{\alpha}\otimes\dots
\otimes  {\cal L}^{\alpha^{ n-1}}\otimes {\cal F})$  vanishes for $q>0$ and
$n>>0$. 
(Notice,  that if $\alpha$ is trivial,  this definition is equivalent to the
usual definition of ample invertible sheaf,  see  [Serre 1955]  \cite{Ser1}.)  
If  $\alpha: X\to X$ is  an automorphism of a projective scheme $X$
over $k$  and ${\cal L}$ is  an  $\alpha$-ample invertible sheaf on $X$,
then  $B(X, {\cal L}, \alpha)$ is called a twisted homogeneous coordinate
ring of $X$.   
\end{dfn}
\begin{thm}\label{thm5.3.2}
{\bf ([Artin \& van den Bergh  1990]  \cite{ArtVdb1})}
\displaymath
\Mod~(B(X, {\cal L}, \alpha)) ~/~\Tors \cong \Coh~(X).  
\enddisplaymath
\end{thm}
\begin{rmk}
\textnormal{
Theorem \ref{thm5.3.2}  extends Theorem \ref{thm5.3.1}  to the 
\linebreak
non-commutative
rings;   hence the name for ring  $B(X, {\cal L}, \alpha)$. 
The question of which invertible sheaves are $\alpha$-ample is fairly
subtle, and there is   no characterization of the automorphisms $\alpha$
for which such an invertible sheaf exists.  However, in many important
special cases this problem is solvable,  see   [Artin \& van den Bergh  1990]  \cite{ArtVdb1},  Corollary 1.6.
}
\end{rmk}
\begin{rmk}\label{rmk5.3.2}
\textnormal{
In practice,  any  twisted homogeneous coordinate ring 
\linebreak
$B(X, {\cal L}, \alpha)$
of a projective scheme $X$ can be  constructed as follows. 
Let $R$ be a commutative graded ring,  such that $X=Spec~(R)$. 
Consider  the ring $B(X, {\cal L}, \alpha):= R[t,t^{-1}; \alpha]$,
where  $R[t,t^{-1}; \alpha]$ is the ring of  skew
Laurent polynomials defined by the commutation relation
$b^{\alpha}t=tb$,
for all $b\in R$;  here   $b^{\alpha}\in R$ is the image of $b$ under automorphism
$\alpha$.   The ring  $B(X, {\cal L}, \alpha)$ satisfies the isomorphism 
$\Mod~(B(X, {\cal L}, \alpha)) / \Tors \cong \Coh~(X)$,  i.e. is the
twisted homogeneous coordinate ring of projective scheme $X$, 
 see Lemma \ref{lem5.3.1}.
 }
\end{rmk}
\begin{exm}
\textnormal{
Let $k$ be a field and $U_{\infty}(k)$ the algebra of polynomials
over $k$ in two non-commuting variables $x_1$ and $x_2$,  and a quadratic relation
$x_1x_2-x_2x_1-x_1^2=0$;  let ${\Bbb P}^1(k)$ be the projective
line over $k$.  Then $B(X, {\cal L}, \alpha)=U_{\infty}(k)$ and $X={\Bbb P}^1(k)$
satisfy equation $\Mod~(B(X, {\cal L}, \alpha)) / \Tors \cong \Coh~(X)$.  
The  ring  $U_{\infty}(k)$  corresponds to the automorphism $\alpha(u)=u+1$ 
of the projective line ${\Bbb P}^1(k)$.   Indeed,  $u=x_2x_1^{-1}=x_1^{-1}x_2$
and,  therefore,  $\alpha$ maps $x_2$ to $x_1+x_2$;  if  one substitutes
$t=x_1,  b=x_2$ and $b^{\alpha}=x_1+x_2$ in equation $b^{\alpha}t=tb$ (see Remark \ref{rmk5.3.2}),   
then  one  gets the defining  relation  $x_1x_2-x_2x_1-x_1^2=0$ for  the algebra    $U_{\infty}(k)$.    
}
\end{exm}
To get a $C^*$-algebra from the ring $B(X, {\cal L}, \alpha)$,   we shall consider  infinite-dimensional 
 representations  of $B(X, {\cal L}, \alpha)$ by bounded linear operators  on a Hilbert space ${\cal H}$;
as usual,   let    ${\cal B}({\cal H})$ stay for  the algebra of   all  bounded linear  operators on  ${\cal H}$.
For a  ring of skew Laurent polynomials $R[t, t^{-1};  \alpha]$  described in Remark \ref{rmk5.3.2},  
we shall consider a homomorphism 
\displaymath
\rho: R[t, t^{-1};  \alpha]\longrightarrow {\cal B}({\cal H}). 
\enddisplaymath
Recall  that algebra ${\cal B}({\cal H})$ is endowed  with a $\ast$-involution;
such an  involution  is the adjoint with respect to   the scalar product on 
the Hilbert space ${\cal H}$. 
\begin{dfn}
The representation $\rho$ will be called  $\ast$-coherent if:

\medskip
(i)  $\rho(t)$ and $\rho(t^{-1})$ are unitary operators,  such that
$\rho^*(t)=\rho(t^{-1})$; 

\smallskip
(ii) for all $b\in R$ it holds $(\rho^*(b))^{\alpha(\rho)}=\rho^*(b^{\alpha})$, 
where $\alpha(\rho)$ is an automorphism of  $\rho(R)$  induced by $\alpha$. 
\end{dfn}
\begin{exm}
\textnormal{
The  ring  $U_{\infty}(k)$   has no    $\ast$-coherent representations.
 Indeed,  involution  acts  on the generators of $U_{\infty}(k)$
 by formula $x_1^*=x_2$;   the latter does not  preserve the defining relation 
 $x_1x_2-x_2x_1-x_1^2=0$.
 }
\end{exm}
\begin{dfn}
By a Serre $C^*$-algebra ${\cal A}_X$ of the projective scheme $X$
one understands   the norm-closure of an  $\ast$-coherent representation 
$\rho(B(X, {\cal L}, \alpha))$ of the twisted homogeneous coordinate ring  
$B(X, {\cal L}, \alpha)\cong R[t, t^{-1};  \alpha]$  of  scheme $X$. 
\end{dfn}
\begin{exm}\label{exm5.3.3}
\textnormal{
For $X\cong {\cal E}_{\tau}$ is an elliptic curve,   the ring  $R[t, t^{-1};  \alpha]$  is isomorphic to 
the Sklyanin algebra,   see Section 5.1.1.   For such algebras there exists a $\ast$-coherent representation
{\it ibid.};   the resulting Serre $C^*$-algebra ${\cal A}_X\cong {\cal A}_{\theta}$,
where ${\cal A}_{\theta}$  is  the noncommutative torus. 
 }
\end{exm}
\begin{rmk}\label{rmk5.3.3}
{\bf (Functor on complex projective varieties)} 
\textnormal{
If {\bf Proj-Alg} is the category of all complex projective varieties $X$ (of dimension $n$)
and {\bf C*-Serre} the category of all Serre $C^*$-algebras ${\cal A}_X$,   then the  
formula $X\mapsto {\cal A}_X$ gives rise to a map   
\displaymath
F:  \hbox{{\bf Proj-Alg}}\longrightarrow \hbox{{\bf C*-Serre}}.
\enddisplaymath
 The map $F$ is actually a functor which takes isomorphisms between projective 
 varieties to the stable isomorphisms (Morita equivalences) between the corresponding
 Serre $C^*$-algebras;   the proof repeats   the argument for elliptic curves given in Section 5.1.1
 and is left to the reader.    
 }
\end{rmk}
Let  $\Spec~(B(X, {\cal L}))$ be the space of all prime ideals 
of the commutative homogeneous coordinate ring $B(X, {\cal L})$
of a complex projective variety $X$,  see Theorem \ref{thm5.3.1}. 
 To get an analog of the classical  formula
 \displaymath 
 X\cong \Spec~(B(X, {\cal L}))
 \enddisplaymath
 for the Serre $C^*$-algebras ${\cal A}_X$,  we shall  recall that for each  
   continuous homomorphism  $\alpha: G\to Aut~({\cal A})$ of a locally compact
group $G$ into the group of automorphisms of a $C^*$-algebra  ${\cal A}$,
there exists a crossed product $C^*$-algebra ${\cal A}\rtimes_{\alpha}G$,
see e.g.  Section 3.2.   Let $G={\Bbb Z}$ and let $\hat {\Bbb Z}\cong S^1$ be
its Pontryagin dual.  
We shall write $\Irred$ for the set of all irreducible representations of  given
$C^*$-algebra. 
\begin{thm}\label{thm5.3.3}
For each Serre $C^*$-algebra  ${\cal A}_X$  there exists 
$\hat\alpha\in Aut~({\cal A}_X)$,  such that:
\displaymath
X\cong  \Irred ~({\cal A}_X\rtimes_{\hat\alpha} \hat {\Bbb Z}). 
\enddisplaymath
\end{thm}
\begin{rmk}
\textnormal{
Note that a naive generalization $X\cong \Spec ~({\cal A}_X)$ is wrong,
because most of the Serre $C^*$-algebras are simple,  i.e. have no ideals
whatsoever. 
}
\end{rmk}

\subsection{Proof of theorem \ref{thm5.3.3}}
\begin{lem}\label{lem5.3.1}
 $B(X, {\cal L}, \alpha)\cong R[t,t^{-1}; \alpha]$,  where $X=\Spec~(R)$.
\end{lem}
{\it Proof.}  
Let us  write the twisted homogeneous coordinate ring  $B(X, {\cal L}, \alpha)$ of projective variety
$X$ in the following  form:
\displaymath
B(X, {\cal L}, \alpha)=\bigoplus_{n\ge 0} H^0(X, {\goth B}_n),
\enddisplaymath
where ${\goth B}_n={\cal L}\otimes {\cal L}^{\alpha}\otimes\dots
\otimes  {\cal L}^{\alpha^{ n}}$ and  $H^0(X, {\goth B}_n)$ is the zero
sheaf cohomology of  $X$,  i.e. the space of sections $\Gamma(X, {\goth B}_n)$;
compare with  [Artin \& van den Bergh  1990]  \cite{ArtVdb1},  formula (3.5). 
If one denotes by ${\cal O}$  the structure sheaf of $X$, then 
\displaymath
{\goth B}_n={\cal O}t^n
\enddisplaymath
can be interpreted  as a free left  ${\cal O}$-module of rank one with basis $\{t^n\}$,
see   [Artin \& van den Bergh  1990]  \cite{ArtVdb1},  p. 252. 
Recall, that spaces  $B_i=H^0(X, {\goth B}_i)$ have  been endowed with the multiplication
rule  between the   sections  $a\in B_m$ and $b\in B_n$,  see Definition \ref{dfn5.3.1};
such a rule translates into the formula
\displaymath
at^mbt^n=ab^{\alpha^m}t^{m+n}. 
\enddisplaymath
One can eliminate $a$ and $t^n$ on   the both sides of the above equation;  
this operation gives us the following equation
\displaymath
t^mb=b^{\alpha^m}t^m. 
\enddisplaymath
First notice, that our ring  $B(X, {\cal L}, \alpha)$ contains a commutative
subring $R$, such that $\Spec~(R)=X$. 
Indeed, let $m=0$ in formula $t^mb=b^{\alpha^m}t^m$;  then $b=b^{Id}$ and, thus, $\alpha=Id$.
We conclude therefore, that $R=B_0$ is a commutative subring of  
$B(X, {\cal L}, \alpha)$,  and $\Spec~(R)=X$.

Let us show that equations  $b^{\alpha}t=tb$ of Remark \ref{rmk5.3.2}  and  
$t^mb=b^{\alpha^m}t^m$  are equivalent. 
First, let us show that  $b^{\alpha}t=tb$  implies $t^mb=b^{\alpha^m}t^m$.   Indeed,  equation 
 $b^{\alpha}t=tb$ can be written as $b^{\alpha}=tbt^{-1}$.   Then: 
\displaymath
\left\{
\begin{array}{ccc}
b^{\alpha^2} &=& tb^{\alpha}t^{-1}= t^2 bt^{-2},\\
b^{\alpha^3} &=& tb^{\alpha^2}t^{-1}= t^3 b t^{-3},\\
&\vdots&\\
b^{\alpha^m} &=& tb^{\alpha^{m-1}}t^{-1}= t^m b t^{-m}. 
\end{array}
\right.
\enddisplaymath
The last equation  of the above system  is  equivalent to equation $t^mb=b^{\alpha^m}t^m$. 
The converse is evident;  one sets $m=1$ in $t^mb=b^{\alpha^m}t^m$ and obtains
equation $b^{\alpha}t=tb$.  Thus,    $b^{\alpha}t=tb$  and  $t^mb=b^{\alpha^m}t^m$  are equivalent
equations.  
 It is easy now to establish an isomorphism  $B(X, {\cal L}, \alpha)\cong R[t,t^{-1}; \alpha]$.
For that,  take  $b\in R\subset B(X, {\cal L}, \alpha)$;  then $B(X, {\cal L}, \alpha)$ 
coincides  with the ring of the skew  Laurent polynomials  $R[t,t^{-1}; \alpha]$,
since the  commutation relation  $b^{\alpha}t=tb$  is equivalent to equation $t^mb=b^{\alpha^m}t^m$.
 Lemma \ref{lem5.3.1} follows. 
 $\square$

\begin{lem}\label{lem5.3.2}
${\cal A}_X\cong C(X)\rtimes_{\alpha} {\Bbb Z}$,  where $C(X)$ is the
$C^*$-algebra of all continuous complex-valued functions on $X$
and $\alpha$  is a $\ast$-coherent  automorphism of  $X$.   
\end{lem}
{\it Proof.} 
By definition of the Serre algebra ${\cal A}_X$,   the ring of skew Laurent 
polynomials $R[t, t^{-1}; \alpha]$  is  dense in ${\cal A}_X$;  roughly
speaking, one has to show that this property defines a crossed product
structure on ${\cal A}_X$.  We shall proceed in the following steps.

\medskip
(i) Recall  that  $R[t, t^{-1}; \alpha]$ consists of the finite sums
\displaymath
\sum  b_k t^k,  \qquad b_k\in R,
\enddisplaymath
subject to the commutation relation
\displaymath
b_k^{\alpha}t=tb_k. 
\enddisplaymath
Because of  the  $\ast$-coherent representation,
there is also an involution on $R[t, t^{-1}; \alpha]$, subject to the
following rules
\displaymath
\left\{
\begin{array}{ccccc}
 &(i)&   t^* &=& t^{-1},\\
 &(ii)&  (b_k^*)^{\alpha} &=& (b_k^{\alpha})^*.
\end{array}
\right.
\enddisplaymath

\bigskip
(ii)  Following [Williams  2007]  \cite{W}, p.47,   we shall consider the set 
$C_c({\Bbb Z}, R)$ of continuous functions from 
${\Bbb Z}$ to $R$ having a compact support;  
then the finite  sums  can be viewed as 
elements of  $C_c({\Bbb Z}, R)$ via  the identification 
\displaymath
k\longmapsto b_k.
\enddisplaymath
It can be verified,   that multiplication operation of the finite sums 
translates into a convolution product of functions $f,g\in C_c({\Bbb Z}, R)$
given by the formula
\displaymath
(f g)(k) = \sum_{l\in {\Bbb Z}} f(l) t^l g(k-l) t^{-l}, 
\enddisplaymath
while involution  translates into an involution on $C_c({\Bbb Z}, R)$
 given by the formula
\displaymath
f^*(k) = t^k f^*(-k) t^{-k}.
\enddisplaymath
It is easy to see, that the multiplication given by the convolution product 
and involution  turn $C_c({\Bbb Z}, R)$ into an $\ast$-algebra,
which is isomorphic to  the algebra $R[t, t^{-1}; \alpha]$.

\bigskip
(iii)  There exists the standard construction of a norm on  
$C_c({\Bbb Z}, R)$;  we omit it here referring the reader 
to [Williams  2007]  \cite{W}, Section 2.3.  The completion of $C_c({\Bbb Z}, R)$
in that norm defines a crossed product $C^*$-algebra 
$R\rtimes_{\alpha}{\Bbb Z}$ [Williams  2007]  \cite{W}, Lemma 2.27.

\bigskip
(iv)  Since $R$ is a commutative $C^*$-algebra
and $X=\Spec~(R)$,  one concludes that $R\cong C(X)$.
Thus, one obtains ${\cal A}_X=C(X)\rtimes_{\alpha}{\Bbb Z}$.
Lemma \ref{lem5.3.2}  follows. 
$\square$

\begin{rmk}\label{rmk5.3.5}
\textnormal{
It is easy to prove,   that  equations $b_k^{\alpha}t=tb_k$ and 
 $t^* = t^{-1}$    imply  equation $(b_k^*)^{\alpha} = (b_k^{\alpha})^*$;
in other words,  if involution does not commute with automorphism $\alpha$,
representation $\rho$  cannot  be  unitary,  i.e. $\rho^*(t)\ne\rho(t^{-1})$.    
 }
\end{rmk}
\begin{lem}\label{lem5.3.3}
There exists  $\hat\alpha\in Aut~({\cal A}_X)$, such that: 
\displaymath
X\cong \Irred~({\cal A}_X\rtimes_{\hat\alpha} \hat {\Bbb Z}).
\enddisplaymath
\end{lem}
{\it Proof.}  The above formula is an implication of the Takai duality 
for the crossed products, see e.g.  [Williams  2007]  \cite{W},  Section 7.1;  for the sake of clarity, 
we shall   repeat  this construction.  
Let $(A,G,\alpha)$ be a $C^*$-dynamical system with $G$ locally compact
abelian group; let $\hat G$ be the dual of $G$. For each $\gamma\in \hat G$,
one can define a map $\hat a_{\gamma}: C_c(G,A)\to C_c(G,A)$
given by the formula:
\displaymath
\hat a_{\gamma}(f)(s)=\bar\gamma(s)f(s), \qquad\forall s\in G.
 \enddisplaymath
In fact, $\hat a_{\gamma}$ is a $\ast$-homomorphism, since it
respects the convolution product and involution on $C_c(G,A)$
[Williams  2007]  \cite{W}.  Because the crossed product $A\rtimes_{\alpha}G$
is the closure of $C_c(G,A)$, one gets an extension of $\hat a_{\gamma}$
to an element of $Aut~(A\rtimes_{\alpha}G)$ and, therefore, a 
homomorphism:
\displaymath
\hat\alpha: \hat G\to Aut~(A\rtimes_{\alpha}G).
 \enddisplaymath
The {\it Takai duality} asserts, that
\displaymath
(A\rtimes_{\alpha} G)\rtimes_{\hat\alpha}\hat G\cong A\otimes {\cal K}(L^2(G)),
 \enddisplaymath
where ${\cal K}(L^2(G))$ is the algebra of compact operators on the
Hilbert space $L^2(G)$. 
Let us substitute $A=C_0(X)$   and   $G={\Bbb Z}$ in the above equation; 
one gets the following isomorphism
\displaymath
(C_0(X)\rtimes_{\alpha} {\Bbb Z})\rtimes_{\hat\alpha}\hat {\Bbb Z}\cong C_0(X)\otimes 
{\cal K}(L^2({\Bbb Z})).
 \enddisplaymath
Lemma \ref{lem5.3.2} asserts that $C_0(X)\rtimes_{\alpha}{\Bbb Z}\cong {\cal A}_X$;
 therefore  one arrives at the following isomorphism 
\displaymath
{\cal A}_X \rtimes_{\hat\alpha} \hat {\Bbb Z}\cong C_0(X)\otimes 
{\cal K}(L^2({\Bbb Z})).
 \enddisplaymath
Consider the set of all irreducible representations of the $C^*$-algebras
in the above equation;  then one gets the following equality of representations
\displaymath
\Irred~({\cal A}_X \rtimes_{\hat\alpha} \hat {\Bbb Z}) =\Irred~(C_0(X)\otimes 
{\cal K}(L^2({\Bbb Z}))).
 \enddisplaymath
Let $\pi$ be a representation of the tensor product 
$C_0(X)\otimes {\cal K}(L^2({\Bbb Z}))$ on the Hilbert
space ${\cal H}\otimes L^2({\Bbb Z})$;  then $\pi=\varphi\otimes\psi$,
where $\varphi: C_0(X)\to {\cal B}({\cal H})$ and $\psi: {\cal K}\to {\cal B}(L^2({\Bbb Z}))$. 
It is known, that the only irreducible representation of the algebra of 
compact operators is the identity representation. Thus,  one gets:
\begin{eqnarray}
\Irred~(C_0(X)\otimes {\cal K}(L^2({\Bbb Z}))) &=&
\Irred~(C_0(X))\otimes \{pt\}=\nonumber\\
       &=& \Irred(C_0(X)).\nonumber
 \end{eqnarray}
Further, the $C^*$-algebra $C_0(X)$ is commutative,  hence the 
following equations are true
\displaymath
\Irred~(C_0(X))=\Spec~(C_0(X))=X.
 \enddisplaymath
 Putting together the last three equations,   one obtains:
\displaymath
\Irred~({\cal A}_X \rtimes_{\hat\alpha} \hat {\Bbb Z})  \cong    X.
 \enddisplaymath
The conclusion of lemma \ref{lem5.3.3} follows from the above equation. 
$\square$

\bigskip
Theorem \ref{thm5.3.3} follows from Lemma \ref{lem5.3.3}.
$\square$

 \index{real multiplication}

\subsection{Real multiplication revisited}
We shall  test Theorem \ref{thm5.3.3} for  ${\cal A}_X\cong {\cal A}_{RM}$,  
i.e  a  noncommutative torus  with real multiplication;   notice that ${\cal A}_{RM}$ 
is the Serre $C^*$-algebra,  see Example \ref{exm5.3.3}. 
\begin{thm}\label{thm5.3.4}
$$\Irred~({\cal A}_{RM}\rtimes_{\hat\alpha} \hat {\Bbb Z})\cong {\cal E}(K),$$
where ${\cal E}(K)$ is  non-singular elliptic curve defined over a field of algebraic 
numbers  $K$. 
\end{thm}
{\it Proof.}  
We shall view the crossed product ${\cal A}_{RM}\rtimes_{\hat\alpha} \hat {\Bbb Z}$
as a $C^*$-dynamical system $({\cal A}_{RM}, \hat {\Bbb Z}, \hat\alpha)$,  
see [Williams  2007]  \cite{W} for the details. 
Recall that the irreducible representations of 
$C^*$-dynamical system  $({\cal A}_{RM}, \hat {\Bbb Z}, \hat\alpha)$
 are in the one-to-one correspondence with the minimal sets of the 
dynamical system (i.e. closed $\hat\alpha$-invariant sub-$C^*$-algebras
of ${\cal A}_{RM}$ not containing a smaller object with the same property).  
 To calculate the minimal sets of $({\cal A}_{RM}, \hat {\Bbb Z}, \hat\alpha)$,
let $\theta$ be quadratic irrationality such that ${\cal A}_{RM}\cong {\cal A}_{\theta}$.
It is known that every non-trivial sub-$C^*$-algebra  of ${\cal A}_{\theta}$ has
the form ${\cal A}_{n\theta}$ for some positive integer $n$, 
see [Rieffel 1981] \cite{Rie2},  p. 419.   It is easy to deduce that the
{\it maximal}  proper sub-$C^*$-algebra of ${\cal A}_{\theta}$ has the form 
${\cal A}_{p\theta}$, where $p$ is   a prime number.  
(Indeed, each composite $n=n_1n_2$ cannot be maximal since
${\cal A}_{n_1n_2\theta}\subset {\cal A}_{n_1\theta}\subset {\cal A}_{\theta}$
or ${\cal A}_{n_1n_2\theta}\subset {\cal A}_{n_2\theta}\subset {\cal A}_{\theta}$,
where all inclusions are strict.) 
We claim that  $({\cal A}_{p\theta}, \hat {\Bbb Z}, \hat\alpha^{\pi(p)})$
is the minimal $C^*$-dynamical system,  where $\pi(p)$ is certain
power of the automorphism $\hat\alpha$.   Indeed,  the automorphism 
$\hat\alpha$ of ${\cal A}_{\theta}$  corresponds to multiplication by
the fundamental unit, $\varepsilon$, of   pseudo-lattice  $\Lambda=
{\Bbb Z}+\theta {\Bbb Z}$.   It is known that certain power, $\pi(p)$,
of $\varepsilon$  coincides with the fundamental unit of pseudo-lattice
${\Bbb Z}+(p\theta){\Bbb Z}$, see e.g. [Hasse 1950]  \cite{HA},  p. 298.  
Thus  one gets the minimal $C^*$-dynamical system 
$({\cal A}_{p\theta}, \hat {\Bbb Z}, \hat\alpha^{\pi(p)})$,  which is defined 
on the sub-$C^*$-algebra ${\cal A}_{p\theta}$ of  ${\cal A}_{\theta}$.  
Therefore we have an isomorphism    
\displaymath
\Irred ~({\cal A}_{RM}\rtimes_{\hat\alpha} \hat {\Bbb Z})\cong
\bigcup_{p\in {\cal P}} \Irred ~( {\cal A}_{p\theta}\rtimes_{\hat\alpha^{\pi(p)}} \hat {\Bbb Z}),
 \enddisplaymath
where ${\cal P}$ is the set of all (but a finite number) of primes.  
To simplify the RHS of the above equation, let us introduce some notation.  
Recall that matrix form of the fundamental unit $\varepsilon$ of
pseudo-lattice $\Lambda$ coincides with the matrix $A$, see above. 
For each prime $p\in {\cal P}$ consider the matrix 
\displaymath
L_p=\left(\matrix{tr~(A^{\pi(p)})-p & p\cr tr~(A^{\pi(p)})-p-1 & p}\right),
 \enddisplaymath
where $tr$ is the trace of matrix.  Let us show, that 
\displaymath
{\cal A}_{p\theta}\rtimes_{\hat\alpha^{\pi(p)}} \hat {\Bbb Z}\cong
{\cal A}_{\theta}\rtimes_{L_p} \hat {\Bbb Z},
 \enddisplaymath
where $L_p$ is an endomorphism of ${\cal A}_{\theta}$ (of degree $p$)
induced by matrix $L_p$.
Indeed,  because $deg~(L_p)=p$ the endomorphism $L_p$ maps
pseudo-lattice $\Lambda={\Bbb Z}+\theta {\Bbb Z}$ to a sub-lattice
of index $p$;  any such can be written in the form 
$\Lambda_p={\Bbb Z}+(p\theta) {\Bbb Z}$, see e.g. 
[Borevich \& Shafarevich  1966]  \cite{BS},  p.131.
Notice that pseudo-lattice $\Lambda_p$ corresponds
to the sub-$C^*$-algebra ${\cal A}_{p\theta}$  of algebra ${\cal A}_{\theta}$  
and $L_p$ induces a shift automorphism of ${\cal A}_{p\theta}$,
see e.g.   [Cuntz 1977]  \cite{Cun1}  beginning of Section 2.1
for terminology and  details of this construction.   
It is not hard to see, that the shift automorphism coincides
with $\hat\alpha^{\pi(p)}$.   Indeed, it is verified directly
that $tr~(\hat\alpha^{\pi(p)})=tr~(A^{\pi(p)})=tr~(L_p)$;
thus one gets a bijection between powers of $\hat\alpha^{\pi(p)}$
and such of $L_p$.   But $\hat\alpha^{\pi(p)}$ corresponds to
the fundamental unit of pseudo-lattice $\Lambda_p$;  therefore 
the shift automorphism induced by $L_p$
must coincide with $\hat\alpha^{\pi(p)}$.   The required isomorphism 
is proved and,  therefore,   our last formula  can be written in the form
\displaymath
\Irred ~({\cal A}_{RM}\rtimes_{\hat\alpha} \hat {\Bbb Z})\cong
\bigcup_{p\in {\cal P}} \Irred ~( {\cal A}_{RM}\rtimes_{L_p} \hat {\Bbb Z}).  
 \enddisplaymath
To calculate  irreducible representations  of the crossed product $C^*$-algebra 
 ${\cal A}_{RM}\rtimes_{L_p}{\hat \Bbb Z}$ at the RHS of  the above equation,
 recall that such are in a one-to-one correspondence with 
 the set of  invariant measures on a subshift of finite type given by the positive integer 
 matrix $L_p$,  see   [Bowen \& Franks 1977]  \cite{BowFra1} and 
  [Cuntz 1977]  \cite{Cun1};  the  measures make an  abelian group under the 
  addition operation.   Such a  group is  isomorphic to  ${\Bbb Z}^2~/~(I-L_p) {\Bbb Z}^2$,
where $I$ is the identity matrix, see  [Bowen \& Franks 1977]  \cite{BowFra1},  Theorem 2.2. 
Therefore our last equation can be written in  the form
\displaymath
\Irred ~({\cal A}_{RM}\rtimes_{\hat\alpha} \hat {\Bbb Z})\cong
\bigcup_{p\in {\cal P}} {{\Bbb Z}^2\over (I-L_p) {\Bbb Z}^2}.
\enddisplaymath
Let ${\cal E}(K)$ be a non-singular elliptic curve defined over 
the algebraic number field $K$;    let ${\cal E}({\Bbb F}_p)$ 
be the reduction of ${\cal E}(K)$ modulo prime ideal over a 
``good'' prime number $p$.  Recall that $|{\cal E}({\Bbb F}_p)|=
det~(I-Fr_p)$,  where  $Fr_p$  is an integer two-by-two matrix
corresponding to the action of Frobenius endomorphism on the
$\ell$-adic cohomology of ${\cal E}(K)$,  see e.g. 
[Tate  1974]  \cite{Tat1}, p. 187.  
Since $|{\Bbb Z}^2 / (I-L_p){\Bbb Z}^2|=det~(I-L_p)$,
one can identify $Fr_p$ and $L_p$  and,  therefore,  
one obtains an isomorphism ${\cal E}({\Bbb F}_p)\cong 
{\Bbb Z}^2 / (I-L_p){\Bbb Z}^2$.  
Thus  our equation  can be written in  the form
\displaymath
\Irred ~({\cal A}_{RM}\rtimes_{\hat\alpha} \hat {\Bbb Z})\cong
\bigcup_{p\in {\cal P}} {\cal E}({\Bbb F}_p).
\enddisplaymath
Finally,   consider an arithmetic scheme, $X$,  corresponding to ${\cal E}(K)$;
the latter fibers over ${\Bbb Z}$,  see [Silverman 1994]  \cite{S},  Example 4.2.2
for the details.   It can be immediately seen,  that the RHS of our last equation 
coincides with the scheme $X$, where  the regular fiber over $p$ corresponds
to ${\cal E}({\Bbb F}_p)$ {\it ibid.}  The argument finishes the proof of 
Theorem \ref{thm5.3.4}. 
$\square$

\vskip1cm\noindent
{\bf Guide to the literature.}
The standard  reference to complex projective varieties is the monograph 
[Hartshorne  1977] \cite{H1}.  Twisted homogeneous coordinate rings of 
projective varieties are covered in the excellent survey  by 
 [Stafford \& van ~den ~Bergh  2001]  \cite{StaVdb1}.  The Serre $C^*$-algebras
 were introduced and studied in \cite{Nik11}.

 \index{mapping class group}

\section{Application:  Mapping class groups}
In the foreword  it was asked:  {\it Why does NCG matter?}  
We  shall answer this question by solving  a problem of classical 
geometry (Harvey's conjecture)  using  invariants attached to the functor    $F:$ {\bf Alg-Gen}  $\to$  {\bf AF-Toric},
see Theorem \ref{thm5.2.1};  the author  is unaware of a ``classical'' proof of this result. 
The invariant in question is the stable  isomorphism group of toric AF-algebra ${\Bbb A}_{\theta}$.

\subsection{Harvey's conjecture}
The mapping class group has  been introduced in the 1920-ies by M.~Dehn [Dehn  1938]   \cite{Deh1}. 
Such a  group, $Mod~(X)$, is defined as the group of isotopy classes of the
orientation-preserving diffeomorphisms of a two-sided closed surface $X$ of
genus $g\ge 1$.  The group is known to be  prominent in  algebraic geometry [Hain \& Looijenga 1997]   \cite{HaLo1}, 
topology [Thurston 1982]   \cite{Thu2}  and dynamics [Thurston 1988]   \cite{Thu1}.  When $X$
is a   torus, the $Mod~(X)$ is isomorphic to the group  $SL(2, {\Bbb Z})$. 
(The $SL(2, {\Bbb Z})$ is called a modular group, hence our notation for the mapping class group.) 
A little is known about the representations of $Mod~(X)$ beyond the case $g=1$.
Recall, that the group is called {\it linear},  if there exists a faithful
representation into the matrix group $GL(m, R)$, where $R$ is a commutative ring.
The braid groups are known  to be linear [Bigelow  2001]   \cite{Big1}.  Using a modification of the argument for 
the braid groups, it is possible to prove, that $Mod~(X)$ is linear in the case $g=2$
[Bigelow \&  Budney  2001]   \cite{BiBu1}. 
\begin{dfn}
{\bf ([Harvey 1979]   \cite{Har1},  p.267)} 
By Harvey's conjecture we understand the claim that 
the mapping class group is linear for  $g\ge 3$. 
\end{dfn}
Recall that  a covariant  functor $F:$ {\bf Alg-Gen}  $\to$  {\bf AF-Toric}
from a category of generic Riemann
surfaces (i.e. complex algebraic curves)  to a category of  toric AF-algebras was 
constructed in Section 5.2;   the functor   maps any  pair of isomorphic Riemann surfaces to a pair of 
 stably isomorphic (Morita equivalent)  toric  AF-algebras.  Since each isomorphism of Riemann surfaces is 
 given by an element of $Mod~(X)$ [Hain \& Looijenga 1997]   \cite{HaLo1},  it is natural to ask about 
 a representation of $Mod~(X)$  by the stable isomorphisms of  toric AF-algebras.  
 Recall that the stable isomorphisms of toric AF-algebras are well understood and
surprisingly simple;   provided  the automorphism group of
the algebra is trivial (this is true for a generic algebra),  its group of stable isomorphism  
admits a faithful representation into the  matrix group $GL(m, {\Bbb Z})$,  see e.g.  [Effros  1981]   \cite{E}.  
This fact, combined  with the properties of functor $F$,  implies a  positive
 solution to the Harvey conjecture.   
\begin{thm}\label{thm5.4.1}
For every surface $X$ of genus $g\ge 2$,  there exists a faithful representation 
$\rho: Mod~(X)\rightarrow GL(6g-6,  {\Bbb Z})$.
\end{thm}

\subsection{Proof of Theorem \ref{thm5.4.1}}
Let {\bf AF-Toric}  denote  the set of all  toric AF-algebras of genus $g\ge 2$. 
Let $G$ be a finitely presented group and  
\displaymath
G\times\hbox{{\bf AF-Toric}} \longrightarrow
\hbox{{\bf AF-Toric}}  
 \enddisplaymath
 be  its action on {\bf AF-Toric}   by  the stable isomorphisms (Morita equivalences)  of toric AF-algebras; 
in other words,   $\gamma ({\Bbb A}_{\theta})\otimes {\cal K}\cong
{\Bbb A}_{\theta}\otimes {\cal  K}$  for all $\gamma\in G$ and all 
${\Bbb A}_{\theta}\in$ {\bf AF-Toric}.   The following preparatory lemma will be 
important.
\begin{lem}\label{lem5.4.1}
For each ${\Bbb A}_{\theta}\in$ {\bf AF-Toric},   there exists a 
representation 
\displaymath
\rho_{{\Bbb A}_{\theta}}: G\to GL(6g-6, {\Bbb Z}).
\enddisplaymath
\end{lem}
{\it Proof.}
The proof of  lemma is based on the following well known criterion
of the stable isomorphism for the (toric)  AF-algebras: a pair of such algebras
${\Bbb A}_{\theta}, {\Bbb A}_{\theta'}$ are stably isomorphic if and only
if their Bratteli diagrams coincide, except (possibly) a finite part
of the diagram, see e.g. [Effros  1981]   \cite{E},   Theorem 2.3. 
\begin{rmk}
\textnormal{
Note  that the order isomorphism  between  the dimension groups {\it ibid.}, 
translates  to the language of the Bratteli diagrams as stated.
}
\end{rmk}
Let $G$ be a finitely presented group on the generators  $\{\gamma_1, \dots, \gamma_m\}$
subject to relations $r_1,\dots,r_n$.  Let  ${\Bbb A}_{\theta}\in$  {\bf AF-Toric}.
Since $G$ acts on the toric AF-algebra ${\Bbb A}_{\theta}$ by stable isomorphisms,
the toric AF-algebras ${\Bbb A}_{\theta_1}:=\gamma_1({\Bbb A}_{\theta}),\dots,   
{\Bbb A}_{\theta_m}:=\gamma_m({\Bbb A}_{\theta})$ are stably isomorphic to 
${\Bbb A}_{\theta}$;  moreover,  by transitivity, they are also pairwise stably isomorphic.
Therefore, the Bratteli diagrams of ${\Bbb A}_{\theta_1},\dots, {\Bbb A}_{\theta_m}$ 
 coincide everywhere except, possibly, some finite parts.  
We shall denote by ${\Bbb A}_{\theta_{\max}}\in$ {\bf AF-Toric}
a toric AF-algebra,  whose Bratteli diagram is the maximal common part
of the Bratteli diagrams of ${\Bbb A}_{\theta_i}$ for $1\le i\le m$;
such a choice is unique and defined correctly because the set $\{{\Bbb A}_{\theta_i}\}$
is a finite set.  By the Definition \ref{dfn5.2.1}  of a toric AF-algebra,  the vectors 
$\theta_i=(1,\theta_1^{(i)},\dots,\theta_{6g-7}^{(i)})$
are related to the vector $\theta_{\max}=(1, \theta_1^{(\max)},\dots,\theta_{6g-7}^{(\max)})$ 
by the formula
\displaymath
\left(\matrix{1\cr \theta_1^{(i)}\cr\vdots\cr\theta_{6g-7}^{(i)}} \right)
=\underbrace{
\left(\matrix{0 &  0 & \dots & 0 & 1\cr
              1 &  0 & \dots & 0 & b_1^{(1)(i)}\cr
              \vdots &\vdots & &\vdots &\vdots\cr
              0 &  0 & \dots & 1 & b_{6g-7}^{(1)(i)}}\right)
\dots 
\left(\matrix{0 &  0 & \dots & 0 & 1\cr
              1 &  0 & \dots & 0 & b_1^{(k)(i)}\cr
              \vdots &\vdots & &\vdots &\vdots\cr
              0 &  0 & \dots & 1 & b_{6g-7}^{(k)(i)}}\right)
}_{A_i}
\left(\matrix{1\cr \theta^{(\max)}_1\cr\vdots\cr\theta^{(\max)}_{6g-7}} \right)
\enddisplaymath
The above expression can be written in the matrix form $\theta_i=A_i\theta_{\max}$, where 
$A_i\in GL(6g-6, {\Bbb Z})$.   Thus, one gets a matrix representation of the
generator $\gamma_i$,  given by the formula 
\displaymath
\rho_{{\Bbb A}_{\theta}}(\gamma_i):=A_i.
\enddisplaymath
The map $\rho_{{\Bbb A}_{\theta}}:  G\to GL(6g-6, {\Bbb Z})$ extends to the rest of the group $G$
 via its values on the generators;   namely,  for every $g\in G$ one sets $\rho_{{\Bbb A}_{\theta}}(g)= A_1^{k_1}\dots A_m^{k_m}$,
whenever $g=\gamma_1^{k_1}\dots \gamma_m^{k_m}$.  Let us verify, that the 
map $\rho_{{\Bbb A}_{\theta}}$ is a well defined homomorphism of groups $G$ and $GL(6g-6, {\Bbb Z})$.
Indeed, let us write $g_1=\gamma_1^{k_1}\dots\gamma_m^{k_m}$ and
$g_2=\gamma_1^{s_1}\dots\gamma_m^{s_m}$ for a pair of elements $g_1,g_2\in G$;
then their product 
$g_1g_2=\gamma_1^{k_1}\dots\gamma_m^{k_m}\gamma_1^{s_1}\dots\gamma_m^{s_m}=
\gamma_1^{l_1}\dots\gamma_m^{l_m}$,
where the last equality is obtained by a reduction of words using the 
relations $r_1,\dots,r_n$.  One can write relations $r_i$ in their matrix
form $\rho_{{\Bbb A}_{\theta}}(r_i)$; thus, one gets the matrix equality
$A_1^{l_1}\dots A_m^{l_m}= A_1^{k_1}\dots A_m^{k_m}A_1^{s_1}\dots A_m^{s_m}$.
It is immediate from the last  equation, that 
$\rho_{{\Bbb A}_{\theta}}(g_1g_2)=
A_1^{l_1}\dots A_m^{l_m}= A_1^{k_1}\dots A_m^{k_m}A_1^{s_1}\dots A_m^{s_m}=
\rho_{{\Bbb A}_{\theta}}(g_1)\rho_{{\Bbb A}_{\theta}}(g_2)$ for  $\forall g_1,g_2\in G$,
i.e.  $\rho_{{\Bbb A}_{\theta}}$ is a homomorphism. Lemma \ref{lem5.4.1}  follows.
$\square$

\bigskip\noindent
Let  {\bf AF-Toric-Aper} $\subset$ {\bf AF-Toric}  be a set consisting  of the toric AF-algebras,  whose  Bratteli diagrams
are {\it not}  periodic;  these are known as  non-stationary toric  AF-algebras (Section 3.5.2) 
and they are generic in the set {\bf AF-Toric}   endowed with the natural topology. 
\begin{dfn}
The action of group  $G$  on the toric AF-algebra  ${\Bbb A}_{\theta}\in$ {\bf AF-Toric}
 will be called free,    if   $\gamma ({\Bbb A}_{\theta})={\Bbb A}_{\theta}$ implies $\gamma=Id$.
\end{dfn}
\begin{lem}\label{lem5.4.2}
If  ${\Bbb A}_{\theta}\in$ {\bf AF-Toric-Aper}   and  the action of group $G$ 
on the ${\Bbb A}_{\theta}$ is free,  then  $\rho_{{\Bbb A}_{\theta}}$ is a 
faithful representation.
\end{lem}
{\it Proof.}
Since the action of $G$ is free, to prove that  $\rho_{{\Bbb A}_{\theta}}$  is faithful, 
it remains  to show, that in the formula  $\theta_i=A_i\theta_{\max}$, it holds
$A_i=I$, if and only if,  $\theta_i=\theta_{\max}$, where $I$
is the unit matrix.   Indeed,  it is immediate that $A_i=I$ implies $\theta_i=\theta_{\max}$.
Suppose now that  $\theta_i=\theta_{\max}$ and, let to the contrary, $A_i\ne I$. 
One gets $\theta_i=A_i \theta_{\max}=\theta_{\max}$.  Such an equation has  a non-trivial solution, 
if and only if,  the vector $\theta_{\max}$ has a periodic 
Jacobi-Perron fraction;  the period of  such a fraction  is given by the matrix $A_i$. This 
is impossible, since it has been assumed, that ${\Bbb A}_{\theta_{\max}}\in$ {\bf AF-Toric-Aper}.
The contradiction proves  Lemma \ref{lem5.4.2}.
$\square$

\bigskip\noindent
Let $G=Mod~(X)$,  where $X$ is a surface of genus $g\ge 2$. The group $G$
is finitely presented, see  [Dehn  1938]   \cite{Deh1}; it  acts on the Teichmueller space $T(g)$
by isomorphisms of the Riemann surfaces.  Moreover, the action of $G$ is free on a 
generic set,   $U\subset T(g)$,  consisting of the Riemann surfaces with the trivial group 
of  automorphisms.   On the other hand,  there exists a functor  
 \displaymath
 F:  ~\hbox{{\bf Alg-Gen}}\longrightarrow \hbox{{\bf AF-Toric}}
 \enddisplaymath
 between the Riemann surfaces (complex algebraic curves)  and 
toric AF-algebras,  see Theorem \ref{thm5.2.1}.
\begin{lem}\label{lem5.4.3}
The pre-image  $F^{-1}(\hbox{{\bf AF-Toric-Aper}})$ is a generic set in the space  $T(g)$.
\end{lem}
{\it Proof.}
Note,  that the set of stationary toric AF-algebras is a countable set.
The functor $F$ is a surjective map,  which is continuous with respect to the natural topology 
on the sets {\bf Alg-Gen}  and {\bf AF-Toric}.  Therefore,  the pre-image of the complement 
of a countable set is a generic set.  Lemma \ref{lem5.4.3} follows.
$\square$

\bigskip\noindent
Consider the set $U\cap F^{-1}(\hbox{{\bf AF-Toric-Aper}})$;
this set  is  non-empty,  since it is the intersection of  two generic subsets of $T(g)$, 
see Lemma \ref{lem5.4.3}.  Let 
\displaymath
S\in U\cap F^{-1}(\hbox{{\bf AF-Toric-Aper}})
\enddisplaymath
be a point   (a Riemann surface) in the above set.  
In  view of Lemma \ref{lem5.4.1},  group $G$ acts on the toric AF-algebra ${\Bbb A}_{\theta}=F(S)$ by the stable
isomorphisms.  By the construction,  the action is free and ${\Bbb A}_{\theta}\in$ {\bf AF-Toric-Aper}.
In view of Lemma \ref{lem5.4.2},   one gets a faithful representation $\rho=\rho_{{\Bbb A}_{\theta}}$
of the group $G\cong Mod~(X)$ into the matrix group $GL(6g-6, {\Bbb Z})$.  Theorem \ref{thm5.4.1}
is proved.
$\square$

\vskip1cm\noindent
{\bf Guide to the literature.}
The mapping class groups were introduced by M.~Dehn [Dehn  1938]   \cite{Deh1}. 
For a  primer  on the mapping class groups we refer the reader to the textbook  [Farb \& Margalit  2011]
\cite{FM}.   The Harvey conjecture was formulated in [Harvey 1979]   \cite{Har1}.
Some infinite-dimensional (asymptotic) faithfulness of the mapping class groups was proved by 
[Anderson  2006]  \cite{And1}.   A faithful representation of $Mod~(X)$ in the matrix group $GL(6g-6, {\Bbb Z})$
was constructed  in \cite{Nik12}.

\section*{Exercises}

\begin{enumerate}

\item
Prove that the skew-symmetric  relations 
\displaymath
\left\{
\begin{array}{ccc}
x_3x_1 &=& q_{13} x_1x_3,\\
x_4x_2 &=&  q_{24}x_2x_4,\\
x_4x_1 &=&  q_{14}x_1x_4,\\
x_3x_2 &=&  q_{23}x_2x_3,\\
x_2x_1&=&  q_{12}x_1x_2,\\
x_4x_3&=&  q_{34}x_3x_4,
\end{array}
\right.
\enddisplaymath
are invariant of the involution  $x_1^*=x_2, x_3^*=x_4$,     if and only if,  the following restrictions on 
the constants $q_{ij}$ hold
\displaymath
\left\{
\begin{array}{ccc}
q_{13} &=&  (\bar q_{24})^{-1},\\
q_{24} &=&  (\bar q_{13})^{-1},\\
q_{14} &= & (\bar q_{23})^{-1},\\
q_{23} &= & (\bar q_{14})^{-1},\\
q_{12} &= & \bar q_{12},\\
q_{34} &= & \bar q_{34},
\end{array}
\right.
\enddisplaymath
where $\bar q_{ij}$ means the complex conjugate of $q_{ij}\in {\Bbb C}\setminus\{0\}$.

\item
Prove that  a family of  free algebras  ${\Bbb C}\langle x_1,x_2,x_3,x_4\rangle$  modulo an ideal 
generated  by  six   skew-symmetric quadratic  relations
\displaymath
\left\{
\begin{array}{cc}
x_3x_1 &= \mu e^{2\pi i\theta}x_1x_3,\\
x_4x_2 &= {1\over \mu} e^{2\pi i\theta}x_2x_4,\\
x_4x_1 &= \mu e^{-2\pi i\theta}x_1x_4,\\
x_3x_2 &= {1\over \mu} e^{-2\pi i\theta}x_2x_3,\\
x_2x_1 &= x_1x_2,\\
x_4x_3 &= x_3x_4,
\end{array}
\right.
\enddisplaymath
consists of the pairwise non-isomorphic algebras for different values of 
$\theta\in S^1$ and $\mu\in (0,\infty)$.

\item
Prove that the system of relations for noncommutative torus ${\cal A}_{\theta}$
\displaymath
\left\{
\begin{array}{cc}
x_3x_1 &=  e^{2\pi i\theta}x_1x_3,\\
x_4x_2 &=  e^{2\pi i\theta}x_2x_4,\\
x_4x_1 &=  e^{-2\pi i\theta}x_1x_4,\\
x_3x_2 &=   e^{-2\pi i\theta}x_2x_3,\\
x_2x_1 &= x_1x_2=e,\\
x_4x_3 &= x_3x_4=e.
\end{array}
\right.
\enddisplaymath
is equivalent to the system of relations
\displaymath
\left\{
\begin{array}{ccc}
x_3x_1x_4 &=&  e^{2\pi i\theta}x_1,\\
x_4 &= & e^{2\pi i\theta}x_2x_4x_1,\\
x_4x_1x_3 &=&  e^{-2\pi i\theta}x_1,\\
x_2 &=&  e^{-2\pi i\theta}x_4x_2x_3,\\
x_1x_2   &=&   x_2x_1   =e,\\
 x_3x_4    &=&  x_4x_3  =e.
\end{array}
\right.
\enddisplaymath
(Hint:  use the last two relations.)

\item
Prove that the system of relations for the Sklyanin  $\ast$-algebra plus 
the scaled unit relation, i.e. 
\displaymath
\left\{
\begin{array}{cc}
x_3x_1 &= \mu e^{2\pi i\theta}x_1x_3,\\
x_4x_2 &= {1\over \mu} e^{2\pi i\theta}x_2x_4,\\
x_4x_1 &= \mu e^{-2\pi i\theta}x_1x_4,\\
x_3x_2 &= {1\over \mu} e^{-2\pi i\theta}x_2x_3,\\
x_2x_1 &= x_1x_2={1\over\mu}e,\\
x_4x_3 &= x_3x_4={1\over\mu}e
\end{array}
\right.
\enddisplaymath
is equivalent to the system
\displaymath
\left\{
\begin{array}{cc}
x_3x_1x_4 &= e^{2\pi i\theta}x_1,\\
x_4 &= e^{2\pi i\theta}x_2x_4x_1,\\
x_4x_1x_3 &= e^{-2\pi i\theta}x_1,\\
x_2 &= e^{-2\pi i\theta}x_4x_2x_3,\\
x_2x_1 &= x_1x_2={1\over\mu}e,\\
x_4x_3 &= x_3x_4={1\over\mu}e.
\end{array}
\right.
\enddisplaymath
(Hint: use   multiplication and cancellation involving the last two equations.)

\item
If {\bf Proj-Alg} is the category of all complex projective varieties $X$ (of dimension $n$)
and {\bf C*-Serre} the category of all Serre $C^*$-algebras ${\cal A}_X$,   then the  
formula $X\mapsto {\cal A}_X$ gives rise to a map   
\displaymath
F:  \hbox{{\bf Proj-Alg}}\longrightarrow \hbox{{\bf C*-Serre}}.
\enddisplaymath
 Prove that the map $F$ is actually a functor which takes isomorphisms between projective 
 varieties to the stable isomorphisms (Morita equivalences) between the corresponding
 Serre $C^*$-algebras.  (Hint:    repeat   the argument for elliptic curves given in Section 5.1.1.)

\item
Prove Remark \ref{rmk5.3.5},  i.e.  that  equations $b_k^{\alpha}t=tb_k$ and 
 $t^* = t^{-1}$    imply  equation $(b_k^*)^{\alpha} = (b_k^{\alpha})^*$.

\end{enumerate}





\chapter{Number Theory}
The most  elegant   functors  (with values in NCG)   are acting  on the arithmetic schemes $X$.   
We  start with the simplest case of  $X$  being  elliptic curve with
 complex multiplication by the number field $k={\Bbb Q}(f\sqrt{-D})$;
in this case $X\cong {\cal E}(K)$,  where $K$ is the Hilbert class field
of $k$  [Serre 1967] \cite{Ser2}.  We prove  in Section 6.1  that functor $F$ sends ${\cal E}(K)$
to  noncommutative torus with real multiplication by the number
 field   ${\Bbb Q}(f\sqrt{D})$.  It is proved in Section 6.2 that   the so-called 
 {\it arithmetic complexity} of such a torus   is  linked  by a simple formula  to the 
 rank of elliptic curve ${\cal E}(K)$ whenever  $D\equiv 3~\mod ~4$ is a prime number
 and $f=1$.  
In Section 6.3  we  introduce an $L$-function $L({\cal A}_{RM}, s)$ 
associated to the noncommutative torus with real multiplication and prove
that any such coincides with the classical Hasse-Weil function $L({\cal E}_{CM},s)$
of an elliptic curve with complex multiplication;  a surprising {\it localization formula}
tells us that the crossed products replace prime (or maximal) ideals familiar 
from the  commutative algebra.   
In Section 6.4 a functor $F:$ {\bf Alg-Num} $\to$ {\bf NC-Tor}  from a category of 
the finite Galois extensions $E$  of the field ${\Bbb Q}$ to the category of even-dimensional
noncommutative tori with real multiplication ${\cal A}_{RM}^{2n}$ is defined.
An $L$-function $L({\cal A}_{RM}^{2n}, s)$ is constructed and it is conjectured
that if ${\cal A}_{RM}^{2n}=F(E)$,  then   $L({\cal A}_{RM}^{2n}, s)\equiv L(\sigma,s)$,
where $L(\sigma, s)$ is the {\it Artin $L$-function} of  $E$  corresponding to 
an irreducible representation $\sigma: Gal~(E|{\Bbb Q})\to GL_n({\Bbb C})$.  
We prove the conjecture for $n=1$ (resp., $n=0$) and $E$ being the Hilbert class field of an
imaginary quadratic field $k$ (resp., field ${\Bbb Q}$).  
Thus we deal with an analog of the 
{\it Langlands program},  where the ``automorphic cuspidal representations of 
group $GL_n$'' are replaced by the noncommutative tori ${\cal A}_{RM}^{2n}$, 
see [Gelbart 1984]  \cite{Gel1} for an introduction to the Langlands program.    
In Section 6.5  we compute  the number of  points of  projective variety $V({\Bbb F}_q)$ 
over a  finite field  ${\Bbb F}_q$  in terms of  invariants  of the 
{\it Serre  $C^*$-algebra}   associated to the complex projective variety  $V({\Bbb C})$, 
 see  Section  5.3.1;   the calculation involves an explicit formula for the traces of 
 {\it Frobenius map} of $V({\Bbb F}_q)$  being  linked  to the {\it Weil Conjectures},
see e.g.  [Hartshorne 1977]  \cite{H1},  Appendix C for an introduction.             
Finally, in Section 6.6 we apply  our  functor  $F:$  {\bf Ell} $\to$  {\bf NC-Tor} to a problem of 
the   {\it transcendental number  theory},  see e.g. [Baker  1975]  \cite{BA} for
an introduction.  Namely,  we  use the formula
$F({\cal E}_{CM}^{(-D,f)})={\cal A}_{RM}^{(D,f)}$ of Section 6.1  to prove that the transcendental function
${\cal J}(\theta,\varepsilon)=e^{2\pi i\theta+\log\log\varepsilon}$
takes algebraic values for the algebraic arguments $\theta$ and $\varepsilon$.
Moreover,  these values of  ${\cal J}(\theta,\varepsilon)$ belong to the Hilbert
class field of the imaginary quadratic field ${\Bbb Q}(\sqrt{-D})$ for all but a finite
set  of values of $D$.

 \index{transcendental number theory}
 \index{real multiplication}
 \index{Frobenius map}
 \index{Langlands program}

 \index{complex multiplication}

\section{Complex multiplication}
We recall that an {\it elliptic curve} is  the subset of the complex projective plane  of the form
${\cal E}({\Bbb C})=\{(x,y,z)\in {\Bbb C}P^2 ~|~ y^2z=4x^3-g_2xz^2-g_3z^3\}$,
where $g_2$ and  $g_3$  are some constant complex numbers.   
The {\it $j$-invariant}  of ${\cal E}({\Bbb C})$ is the complex number
\displaymath
j({\cal E}({\Bbb C}))={1728 g_2^3\over g_2^3-27g_3^2},
\enddisplaymath
which is constant only on isomorphic elliptic curves.  The Weierstrass function $\wp(z)$
defines an isomorphism ${\cal E}({\Bbb C})\cong {\Bbb C}/({\Bbb Z}+{\Bbb Z}\tau)$
between elliptic curves and complex tori of modulus $\tau\in {\Bbb H}:=\{z=x+iy\in {\Bbb C}~|~y>0\}$,
see Theorem \ref{thm5.1.1};  by ${\cal E}_{\tau}$ we understand an elliptic curve of 
complex modulus $\tau$. 
\begin{dfn}
By an isogeny between elliptic curves ${\cal E}_{\tau}$ and ${\cal E}_{\tau'}$
one understands an analytic map  $\varphi:  {\cal E}_{\tau}\to {\cal E}_{\tau'}$,
such that $\varphi(0)=0$.  Clearly,  the invertible isogeny corresponds to an  
isomorphism between elliptic curves. 
\end{dfn}
 \index{isogeny}
\begin{rmk}
\textnormal{
The elliptic curves ${\cal E}_{\tau}$ and  ${\cal E}_{\tau'}$ are isogenous  
 if and only if 
\displaymath
\tau'={a\tau+b\over c\tau+d} \quad \hbox{for some matrix} 
 \quad\left(\matrix{a & b\cr c & d}\right)  \in M_2({\Bbb Z})
\quad \hbox{with} \quad ad-bc>0.
 \enddisplaymath
 The case of  an invertible matrix  (i.e. $ad-bc=1$)   corresponds to 
 an isomorphism between elliptic curves.  
 (We leave the proof  to the reader.  Hint:  notice that $z\mapsto\alpha z$
is an invertible holomorphic map for each $\alpha\in {\Bbb C}-\{0\}$.)
}
\end{rmk}
 An {\it endomorphism}  of ${\cal E}_{\tau}$ is a multiplication 
of the lattice $L_{\tau}:={\Bbb Z}+{\Bbb Z}\tau\subset {\Bbb C}$ by complex number $z$
such that 
\displaymath
zL_{\tau}\subseteq L_{\tau}. 
\enddisplaymath
In other words,  the endomorphism is an isogeny of the elliptic curve into itself.   
The sum and product of two endomorphisms is an endomorphism of ${\cal E}_{\tau}$;
thus one gets  a commutative ring of all endomorphisms of ${\cal E}_{\tau}$ 
denoted by $End~({\cal E}_{\tau})$.  Typically  $End~({\cal E}_{\tau})\cong {\Bbb Z}$,
i.e. the only endomorphisms of ${\cal E}_{\tau}$ are the multiplication-by-$m$ 
endomorphisms;  however,   for a countable set of $\tau$
\displaymath
End~({\cal E}_{\tau})\cong {\Bbb Z}+ fO_k,
\enddisplaymath
where $k={\Bbb Q}(\sqrt{-D})$ is an imaginary quadratic field,  $O_k$ its ring
of integers and $f\ge 1$ is the conductor of a finite index subring of $O_k$.
(The proof of this simple but fundamental fact is left to the reader.)  It is easy 
to see that in such a case $\tau\in End~({\cal E}_{\tau})$,  i.e. complex modulus 
itself is an imaginary quadratic number.   
\begin{dfn}
Elliptic curve ${\cal E}_{\tau}$ is said to have complex multiplication if 
$End~({\cal E}_{\tau})\cong {\Bbb Z}+ fO_k$,  i.e. $\tau$ is an imaginary
quadratic number;  such a curve will be denoted by ${\cal E}_{CM}^{(-D, f)}$.
\end{dfn}
\begin{rmk}
\textnormal{
There is a finite number of pairwise non-isomorphic elliptic curves
with the same ring of non-trivial endomorphisms  $R:=End~({\cal E}_{\tau})$;
such a number is equal to $|Cl~(R)|$, where $Cl~(R)$ is the class group of
ring $R$.   This fact is extremely important,  because the $j$-invariant
$j({\cal E}_{CM}^{(-D,f)})$ is known to be an algebraic number and, 
therefore,  $Gal~(K|k)\cong Cl~(R)$,
where $K=k(j({\cal E}_{CM}^{(-D,f)}))$ and $Gal~(K|k)$ is the Galois 
group of the field extension $K|k$.  In other words, the number field $K$
is the {\it Hilbert class field}  of   imaginary quadratic field $k$.   
Moreover,
\displaymath
{\cal E}_{CM}^{(-D,f)}\cong {\cal E}(K),
\enddisplaymath
i.e. the complex constants $g_2$ and $g_3$ in the cubic equation for
${\cal E}_{CM}^{(-D,f)}$ must belong to the number field $K$.    
}
\end{rmk}

\subsection{Functor on elliptic curves with complex multiplication}
\begin{dfn}
 By {\bf Ell-Isgn}  we shall mean the category of all elliptic curves   ${\cal E}_{\tau}$;
 the arrows of {\bf Ell-Isgn}  are identified with the isogenies  between 
 elliptic curves ${\cal E}_{\tau}$.  We shall write {\bf NC-Tor-Homo} to denote the
  category of all noncommutative tori  ${\cal A}_{\theta}$;
 the arrows of {\bf NC-Tor-Homo}  are identified with the stable homomorphisms
  between  noncommutative tori  ${\cal A}_{\theta}$.  
 \end{dfn}
\begin{rmk}\label{rmk6.1.3}
\textnormal{
The noncommutative tori ${\cal A}_{\theta}$ and  ${\cal A}_{\theta'}$ are stably
homomorphic   if and only if 
\displaymath
\theta'={a\theta+b\over c\theta+d} \quad \hbox{for some matrix} 
 \quad\left(\matrix{a & b\cr c & d}\right)  \in M_2({\Bbb Z})
\quad \hbox{with} \quad ad-bc>0.
 \enddisplaymath
 The case of  an invertible matrix  (i.e. $ad-bc=1$)   corresponds to 
 a stable isomorphism (Morita equivalence) between noncommutative tori.  
 (We leave the proof  to the reader.  Hint:  follow and modify the argument of 
 [Rieffel 1981] \cite{Rie2}.)
}
\end{rmk}
\begin{figure}[here]
\begin{picture}(300,110)(-120,-5)
\put(20,70){\vector(0,-1){35}}
\put(130,70){\vector(0,-1){35}}
\put(45,23){\vector(1,0){60}}
\put(45,83){\vector(1,0){60}}
\put(15,20){${\cal A}_{\theta}$}
\put(0,50){$F$}
\put(145,50){$F$}
\put(123,20){${\cal A}_{\theta'={a\theta+b\over c\theta+d}}$}
\put(17,80){${\cal E}_{\tau}$}
\put(122,80){${\cal E}_{\tau'={a\tau+b\over c\tau+d}}$}
\put(60,30){\sf stably}
\put(45,10){\sf  homomorphic}
\put(50,90){\sf isogenous}
\end{picture}
\caption{Functor on isogenous elliptic curves.}
\end{figure}
\begin{thm}\label{thm6.1.1}
{\bf (Functor on isogenous elliptic curves)}
There exists a covariant functor 
\displaymath
F:  \hbox{{\bf Ell-Isgn}} \longrightarrow \hbox{{\bf NC-Tor-Homo}},
\enddisplaymath
which maps isogenous  elliptic curves ${\cal E}_{\tau}$ to the 
stably homomorphic  noncommutative tori ${\cal A}_{\theta}$,
see Fig. 6.1;    the functor $F$ is non-injective and $Ker~F\cong (0,\infty)$.
In particular, $F$ maps  isomorphic elliptic curves  to the stably isomorphic 
(Morita equivalent) noncommutative tori. 
 \end{thm}
\begin{thm}\label{thm6.1.2}
{\bf (Functor on  elliptic curves with complex multiplication)}
If  $Isom~({\cal E}_{CM}):=\{{\cal E}_{\tau}\in \hbox{{\bf Ell-Isgn}}~|~{\cal E}_{\tau}\cong {\cal E}_{CM}\}$
is the isomorphism class of an elliptic curve with complex multiplication and 
${\goth m}_{CM}:=\mu_{CM}({\Bbb Z}+{\Bbb Z}\theta_{CM})\subset {\Bbb R}$
is a ${\Bbb Z}$-module such that ${\cal A}_{\theta_{CM}}=F({\cal E}_{CM})$
and $\mu_{CM}\in Ker~F$,   then:

\medskip
(i) ${\goth m}_{CM}$ is an invariant of  $Isom~({\cal E}_{CM})$;

\smallskip
(ii) ${\goth m}_{CM}$ is a full module in the  real quadratic number field.

\medskip\noindent
In particular, ${\cal A}_{\theta_{CM}}$ is a noncommutative torus with real multiplication.
\end{thm}
\begin{dfn}
If  ${\cal A}_{RM}^{(D,f)}$ is a noncommutative torus with real multiplication,
then the Riemann surface $X({\cal A}_{RM}^{(D,f)})$ is called associated
to ${\cal A}_{RM}^{(D,f)}$ whenever the covering of geodesic spectrum 
of $X({\cal A}_{RM}^{(D,f)})$ on the half-plane ${\Bbb H}$ contains the set
 $\{\tilde\gamma (x,\bar x) ~: ~\forall x\in {\goth m}_{CM}\}$,  where
\displaymath
\tilde \gamma (x,\bar x)={xe^{t\over 2}+i\bar x e^{-{t\over 2}}\over e^{t\over 2}+ie^{-{t\over 2}}},
\qquad -\infty\le t\le \infty
\enddisplaymath
is the geodesic half-circle through the pair of conjugate quadratic irrationalities
$x,\bar x\in {\goth m}_{CM}\subset\partial {\Bbb H}$,   see Definition \ref{dfn6.1.5}.  
\end{dfn}
\begin{thm}\label{thm6.1.3}
{\bf (Functor on noncommutative tori with real multiplication)}
For each  square-free  integer $D>1$ and integer $f\ge 1$  there exists a holomorphic 
map $F^{-1}: X({\cal A}_{RM}^{(D,f)})\to {\cal E}_{CM}^{(-D,f)}$,  where
$F({\cal E}_{CM}^{(-D,f)})={\cal A}_{RM}^{(D,f)}$.
\end{thm}
 \index{real multiplication}
\begin{rmk}
\textnormal{
Roughly speaking,  Theorem \ref{thm6.1.3}  is an explicit form of functor $F$ 
constructed in Theorem \ref{thm6.1.2};
moreover,  Theorem \ref{thm6.1.3}  says the $F$ is a bijection by constructing  an
explicit inverse functor $F^{-1}$. 
}
\end{rmk}

\subsection{Proof of Theorem \ref{thm6.1.1}}
The proof is a modification of the one for Theorem \ref{thm5.1.2};   we freely use the notation 
and facts of the Teichm\"uller theory  introduced  in  Section 5.1.2. 
Let $\phi=Re~\omega$ be a 1-form defined by a holomorphic form $\omega$ on the 
complex torus $S$.   Since $\omega$
is holomorphic, $\phi$ is a closed $1$-form on topological torus $T^2$. The ${\Bbb R}$-isomorphism
$h_q: H^0(S,\Omega)\to Hom~(H_1(T^2); {\Bbb R})$, as explained,
is given by the formulas:
\displaymath
\left\{
\begin{array}{cc}
\lambda_1 &= \int_{\gamma_1}\phi\\
\lambda_2 &= \int_{\gamma_2}\phi,
\end{array}
\right.
\enddisplaymath
where $\{\gamma_1,\gamma_2\}$ is a basis in the first homology group of $T^2$. 
We further assume that, after a proper choice of the basis,  $\lambda_1,\lambda_2$ are positive real
numbers.  Denote by $\Phi_{T^2}$ the space of measured foliations on $T^2$.
Each ${\cal F}\in \Phi_{T^2}$ is measure equivalent to a foliation
by a family of the parallel lines of a slope $\theta$ and the invariant  transverse
measure $\mu$,  see Fig. 6.2.

\begin{figure}[here]
\begin{picture}(300,60)(-30,0)

\put(130,10){\line(1,0){40}}
\put(130,10){\line(0,1){40}}
\put(130,50){\line(1,0){40}}
\put(170,10){\line(0,1){40}}

\put(130,40){\line(2,1){20}}
\put(130,30){\line(2,1){40}}
\put(130,20){\line(2,1){40}}
\put(130,10){\line(2,1){40}}

\put(150,10){\line(2,1){20}}

\end{picture}

\caption{Measured foliation ${\cal F}$ on $T^2={\Bbb R}^2/{\Bbb Z}^2$.}
\end{figure}

\medskip\noindent
We use the notation ${\cal F}^{\mu}_{\theta}$ for such a foliation. 
There exists a simple  relationship between the reals  $(\lambda_1,\lambda_2)$
and $(\theta,\mu)$. Indeed, the closed $1$-form $\phi=Const$ defines a measured foliation,
${\cal F}^{\mu}_{\theta}$, so that
\displaymath
\left\{
\begin{array}{ccc}
\lambda_1  &= \int_{\gamma_1}\phi &= \int_0^1\mu dx\\
\lambda_2  &= \int_{\gamma_2}\phi &= \int_0^1\mu dy
\end{array}
\right.
\hbox{,  where}
\quad
{dy\over dx}=\theta.
\enddisplaymath
By the integration:
\displaymath
\left\{
\begin{array}{ccc}
\lambda_1  &= \int_0^1\mu dx &= \mu\\
\lambda_2  &= \int_0^1\mu\theta  dx &= \mu\theta.
\end{array}
\right.
\enddisplaymath
Thus, one gets $\mu=\lambda_1$ and $\theta={\lambda_2\over\lambda_1}$.  
Recall that the Hubbard-Masur theory  establishes a homeomorphism
$h: T_S(1)\to \Phi_{T^2}$, where $T_S(1)\cong {\Bbb H}=\{\tau: Im~\tau>0\}$
is the Teichm\"uller space of the  torus,   see Corollary \ref{cor5.1.1}.  
Denote by $\omega_N$ an invariant
(N\'eron) differential of the complex torus ${\Bbb C}/(\omega_1{\Bbb Z}+\omega_2{\Bbb Z})$.
It is well known that $\omega_1=\int_{\gamma_1}\omega_N$ and 
$\omega_2=\int_{\gamma_2}\omega_N$, where $\gamma_1$ and $\gamma_2$ are the meridians of the torus.
Let $\pi$ be a projection acting by the formula $(\theta,\mu)\mapsto \theta$. 
An explicit formula for the functor  $F:$  {\bf Ell-Isgn} $\to$ {\bf NC-Tor-Homo}  is given by the composition
$F=\pi\circ h$,   where $h$ is the Hubbard-Masur homeomorphism. 
In other words, one gets the following explicit correspondence between the complex 
and noncommutative tori: 
\displaymath
{\cal E}_{\tau}={\cal E}_{(\int_{\gamma_2}\omega_N) / (\int_{\gamma_1}\omega_N)}
\buildrel\rm h\over
\longmapsto
{\cal F}^{\int_{\gamma_1}\phi}_{(\int_{\gamma_2}\phi)/(\int_{\gamma_1}\phi)}
\buildrel\rm\pi \over
\longmapsto
{\cal A}_{(\int_{\gamma_2}\phi)/(\int_{\gamma_1}\phi)}= {\cal A}_{\theta},
\enddisplaymath
where ${\cal E}_{\tau}={\Bbb C}/({\Bbb Z}+{\Bbb Z}\tau)$.   Let 
\displaymath
\varphi: {\cal E}_{\tau}\longrightarrow  {\cal E}_{\tau'}
\enddisplaymath
be an isogeny of the elliptic curves.  The action of $\varphi$ on the homology basis $\{\gamma_1,\gamma_2\}$
of $T^2$ is given by the formulas
\displaymath
\left\{
\begin{array}{cc}
\gamma_1' &= a\gamma_1+b\gamma_2\nonumber\\
\gamma_2' &= c\gamma_1+d\gamma_2
\end{array}
\right.
\hbox{,  where}
\left(\matrix{a & b\cr c & d}\small\right)\in M_2({\Bbb Z}). 
\enddisplaymath
Recall that the functor  $F:$  {\bf Ell-Isgn} $\to$ {\bf NC-Tor-Homo}
 is given by the formula
\displaymath
\tau={\int_{\gamma_2}\omega_N\over\int_{\gamma_1}\omega_N}
\longmapsto
\theta={\int_{\gamma_2}\phi\over\int_{\gamma_1}\phi},
\enddisplaymath
where $\omega_N$ is an invariant differential on ${\cal E}_{\tau}$
and $\phi=Re~\omega$ is a closed 1-form on $T^2$.

\bigskip
(i) From the left-hand side of the above equation,  one obtains 
\displaymath
\left\{
\begin{array}{ccccc}
\omega_1' &= \int_{\gamma_1'}\omega_N &=  \int_{a\gamma_1+b\gamma_2}\omega_N  &= 
a\int_{\gamma_1}\omega_N +b\int_{\gamma_2}\omega_N &= a\omega_1+b\omega_2\nonumber\\
\omega_2' &= \int_{\gamma_2'}\omega_N &=  \int_{c\gamma_1+d\gamma_2}\omega_N  &= 
c\int_{\gamma_1}\omega_N +d\int_{\gamma_2}\omega_N &= c\omega_1+d\omega_2,
\end{array}
\right.
\enddisplaymath
and therefore $\tau'={\int_{\gamma_2'}\omega_N\over\int_{\gamma_1'}\omega_N}=
{c+d\tau\over a+b\tau}$.

\bigskip
(ii) From the right-hand side, one obtains
\displaymath
\left\{
\begin{array}{ccccc}
\lambda_1' &= \int_{\gamma_1'}\phi &=  \int_{a\gamma_1+b\gamma_2}\phi  &= 
a\int_{\gamma_1}\phi +b\int_{\gamma_2}\phi &= a\lambda_1+b\lambda_2\nonumber\\
\lambda_2' &= \int_{\gamma_2'}\phi &=  \int_{c\gamma_1+d\gamma_2}\phi  &= 
c\int_{\gamma_1}\phi +d\int_{\gamma_2}\phi &= c\lambda_1+d\lambda_2,
\end{array}
\right.
\enddisplaymath
and therefore $\theta'={\int_{\gamma_2'}\phi\over\int_{\gamma_1'}\phi}=
{c+d\theta\over a+b\theta}$.   Comparing (i) and (ii),  one gets the conclusion
of the first part of Theorem  \ref{thm6.1.1}.   To prove the second part,  recall that
the invertible isogeny is an isomorphism of the elliptic curves.
In this case  $\small\left(\matrix{a & b\cr c & d}\small\right)\in SL_2({\Bbb Z})$
and $\theta'=\theta ~mod~SL_2({\Bbb Z})$.  Therefore $F$ sends the isomorphic
elliptic curves to the stably isomorphic noncommutative tori. The second part 
of Theorem  \ref{thm6.1.1} is proved.    It follows from the proof that 
$F:$  {\bf Ell-Isgn} $\to$ {\bf NC-Tor-Homo}  is a covariant
functor.   Indeed, $F$ preserves the morphisms and does not reverse the arrows:
$F(\varphi_1\varphi_2)=\varphi_1\varphi_2=F(\varphi_1)F(\varphi_2)$
for any pair of the isogenies $\varphi_1,\varphi_2\in Mor~(\hbox{{\bf Ell-Isgn}})$. 
Theorem \ref{thm6.1.1} follows.
$\square$

\subsection{Proof of Theorem \ref{thm6.1.2}}
\begin{lem}\label{lem6.1.1}
Let ${\goth m}\subset {\Bbb R}$ be a module of the rank 2, i.e
${\goth m}={\Bbb Z}\lambda_1+{\Bbb Z}\lambda_2$, where
$\theta={\lambda_2\over\lambda_1}\not\in {\Bbb Q}$. 
If ${\goth m}'\subseteq {\goth m}$ is a submodule of the rank 2,
then ${\goth m}'=k {\goth m}$, where either:

\medskip
(i) $k\in {\Bbb Z}-\{0\}$ and $\theta\in {\Bbb R}-{\Bbb Q}$, or 

\smallskip
(ii) $k$ and $\theta$ are the irrational numbers of  
a quadratic number field. 
\end{lem}
{\it Proof.} Any rank 2 submodule of $m$ can be written as  
${\goth m}'=\lambda_1'{\Bbb Z}+{\lambda_2'}{\Bbb Z}$, where
\displaymath
\left\{
\begin{array}{cc}
\lambda_1'  &= a\lambda_1 +b\lambda_2\\
\lambda_2'  &= c\lambda_1 +d\lambda_2
\end{array}
\right.
\qquad \hbox{and} \qquad 
\small\left(\matrix{a & b\cr c & d}\small\right)\in M_2({\Bbb Z}).
\enddisplaymath

\bigskip
(i) Let us assume that $b\ne 0$.  Let $\Delta= (a+d)^2-4(ad-bc)$ and $\Delta'=(a+d)^2-4bc$. 
We shall consider the following cases.

\medskip
{\bf Case 1: $\Delta>0$ and $\Delta\ne m^2$, $m\in {\Bbb Z}-\{0\}$.} 
The real number  $k$ can be determined  from the equations: 
\displaymath
\left\{
\begin{array}{ccc}
\lambda_1'  &= k\lambda_1 &= a\lambda_1 +b\lambda_2\\
\lambda_2'  &= k\lambda_2 &= c\lambda_1 +d\lambda_2.
\end{array}
\right.
\enddisplaymath
Since  $\theta={\lambda_2\over\lambda_1}$, one gets the  equation $\theta={c+d\theta\over a+b\theta}$
by taking the ratio of  two  equations above.  A quadratic equation for $\theta$ writes
as $b\theta^2+(a-d)\theta-c=0$. The discriminant of the equation coincides with $\Delta$
and therefore there exist  real roots $\theta_{1,2}={d-a\pm\sqrt{\Delta}\over 2b}$.
Moreover, $k=a+b\theta={1\over 2}(a+d\pm\sqrt{\Delta})$.   Since $\Delta$ is
not the square of an integer, $k$ and $\theta$ are  irrationalities   of the quadratic number field
${\Bbb Q}(\sqrt{\Delta})$.

\medskip
{\bf Case 2: $\Delta>0$ and $\Delta=m^2$, $m\in {\Bbb Z}-\{0\}$.}
Note that $\theta={a-d\pm |m|\over 2c}$ is a rational number. 
Since $\theta$ does not satisfy the rank assumption of the lemma,
the case should be omitted.

\medskip
{\bf Case 3: $\Delta=0$.}
The quadratic equation has a double root $\theta={a-d\over 2c}\in {\Bbb Q}$.
This case leads to a module of the rank $1$, which is contrary to an assumption of the
lemma.

\medskip
{\bf Case 4: $\Delta<0$ and $\Delta'\ne m^2$, $m\in {\Bbb Z}-\{0\}$.}
Let us  define  a new basis
$\{\lambda_1'',\lambda_2''\}$ in ${\goth m}'$ so that 
\displaymath
\left\{
\begin{array}{cc}
\lambda_1''  &=  \lambda_1'\\
\lambda_2''  &= -\lambda_2'.
\end{array}
\right.
\enddisplaymath
Then:
\displaymath
\left\{
\begin{array}{cc}
\lambda_1''  &= a\lambda_1 +b\lambda_2\\
\lambda_2''  &= -c\lambda_1 -d\lambda_2,
\end{array}
\right.
\enddisplaymath
and $\theta={\lambda_2''\over\lambda_1''}={-c-d\theta\over a+b\theta}$.
The quadratic equation for $\theta$ has the form $b\theta^2+(a+d)\theta+c=0$,
whose discrimimant is $\Delta'=(a+d)^2-4bc$.
Let us show that $\Delta'>0$. Indeed, $\Delta= (a+d)^2-4(ad-bc)<0$ 
and the evident inequality $-(a-d)^2\le 0$ have the same sign, and
we shall add them up.  After an obvious elimination, one gets $bc<0$.
Therefore $\Delta'$ is a  sum of the two positive integers, which is   
always a positive integer. 
Thus,  there exist the real roots $\theta_{1,2}={-a-d\pm\sqrt{\Delta'}\over 2b}$.
Moreover, $k=a+b\theta={1\over 2}(a-d\pm\sqrt{\Delta'})$. Since $\Delta'$ is
not the square of an integer, $k$ and $\theta$ are the irrational numbers in the 
quadratic field ${\Bbb Q}(\sqrt{\Delta'})$.

\medskip
{\bf Case 5: $\Delta<0$ and $\Delta'=m^2$, $m\in {\Bbb Z}-\{0\}$.}
Note that $\theta={-a-d\pm |m|\over 2b}$ is a rational number.
Since $\theta$ does not satisfy the rank assumption of the lemma,
the case should be omitted.

\bigskip
(ii) Assume that $b=0$.

\medskip
{\bf Case 1: $a-d\ne 0$.}
The quadratic equation for $\theta$ degenerates to a 
linear eauation $(a-d)\theta+c=0$. The root $\theta={c\over d-a}\in {\Bbb Q}$
does not satisfy the rank assumption again, and we omit the case.

\medskip
{\bf Case 2: $a=d$ and $c\ne 0$.}
It is easy to see, that the set of the solutions for $\theta$ is an empty set.

\medskip
{\bf Case 3: $a=d$ and $c=0$.}
Finally, in this case all coefficients of the quadratic equation vanish,
so that any $\theta\in {\Bbb R}-{\Bbb Q}$ is a solution. Note that 
$k=a=d\in {\Bbb Z}$.  Thus, one gets case (i) of the lemma. 
Since there are no other possiblities left,  Lemma \ref{lem6.1.1}  is proved.
$\square$

\bigskip   
\begin{lem}\label{lem6.1.2}
Let ${\cal E}_{CM}$ be an elliptic curve with  complex multiplication
and  consider a ${\Bbb Z}$-module  $F(Isom~({\cal E}_{CM}))=\mu_{CM}({\Bbb Z}+{\Bbb Z}\theta_{CM}):= {\goth m}_{CM}$.
Then:

\medskip
(i) $\theta_{CM}$ is a quadratic irrationality, 

\smallskip
(ii) $\mu_{CM}\in {\Bbb Q}$ (up to a choice of map $F$).  
\end{lem}
{\it Proof.} 
(i) Since ${\cal E}_{CM}$  has  complex multiplication,  one gets $End~({\cal E}_{CM})> {\Bbb Z}$.
In particular,  there exists a non-trivial  isogeny  
\displaymath
\varphi:  ~{\cal E}_{CM}\longrightarrow {\cal E}_{CM}, 
\enddisplaymath
i.e an endomorphism which is {\it not}   the multiplication by $k\in {\Bbb Z}$.  
By Theorem  \ref{thm6.1.1} and Remark \ref{rmk6.1.3},   each isogeny $\varphi$ defines a rank $2$ submodule ${\goth m}'$
 of module  ${\goth m}_{CM}$.    By Lemma \ref{lem6.1.1},  ${\goth m}'=k {\goth m}_{CM}$  for a
$k\in {\Bbb R}$.  Because  $\varphi$ is a non-trivial endomorphism,  we get  $k\not\in {\Bbb Z}$;
thus,   option (i) of Lemma \ref{lem6.1.1}  is excluded.  Therefore,   
by the item (ii) of  Lemma \ref{lem6.1.1},  real number $\theta_{CM}$ must be  a quadratic irrationality.

\smallskip
(ii) Recall that $E_{\cal F}\subset Q- \{0\}$ is the space of holomorphic
differentials on the complex torus, whose  horizontal trajectory
structure is equivalent to given  measured foliation ${\cal F}={\cal F}^{\mu}_{\theta}$. 
We shall vary ${\cal F}_{\theta}^{\mu}$, thus varying the Hubbard-Masur homeomorphism
$h=h({\cal F}^{\mu}_{\theta}): E_{\cal F}\to T(1)$,  see Section 6.1.2.  Namely, consider a 1-parameter
continuous family of such maps  $h=h_{\mu}$, where $\theta=Const$ and $\mu\in {\Bbb R}$.  
Recall that $\mu_{CM}=\lambda_1=\int_{\gamma_1}\phi$, where $\phi=Re~\omega$ and $\omega\in E_{\cal F}$. 
The family $h_{\mu}$ generates  a family $\omega_{\mu}=h^{-1}_{\mu}(C)$, where $C$ is a fixed
point in $T(1)$. Denote by $\phi_{\mu}$ and $\lambda_1^{\mu}$ the corresponding families
of the closed 1-forms and their periods, respectively. By the continuity, $\lambda_1^{\mu}$
takes on a rational value for a $\mu=\mu'$.  (Actually, every  neighborhood of $\mu_0$
contains such a $\mu'$.) Thus, $\mu_{CM}\in {\Bbb Q}$ for the Hubbard-Masur homeomorphism $h=h_{\mu'}$.  
Lemma \ref{lem6.1.2} follows. 
$\square$

\bigskip\noindent
The claim (ii) of Theorem \ref{thm6.1.2}  follows from (i) of  Lemma \ref{lem6.1.2} and claim (i)
of  Theorem \ref{thm6.1.2}.  To prove claim (i) of  Theorem \ref{thm6.1.2},  notice that  whenever   
${\cal E}_1, {\cal E}_2\in Isom ({\cal E}_{CM})$  the  respective ${\Bbb Z}$-modules coincide, 
i.e.  ${\goth m}_1={\goth m}_2$;  this happens because  an isomorphism  between  elliptic curves  
corresponds to a change of basis in the module ${\goth m}$,  see Theorem  \ref{thm6.1.1} 
and Remark \ref{rmk6.1.3}.   Theorem \ref{thm6.1.2}  is proved.
$\square$

\subsection{Proof of Theorem \ref{thm6.1.3}}
Let us  recall some classical facts and notation,  and give an exact definition
of the Riemann surface    $X({\cal A}_{RM}^{(D,f)})$.
Let $N\ge 1$ be an integer;  recall that $\Gamma_1(N)$ is a subgroup 
of the modular group $SL_2({\Bbb Z})$ consisting of matrices of the form
\displaymath
\left\{\left(\matrix{a & b\cr c & d}\right)\in SL_2({\Bbb Z}) ~| ~ a,d\equiv~1~mod~N, 
~c\equiv~0~mod~N\right\};
\enddisplaymath
the corresponding Riemann surface ${\Bbb H}/\Gamma_1(N)$ will be denoted by $X_1(N)$. 
Consider the geodesic spectrum of $X_1(N)$,  i.e the set $Spec~X_1(N)$
consisting of all closed geodesics of the surface $X_1(N)$;   each geodesic $\gamma\in Spec~X_1(N)$
is the image under the covering map ${\Bbb H}\to {\Bbb H}/\Gamma_1(N)$  of a geodesic 
half-circle $\tilde\gamma\in {\Bbb H}$ passing through the
  points $x$ and $\bar x$  fixed by the linear fractional transformation
$x\mapsto {ax+b\over cx+d}$, where matrix  $(a, b, c, d)\in \Gamma_1(N)$.  
It is not hard to see,  that $x$ and $\bar x$ are quadratic irrational numbers;
the numbers are real when $|a+d|>2$.      
\begin{dfn}\label{dfn6.1.5}
We shall say that the Riemann surface $X$ is associated to the noncommutative
torus  ${\cal A}_{RM}^{(D,f)}$, if $\{\tilde\gamma (x,\bar x) ~: ~\forall x\in {\goth m}_{RM}^{(D,f)}\}
\subset \widetilde{Spec}~X$,  where  $\widetilde{Spec}~X\subset {\Bbb H}$ is the set of 
geodesic half-circles covering the geodesic spectrum of $X$ and  ${\goth m}_{RM}^{(D,f)}$
is a ${\Bbb Z}$-module (a pseudo-lattice) in ${\Bbb R}$ generated by torus ${\cal A}_{RM}^{(D,f)}$; 
 the associated  Riemann surface will be denoted by  $X({\cal A}_{RM}^{(D,f)})$.
 \end{dfn}
 \index{geodesic spectrum}
\begin{lem}\label{lem6.1.3}
$X({\cal A}_{RM}^{(D,f)})\cong X_1(fD)$.
\end{lem}
{\it Proof.}
Recall that ${\goth m}_{RM}^{(D,f)}$ is a ${\Bbb Z}$-module  (a pseudo-lattice) 
with  real  multiplication by an order $R$ in the real quadratic number field
${\Bbb Q}(\sqrt{D})$;  it is known, that ${\goth m}_{RM}^{(D,f)}\subseteq R$
and $R={\Bbb Z}+(f\omega){\Bbb Z}$, where $f\ge 1$ is the conductor of $R$ and  
\displaymath
\omega=\cases{{1+\sqrt{D}\over 2} & if $D\equiv 1 ~mod~4$,\cr
               \sqrt{D} & if $D\equiv 2,3 ~mod~4$,}
\enddisplaymath
see e.g. [Borevich \&  Shafarevich  1988]   \cite{BS},  pp. 130-131.  
Recall that matrix $(a,b,c,d)\in SL_2({\Bbb Z})$ has a pair
of real fixed points $x$ and $\bar x$ if and only if $|a+d|>2$ (the hyperbolic matrix);
the fixed points can be found from the equation $x=(ax+b)(cx+d)^{-1}$  by the
formulas
\displaymath
x={a-d\over 2c}+\sqrt{{(a+d)^2-4\over 4c^2}}, \qquad
\bar x={a-d\over 2c}-\sqrt{{(a+d)^2-4\over 4c^2}}. 
\enddisplaymath

\bigskip
{\sf Case I.} If $D\equiv 1~mod~4$, then the above formulas  imply  that
$R=(1+{f\over 2}){\Bbb Z}+{\sqrt{f^2D}\over 2}{\Bbb Z}$. If $x\in {\goth m}_{RM}^{(D,f)}$
is fixed point of a transformation $(a,b,c,d)\in SL_2({\Bbb Z})$,   then
\displaymath
\left\{
\begin{array}{ccc}
{a-d\over 2c} &=& (1+{f\over 2})z_1\\
{(a+d)^2-4\over 4c^2}  &=& {f^2D\over 4}z_2^2
\end{array}
\right.
\enddisplaymath
for some integer numbers $z_1$ and $z_2$.  
The second equation can be written in the form $(a+d)^2-4=c^2f^2Dz_2^2$;
we have therefore $(a+d)^2\equiv 4~mod~(fD)$ and $a+d\equiv\pm 2~mod~(fD)$.
Without loss of generality we assume $a+d\equiv 2~mod~(fD)$ since
matrix $(a,b,c,d)\in SL_2({\Bbb Z})$ can be multiplied by $-1$.  Notice that
the last equation admits a solution $a=d\equiv 1~mod~(fD)$.
The first equation yields us ${a-d\over c}=(2+f)z_1$, where $c\ne0$ since
the matrix $(a,b,c,d)$ is hyperbolic.  Notice that  $a-d\equiv 0~mod~(fD)$;
since the ratio ${a-d\over c}$ must be integer, we conclude that $c\equiv 0~mod~(fD)$.  
Summing up,   we get: 
\displaymath
a\equiv 1~mod~(fD), \quad d\equiv  1~mod~(fD), \quad  c\equiv 0~mod~(fD).
\enddisplaymath

\bigskip
{\sf Case II.}   If $D\equiv 2$ or $3~mod~4$, then
$R={\Bbb Z}+(\sqrt{f^2D})~{\Bbb Z}$.  If $x\in{\goth m}_{RM}^{(D,f)}$
is  fixed point of a  transformation $(a,b,c,d)\in SL_2({\Bbb Z})$, then
\displaymath
\left\{
\begin{array}{ccc}
{a-d\over 2c} &=& z_1\\
{(a+d)^2-4\over 4c^2}  &=& f^2Dz_2^2
\end{array}
\right.
\enddisplaymath
for some integer numbers $z_1$ and $z_2$.   The second equation
gives  $(a+d)^2-4=4c^2f^2Dz_2^2$;  therefore  $(a+d)^2\equiv 4~mod~(fD)$ and $a+d\equiv\pm 2~mod~(fD)$.
Again without loss of generality we assume $a+d\equiv 2~mod~(fD)$ since
matrix $(a,b,c,d)\in SL_2({\Bbb Z})$ can be  multiplied by $-1$.  
The last equation admits a solution $a=d\equiv 1~mod~(fD)$.
The first equation is ${a-d\over c}=2z_1$, where $c\ne0$.
Since  $a-d\equiv 0~mod~(fD)$ and the ratio ${a-d\over c}$ must be integer,
one concludes  that $c\equiv 0~mod~(fD)$.
 All together,  one  gets
\displaymath
a\equiv 1~mod~(fD), \quad d\equiv  1~mod~(fD), \quad  c\equiv 0~mod~(fD).
\enddisplaymath
 Since all possible cases are exhausted,  Lemma \ref{lem6.1.3}  follows.
$\square$

\begin{rmk}
\textnormal{
There exist other finite index subgroups of $SL_2({\Bbb Z})$ whose
geodesic spectrum contains the set 
 $\{\tilde\gamma (x,\bar x) ~: ~\forall x\in {\goth m}_{RM}^{(D,f)}\}$;
however $\Gamma_1(fD)$ is a unique group with such a property 
among subgroups of the principal congruence group.
}
\end{rmk}
\begin{rmk}
\textnormal{
Not all geodesics of $X_1(fD)$  have  the above form;  thus the set \linebreak
 $\{\tilde\gamma (x,\bar x) ~: ~\forall x\in {\goth m}_{RM}^{(D,f)}\}$ is strictly
included in the geodesic spectrum of modular curve $X_1(fD)$.
}
\end{rmk}
\begin{dfn}
The group
\displaymath
\Gamma(N):=\left\{\left(\matrix{a & b\cr c & d}\right)\in SL_2({\Bbb Z})~|~a,d\equiv~1~mod~N,~b,c\equiv~0~mod~N\right\}
\enddisplaymath
is called a  principal congruence group  of level $N$;  the corresponding compact  
modular curve will be denoted by $X(N)={\Bbb H}/\Gamma(N)$. 
\end{dfn}
\begin{lem}\label{lem6.1.4}
 {\bf (Hecke)}
There exists a holomorphic map $X(fD)\to {\cal E}_{CM}^{(-D,f)}.$ 
\end{lem}
{\it Proof.}  A detailed proof of this beautiful fact is given in
[Hecke  1928]  \cite{Hec1}.   For the sake of clarity,  we shall  
give an idea of the proof.
Let  ${\goth R}$ be an order of conductor $f\ge 1$
in the imaginary quadratic number field ${\Bbb Q}(\sqrt{-D})$;  consider
an $L$-function attached to ${\goth R}$
\displaymath
L(s, \psi)=\prod_{{\goth P}\subset {\goth R}}{1\over 1-{\psi({\goth P})\over N({\goth P})^s}},
\quad s\in {\Bbb C}, 
\enddisplaymath
where ${\goth P}$ is a prime ideal in ${\goth R}$, $N({\goth P})$ its norm
and $\psi$ a  Gr\"ossencharacter. 
A crucial observation of Hecke  says that the series $L(s, \psi)$ converges to a cusp form $w(s)$
of the principal congruence group $\Gamma(fD)$. 
By the Deuring Theorem,   $L({\cal E}_{CM}^{(-D,f)},s)=L(s,\psi)L(s, \bar\psi)$, where $L({\cal E}_{CM}^{(-D,f)},s)$ is the
Hasse-Weil $L$-function of the elliptic curve and $\bar\psi$ a conjugate of the 
Gr\"ossencharacter,  see  e.g.   [Silverman  1994]  \cite{S2},   p. 175;
moreover $L({\cal E}_{CM}^{(-D,f)},s)=L(w, s)$,  where $L(w, s):=\sum_{n=1}^{\infty}{c_n\over n^s}$ and $c_n$
the Fourier coefficients of the cusp form $w(s)$.  In other words,   
${\cal E}_{CM}^{(-D,f)}$ is a modular elliptic curve.
One can now apply the modularity principle:  if  $A_w$ is an abelian variety given by 
the periods of holomorphic differential $w(s)ds$ (and its conjugates)  on  
$X(fD)$, then the diagram in Fig. 6.3 is  commutative.  
The holomorphic map  $X(fD)\to {\cal E}_{CM}^{(-D,f)}$ is obtained as  a composition of the canonical
embedding $X(fD)\to A_w$ with the subsequent holomorphic projection $A_w\to {\cal E}_{CM}^{(-D,f)}$.
Lemma \ref{lem6.1.4} is proved.
$\square$

 \index{Gr\"ossencharacter}

\bigskip
\begin{figure}[here]
\begin{picture}(300,110)(-100,-5)
\put(130,70){\vector(0,-1){35}}
\put(55,70){\vector(2,-1){53}}
\put(55,83){\vector(1,0){53}}
\put(128,20){${\cal E}_{CM}^{(-D,f)}$}
\put(10,80){$X(fD)$}
\put(125,80){$A_w$}
\put(60,100){\sf canonical}
\put(60,90){\sf embedding}
\put(140,60){\sf holomorphic}
\put(140,50){\sf projection}
\end{picture}
\caption{Hecke lemma.}
\end{figure}
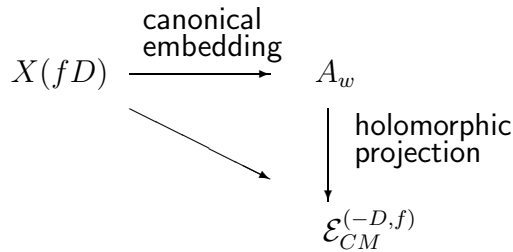

 \index{Hecke lemma}

\begin{lem}\label{lem6.1.5}
The functor $F$ acts by the formula  ${\cal E}_{CM}^{(-D,f)}\mapsto {\cal A}_{RM}^{(D,f)}$. 
\end{lem}
{\it Proof.} Let $L_{CM}$ be a lattice with complex multiplication by an order ${\goth R}={\Bbb Z}+(f\omega){\Bbb Z}$
in the imaginary quadatic field ${\Bbb Q}(\sqrt{-D})$; the multiplication by $\alpha\in {\goth R}$ generates an
endomorphism  $(a,b,c,d)\in M_2({\Bbb Z})$ of the lattice $L_{CM}$. 
It is known from  Section 6.1.3,  {\bf Case 4},  that the endomorphisms of 
lattice $L_{CM}$ and endomorphisms of the pseudo-lattice ${\goth m}_{RM}=F(L_{CM})$ 
are related by the following explicit  map
\displaymath
\left(\matrix{a & b\cr c & d}\right)\in End~(L_{CM})
\longmapsto
 \left(\matrix{a & b\cr -c & -d}\right)\in End~({\goth m}_{RM}),
\enddisplaymath
Moreover, one can always assume  $d=0$ in a proper basis of $L_{CM}$. 
We shall consider the following two cases.

\medskip
{\sf Case I.} If $D\equiv 1~mod~4$ then we have 
${\goth R}={\Bbb Z}+({f+\sqrt{-f^2D}\over 2}){\Bbb Z}$;
thus  $\alpha={2m+fn\over 2}+\sqrt{{-f^2Dn^2\over 4}}$ for some 
$m, n\in {\Bbb Z}$.  Therefore  multiplication by $\alpha$ corresponds to an  endomorphism $(a,b,c,0)\in M_2({\Bbb Z})$,
where 
\displaymath
\left\{
\begin{array}{ccc}
a &=& Tr (\alpha)=\alpha+\bar\alpha=2m+fn\\
b &=& -1\\
c &=& N (\alpha)=\alpha\bar\alpha=\left({2m+fn\over 2}\right)^2+{f^2Dn^2\over 4}. 
\end{array}
\right.
\enddisplaymath
To calculate a primitive generator of endomorphisms of the lattice $L_{CM}$
one should  find a multiplier $\alpha_0\ne 0$ such that
\displaymath
|\alpha_0|=\min_{m.n\in {\Bbb Z}}|\alpha|=\min_{m.n\in {\Bbb Z}}\sqrt{N(\alpha)}.
\enddisplaymath
From the  equation for $c$  the minimum is attained at  $m=-{f\over 2}$ and
$n=1$ if $f$ is even or  $m=-f$ and $n=2$ if $f$ is odd. Thus
\displaymath
\alpha_0=\cases{\pm{f\over 2}\sqrt{-D}, & if $f$ is even\cr
                \pm f\sqrt{-D}, & if $f$ is odd.}
\enddisplaymath
To find the matrix form of the endomorphism $\alpha_0$,  we shall substitute
in the corresponding formula  $a=d=0, b=-1$ and $c={f^2D\over 4}$ if $f$ is even or
$c=f^2D$ if $f$ is odd.  Thus functor $F$  maps the
multiplier $\alpha_0$ into
\displaymath
F(\alpha_0)=\cases{\pm{f\over 2}\sqrt{D}, & if $f$ is even\cr
                \pm f\sqrt{D}, & if $f$ is odd.}
\enddisplaymath
Comparing the above equations,   one verifies
that formula  $F({\cal E}_{CM}^{(-D,f)})={\cal A}_{RM}^{(D,f)}$ is true in this case.

\medskip
{\sf Case II.} If $D\equiv 2$ or $3 ~mod~4$ then 
${\goth R}={\Bbb Z}+(\sqrt{-f^2D})~{\Bbb Z}$;
thus the multiplier $\alpha=m+\sqrt{-f^2Dn^2}$ for some 
$m, n\in {\Bbb Z}$.  A  multiplication by $\alpha$ corresponds to an  
endomorphism $(a,b,c,0)\in M_2({\Bbb Z})$,  where 
\displaymath
\left\{
\begin{array}{ccc}
a &=& Tr (\alpha)=\alpha+\bar\alpha=2m\\
b &=& -1\\
c &=& N (\alpha)=\alpha\bar\alpha=m^2+f^2Dn^2. 
\end{array}
\right.
\enddisplaymath
We shall repeat the argument of {\sf Case I};  then
from the  equation for $c$  the minimum of $|\alpha|$ is attained at  $m=0$ and
$n=\pm 1$. Thus $\alpha_0=\pm f\sqrt{-D}$. 
To find the matrix form of the endomorphism $\alpha_0$ we substitute
in the corresponding equation   $a=d=0, b=-1$ and $c=f^2D$.  Thus  functor $F$  maps the
multiplier $\alpha_0=\pm f\sqrt{-D}$  into $F(\alpha_0)=\pm f\sqrt{D}$.  
 In other words, formula  $F({\cal E}_{CM}^{(-D,f)})={\cal A}_{RM}^{(D,f)}$ is true in this case
as well.   Since all possible cases are exhausted,  Lemma \ref{lem6.1.5} is proved.
$\square$

\begin{lem}\label{lem6.1.6}
For every $N\ge 1$  there exists a holomorphic map $X_1(N)\to X(N)$.
\end{lem}
{\it Proof.} Indeed, $\Gamma(N)$ is a normal subgroup of index $N$  of the
group $\Gamma_1(N)$;  therefore there exists a degree $N$ holomorphic 
map $X_1(N)\to X(N)$.   Lemma \ref{lem6.1.6} follows. 
$\square$

\bigskip
Theorem \ref{thm6.1.3}  follows from Lemmas \ref{lem6.1.3}-\ref{lem6.1.5} and  Lemma \ref{lem6.1.6} for
$N=fD$. 
$\square$

\vskip1.5cm\noindent
{\bf Guide to the literature.}
D.~Hilbert counted complex multiplication as not only the most beautiful
part of mathematics but also of entire science;  it surely  does as it links
complex analysis and number theory.  One cannot beat [Serre 1967] \cite{Ser2}
for an introduction,  but more comprehensive   [Silverman 1994]   \cite{S2},  Chapter 2
is the must.  Real multiplication has been introduced in [Manin 2004]  \cite{Man1}.  
The link between the two was  the subject of \cite{Nik4}  and the inverse functor
$F^{-1}$ was constructed in \cite{Nik5}.

 \index{$K$-rational elliptic curve}

\section{Ranks of the $K$-rational elliptic curves}
We are working in the category of elliptic curves with complex multiplication;
such curves were  denoted by ${\cal E}_{CM}^{(-D,f)}$,  where $f\ge 1$ is
the conductor of an order in the imaginary quadratic field $k={\Bbb Q}(\sqrt{-D})$. 
Recall that  ${\cal E}_{CM}^{(-D,f)}\cong {\cal E}(K)$,  where $K=k(j({\cal E}_{CM}^{(-D,f)}))$
is the Hilbert class field of $k$,  see e.g.   [Serre 1967] \cite{Ser2}.  In other words,
we deal with a $K$-rational  projective  curve
\displaymath
{\cal E}(K)=\{(x,y,z)\in {\Bbb C}P^2 ~|~ y^2z=4x^3-g_2xz^2-g_3z^3\},
\enddisplaymath
where constants $g_2$ and  $g_3$  belong to the number field $K$.  
It is well known,  that any pair of points $p,p'\in {\cal E}(K)$ defines a
sum $p+p'\in {\cal E}(K)$  and an inverse $-p\in {\cal E}(K)$ so that ${\cal E}(K)$ 
has the structure of an  abelian group;  the next result is now a standard fact,
see e.g.  [Tate 1974]  \cite{Tat1}, p. 192.  
\begin{thm}
{\bf (Mordell-N\'eron)}
The ${\cal E}(K)$ is finitely generated abelian group.
\end{thm}
 \index{Mordell-N\'eron Theorem}
 \index{rank of elliptic curve}
\begin{dfn}
By $rk~({\cal E}_{CM}^{(-D,f)})$ we understand the integer number equal to the
rank of abelian group ${\cal E}(K)$.  
\end{dfn}
\begin{rmk}
\textnormal{
The  $rk~({\cal E}_{CM}^{(-D,f)})$ is an invariant of the $K$-isomorphism class
of  ${\cal E}_{CM}^{(-D,f)}$ but not of the general isomorphism class;  those variations
of the rank are known as {\it twists} of ${\cal E}_{CM}^{(-D,f)}$.   Of course,  if two curves
are $K$-isomorphic,  they are also isomorphic over ${\Bbb C}$.  
}
\end{rmk}
In what remains,  we calculate $rk~({\cal E}_{CM}^{(-D,f)})$ in terms of 
 invariants of the noncommutative  torus  ${\cal A}_{RM}^{(D,f)}=F({\cal E}_{CM}^{(-D,f)})$;
one such invariant called an {\it arithmetic complexity} will be introduced below.

 \index{arithmetic complexity}

\subsection{Arithmetic complexity of noncommutative tori}
Let $\theta$ be a quadratic irrationality,  i.e.  irrational root of a quadratic
polynomial  $ax^2+bx+c=0$, where $a,b,c\in {\Bbb Z}$;   denote by 
$Per~(\theta):=(\overline{a_1, a_2, \dots, a_P})$ the minimal period of continued fraction
of $\theta$ taken up to a cyclic permutation.  
Fix  $P$  and suppose for a moment that $\theta$ is  a function  of
its period
\displaymath
\theta(x_0, x_1,\dots,x_P)=[x_0, \overline{x_1,\dots,x_P}], 
\enddisplaymath
where $x_i\ge 1$ are integer variables;  then
$\theta(x_0,\dots, x_P)\in {\Bbb Q}+\sqrt{{\Bbb Q}}$,
where  $\sqrt{{\Bbb Q}}$ are square roots of positive rationals. 
Consider a  constraint  (a restriction)   $x_1=x_{P-1}, x_2=x_{P-2},\dots, x_P=2x_0$; 
 then $\theta(x_0, x_1, x_2, \dots,x_2, x_1,  2x_0)\in \sqrt{{\Bbb Q}}$,
see e.g.  [Perron 1954]  \cite{P}, p. 79.  Notice, that in this case there are   
 ${1\over 2}P+1$ independent variables, if $P$ is even
and ${1\over 2}(P+1)$,   if $P$ is odd.  
The number of independent variables  will further decrease,  if $\theta$ is square root of an integer;
let us introduce some notation. 
For a regular fraction $[a_0,a_1,\dots]$ one associates the linear equations
\displaymath
\left\{
\begin{array}{ccc}
y_0 &=& a_0y_1+y_2\\
y_1 &=& a_1y_2+y_3\\
y_2 &=& a_2y_3+y_4\\
       &\vdots& 
\end{array}
\right.
\enddisplaymath
One can  put  above equations in the form
\displaymath
\left\{
\begin{array}{ccc}
y_j &=& A_{i-1,j}y_{i+j}+a_{i+j}A_{i-2,j}y_{i+j+1}\\
y_{j+1} &=& B_{i-1,j}y_{i+j}+a_{i+j}B_{i-2,j}y_{i+j+1},
\end{array}
\right.
\enddisplaymath
where the polynomials $A_{i,j}, B_{i,j}\in {\Bbb Z}[a_0,a_1,\dots]$ are called
{\it Muir's symbols}, see   [Perron 1954]  \cite{P}, p.10.  The following lemma will play an
important r\^ole.    
 \index{Muir symbols}
\begin{lem}\label{lem6.2.1}
{\bf  ([Perron 1954]  \cite{P},  pp. 88 and 107)} 
There exists  a square-free integer $D>0$,  such that
\displaymath
 [x_0, \overline{x_1, \dots, x_1, x_P}]=
 \cases{\sqrt{D}, & if  $x_P=2x_0$  and $D=2,3~\mod~4,$\cr
                {\sqrt{D}+1\over 2}, & if   $x_P=2x_0-1$ and  $D=1~\mod~4,$}
\enddisplaymath
if and only if   $x_P$ satisfies  the diophantine equation
\displaymath
x_P=mA_{P-2, 1}-(-1)^PA_{P-3,1}B_{P-3,1}, 
\enddisplaymath
for an integer  $m>0$;  moreover, in this case $D= {1\over 4}x_P^2+mA_{P-3,1}-(-1)^PB_{P-3,1}^2$.
\end{lem}
Let $(x_0^*,\dots, x_P^*)$ be a solution of the diophantine equation of Lemma \ref{lem6.2.1}.
By  {\it dimension},  $d$,   of  this  solution one understands the maximal number
of variables $x_i$,  such that for every $s\in {\Bbb Z}$  there exists a solution
of  the above diophantine equation  of the form   $(x_0,\dots, x_i^*+s,\dots, x_P)$. 
In geometric terms,  $d$ is equal  to  dimension of a  connected component through the 
 point $(x_0^*,\dots,x_P^*)$  of an  affine variety $V_m$ (i.e. depending on $m$) 
  defined by the diophantine  equation.   For the sake of clarity,  let us consider a simple example.  
\begin{exm}\label{exm6.2.1}
{\bf ([Perron 1954]  \cite{P}, p. 90)}
{\normalfont
If  $P=4$,   then Muir's symbols are: $A_{P-3,1}=A_{1,1}=x_1x_2+1$, ~$B_{P-3,1}=B_{1,1}=x_2$ and 
$A_{P-2,1}=A_{2,1}=x_1x_2x_3+x_1+x_3=x_1^2x_2+2x_1$, since $x_3=x_1$.  Thus, our
diophantine equation takes the form 
\displaymath
2x_0= m(x_1^2x_2+2x_1)-x_2(x_1x_2+1),
\enddisplaymath
and, therefore,  $\sqrt{x_0^2+m(x_1x_2+1)-x_2^2}=[x_0,\overline{x_1,x_2,x_1,2x_0}]$. 
First,  let us show that the affine variety defined by the last equation is not connected. 
Indeed,  by Lemma \ref{lem6.2.1},  parameter $m$ must be integer for all (integer) values
of $x_0,x_1$ and $x_2$.  This is not possible in general,  since from our last  equation 
one obtains  $m=(2x_0+x_2(x_1x_2+1))(x_1^2x_2+2x_1)^{-1}$ is a rational number. 
However,  a restriction to  $x_1=1, ~x_2=x_0-1$ defines a  (maximal)
connected component of the variety corresponding to our equation,  since in this case $m=x_0$ is always an integer.    
Thus,  one gets a family of solutions of of the form
 $\sqrt{(x_0+1)^2-2}=[x_0,\overline{1, x_0-1, 1, 2x_0}]$, 
 where   each solution has dimension $d=1$. 
 }
 \end{exm}
\begin{dfn}\label{dfn6.2.2}
By an arithmetic complexity $c({\cal A}_{RM}^{(D,f)})$ of the noncommutative torus 
 ${\cal A}_{RM}^{(D,f)}$  one  understands  an integer equal to  dimension $d$ of solution 
 $(x_0^*,\dots, x_P^*)$ of  diophantine equation 
$x_P=mA_{P-2, 1}-(-1)^PA_{P-3,1}B_{P-3,1}$;
if ${\cal A}_{\theta}$ has no real multiplication,  then  
 the arithmetic  complexity is assumed to be infinite.
 \end{dfn}

 \index{${\Bbb Q}$-curve}

\subsection{${\Bbb Q}$-curves}
For the sake of simplicity,  we shall restrict our considerations to a family of elliptic curves ${\cal E}_{CM}^{(-D,f)}$
known as the {\it ${\Bbb Q}$-curves};   a general result exists only in a conjectural form so far,
see {\it Exercises, problems and conjectures}.    
\begin{dfn}
{\bf ([Gross 1980]   \cite{G})}
Let $({\cal E}_{CM}^{(-D,f)})^{\sigma},  ~\sigma\in Gal~(k|{\Bbb Q})$ be 
the  Galois conjugate of the curve  ${\cal E}_{CM}^{(-D,f)}$; 
by a {\it ${\Bbb Q}$-curve} one understands   ${\cal E}_{CM}^{(-D,f)}$,
such that  there exists  an  isogeny between $({\cal E}_{CM}^{(-D,f)})^{\sigma}$ and   ${\cal E}_{CM}^{(-D,f)}$
for  each    $\sigma\in Gal~(k|{\Bbb Q})$. 
 \end{dfn}
\begin{rmk}
\textnormal{
The curve ${\cal E}_{CM}^{(-p,1)}$ is a ${\Bbb Q}$-curve,  
whenever $p\equiv  3 ~\mod ~4$ is a prime number,  see [Gross 1980]  \cite{G},  p. 33;
we shall write ${\goth P}_{3 ~\mod ~4}$ to denote the set of all such primes. 
}
\end{rmk}
\begin{rmk}
\textnormal{
The rank of  ${\cal E}_{CM}^{(-p,1)}$  is always divisible by $2h_k$,  where $h_k$ is the 
 class number of number  field $k:={\Bbb Q}(\sqrt{-p})$,  see [Gross 1980]   \cite{G},  p. 49.  
}
\end{rmk}
 \index{${\Bbb Q}$-rank}
\begin{dfn}
By a ${\Bbb Q}$-rank   of  ${\cal E}_{CM}^{(-p,1)}$
one  understands  the integer
\displaymath 
 rk_{\Bbb Q}({\cal E}_{CM}^{(-p,1)}):={1\over 2h_k}~rk~({\cal E}_{CM}^{(-p,1)}).   
 \enddisplaymath
 \end{dfn}
The following result links invariants of   noncommutative tori  and   geometry of  the 
$K$-rational elliptic curves;   namely,  the arithmetic complexity plus one is equal to the ${\Bbb Q}$-rank
of the corresponding elliptic curve. 
\begin{thm}\label{thm6.2.1}
$rk_{\Bbb Q}~({\cal E}_{CM}^{(-p,1)}) +1=c({\cal A}_{RM}^{(p,1)})$
whenever $p\equiv 3 ~\mod ~4$. 
\end{thm}
\begin{rmk}
\textnormal{
The general formula $rk~({\cal E}_{CM}^{(-D,f)}) +1=c({\cal A}_{RM}^{(D,f)})$
for all $D\ge 2$ and $f\ge 1$ is known as the {\it rank conjecture},  see
 {\it Exercises, problems and conjectures}. 
}
\end{rmk}

\subsection{Proof of  Theorem  \ref{thm6.2.1}}
\begin{lem}\label{lem6.2.2}
If $[x_0,\overline{x_1,\dots, x_k,\dots, x_1, 2x_0}]\in\sqrt{{\goth P}_{3~\mod~4}}$  ,
then:

\smallskip
(i) $P=2k$ is an even number,  such that:

\smallskip
\hskip1cm
(a) $P\equiv 2~\mod~4$,  if $p\equiv 3~\mod~8$;

\smallskip
\hskip1cm
(b) $P\equiv 0~\mod~4$,  if $p\equiv 7~\mod~8$;

\medskip
(ii) either of two is true:

\smallskip
\hskip1cm
(a) $x_k=x_0$ (a culminating period);

 \index{culminating period}
 \index{almost-culminating period}

\smallskip
\hskip1cm
(b) $x_k=x_0-1$ and $x_{k-1}=1$ (an almost-culminating period).
\end{lem}
{\it Proof.} (i)  Recall that  if $p\ne 2$ is a prime,  then one and only one of the following
diophantine equations is solvable:
\displaymath
\left\{
\begin{array}{ccc}
x^2-py^2 &=& -1,\\
x^2-py^2 &=& 2,\\
x^2-py^2 &=& -2,
\end{array}
\right.
\enddisplaymath
see e.g. [Perron 1954]  \cite{P}, Satz 3.21.  Since  $p\equiv 3~\mod~4$,  one concludes
that $x^2-py^2=-1$ is not solvable [Perron 1954]  \cite{P}, Satz 3.23-24;  this happens if and only if $P=2k$ is
even (for otherwise the continued fraction of $\sqrt{p}$ would provide a
solution).

It is known,  that for even periods $P=2k$ the convergents  $A_i/B_i$ satisfy the
diophantine equation $A_{k-1}^2-pB_{k-1}^2=(-1)^k ~2$,  see [Perron 1954]  \cite{P},  p.103;
thus if $P\equiv 0~\mod~4$,  the equation $x^2-py^2=2$ is solvable and if
$P\equiv 2~\mod~4$,  then the equation $x^2-py^2=-2$ is sovable.
But  equation   $x^2-py^2=2$  (equation $x^2-py^2=-2$,  resp.) is solvable if and only if 
$p\equiv 7~\mod~8$  ($p\equiv 3~\mod~8$, resp.),  see [Perron 1954]  \cite{P}, Satz 3.23
(Satz 3.24, resp.).  Item (i) follows.

\medskip
(ii)   The equation $A_{k-1}^2-pB_{k-1}^2=(-1)^k ~2$ is a special case of 
equation $A_{k-1}^2-pB_{k-1}^2=(-1)^k ~Q_k$,  where $Q_k$ is the full
quotient of continued fraction [Perron 1954]  \cite{P}, p.92;  therefore, $Q_k=2$.
One  can now apply  [Perron 1954]  \cite{P},   Satz 3.15,  which says that for $P=2k$ and 
$Q_k=2$ the continued fraction of    $\sqrt{{\goth P}_{3~\mod~4}}$ is either
culminating (i.e.   $x_k=x_0$) or almost-culminating  (i.e. $x_k=x_0-1$ and $x_{k-1}=1$).
Lemma \ref{lem6.2.2}  follows.
$\square$

\medskip
\begin{lem}\label{lem6.2.3}
If $p\equiv 3~\mod~8$,   then  $c({\cal A}_{RM}^{(p,1)})=2$.
\end{lem}
{\it Proof.}  The proof proceeds by induction in period $P$, which is
in this case   $P\equiv 2~\mod~4$ by Lemma \ref{lem6.2.2}.    We shall start
with $P=6$, since $P=2$ reduces to it,   see item (i) below.

\smallskip
(i)   Let $P=6$ be a culminating period;  then the diophantine equation in Definition \ref{dfn6.2.2}  admits 
a general solution $[x_0,\overline{x_1,  2x_1,  x_0, ,  2x_1, x_1,  2x_0}]=\sqrt{x_0^2+4nx_1+2}$,
where $x_0=n(2x_1^2+1)+x_1$,  see   [Perron 1954]  \cite{P},  p. 101.   The solution depends on two
integer variables $x_1$ and $n$,  which is the maximal  possible number of  variables in this case;
therefore,  the dimension of the solution is $d=2$,  so as  complexity of the corresponding torus.   
Notice that the case  $P=2$ is obtained from $P=6$ by restriction to  $n=0$;  thus the complexity   
for $P=2$ is equal to $2$.

\smallskip
(ii)   Let $P=6$ be an almost-culminating period;  
 then the diophantine equation in Definition \ref{dfn6.2.2} has 
a solution $[3s+1,\overline{2,  1,  3s ,  1,  2,  6s+2}]=\sqrt{(3s+1)^2+2s+1}$,
where $s$ is an integer variable  [Perron 1954]  \cite{P},  p. 103.  We encourage the
reader to verify,  that this solution is a restriction of solution (i) to $x_1=-1$
and $n=s+1$;   thus,  the dimension of our solution is $d=2$, so as the
complexity of the corresponding torus.

\smallskip
(iii)  Suppose  a solution $[x_0,\overline{x_1,\dots, x_{k-1}, x_k, x_{k-1}, \dots, x_1, 2x_0}]$
with the (culminating or almost-culminating) period $P_0\equiv 3~\mod~8$ has dimension $d=2$;
let us show that a solution
\displaymath
[x_0,\overline{y_1, x_1,\dots, x_{k-1}, y_{k-1},  x_k, y_{k-1},  x_{k-1}, \dots, x_1, y_1, 2x_0}]
\enddisplaymath
with period $P_0+4$ has also dimension $d=2$.  According to  [Weber 1926]  \cite{Web1}, 
if above fraction  is a solution to the diophantine equation in Definition \ref{dfn6.2.2},  then either
(i) $y_{k-1}=2y_1$ or (ii) $y_{k-1}=2y_1+1$ and $x_1=1$.  We proceed by showing 
that case (i)  is not possible for the square roots of prime numbers.

Indeed, let to the contrary $y_{k-1}=2y_1$;  then the following system of equations
must be compatible:
\displaymath
\left\{
\begin{array}{ccc}
A_{k-1}^2 &-&   pB_{k-1}^2=-2,\\
A_{k-1} &=& 2y_1A_{k-2}+A_{k-3},\\
B_{k-1} &=& 2y_1B_{k-2}+B_{k-3},
\end{array}
\right.
\enddisplaymath
 where $A_i, B_i$  are convergents and the first equation is solvable
 since $p\equiv 3~\mod~8$.  From the first equation, both convergents 
 $A_{k-1}$ and $B_{k-1}$ are odd numbers. (They are both odd or even,
 but even excluded,  since $A_{k-1}$ and $B_{k-1}$ are relatively prime.)
 From the last two equations, the convergents $A_{k-3}$ and $B_{k-3}$
 are also odd.  Then the convergents $A_{k-2}$ and $B_{k-2}$ must
 be even,  since among six consequent  convergents $A_{k-1}, B_{k-1}, A_{k-2}, B_{k-2},
 A_{k-3}, B_{k-3}$ there are always two even;   but this is not possible,  because
 $A_{k-2}$ and $B_{k-2}$ are relatively prime. Thus, $y_{k-1}\ne 2y_1$.

 Therefore the above equations give  a solution of
 the diophantine equation in Definition \ref{dfn6.2.2}
  if and only if  $y_{k-1}=2y_1+1$ and $x_1=1$;  the dimension of such a solution
 coincides with the dimension of solution 
 $$[x_0,\overline{x_1,\dots, x_{k-1}, x_k, x_{k-1}, \dots, x_1, 2x_0}],$$
 since for two new integer variables $y_1$ and $y_{k-1}$ one gets
 two new constraints.  Thus,  the dimension of the above solution 
 is $d=2$,  so as the  complexity of the corresponding torus. 
 Lemma \ref{lem6.2.3} follows.
 $\square$

 \medskip
\begin{lem}\label{lem6.2.4}
If $p\equiv 7~\mod~8$,   then  $c({\cal A}_{RM}^{(p,1)})=1$.
\end{lem}
{\it Proof.} 
 The proof proceeds by induction in period $P\equiv 0~\mod~4$, see 
  Lemma \ref{lem6.2.2};   we  start with $P=4$.

\smallskip
(i)   Let $P=4$ be a culminating period;  then equation in Definition \ref{dfn6.2.2} admits 
a solution $[x_0,\overline{x_1,x_2,x_1,2x_0}]= \sqrt{x_0^2+m(x_1x_2+1)-x_2^2}$,
where $x_2=x_0$,   see  Example  \ref{exm6.2.1}  for the details.   
Since  the polynomial  $m(x_0x_1+1)$  under   the square root
represents   a prime number,  we have   $m=1$;   the latter equation
is not solvable in integers $x_0$ and $x_1$,   since 
$m=x_0(x_0x_1+3)x_1^{-1}(x_0x_1+2)^{-1}$.   Thus,  there are no solutions
of the diophantine equation in Definition \ref{dfn6.2.2}  with the culminating period $P=4$.

\smallskip
(ii)   Let $P=4$ be an almost-culminating period;  then equation  in Definition \ref{dfn6.2.2} 
admits  a solution  $[x_0,\overline{1, x_0-1, 1, 2x_0}]= \sqrt{(x_0+1)^2-2}$.
The dimension of this solution was proved to be $d=1$,  see Example \ref{exm6.2.1};
thus,  the complexity of the corresponding torus is equal  to $1$.

\smallskip
(iii)  Suppose  a solution $[x_0,\overline{x_1,\dots, x_{k-1}, x_k, x_{k-1}, \dots, x_1, 2x_0}]$
with the (culminating or almost-culminating) period $P_0\equiv 7~\mod~8$ has dimension $d=1$.
It can be shown by the same argument as in Lemma \ref{lem6.2.3},  that for 
 a solution of the form  
 $$[x_0,\overline{y_1, x_1,\dots, x_{k-1}, y_{k-1},  x_k, y_{k-1},  x_{k-1}, \dots, x_1, y_1, 2x_0}]$$
  having the  period $P_0+4$ the  dimension remains
 the same, i.e. $d=1$;  we leave details to the reader.  Thus, complexity of the corresponding
 torus is equal to $1$. Lemma \ref{lem6.2.4} follows.
 $\square$

\begin{lem}\label{lem6.2.5}
{\bf ([Gross 1980]  \cite{G},  p. 78)}
\begin{equation}\label{eq17}
 rk_{\Bbb Q} ~({\cal E}_{CM}^{(-p,1)})=
 \cases{1, & if  ~$p\equiv 3~\mod~8$\cr
              0, & if   ~$p\equiv 7~\mod~8.$}
\end{equation}
\end{lem}

\bigskip
Theorem \ref{thm6.2.1} follows from  Lemmas \ref{lem6.2.3}-\ref{lem6.2.5}.
$\square$

\subsection{Numerical examples}
To illustrate Theorem \ref{thm6.2.1} by numerical examples,  we refer 
the reader to Fig. 6.4  with  all ${\Bbb Q}$-curves ${\cal E}_{CM}^{(-p,1)}$ for $p<100$;
notice,  that   there are infinitely many pairwise  non-isomorphic ${\Bbb Q}$-curves  
[Gross 1980]   \cite{G}.

\begin{figure}[here]
\begin{tabular}{c|c|c|c}
\hline
&&&\\
$p\equiv 3~\mod~4$ & $rk_{\Bbb Q}({\cal E}_{CM}^{(-p,1)})$ & $\sqrt{p}$ & $c({\cal A}_{RM}^{(p,1)})$\\
&&&\\
\hline
$3$ & $1$ & $[1,\overline{1,2}]$ & $2$\\
\hline
$7$ & $0$ & $[2,\overline{1,1,1,4}]$ & $1$\\
\hline
$11$ & $1$ & $[3,\overline{3,6}]$ & $2$\\
\hline
$19$ & $1$ & $[4,\overline{2,1,3,1,2,8}]$ & $2$\\
\hline
$23$ & $0$ & $[4,\overline{1,3,1,8}]$ & $1$\\
\hline
$31$ & $0$ & $[5,\overline{1,1,3,5,3,1,1,10}]$ & $1$\\
\hline
$43$ & $1$ & $[6,\overline{1,1,3,1,5,1,3,1,1,12}]$ & $2$\\
\hline
$47$ & $0$ & $[6,\overline{1,5,1,12}]$ & $1$\\
\hline
$59$ & $1$ & $[7,\overline{1,2,7,2,1,14}]$ & $2$\\
\hline
$67$ & $1$ & $[8,\overline{5,2,1,1,7,1,1,2,5,16}]$ & $2$\\
\hline
$71$ & $0$ & $[8,\overline{2,2,1,7,1,2,2,16}]$ & $1$\\
\hline
$79$ & $0$ & $[8,\overline{1,7,1,16}]$ & $1$\\
\hline
$83$ & $1$ & $[9,\overline{9,18}]$ & $2$\\
\hline
\end{tabular}
\caption{The ${\Bbb Q}$-curves ${\cal E}_{CM}^{(-p,1)}$  with  $p<100$.}
\end{figure}

\vskip1.5cm\noindent
{\bf Guide to the literature.}
The problem of ranks of rational elliptic curves was raised by  [Poincar\'e  1901]  \cite{Poi1},
p. 493.   It was proved by [Mordell 1922]  \cite{Mor1}  that ranks of rational and 
by [N\'eron  1952]  \cite{Ner1} that ranks of  the $K$-rational elliptic curves are always
finite.  The ranks of individual elliptic curves  are calculated by the method of 
{\it descent},  see e.g.  [Cassels 1966]  \cite{Cas1},  p.205.  More conceptual 
approach uses an analytic object called the {\it Hasse-Weil L-function} $L({\cal E}_{\tau}, s)$;
it was conjectured by B.~J.~Birch and  H.~P.~F.~Swinnerton-Dyer that the
order of zero of such a function at $s=1$  is equal to the rank of ${\cal E}_{\tau}$,
see e.g. [Tate 1974] \cite{Tat1}, p. 198.  The {\it rank conjecture}  involving 
invariants of noncommutative tori was formulated in \cite{Nik4}  and proved 
in \cite{Nik6} for the ${\Bbb Q}$-curves ${\cal E}_{CM}^{(-p,1)}$ with
prime $p\equiv 3~\mod~4$.

 \index{Birch and Swinnerton-Dyer Conjecture}

 \index{non-commutative reciprocity}

\section{Non-commutative reciprocity}
In the world of $L$-functions each equivalence between two $L$-functions
is called a {\it reciprocity},  see e.g.  [Gelbart 1984]  \cite{Gel1}. 
In this section we shall introduce an $L$-function $L({\cal A}_{RM}, s)$ 
associated to the noncommutative torus with real multiplication and prove
that any such coincides with the classical Hasse-Weil function $L({\cal E}_{CM},s)$
of an elliptic curve with complex multiplication.  The necessary and sufficient
condition for such a reciprocity is the relation
\displaymath
{\cal A}_{RM}=F({\cal E}_{CM}),
\enddisplaymath
where $F:$ {\bf Ell} $\to$ {\bf NC-Tor}  is the functor introduced 
in Section 6.1.1;    we shall call  such a relation a {\it non-commutative reciprocity},
because it involves invariants of  the  non-commutative algebra ${\cal A}_{RM}$.  
The reciprocity provides us with  explicit  {\it localization formulas}  
for the torus ${\cal A}_{RM}$ at each prime number $p$;  we shall use 
these formulas in the sequel.      
\begin{rmk}
\textnormal{
The reader can think of the non-commutative reciprocity
as an analog of the  Eichler-Shimura theory;  recall that
such a theory  identifies the $L$-function coming from
certain cusp form (of weight two) and the Hasse-Weil function 
  of a  rational elliptic curve,   see e.g.  [Knapp 1992]  \cite{K1}, 
  Chapter XI.  
}
\end{rmk}

 \index{Eichler-Shimura theory}
 \index{$L$-function of noncommutative torus}
 \index{Hasse-Weil $L$-function}

\subsection{$L$-function of noncommutative tori}
Let $p$ be a prime number and ${\cal E}_{CM}$ be an elliptic curve with complex multiplication;   
denote by  ${\cal E}_{CM}({\Bbb F}_p)$   localization  of 
the  ${\cal E}_{CM}$  at the prime ideal ${\goth P}$ over $p$, see e.g.  [Silverman 1994]  \cite{S2},  p.171.  
We are looking for a proper concept of   localization  of the  algebra 
${\cal A}_{RM}=F({\cal E}_{CM})$  corresponding   to the localization   ${\cal E}_{CM}({\Bbb F}_p)$  
of  elliptic curve ${\cal E}_{CM}$  at prime $p$.   To attain the goal,  recall that 
the cardinals $|{\cal E}_{CM}({\Bbb F}_p)|$ generate  the Hasse-Weil function $L({\cal E}_{CM},s)$ of 
the curve ${\cal E}_{CM}$,  see e.g.  [Silverman 1994]  \cite{S2},  p.172;  
thus,  we have to  define   an $L$-function of the noncommutative  torus ${\cal A}_{RM}=F({\cal E}_{CM})$
equal to  the Hasse-Weil function of  the curve ${\cal E}_{CM}$.  
\begin{dfn}
If ${\cal A}_{RM}$ is a noncommutative torus with real multiplication,  consider
an integer matrix
\displaymath
A=\left(\matrix{a_1 & 1\cr 1 & 0}\right)\dots \left(\matrix{a_k & 1\cr 1 & 0}\right),
\enddisplaymath
where $\overline{(a_1,\dots,a_k)}$  is the minimal period of continued fraction
of a quadratic irrationality  $\theta$ corresponding to ${\cal A}_{RM}$.  For each prime $p$
consider an integer matrix
\displaymath
L_p:=\left(\matrix{tr~(A^{\pi(p)}) & p\cr -1 & 0}\right),
\enddisplaymath
where $tr~(\bullet)$ is the trace of a matrix and $\pi(n)$ is an integer-valued function
defined in the Supplement 6.3.3.   By a local  zeta function  of  torus  ${\cal A}_{RM}$
one understands the analytic function
\displaymath
\zeta_p({\cal A}_{RM}, z):=
\exp\left(
\sum_{n=1}^{\infty}{|K_0({\cal O}_{\varepsilon_n})|\over n} ~z^n
\right),
~\varepsilon_n=
\cases{L_p^n, & \mbox{if}  $p\nmid tr^2(A)-4$   \cr
1-\alpha^n,  & \mbox{if}  $p ~| ~tr^2(A)-4$,}
\enddisplaymath
where $\alpha\in\{-1, 0, 1\}$,  
${\cal O}_{\varepsilon_n}={\cal A}_{RM}\rtimes_{\varepsilon_n}{\Bbb Z}$
is the Cuntz-Krieger algebra and $K_0(\bullet)$ its $K_0$-group,   see Section 3.7.   
By an $L$-function of the  noncommutative torus  ${\cal A}_{RM}$
we understand  the analytic function
\displaymath
L({\cal A}_{RM}, s):=\prod_{p}  ~\zeta_p({\cal A}_{RM}, p^{-s}), \quad  s\in {\Bbb C},
\enddisplaymath
where $p$ runs through the set of all prime numbers.    
\end{dfn}
\begin{thm}\label{thm6.3.1}
The following conditions are equivalent:

\medskip
(i)  \hskip3cm   ${\cal A}_{RM}=F({\cal E}_{CM})$;

\smallskip
(ii) 
$$
\left\{
\begin{array}{ccc}
L({\cal A}_{RM}, s) &\equiv& L({\cal E}_{CM}, s),\\
K_0({\cal O}_{\varepsilon_n}) &\cong&  {\cal E}_{CM}({\Bbb F}_{p^n}),
\end{array}
\right.
$$
where $F:$ {\bf Ell} $\to$ {\bf NC-Tor} is the functor defined in Section 6.1.1 and  $L({\cal E}_{CM}, s)$
is the Hasse-Weil $L$-function of elliptic curve ${\cal E}_{CM}$. 
\end{thm}
\begin{rmk}
{\bf (Non-commutative localization)}
\textnormal{
Theorem \ref{thm6.3.1} implies  a localization formula for the  torus 
${\cal A}_{RM}$ at  a prime $p$,  since   the Cuntz-Krieger algebra 
${\cal O}_{\varepsilon_n}\cong {\cal A}_{RM}\rtimes_{\varepsilon_n}{\Bbb Z}$ 
can be viewed as  a  non-commutative coordinate ring  of  elliptic curve  
${\cal E}_{CM}({\Bbb F}_{p^n})$.  Thus,  to localize a non-commutative ring
one takes its crossed product  rather than taking its prime (or maximal) ideal as prescribed by the
familiar  formula  for  commutative rings.  
}
\end{rmk}

 \index{localization formula}

\subsection{Proof of  Theorem \ref{thm6.3.1}}
Let $p$ be such, that ${\cal E}_{CM}$ has a good reduction at ${\goth P}$;  
the corresponding local zeta function   $\zeta_p({\cal E}_{CM},z)=(1-tr~(\psi_{{\cal E}(K)}({\goth P}))z+pz^2)^{-1}$,
where $\psi_{{\cal E}(K)}$ is the Gr\"ossencharacter on $K$ and  $tr$ is the trace of algebraic 
number.  We have to prove,  that $\zeta_p({\cal E}_{CM},z)=\zeta_p({\cal   A}_{RM},z):=(1-tr~(A^{\pi(p)})z+pz^2)^{-1}$;
the last equality is a consequence of definition of $\zeta_p({\cal  A}_{RM},z)$.
Let ${\cal E}_{CM}\cong {\Bbb C}/L_{CM}$,  where $L_{CM}={\Bbb Z}+{\Bbb Z}\tau$
is a lattice in the complex plane  [Silverman 1994]  \cite{S},  pp. 95-96;  
let $K_0({\cal A}_{RM})\cong{\goth m}_{RM}$,  where 
${\goth m}_{RM}={\Bbb Z}+{\Bbb Z}\theta$ is a {\it pseudo-lattice} in 
${\Bbb R}$,  see  [Manin 2004]  \cite{Man1}. 
Roughly speaking,  we construct an invertible element (a unit) $u$ of
the ring $End~({\goth m}_{RM})$ attached to pseudo-lattice ${\goth m}_{RM}=F(L_{CM})$,
  such  that
\displaymath
tr~(\psi_{{\cal E}(K)}({\goth P}))=tr~(u)=tr~(A^{\pi(p)}).
\enddisplaymath
The latter  will be achieved with the help of an explicit formula 
connecting endomorphisms of lattice $L_{CM}$ with such of 
the pseudo-lattice ${\goth m}_{RM}$
\displaymath
\left(\matrix{a & b\cr c & d}\right)\in End~(L_{CM})
\longmapsto
 \left(\matrix{a & b\cr -c & -d}\right)\in End~({\goth m}_{RM}),
\enddisplaymath
see  proof of Lemma \ref{lem6.1.5} for the details.  
We shall split the proof into a series of lemmas,
 starting with the following elementary lemma.  
\begin{lem}\label{lem6.3.1}
Let $A=(a,b,c,d)$ be an integer matrix with $ad-bc\ne 0$
and $b=1$. Then $A$ is similar to the matrix
$(a+d, 1, c-ad, 0)$.
\end{lem}
{\it Proof.} Indeed, take a matrix $(1,0,d,1)\in SL_2({\Bbb Z})$.
The matrix realizes the similarity, i.e.
\displaymath
\left(\matrix{1 & 0\cr -d & 1}\right)
\left(\matrix{a & 1\cr  c & d}\right)
\left(\matrix{1 & 0\cr d & 1}\right)=
\left(\matrix{a+d & 1\cr c-ad & 0}\right).
\enddisplaymath
Lemma \ref{lem6.3.1} follows.
$\square$

\begin{lem}\label{lem6.3.2}
The matrix $A=(a+d, 1, c-ad, 0)$ is similar to its
transpose $A^t=(a+d, c-ad, 1, 0)$. 
\end{lem}
{\it Proof.} We shall use the following criterion: the
(integer) matrices $A$ and $B$ are similar, if and only if 
the characteristic matrices $xI-A$ and $xI-B$ have the same Smith normal
form.  The calculation for the matrix $xI-A$ gives
$$
\left(\matrix{x-a-d & -1\cr ad-c & x}\right)\sim
\left(\matrix{x-a-d & -1\cr  x^2-(a+d)x+ad-c & 0}\right)\sim
$$
$$
\sim \left(\matrix{1 & 0\cr 0 & x^2-(a+d)x+ad-c}\right),
$$
where $\sim$ are the elementary operations between the rows (columns)  
of the matrix.  Similarly, a calculation for the matrix $xI-A^t$
gives
$$
\left(\matrix{x-a-d & ad-c\cr -1 & x}\right)\sim
\left(\matrix{x-a-d & x^2-(a+d)x+ad-c\cr -1& 0}\right)\sim
$$
$$
\sim\left(\matrix{1 & 0\cr 0 & x^2-(a+d)x+ad-c}\right).
$$
Thus, $(xI-A)\sim (xI-A^t)$ and  Lemma \ref{lem6.3.2} follows.
$\square$

\begin{cor}\label{cor6.3.1}
The matrices $(a, 1, c, d)$ and $(a+d, c-ad, 1, 0)$
are similar. 
\end{cor}
{\it Proof.} The Corollary \ref{cor6.3.1}  follows from Lemmas \ref{lem6.3.1}--\ref{lem6.3.2}.
$\square$

\bigskip\noindent
Recall  that if ${\cal E}_{CM}^{(-D,f)}$ is an elliptic curve with complex multiplication
by order $R={\Bbb Z}+fO_k$ in imaginary quadratic field $k={\Bbb Q}(\sqrt{-D})$, 
then ${\cal A}_{RM}^{(D,f)}=F({\cal E}_{CM}^{(-D,f)})$ is the noncommutative torus
with real multiplication by the order ${\goth R}={\Bbb Z}+fO_{\goth k}$ in real quadratic
field ${\goth k}={\Bbb Q}(\sqrt{D})$.  
\begin{lem}\label{lem6.3.3}
Each $\alpha\in R$ goes under $F$ into an $\omega\in {\goth R}$,
such that $tr~(\alpha)=tr~(\omega)$,  where $tr~(x)=x+\bar x$ is the trace
of an algebraic number $x$.    
\end{lem}
{\it Proof.} Recall that each $\alpha\in R$ can be written  in a matrix
form for a given base $\{\omega_1,\omega_2\}$ of the lattice
$L_{CM}$.  Namely,
\displaymath
\left\{
\begin{array}{cc}
\alpha\omega_1 &= a\omega_1 +b\omega_2\\
\alpha\omega_2 &= c\omega_1 +d\omega_2,
\end{array}
\right.
\enddisplaymath
where $(a,b,c,d)$ is an integer matrix with $ad-bc\ne 0$. 
and $tr~(\alpha)=a+d$.  
The first equation implies $\alpha=a+b\tau$;  since both $\alpha$ and
$\tau$ are algebraic integers, one concludes that $b=1$. 
In view of Corollary \ref{cor6.3.1},  in a base $\{\omega_1',\omega_2'\}$,
the $\alpha$ has a matrix form $(a+d, c-ad, 1, 0)$. 
To calculate a real quadratic $\omega\in {\goth R}$ corresponding to $\alpha$,
recall an explicit formula obtained  in the proof of  Lemma \ref{lem6.1.5};
namely,   each endomorphism $(a,b,c,d)$ of the lattice $L_{CM}$  gives rise to  the 
endomorphism  $(a,b,-c,-d)$ of  pseudo-lattice ${\goth m}_{RM}=F(L_{CM})$.
Thus,   one gets a map:
\displaymath
F:  ~\left(\matrix{a+d & c-ad\cr 1 & 0}\right)
\longmapsto
 \left(\matrix{a+d & c-ad\cr -1 & 0}\right). 
\enddisplaymath
In other words, for a given base $\{\lambda_1,\lambda_2\}$ of
the pseudo-lattice ${\Bbb Z}+{\Bbb Z}\theta$ one can write
\displaymath
\left\{
\begin{array}{cc}
\omega\lambda_1 &= (a+d)\lambda_1 +(c-ad)\lambda_2\\
\omega\lambda_2 &= -\lambda_1.
\end{array}
\right.
\enddisplaymath
It is  easy  to verify,  that $\omega$ is a real  quadratic integer with 
$tr~(\omega)=a+d$. The latter coincides  with the $tr~(\alpha)$.
Lemma \ref{lem6.3.3} follows.
$\square$

 \index{pseudo-lattice}

\bigskip\noindent
Let $\omega\in {\goth R}$ be an endomorphism of the pseudo-lattice
${\goth m}_{RM}={\Bbb Z}+{\Bbb Z}\theta$ of degree $deg~(\omega):=\omega\bar\omega=n$.
The endomorphism maps ${\goth m}_{RM}$  to a sub-lattice   ${\goth m}_0\subset {\goth m}_{RM}$  of index $n$;
any such has the form ${\goth m}_0={\Bbb Z}+(n\theta){\Bbb Z}$, see e.g.  
[Borevich \& Shafarevich 1966]  \cite{BS},  p.131. 
Moreover,  $\omega$ generates an automorphism, $u$,  of the pseudo-lattice
${\goth m}_0$;  the traces of $\omega$ and $u$  are related.
\begin{lem}\label{lem6.3.4}
$tr~(u)=tr~(\omega)$. 
\end{lem}
{\it Proof.}
Let us calculate the action of endomorphism 
$\omega=(a+d, c-ad, -1, 0)$ on the pseudo-lattice 
${\goth m}_0={\Bbb Z}+(n\theta){\Bbb Z}$.
Since $deg~(\omega)=c-ad=n$,  one gets
\displaymath
\left(\matrix{a+d & n\cr -1 & 0}\right)
\left(\matrix{1\cr \theta}\right)=
\left(\matrix{a+d & 1\cr -1 & 0}\right)
\left(\matrix{1\cr n\theta}\right),
\enddisplaymath
where $\{1,\theta\}$ and $\{1,n\theta\}$ are bases of the pseudo-lattices 
${\goth m}_{RM}$ and ${\goth m}_0$, respectively,  and $u=(a+d, 1, -1, 0)$
is an automorphism of ${\goth m}_0$.   It is easy to see,  that
$tr~(u)=a+d=tr~(\omega)$.   
Lemma \ref{lem6.3.4} follows.
$\square$

\begin{rmk}\label{rm1}
{\bf (Second proof of Lemma \ref{lem6.3.4})}
{\normalfont
There exists a canonical proof of Lemma \ref{lem6.3.4}  based on the notion of
a subshift of finite type [Wagoner 1999]  \cite{Wag1};  we shall give such a proof  below,  since
it generalizes to   pseudo-lattices of any rank.   
Consider a dimension group ([Blackadar 1986]  \cite{B},  p.55) corresponding to the endomorphism $\omega$
of lattice ${\Bbb Z}^2$,  i.e. the limit $G(\omega)$:
\displaymath
{\Bbb Z}^2 \buildrel\omega \over\to 
{\Bbb Z}^2 \buildrel\omega \over\to 
{\Bbb Z}^2 \buildrel\omega \over\to
\dots
\enddisplaymath
 It is known that $G(\omega)\cong {\Bbb Z}[{1\over\lambda}]$, where $\lambda>1$ is the Perron-Frobenius
 eigenvalue of $\omega$.  We shall write $\hat\omega$ to denote the shift automorphism of dimension
 group $G(\omega)$, ([Effros  1981]  \cite{E}, p. 37)  and 
 $\zeta_{\omega}(t)=\exp\left(\sum_{k=1}^{\infty}{tr~(\omega^k)\over k}t^k\right)$ and 
 $\zeta_{\hat\omega}(t)=\exp\left(\sum_{k=1}^{\infty}{tr~(\hat\omega^k)\over k}t^k\right)$ the 
 corresponding Artin-Mazur zeta functions \cite{Wag1},  p. 273.     
 Since the Artin-Mazur zeta function of the  subshift of finite type  is an invariant of shift equivalence,
 we conclude that  $\zeta_{\omega}(t)\equiv\zeta_{\hat\omega}(t)$;  in particular, $tr~(\omega)=tr~(\hat\omega)$.  
 Hence the matrix form of $\hat\omega=(a+d, 1, -1, 0)=u$ and, therefore, $tr~(u)=tr~(\omega)$. 
Lemma \ref{lem6.3.4} is proved by a different method. 
$\square$
}
\end{rmk}
 \index{Gr\"ossencharacter}
\begin{lem}\label{lem6.3.5}
The automorphism $u$ is a unit of the ring ${\goth R}_0:=End~({\goth m}_0)$;
it is the fundamental unit of ${\goth R}_0$,  whenever $n=p$ is a prime number
and  $tr~(u)=tr~(\psi_{{\cal E}(K)}({\goth P}))$,  where  $(\psi_{{\cal E}(K)}({\goth P}))$
is the Gr\"ossencharacter associated to prime $p$, see  Supplement 6.3.3.  
\end{lem}
{\it Proof.}
(i)  Since $deg~(u)=1$,  the element $u$ is invertible and,  therefore,
a unit of the ring ${\goth R}_0$;  in general, unit  $u$ is not the fundamental
unit of ${\goth R}_0$,  since it is possible that  $u=\varepsilon^a$,  where
$\varepsilon$ is another unit of ${\goth R}_0$ and $a\ge 1$.

\medskip
(ii)  When $n=p$ is a prime number,  then we let   $\psi_{{\cal E}(K)}({\goth P})$
be the corresponding Gr\"ossencharacter on $K$ attached to an elliptic
curve ${\cal E}_{CM}\cong {\cal E}(K)$,  see Supplement 6.3.3  for the notation. 
The Gr\"ossencharacter can be identified with a complex number $\alpha\in k$
of the imaginary quadratic field $k$ associated to the complex multiplication. 
Let  $tr~(u)=tr~(\psi_{{\cal E}(K)}({\goth P}))$ and suppose to the contrary, that
$u$ is not the fundamental unit of ${\goth R}_0$,  i.e.  $u=\varepsilon^a$
for a unit $\varepsilon\in {\goth R}_0$ and an integer $a\ge 1$. 
Then there exists a Gr\"ossencharacter $\psi^{\prime}_{{\cal E}(K)}({\goth P})$,
such that 
\displaymath
tr~(\psi^{\prime}_{{\cal E}(K)}({\goth P}))< tr~(\psi_{{\cal E}(K)}({\goth P})). 
\enddisplaymath
Since $tr~(\psi_{{\cal E}(K)}({\goth P}))=q_{\goth P}+1-\# \tilde E({\Bbb F}_{\goth P})$,
one concludes that  $\# \tilde E({\Bbb F}^{\prime}_{\goth P})  > \# \tilde E({\Bbb F}_{\goth P})$;
in other words,  there exists a non-trivial extension 
${\Bbb F}^{\prime}_{\goth P}\supset {\Bbb F}_{\goth P}$
of the finite field  ${\Bbb F}_{\goth P}$.  The latter is impossible,
since any extension of  ${\Bbb F}_{\goth P}$ has the form ${\Bbb F}_{{\goth P}^n}$
for some $n\ge 1$;  thus $a=1$, i.e. unit $u$ is the fundamental unit 
of the ring ${\goth R}_0$.  Lemma \ref{lem6.3.5} is proved. 
$\square$

 \index{Hasse lemma}

\begin{lem}\label{lem6.3.6}
$tr~(\psi_{{\cal E}(K)}({\goth P}))=tr~(A^{\pi(p)})$. 
\end{lem}
{\it Proof.}
Recall that  the fundamental unit of the order ${\goth R}_0$ is given by 
the formula $\varepsilon_p=\varepsilon^{\pi(p)}$,  where $\varepsilon$ is 
the fundamental unit of the ring $O_{\goth k}$  and $\pi(p)$ an integer number,
see  Hasse's  Lemma \ref{lem6.3.10}  of  Supplement 6.3.3.    
On the other hand,  matrix $A=\prod_{i=1}^n(a_i, 1, 1, 0)$,
where $\theta=\overline{(a_1,\dots,a_n)}$ is a purely periodic
continued fraction.  Therefore
\displaymath
A\left(\matrix{1\cr\theta}\right)=\varepsilon\left(\matrix{1\cr\theta}\right), 
\enddisplaymath
where $\varepsilon>1$ is the fundamental unit of the real
quadratic field ${\goth k}={\Bbb Q}(\theta)$.  In other words,
$A$ is the matrix form of  the fundamental unit $\varepsilon$. 
Therefore  the matrix form of the fundamental unit  $\varepsilon_p=\varepsilon^{\pi(p)}$ 
of ${\goth R}_0$ is given by matrix   $A^{\pi(p)}$.
One can apply  Lemma \ref{lem6.3.5} and get  
\displaymath
tr~(\psi_{{\cal E}(K)}({\goth P}))=tr~(\varepsilon_p)=tr~(A^{\pi(p)}).
\enddisplaymath
Lemma \ref{lem6.3.6} follows.
$\square$

\bigskip\noindent
One can finish the proof of Theorem \ref{thm6.3.1}  by comparing the local $L$-series 
of the Hasse-Weil $L$-function for the ${\cal E}_{CM}$ with  that of the local zeta  for
the ${\cal A}_{RM}$. 
The local $L$-series for ${\cal E}_{CM}$ are
$L_{\goth P}({\cal E}(K),T)=1-a_{\goth P}T+q_{\goth P}T^2$ if the ${\cal E}_{CM}$ has a good reduction at ${\goth P}$
and $L_{\goth P}({\cal E}(K),T)=1-\alpha T$ otherwise; here 
\displaymath
\left\{
\begin{array}{lll}
q_{\goth P} &=& N^K_{\Bbb Q}{\goth P}=\# {\Bbb F}_{\goth P}=p,\\
a_{\goth P} &=& q_{\goth P}+1-\# \tilde E({\Bbb F}_{\goth P})=tr~(\psi_{{\cal E}(K)}({\goth P})),\\
\alpha &\in & \{-1,0,1\}.
\end{array}
\right.
\enddisplaymath
Therefore, 
\displaymath
L_{\goth P}({\cal E}_{CM},T)=
\cases{
1-tr~(\psi_{{\cal E}(K)}({\goth P}))T+pT^2, & \mbox{for good reduction }  \cr
1-\alpha T,  & \mbox{for bad reduction}.
}
\enddisplaymath
\begin{lem}\label{lem6.3.7}
For ${\cal A}_{RM}=F({\cal E}_{CM})$,  it holds 
$$\zeta_p^{-1}({\cal A}_{RM}, T)=1-tr~(A^{\pi(p)})T+pT^2,$$ 
 whenever $p\nmid tr^2(A)-4$.
\end{lem}
{\it Proof.}   By the formula   $K_0({\cal O}_{B})={\Bbb Z}^2/(I-B^t){\Bbb Z}^2$,
one gets
\displaymath
|K_0({\cal O}_{L_p^n})|=\left|{{\Bbb Z}^2\over (I-(L_p^n)^t){\Bbb Z}^2}\right|=
|det(I-(L_p^n)^t)|=|Fix ~(L_p^n)|,
\enddisplaymath
where $Fix~(L_p^n)$ is the set of (geometric) fixed points of the
endomorphism $L_p^n: {\Bbb Z}^2\to {\Bbb Z}^2$.    Thus,
\displaymath
\zeta_p({\cal A}_{RM}, z)= \exp \left(\sum_{n=1}^{\infty}{|Fix~(L_p^n)|\over n} ~z^n\right), 
\quad z\in {\Bbb C}. 
\enddisplaymath
But the latter series is an Artin-Mazur zeta function of the endomorphism $L_p$;
it converges to a rational function $det^{-1}(I-zL_p)$, see e.g.    [Hartshorn  1977]  \cite{H1},  p.455. 
Thus,  $\zeta_p({\cal A}_{RM}, z)= det^{-1}(I-zL_p)$. 
The substitution $L_p=(tr~(A^{\pi(p)}), p, -1, 0)$ gives us
\displaymath
det~(I-zL_p)=
det~\left(\matrix{1-tr~(A^{\pi(p)})z & -pz\cr z & 1}\right)=1-tr~(A^{\pi(p)})z+pz^2.
\enddisplaymath
Put $z=T$ and get $\zeta_p({\cal A}_{RM}, T)=(1-tr~(A^{\pi(p)})T+pT^2)^{-1}$,
which is a conclusion of  Lemma \ref{lem6.3.7}.
$\square$

\begin{lem}\label{lem6.3.8}
For ${\cal A}_{RM}=F({\cal E}_{CM})$,  it holds 
$$\zeta_p^{-1}({\cal A}_{RM}, T)=1-\alpha T,$$ 
whenever $p ~| ~tr^2(A)-4$.
\end{lem}
{\it Proof.}
Indeed, $K_0({\cal O}_{1-\alpha^n})={\Bbb Z}/(1-1+\alpha^n){\Bbb Z}={\Bbb Z}/\alpha^n{\Bbb Z}$.
Thus, $|K_0({\cal O}_{1-\alpha^n})|=det~(\alpha^n)=\alpha^n$. By the definition,
\displaymath
\zeta_p({\cal A}_{RM},z)=\exp\left(\sum_{n=1}^{\infty} {\alpha^n\over n}z^n\right)=
\exp\left(\sum_{n=1}^{\infty} {(\alpha z)^n\over n}\right)=
{1\over 1-\alpha z}.
\enddisplaymath
The substitution $z=T$ gives the conclusion of  Lemma \ref{lem6.3.8}.
$\square$

\begin{lem}\label{lem6.3.9}
Let ${\goth P}\subset K$ be a prime ideal over $p$; 
then ${\cal E}_{CM}={\cal E}(K)$ has a bad reduction at ${\goth P}$
if and only if  $p~|~tr^2(A)-4$. 
\end{lem}
{\it Proof.}
Let $k$ be a  field of complex 
multiplication of the ${\cal E}_{CM}$; its discriminant we
shall write as $\Delta_k<0$. It is known, that whenever $p~|~\Delta_k$,
the ${\cal E}_{CM}$ has a bad reduction at the prime ideal ${\goth P}$ over $p$.
On the other hand,   the explicit formula for functor $F$ applied to the matrix $L_p$
gives us $F: (tr~(A^{\pi(p)}),  p,  -1, 0) \mapsto (tr~(A^{\pi(p)}), p,  1, 0)$,
see proof of Lemma \ref{lem6.1.5}. 
The characteristic polynomials of the above matrices
are  $x^2-tr~(A^{\pi(p)})x+p$ and $x^2-tr~(A^{\pi(p)})x-p$, respectively.   
They generate an imaginary (resp., a real) quadratic field 
$k$ (resp., ${\goth k}$) with the discriminant $\Delta_k=tr^2 (A^{\pi(p)})-4p<0$
(resp., $\Delta_{\goth k}=tr^2 (A^{\pi(p)})+4p>0$). Thus,
$\Delta_{\goth k}-\Delta_k=8p$. It is easy to see, that 
$p~|~\Delta_{\goth k}$ if and only if  $p~|~\Delta_k$.
It remains to express the discriminant $\Delta_{\goth k}$ in terms
of the matrix $A$. Since the characteristic polynomial for $A$
is $x^2-tr~(A)x+1$, it follows  that $\Delta_{\goth k}=tr^2(A)-4$.
Lemma \ref{lem6.3.9} follows.  
$\square$

\bigskip\noindent
Let us   prove that the  first part of Theorem \ref{thm6.3.1} implies the first claim of its second part;
notice,  that the critical piece of information is   provided by 
Lemma \ref{lem6.3.6},  which says that $tr~(\psi_{{\cal E}(K)}({\goth P}))=tr~(A^{\pi(p)})$.
Thus, Lemmas \ref{lem6.3.7}--\ref{lem6.3.9} imply that  $L_{\goth P}({\cal E}_{CM},T)\equiv \zeta_p^{-1}({\cal A}_{RM}, T)$.
The first claim of part (ii)  of  Theorem \ref{thm6.3.1} follows.

\bigskip\noindent
\underline{{\sf  A. Let $p$ be a good prime.}}  Let us prove the second claim of part (ii)  of 
Theorem \ref{thm6.3.1} in the case  $n=1$.   From the left side: 
$K_0({\cal A}_{RM}\rtimes_{L_p}{\Bbb Z})\cong K_0({\cal O}_{L_p})\cong
{\Bbb Z}^2/(I-L_p^t){\Bbb Z}^2$,  where $L_p=(tr~(A^{\pi(p)}), p, -1, 0)$. To
calculate the abelian group ${\Bbb Z}^2/(I-L_p^t){\Bbb Z}^2$, we shall 
use a reduction of the matrix $I-L_p^t$ to the Smith normal form:
\displaymath
I-L_p^t=
\left(\matrix{1-tr~(A^{\pi(p)}) & 1\cr -p & 1}\right)\sim
\left(\matrix{1+p-tr~(A^{\pi(p)}) & 0\cr -p & 1}\right)\sim
$$
$$
\sim\left(\matrix{1 & 0\cr 0 & 1+p-tr~(A^{\pi(p)})}\right).  
\enddisplaymath
Therefore, $K_0({\cal O}_{L_p})\cong {\Bbb Z}_{1+p-tr~(A^{\pi(p)})}$. 
From the right side, the ${\cal E}_{CM}({\Bbb F}_{\goth P})$ is an elliptic
curve over the field of characteristic $p$. Recall,  that the chord and tangent law
turns the ${\cal E}_{CM}({\Bbb F}_{\goth P})$ into a finite abelian group. The group
is cyclic and has the order $1+q_{\goth P}-a_{\goth P}$.
But $q_{\goth P}=p$ and $a_{\goth P}= tr~(\psi_{{\cal E}(K)}({\goth P}))=tr~(A^{\pi(p)})$,
see Lemma \ref{lem6.3.6}.   Thus, ${\cal E}_{CM}({\Bbb F}_{\goth P})\cong {\Bbb Z}_{1+p-tr~(A^{\pi(p)})}$;
therefore  $K_0({\cal O}_{L_p})\cong {\cal E}_{CM}({\Bbb F}_p)$. 
The general case $n\ge 1$ is treated likewise. Repeating  the argument of Lemmas \ref{lem6.3.1}--\ref{lem6.3.2},
it follows that 
\displaymath
L_p^n=\left(\matrix{tr~ (A^{n\pi(p)}) & p^n\cr  -1 & 0}\right).
\enddisplaymath
Then one gets  $K_0({\cal O}_{L_p^n})\cong {\Bbb Z}_{1+p^n-tr~(A^{n\pi(p)})}$
on the left side. From the right side, $|{\cal E}_{CM}({\Bbb F}_{p^n})|=1+p^n-tr~(\psi^n_{{\cal E}(K)}({\goth P}))$;
but a repetition of the argument of  Lemma \ref{lem6.3.6}  yields us $tr~(\psi^n_{{\cal E}(K)}({\goth P}))=tr~(A^{n\pi(p)})$.
Comparing the left and right sides, one gets that $K_0({\cal O}_{L_p^n})\cong {\cal E}_{CM}({\Bbb F}_{p^n})$.
This argument finishes the proof of the second claim of part (ii)  of  Theorem \ref{thm6.3.1} for the good primes.

\bigskip\noindent
\underline{{\sf B. Let $p$ be a bad prime.}}
From the proof of Lemma \ref{lem6.3.8},  one gets for the left side  $K_0({\cal O}_{\varepsilon_n})\cong {\Bbb Z}_{\alpha^n}$.
From the right side, it holds    $|{\cal E}_{CM}({\Bbb F}_{p^n})|=1+q_{\goth P}-a_{\goth P}$,
where $q_{\goth P}=0$ and $a_{\goth P}=tr~(\varepsilon_n)=\varepsilon_n$. 
Thus,   $|{\cal E}_{CM}({\Bbb F}_{p^n})|=1-\varepsilon_n=1-(1-\alpha^n)=\alpha^n$. 
Comparing the left and right sides, we conclude that
$K_0({\cal O}_{\varepsilon_n})\cong {\cal E}_{CM}({\Bbb F}_{p^n})$ at the bad primes.

\bigskip
All  cases are exhausted;  thus part (i) of    Theorem \ref{thm6.3.1}  implies its part (ii). 
The proof of  converse consists in a step by step claims   similar to just     proved
and is left to the reader.  Theorem \ref{thm6.3.1} is proved. 
$\square$

 \index{Gr\"ossencharacter}
 \index{unit of algebraic number field}

\subsection{Supplement: Gr\"ossencharacters, units and $\pi(n)$ }
We shall briefly review  the well known facts about  complex multiplication
and  units in  subrings of the ring of integers  in algebraic number fields;
for the detailed account,  we refer the reader to [Silverman 1994]  \cite{S2} and 
[Hasse 1950] \cite{HA},   respectively.

\subsubsection{Gr\"ossencharacters}
Let ${\cal E}_{CM}\cong {\cal E}(K)$ be elliptic curve with complex
multiplication and $K\cong k(j({\cal E}_{CM}))$ the Hilbert class field  attached to ${\cal E}_{CM}$.
For each prime ideal ${\goth P}$ of $K$,  let ${\Bbb F}_{\goth P}$
be a residue field of $K$ at ${\goth P}$ and 
$q_{\goth P}=N^K_{\Bbb Q}{\goth P}=\# {\Bbb F}_{\goth P}$,
where $N^K_{\Bbb Q}$ is the norm of the ideal ${\goth P}$.
If ${\cal E}(K)$ has a good reduction at ${\goth P}$, 
one defines  $a_{\goth P}=q_{\goth P}+1-\# \tilde {\cal E}({\Bbb F}_{\goth P})$,
where $\tilde {\cal E}({\Bbb F}_{\goth P})$ is a reduction of ${\cal E}(K)$ 
modulo the prime ideal ${\goth P}$.  If ${\cal E}(K)$ has good reduction at ${\goth P}$,  the polynomial
\displaymath
L_{\goth P}({\cal E}(K),T)=1-a_{\goth P}T+q_{\goth P}T^2,
\enddisplaymath
is called the {\it local $L$-series} of ${\cal E}(K)$ at ${\goth P}$.
If ${\cal E}(K)$ has bad reduction at ${\goth P}$, the local $L$-series
are $L_{\goth P}({\cal E}(K),T)=1-T$ (resp. $L_{\goth P}({\cal E}(K),T)=1+T$; $L_{\goth P}({\cal E}(K),T)=1$)
if ${\cal E}(K)$ has split multiplicative reduction at ${\goth P}$ (if ${\cal E}(K)$ has non-split
multiplicative reduction at ${\goth P}$; if ${\cal E}(K)$ has additive reduction at ${\goth P}$). 
\begin{dfn}
By the Hasse-Weil $L$-function of elliptic curve ${\cal E}(K)$ one
understands  the global $L$-series defined by the Euler product
\displaymath
L({\cal E}(K),s)=\prod_{\goth P} ~[L_{\goth P}({\cal E}(K), q_{\goth P}^{-s})]^{-1}.
\enddisplaymath
\end{dfn}
\begin{dfn}
If  $A_K^*$ be the idele group of the number field $K$,   then by a Gr\"ossen\-char\-ac\-ter  
on $K$ one understands   a continuous homomorphism 
\displaymath
\psi: A_K^*\longrightarrow  {\Bbb C}^*
\enddisplaymath
 with the property $\psi(K^*)=1$;  the asterisk denotes the group of invertible elements
of the corresponding ring.  The  Hecke $L$-series attached
to the Gr\"ossencharacter $\psi: A_K^*\to {\Bbb C}^*$ is defined 
by the Euler product
\displaymath
L(s,\psi)=\prod_{\goth P} (1-\psi({\goth P})q_{\goth P}^{-s})^{-1},
\enddisplaymath
where the product is taken over all prime ideals of $K$.
\end{dfn}
\begin{rmk}
\textnormal{
For a prime ideal  ${\goth P}$  of field  $K$ at which ${\cal E}(K)$
has  good reduction  and  $\tilde {\cal E}({\Bbb F}_{\goth P})$ being  the reduction 
of ${\cal E}(K)$  at ${\goth P}$,  we let 
\displaymath
\phi_{\goth P}: \tilde {\cal E}({\Bbb F}_{\goth P})\longrightarrow \tilde {\cal E}({\Bbb F}_{\goth P}) 
\enddisplaymath
denote  the associated {\it Frobenius map};    if  $\psi_{{\cal E}(K)}: A_K^*\to k^*$ is  the 
Gr\"ossen\-char\-ac\-ter  attached  to the ${\cal E}_{CM}$,  then
the diagram in Fig. 6.5 is known to be commutative,   see [Silverman 1994]  \cite{S2},  p.174. 
In particular, $\psi_{{\cal E}(K)}({\goth P})$  is an endomorphism of the ${\cal E}(K)$ given by the complex number
$\alpha_{{\cal E}(K)}({\goth P})\in R$,  where $R={\Bbb Z}+fO_k$ is an order in imaginary quadratic field $k$ . 
If  $\overline{\psi}_{{\cal E}(K)}({\goth P})$  is  the conjugate Gr\"ossencharacter viewed as a complex number,
then  the {\it Deuring Theorem} says that the Hasse-Weil $L$-function of the ${\cal E}(K)$
is related to the Hecke $L$-series of the $\psi_{{\cal E}(K)}$ by the formula
\displaymath
L({\cal E}(K),s)\equiv L(s,\psi_{{\cal E}(K)})L(s, \overline{\psi}_{{\cal E}(K)}).
\enddisplaymath
}
\end{rmk}

 \index{Deuring Theorem}

\bigskip
\begin{figure}[here]
\begin{picture}(300,110)(-120,-5)
\put(20,70){\vector(0,-1){35}}
\put(130,70){\vector(0,-1){35}}
\put(50,23){\vector(1,0){53}}
\put(50,83){\vector(1,0){53}}
\put(15,20){$\tilde {\cal E}({\Bbb F}_{\goth P})$}
\put(123,20){$\tilde {\cal E}({\Bbb F}_{\goth P})$}
\put(15,80){${\cal E}(K)$}
\put(115,80){${\cal E}(K)$}
\put(70,30){$\phi_{\goth P}$}
\put(55,90){$\psi_{{\cal E}(K)}({\goth P})$}
\end{picture}
\caption{The Gr\"ossencharacter  $\psi_{{\cal E}(K)}({\goth P})$.}
\end{figure}

 \index{function $\pi(n)$}

\subsubsection{Units and function $\pi(n)$}
Let ${\goth k}={\Bbb Q}(\sqrt{D})$ be a real quadratic number field
and $O_{\goth k}$ its ring of integers.  For rational integer $n\ge 1$
we shall write ${\goth R}_n\subseteq O_{\goth k}$  to denote an
order (i.e. a subring containing  $1$) of $O_{\goth k}$.  The order
${\goth R}_n$ has a basis $\{1, n\omega\}$, where 
\displaymath
\omega=\cases{{\sqrt{D}+1\over 2} & if $D\equiv 1 ~ mod~4$,\cr
               \sqrt{D} & if $D\equiv 2,3 ~ mod~4$.}
\enddisplaymath
In other words, ${\goth R}_n={\Bbb Z}+(n\omega){\Bbb Z}$. 
It is clear, that ${\goth R}_1=O_{\goth k}$ and the fundamental 
unit of $O_{\goth k}$ we shall denote by $\varepsilon$. 
Each ${\goth R}_n$ has its own fundamental unit, which
we shall write as $\varepsilon_n$;  notice that $\varepsilon_n\ne
\varepsilon$  unless $n= 1$. 
There exists the well-known  formula,  which relates $\varepsilon_n$
to the fundamental unit $\varepsilon$,  see e.g. [Hasse  1950]  \cite{HA}, p.297.
Denote by ${\goth G}_n:=U(O_{\goth k}/nO_{\goth k})$ 
the multiplicative group of invertible elements (units) of 
the residue ring  $O_{\goth k}/nO_{\goth k}$;   clearly, all units 
of $O_{\goth k}$ map (under the natural $mod~n$ homomorphism) to ${\goth G}_n$. 
Likewise  let ${\goth g}_n:=U({\goth R}_n/n{\goth R}_n)$
be the group of units of the residue ring ${\goth R}_n/n{\goth R}_n$;
it is not hard to prove ([Hasse  1950]  \cite{HA}, p.296), that ${\goth g}_n\cong U({\Bbb Z}/n{\Bbb Z})$
the ``rational'' unit group of the residue ring ${\Bbb Z}/n{\Bbb Z}$. 
Similarly, all units of the order ${\goth R}_n$ map to ${\goth g}_n$.  
Since units of ${\goth R}_n$ are also units of $O_{\goth k}$
(but not vice versa),  ${\goth g}_n$ is a subgroup of ${\goth G}_n$;
in particular,   $|{\goth G}_n|/|{\goth g}_n|$ is an integer number
and $|{\goth g}_n|=\varphi(n)$, where $\varphi(n)$ is the Euler totient function. 
In general,  the following formula is true
\displaymath
{|{\goth G}_n|\over |{\goth g}_n|}=n\prod_{p_i | n}
\left(1-\left({D\over p_i}\right){1\over p_i}\right),
\enddisplaymath
where $\left({D\over p_i}\right)$ is the Legendre symbol, 
see [Hasse  1950]  \cite{HA}, p. 351. 
\begin{dfn}
By the  function $\pi(n)$ one understands   the least integer number  dividing 
$|{\goth G}_n|/|{\goth g}_n|$ and   such that
$\varepsilon^{\pi(n)}$ is  a unit of ${\goth R}_n$, i.e. belongs to ${\goth g}_n$.
\end{dfn}
\begin{lem}\label{lem6.3.10}
{\bf ([Hasse  1950]  \cite{HA}, p.298)}
$\varepsilon_n=\varepsilon^{\pi(n)}.$
\end{lem}
\begin{rmk}
\textnormal{
 Lemma \ref{lem6.3.10} asserts  existence  of the number $\pi(n)$ 
as one of the divisors of  $|{\goth G}_n|/|{\goth g}_n|$,  yet no analytic 
formula for $\pi(n)$ is known;   it would be rather interesting to have such
a formula. 
}
\end{rmk}
\begin{rmk}
\textnormal{
In the special case $n=p$  is a prime number, the following formula is true
\displaymath
{|{\goth G}_p|\over |{\goth g}_p|}=p-\left({D\over p}\right).
\enddisplaymath
}
\end{rmk}

\vskip1.5cm\noindent
{\bf Guide to the literature.} 
The Hasse-Weil $L$-functions $L({\cal E}(K), s)$ of the $K$-rational 
elliptic curves are covered in the textbooks by   [Husem\"oller 1986]  \cite{H2},
[Knapp 1992] \cite{K1} and  [Silverman 1994]   \cite{S2};  see also the survey
[Tate 1974]  \cite{Tat1}.  The reciprocity of  $L({\cal E}(K), s)$ with an $L$-function
obtained from certain cusp form of weight two is  subject of the {\it Eichler-Shimura theory},
see e.g.    [Knapp 1992]  \cite{K1},  Chapter XI;  such a reciprocity coupled with 
the {\it Shimura-Taniyama Conjecture}  was critical to solution of the Fermat Last 
Theorem by A.~Wiles.  The {\it non-commutative reciprocity}  of   $L({\cal E}(K), s)$
with an $L$-function obtained from a noncommutative torus with real multiplication
was proved in \cite{Nik13}.    

 \index{Shimura-Taniyama Conjecture}

 \index{Langlands program}

\section{Langlands program  for noncommutative tori}
We dealt with functors on the arithmetic schemes $X$ so far.   In this section 
we shall define a functor $F$ on the category of all finite Galois extensions $E$
of the field ${\Bbb Q}$;  the functor ranges in a category of the even-dimensional 
noncommutative tori with real multiplication.  
For such a torus,  ${\cal A}_{RM}^{2n}$, we construct an $L$-function 
$L({\cal A}_{RM}^{2n}, s)$;  it is conjectured that for each $n\ge 1$ and each
irreducible representation 
\displaymath
\sigma: Gal~(E|{\Bbb Q})\longrightarrow GL_n({\Bbb C}),
\enddisplaymath
the corresponding {\it Artin $L$-function} $L(\sigma, s)$ coincides with   
$L({\cal A}_{RM}^{2n}, s)$,   whenever ${\cal A}_{RM}^{2n}=F(E)$.   
Our main result Theorem \ref{thm6.4.1} says  that the conjecture is true  
for $n=1$ (resp., $n=0$) and $E$ being the Hilbert class field of an
imaginary quadratic field $k$ (resp., the rational field ${\Bbb Q}$) .   
Thus we are dealing  with an analog of the 
{\it Langlands program},   where the ``automorphic cuspidal representations of group
$GL_n$'' are replaced by the noncommutative tori ${\cal A}_{RM}^{2n}$, 
see [Gelbart 1984]  \cite{Gel1} for an  introduction.

 \index{function $L({\cal A}_{RM}^{2n},s)$}

\subsection{$L({\cal A}_{RM}^{2n}, s)$}
The higher-dimensional noncommutative tori were introduced in Section 3.4.1;
let us recall some notation.  Let $\Theta=(\theta_{ij})$ be a real skew symmetric  matrix of even dimension $2n$;
by ${\cal A}_{\Theta}^{2n}$  we shall mean  the even-dimensional noncommutaive
torus defined by matrix $\Theta$,   i.e.  a  universal $C^*$-algebra on the unitary generators $u_1,\dots, u_{2n}$
and  relations  
\displaymath
u_ju_i=e^{2\pi i\theta_{ij}}u_iu_j, \quad 1\le i,j\le 2n.    
\enddisplaymath
It is known, that by the orthogonal linear transformations 
every (generic)  real even-dimensional skew symmetric matrix can be brought to
 the  normal form
\displaymath
\Theta_0=
\left(
\matrix{      0 & \theta_1 &          &            &   \cr
      -\theta_1 &   0      &          &            &   \cr
                &          &  \ddots  &            &   \cr
                &          &          & 0          & \theta_n  \cr
                &          &          & -\theta_n  &    0
}\right)
\enddisplaymath
where $\theta_i>0$ are linearly independent over ${\Bbb Q}$. 
We shall consider the noncommutative tori ${\cal A}_{\Theta_0}^{2n}$,
given by matrix in the above normal form;   we refer to the family 
as  a {\it normal family}. 
Recall  that $K_0({\cal A}_{\Theta_0}^{2n})\cong {\Bbb Z}^{2^{2n-1}}$
and  the positive cone $K_0^+({\cal A}_{\Theta_0}^{2n})$ is
given by the pseudo-lattice 
\displaymath
{\Bbb Z}+\theta_1{\Bbb Z}+\dots+\theta_n{\Bbb Z}
+\sum_{i=n+1}^{2^{2n-1}}p_i(\theta){\Bbb Z}  ~\subset ~{\Bbb R},
\enddisplaymath
where $p_i(\theta)\in {\Bbb Z}[1, \theta_1,\dots,\theta_n]$,   see e.g.  
[Elliott  1982]  \cite{Ell2}. 
\begin{dfn}
The noncommutative torus ${\cal A}_{\Theta_0}^{2n}$ is said to
have real multiplication  if the endomorphism  ring $End~(K_0^+({\cal A}_{\Theta_0}^{2n}))$
is non-trivial,  i.e.  exceeds  the  ring ${\Bbb Z}$;   we shall denote such a torus by ${\cal A}_{RM}^{2n}$. 
\end{dfn}
\begin{rmk}\label{rmk6.4.1}
\textnormal{
It is easy to see that  if ${\cal A}_{\Theta_0}^{2n}$  has real multiplication,  then $\theta_i$
are algebraic integers;  we leave the proof to the reader.  (Hint:  each endomorphism
of   $K_0^+({\cal A}_{\Theta_0}^{2n})\cong {\Bbb Z}+\theta_1{\Bbb Z}+\dots+\theta_n{\Bbb Z}
+\sum_{i=n+1}^{2^{2n-1}}p_i(\theta){\Bbb Z}$ is multiplication by a real number;  
thus the endomorphism is described by an integer matrix,  which 
defines a polynomial equation involving $\theta_i$.)
}
\end{rmk}
\begin{rmk}\label{rmk6.4.2}
\textnormal{
Remark \ref{rmk6.4.1} says that $\theta_i$ are algebraic integers whenever
${\cal A}_{\Theta_0}^{2n}$ has real multiplication;   so will be the values 
of polynomials $p_i(\theta)$ in this case.  Since such values belong to the 
number field ${\Bbb Q}(\theta_1,\dots,\theta_n)$,  one concludes that 
\displaymath
K_0^+({\cal A}_{RM}^{2n})\cong
{\Bbb Z}+\theta_1{\Bbb Z}+\dots+\theta_n{\Bbb Z} ~\subset ~{\Bbb R}.
\enddisplaymath
}
\end{rmk}
Let $A\in GL_{n+1}({\Bbb Z})$ be a positive matrix such that
\displaymath
A\left(\matrix{1\cr\theta_1\cr\vdots\cr\theta_n}\right)=
\lambda_A \left(\matrix{1\cr\theta_1\cr\vdots\cr\theta_n}\right),
\enddisplaymath
where $\lambda_A$ is the Perron-Frobenius eigenvalue of $A$;  
in other words,  $A$ is a  matrix corresponding to the 
{\it shift automorphism $\sigma_A$}  of   $K_0^+({\cal A}_{RM}^{2n})$
regarded as a stationary dimension group,  see Definition \ref{dfn3.5.4}.   
For each prime number $p$,  consider the characteristic polynomial 
of matrix $A^{\pi(p)}$,  where $\pi(n)$ is the integer-valued function 
introduced in Section 6.3.3;  in other words,
\displaymath
\hbox{{\bf Char}} (A^{\pi(p)}):=\det (xI-A^{\pi(p)})=x^{n+1}-a_1x^n-\dots-a_nx-1 ~\in ~{\Bbb Z}[x]. 
\enddisplaymath
\begin{dfn}
By a local zeta function of  the noncommutative torus  ${\cal A}_{RM}^{2n}$ we 
understand the function
\displaymath
\zeta_p({\cal A}_{RM}^{2n}, z):= {1\over 1-a_1z+a_2z^2-\dots-a_n z^n +pz^{n+1}}, \quad z\in {\Bbb C}.
\enddisplaymath
\end{dfn}
\begin{rmk}\label{rmk6.4.3}
\textnormal{
To explain the structure of  $\zeta_p({\cal A}_{RM}^{2n}, z)$,  consider the 
{\it companion matrix}
\displaymath
J=\left(
\matrix{
a_1     &  1    & \dots  & 0 &  0\cr
a_2     & 0      & 1  &  0      & 0\cr
\vdots  & \vdots  & \ddots & \vdots  & \vdots\cr
a_n & 0  & \dots &  0 & 1 \cr
1      &  0     & \dots  &    0   & 0
}\right) 
\enddisplaymath
of polynomial {\bf Char}  $(A^{\pi(p)})=x^{n+1}-a_1x^n -\dots-a_nx-1$,
i.e. the matrix $J$ such that  $\det(xI-J)=x^{n+1}-a_1x^n-\dots-a_nx-1$.  
It is not hard to see, that the non-negative  integer matrix $J$ corresponds
to the {\it shift automorphism}  of a stationary dimension group
\displaymath
{\Bbb Z}^{n+1}\buildrel\rm
J_p
\over\longrightarrow {\Bbb Z}^{n+1}
\buildrel\rm
J_p
\over\longrightarrow
{\Bbb Z}^{n+1}\buildrel\rm
J_p
\over\longrightarrow \dots
\enddisplaymath
where 
\displaymath
J_p=\left(
\matrix{
a_1     &  1    & \dots  & 0 &  0\cr
a_2     & 0      & 1  &  0      & 0\cr
\vdots  & \vdots  & \ddots & \vdots  & \vdots\cr
a_n & 0  & \dots &  0 & 1 \cr
p      &  0     & \dots  &    0   & 0
}\right). 
\enddisplaymath
On the other hand,  the companion matrix of polynomial 
{\bf Char} $(\sigma(Fr_p))=\det(xI-\sigma(Fr_p))=x^{n+1}-a_1x^n+\dots-a_nx+p$
has the form
\displaymath
W_p=\left(
\matrix{
a_1     &  1    & \dots  & 0 &  0\cr
-a_2     & 0      & 1  &  0      & 0\cr
\vdots  & \vdots  & \ddots & \vdots  & \vdots\cr
a_n  & 0  & \dots &  0 & 1 \cr
-p      &  0     & \dots  &    0   & 0
}\right),
\enddisplaymath
see Section 6.4.3 for the meaning of $\sigma(Fr_p)$.  
Thus the action of functor $F:$ {\bf Alg-Num} $\to$ {\bf NC-Tor} on the 
corresponding companion matrices $W_p$ and $J_p$ is given by the formula
\displaymath
F:  ~\left(
\matrix{
a_1     &  1    & \dots  & 0 &  0\cr
-a_2     & 0      & 1  &  0      & 0\cr
\vdots  & \vdots  & \ddots & \vdots  & \vdots\cr
a_n & 0  & \dots &  0 & 1 \cr
-p      &  0     & \dots  &    0   & 0
}\right)
\mapsto
\left(
\matrix{
a_1     &  1    & \dots  & 0 &  0\cr
a_2     & 0      & 1  &  0      & 0\cr
\vdots  & \vdots  & \ddots & \vdots  & \vdots\cr
a_n  & 0  & \dots &  0 & 1 \cr
p      &  0     & \dots  &    0   & 0
}\right). 
\enddisplaymath
It remains to compare our formula for $\zeta_p({\cal A}_{RM}^{2n}, z)$ with the well-known 
formula for  the Artin zeta function 
\displaymath
\zeta_p(\sigma_n, z)={1\over\det(I_n-\sigma_n(Fr_p)z)},
\enddisplaymath
where $z=x^{-1}$,  see [Gelbart 1984] \cite{Gel1}, p. 181.      
}
\end{rmk}
\begin{dfn}
By an  $L$-function of  the noncommutative torus  ${\cal A}_{RM}^{2n}$ one  
understand the product 
\displaymath
L({\cal A}_{RM}^{2n},s):=\prod_p
~\zeta_p({\cal A}_{RM}^{2n}, p^{-s}), \quad s\in {\Bbb C},
\enddisplaymath
over all but a finite number of primes $p$.  
\end{dfn}
 \index{Langalands Conjecture for noncommutative torus}
\begin{cnj}\label{cnj6.4.1}
{\bf (Langlands conjecture for noncommutative tori)}
For each  finite extension $E$  of the field of rational numbers ${\Bbb Q}$
with the  Galois group $Gal~(E |{\Bbb Q})$ and each  
 irreducible representation 
 \displaymath
 \sigma_{n+1}: Gal~(E | {\Bbb Q})\to GL_{n+1}({\Bbb C}),
\enddisplaymath
there exists a $2n$-dimensional  noncommutative torus with 
real  multiplication  ${\cal A}_{RM}^{2n}$,   such that
\displaymath
L(\sigma_{n+1},s)\equiv L({\cal A}_{RM}^{2n}, s),
\enddisplaymath
where $L(\sigma_{n+1}, s)$ is the Artin $L$-function attached to representation 
$\sigma_{n+1}$   and  $L({\cal A}_{RM}^{2n}, s)$ is  the  $L$-function of  the
noncommutative torus  ${\cal A}_{RM}^{2n}$.
\end{cnj}
\begin{rmk}
\textnormal{
Roughly speaking,   Conjecture \ref{cnj6.4.1}  says  that the Galois extensions (abelian or not) 
of the field  ${\Bbb Q}$ are in a one-to-one correspondence with  the even-dimensional 
noncommutative  tori with real multiplication.   In the context of the Langlands program,   
the noncommutative torus ${\cal A}_{RM}^{2n}$  can  be regarded as 
an  analog   of  the ``automorphic cuspidal representation  $\pi_{\sigma_{n+1}}$ of the group 
$GL(n+1)$''.    This  appearance of  ${\cal A}_{RM}^{2n}$ is {\it not}  random
because the noncommutative tori classify the irreducible infinite-dimensional representations of 
the  Lie group $GL(n+1)$,  see the remarkable paper  by  [Poguntke 1983]   \cite{Pog1};    
such representations are known to be at the heart  of the Langlands  philosophy,  
see  [Gelbart 1984] \cite{Gel1} . 
}
\end{rmk}
\begin{thm}\label{thm6.4.1}
Conjecture \ref{cnj6.4.1} is true for $n=1$ (resp., $n=0$) and $E$ abelian 
extension of an imaginary quadratic field $k$ (resp., the rational field ${\Bbb Q}$). 
\end{thm}

\subsection{Proof of  Theorem \ref{thm6.4.1}}
\subsubsection{Case $n=1$}
Roughly speaking,  this case is equivalent to  Theorem \ref{thm6.3.1};  
it was a model example for  Conjecture \ref{cnj6.4.1}.   
Using the Gr\"ossencharacters,  one can identify  the Artin $L$-function for 
abelian extensions of   the imaginary quadratic fields $k$ with the 
Hasse-Weil $L$-function $L({\cal E}_{CM}, s)$,
where ${\cal E}_{CM}$ is an elliptic curve with complex multiplication by $k$;
but Theorem \ref{thm6.3.1} says  that $L({\cal E}_{CM}, s)\equiv L({\cal A}_{RM}, s)$,  
where  $L({\cal A}_{RM}, s)$ is  the special case $n=1$ of our function  $L({\cal A}_{RM}^{2n}, s)$.

To  give  the  details,  let $k$ be an imaginary quadratic field and let ${\cal E}_{CM}$
be an elliptic curve with complex multiplication by (an order) in $k$. 
By the  theory of complex multiplication,  the Hilbert class field $K$ of $k$ 
is given by  the  $j$-invariant of ${\cal E}_{CM}$,  i.e.
\displaymath
K\cong k(j({\cal E}_{CM})), 
\enddisplaymath
and  $Gal~(K|k)\cong Cl~(k)$,  where $Cl~(k)$ is the ideal class group of $k$;  moreover,
\displaymath
{\cal E}_{CM}\cong {\cal E}(K),
\enddisplaymath
see e.g.   [Silverman 1994]  \cite{S2}.  Recall that functor $F:$ {\bf Ell} $\to$ {\bf NC-Tor} maps 
${\cal E}_{CM}$ to a two-dimensional noncommutive torus with real multiplication  ${\cal A}_{RM}^2$. 
To  calculate $L({\cal A}_{RM}^2,s)$,   let  $A\in GL_2({\Bbb Z})$  be  positive  matrix
corresponding the {\it shift automorphism} of ${\cal A}_{RM}^2$,  i.e. 
\displaymath
A\left(\matrix{1\cr\theta}\right)=
\lambda_A \left(\matrix{1\cr\theta}\right),
\enddisplaymath
where $\theta$ is a quadratic irrationality  and $\lambda_A$ the Perron-Frobenius eigenvalue 
of $A$.  If $p$ is a prime,  then the characteristic  polynomial of matrix $A^{\pi(p)}$ 
can be written as   {\bf Char} $A^{\pi(p)}=x^2-tr~(A^{\pi(p)})x-1$;
therefore,  the local zeta function of  torus ${\cal A}_{RM}^2$ has the form
\displaymath
\zeta_p({\cal A}_{RM}^2,z)={1\over 1-tr~(A^{\pi(p)})z+pz^2}. 
\enddisplaymath
On the other hand,   Lemma \ref{lem6.3.6}  says that
\displaymath
tr~(A^{\pi(p)})= tr~(\psi_{{\cal E}(K)}({\goth P})),
\enddisplaymath
where  $\psi_{{\cal E}(K)}$ is the Gr\"ossencharacter on $K$ and  ${\goth P}$ the
prime ideal of $K$ over $p$.   But we know that the local zeta function of ${\cal E}_{CM}$
has the form
\displaymath
\zeta_p({\cal E}_{CM}, z)={1\over 1-tr~(\psi_{{\cal E}(K)}({\goth P}))z+pz^2};
\enddisplaymath
thus for each prime $p$ it holds  $\zeta_p({\cal A}_{RM}^2,z)=\zeta_p({\cal E}_{CM},z)$.
Leaving aside the  bad primes,  one derives  the following important equality of the $L$-functions
\displaymath
L({\cal A}_{RM}^2,s)\quad\equiv\quad L({\cal E}_{CM}, s),
\enddisplaymath
where $L({\cal E}_{CM}, s)$ is the Hasse-Weil $L$-function of elliptic curve ${\cal E}_{CM}$. 
Case $n=1$ of Theorem \ref{thm6.4.1} becomes  an implication of the following lemma.
\begin{lem}\label{lem6.4.1}
$L({\cal E}_{CM},s)\equiv L(\sigma_2, s)$,
where  $L(\sigma_2, s)$ the Artin $L$-function for an irreducible  representation 
$\sigma_2: Gal~(K | k)\to GL_2({\Bbb C})$.
\end{lem}
{\it Proof.}  The Deuring theorem says that 
\displaymath
L({\cal E}_{CM}, s)=L(\psi_K, s)L(\overline{\psi}_K, s),
\enddisplaymath
where  $L(\psi_K,s)$ is the Hecke $L$-series attached to the 
Gr\"ossencharacter  $\psi: {\Bbb A}_K^*\to {\Bbb C}^*$;   here 
${\Bbb A}_K^*$ denotes the adele ring of the field $K$ and the bar
means a complex conjugation,  see e.g.  [Silverman 1994]  \cite{S2},  p.175. 
 Because  our elliptic curve has complex   multiplication,  the group $Gal~(K | k)$ is abelian;  
 one can apply  the result of   [Knapp  1997]  \cite{Kna1},  
Theorem 5.1,   which says that the Hecke $L$-series 
$L(\sigma_1\circ \theta_{K|k}, s)$ equals the Artin $L$-function $L(\sigma_1, s)$,
where $\psi_K=\sigma\circ \theta_{K|k}$ is the Gr\"ossencharacter and 
 $\theta_{K|k}: {\Bbb A}_K^*\to Gal~(K | k)$  the canonical homomorphism.
Thus one gets 
\displaymath
L({\cal E}_{CM}, s)\equiv L(\sigma_1,s)L(\overline{\sigma}_1,s),
\enddisplaymath
where $\overline{\sigma}_1: Gal~(K | k)\to {\Bbb C}$ means a (complex) conjugate
representation of the Galois group. 
Consider the local factors of the Artin $L$-functions $L(\sigma_1,s)$ and $L(\overline{\sigma}_1,s)$;
it is immediate, that they are $(1-\sigma_1(Fr_p)p^{-s})^{-1}$ and  $(1-\overline{\sigma}_1(Fr_p)p^{-s})^{-1}$,
respectively.  Let us consider a representation $\sigma_2: Gal~(K | k)\to GL_2({\Bbb C})$, 
such that 
\displaymath
\sigma_2(Fr_p)=
\left(\matrix{\sigma_1(Fr_p) & 0\cr 0 & \overline{\sigma}_1(Fr_p)}\right). 
\enddisplaymath
It can be verified, that $det^{-1}(I_2-\sigma_2(Fr_p)p^{-s})=
 (1-\sigma_1(Fr_p)p^{-s})^{-1}(1-\overline{\sigma}_1(Fr_p)p^{-s})^{-1}$,
i.e. $L(\sigma_2,s)=L(\sigma_1,s)L(\overline{\sigma}_1,s)$. 
Lemma \ref{lem6.4.1} follows.
$\square$

\bigskip
From  lemma \ref{lem6.4.1} and $L({\cal A}_{RM}^2,s)\equiv  L({\cal E}_{CM}, s)$,  one gets
\displaymath 
L({\cal A}_{RM}^2,s)\equiv L(\sigma_2,s)
\enddisplaymath
 for an irreducible   representation $\sigma_2: Gal~(K|k)\to GL_2({\Bbb C})$.
It remains to notice that  $L(\sigma_2,s)= L(\sigma_2',s)$,
where $\sigma_2': Gal~(K|{\Bbb Q})\to GL_2({\Bbb C})$, see e.g. 
[Artin 1924]  \cite{Art1},  Section 3.  Case $n=1$ of  Theorem \ref{thm6.4.1} is proved
$\square$

\subsubsection{Case $n=0$}
When $n=0$, one gets a one-dimensional (degenerate) noncommutative
torus; such an object, ${\cal A}_{\Bbb Q}$, can be obtained from the 
$2$-dimensional torus ${\cal A}^2_{\theta}$ by forcing $\theta=p/q\in {\Bbb Q}$ be a rational
number (hence our notation). 
One can always assume $\theta=0$ and, thus
\displaymath   
K_0^+({\cal A}_{\Bbb Q})\cong {\Bbb Z}. 
\enddisplaymath
The group of automorphisms of  ${\Bbb Z}$-module  $K_0^+({\cal A}_{\Bbb Q})\cong {\Bbb Z}$
is trivial,  i.e. the  multiplication by $\pm 1$;   hence matrix $A$ corresponding to the shift automorphisms 
is either $1$ or $-1$.  Since $A$ must be positive,  one  gets $A=1$.
 However,  $A=1$ is not a primitive;
indeed,  for any $N>1$  matrix $A'=\zeta_N$ gives us $A=(A')^N$,  
where $\zeta_N=e^{2\pi i\over N}$ is the $N$-th root of unity.
 Therefore, one gets 
 \displaymath
 A=\zeta_N.
 \enddisplaymath
 Since for the field ${\Bbb Q}$ it holds $\pi(n)=n$,  one obtains 
 $tr~(A^{\pi(p)})=tr~(A^p)=\zeta_N^p$.
A degenerate noncommutative torus, corresponding to the matrix $A=\zeta_N$,
we shall write as ${\cal A}_{\Bbb Q}^N$.

Suppose that $Gal~(K|{\Bbb Q})$ is abelian and   
let $\sigma: Gal~(K|{\Bbb Q})\to {\Bbb C}^{\times}$ be a
homomorphism.   By the {\it Artin reciprocity},   there
exists an integer  $N_{\sigma}$ and the Dirichlet character
\displaymath
 \chi_{\sigma}: ({\Bbb Z}/N_{\sigma} {\Bbb Z})^{\times}\to {\Bbb C}^{\times},
\enddisplaymath
such that $\sigma(Fr_p)=\chi_{\sigma}(p)$,   see e.g.  [Gelbart 1984]  \cite{Gel1}.
On the other hand,  it is verified directly, that $\zeta_{N_{\sigma}}^p=e^{{2\pi i\over N_{\sigma}}p}=
\chi_{\sigma}(p)$.   Therefore  {\bf Char} $(A^p)= \chi_{\sigma}(p)x-1$
and one gets 
\displaymath
\zeta_p({\cal A}_{\Bbb Q}^{N_{\sigma}},z)= {1\over 1-\chi_{\sigma}(p)z},
\enddisplaymath
where $\chi_{\sigma}(p)$ is the Dirichlet character. 
Therefore, $L({\cal A}^{N_{\sigma}}_{\Bbb Q}, s)\equiv L(s,\chi_{\sigma})$ 
is the Dirichlet $L$-series;    such a series, by construction,  coincides 
with the Artin $L$-series of the representation $\sigma: Gal~(K|{\Bbb Q})\to {\Bbb C}^{\times}$.
Case $n=0$ of Theorem \ref{thm6.4.1} is proved.
$\square$

 \index{Dirichlet $L$-series}

 \index{Artin $L$-function}
 \index{class field theory}

\subsection{Supplement:   Artin $L$-function}
The {\it Class Field Theory} (CFT) studies algebraic extensions of the number fields;
the objective of CFT is a description of arithmetic of the extension $E$ in terms of arithmetic
of the ground field $k$ and the Galois group $Gal~(E|k)$ of the extension.  
Unless $Gal~(E|k)$ is abelian,  the CFT is out of reach so far;   yet a series of conjectures
called the {\it Langlands program} (LP) are designed to attain the goals of CFT.
We refer the interested reader to [Gelbart 1984]   \cite{Gel1}  for an introduction 
to the CFT and LP;   roughly speaking,  the LP consists in an  $n$-dimensional  
generalization of the Artin reciprocity  based on the ideas and methods of representation
theory of the locally compact Lie groups.  The centerpiece of LP is the {\it Artin $L$-function}
attached to representation $\sigma: Gal~(E|k)\to GL_n({\Bbb C})$ of the Galois group 
of $E$;  we shall give a brief account of this $L$-function following the survey by 
 [Gelbart 1984]   \cite{Gel1}.

 \index{isotropy subgroup}
 \index{Frobenius element}

The fundamental problem in algebraic number theory is to describe how an ordinary
prime $p$ factors into prime ideals ${\goth P}$ in the ring of integers of an arbitrary finite extensions
$E$ of the rational field ${\Bbb Q}$.    Let $O_E$ be the ring of integers of the extension $E$
and $pO_E$ a principal ideal;  it is known that
\displaymath
pO_E=\prod {\goth P}_i,
\enddisplaymath
where ${\goth P}_i$ are prime ideals of $O_E$.  If $E$ is the Galois extension of ${\Bbb Q}$
and $Gal~(E|{\Bbb Q})$ is the corresponding Galois group,  then each automorphism 
$g\in Gal~(E|{\Bbb Q})$ ``moves around''  the ideals ${\goth P}_i$ in the prime decomposition
of $p$ over $E$.  An {\it isotropy subgroup} of $Gal~(E|{\Bbb Q})$ (for given $p$)  consists of 
the  elements of    $Gal~(E|{\Bbb Q})$ which fix all the ideals ${\goth P}_i$.   For simplicity,
we shall assume that $p$ is {\it unramified} in $E$, i.e. all ${\goth P}_i$ are distinct;
in this case the isotropy  subgroup are {\it cyclic}.   The (conjugacy class of)  generator
in the cyclic isotropy subgroup of $Gal~(E|{\Bbb Q})$ corresponding to $p$ is called 
the {\it Frobenius element} and denoted by $Fr_p$.   The element $Fr_p\in Gal~(E|{\Bbb Q})$ 
describes completely the factorization of $p$ over $E$ and the major goal of the CFT 
is to express $Fr_p$ in terms of  arithmetic of the ground field ${\Bbb Q}$.  
To handle this hard problem,  it was suggested by E.~Artin to consider the $n$-dimensional
irreducible representations 
\displaymath
\sigma_n:  Gal~(E|{\Bbb Q})\longrightarrow GL_n({\Bbb C}),
\enddisplaymath
of the Galois group $Gal~(E|{\Bbb Q})$,  see [Artin 1924]  \cite{Art1}.   The idea
was to use the characteristic polynomial {\bf Char}  
$(\sigma_n(Fr_p)):=\det (I_n-\sigma_n(Fr_p)z)$ of the 
matrix $\sigma_n(Fr_p)$;  the polynomial is independent of the similarity class 
of   $\sigma_n(Fr_p)$ in the group $GL_n({\Bbb C})$ and provides
an {\it intrinsic} description of the Frobenius element $Fr_p$.  
\begin{dfn}
By an Artin zeta function of representation $\sigma_n$ 
one understands the function
\displaymath
\zeta_p(\sigma_n, z):={1\over \det(I_n-\sigma_n(Fr_p)z)}, \quad z\in {\Bbb C}.  
\enddisplaymath
By an Artin $L$-function of representation $\sigma_n$  one understands the product 
\displaymath
L(\sigma_n, s):= \prod_p ~\zeta_p(\sigma_n, p^{-s}), \quad s\in {\Bbb C}, 
\enddisplaymath
over all but a finite set of primes $p$. 
\end{dfn}
 \index{Artin reciprocity}
\begin{rmk}
{\bf (Artin reciprocity)}
\textnormal{
If $n=1$ and $Gal~(E|{\Bbb Q})\cong {\Bbb Z}/N{\Bbb Z}$ is abelian,
then the calculation of the Artin $L$-function gives the  equality
\displaymath
L(\chi, s)= \prod_p {1\over 1-\chi(p)p^{-s}}, 
\enddisplaymath
where $\chi:  ({\Bbb Z}/N{\Bbb Z})\to {\Bbb C}^{\times}$ is the Dirichlet character;
the RHS of the equality is known as the {\it Dirichlet $L$-series} for $\chi$. Thus one gets
a  formula
\displaymath
\sigma(Fr_p)=\chi(p)
\enddisplaymath
called the {\it Artin reciprocity law};    the formula generalizes many classical reciprocity
results known for the particular values of $N$.  
}
\end{rmk}

\vskip1.5cm\noindent
{\bf Guide to the literature.} 
The Artin $L$-function first appeared  in [Artin 1924]  \cite{Art1}.
The origins of the Langlands program (and philosophy) can be found 
in his  letter to Andr\'e Weil,   see  [Langlands  1960's]   \cite{Lan1}.    An excellent introduction 
to the Langlands program has been written by [Gelbart 1984] \cite{Gel1}.
The Langlands program for the even-dimensional noncommutative tori
was the subject of \cite{Nik14}.

 \index{finite field}
 \index{Serre $C^*$-algebra}
 \index{Weil Conjectures}

\section{Projective varieties over finite fields}
In Section 5.3  we constructed a covariant functor 
\displaymath
F:  \hbox{{\bf Proj-Alg}}\longrightarrow \hbox{{\bf C*-Serre}}
\enddisplaymath
from the category of complex projective varieties $V({\Bbb C})$ to a category
of the {\it Serre $C^*$-algebras}  ${\cal A}_V$.   Provided $V({\Bbb C})\cong V(K)$
for a number field $K\subset {\Bbb C}$,  one can reduce variety $V(K)$ modulo
a prime ideal ${\goth P}\subset K$ over the prime product $q=p^r$;  the reduction
corresponds to  a projective variety $V({\Bbb F}_q)$ defined  over the {\it finite
 field} ${\Bbb F}_q$.   In this section we express the geometric invariant 
  $|V({\Bbb F}_q)|$ of  $V({\Bbb F}_q)$ (the number of points of  variety  $V({\Bbb F}_q)$) 
   in terms of the  noncommutative invariants of the Serre $C^*$-algebra ${\cal A}_V$;
  the obtained formula shows  interesting links   to the {\it Weil Conjectures},   
see e.g.  [Hartshorne 1977]  \cite{H1},  Appendix C for an introduction. 
We test our formula on  the concrete families  of complex multiplication and   
rational elliptic curves.

 \index{trace of Frobenius endomorphism}
 \index{$\ell$-adic cohomology}

\subsection{Traces of Frobenius endomorphisms}
The number of solutions of a system of polynomial equations over
a finite field is an important invariant of the system and 
an old problem dating back to Gauss.    Recall that if ${\Bbb F}_q$
is a field with $q=p^r$ elements and $V({\Bbb F}_q)$ a smooth $n$-dimensional 
projective variety over ${\Bbb F}_q$,  then  one can define a zeta
function $Z(V; t):= \exp~\left(\sum_{r=1}^{\infty} |V({\Bbb F}_{q^r})|{t^r\over r}\right)$;
the function is rational,  i.e.
\displaymath
Z(V; t)={P_1(t)P_3(t)\dots P_{2n-1}(t)\over P_0(t)P_2(t)\dots P_{2n}(t)},
\enddisplaymath
where   $P_0(t)=1-t$, $P_{2n}(t)=1-q^nt$ and for each $1\le i\le 2n-1$ the polynomial 
$P_i(t)\in {\Bbb Z}[t]$ can be written as  $P_i(t)=\prod_{j=1}^{deg~P_i(t)} (1-\alpha_{ij}t)$
so that   $\alpha_{ij}$ are algebraic integers with
$|\alpha_{ij}|=q^{i\over 2}$, see e.g.  [Hartshorne 1977]  \cite{H1},  pp. 454-457. 
The  $P_i(t)$ can be viewed as characteristic polynomial of the Frobenius 
endomorphism $Fr_q^i$  of  the $i$-th $\ell$-adic cohomology
group  $H^i(V)$;  such an endomorphism is induced 
by the map acting on  points  of variety $V({\Bbb F}_q)$
according to the formula  $(a_1,\dots, a_n)\mapsto (a_1^q,\dots, a_n^q)$;
we assume throughout  the  {\it Standard Conjectures},  see [Grothendieck 1968]  \cite{Gro1}. 
If   $V({\Bbb F}_q)$ is  defined  by a  system of polynomial equations,
 then the number of solutions of the system  is given by the formula 
\displaymath
|V({\Bbb F}_q)|=\sum_{i=0}^{2n}(-1)^i  ~tr~(Fr^i_q), 
\enddisplaymath
where  $tr$ is the trace of Frobenius  endomorphism,  see [Hartshorne 1977]  \cite{H1},   {\it loc. cit}.

 Let $V(K)$ be a complex projective variety defined over an algebraic number field $K\subset {\Bbb C}$;
 suppose that  projective variety  $V({\Bbb F}_q)$ is the reduction of $V(K)$ modulo  the prime ideal 
 ${\goth P}\subset K$ corresponding to  $q=p^r$.   Denote by ${\cal A}_V$ the {\it Serre $C^*$-algebra} 
of projective variety $V(K)$,  see Section 5.3.1.  
Consider the  stable $C^*$-algebra of ${\cal A}_V$,  i.e. the $C^*$-algebra  ${\cal A}_V\otimes {\cal K}$,
 where ${\cal K}$ is the $C^*$-algebra of compact operators on ${\cal H}$. 
 Let $\tau: {\cal A}_V\otimes {\cal K}\to {\Bbb R}$   be the unique normalized trace (tracial state) on  ${\cal A}_V\otimes {\cal K}$, 
  i.e. a positive linear functional   of norm $1$  such that $\tau(yx)=\tau(xy)$ for all $x,y\in {\cal A}_V\otimes {\cal K}$,  see  
 [Blackadar 1986] \cite{B},  p. 31.    
Recall that ${\cal A}_V$ is the crossed product $C^*$-algebra of the form
${\cal A}_V\cong C(V)\rtimes {\Bbb Z}$,  where $C(V)$ is the 
commutative $C^*$-algebra of complex valued functions on $V$ 
and the product is taken by an automorphism of algebra $C(V)$ 
induced by the map $\sigma: V\to V$,  see Lemma \ref{lem5.3.2}.   
From the Pimsner-Voiculescu six term exact sequence for
crossed products,  one gets the  short exact sequence of algebraic $K$-groups 
\displaymath
0\to K_0(C(V))\buildrel  i_*\over\to  K_0({\cal A}_V)\to K_1(C(V))\to 0, 
\enddisplaymath
where   map  $i_*$  is induced by an  embedding of $C(V)$ 
into ${\cal A}_V$,   see   [Blackadar 1986]  \cite{B}, p. 83 for the details.  
 We  have $K_0(C(V))\cong K^0(V)$ and 
$K_1(C(V))\cong K^{-1}(V)$,  where $K^0$ and $K^{-1}$  are  the topological
$K$-groups of variety $V$, see  [Blackadar 1986]  \cite{B}, p. 80. 
By  the Chern character formula,  one gets
 \index{Chern character formula}
\displaymath
\left\{
\begin{array}{ccc}
K^0(V)\otimes {\Bbb Q} &\cong&  H^{even}(V; {\Bbb Q})\\
K^{-1}(V)\otimes {\Bbb Q} &\cong&   H^{odd}(V; {\Bbb Q}),
\end{array}
\right.
\enddisplaymath
where $H^{even}$  ($H^{odd}$)  is the direct sum of even (odd, resp.) 
cohomology groups of $V$. 
\begin{rmk}
\textnormal{
 It is known,  that  $K_0({\cal A}_V\otimes {\cal K})\cong K_0({\cal A}_V)$  because
of  stability of the $K_0$-group with respect to tensor products by the algebra 
${\cal K}$,  see e.g.   [Blackadar 1986]  \cite{B}, p. 32.
}
\end{rmk}
Thus one gets the  commutative diagram shown in Fig. 6.6, 
where $\tau_*$ denotes  a homomorphism  induced on $K_0$ by  the canonical  trace 
$\tau$ on the $C^*$-algebra  ${\cal A}_V\otimes {\cal K}$.

\bigskip
\begin{figure}[here]
\begin{picture}(300,100)(-50,5)
\put(160,72){\vector(0,-1){35}}
\put(80,65){\vector(2,-1){45}}
\put(240,65){\vector(-2,-1){45}}
\put(10,80){$ H^{even}(V)\otimes {\Bbb Q} 
\buildrel  i_*\over\longrightarrow  K_0({\cal A}_V\otimes{\cal K})\otimes {\Bbb Q} 
\longrightarrow H^{odd}(V)\otimes {\Bbb Q}$}
\put(167,55){$\tau_*$}
\put(157,20){${\Bbb R}$}
\end{picture}
\caption{$K$-theory of the Serre $C^*$-algebra ${\cal A}_V$.}
\end{figure}

\noindent
Because   $H^{even}(V):=\oplus_{i=0}^n H^{2i}(V)$ and  
$H^{odd}(V):=\oplus_{i=1}^n H^{2i-1}(V)$,   one gets  for each  $0\le i\le 2n$ 
 an injective  homomorphism 
 \displaymath
 H^i(V)\to  {\Bbb R}
 \enddisplaymath
 and we shall denote by $\Lambda_i$  an  additive abelian subgroup of real numbers
 defined by the homomorphism.    
\begin{rmk}
\textnormal{
The $\Lambda_i$ is  called  a {\it pseudo-lattice}  [Manin 2004]  \cite{Man1},  Section 1. 
}
\end{rmk}
Recall that  endomorphisms  of a pseudo-lattice are given as 
multiplication of points of $\Lambda_i$ by the real numbers $\alpha$
such that $\alpha\Lambda_i\subseteq\Lambda_i$.   It is known that
$End~(\Lambda_i)\cong {\Bbb Z}$ or $End~(\Lambda_i)\otimes {\Bbb Q}$
is a real algebraic number field such that $\Lambda_i\subset  End~(\Lambda_i)\otimes {\Bbb Q}$,
see e.g.  [Manin 2004]  \cite{Man1},  Lemma 1.1.1 for the case of quadratic fields. 
We shall write $\varepsilon_i$ to denote the unit of the order in the field $K_i:=End~(\Lambda_i)\otimes {\Bbb Q}$,  
which induces the  shift automorphism  of $\Lambda_i$,  
 see [Effros 1981]  \cite{E},  p. 38   for the details and terminology. 
  Let $p$ be a ``good prime''   in  the reduction $V({\Bbb F}_q)$ of 
  complex projective variety $V(K)$ modulo a prime ideal over  $q=p^r$.  
   Consider a sub-lattice $\Lambda_i^{q}$ of $\Lambda_i$ of the index $q$; 
  by an  index of  the sub-lattice we understand  its  index as an abelian subgroup of $\Lambda_i$.
We shall write  $\pi_i(q)$ to  denote  an
 integer,   such that  multiplication by $\varepsilon_i^{\pi_i(q)}$ 
 induces  the shift automorphism of   $\Lambda_i^q$. 
  The trace of an algebraic number will be written as $tr~(\bullet)$.    
  The following result relates invariants $\varepsilon_i$ and $\pi_i(q)$ of the $C^*$-algebra
 ${\cal A}_V$ to the cardinality of the set $V({\Bbb F}_q)$.   
\begin{thm}\label{thm6.5.1}
{\bf (Noncommutative invariant of projective varieties over finite fields)}
$$|V({\Bbb F}_q)|=\sum_{i=0}^{2n}(-1)^i  ~tr~\left(\varepsilon_i^{\pi_i(q)}\right)$$.  
\end{thm}

\subsection{Proof of  Theorem \ref{thm6.5.1} }
\begin{lem}\label{lem5.5.1}
There exists a symplectic unitary matrix $\Theta_q^i\in Sp ~(deg~P_i; ~{\Bbb R})$,  such that
\displaymath
Fr^i_q=q^{i\over 2}\Theta_q^i.  
\enddisplaymath
\end{lem}
{\it Proof.}  Recall that the eigenvalues of $Fr_q^i$ have absolute value $q^{i\over 2}$;
they come in the complex conjugate pairs.  On the other hand,  symplectic unitary matrices in
group   $Sp ~(deg~P_i; ~{\Bbb R})$ are known to have eigenvalues of absolute value $1$  
coming in complex conjugate pairs.   Since the spectrum of a matrix defines the  similarity class of 
matrix,  one can write the characteristic polynomial of $Fr_q^i$ in the form
\displaymath
 P_i(t)=det~(I-q^{i\over 2}\Theta_q^i t), 
\enddisplaymath
where matrix $\Theta_q^i\in Sp~(deg~P_i; ~{\Bbb Z})$ and its  eigenvalues have  absolute value $1$. 
It remains to compare the above equation with the formula
\displaymath
 P_i(t)=det~(I-Fr_q^i t), 
\enddisplaymath
i.e. $Fr_q^i=q^{i\over 2}\Theta_q^i$.   Lemma \ref{lem5.5.1} follows.
$\square$

\begin{lem}\label{lem5.5.2}
Using a symplectic transformation one can bring matrix $\Theta_q^i$ to the block form
\displaymath
\Theta_q^i=\left(\matrix{A & I\cr -I & 0}\right), 
\enddisplaymath
where $A$ is a positive symmetric and $I$  the identity matrix.  
\end{lem}
{\it Proof.} Let  us write $\Theta_q^i$ in the block form
\displaymath
\Theta_q^i=\left(\matrix{A & B\cr C & D}\right), 
\enddisplaymath
where matrices $A,B,C, D$ are invertible and their transpose $A^T,B^T,C^T, D^T$
satisfy the symplectic equations 
\displaymath
\left\{
\begin{array}{cc}
A^TD-C^TB &= I,\\
A^TC-C^TA &= 0,\\
B^TD-D^TB &= 0. 
\end{array}
\right.
\enddisplaymath
Recall that symplectic matrices correspond to the linear
fractional  transformations  $\tau\mapsto {A\tau+B\over C\tau +D}$
of the Siegel half-space ${\Bbb H}_n=\{\tau=(\tau_j)\in {\Bbb C}^{{n(n+1)\over 2}}
~|~\Im (\tau_j)>0\}$ consisting of symmetric $n\times n$ matrices, see e.g.  [Mumford 1983]
\cite{MU},   p. 173.   One can always multiply    the nominator and denominator
of such a transformation  by $B^{-1}$  without affecting the transformation;
thus with no loss of generality,  we can assume that $B=I$.  
We shall consider the symplectic matrix $T$ and its inverse $T^{-1}$
given by the formulas
\displaymath
T=\left(\matrix{I & 0\cr D & I}\right) \quad\hbox{and}
\quad T^{-1}=\left(\matrix{I & 0\cr -D & I}\right). 
\enddisplaymath
It is verified directly, that  
\displaymath
T^{-1}\Theta_q^i T=
\left(\matrix{I & 0\cr -D & I}\right)
\left(\matrix{A & I\cr C & D}\right)
\left(\matrix{I & 0\cr D & I}\right)
=\left(\matrix{A+D & I\cr C-DA & 0}\right). 
\enddisplaymath
The system of symplectic equations with $B=I$ implies the 
following two equations
\displaymath
A^TD-C^T = I
 \quad\hbox{and}\quad
D=D^T. 
\enddisplaymath
Applying transposition to the both parts  of the first equation of  the above equations,  
 one gets  $(A^TD-C^T)^T = I^T$ and,  therefore, $D^TA-C=I$.    
But the second equation  says that $D^T=D$;  thus one arrives 
at the equation $DA-C=I$.  The latter gives us $C-DA=-I$, 
which we substitute  in the above equations  and get (in a new notation) the
conclusion of Lemma \ref{lem5.5.2}. 
Finally, the middle of the symplectic equations  with $C=-I$ implies $A=A^T$,
i.e. $A$ is a symmetric matrix.   Since the eigenvalues of symmetric matrix
are always real and in view of  $tr~(A)>0$ (because $tr~(Fr_q^i)>0$),  one concludes that 
$A$ is similar to a positive matrix, see e.g.  [Handelman 1981]  \cite{Han2}, Theorem 1.   
Lemma \ref{lem5.5.2} follows. 
$\square$

\begin{lem}\label{lem5.5.3}
The symplectic unitary transformation $\Theta_q^i$ of $H^i(V; {\Bbb Z})$
descends to  an automorphism of $\Lambda_i$  given by the matrix
\displaymath
M_q^i=\left(\matrix{A & I\cr I & 0}\right).  
\enddisplaymath
\end{lem}
\begin{rmk}
\textnormal{
In other words,  Lemma \ref{lem5.5.3} says that functor $F:$  {\bf Proj-Alg} $\to$ {\bf C*-Serre}
acts between matrices  $\Theta_q^i$ and $M_q^i$  according to the formula
\displaymath
F: \left(\matrix{A & I\cr -I & 0}\right)\longmapsto\left(\matrix{A & I\cr I & 0}\right).
\enddisplaymath
}
\end{rmk}
{\it Proof.}
Since $\Lambda_i\subset K_i$ there exists a basis of $\Lambda_i$
consisting of algebraic numbers;  denote by $(\mu_1,\dots,\mu_k;   ~\nu_1,\dots,\nu_k)$ a basis
of  $\Lambda_i$ consisting of positive algebraic numbers  $\mu_i>0$ and 
$\nu_i>0$.  Using the injective homomorphism $\tau_*$,  
one can descend  $\Theta^i_q$   to an  automorphism of $\Lambda_i$   so that 
\displaymath
\left(\matrix{\mu'\cr \nu'}\right)
=\left(\matrix{A & I\cr -I & 0}\right)
\left(\matrix{\mu\cr \nu}\right)=
\left(\matrix{A\mu+\nu\cr -\mu}\right),
\enddisplaymath
where $\mu=(\mu_1,\dots,\mu_k)$ and $\nu=(\nu_1,\dots,\nu_k)$.
Because vectors $\mu$ and $\nu$ consist of positive entries and 
$A$ is a positive matrix,  it is immediate that  $\mu'=A\mu+\nu>0$ while
$\nu'=-\mu<0$.   
 \index{Markov category}
\begin{rmk}
\textnormal{
All automorphisms in the (Markov) category of pseudo-lattices 
come from multiplication of the basis vector 
$(\mu_1,\dots,\mu_k;   ~\nu_1,\dots,\nu_k)$ of $\Lambda_i$ 
by an algebraic unit $\lambda>0$ of field $K_i$;   in particular, 
any such an automorphism must be given by a   
non-negative matrix,  whose  Perron-Frobenius eigenvalue coincides with $\lambda$.
Thus for any automorphism of $\Lambda_i$ it must hold 
$\mu'>0$ and $\nu'>0$.  
}
\end{rmk}
In view of the above,  we shall  consider an automorphism of $\Lambda_i$
given by  matrix $M_q^i=(A, I, I, 0)$;   clearly,  for $M_q^i$ it holds
 $\mu'=A\mu+\nu>0$ and  $\nu'=\mu>0$.   
Therefore  $M_q^i$ is a non-negative matrix satisfying  
  the   necessary condition  to belong to the  Markov category.  
   It is also a sufficient one,  because the similarity class of $M_q^i$ 
contains a representative whose Perron-Frobenius eigenvector can be taken for a
basis $(\mu, \nu)$ of $\Lambda_i$.  This argument finishes the proof of Lemma \ref{lem5.5.3}.
$\square$

\begin{cor}\label{cor5.5.1}
$tr~(M_q^i)=tr~(\Theta_q^i)$. 
\end{cor}
{\it Proof.}  This fact is an implication of the above formulas 
and  a direct computation  $tr~(M_q^i)=tr~(A)=tr~(\Theta_q^i)$.
$\square$

\begin{dfn}
We shall call $q^{i\over 2}M_q^i$ a Markov endomorphism of $\Lambda_i$
and denote it by $Mk_q^i$.  
\end{dfn}
\begin{lem}\label{lem5.5.4}
$tr~(Mk_q^i)=tr~(Fr_q^i)$.
\end{lem}
{\it Proof.}
Corollary \ref{cor5.5.1}  says that $tr~(M_q^i)=tr~(\Theta_q^i)$,  and  
therefore
\displaymath
\begin{array}{ccc}
tr~(Mk_q^i) &=&  tr~(q^{i\over 2}M_q^i)= q^{i\over 2}~tr~(M_q^i)=\\
                     &=& q^{i\over 2} ~tr~(\Theta_q^i) =  tr~(q^{i\over 2}\Theta_q^i)=tr~(Fr_q^i). 
\end{array}
\enddisplaymath
In words,    Frobenius and Markov endomorphisms have the same trace,
i.e.  $tr~(Mk_q^i)=tr~(Fr_q^i)$.
Lemma \ref{lem5.5.4} follows. 
$\square$

\begin{rmk}
\textnormal{
Notice  that,  unless $i$ or $r$  are even,  neither $\Theta_q^i$ nor $M_q^i$ are 
integer matrices;   yet    $Fr_q^i$ and $Mk_q^i$ are always  integer matrices.    
 }
\end{rmk}
\begin{lem}\label{lem5.5.5}
There exists an algebraic unit $\omega_i\in K_i$ such that:

\medskip
(i) $\omega_i$ corresponds to the shift automorphism of an index
$q$ sub-lattice of pseudo-lattice $\Lambda_i$; 

\smallskip
(ii)  $tr~(\omega_i)=tr~(Mk_q^i)$.
\end{lem}
{\it Proof.}
(i) To prove Lemma \ref{lem5.5.5},   we shall  use the notion of a stationary dimension group and 
the corresponding shift  automorphism;  we refer the reader to  [Effros 1981]  \cite{E}, p. 37 and 
[Handelman  1981]  \cite{Han2},  p.57 for the notation and details  on  stationary dimension groups
and a survey of [Wagoner 1999]   \cite{Wag1}  for the general theory of subshifts of finite type.   
Consider a stationary dimension group, $G(Mk_q^i)$, generated by  the Markov
endomorphism $Mk_q^i$
\displaymath
{\Bbb Z}^{b_i} \buildrel Mk_q^i \over\to 
{\Bbb Z}^{b_i} \buildrel Mk_q^i \over\to 
{\Bbb Z}^{b_i} \buildrel Mk_q^i \over\to
\dots,
\enddisplaymath
where $b_i=deg~P_i(t)$.  
Let $\lambda_M$ be the Perron-Frobenius eigenvalue of matrix $M_q^i$. 
It is known, that $G(Mk_q^i)$ is  order-isomorphic to a  dense additive  abelian subgroup 
${\Bbb Z}[{1\over\lambda_M}]$ of  ${\Bbb R}$;  
here ${\Bbb Z}[x]$ is the set of all polynomials in one variable with the integer
coefficients.
Let $\widehat{Mk_q^i}$  be a  shift automorphism of  $G(Mk_q^i)$ [Effros 1981]  \cite{E},
p. 37.    To calculate the automorphism,  notice that multiplication 
of ${\Bbb Z}[{1\over\lambda_M}]$ by $\lambda_M$  induces an automorphism
of dimension group ${\Bbb Z}[{1\over\lambda_M}]$.   Since the determinant of matrix $M_q^i$ (i.e. the degree of 
Markov endomorphism)   is equal to $q^n$,   one concludes that  such an  automorphism corresponds to  a unit 
of the  endomorphism ring  of  a  sub-lattice of  $\Lambda_i$ of index $q^n$.
We shall denote such a unit by $\omega_i$.    Clearly,  $\omega_i$ generates 
the required shift automorphism $\widehat{Mk_q^i}$  through  multiplication 
of  dimension group  ${\Bbb Z}[{1\over\lambda_M}]$  by the algebraic number  $\omega_i$. 
Item (i) of Lemma \ref{lem5.5.5} follows.

 \index{Artin-Mazur zeta function}

\bigskip
(ii)  Consider the Artin-Mazur zeta function of $Mk_q^i$
\displaymath
\zeta_{Mk_q^i}(t)=\exp\left(\sum_{k=1}^{\infty}{tr~\left[(Mk_q^i)^k\right]\over k}t^k\right)
\enddisplaymath
and such of $\widehat{Mk_q^i}$
\displaymath
\zeta_{\widehat{Mk_q^i}}(t)=\exp\left(\sum_{k=1}^{\infty}{tr~\left[(\widehat{Mk_q^i})^k\right]\over k}t^k\right).
\enddisplaymath
Since  $Mk_q^i$ and $\widehat{Mk_q^i}$ are shift equivalent matrices,  one concludes 
 that  $\zeta_{Mk_q^i}(t)\equiv \zeta_{\widehat{Mk_q^i}}(t)$,    see [Wagoner 1999]  \cite{Wag1}, p. 273.    
 In particular, 
\displaymath 
tr~(Mk_q^i)=tr~(\widehat{Mk_q^i}).
\enddisplaymath
But   $tr~(\widehat{Mk_q^i})=tr~(\omega_i)$,  where on the right hand side is the
trace of an algebraic number.  In view of the above,  one gets the conclusion of 
item (ii)  of  Lemma \ref{lem5.5.5}. 
$\square$

\begin{lem}\label{lem5.5.6}
There exists a positive integer $\pi_i(q)$,  such that
\displaymath
\omega_i=\varepsilon_i^{\pi_i(q)},
\enddisplaymath
where $\varepsilon_i\in End~(\Lambda_i)$ is the fundamental unit  
corresponding to  the shift automorphism of pseudo-lattice $\Lambda_i$. 
\end{lem}
{\it Proof.}  
Given an automorphism $\omega_i$ of a finite-index sub-lattice of $\Lambda_i$
one can extend $\omega_i$ to an automorphism of entire $\Lambda_i$,
since $\omega_i\Lambda_i=\Lambda_i$.  Therefore each unit of (endomorphism ring of)
a sub-lattice is also a unit of the host pseudo-lattice.  Notice that the converse statement 
is false in general. 
On the other hand,   by virtue of the Dirichlet Unit Theorem  each  unit of $End~(\Lambda_i)$ 
is a product of a finite number of (powers of)  fundamental units of   $End~(\Lambda_i)$.
We shall denote by $\pi_i(q)$ the least positive integer,  such that  $\varepsilon_i^{\pi_i(q)}$
is the shift automorphism of a sub-lattice of index $q$ of pseudo-lattice $\Lambda_i$. 
The number $\pi_i(q)$ exists and uniquely defined,  albeit no general formula for its calculation 
 is known,   see Remark \ref{rmk5.5.6}.   It is clear from construction,  that $\pi_i(q)$ 
 satisfies the claim of Lemma \ref{lem5.5.6}.
 $\square$

 \index{Dirichlet Unit Theorem}

\begin{rmk}\label{rmk5.5.6}
\textnormal{
No general formula for the number $\pi_i(q)$ as a function of $q$ 
is known;  however,  if the rank of $\Lambda_i$ is two (i.e. $n=1$), 
then there are classical results recorded in  e.g.   [Hasse 1950]  \cite{HA},  p.298;
see also Section 5.5.3.     
}
\end{rmk}

\bigskip
Theorem \ref{thm6.5.1} follows from  Lemmas \ref{lem5.5.4}-\ref{lem5.5.6}
and the known formula $|V({\Bbb F}_q)|=\sum_{i=0}^{2n}(-1)^i  ~tr~(Fr^i_q)$.
$\square$

\subsection{Examples}
Let $V({\Bbb C})\cong {\cal E}_{\tau}$ be an elliptic curve;   it is well known that 
its Serre $C^*$-algebra ${\cal A}_{\cal E_{\tau}}$  is isomorphic to the  noncommutative torus ${\cal A}_{\theta}$
with the   unit   scaled by  a constant $0<\log\mu<\infty$.  
Furthermore,  $K_0({\cal A}_{\theta})\cong K_1({\cal A}_{\theta})\cong {\Bbb Z}^2$
and  the canonical trace $\tau$ on ${\cal A}_{\theta}$  gives us the following 
formula
\displaymath
\tau_*(K_0({\cal A}_{\cal E_{\tau}}\otimes {\cal K}))=\mu({\Bbb Z}+{\Bbb Z}\theta).
\enddisplaymath
 Because  $H^0({\cal E}_{\tau}; {\Bbb Z})=H^2({\cal E}_{\tau};  {\Bbb Z})\cong {\Bbb Z}$ while  
 $H^1({\cal E}_{\tau};  {\Bbb Z})\cong {\Bbb Z}^2$,  one gets 
  the following pseudo-lattices
\displaymath
\Lambda_0=\Lambda_2\cong {\Bbb Z} \quad\hbox{and}\quad  
\Lambda_1\cong\mu({\Bbb Z}+{\Bbb Z}\theta).
\enddisplaymath
For the sake of simplicity,  we shall focus on the following families of elliptic curves.

\subsubsection{Complex multiplication}
Suppose  that ${\cal E}_{\tau}$ has  complex multiplication;  recall that such a
curve was denoted by ${\cal E}_{CM}^{(-D,f)}$,  i.e.  the endomorphism
ring of ${\cal E}_{\tau}$   is  an order of conductor $f\ge 1$ in the imaginary 
quadratic field ${\Bbb Q}(\sqrt{-D})$.  By the  results of Section 5.1 on   elliptic curve 
${\cal E}_{CM}^{(-D,f)}$,   the formulas for $\Lambda_i$  are as follows
\displaymath
\Lambda_0=\Lambda_2\cong {\Bbb Z} \quad\hbox{and}\quad
\Lambda_1 = \varepsilon[{\Bbb Z}+(f\omega){\Bbb Z}], 
\enddisplaymath
where  $\omega={1\over 2}(1+\sqrt{D})$ if $D\equiv 1~mod~4$ and $D\ne 1$  or $\omega=\sqrt{D}$ if $D\equiv 2,3~mod~4$
and  $\varepsilon>1$ is the fundamental unit of  order ${\Bbb Z}+(f\omega){\Bbb Z}$.    
\begin{rmk}
\textnormal{
The reader can verify,  that  $\Lambda_1\subset K_1$,   where $K_1\cong {\Bbb Q}(\sqrt{D})$.  
}
\end{rmk}
Let $p$ be a good prime.  Consider a localization ${\cal E}({\Bbb F}_p)$
of curve  ${\cal E}_{CM}^{(-D,f)}\cong {\cal E}(K)$ at the prime ideal ${\goth P}$ over $p$. 
It is well known, that the Frobenius endomorphism  of elliptic curve
with complex multiplication  is defined  by the  Gr\"ossencharacter;
the latter is  a complex number  $\alpha_{{\goth P}}\in {\Bbb Q}(\sqrt{-D})$
of absolute value $\sqrt{p}$.  Moreover,
multiplication of the lattice $L_{CM}={\Bbb Z}+{\Bbb Z}\tau$ by $\alpha_{{\goth P}}$
induces the Frobenius endomorphism $Fr_p^1$ on $H^1({\cal E}(K); {\Bbb Z})$,
see e.g.  [Silverman 1994] \cite{S2},  p. 174.  
Thus   one arrives  at the following matrix form for the Frobenius \&  Markov endomorphisms 
and the shift automorphism,  respectively: 
\displaymath
\left\{
\begin{array}{cc}
Fr_p^1 &=  \left(\matrix{tr~(\alpha_{{\goth P}}) & p\cr -1 & 0}\right),\\
Mk_p^1 &=  \left(\matrix{tr~(\alpha_{{\goth P}}) & p\cr 1 & 0}\right), \\
\widehat{Mk_p^1} &=  \left(\matrix{tr~(\alpha_{{\goth P}}) & 1\cr  1 & 0}\right) . 
\end{array}
\right.
\enddisplaymath
To calculate  positive integer $\pi_1(p)$ appearing in Theorem \ref{thm6.5.1},
  denote by $\left({D\over p}\right)$  the Legendre symbol of $D$ and $p$. 
    A classical result of the theory of real quadratic fields asserts that     
$\pi_1(p)$  must be one of the divisors of the integer number
\displaymath
p-\left({D\over p}\right),
\enddisplaymath
see e.g. [Hasse  1950]  \cite{HA}, p. 298.   Thus the trace of Frobenius 
endomorphism on $H^1({\cal E}(K); {\Bbb Z})$ is given by the formula
\displaymath
tr~(\alpha_{{\goth P}})=tr~(\varepsilon^{\pi_1(p)}).
\enddisplaymath
 The right hand side of 
the above equation  can be further simplified,   since 
\displaymath 
tr~(\varepsilon^{\pi_1(p)})= 2T_{\pi_1(p)}\left[~{1\over 2}~ tr~(\varepsilon)\right],
\enddisplaymath
where $T_{\pi_1(p)}(x)$ is the Chebyshev polynomial (of the first kind) 
of degree $\pi_1(p)$.  Thus  one obtains a  formula for the number of (projective) solutions
of a cubic equation over field  ${\Bbb F}_p$  in terms of invariants of
pseudo-lattice $\Lambda_1$         
\displaymath
|{\cal E}({\Bbb F}_p)|=1 + p - 2T_{\pi_1(p)}\left[ ~{1\over 2}~ tr~(\varepsilon)\right] . 
\enddisplaymath

 \index{rational elliptic curve}

\subsubsection{Rational elliptic curve}
Let  $b\ge 3$ be an integer and  consider a rational elliptic curve ${\cal E}({\Bbb Q})\subset {\Bbb C}P^2$ 
 given by the homogeneous Legendre equation
\displaymath 
y^2z=x(x-z)\left(x-{b-2\over b+2}z\right).  
\enddisplaymath
The Serre $C^*$-algebra of projective variety $V\cong {\cal E}({\Bbb Q})$  is isomorphic
(modulo an ideal) to the   Cuntz-Krieger algebra ${\cal O}_B$,  where 
\displaymath
B=\left(\matrix{b-1 & 1\cr b-2 & 1}\right),  
\enddisplaymath
see \cite{Nik15}.  Recall that ${\cal O}_B\otimes {\cal K}$ is the   crossed product 
$C^*$-algebra  of a stationary AF $C^*$-algebra by its shift automorphism, 
see [Blackadar 1986]  \cite{B},  p. 104; 
the AF $C^*$-algebra has the following dimension group
\displaymath
{\Bbb Z}^2 \buildrel B^T \over\to 
{\Bbb Z}^2\buildrel B^T \over\to 
{\Bbb Z}^2 \buildrel B^T\over\to
\dots,
\enddisplaymath
where $B^T$ is the transpose of matrix $B$.   
Because $\mu$ must be a positive eigenvalue of matrix $B^T$,   one gets  
\displaymath
\mu  =  {2-b+\sqrt{b^2-4}\over 2}. 
\enddisplaymath
Likewise,  since $\theta$ must be the corresponding 
positive eigenvector $(1,\theta)$ of the same matrix,  one gets
\displaymath
\theta =   {1\over 2}\left(\sqrt{{b+2\over b-2}}-1\right).
\enddisplaymath
Therefore,  pseudo-lattices $\Lambda_i$ are  $\Lambda_0=\Lambda_2\cong {\Bbb Z}$
and 
\displaymath
\Lambda_1\cong 
{2-b+\sqrt{b^2-4}\over 2} 
\left[{\Bbb Z}+
{1\over 2}\left(\sqrt{{b+2\over b-2}}-1\right)
{\Bbb Z}\right].
\enddisplaymath
\begin{rmk}
\textnormal{
The pseudo-lattice $\Lambda_1\subset K_1$,  where $K_1={\Bbb Q}(\sqrt{b^2-4})$.  
}
\end{rmk}
Let $p$ be a good prime and let ${\cal E}({\Bbb F}_p)$ be the reduction
of our rational  elliptic curve   modulo  $p$.  It follows from Section 6.3.3,
that  $\pi_1(p)$ as one of the divisors of integer number
\displaymath
p-\left({b^2-4\over p}\right).
\enddisplaymath
Unlike the case of complex multiplication, 
the Gr\"ossencharacter is no longer  available for ${\cal E}({\Bbb Q})$;   
yet  the trace of Frobenius endomorphism  can be computed using 
Theorem \ref{thm6.5.1}, i.e.  
\displaymath 
tr~(Fr_p^1)=tr~\left[(B^T)^{\pi_1(p)}\right].
\enddisplaymath
Using the Chebyshev polynomials,  one can write the last equation in the 
form
\displaymath
tr~(Fr_p^1)=2T_{\pi_1(p)}\left[~{1\over 2}~ tr~(B^T)\right]. 
\enddisplaymath
Since  $tr~(B^T)=b$,  one gets
\displaymath 
tr~(Fr_p^1)=2T_{\pi_1(p)}\left({b\over 2}\right). 
\enddisplaymath
Thus  one obtains a  formula for the number of  solutions
of equation  $y^2z=x(x-z)\left(x-{b-2\over b+2}z\right)$
over field  ${\Bbb F}_p$  in terms of the {\it noncommutative invariants}  of
pseudo-lattice $\Lambda_1$  of the form    
\displaymath 
|{\cal E}({\Bbb F}_p)|=1 + p - 2T_{\pi_1(p)}\left({b\over 2}\right).
\enddisplaymath
We shall conclude by a concrete  example comparing the obtained  formula with
the known results for rational elliptic curves   in the Legendre form, 
see e.g. [Hartshorne  1977]  \cite{H1},  p.  333  and  [Kirwan 1992]  \cite{KI}, 
pp. 49-50.

\bigskip
\begin{exm}\label{ex1}
{\bf (Comparison to classical invariants)}
{\normalfont
Suppose that $b\equiv 2~mod~4$.   Recall that the $j$-invariant takes the same value on 
$\lambda$, $1-\lambda$ and ${1\over\lambda}$,   see   e.g.  [Hartshorne  1977]  \cite{H1},  p.  320.
Therefore,   one can bring equation $y^2z=x(x-z)\left(x-{b-2\over b+2}z\right)$ to the form
\displaymath
y^2z=x(x-z)(x-\lambda z),    
\enddisplaymath
 where $\lambda={1\over 4}(b+2)\in \{2,3,4,\dots\}$.      Notice that for the above  curve 
\displaymath
tr~(B^T)=b=2(2\lambda-1).     
\enddisplaymath
 To calculate  $tr~(Fr_p^1)$ for our elliptic curve,   recall that in view of last equality, 
 one gets
\displaymath
  tr~(Fr_p^1)=2 ~T_{\pi_1(p)} (2\lambda-1).  
\enddisplaymath
 \index{Chebyshev polynomial}
It will be useful to express  Chebyshev polynomial   $T_{\pi_1(p)} (2\lambda-1)$   
in terms of the hypergeometric function $_2F_1(a, b; c; z)$;  the standard 
formula brings our last equation  to the form 
\displaymath
  tr~(Fr_p^1)=2 ~_2F_1(-\pi_1(p), ~\pi_1(p); ~{1\over 2}; ~1-\lambda).  
\enddisplaymath
We leave  to the reader to prove the identity
\displaymath
\begin{array}{cc}
  2 ~_2F_1(-\pi_1(p), ~\pi_1(p); ~{1\over 2}; ~1-\lambda) =&\\
   &\\
=  (-1)^{\pi_1(p)} ~_2F_1(\pi_1(p)+1,~\pi_1(p)+1; ~1; ~\lambda). &
\end{array}
\enddisplaymath
In the last formula
\displaymath 
 _2F_1(\pi_1(p)+1,~\pi_1(p)+1; ~1; ~\lambda)=
 \sum_{r=0}^{\pi_1(p)}\left(\matrix{\pi_1(p)\cr r}\right)^2 \lambda^r,  
\enddisplaymath
see  [Carlitz 1966]  \cite{Car1},  p. 328.  
  Recall that $\pi_1(p)$ is a divisor of 
 $p-\left({b^2-4\over p}\right)$,  which in our case 
 takes the value ${p-1\over 2}$.  Bringing together the above formulas,  one gets
\displaymath
|{\cal E}({\Bbb F}_p)|=1+p + (-1)^{{p-1\over 2}} 
\sum_{r=0}^{{p-1\over 2}} \left(\matrix{{p-1\over 2} \cr r}\right)^2 \lambda^r.  
\enddisplaymath
The reader is encouraged to compare  the obtained formula with the classical result
in    [Hartshorne  1977]  \cite{H1},  p.  333   and [Kirwan 1992]  \cite{KI},   pp. 49-50;  
notice also   an intriguing  relation with the   {\it Hasse  invariant}.  
 }
\end{exm}

\vskip1.5cm\noindent
{\bf Guide to the literature.}
The {\it Weil Conjectures} (WC)  were formulated in [Weil 1949]  \cite{Wei1};
along with the Langalands Program,  the WC shaped the modern look of
number theory.   The theory of {\it motives}  was elaborated by 
[Grothendieck 1968] \cite{Gro1} to solve the WC.  An excellent introduction to the WC
 can be found in   [Hartshorne 1977]  \cite{H1},  Appendix C.   
The related noncommutative invariants   were 
calculated in \cite{Nik16}.

 \index{transcendental number theory}

\section{Transcendental number  theory}
The functor
\displaymath
F: \hbox{{\bf Ell}} \longrightarrow \hbox{{\bf NC-Tor}}
\enddisplaymath
constructed in Section 5.1 has an amazing application in the 
{\it transcendental number  theory},  see e.g. [Baker  1975]  \cite{BA} for
an introduction.  Namely,  we shall use the formula
$F({\cal E}_{CM}^{(-D,f)})={\cal A}_{RM}^{(D,f)}$ obtained 
in Section 6.1 to prove that the transcendental function
\displaymath
{\cal J}(\theta,\varepsilon):=e^{2\pi i\theta+\log\log\varepsilon}
\enddisplaymath
takes algebraic values for the algebraic arguments $\theta$ and $\varepsilon$.
Moreover,  these values of  ${\cal J}(\theta,\varepsilon)$ belong to the Hilbert
class field of the imaginary quadratic field ${\Bbb Q}(\sqrt{-D})$ for all but a finite
set  of values of $D$.

 \index{Gelfond-Schneider Theorem}

\subsection{Algebraic values of transcendental functions}
Recall that   an old and difficult  problem of number theory is 
to determine  if  given irrational value of  a  transcendental function is 
algebraic or transcendental for certain algebraic arguments;  the algebraic values 
are particularly remarkable and worthy of thorough investigation,  see   [Hilbert  1902]  \cite{Hil1},  p. 456.
Only few general results are known,   see e.g.  [Baker  1975]  \cite{BA}.   We shall mention the famous 
Gelfond-Schneider Theorem saying  that $e^{\beta\log\alpha}$ is a transcendental
number,  whenever $\alpha\not\in \{0, 1\}$ is an algebraic and $\beta$ an irrational
algebraic number.  
In  contrast,   Klein's  invariant  $j(\tau)$   is known to take
algebraic values  whenever  $\tau\in {\Bbb H}:=\{x+iy\in {\Bbb C}~|~y>0\}$ is 
an imaginary quadratic number.    
In follows we shall focus on  algebraic  values  of  the  transcendental  function
\displaymath
{\cal J}(\theta,\varepsilon):=\{e^{2\pi i\theta+\log\log\varepsilon} ~|-\infty<\theta<\infty, ~1<\varepsilon<\infty\}  
\enddisplaymath
for  real   arguments  $\theta$ and $\varepsilon$. 
\begin{rmk}
\textnormal{
The  ${\cal J}(\theta,\varepsilon)$ can be viewed as an extension 
of  Klein's invariant $j(\tau)$ to the boundary of  half-plane ${\Bbb H}$;
hence  our  notation. 
}
\end{rmk}
Let  $K={\Bbb Q}(\sqrt{-D})$  be an imaginary  quadratic field
of class number  $h\ge 1$;    let  $\{{\cal E}_1,\dots {\cal E}_h\}$ be  pairwise non-isomorphic elliptic
curves with complex multiplication by the ring of integers of  field   $K$. 
For  $1\le i\le h$  we shall write   ${\cal A}_{\theta_i}=F({\cal E}_i)$ to denote
the  noncommutative torus,   where $F:$  {\bf Ell} $\to$  {\bf NC-Tor}
is the functor defined in Section 5.1.  It follows from Theorem \ref{thm6.1.2}
that each  $\theta_i$ is a quadratic irrationality of  the field ${\Bbb Q}(\sqrt{D})$.   
We shall write $\overline{(a_1^{(i)},\dots, a_n^{(i)})}$ to denote the  period of continued 
fraction for $\theta_i$ and   for each $\theta_i$ we shall consider the matrix
\displaymath
A_i=\left(\matrix{a_1^{(i)} & 1\cr 1 & 0}\right)\dots
\left(\matrix{a_n^{(i)} & 1\cr 1 & 0}\right). 
\enddisplaymath
\begin{rmk}\label{rmk6.6.2}
\textnormal{
In other words, matrix $A_i$ corresponds to the {\it shift automorphism} 
$\sigma_{A_i}$ of the dimension group 
$K_0^+({\cal A}_{\theta_i})\cong {\Bbb Z}+{\Bbb Z}\theta_i$,  see Section 3.5.2.
}
\end{rmk}
Let $\varepsilon_i>1$  be  the Perron-Frobenius eigenvalue of  matrix $A_i$;
it is easy to see, that  $\varepsilon_i$ is a quadratic irrationality of the field ${\Bbb Q}(\sqrt{D})$. 
\begin{thm}\label{thm6.6.1}
 If $D\not\in\{1,2,3,7,11,19,43,67,163\}$ is a  square-free positive integer,
 then $\{{\cal J}(\theta_i,\varepsilon_i) ~|~ 1\le i\le h\}$ are conjugate algebraic numbers.
 Moreover,   such  numbers are  generators of  the Hilbert class field  of    
 imaginary quadratic field  ${\Bbb Q}(\sqrt{-D})$.
\end{thm}

\subsection{Proof of Theorem \ref{thm6.6.1}}
The idea of proof is remarkably simple;  indeed,  recall that the system of
defining relations 
\displaymath
\left\{
\begin{array}{cc}
x_3x_1 &=  e^{2\pi i\theta}x_1x_3,\\
x_4x_2 &=  e^{2\pi i\theta}x_2x_4,\\
x_4x_1 &=  e^{-2\pi i\theta}x_1x_4,\\
x_3x_2 &=   e^{-2\pi i\theta}x_2x_3,\\
x_2x_1 &= x_1x_2=e,\\
x_4x_3 &= x_3x_4=e,
\end{array}
\right.
\enddisplaymath
for the noncommutative torus ${\cal A}_{\theta}$ and  defining relations
\displaymath
\left\{
\begin{array}{cc}
x_3x_1 &= \mu e^{2\pi i\theta}x_1x_3,\\
x_4x_2 &= {1\over \mu} e^{2\pi i\theta}x_2x_4,\\
x_4x_1 &= \mu e^{-2\pi i\theta}x_1x_4,\\
x_3x_2 &= {1\over \mu} e^{-2\pi i\theta}x_2x_3,\\
x_2x_1 &= x_1x_2,\\
x_4x_3 &= x_3x_4,
\end{array}
\right.
\enddisplaymath
for the Sklyanin $\ast$-algebra $S(q_{13})$ with $q_{13}=\mu e^{2\pi i\theta}\in {\Bbb C}$
are identical modulo the ``scaled unit relation''
\displaymath
x_1x_2=x_3x_4={1\over\mu}e,
\enddisplaymath
see Lemma \ref{lem5.1.3}.   On the other hand,  we know that if our elliptic curve
is isomorphic to ${\cal E}_{CM}^{(-D,f)}$,  then:

\medskip
(i)  $q_{13}=\mu e^{2\pi i\theta}\in K$,  where $K=k(j({\cal E}_{CM}^{(-D,f)}))$ is the
Hilbert class field of the imaginary quadratic field $k={\Bbb Q}(\sqrt{-D})$;

\smallskip
(ii) $\theta\in {\Bbb Q}(\sqrt{D})$,  since Theorem \ref{thm6.1.3} says that
$F({\cal E}_{CM}^{(-D,f)})={\cal A}_{RM}^{(D,f)}$.

\bigskip\noindent
Thus one gets the following inclusion
\displaymath
K\ni \mu  e^{2\pi i\theta},  \quad\hbox{where} ~ \theta\in {\Bbb Q}(\sqrt{D}). 
\enddisplaymath
The only missing piece of data is the constant $\mu$;  however, the following
lemma fills in the gap. 
\begin{lem}\label{lem6.6.1}
For each noncommutative torus ${\cal A}_{RM}^{(D,f)}$ the constant  $\mu=\log\varepsilon$,
where $\varepsilon>1$ is the Perron-Frobenius eigenvalue of positive integer
matrix
\displaymath
A=\left(\matrix{a_1 & 1\cr 1 & 0}\right)\dots \left(\matrix{a_n & 1\cr 1 & 0}\right)
\enddisplaymath
 and  $\overline{(a_1,\dots,a_n)}$ is the period of continued fraction
 of  the corresponding quad\-ratic irrationality $\theta$,  see also Remark \ref{rmk6.6.2}.    
\end{lem}
{\it Proof.} 
(i) Recall that the range of the canonical trace $\tau$ on projections 
of algebra ${\cal A}_{\theta}\otimes {\cal K}$ is given by  pseudo-lattice
$\Lambda={\Bbb Z}+{\Bbb Z}\theta$ [Rieffel 1990] \cite{Rie1}, p. 195.  
Because $\tau({1\over\mu}e)={1\over\mu}\tau(e)={1\over\mu}$, the pseudo-lattice
corresponding to the algebra ${\cal A}_{\theta}$ with a scaled unit can be written
as $\Lambda_{\mu}=\mu({\Bbb Z}+{\Bbb Z}\theta)$.

\bigskip
(ii)  To express $\mu$ in terms of the inner  invariants of pseudo-lattice $\Lambda$,
denote by $R$ the ring of endomorphisms of $\Lambda$ and by $U_R\subset R$
the multiplicative group of automorphisms (units) of $\Lambda$.  For each 
$\varepsilon, \varepsilon'\in U_R$ it must hold
\displaymath
\mu(\varepsilon\varepsilon' \Lambda)=\mu(\varepsilon\varepsilon') \Lambda=
\mu(\varepsilon)\Lambda+\mu(\varepsilon')\Lambda,
\enddisplaymath
since $\mu$ is an additive functional on the pseudo-lattice $\Lambda$.    
Canceling $\Lambda$ in the above equation,  one gets
\displaymath
\mu(\varepsilon\varepsilon')=\mu(\varepsilon)+\mu(\varepsilon'),
\qquad \forall \varepsilon, \varepsilon'\in U_R.
\enddisplaymath
The only real-valued function on $U_R$ with such a property
is the logarithmic function (a {\it regulator} of $U_R$);  thus $\mu(\varepsilon)=\log\varepsilon$.

\bigskip
(iii)   Notice  that $U_R$ is generated by a single element $\varepsilon$.  To calculate the 
generator,  recall that pseudo-lattice $\Lambda={\Bbb Z}+{\Bbb Z}\theta$ is isomorphic 
to a pseudo-lattice $\Lambda'={\Bbb Z}+{\Bbb Z}\theta'$, where $\theta'=\overline{(a_1,\dots,a_n)}$
is purely periodic continued fraction and $(a_1,\dots, a_n)$ is the  period of 
continued fraction of $\theta$.  From the standard facts of the theory of 
continued fractions,  one gets that $\varepsilon$ coincides with the Perron-Frobenius 
eigenvalue of matrix  
\displaymath
A=\left(\matrix{a_1 & 1\cr 1 & 0}\right)\dots \left(\matrix{a_n & 1\cr 1 & 0}\right).
\enddisplaymath
Clearly,  $\varepsilon>1$ and it is an invariant of the stable isomorphism class 
of algebra  ${\cal A}_{\theta}$.     Lemma \ref{lem6.6.1} is proved.
 $\square$

\begin{rmk}
{\bf (Second proof of Lemma \ref{lem6.6.1})}
\textnormal{
Lemma \ref{lem6.6.1}  follows from a purely measure-theoretic argument.  
Indeed,  if $h_x: {\Bbb R}\to {\Bbb R}$ is a ``stretch-out'' automorphism 
of real line ${\Bbb R}$ given by the formula $t\mapsto tx,~\forall t\in {\Bbb R}$,
then the only $h_x$-invariant measure $\mu$ on ${\Bbb R}$ is the ``scale-back''
measure $d\mu={1\over t} dt$.  Taking the antiderivative and integrating 
between $t_0=1$ and $t_1=x$,  one gets
\displaymath
\mu=\log x.
\enddisplaymath
It remains to notice that for pseudo-lattice 
$K_0^+({\cal A}_{RM}^{(D,f)})\cong {\Bbb Z}+{\Bbb Z}\theta\subset {\Bbb R}$,
the automorphism $h_x$ corresponds to $x=\varepsilon$,  where $\varepsilon>1$
is the Perron-Frobenius eigenvalue of matrix $A$.  
Lemma \ref{lem6.6.1} follows. $\square$. 
}
\end{rmk}

\bigskip
Theorem \ref{thm6.6.1} is an implication of the following  argument.

\bigskip
(i)  Let  $D\not\in\{1,2,3,7,11,19,43, 67,163\}$ be a positive square-free integer;  for the sake of simplicity,  let  
$f=1$,  i.e. the endomorphism ring of ${\cal E}_{CM}^{(-D,f)}$ coincides with the ring of integers $O_k$ of the 
imaginary quadratic field $k={\Bbb Q}(\sqrt{-D})$.    In  this case  ${\cal E}_{CM}^{(-D, f)}\cong {\cal E}(K)$, 
 where $K=k(j({\cal E}(K)))$ is the Hilbert class field  of $k$.  
It follows from the well-known facts of complex multiplication, that condition  
$D\not\in\{1,2,3,7,11,19,43, 67,163\}$ guarantees  that $K\not\cong {\Bbb Q}$,
i.e. the field $K$ has complex embedding.

\bigskip
(ii) Let ${\cal E}_1(K),\dots, {\cal E}_h(K)$ be pairwise  non-isomorphic curves  with
\linebreak
$End~({\cal E}_i(K))\cong O_k$,  where $h$ is the class number of  $O_k$.  
Repeating for each ${\cal E} _i (K)$  the argument at the beginning of proof of  
Theorem \ref{thm6.6.1},  we conclude   that $\mu_i e^{2\pi i\theta_i}\in K$.

\bigskip
(iii)  But Lemma \ref{lem6.6.1} says that  $\mu_i=\log\varepsilon_i$; 
thus for each $1\le i\le h$,   one gets 
\displaymath
(\log\varepsilon_i) ~e^{2\pi i\theta_i}=e^{2\pi i\theta_i + \log\log\varepsilon_i}={\cal J}(\theta_i, \varepsilon_i) \in K.
\enddisplaymath

\bigskip
(iv)  The transitive action of the ideal class group $Cl~(k)\cong Gal~(K|k)$ on the elliptic curves ${\cal E}_i(K)$
extends to the algebraic numbers ${\cal J}(\theta_i, \varepsilon_i)$;   thus ${\cal J}(\theta_i, \varepsilon_i)\in K$
are algebraically conjugate.

\bigskip
Theorem \ref{thm6.6.1} is proved.
$\square$

\vskip1cm\noindent
{\bf Guide to the literature.}
The complex number is {\it algebraic} whenever it is the root of a 
polynomial with integer coefficients;  it is notoriously hard to tell 
if given complex number is algebraic or not.  Thanks to Ch.~Hermite 
and C.~L.~Lindemann   the $e$ and $\pi$ are not,  but even for $e\pm\pi$
the answer is unknown.  The Seventh Hilbert Problem deals with such type
of questions,  see   [Hilbert  1902]  \cite{Hil1},  p. 456.
The famous  Gelfond-Schneider Theorem says  that $e^{\beta\log\alpha}$ is a transcendental
number,  whenever $\alpha\not\in \{0, 1\}$ is an algebraic and $\beta$ an irrational
algebraic number.   An excellent introduction to the theory of transcendental
numbers in the book by [Baker  1975]  \cite{BA}.  
The noncommutative invariants in transcendence theory were the subject of 
\cite{Nik17}.

 \index{Seventh Hilbert Problem}

\section*{Exercises, problems and conjectures}

\begin{enumerate}

\item
Prove that   elliptic curves ${\cal E}_{\tau}$ and  ${\cal E}_{\tau'}$ are isogenous  
 if and only if 
\displaymath
\tau'={a\tau+b\over c\tau+d} \quad \hbox{for some matrix} 
 \quad\left(\matrix{a & b\cr c & d}\right)  \in M_2({\Bbb Z})
 \quad \hbox{with} \quad ad-bc>0.
 \enddisplaymath
 (Hint:  notice that $z\mapsto\alpha z$
is an invertible holomorphic map for each $\alpha\in {\Bbb C}-\{0\}$.)

\item
Prove that typically  $End~({\cal E}_{\tau})\cong {\Bbb Z}$,
i.e. the only endomorphisms of ${\cal E}_{\tau}$ are the multiplication-by-$m$ 
endomorphisms.

\item 
Prove that  for a countable set of $\tau$
\displaymath
End~({\cal E}_{\tau})\cong {\Bbb Z}+ fO_k,
\enddisplaymath
where $k={\Bbb Q}(\sqrt{-D})$ is an imaginary quadratic field,  $O_k$ its ring
of integers and $f\ge 1$ is the conductor of a finite index subring of $O_k$;
prove that in such a case $\tau\in End~({\cal E}_{\tau})$,   i.e. complex modulus 
itself is an imaginary quadratic number.

\item
Show that 
the noncommutative tori ${\cal A}_{\theta}$ and  ${\cal A}_{\theta'}$ are stably
homomorphic   if and only if 
\displaymath
\theta'={a\theta+b\over c\theta+d} \quad \hbox{for some matrix} 
 \quad\left(\matrix{a & b\cr c & d}\right)  \in M_2({\Bbb Z})
\quad \hbox{with} \quad ad-bc>0.
 \enddisplaymath
(Hint:  follow and modify the argument of  [Rieffel 1981] \cite{Rie2}.)

 \index{rank conjecture}

\item
{\bf (The rank conjecture)}
Prove that the formula  
\displaymath
rk~({\cal E}_{CM}^{(-D,f)}) +1=c({\cal A}_{RM}^{(D,f)})
\enddisplaymath
is true in general,  i.e. for all  $D\ge 2$ and $f\ge 1$.

\item
Verify that the continued fraction 
$$[3s+1,\overline{2,  1,  3s ,  1,  2,  6s+2}]=\sqrt{(3s+1)^2+2s+1}$$
is a restriction of the continued fraction  
$$[x_0,\overline{x_1,  2x_1,  x_0, ,  2x_1, x_1,  2x_0}]=\sqrt{x_0^2+4nx_1+2}$$
with  $x_0=n(2x_1^2+1)+x_1$   to the case $x_1=-1$  and $n=s+1$.

\item
Assume that  a solution 
$$[x_0,\overline{x_1,\dots, x_{k-1}, x_k, x_{k-1}, \dots, x_1, 2x_0}]$$
of the diophantine equation  in Definition \ref{dfn6.2.2}
with the (culminating or almost-culminating) period $P_0\equiv 7~\mod~8$ has dimension $d=1$;
prove that for the solution 
$$[x_0,\overline{y_1, x_1,\dots, x_{k-1}, y_{k-1},  x_k, y_{k-1},  x_{k-1}, \dots, x_1, y_1, 2x_0}]$$
  having the  period $P_0+4$ the  dimension remains
 the same, i.e. $d=1$.   (Hint:  use  the same argument as in Lemma \ref{lem6.2.3}.)

\item
Prove that part (ii) of Theorem \ref{thm6.3.1} implies its part (i). 
(Hint:  repeat the step by step argument of Section 6.3.2.)

\item
Prove Remark \ref{rmk6.4.1}, i.e.   that  if  noncommutative torus 
${\cal A}_{\Theta_0}^{2n}$  has real multiplication,  then $\theta_i$
are algebraic integers.  (Hint:  each endomorphism
of   $K_0^+({\cal A}_{\Theta_0}^{2n})\cong {\Bbb Z}+\theta_1{\Bbb Z}+\dots+\theta_n{\Bbb Z}
+\sum_{i=n+1}^{2^{2n-1}}p_i(\theta){\Bbb Z}$ is multiplication by a real number;  
thus the endomorphism is described by an integer matrix,  which 
defines a polynomial equation involving $\theta_i$.)

\item
{\bf (Langlands conjecture for noncommutative tori)}
Prove Conjecture \ref{cnj6.4.1}, i.e. that 
for each  finite extension $E$  of the field of rational numbers ${\Bbb Q}$
with the  Galois group $Gal~(E |{\Bbb Q})$ and each  
 irreducible representation 
 \displaymath
 \sigma_{n+1}: Gal~(E | {\Bbb Q})\to GL_{n+1}({\Bbb C}),
\enddisplaymath
there exists a $2n$-dimensional  noncommutative torus with 
real  multiplication  ${\cal A}_{RM}^{2n}$,   such that
\displaymath
L(\sigma_{n+1},s)\equiv L({\cal A}_{RM}^{2n}, s),
\enddisplaymath
where $L(\sigma_{n+1}, s)$ is the Artin $L$-function attached to representation 
$\sigma_{n+1}$   and  $L({\cal A}_{RM}^{2n}, s)$ is  the  $L$-function of  the
noncommutative torus  ${\cal A}_{RM}^{2n}$.

\item
Prove the identity
\displaymath
\begin{array}{cc}
  2 ~_2F_1(-\pi_1(p), ~\pi_1(p); ~{1\over 2}; ~1-\lambda) =&\\
   &\\
=  (-1)^{\pi_1(p)} ~_2F_1(\pi_1(p)+1,~\pi_1(p)+1; ~1; ~\lambda), &
\end{array}
\enddisplaymath
where $_2F_1(a, b; c; z)$ is the hypergeometric function.

\end{enumerate}






\part{BRIEF SURVEY OF NCG}

\chapter{Finite Geometries}
Let $D$ be a {\it skew field},  i.e. the associative division ring. 
It is long known that  some results of classical projective geometry (e.g. the Desargues Theorem)
are true in projective spaces ${\Bbb P}^n(D)$ for any  $D$ provided $n\ge 3$; 
the other (e.g. the Pappus Theorem) hold only when $D$ is commutative,
i.e. a field.    In other words,  the algebraic structure of $D$ is reflected in geometry 
of the corresponding projective space ${\Bbb P}^n(D)$. 
In general,  there exists a bijection between the so-called {\it rational identities}
in the skew field $D$ and {\it configurations} (of  points, lines, etc.) in the corresponding 
projective space ${\Bbb P}^n(D)$.     We shall call projective space
 ${\Bbb P}^n(D)$ a {\it finite geometry},    since there exist concrete examples
of such spaces consisting of only a finite set of points, lines, etc.;  this case 
corresponds to  $D\cong {\Bbb F}_p$ being a finite field.  
We refer the reader to the classical book  [Hilbert 1930]  \cite{HIL} for a detailed
account of this field of mathematics.  For a more algebraic approach to the finite
geometries,  see the monographs  [Artin 1957]  \cite{ART}  and  [Baer 1952]  \cite{BAE}.
Our exposition follows the guidelines of  [Shafarevich 1990]  \cite{SHA},   pp. 84-87.  

 \index{skew field}
 \index{finite geometry}
 \index{projective space ${\Bbb P}^n(D)$}
 \index{rational identity}

 \index{configurations}

\section{Axioms of projective geometry}
It is well known that all statements of the elementary geometry 
can be deduced from  a finite set of the axioms;  for the sake of
clarity,  we shall focus on the plane geometry and the following 
set of axioms reflecting the incidence and parallelism 
of points and lines in the plane:

\medskip
(a) Through any two distinct points there is one and only
one line;

\smallskip
(b)  Given any line and a point not on it,  there exists one and only
one line through the point and not intersecting the line, i.e. 
parallel to it;  

\smallskip
(c) There exist three points not on any line.

\bigskip\noindent
One can wonder what realizations of the chosen set of axioms are 
possible;  in other words, do there exists a system of objects, different
from the ordinary plane geometry of points and lines,  for which the
set of axioms (a)-(c) are  satisfied?  The following example indicates 
that such realizations exist and can be in a stark contrast to our intuition.     
\begin{exm}
\textnormal{
Consider a model with four points $A,B,C$ and $D$ and six lines
$AB, BC, CD, DA, AC$ and $BD$,  see Fig. 7.1  The reader can verify that the
axioms (a)-(c) are satisfied.   The parallel lines are:  $AB ~|| ~CD$, $BC ~||~ AD$ and,
counter-intuitively, $AC ~|| ~BD$.     
}
\end{exm}
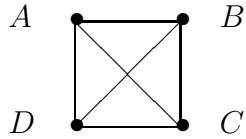
\begin{figure}[here]
\begin{picture}(500,110)(50,-5)

\put(225,80){\line(1,0){40}}
\put(265,80){\line(0,-1){40}}
\put(265,40){\line(-1,0){40}}
\put(225,40){\line(0,1){40}}
\put(225,80){\line(1,-1){40}}
\put(225,40){\line(1,1){40}}

\put(223,78){$\bullet$}
\put(263,78){$\bullet$}
\put(263,38){$\bullet$}
\put(223,38){$\bullet$}

\put(200,78){$A$}
\put(280,78){$B$}
\put(280,38){$C$}
\put(200,38){$D$}

\end{picture}

\caption{Finite geometry on  four points and six lines.}
\end{figure}
Following the seminal idea of  R\'en\'e Descartes,  one can introduce a coordinate system $(X,Y)$ in the plane
and  write the points of a line as solutions to the linear equation 
\displaymath
aX+bY=c,
\enddisplaymath
for some constants $a,b,c$ in a skew field $D$;  the axioms (a)-(c)  will reduce to a system of algebraic 
equations over $D$.  
\begin{exm}
\textnormal{
Let $D\cong {\Bbb F}_2$ be the field of characteristic $2$. 
The reader can verify that finite geometry on four points and six lines 
can be written  as 
\displaymath
\left\{
\begin{array}{ccc}
A &=& (0,0)\\ 
B &=& (0,1)\\ 
C&=&  (1,0)\\
D &=& (1,1),                   
\end{array}
\right.
\enddisplaymath
and 
\displaymath
\left\{
\begin{array}{ccc}
AB  &:& X=0\\ 
CD &: & X=1\\ 
AD &:&  X+Y=0\\
BC &: & X+Y=1\\ 
AC &:&  Y=0\\
BD &:& Y=1.                   
\end{array}
\right.
\enddisplaymath
}
\end{exm}

 \index{projective space ${\Bbb P}^n(D)$}

\section{Projective spaces over skew fields}
Let $D$ be a skew field;  consider its matrix ring $M_n(D)$.  It is known that the 
left ideals of $M_n(D)$ are in a one-to-one correspondence with subspaces $V_i$
of the $n$-dimensional space $M_n(D^{op})$,  where $D^{op}$ is the opposite 
skew field of $D$.  The set of all $V_i$ has a natural (partial) order structure coming
from the inclusions of $V_i$ into each other;  one can wonder if it is possible to recover 
the matrix ring $M_n(D)$, its dimension $n$ and the skew field $D$ itself from the 
partially ordered set $\{V_i\}$?  Note that the set $\{V_i\}$ 
corresponds to the set of all linear subspaces of the projective space ${\Bbb P}^{n-1} (D)$
and, therefore, our question  reduces  to the one about  the  ``axiomatic structure''
of the projective space   ${\Bbb P}^{n-1} (D)$.  Thus one arrives at the following 
\begin{dfn}\label{dfn7.2.1} 
By a projective space one understands a partially ordered set $({\cal P}, \ge)$,  which satisfies
the following conditions:

\medskip
(i)  For any set of elements $x_{\alpha}\in {\cal P}$ there exists an element $y$ such that
$y\ge x_{\alpha}$ for all $\alpha$;  if $z\ge x_{\alpha}$,  then $z\ge y$.  The element $y$ 
is called the sum of the $x_{\alpha}$ and is denoted by $\cup x_{\alpha}$.  The sum of 
all $x\in {\cal P}$ exists and is called the whole projective space $I({\cal P})$.

\smallskip
(ii)  For any set of elements $x_{\alpha}\in {\cal P}$ there exists an element $y'$ such
that $y'\le x_{\alpha}$ for all $\alpha$;   if $z'\le x_{\alpha}$, then $z'\le y'$.  The element
$y'$ is called the intersection of the $x_{\alpha}$ and is denoted by    $\cap x_{\alpha}$.
The intersection of all $x\in {\cal P}$ exists and is called the empty set $\emptyset ({\cal P})$.

\smallskip
(iii)  For any $x,y\in {\cal P}$ and $a\in x/y$ there exists an element $b\in x/y$ such that
$a\cup b=I(x/y)$ and $a\cap b=\emptyset (x/y)$, where $x/y$ is the partially ordered set 
of all $z\in {\cal P}$, such that $y\le z\le x$.  If $b'\in x/y$ is another element with the same
properties and if $b\le b'$,  then $b=b'$.

\smallskip
(iv)  The length of all chains $a_1\le a_2\le\dots\le a_r$ with $a_1\ne a_2\ne a_3\ne\dots\ne a_r$
is bounded.  

\smallskip
(v) The element $a\in {\cal P}$ is called a point if $b\le a$ and $b\ne a$ implies that
$b=\emptyset ({\cal P})$.     For any two points $a$ and $b$ there exists a point $c$,
such that $c\ne a, c\ne b$ and $c\le a\cup b$.       
 \end{dfn}
 \index{dimension function}
\begin{dfn}\label{dfn7.2.2}
By a dimension function on the projective space $({\cal P}, \ge)$ one understands 
the function $d(a)$ equal to the maximum length of a chain $\emptyset\le\dots\le a$;
such a function satisfies the equality 
\displaymath
d(a\cap b)+d(a\cup b)=d(a)+d(b).
\enddisplaymath
for all $a,b\in {\cal P}$.  The number $d(I({\cal P}))$ is called the dimension 
of projective space $({\cal P}, \ge)$. 
\end{dfn}
\begin{exm}
\textnormal{
If  $({\cal P},\ge)$ is the projective space of all linear subspaces of the space 
${\Bbb P}^n(D)$,  then $n=d(I({\cal P}))$.   
}
\end{exm}
The following basic result says that the skew field $D$ can be recovered from the
corresponding partially ordered set ${\Bbb P}^n(D)$;  in other words,  given projective geometry
determines the  skew field $D$. 
\begin{thm}\label{thm7.2.1}
{\bf (Fundamental Theorem of Projective Geometry)}

\medskip
(i) If  $n\ge 2$ and $({\cal P},\ge)$ is a  projective space
 coming from the linear subspaces of the space ${\Bbb P}^n(D)$
for a skew field $D$,   then  $({\cal P},\ge)$  determines the number $n$ and 
the skew field $D$;

\smallskip
(ii)  If  $({\cal P},\ge)$ is an arbitrary projective space of dimension $n\ge 3$,
then it is isomorphic to the space ${\Bbb P}^n(D)$ over some skew field $D$.  
\end{thm}
 \index{fundamental theorem of projective geometry}
\begin{rmk}
\textnormal{
Not every projective space $({\cal P},\ge)$ of dimension $2$ is isomorphic to ${\Bbb P}^2(D)$, 
unless one adds an extra axiom (called the Desargues axiom)  to the list
defining the projective space. 
}
\end{rmk}

\section{Desargues and Pappus axioms}
We shall remind the following two  results of elementary geometry in the plane. 
\begin{thm}
{\bf (Desargues)}
If the three lines $AA', BB'$ and $CC'$ joining corresponding vertices of two
triangles $ABC$ and $A'B'C'$ intersect in a point $O$,  then the points of
intersection of the corresponding sides are collinear.   
\end{thm}
 \index{Desargues axiom}
\begin{rmk}
\textnormal{
We encourage the reader to graph  the Deasargues Theorem. 
 }
\end{rmk}
\begin{thm}
{\bf (Pappus)}
If the vertices of a hexagon $P_1,\dots, P_6$ lie three by three on two lines,
then the points of intersection of the opposite sides $P_1P_2$ and $P_4P_5$,
$P_2P_3$ and $P_5P_6$,  $P_3P_4$ and $P_6P_1$ are collinear, see Fig. 7.2. 
\end{thm}
\begin{figure}[here]
\begin{picture}(500,110)(50,15)

\put(225,80){\line(2,1){65}}
\put(225,80){\line(2,-1){70}}

\put(235,85){\line(1,-1){21}}
\put(245,90){\line(1,-1){39}}

\put(245,90){\line(-1,-2){8}}

\put(252,65){\line(1,2){18}}

\put(223,78){$\bullet$}
\put(232,82){$\bullet$}
\put(242,87){$\bullet$}
\put(250,63){$\bullet$}
\put(280,49){$\bullet$}
\put(234,70){$\bullet$}
\put(267,98){$\bullet$}

\put(223,93){$P_1$}
\put(235,100){$P_5$}
\put(255,110){$P_3$}
\put(250,48){$P_2$}
\put(280,34){$P_6$}
\put(234,55){$P_4$}

\end{picture}

\caption{Pappus Theorem.}
\end{figure}
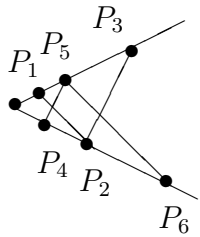
 \index{Desargues projective plane}
\begin{dfn}
By a Desargues projective plane one understands a two-dimensional projective plane 
$({\cal P},\ge)$  which satisfies the Desargues Theorem.    
 \end{dfn}
 \index{Pappus axiom}
 \index{Pappus projective space}
\begin{thm}
The Desargues  projective plane   $({\cal P},\ge)$ is isomorphic to ${\Bbb P}^2(D)$ for 
some skew field $D$.  
\end{thm}
\begin{dfn}
By a Pappus projective space  one understands an $n$-di\-men\-sional Desargues projective space
$({\cal P},\ge)$  which satisfies the Pappus Theorem.    
 \end{dfn}
The following beautiful result links geometry of projective spaces  to algebraic properties
of  the skew field $D$.  
\begin{thm}
 The $n$-dimensional Desargues projective space   $({\cal P},\ge)\cong {\Bbb P}^n(D)$ is
isomorphic to a Pappus projective space if and only if   $D$ is commutative,  i.e. a field.  
\end{thm}
 \index{commutative Desargues space}

\vskip1cm\noindent
{\bf Guide to the literature.}
The subject of finite geometries is  covered in the classical book  [Hilbert 1930]  \cite{HIL}.
We refer the reader to the monographs  [Artin 1957]  \cite{ART}  and  [Baer 1952]  \cite{BAE}
for a more algebraic treatment.  Our exposition follows the guidelines of  
[Shafarevich 1990]  \cite{SHA},   pp. 84-87.






\chapter{Continuous Geometries}
Roughly  speaking, {\it  continuous geometry}  is the generalization of  finite geometry ${\Bbb P}^n(D)$
to the case $n=\infty$ and $D\cong {\Bbb C}$;  the ${\Bbb P}^{\infty}({\Bbb C})$ is identified 
with the ring of bounded linear on a Hilbert space ${\cal H}$ known as a {\it factor} of the
{\it von Neumann algebra} ${\cal A}$ (or {\it $W^*$-algebra} for brevity).  
The latter  is a $\ast$-subalgebra  of the algebra $B({\cal H})$ of all bounded linear
operators on  ${\cal H}$,  such that its topological closure
can be expressed in purely algebraic terms (as a double commutant of the
$\ast$-subalgebra).  The topology is the weak topology;  therefore 
$W^*$-algebras can be viewed as the $C^*$-algebras with additional properties.
The  $W^*$- algebras and factors are reviewed in Section 8.1,    and  
corresponding continuous geometries are reviewed in Section 8.2. 
The $W^*$-algebras were introduced in the seminal paper by [Murray \& von Neumann 1936]   \cite{MuNeu1}.
We refer the reader to the textbooks  [Schwartz 1967]   \cite{SCH}   and [Fillmore 1996] \cite{F},  Chapter 4
for a detailed account  of the area.  
Continuous geometries were introduced and studied in the monograph [von Neumann 1960]  \cite{NEU}. 

 \index{continuous geometry}
 \index{$W^*$-algebra}
 \index{weak topology}

\section{$W^*$-algebras}
\begin{dfn}
The weak topology on the space $B({\cal H})$ is defined by the semi-norm
$||\bullet ||$ on each $S\in B({\cal H})$ given by the formula 
\displaymath
||S||=|(Sx, y)|
\enddisplaymath
 as $x$ and $y$ run through all of ${\cal H}$.  
\end{dfn}
\begin{dfn}
By a commutant of a subset ${\cal S}\subset B({\cal H})$ one understands  a unital
subalgebra ${\cal S}'$ of $B({\cal H})$ given by the formula
\displaymath
{\cal S}'=\{T\in B({\cal H}) ~| ~ TS=ST ~\hbox{for all} ~S\in {\cal S}\}. 
\enddisplaymath
The double commutant  of ${\cal S}$ is denoted by ${\cal S}'':=({\cal S}')'$.   
\end{dfn}
 \index{double commutant}
\begin{rmk}
\textnormal{
The reader can verify that ${\cal S}'$ is closed in the weak operator topology. 
}
\end{rmk}
\begin{dfn}
Let ${\cal S}$ be a subset of $B({\cal H})$.  Then:

\medskip
(i)  {\bf Alg} ${\cal S}$  denotes the subalgebra of $B({\cal H})$ generated by finite linear
combinations of the finite products of elements of ${\cal S}$;

\smallskip
(ii)  ${\cal S}$ is called self-adjoint if for each $S\in {\cal S}$ we have $S^*\in {\cal S}$;

\smallskip
(iii) ${\cal S}$ is called non-degenerate if the only $x\in {\cal H}$ for which $Sx=0$
for all $S\in {\cal S}$ is $x=0$.
\end{dfn}
\begin{thm}
{\bf (Murray-von Neumann)}
If  ${\cal S}\subset B({\cal H})$ is  self-adjoint and non-degenerate,  then {\bf Alg} ${\cal S}$
is dense in ${\cal S}''$.  
\end{thm}
 \index{Double Commutant Theorem}
\begin{cor}
{\bf (Double Commutant Theorem)}
For any unital $\ast$-algebra ${\cal A}\subset B({\cal H})$ the double commutant 
${\cal A}''$ coincides with the weak closure {\bf Clos} ${\cal A}$ of the algebra ${\cal A}$. 
\end{cor}
\begin{dfn}
By a (concrete) $W^*$-algebra one understands a unital $\ast$-subalgebra ${\cal A}$ of $B({\cal H})$,
such that 
\displaymath
{\cal A}\cong {\cal A}''\cong  \hbox{{\bf Clos}} ~{\cal A}. 
\enddisplaymath
 \end{dfn}
\begin{exm}
\textnormal{
The algebra $B({\cal H})$ is closed in the weak topology, thus $B({\cal H})$ is a 
$W^*$-algebra.  
}
\end{exm}
\begin{exm}
\textnormal{
$M_n({\Bbb C})$ is a finite-dimensional $W^*$-algebra,  because it is a $C^*$-algebra and 
the norm topology and the weak topology coincide for the linear operators on ${\cal H}\cong {\Bbb C}^n$. 
}
\end{exm}
\begin{exm}
\textnormal{
If $G$ is a group and $g\mapsto u_g$ is a unitary representation of $G$,
then the commutant $\{u_g\}'$ is a $W^*$-algebra.  
}
\end{exm}
\begin{exm}
\textnormal{
Let $S^1$ be the unit circle and $L^{\infty}(S^1)$ the algebra of measurable
functions on $S^1$;  consider the crossed product $L^{\infty}(S^1)\rtimes {\Bbb Z}$ 
by the automorphisms of  $L^{\infty}(S^1)$ induced by rotation of $S^1$ through an
angle $\alpha$.  The $L^{\infty}(S^1)\rtimes {\Bbb Z}$ is a $W^*$-algebra.  
}
\end{exm}
An analogue of the AF-algebras for the $W^*$-algebras is given by the
following definition.
\begin{dfn}
A $W^*$-algebra ${\cal A}$ is called hyper-finite if there exists an increasing sequence 
${\cal A}_n$ of finite-dimensional $W^*$-subalgebras of ${\cal A}$ whose union 
is weakly dense in ${\cal A}$.  
\end{dfn}
 \index{hyper-finite $W^*$-algebra}
\begin{exm}
\textnormal{
The $W^*$-algebra $L^{\infty}(S^1)\rtimes {\Bbb Z}$ is hyper-finite. 
 }
\end{exm}
\begin{dfn}
A $W^*$-algebra ${\cal A}$ is called a factor if the center of ${\cal A}$ (i.e. abelian algebra $Z({\cal A})\subset {\cal A}$
which commutes with everything in ${\cal A}$)   is isomorphic to
the scalar multiple of the identity operator. 
 \end{dfn}
 \index{factor}
\begin{exm}
\textnormal{
The $W^*$-algebra $B({\cal H})$ is a factor.
}
\end{exm}
\begin{exm}
\textnormal{
The $W^*$-algebra $L^{\infty}(S^1)\rtimes {\Bbb Z}$ is a factor.  
}
\end{exm}
Each $W^*$-algebra decomposes into a ``direct integral'' of factors;  to classify the latter,
one needs the following definition.  
\begin{dfn}
A pair of projections $p,q\in {\cal A}$ is written $p\le q$ if there exists a partial isometry $u\in {\cal A}$
such that $p=uu^*$ and $u^*uq=u^*u$;  they say that $p\sim q$ are equivalent if $p=uu^*$ and
$q=u^*u$.  
  \end{dfn}
\begin{thm}\label{thm8.1.2}
{\bf (Murray-von Neumann)}
If ${\cal A}$ is a factor and $p, q$ are projections in ${\cal A}$,  then
either $p\le q$ or $q\le p$.   
\end{thm}
\begin{dfn}
If $p\in {\cal A}$ is a projection,  then 

\medskip
(i) $p$ is called finite if $p\sim q\le  p$ implies $q=p$ for all projections $q\in {\cal A}$;

\smallskip
(ii) $p$ is called infinite if there is projection $q\sim p$ such that $p<q$;

\smallskip
(iii) $p\ne 0$ is called minimal if $p$ dominates no other projection in ${\cal A}$.    
\end{dfn}
 \index{type {\bf I}, {\bf II}$_1$, {\bf II}$_{\infty}$,  {\bf III}-factors}
\begin{dfn}
The factor ${\cal A}$ is said belong to the type:

\medskip
(i) {\bf I} if ${\cal A}$ has a minimal projection;

\smallskip
(ii) {\bf II}$_1$ if ${\cal A}$ has no minimal projections and every projection is finite;

\smallskip
(iii)  {\bf II}$_{\infty}$ if ${\cal A}$ has no minimal projections but ${\cal A}$ has both finite
and infinite projections;

 \smallskip
 (iv)  {\bf III} if ${\cal A}$ has no finite projections except $0$.  
\end{dfn}
\begin{rmk}
\textnormal{
The factors of types  {\bf I}, {\bf II}$_1$,  {\bf II}$_{\infty}$ and  {\bf III}
can be realized  as the concrete $W^*$-algebras. 
}
\end{rmk}
\begin{exm}
\textnormal{
The $W^*$-algebra $L^{\infty}(S^1)\rtimes {\Bbb Z}$ is a  type {\bf II}$_1$
 factor whenever the angle $\alpha$ is irrational.  
}
\end{exm}
\begin{thm}
{\bf (Murray-von Neumann)}
All hyper-finite type {\bf II}$_1$ factors are isomorphic.  
\end{thm}
\begin{rmk}
\textnormal{
Unlike the irrational rotation $C^*$-algebras,  the corresponding 
$W^*$-algebras $L^{\infty}(S^1)\rtimes {\Bbb Z}$ are isomorphic to each other;  
it follows from the Murray-von Neumann Theorem. 
}
\end{rmk}

 \index{von Neumann geometry}

\section{Von Neumann  geometry}
Murray and von ~Neumann observed that projections in the type {\bf I}$_n$,  {\bf I}$_{\infty}$,
{\bf II}$_1$,  {\bf II}$_{\infty}$ and {\bf III}-factors in $W^*$-algebras are partially ordered
sets $({\cal P}, \ge)$,  see Theorem \ref{thm8.1.2};   the set   $({\cal P}, \ge)$ satisfies 
the  axioms (i)-(v)  of projective space ${\Bbb P}^n(D)$,  see  Definition \ref{dfn7.2.1}.   
However,  the range of {\it dimension function} $d$ introduced in Definition \ref{dfn7.2.2}  
is no longer a discrete but a {\it continuous}   subset of ${\Bbb R}$ unless ${\Bbb P}^n(D)$
is the type {\bf I}$_n$ factor. 
According to  Theorem \ref{thm7.2.1}   the  type {\bf I}$_{\infty}$
{\bf II}$_1$,  {\bf II}$_{\infty}$ and {\bf III}-factors  correspond to the infinite-dimensional 
projective spaces  ${\Bbb P}^{\infty}({\Bbb C})$;  the ${\Bbb P}^{\infty}({\Bbb C})$ will be
called a {\it von Neumann geometry}. 
\begin{thm}
{\bf (J.~von ~Neumann)}
The von Neumann  geometry  belong to one of the following types:

\medskip
(i)  {\bf I}$_n$ if the dimension function takes values $0,1, 2, \dots, n$ and the corresponding
factor is isomorphic to the ring $M_n({\Bbb C})$;

\smallskip
(ii) {\bf I}$_{\infty}$ if the dimension function takes values $0,1, 2, \dots, \infty$ and the 
corresponding factor is isomorphic to the ring ${\cal B}({\cal H})$ of all bounded 
operators on a Hilbert space ${\cal H}$;

\smallskip
(iii)  {\bf II}$_1$  if the dimension function takes values in the interval $[0,1]$;

\smallskip
(iv)  {\bf II}$_{\infty}$  if the dimension function takes values in the interval $[0,\infty]$;

\smallskip
(v)  {\bf III}  if the dimension function takes values $0$ and $\infty$. 
\end{thm}

\vskip1cm\noindent
{\bf Guide to the literature.}
The $W^*$-algebras were introduced in the seminal paper by [Murray \& von Neumann 1936]   \cite{MuNeu1}.
The textbooks  [Schwartz 1967]   \cite{SCH}   and [Fillmore 1996] \cite{F},  Chapter 4  
cover    the $W^*$-algebras.   Continuous geometries were introduced and studied in the monograph 
[von Neumann 1960]  \cite{NEU}.






\chapter{Connes Geometries}
The main thrust of  Connes'  method    is to exploit   geometric ideas 
and constructions  to study and classify algebras;  similar  approach prevails
 in non-commutative algebraic geometry of  M.~Artin, J.~Tate and M.~van den Bergh, 
see Chapter  13.    The early works of A.~Connes indicated that purely $W^*$-algebraic problems
 can be solved  using  ideas and methods coming  from flow dynamics and geometry of 
 foliations;   subsequently,  this approach  became  a vast program where   one strives to reformulate geometry 
 of the space $X$ in terms of  its commutative $C^*$-algebra $C(X)$ and  to use the acquired 
formal language to  deal  with the non-commutative $C^*$-algebras.  
The program is known  as  noncommutative geometry  in the sense of A.~Connes,  or 
{\it Connes geometry} for brevity.

 \index{Connes geometry}

\section{Classification of type III  factors}
The first successful application of this method was classification 
of the type III factors of  $W^*$-algebras in purely dynamical terms (flow of weights); 
a precursor of the classification was the Tomita-Takesaki Theory of the modular
automorphism group connected to a $W^*$-algebra ${\cal M}$.

 \index{Tomita-Takesaki Theory}

\subsection{Tomita-Takesaki Theory}
 \index{weight}
\begin{dfn}
Let ${\cal M}$  be a $W^*$-algebra.  By a weight $\varphi$ on ${\cal M}$ one understands
a linear map
\displaymath
\varphi: {\cal M}_+\longrightarrow [0,\infty],
\enddisplaymath
where ${\cal M}_+$ is the positive cone of ${\cal M}$.  The trace on ${\cal M}$ is 
a weight  $\varphi$  which is invariant under the action of the unitary group $U({\cal M})$
of ${\cal M}$, i.e.
\displaymath
\varphi(uxu^{-1})=\varphi(x), \quad\forall u\in U({\cal M}). 
\enddisplaymath
 \end{dfn}
\begin{rmk}
\textnormal{
Every $W^*$-algebra on a separable Hilbert space ${\cal H}$ has a weight 
but not a trace;  for instance, type III factors have a weight but no traces. 
}
\end{rmk}
\begin{rmk}
\textnormal{
Recall that the GNS construction is well defined for the traces.  In the case of a
weight $\varphi$ on ${\cal M}$,  the representation of ${\cal M}$ using the GNS 
construction does not yield the unique scalar product on the Hilbert space ${\cal H}_{\varphi}$
because in general $\varphi(x^*x)\ne\varphi(xx^*)$ for $x\in {\cal M}$;  thus one gets 
an unbounded operator 
\displaymath
S_{\varphi}: {\cal H}_{\varphi}\to {\cal H}_{\varphi}
\enddisplaymath
which measures the discrepancy between the scalar products $\varphi(x^*x)$ and $\varphi(xx^*)$.
While dealing with the unbounded operator $S_{\varphi}$ it was suggested by Tomita and Takesaki 
to consider the corresponding polar decomposition
\displaymath
S_{\varphi} =J\Delta^{1\over 2}_{\varphi},
\enddisplaymath
where $J^2=I$ is an isometric involution  on ${\cal H}_{\varphi}$ 
such that $J{\cal M}J={\cal M}'$ replaces 
the $\ast$-involution and $\Delta_{\varphi}$ a positive operator. 
}
\end{rmk}
 \index{Tomita-Takesaki Theorem}
\begin{thm}
{\bf (Tomita and Takesaki)}
For every $t\in {\Bbb R}$ the positive operator $\Delta_{\varphi}$
satisfies an isomorphism
\displaymath
\Delta_{\varphi}^{it}~{\cal M} ~\Delta^{-it}_{\varphi}\cong {\cal M}.
\enddisplaymath
\end{thm}
\begin{rmk}
\textnormal{
The Tomita-Takesaki formula defines a one-parameter group $\sigma_{\varphi}^t$
of automorphisms (a flow or a time evolution) of the $W^*$-algebra ${\cal M}$;
it is called the {\it modular automorphism group} of ${\cal M}$.
}
\end{rmk}
\begin{rmk}
\textnormal{
The lack of trace property $\varphi(xy)=\varphi(yx)$ for general weights $\varphi$
on ${\cal M}$ has  an elegant phrasing as a Kubo-Martin-Schwinger {\it (KMS-) condition}
\displaymath
\varphi(x\sigma^{-i\beta}_{\varphi}(y))=\varphi(yx),
\enddisplaymath
where $t=-i\beta$ and $\beta$ is known in the quantum statistical physics 
as an {\it inverse temperature}. 
}
\end{rmk}
 \index{KMS condition}
 \index{inverse temperature}

 \index{Connes Invariant}
\subsection{Connes Invariants}
\begin{thm}
{\bf (Connes)}
The modular automorphism group $\sigma_{\varphi}^t$ of the $W^*$-algebra
${\cal M}$ is independent of the weight $\varphi$ modulo an inner automorphism
of ${\cal M}$. 
\end{thm}
 \index{flow of weights}
\begin{rmk}
\textnormal{
The canonical modular automorphism group $\sigma^t: {\cal M}\to {\cal M}$
is called a {\it flow of weights} on ${\cal M}$.
}
\end{rmk}
\begin{dfn}
By the Connes invariants $T({\cal M})$ and $S({\cal M})$ of the $W^*$-algabra ${\cal M}$
one understands the following subsets of the real line ${\Bbb R}$:
\displaymath
\left\{
\begin{array}{ccc}
T({\cal M}) &:= & \{t\in {\Bbb R} ~|~ \sigma_t\in Aut~({\cal M}) ~\hbox{\sf is inner}\} \\ 
S({\cal M}) &:= & \{\cap_{\varphi} ~\hbox{{\bf Spec}} 
~(\Delta_{\varphi})\subset {\Bbb R}  ~|~ \varphi ~\hbox{\sf a weight on}~ {\cal M}\},
 \end{array}
\right.
\enddisplaymath
where {\bf Spec} ~$(\Delta_{\varphi})$ is the spectrum of the self-adjoint linear
operator $\Delta_{\varphi}$.  
 \end{dfn}
\begin{thm}
{\bf (Connes)}
If ${\cal M}$ is a type III factor,  then $T({\cal M})\cong {\Bbb R}$ and 
$S({\cal M})$ is a closed multiplicative semigroup of the interval $[0,\infty]$.  
\end{thm}
\begin{cor}
The type III factors ${\cal M}$ are isomorphic to one of the following
\displaymath
S({\cal M})=
\left\{
\begin{array}{cc}
[0,\infty)  &  ~\hbox{\sf type III}_0,\\ 
\{\lambda^{\Bbb Z} \cup\{0\} ~|~ 0<\lambda< 1\}  & ~\hbox{\sf type III$_{\lambda}$},\\
\{0, 1\}   &  ~\hbox{\sf type III$_1$}.
 \end{array}
\right.
\enddisplaymath
\end{cor}
 \index{type {\bf III}$_0$,  {\bf III}$_{\lambda}$, {\bf III}$_1$-factors}
\begin{rmk}
{\bf (Connes and Takesaki)}
\textnormal{
The type III factors have a crossed product structure as follows
\displaymath
\left\{
\begin{array}{ccc}
\hbox{{\bf III}}_0  &\cong&  \hbox{{\bf II}}_{\infty}\rtimes {\Bbb Z},\\ 
\hbox{{\bf III}}_{\lambda} &\cong  & \hbox{{\bf II}}_{\infty}\rtimes_{\lambda} {\Bbb Z},\\
\hbox{{\bf III}}_1 &\cong &  \hbox{{\bf II}}_{\infty}\rtimes {\Bbb R},
 \end{array}
\right.
\enddisplaymath
where the automorphism of  type {\bf II}$_{\infty}$ factor is given by the scaling of trace by $0<\lambda<1$.  
}
\end{rmk}

\section{Noncommutative differential geometry}
Let ${\cal A}\cong C(X)$ be a commutative $C^*$-algebra;  the Gelfand-Naimark Theorem 
says that such algebras are one-to-one with the Hausdorff topological spaces $X$.  
If $X$ is an algebraic variety,  then its coordinate ring $A(X)$ consists of smooth
(even meromorphic!) complex-valued functions on $X$;  the Nullstellensatz says that rings  $A(X)$ are one-to-one with the algebraic varieties $X$.   
Therefore,  to study the  Hausdorff spaces with an extra structure (like the differentiable structure), 
 it makes sense to look at a  dense sub-algebra of ${\cal A}$ consisting of all smooth complex-valued functions
 on $X$;  an  analog of the de Rham theory for such sub-algebras is known as  a {\it cyclic homology}.    
The cyclic homology is a modification of a {\it Hochschild homology} known in algebraic geometry; 
unlike the Hochschild homology, the cyclic homology  is  linked to the $K$-theory by a homomorphism
called   the {\it Chern-Connes character}.  

 \index{Hochschild homology}
 \index{cyclic homology}
 \index{Chern-Connes character}

\subsection{Hochschild homology}
\begin{dfn}
For a commutative algebra $A$ by the Hochschild homology $HH_*(A)$  of the
complex $C_*(A)$:
\displaymath
\dots
\buildrel\rm b_{n+1}\over\longrightarrow
A^{\otimes n +1}
\buildrel\rm  b_n \over\longrightarrow
A^{\otimes n} 
\buildrel\rm  b_{n-1} \over\longrightarrow 
 \dots
\buildrel\rm b_1\over\longrightarrow 
 A,
\enddisplaymath
where 
$$A^{\otimes n}:=\underbrace{A\otimes\dots\otimes A}_{n ~times},$$
and the boundary map $b$ is defined by the formula
\displaymath
\begin{array}{lll}
b_n(a_0\otimes a_1\otimes\dots\otimes a_n) &=&  
\sum_{i=0}^{n-1} (-1)^i a_0\otimes\dots\otimes a_ia_{i+1}\otimes\dots\otimes a_n+\\ 
&&\\
&+& 
(-1)^n a_na_0\otimes a_1\otimes\dots\otimes a_{n-1}.
 \end{array}
\enddisplaymath
 \end{dfn}
 \index{Hochschild-Kostant-Rosenberg Theorem}
\begin{exm}\label{exm9.2.1}
{\bf (Hochschild-Kostant-Rosenberg)}
\textnormal{
Let $A$ be the algebra of regular functions on a smooth affine variety $X$;   denote
by $\Omega^n(A)$ the module of  differential  $n$-forms over an  algebraic de Rham complex of $A$.
For each $n\ge 0$ there exists an isomorphism:
\displaymath
\Omega^n(A)\cong HH_n(A).
\enddisplaymath
}
\end{exm}
\begin{rmk}
\textnormal{
It follows from the above example, that for a commutative algebra the Hochschild homology 
is an analog of the space of differential forms over the algebra;  note that the groups $HH_n(A)$
are well defined when $A$ is a non-commutative associative algebra.  
}
\end{rmk}
\begin{rmk}
\textnormal{
Let $K_*(A)$ be the $K$-theory of  algebra $A$.   Although there exist natural maps $$K_i(A)\to HH_i(A),$$
they are rarely  homomorphisms;   one needs to modify the Hochschild homology to get the desired homomorphisms
(an analog of the Chern character).  Thus one arrives at the notion of cyclic homology, see next section.   
}
\end{rmk}

 \index{cyclic homology}

\subsection{Cyclic homology}
Roughly speaking,  the cyclic homology comes from a sub-complex of the Hochschild complex
$C_*(A)$ of algebra $A$;  the sub-complex is known as a {\it cyclic complex}. 
\begin{dfn}
An $n$-chain of the Hochschild complex $C_*(A)$ of algebra $A$ is called cyclic if 
\displaymath
a_0\otimes a_1\otimes\dots\otimes a_n=(-1)^n a_n\otimes a_0\otimes a_1\otimes\dots\otimes a_{n-1}
\enddisplaymath
for all $a_0,\dots, a_n$ in $A$;  in other words, the cyclic $n$-chains are invariant of a map
$\lambda: C_*(A)\to C_*(A)$ such that $\lambda^n=Id$.  We shall denote the space of all 
cyclic $n$-chains by $C_n^{\lambda}(A)$ .  The sub-complex  of the Hochschild complex $C_*(A)$
of the form
\displaymath
\dots
\buildrel\rm b_{n+1}\over\longrightarrow
C_{n +1}^{\lambda}(A)
\buildrel\rm  b_n \over\longrightarrow
C_{n}^{\lambda}(A) 
\buildrel\rm  b_{n-1} \over\longrightarrow 
 \dots
\buildrel\rm b_1\over\longrightarrow 
 C_1^{\lambda}(A)
\enddisplaymath
is called the cyclic complex of algebra $A$.   The homology of the cyclic complex 
is denoted $HC_*(A)$ and is called a cyclic homology.   
 \end{dfn}
\begin{rmk}
{\bf (Connes)}
\textnormal{
The inclusion map of the cyclic complex $C^{\lambda}_*$ into the Hochschild complex $C_*(A)$
\displaymath
\iota:  C^{\lambda}_*(A)\longrightarrow C_*(A)
\enddisplaymath
induces a  map 
\displaymath
\iota_*:  HH_n(A)\longrightarrow HC_n(A).
\enddisplaymath
The map is an isomorphism for $n=0$ and a surjection for $n=1$. 
}
\end{rmk}
 \index{cyclic complex}
\begin{exm}
\textnormal{
Let ${\cal A}_{\theta}^0\subset {\cal A}_{\theta}$ be the dense sub-algebra of
smooth functions of noncommutative torus ${\cal A}_{\theta}$,  where 
$\theta\in {\Bbb R}\backslash {\Bbb Q}$;   let $HC^*({\cal A}_{\theta}^0)$ denote
the cyclic cohomology of  ${\cal A}_{\theta}^0$,  i.e. a dual of the cyclic homology
$HC_*({\cal A}_{\theta}^0)$.  Then
\displaymath
HC^0({\cal A}_{\theta}^0)\cong {\Bbb C}
\enddisplaymath
and the map
\displaymath
\iota^*:  HC^1({\cal A}_{\theta}^0)\longrightarrow HH^1({\cal A}_{\theta}^0)
\enddisplaymath
is an isomorphism.  
}
\end{exm}
 \index{Chern-Connes character}
\begin{thm}
{\bf (Chern-Connes character)}
For an associative algebra $A$ and each integer $n\ge 0$ there exist natural maps
\displaymath
\left\{
\begin{array}{ccc}
\chi^{2n} & : &  K_0(A)\longrightarrow HC_{2n}(A)\\ 
\chi^{2n+1} & : &  K_1(A)\longrightarrow HC_{2n+1}(A),
 \end{array}
\right.
\enddisplaymath
Moreover,  in case $n=0$  the map $\chi^0$ is a homomorphism.  
\end{thm}

 \index{Novikov Conjecture for hyperbolic\linebreak 
 groups}

\subsection{Novikov Conjecture for hyperbolic groups}
Similar to Kasparov's $KK$-theory for $C^*$-algebras,  the cyclic homology 
can  be used to prove certain cases of the Novikov  Conjecture on the higher signatures
of smooth manifolds, see Ch.10;   thus the methods  of noncommutative differential geometry solve  open 
problems in topology unsolved by other methods.   
\begin{dfn}
Let $G=\langle g_1,\dots, g_n ~|~ r_1,\dots, r_s\rangle$  be a group on $n$ generators 
and $s$ relations.  By the length $l(w)$ of a word $w\in G$ one understands the minimal
number of $g_i$ (counted with powers) which is necessary to write $w$;  the function
\displaymath
l: G\to {\Bbb R}
\enddisplaymath
 turns $G$ into a metric space $(G, d)$,  where $d$ is the distance function. 
 The group $G$ is called hyperbolic whenever $(G,d)$ is a hyperbolic metric 
 space, i.e. for some $w_0\in G$ there exists $\delta_0>0$ such that 
\displaymath
\delta(w_1,w_3)\ge \inf ~(\delta(w_1,w_2), \delta(w_2,w_3))-\delta_0,\quad \forall w_1,w_2,w_3\in G,
\enddisplaymath
 where 
\displaymath
\delta(w,w'):={d(w,w_0)+d(w',w_0)-d(w,w')\over 2}.
\enddisplaymath
 \end{dfn}
 \index{Connes-Moscovici Theorem}
\begin{thm}
{\bf (Connes and Moscovici)}
If $G$ is a hyperbolic group, then $G$ satisfies the Strong Novikov Conjecture. 
 \end{thm}
{\it Proof.}  The proof uses homotopy invariance of the cyclic homology.
$\square$

\section{Connes'  Index  Theorem}
Index Theory says that index of the Fredholm operator on a closed smooth
manifold $M$ is an integer number,  see Ch.10.    A non-commutative analog
of the Atiyah-Singer Index Theorem  (due to A.~Connes)  says that  the index 
can be any real number;   the idea is to consider Fredholm operators not
on $M$ but on a ``foliated $M$''  (i.e. on the leaves of a foliation ${\cal F}$ 
of $M$) and the values of index equal to the trace of elements of the group
$K_0(C_{red}^*({\cal F}))$,   where $C^*_{red}({\cal F})$ is a (reduced)  $C^*$-algebra
of foliation  ${\cal F}$.  A precursor of  Connes' Index Theorem was the 
Atiyah-Singer Index Theorem for families of elliptic operators on closed 
manifolds.

\subsection{Atiyah-Singer Theorem for families of elliptic operators}
We refer the reader to Ch.10 for a brief introduction to the Index Theory.
Unlike the classical Atiyah-Singer Theorem,  the analytic index of a family
of elliptic operators on a compact manifold $M$ no longer belongs to ${\Bbb Z}$
but to the $K$-homology group of the corresponding parameter space.
(It was the way M.~Atiyah constructed his realization of the $K$-homology,
which is a theory dual to the topological $K$-theory.)  Below we briefly review 
the construction.

Let $M$ be an oriented compact smooth manifold and let 
\displaymath
D: \Gamma^{\infty}(E)\to \Gamma^{\infty}(F)
\enddisplaymath
be an elliptic operator defined on the cross-sections of the vector bundles $E, F$
on $M$.  Let $B$ be a locally compact space called a {\it parameter space} 
and consider a continuous family of elliptic operators on $M$
\displaymath
{\cal D}:=\{D_b~|~b\in B\}
\enddisplaymath
parametrized by $B$.  It is not hard to see that the maps 
\displaymath
\left\{
\begin{array}{ccc}
b & \longrightarrow &   Ker ~D_b\\ 
b & \longrightarrow &  Ker ~D_b^*
 \end{array}
\right.
\enddisplaymath
 define two vector bundles $E,F$ over $B$ with a compact support. 
\begin{dfn}
By an analytic index {\bf Ind} ~$({\cal D})$ of the family ${\cal D}$ one understands 
the difference of equivalence classes $[E], [F]$ of the vector bundles $E, F$, i.e.
\displaymath
\hbox{{\bf Ind}} ~({\cal D}):=[E]-[F]\in K^0(B),
\enddisplaymath
where $K^0(B)$ is the $K$-homology group of the topological space $B$.  
 \end{dfn}
 \index{Atiyah-Singer Theorem for family of elliptic operators}
\begin{thm}
{\bf (Atiyah and Singer)}
For any continuous family ${\cal D}$ of elliptic  operators  on a compact
oriented differentiable manifold $M$ the analytic index {\bf Ind} ~$({\cal D}))$
is given by the formula
\displaymath
\hbox{{\bf Ind}} ~({\cal D})=
\langle \hbox{{\bf ch}}~({\cal D}) ~\hbox{{\bf Td}}~(M)\rangle  [M]\in K_0(B),
\enddisplaymath
where  {\bf ch} ~$({\cal D})$ is the Chern class of family ${\cal D}$, 
{\bf Td}~$(M)$ is the Todd genus of $M$,  $[M]$ the fundamental homology class of $M$
and $K_0(B)$ the $K_0$-group of the parameter space $B$. 
  \end{thm}

 \index{foliated space}

\subsection{Foliated spaces}
Roughly speaking,  foliated space (or a foliation) is a decomposition  of
the space $X$ into a disjoint union of sub-spaces of $X$ called {\it leaves}
such that the neighborhood of each point $x\in X$ is a trivial fibration;
the global behavior of leaves is unknown and can be rather complicated.  
Each (global) fiber bundle $(E, p, B)$ is a foliation of the total space $E$ yet
the converse is false;  thus foliations is a generalization of the notion of
the locally trivial fiber bundles.    
\begin{dfn}
By a $p$-dimensional class $C^r$ foliation ${\cal F}$  of an
\linebreak
 $m$-dimensional manifold
$M$ one understands a decomposition of $M$ into a union of disjoint connected subsets
$\{{\cal L}_{\alpha}\}_{\alpha\in A}$ called leaves of the foliation,  so that every point
of $M$ has neighborhood $U$ and a system of local class $C^r$ coordinates 
$x=(x_1,\dots, x_m):U\to {\Bbb R}^m$  such that for each leaf ${\cal L}_{\alpha}$
the components $U\cap {\cal L}_{\alpha}$ are described by the equations
\displaymath
\left\{
\begin{array}{ccc}
x^{p+1}& =&   Const\\ 
x^{p+2}& = &  Const\\
 \vdots &&\\
x^m &=& Const. 
 \end{array}
\right.
\enddisplaymath
The integer number $m-p$ is called a co-dimension of foliation ${\cal F}$.    
If $x, y$ are two points on the leaf ${\cal L}_{\alpha}$ of a foliation ${\cal F}$,
then one can consider  $(m-p)$-dimensional planes $T_x$ and $T_y$ through $x$ and $y$ transversal 
to ${\cal F}$;  a map $H:  T_x\to T_y$ defined by the shift of points along the nearby leaves of ${\cal F}$
is called a holonomy of the foliation ${\cal F}$.  Whenever the holonomy map $H: T_x\to T_y$ 
preserves a measure on $T_x$ and $T_y$ for all $x, y\in M$,  foliation will be called a foliation with the
measure-preserving holonomy, or a measured foliation for short.  
 \end{dfn}
\begin{exm}
\textnormal{
Let $M$ and $M'$ be manifolds of dimension $m$ and $m'\le m$ respectively;
let $f: M\to M'$ be a submersion of $rank~(df)=m'$.  It follows from the Implicit Function
Theorem that $f$ induces a co-dimension $m'$ foliation on $M$ whose leaves are
defined to be the components of $f^{-1}(x)$,  where $x\in M'$.    
}
\end{exm}
\begin{exm}
\textnormal{
The dimension one foliation ${\cal F}$ is given by the orbits of a non-singular
flow $\varphi^t: M\to M$ on the manifold $M$;  the holonomy $H$ of ${\cal F}$ 
coincides with the Poincar\'e ``map of the first return'' along  the flow. Whenever 
 $\varphi^t$  admits an invariant  measure, foliation ${\cal F}$ will be 
a measured foliation.  For instance,  for the flow ${dx\over dt}=\theta=Const$ on the 
torus $T^2$ the corresponding foliation ${\cal F}_{\theta}$ is  a measured foliation
for any $\theta\in {\Bbb R}$.     
}
\end{exm}

 \index{Connes'  Index Theorem}
\subsection{Index Theorem for foliated spaces}
Notice that the Atiyah-Singer Theorem for the families of elliptic operators
can be rephrased as an index theorem for a foliation ${\cal F}$ which is defined
by the fiber bundle $(E, p, B)$ whose base space coincides with the parameter
space $B$;   in this case each elliptic operator $D_b\in {\cal D}$ acts on the leaf $M$,
where $M$ is the manifold associated to the family ${\cal D}$.  What happens if ${\cal F}$
is no longer a fiber bundle but a general foliation?  An analog of the base space $B$ in 
this case will be the leaf space of foliation ${\cal F}$,  i.e. $E/{\cal F}$;  such a space is far
from being a Hausdorff topological space.  The innovative idea of A.~Connes  is to replace
$B$ by a $C^*$-algebra $C^*_{red}({\cal F})$ attached naturally to a foliation ${\cal F}_{\mu}$
endowed with measure $\mu$;  the analytic index of the family ${\cal D}$ is given by the canonical trace
$\tau$ (defined by measure $\mu$)  of  an element of group $K_0(C^*_{red}({\cal F}))$
described by family ${\cal D}$.  The topological index is obtained from the Ruelle-Sullivan current $C_{\mu}$
attached to the measure $\mu$.  The so-defined index is no longer an integer number or an element of the group
$K^0(B)$,  but an arbitrary real number.  
\begin{dfn}
By a $C^*$-algebra $C^*_{red}({\cal F})$ of measured foliation ${\cal F}$ on a manifold $M$
one understands the reduced groupoid $C^*$-algebra $C_{red}^*(H_{\mu})$, where $H_{\mu}$ 
is the holonomy  groupoid of foliation ${\cal F}$.   
 \end{dfn}
\begin{exm}
\textnormal{
Let ${\cal F}_{\theta}$ be a foliation given by the flow ${dx\over dt}=\theta=Const$ on the torus $T^2$. 
Then $C^*_{red}({\cal F})\cong {\cal A}_{\theta}$,   where ${\cal A}_{\theta}$ is the noncommutative torus.}
\end{exm}
\begin{rmk}
{\bf (Connes)}
\textnormal{
If ${\cal D}$ is a collection of elliptic operators defined on each leaf ${\cal F}_{\alpha}$ of foliation ${\cal F}$,
then  ${\cal D}$ gives rise to an element $p_{\cal D}\in K_0(C^*_{red}({\cal F}))$;   moreover,  if $\tau$
is the canonical trace on $C^*_{red}({\cal F})$ coming from the invariant measure $\mu$ on ${\cal F}$,
then the analytic index of ${\cal D}$ is defined by the formula
\displaymath
\hbox{{\bf Ind}} ~({\cal D}):= \tau(p_{\cal D})\in {\Bbb R}. 
\enddisplaymath
 }
\end{rmk}
\begin{thm}
{\bf (Connes)}
Let $M$ be a compact smooth manifold and let ${\cal F}$ be an oriented foliation
of $M$ endowed with invariant transverse measure $\mu$;  suppose that 
$C_{\mu}\in H_*(M; {\Bbb R})$ is the Ruelle-Sullivan current associated to ${\cal F}$.
Then 
\displaymath
\hbox{{\bf Ind}} ~({\cal D})=
\langle \hbox{{\bf ch}}~({\cal D}) ~\hbox{{\bf Td}}~(M), [C_{\mu}] \rangle \in {\Bbb R},
\enddisplaymath
where  {\bf ch} ~$({\cal D})$ is the Chern class of family ${\cal D}$ and 
{\bf Td}~$(M)$ is the Todd genus of $M$. 
  \end{thm}

 \index{Frobenius endomorphism}
 \index{Riemann Hypothesis}
 \index{Bost-Connes system}
\section{Bost-Connes dynamical system}
Weil's Conjectures and Delignes' proof of the Riemann Hypothesis for
the  zeta function  of  projective varieties over  finite fields  rely heavily 
on the spectral data of a linear operator known as the {\it Frobenius endomorphism}.
J.~P.~Bost and A.~Connes considered the following related problem:  
To characterize an operator $T: {\cal H}\to {\cal H}$ on a Hilbert space ${\cal H}$
such that      
\displaymath
\hbox{{\bf Spec}} ~(T)=\{2, 3, 5, 7,\dots\}:={\cal P},
\enddisplaymath
where ${\cal P}$ is the set of prime numbers.  The $C^*$-algebra $C^*({\Bbb N}^*)$
generated by $T$ has  a spate of remarkable properties,  e.g
the ``flow of weights'' $(C^*({\Bbb N}^*),\sigma_t)$ on  $C^*({\Bbb N}^*)$
gives  the partition function $Tr~(e^{-sH})=\zeta(s)$,  where $H$ is a self-adjoint
element of   $C^*({\Bbb N}^*)$ and $\zeta(s)$ is the Riemann zeta function.  
Below we give a brief review of a more general $C^*$-algebra with the same properties, 
which is known as the Hecke $C^*$-algebra.

 \index{Hecke $C^*$-algebra}
\subsection{Hecke $C^*$-algebra}
\begin{dfn}
Let $\Gamma$ be a discrete group and $\Gamma_0\subset\Gamma$ a subgroup 
such that $\Gamma/\Gamma_0$ is a finite set.  By the Hecke $C^*$-algebra 
$C^*_{red}(\Gamma,\Gamma_0)$ one understands the norm-closure of a convolution
algebra on the Hilbert space $\ell^2(\Gamma_0\backslash\Gamma)$,  where the
convolution of $f,f'\in C^*_{red}(\Gamma,\Gamma_0)$ is given by the formula
\displaymath
(f\ast f')(\gamma)=\sum_{\gamma_1\in\Gamma_0\backslash\Gamma} f(\gamma\gamma_1^{-1})f'(\gamma_1),
\quad\forall\gamma\in\Gamma.
\enddisplaymath
\end{dfn}
\begin{rmk}
\textnormal{
It is not hard to see  that $C^*_{red}(\Gamma,\Gamma_0)$ is a reduced group $C^*$-algebra;  hence the notation.  
 }
\end{rmk}
 \index{Hecke operator}
\begin{exm}
\textnormal{
Let $\Gamma\cong GL(2, {\Bbb Q})$ and $\Gamma_0\cong GL(2, {\Bbb Z})$,
Then $C^*_{red}(\Gamma,\Gamma_0)$ contains a (commutative) sub-algebra
of the classical Hecke operators acting on the space of automorphic cusp 
forms.
}  
\end{exm}
\begin{rmk}
{\bf (Bost and Connes)}
\textnormal{
There exists a unique one-parameter group of automorphisms $\sigma_t\in Aut~
(C^*_{red}(\Gamma,\Gamma_0))$ such that
\displaymath
(\sigma_t(f))(\gamma)=\left({L(\gamma)\over R(\gamma)}\right)^{-it}f(\gamma),
\quad \forall\gamma\in \Gamma_0\backslash\Gamma/\Gamma_0,
\enddisplaymath
where 
\displaymath
\left\{
\begin{array}{ccc}
L(\gamma) & =&  |\Gamma_0\gamma\Gamma_0| ~\hbox{in} ~\Gamma/\Gamma_0\\ 
R(\gamma) & =&  |\Gamma_0\gamma\Gamma_0| ~\hbox{in} ~\Gamma_0\backslash\Gamma.
 \end{array}
\right.
\enddisplaymath
 }
\end{rmk}

\subsection{Bost-Connes Theorem}
\begin{dfn}
Consider the matrix group over a ring $R$
\displaymath
P_R:=\left\{\left(\matrix{1 & b\cr  0 & a}\right) ~|~ aa^{-1}=a^{-1}a=1, ~a\in R\right\}
\enddisplaymath
and let $P^+_R$ denote a restriction of the group to the case $a>0$.  
\end{dfn}
 \index{Bost-Connes Theorem}
 \index{KMS state}
\begin{thm}
{\bf (Bost and Connes)}
Let $(A,\sigma_t)$ be the $C^*$-dynamical system associated to the Hecke $C^*$-algebra
\displaymath
A:= C^*_{red}(P_{\Bbb Q}^+, P_{\Bbb Z}^+). 
\enddisplaymath
Then:

\medskip
(i)  for $0<\beta\le 1$ (with $\beta=it$) there exists a unique KMS state $\varphi_{\beta}$ on 
 $(A,\sigma_t)$ and each $\varphi_{\beta}$ is a factor state so that the associated factor is the
 hyper-finite factor of type {\bf III}$_1$;
 
 \smallskip
 (ii) for $\beta>1$ the KMS states on   $(A,\sigma_t)$ form a simplex whose extreme points
 $\varphi_{\beta}^{\chi}$ are parameterized by the complex embeddings
\displaymath
\chi:  {\Bbb Q}_{ab}\to {\Bbb C},
\enddisplaymath
 where ${\Bbb Q}_{ab}$ the  abelian extension of the field ${\Bbb Q}$ by the roots of unity, 
 and such  states correspond to the type {\bf I}$_{\infty}$ factors;
 
 \smallskip
 (iii) the partition function $Tr~(e^{-\beta H})$ of the $C^*$-dynamical system   $(A,\sigma_t)$
 is given by the formula
\displaymath
Tr~(e^{-s H})=\zeta(s),
\enddisplaymath
 where 
\displaymath
\zeta(s)=\sum_{n=1}^{\infty} {1\over n^s}
\enddisplaymath
 is the Riemann zeta function. 
 \end{thm}
 \index{Riemann zeta function}

\vskip1cm\noindent
{\bf Guide to the literature.}
The  Tomita-Takesaki Modular Theory was developed in [Takesaki 1970]  \cite{T1}.  
The Connes Invariants  $T({\cal M})$ and $S({\cal M})$ were introduced and studied 
in  [Connes 1973]  \cite{Con1};  see also an excellent survey [Connes 1978]  \cite{Con2}. 
The Hochschild homology was introduced in [Hochschild 1945]  \cite{Hoc1};  
the Hochschild-Kostant-Rosenberg Theorem (Example \ref{exm9.2.1}) was proved by
[Hochschild,  Kostant \& Rosenberg 1962]  \cite{HoKoRo1}. 
The cyclic homology was introduced by [Connes 1985]  \cite{Con3}  and,  independently,
by [Tsygan 1983]  \cite{Tsy1}  in the context of matrix  Lie algebras.  
The Novikov Conjecture for hyperbolic groups was proved by [Connes \& Moscovici 1990] 
\cite{CoMo1}.  The Index Theorem for families of elliptic operators was proved by
[Atiyah \& Singer 1968]  \cite{AtSi1}.  The Index Theorem for foliation is due to  [Connes  1981] 
\cite{Con4};  see also the monograph  [Moore \& Schochet  2006]  \cite{MS} for a detailed account.    
The Bost-Connes dynamical system was introduced and studied by  [Bost \& Connes 1995] \cite{BoCo1}.






\chapter{Index Theory}
The Index Theory can be viewed as is a covariant version of the topological $K$-theory;
it grew from the Atiyah-Singer Theorem which says the index
of a Fredholm operator on a manifold $M$  can be expressed in purely topological
terms.  M.~Atiayh himself and later L.~Brown, R.~Douglas and P.~Fillmore
elaborated an analytic (i.e. the $C^*$-algebras) realization of such a theory 
which became known as the  {\it $K$-homology}.  The $KK$-theory of 
G.~Kasparov merges together  $K$-theory with the $K$-homology
into a spectacular  bi-functor $KK({\cal A}, {\cal B})$ on the category of pairs of the $C^*$-algebras
${\cal A}$ and ${\cal B}$.  The so far topological applications of the Index Theory include
the proof of certain cases of the  Novikov Conjecture which would 
remain out of reach otherwise.

 \index{Kasparov's KK-theory}
 \index{K-homology}
 \index{Fredholm operator}

\section{Atiyah-Singer Theorem}
\subsection{Fredholm operators}
\begin{dfn}
An operator $F\in B({\cal H})$ is called Fredholm if $F({\cal H})$ is a closed 
subspace of the Hilbert space ${\cal H}$ and the subspaces $Ker~F$ and $Ker~F^*$
are finite-dimensional;  the integer number
\displaymath
\hbox{{\bf Ind}} ~(F):= \dim ~(Ker ~F)- \dim ~(Ker~F^*)
\enddisplaymath
is called an index of the Fredholm operator $F$.  
\end{dfn}
\begin{exm}
\textnormal{
If $\{e_i\}$ is a basis of the Hilbert space ${\cal H}$,  then the unilateral shift
 $S(e_n)=e_{n+1}$ is a Fredholm operator;  the reader can verify that 
 {\bf Ind}~$(S)=-1$. 
}
\end{exm}
 \index{Atkinson Theorem}
\begin{thm}
{\bf (Atkinson)}
The following conditions on an operator $F\in B({\cal H})$ are equivalent:

\medskip
(i)  $F$ is Fredholm;

\smallskip
(ii) there exists an operator $G\in B({\cal H})$ such that $GF-I$  and $FG-I$ are compact operators;

\smallskip
(iii)   the image of $F$ in the Calkin algebra $B({\cal H})/{\cal K}({\cal H})$ is invertible.   
\end{thm}
\begin{cor}
If $F$ is  the Fredholm operator,  then $F^*$ 
is a Fredholm operator;   moreover, 
\displaymath
\hbox{{\bf Ind}} ~(F^*)= -\hbox{{\bf Ind}} ~(F).  
\enddisplaymath
\end{cor}

\begin{cor}
If $F_1$ and $F_2$ are the Fredholm operators,  then their  product $F_1F_2$
is a Fredholm operator;  moreover, 
\displaymath
\hbox{{\bf Ind}} ~(F_1F_2) =  \hbox{{\bf Ind}} ~(F_1) + \hbox{{\bf Ind}} ~(F_2) .  
\enddisplaymath
\end{cor}
\begin{rmk}
\textnormal{
If ${\cal F}_n$ is the set of all Fredholm operators on a Hilbert space ${\cal H}$ having index $n$,
then the map
\displaymath
\hbox{{\bf Ind}}:  \bigcup_{n=-\infty}^{n=\infty} {\cal F}_n:={\cal F}\longrightarrow {\Bbb Z} 
\enddisplaymath
is locally constant,  i.e. ${\cal F}_n$ is an open subset of $B({\cal H})$.   
Two Fredholm operators are connected by a norm-continuous path of Fredholm operators
in the space $B({\cal H})$ 
if an only if they have the same index;  in particular, the connected components of the space
${\cal F}$ coincide with ${\cal F}_n$ for $n\in {\Bbb Z}$.  
}
\end{rmk}
 \index{essential spectrum}
\begin{dfn}
By an essential spectrum {\bf SpEss} ~$(T)$ of an operator $T\in B({\cal H})$ one 
understands a subset of complex plane given by the formula
\displaymath
\hbox{{\bf SpEss}}~(T):=\{\lambda\in {\Bbb C} ~|~ T-\lambda I\not\in {\cal F}\}.
\enddisplaymath
\end{dfn}
\begin{rmk}
\textnormal{
It is not hard to see that
\displaymath
\hbox{{\bf SpEss}}~(T)\subseteq \hbox{{\bf Sp}}~(T),
\enddisplaymath
where {\bf Sp}~$(T)$ is the usual spectrum of $T$. 
}
\end{rmk}

\subsection{Elliptic operators on manifolds}
Let $M$ be a compact oriented smooth $n$-dimensional manifold. 
Let $E$ be a smooth complex vector bundle over $M$.  Denote by
$\Gamma^{\infty}(E)$ the space of all smooth sections of $E$.  We shall
consider linear differential operators
\displaymath
D: \Gamma^{\infty}(E)\to\Gamma^{\infty}(F)
\enddisplaymath
for a pair $E,F$ of vector bundles over $M$.  Denote by $T^* M$ the cotangent 
vector bundle of $M$, by $S(M)$ the unit vector sub-bundle of $T^*M$ and 
by $\pi: S(M)\to M$ the natural projection.    
\begin{dfn}
By a symbol $\sigma(D)$ of the linear differential operator $D$ one understands 
a  vector bundle homomorphism 
\displaymath
\sigma(D): \pi^*E\to \pi^*F.
\enddisplaymath
The operator $D$ is called elliptic if $\sigma(D)$ is an isomorphism.
\end{dfn}
\begin{rmk}
\textnormal{
If the linear differential operator $D: \Gamma^{\infty}(E)\to\Gamma^{\infty}(F)$
is elliptic,  then $D$ is a Fredholm operator.  We shall denote by {\bf Ind} $(D)$
the corresponding index.
}
\end{rmk}
Denote by $B(M)$ the unit ball sub-bundle of $T^*M$;  clearly, $S(M)\subset B(M)$.
Recall that two vector bundles $E,F$ on a topological space $X$ with an isomorphism 
$\sigma$ on a subspace $X_0\subset X$ define an element
\displaymath
d(E,F,\sigma)\in K(X\backslash X_0),
\enddisplaymath
where  $K(X\backslash X_0)$ is the commutative ring (under the tensor product of vector
bundles) coming from the topological $K$-theory of the space $X$ minus $X_0$. 
Thus if $p: B(M)\to M$ is the natural projection,   then the elliptic operator $D$ defines 
an element 
\displaymath
d(p^*E,p^*F,\sigma(D))\in K(B(M)\backslash S(M)). 
\enddisplaymath
But the Chern character provides us with a ring homomorphism
\displaymath
\hbox{{\bf ch}} : K(B(M)\backslash S(M))\to H^*(B(M)\backslash S(M); {\Bbb Q}),
\enddisplaymath
where $H^*(B(M)\backslash S(M); {\Bbb Q})$ is the rational cohomology ring of the 
topological space $B(M)\backslash S(M)$;  thus one gets an element 
\displaymath
\hbox{{\bf ch}} ~d(p^*E,p^*F,\sigma(D))\in H^*(B(M)\backslash S(M); {\Bbb Q}).
\enddisplaymath
 \index{Chern character of elliptic operator}
\begin{dfn}
By a Chern character {\bf ch} $(D)$ of the elliptic operator $D$ on manifold $M$ 
one understands the element
\displaymath
\hbox{{\bf ch}} ~(D):=
\phi_*^{-1}~\hbox{{\bf ch}} ~d(p^*E,p^*F,\sigma(D))\in H^*(M; {\Bbb Q}),
\enddisplaymath
where $\phi_*:  H^k(M, {\Bbb Q})\to H^{n+k}(B(M)\backslash S(M); {\Bbb Q})$
is the Thom isomorphism of the cohomology rings.  
\end{dfn}

 \index{Atiyah-Singer Index Theorem}
\subsection{Index Theorem}
\begin{thm}
{\bf (Atiyah and Singer)}
For any elliptic differential operator $D$ on a compact
oriented differentiable manifold $M$ the index {\bf Ind} ~$(D)$
is given by the formula
\displaymath
\hbox{{\bf Ind}} ~(D)=
\langle \hbox{{\bf ch}}~(D) ~\hbox{{\bf Td}}~(M)\rangle  [M],
\enddisplaymath
where  {\bf Td}~$(M)\in H^*(M; {\Bbb Q})$ is the Todd genus of the complexification of the tangent bundle 
over $M$ and $\langle\alpha\rangle[M]$ is the value of the top-dimensional component 
of an element $\alpha\in  H^*(M; {\Bbb Q})$ on the fundamental homology class $[M]$ of $M$. 
  \end{thm}
\begin{exm}
\textnormal{
Let $M\cong S^1$ be the circle and let $E,F$ be the one-dimensional trivial vector bundles
over $S^1$;  in this case $\Gamma^{\infty}(E)\cong\Gamma^{\infty}(F)\cong C^{\infty}(S^1)$. 
It is not hard to see that elliptic operators $D: C^{\infty}(S^1)\to C^{\infty}(S^1)$ correspond 
to the multiplication operators $M_f$ acting by the pointwise multiplication of functions of
$C^{\infty}(S^1)$ by a function $f\in C^{\infty}(S^1)$;  therefore 
\displaymath
\hbox{{\bf Ind}} ~(D)=-w(f),
\enddisplaymath
where $w(f)$ is the ``winding number'' of $f$,  i.e. the degree of map $f$ taken with the
plus or minus sign.  
}
\end{exm}
 \index{Hirzebruch-Riemann-Roch\linebreak
  Formula}
\begin{exm}
{\bf (Hirzebruch-Riemann-Roch Formula)}
\textnormal{
Let $M$ be a complex manifold of dimension $l$ endowed with a holomorphic vector bundle
$W$ of dimension $n$. Then the spaces $\Gamma^{\infty}(E)$ and $\Gamma^{\infty}(F)$
can be interpreted as the differential harmonic forms of type $(0,p)$ on $M$, where $0\le p\le l$. 
 Using the Dolbeault isomorphism, one can show that {\bf Ind} $(D)=\sum_{p=0}^l (-1)^p \dim~H^p(M, W)$,
 where $H^p(M, W)$ denotes the $p$-dimensional cohomology group of $M$ with coefficients
 in the sheaf of germs of holomorphic sections of $W$. Thus for any compact complex manifold $M$
 and any holomorphic bundle $W$ one gets
\displaymath
\sum_{p=0}^l (-1)^p \dim~H^p(M, W)=\langle \hbox{{\bf ch}}~(W) ~\hbox{{\bf Td}}~(M)\rangle  [M].  
\enddisplaymath
}
\end{exm}

\section{K-homology}
Recall that the $K$-theory is a covariant functor on algebras and a contravariant functor on spaces (to be constructed below).
The category theory tells us that there exists an abstract  dual functor,  which is contravariant on algebras and covariant on spaces;
such a functor will be called a {\it $K$-homology}.  Surprisingly (or not)  all  realizations of    the $K$-homology involve   Fredholm
operators and their indexes;  in particular,  the Atiyah-Singer Index Theorem becomes a pairing statement between $K$-homology
and $K$-theory.

 \index{topological K-theory}

\subsection{Topological K-theory}
\begin{dfn}
For  a compact Hausdorff topological space $X$ a vector bundle over $X$
is a topological space $E$ ,  a continuous map $p: E\to X$ and a finite-dimensinal
vector space $E_x=p^{-1}(X)$ ,  such that $E$ is locally trivial;  usually,  the vector space is taken over
the field of complex numbers ${\Bbb C}$.   An isomorphism
of vector bundles $E$ and $F$ is a homeomorphism from $E$ to $F$ which takes 
fiber $E_x$ to fiber $F_x$ for each $x\in X$ and which is linear on the fibers.  
A trivial bundle over $X$ is a bundle of the form $X\times V$, where $V$ is a fixed 
finite-dimensional vector space and $p$ is projection onto the first coordinate. 
The Whitney sum $E\oplus F$ of the vector bundles $E$ and $F$ is a vector bundle
of the form 
\displaymath
E\oplus F=\{(e,f)\in E\times F ~|~ p(e)=q(f)\},
\enddisplaymath
where $p: E\to X$ and $q: F\to X$ are the corresponding projections.  
\end{dfn}
\begin{rmk}
\textnormal{
The Whitney sum makes the set of isomorphism classes of (complex) vector bundles 
over $X$ into an abelian semigroup  $V_{\Bbb C}(X)$ with an identity;
the identity is the isomorphism class of the $0$-dimensional  trivial bundle.
A continuous map $\phi: Y\to X$ induces a homomorphism $\phi^*: V_{\Bbb C}(X)\to
V_{\Bbb C}(Y)$; thus we have a contravariant functor from topological spaces to 
abelian semigroups.     
}
\end{rmk}
\begin{dfn}
By a $K$-group $K(X)$ of the compact Hausdorff topological space $X$ one understands
the Grothendieck group of the abelian semigroup $V_{\Bbb C}(X)$. 
\end{dfn}
\begin{rmk}
\textnormal{
A continuous map $\phi: Y\to X$ induces a homomorphism $\phi^*: K(X)\to
K(Y)$;  thus we have a contravariant functor from topological spaces to the
abelian groups,  see  Fig. 10.1.  
}
\end{rmk}
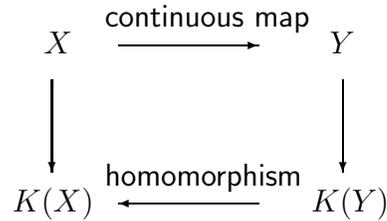
\begin{figure}[here]
\begin{picture}(300,100)(-120,0)
\begin{picture}(300,110)(0,0)
\put(20,70){\vector(0,-1){35}}
\put(130,70){\vector(0,-1){35}}
\put(98,23){\vector(-1,0){53}}
\put(45,83){\vector(1,0){53}}
\put(5,20){$K(X)$}
\put(118,20){$K(Y)$}
\put(17,80){$X$}
\put(125,80){$Y$}
\put(40,30){\sf homomorphism}
\put(40,90){\sf continuous map}
\end{picture}
\end{picture}
\caption{Contravariant functor $K(\bullet)$.}
\end{figure}
\begin{exm}
\textnormal{
If $X\cong S^2$ is a two-dimensional sphere, then $K(S^2)\cong {\Bbb Z}^2$. 
}
\end{exm}
\begin{exm}
\textnormal{
If $X$ is a compact Hausdorff topological space, then
\displaymath
K(X)\cong K_0(C(X)),
\enddisplaymath
where $C(X)$ is the abelian $C^*$-algebra of continuous  complex-valued function on $X$ 
and $K_0(C(X))$ is the respective $K_0$-group,  see Chapter 3.  
}
\end{exm}
\begin{dfn}
If $X$ is a compact Hausdorff space, then by a suspension of $X$ one understands
a topological space  $SX$ given by the formula $SX:= X\times {\Bbb R}$.  The 
higher-order $K$-theory of $X$ can be defined  according to the formulas: 
\displaymath
\left\{
\begin{array}{ccc}
K^0(X) &=& K(X)\\ 
K^{-1}(X)  &=& K(SX)\cong K(X\times {\Bbb R})\\
\vdots       & &\vdots\\
K^{-n}(X) &=& K(S^nX)\cong K(X\times {\Bbb R}^n).        
\end{array}
\right.
\enddisplaymath
\end{dfn}
 \index{Bott periodicity}
\begin{thm}
{\bf (Bott Periodicity)}
There exists an isomorphism (a natural transformation, resp.) between the the groups
(the functors, resp.) $K(X)$ and $K^{-2}(X)$ so that one gets a general isomorphism: 
\displaymath
K^{-n}\cong K^{-n-2}(X).
\enddisplaymath
 \end{thm}
\begin{rmk}
\textnormal{
The topological $K$-theory is an extraordinary cohomology theory,  i.e.
a sequence of the homotopy invariant contravariant functors from compact 
Hausdorff spaces to abelian groups satisfying all the Eilenberg-Steenrod axioms
but  the dimension axiom.    
}
\end{rmk}
 \index{Chern character formula}
\begin{thm}
{\bf (Chern Character)}
If  $X$ is  a compact Hausdorff topological space, then
\displaymath
\left\{
\begin{array}{ccc}
K^0(X)\rtimes {\Bbb Q} &\cong& \bigoplus_k  H^{2k}(X; {\Bbb Q})\\ 
K^{-1}(X)\rtimes {\Bbb Q} &\cong& \bigoplus_k  H^{2k-1}(X; {\Bbb Q}),      
\end{array}
\right.
\enddisplaymath
where $H^k(X; {\Bbb Q})$ denotes the $k$-th cohomology group of $X$
with coefficients in ${\Bbb Q}$. 
\end{thm}

 \index{K-homology}
 \index{Atiyah's realization of K-homology}
\subsection{Atiayh's re\-ali\-zation of K-ho\-mo\-lo\-gy}
An {\it abstract} covariant functor on compact Hausdorff topological spaces $X$
will be denoted by $K_0(X)$.  While looking for a {\it realization} of $K_0(X)$
(i.e. a specific construction yielding the required properties of the functor)
it was observed by M.~Atiyah that the index map 
\displaymath
\hbox{{\bf Ind}} : {\cal F}\to {\Bbb Z}
\enddisplaymath
behaves much like a ``homology'' on the space ${\cal F}$ of all Fredholm 
operators on a Hilbert space ${\cal H}$;  to substantiate such a homology
the following definition was introduced.  
\begin{dfn}
Let $X$ is a compact Hausdorff topological space and let $C(X)$ be the $C^*$-algebra
of continuous complex-valued functions on $X$.  Let
\displaymath
\left\{
\begin{array}{ccc}
\sigma_1 &:& C(X)\to B({\cal H}_1)\\ 
\sigma_2 &:& C(X)\to B({\cal H}_2)
.        
\end{array}
\right.
\enddisplaymath
be a pair of representations of $C(X)$ by bounded linear operators on the Hilbert 
spaces ${\cal H}_1$ and ${\cal H}_2$,  respectively.  Denote by $F: {\cal H}_1\to {\cal H}_2$
a Fredholm operator such that the operator 
\displaymath
F\circ\sigma_1(f)-\sigma_2(f)\circ F 
\enddisplaymath
is a compact operator for all $f\in C(X)$.  By {\bf Ell} $(X)$ one understands the set
of all triples $(\sigma_1,\sigma_2, F)$.  A binary operation of addition on the set
{\bf Ell} $(X)$ is defined as the orthogonal direct sum of the Hilbert spaces ${\cal H}_i$;
the operation turns {\bf Ell} $(X)$ into an abelian semigroup and the Grothendieck 
completion of the semigroup turns  {\bf Ell} $(X)$ into  an abelian group.  
  \end{dfn}
\begin{exm}
\textnormal{
The basic  example is $F$ being an elliptic differential operator
between two smooth vector bundles $E$ and $F$ on a smooth manifold $M$;
it can be proved that $F$ is a Fredholm operator which intertwines the action of
$C(M)$ modulo a compact operator.   
}
\end{exm}
\begin{thm}
{\bf (Atiyah)}
There exists a surjective map
\displaymath
\hbox{{\bf Ell}} ~(X)\longrightarrow K_0(X)
\enddisplaymath
whenever $X$ is a finite $CW$-complex. 
 \end{thm}
\begin{rmk}
\textnormal{
To finalize Atiyah's  realization of the $K$-homology by the group {\bf Ell} $(X)$,
one needs a proper notion of the equivalence relation on {\bf Ell} $(X)$
to make the quotient equal to $K_0(X)$;  it was an open problem solved 
by   L.~Brown, R.~Douglas and P.~Fillmore, see next Section. 
}
\end{rmk}

 \index{Brown-Douglas-Fillmore Theory}

\subsection{Brown-Douglas-Fillmore Theory}
Larry Brown, Ron Douglas and Peter Fillmore  introduced a realization of
the $K$-homology inadvertently (or not)  while working on the problem of
classification of the essentially normal operators on a Hilbert space ${\cal H}$
raised by J. von ~Neumann; solution of  the problem depends on the classification of 
the extensions of the $C^*$-algebra ${\cal K}$ of all compact operators by 
a $C^*$-algebra $C(X)$ for some compact metrizable space $X$.  
\begin{dfn}
By an extension of the $C^*$-algebra ${\cal A}$ by a $C^*$-algebra ${\cal B}$
one understands a $C^*$-algebra ${\cal E}$ which fits into the short exact 
sequence
\displaymath
0\to {\cal A}\to {\cal E}\to {\cal B}\to 0,
\enddisplaymath
i.e. such that ${\cal A}\cong {\cal E}/{\cal B}$.  Two extensions $\tau$ and $\tau'$
are said to be stably equivalent if they differ by a trivial extension. 
If $\tau_1$ and $\tau_2$ are two extensions of ${\cal A}$ by ${\cal B}$, then
there is a naturally defined extension $\tau_1\oplus\tau_2$ called the sum 
of $\tau_1$ and $\tau_2$.   With the addition operation the set of all extensions of 
${\cal A}$ by ${\cal B}$ is an abelian semigroup;    the Grothendieck 
completion of the semigroup yields   an abelian group denoted by {\bf Ext}~ $({\cal A}, {\cal B})$.  
 \end{dfn}
 \index{Brown-Douglas-Fillmore Theorems}
\begin{thm}
{\bf (Brown, Douglas and Fillmore)}
If $X$ is a compact metrizable topological space, then the map
\displaymath
X\mapsto \hbox{{\bf Ext}} ~({\cal K}, C(X)) 
\enddisplaymath
is a homotopy invariant covariant functor from the category of compact metrizable
spaces to the category of abelian groups;  the functor is denoted by $K_1(X)$.
 \end{thm}
\begin{cor}
The suspension $SX$ of $X$  gives a realization of the $K$-homology according to
 the formula
\displaymath
K_0(X)\cong K_1(SX).
\enddisplaymath
 \end{cor}
\begin{rmk}
\textnormal{
Existence of the functor $K_1(X)$ yields a solution to the following problem of J. ~von~Neumann:
Given two operators $T_1,T_2\in B({\cal H})$,  when is it true 
\displaymath
T_1=UT_2U^*+K,
\enddisplaymath
where $U$ is a unitary and $K$ is a compact operator?  Below we give a brief account
of the solution, which was  the first successful application of the $K$-theory to an open
problem in analysis. 
}
\end{rmk}
\begin{dfn}
An operator $T\in B({\cal H})$ is called essentially normal whenever
\displaymath
TT^*-T^*T\in {\cal K}.
\enddisplaymath
\end{dfn}
Let $T\in B({\cal H})$ be an operator and denote by ${\cal E}_T:=C^*(I,T, {\cal K})$ a 
$C^*$-algebra  generated by the identity operator $I$, by  operator $T$ and by the $C^*$-algebra
${\cal K}$.  It follows from the definition that ${\cal E}_T/{\cal K}$ is a commutative $C^*$-algebra
if and only if $T$ is an essentially normal operator.  Thus one gets a short exact sequence
\displaymath
0\to {\cal K}\to {\cal E}_T\to C(\hbox{{\bf SpEss}}~(T)) \to 0,
\enddisplaymath
where $T$ is an essentially normal operator and {\bf SpEss} ~$(T)$ its essential 
spectrum.  Therefore the von ~Neumann  problem for the essentially normal
operators reduces to the classification  (modulo stable equivalence)  of the extensions of 
${\cal K}$ by $C(X)$,  where $X$ is a subset of the complex plane.  A
precise statement is this.    
\begin{thm}
{\bf (Brown, Douglas and Fillmore)}
Two essentially normal operators $T_1$ and $T_2$ are related 
by the formula $T_1=UT_2U^*+K$ if and only if
\displaymath
\left\{
\begin{array}{ccc}
\hbox{{\bf SpEss}} ~(T_1) &= & \hbox{{\bf SpEss}} ~(T_1) \\ 
\hbox{{\bf Ind}} ~(T_1-\lambda I) &= & \hbox{{\bf Ind}} ~(T_2-\lambda I)
 \end{array}
\right.
\enddisplaymath
for all  $\lambda\in {\Bbb C}\backslash\hbox{{\bf SpEss}} ~(T_i)$.
\end{thm}

\section{Kasparov's  KK-theory}
Kasparov's bi-functor $KK({\cal A},{\cal B})$  on pairs of the $C^*$-algebras
${\cal A}$ and ${\cal B}$ blends $K$-homology with the $K$-theory.   Such a  functor
was designed to solve a concrete open problem of topology -- the higher signatures 
hypothesis of S.~P.~Novikov;  it does achieve the goal, albeit in certain special cases. 
The original definition of the $KK$-groups uses the notion of a Hilbert module.    
 \index{Kasparov's KK-theory}
 \index{Hilbert module}

\subsection{Hilbert modules}
Roughly speaking,  the Hilbert module is a Hilbert space  whose inner product takes 
values not in ${\Bbb C}$ but in an arbitrary  $C^*$-algebra ${\cal A}$;  one can define  
bounded and compact operators acting on such modules and the construction appears
to be  very useful. 
\begin{dfn}
For a $C^*$-algebra ${\cal A}$ by a Hilbert module over ${\cal A}$ one understands 
a right ${\cal A}$-module $E$ endowed with a ${\cal A}$-valued inner product
\displaymath
\langle\bullet, \bullet\rangle: E\times E\to {\cal A}, 
\enddisplaymath
which satisfies the following properties:

\medskip
(i) $\langle\bullet, \bullet\rangle$ is sesquilinear; 

\smallskip
(ii) $\langle x, ya\rangle=\langle x, y\rangle a$ for all $x,y\in E$ and $a\in {\cal A}$;

\smallskip
(iii)  $\langle y, x\rangle=\langle x, y\rangle^*$ for all $x,y\in E$; 

\smallskip
(iv)  $\langle x, x\rangle\ge 0$ and if  $\langle x, x\rangle=0$ then $x=0$.

\medskip\noindent
and the norm $||x||=\sqrt{|| \langle x, x\rangle ||}$  under which $E$ is complete. 
 \end{dfn}
\begin{exm}
\textnormal{
The $C^*$-algebra ${\cal A}$ is itself a Hilbert ${\cal A}$-module with  the
inner product  $\langle a, b\rangle:=a^*b$.
}
\end{exm}
\begin{exm}
\textnormal{
Any closed right ideal of the $C^*$-algebra ${\cal A}$ is a Hilbert ${\cal A}$-module. 
}
\end{exm}
\begin{dfn}
For a Hilbert ${\cal A}$-module $E$ by $B(E)$ one understands the set of all module
homomorphisms 
\displaymath
T: E\longrightarrow E
\enddisplaymath
 for which there is an adjoint module homomorphism $T^*: E\to E$ such that
 $\langle Tx,y\rangle=\langle x,T^*y\rangle$ for all $x,y\in E$.  
 \end{dfn}
\begin{rmk}
\textnormal{
Each homomorphism (``operator'') in $B(E)$ is bounded and $B(E)$ itself is a $C^*$-algebra
with respect to the operator norm.
}
\end{rmk}
\begin{dfn}
By the set ${\cal K}(E)\subset B(E)$ one understands the closure of the linear span 
of all  bounded operators $T(x,y)$ acting by the formula
\displaymath
z\longmapsto x\langle y,z\rangle, \qquad z\in E,
\enddisplaymath
where $x,y\in E$ and $T^*(x,y)=T(y,x)$.    
 \end{dfn}
\begin{rmk}
\textnormal{
The set ${\cal K}(E)$ is a closed ideal in $B(E)$.  
}
\end{rmk}

 \index{KK-groups}

\subsection{KK-groups}
\begin{dfn}
For a pair of the $C^*$-algebras ${\cal A}$ and ${\cal B}$  by a Kasparov module 
{\bf E}$({\cal A}, {\cal B})$ one understands the set of all triples $(E, \phi, F)$, where $E$
is a countably generated (graded) Hilbert module over ${\cal B}$, 
\displaymath
 \phi:  {\cal A}\to B(E)
\enddisplaymath
is an $\ast$-homomorphism and $F\in B(E)$ is a bounded operator  such that 
\displaymath
\left\{
\begin{array}{ccc}
F\circ\phi(a) -\phi(a)\circ F&\in & K(E) \\ 
(F^2-I)\circ \phi(a) &\in & K(E)\\
(F-F^*)\circ \phi(a) &\in & K(E),
 \end{array}
\right.
\enddisplaymath
for all $a\in {\cal A}$. 
\end{dfn}
\begin{exm}
{\bf (Basic example)}
\textnormal{
Let $M$ be a compact oriented smooth manifold endowed with a riemannian metric.
Consider a pair of vector bundles  $V_1$ and $V_2$ over $M$ and let
\displaymath
P:  \Gamma^{\infty} (V_1) \to \Gamma^{\infty} (V_2)
\enddisplaymath
be an elliptic operator on $M$ which extends to a Fredholm operator from $L^2(V_1)$ 
to $L^2(V_2)$. Denote by $Q$ the parametrix for operator $P$.    Let $E= L^2(V_1)\oplus L^2(V_2)$ be the Hilbert module (over ${\Bbb C}$)
and consider a $\ast$-homomorphism    
\displaymath
\phi: C(M)\to B(E)
\enddisplaymath
realized by the multiplication operators on $B(E)$.     The triple
\displaymath
(E, \phi, \left(\small\matrix{0 & Q\cr P & 0}\right)) 
\enddisplaymath
is an element of the Kasparov module {\bf E}$(C(M), {\Bbb C})$.    
}
\end{exm}
 \index{Kasparov module}
\begin{dfn}
Let ${\cal A}$ and ${\cal B}$  be a pair of the $C^*$-algebras and consider the 
 set of equivalence classes of the Kasparov modules {\bf E}$({\cal A}, {\cal B})$ 
 under a homotopy equivalence relation.  The direct sum of the Kasparov modules
 turns  the equivalence classes of  {\bf E}$({\cal A}, {\cal B})$ into an abelian semigroup;
 the Grothendieck completion of the semigroup is an abelian group denoted 
 by $KK({\cal A}, {\cal B})$.   
\end{dfn}
\begin{exm}
\textnormal{
If ${\cal A}\cong {\cal B}\cong {\Bbb C}$, then $KK({\Bbb C}, {\Bbb C})\cong {\Bbb Z}$. 
}
\end{exm}
\begin{rmk}
\textnormal{
The $KK$-groups is a bi-functor defined on the pairs of $C^*$-algebras and with values in the 
abelian groups;   the functor is homotopy-invariant in each variable.  A powerful new feature 
of the $KK$-functor is an {\it intersection product}  acting  by the formula
\displaymath
KK({\cal A}, {\cal B}) \times KK({\cal B}, {\cal C})\to KK({\cal A}, {\cal C}). 
\enddisplaymath
The intersection product generalizes the cap and cup products of algebraic topology
and turns $KK({\cal A}, {\cal B})$ into a ring with the multiplication operation defined 
by the intersection product.  
}
\end{rmk}
\begin{exm}
\textnormal{
If ${\cal A}$ is an AF-algebra,  then 
\displaymath
KK({\cal A}, {\cal A}) \cong_{ring}  End~(K_0({\cal A})),  
\enddisplaymath
where $\cong_{ring}$ is a ring isomorphism and $End~(K_0({\cal A}))$ is the endomorphism ring 
of $K_0$-group of the AF-algebra ${\cal A}$.  
}
\end{exm}

 \index{Novikov Conjecture}

\section{Applications of Index Theory}
The Index Theory is known to have successful applications on topology (higher signatures) 
and geometry (positive scalar curvature);  the Baum-Connes Conjecture can be viewed as a 
generalization of the Index Theory in topology.   
\subsection{Novikov Conjecture}
It is well known that the homology and cohomology groups are invariants of 
homeomorphisms and also  more general continuous maps known as homotopies 
of the topological space;   for instance,   homotopy affects  the dimension of topological
space while homeomorphism always preserves it.  The  numerical invariants obtained 
by  pairing of homology and cohomology groups  (i.e integration of cycles against the co-cylcles)   are 
invariants of  homeomorphisms but  not of homotopies.  It is an interesting and difficult problem 
to find   all numerical invariants preserved by  homotopies of the topological space.         
\begin{exm}
{\bf (Hirzebruch Signature)}
\textnormal{
Let $M^{4k}$ be a smooth
\linebreak
 $4k$-dimensional manifold.  Denote by $H^*(M^{4k}; {\Bbb Q})$ its
cohomology ring with the coefficients in ${\Bbb Q}$ and let $L(M^{4k})\in H^*(M^{4k}; {\Bbb Q})$
be the Hirzebruch characteristic class (an $L$-class) of the manifold $M^{4k}$.  If 
$[M^{4k}]\in H_{4k}(M^{4k}; {\Bbb Z})$ is the fundamental class of $M^{4k}$,  then one 
can consider a (rational) number
\displaymath
\hbox{{\bf Sign}} ~(M^{4k})=\langle L(M^{4k}), [M^{4k}]\rangle\in {\Bbb Q}  
\enddisplaymath
obtained by integration of the Hirzebruch characteristic class against the fundamental
class of manifold $M^{4k}$.  The number coincides with the signature of a non-degenerate
symmetric bilinear form on the space $H^{2k}(M^{4k}; {\Bbb Q})$ and for this reason 
is called a {\it signature}.  The signature {\bf Sign} $(M^{4k})$ is a homotopy invariant 
of the manifold $M^{4k}$.  
}
\end{exm}
 \index{Hirzebruch Signature}
\begin{exm}
\textnormal{
Let $M^n$ be an $n$-dimensional manifold and let $p_k(M^n)\in H^{4k}(M^n; {\Bbb Q})$
be its Pontryagin characteristic classes, i.e. such classes for the tangent bundle over  $M^n$.  
The integrals
\displaymath
\langle p_k(M^n), [c_{n-4k}]\rangle\in {\Bbb Q}  
\enddisplaymath
of the Pontryagin classes against the cycles $[c_{n-4k}]\in H_{n-4k}(M^n; {\Bbb Z})$
are topological (i.e. homeomorphism) invariant but not homotopy invariant.  
}
\end{exm}
 \index{Pontryagin characteristic classes}
\begin{rmk}
\textnormal{
Let $K$ denote  a CW-complex.   It is known that $K$ has a distinguished set 
of cohomology classes determined by its  fundamental group $\pi_1(K)$.  
Namely, such classes come from the canonical continuous  map 
\displaymath
f: K\to K(\pi, 1),
\enddisplaymath
where $K(\pi, 1)$ is the classifying space of the group $\pi_1(K)$;
the distinguished set is given by the induced homomorphism 
\displaymath
f^*: H^*(K(\pi, 1))\to H^*(K)
\enddisplaymath
with any coefficients. In the case of the rational coefficients,  the corresponding class of 
distinguished cycles will be denoted by 
\displaymath
Df^* H^*(\pi; {\Bbb Q})\subset H_*(K; {\Bbb Q}). 
\enddisplaymath
}
\end{rmk}
 \index{Novikov Conjecture}
\begin{cnj}
{\bf (Novikov Conjecture on the Higher Signatures)} 
For each cycle $z\in Df^* H^*(\pi; {\Bbb Q})$ the integral 
\displaymath
\langle L_k(p_1,\dots, p_k), z\rangle\in {\Bbb Q}
\enddisplaymath
of the Pontryagin-Hirzebruch 
characteristic class $L_k(p_1,\dots, p_k)\in H^*(M^n; {\Bbb Q})$ against $z$
is a homotopy invariant.  
\end{cnj}
 \index{Kasparov Theorem}
\begin{thm}
{\bf (Kasparov)} 
Whenever $\pi_1(M^n)$ is isomorphic to a discrete subgroup of a connected Lie group,
the Novikov Conjecture for the manifold $M^n$ is true. 
\end{thm}
{\it Outline of  proof.}  The proof is based on  the $KK$-theory.  The idea is to introduce a pair of
the $C^*$-algebras ${\cal A}\cong C^*(\pi_1(M^n))$ and ${\cal B}\cong C(M^n)$,
where  $C^*(\pi_1(M^n))$ is a group $C^*$-algebra attached to the group $\pi_1(M^n)$. 
(Such an  approach can recapitulated in terms of  representation of $\pi_1(M^n)$  
by linear operators on a Hilbert space ${\cal H}$;  thus one gets a functor from topological
spaces to the $C^*$-algebras.)  The bilinear form $\langle L_k(p_1,\dots, p_k), z\rangle$
can be written in terms of the intersection product on the Kasparov's $KK$-groups.  
Since the $KK$-functor is a homotopy invariant  so will be the bilinear form and its signature,
whenever $\pi_1(M^n)$  satisfies assumptions of the theorem.
$\square$  
\begin{rmk}
\textnormal{
Note that Novikov's Conjecture can be rephrased as a question about general discrete
groups and their representations by the $C^*$-algebras;  positive answer to  such a question
would imply the truth of Novikov's Conjecture.  Namely, one has the following  
}
\end{rmk}
 \index{Strong Novikov Conjecture}
\begin{cnj}
{\bf (Strong Novikov Conjecture)} 
 Let $B\pi$ be the classifying space of a discrete  group $\pi$;  then the map 
\displaymath
\beta:  K_*(B\pi)\to K_*(C^*(\pi))
\enddisplaymath
is rationally injective.
\end{cnj}
\begin{rmk}
\textnormal{
Kasparov's proof of Novikov's Conjecture is valid for the Strong Novikov Conjecture. 
}
\end{rmk}

 \index{Baum-Connes Conjecture}

\subsection{Baum-Connes  Conjecture}
The Baum-Connes Conjecture can be viewed as strengthening
(to an isomorphism) and generalization (to the Lie groups)  of the 
Strong Novikov Conjecture.   
\begin{cnj}
{\bf (Baum-Connes Conjecture)} 
 Let $G$ be a Lie group or a countable discrete group;  let $BG$ be the 
 corresponding classifying space.  Let $C^*_{red}(G)$ be the reduced group
 $C^*$-algebra of $G$. Then there exists an isomorphism 
\displaymath
\mu:  K^0(BG)\to K_0(C^*_{red}(G)),
\enddisplaymath
where $K^0(BG)$ is the K-homology of the topological space $BG$.  
\end{cnj}
\begin{rmk}
\textnormal{
The Baum-Connes Conjecture is proved in many cases, e.g. for the hyperbolic, amenable, etc. groups;
we refer the reader to the respective literature.  It seems that the cyclic cohomology is an appropriate
replacement for the $KK$-theory in this case.  
}
\end{rmk}

 \index{positive scalar curvature}
\subsection{Positive scalar curvature}
It is long known that topology of a manifold imposes constraints on the
type of metric one can realize on the manifold;  for instance,  there are no
zero scalar curvature riemannian metrics on any surface of genus $g>1$.
It is remarkable that the Index Theory detects the topological obstructions
for having a metric of positive scalar curvature on a given manifold;  below
we give a brief account of this  theory.  

Let $M$ be a smooth manifold having an oriented spin structure;  let $D$ be the
canonical Dirac operator corresponding to the spin structure.  Suppose that $M$
has even dimension and let $D=D^+\oplus D^-$ be the corresponding canonical 
decomposition of the Dirac operator on positive and negative components.  
For a unital $C^*$-algebra $B$ one can define a flat $B$-vector bundle $V$ over
$M$ and consider the Dirac operators $D_V^+$ and $D_V^-$;  the Dirac operator
$D_V^+$ is a Fredholm operator and 
\displaymath
\hbox{{\bf Ind}}_B(D_V^+)\in K_0(B)\otimes {\Bbb Q}.  
\enddisplaymath
\begin{thm}
{\bf (Rosenberg)} 
If $M$ admits a metric of positive scalar curvature, then 
\displaymath
\hbox{{\bf Ind}}_B(D_V^+)=0
\enddisplaymath
for any flat $B$-vector bundle $V$.
\end{thm}
 \index{Mischenko-Fomenko Theorem}
\begin{rmk}
{\bf (Mischenko-Fomenko)}
\textnormal{
One can express the left-hand side of the last equality in purely topological terms;
namely,
\displaymath
\hbox{{\bf Ind}}_B(D_V^+)=\langle \hbox{{\bf A}}(M) \cup \hbox{{\bf ch}} (V),  [M]\rangle,
\enddisplaymath
where {\bf A}$(M)\in H^*(M, {\Bbb Q})$,  {\bf ch} $(V)$ is the Chern class of $V$ and $[M]$ is
the fundamental class of $M$.  
}
\end{rmk}
\begin{cor}
If $M$ admits a metric of positive scalar curvature, then 
\displaymath
\langle \hbox{{\bf A}}(M) \cup \hbox{{\bf ch}} (V),  [M]\rangle=0
\enddisplaymath
for any flat $B$-vector bundle $V$.
\end{cor}
 \index{Gromov-Lawson Conjecture}
\begin{cnj}
{\bf (Gromov-Lawson Conjecture)} 
 If $M$ has positive scalar curvature,  then for all $x\in H^*(B\pi; {\Bbb Q})$ 
 it holds   
\displaymath
\langle \hbox{{\bf A}}(M) \cup f^*(x),  [M]\rangle=0
\enddisplaymath
where $f: M\to B\pi$ is the classifying map.  
\end{cnj}
\begin{rmk}
\textnormal{
The Strong Novikov Conjecture implies the Gromov-Lawson Conjecture. 
}
\end{rmk}

 \index{coarse geometry}
\section{Coarse geometry}
Index Theory gives only rough topological  invariants of manifolds; 
the idea of {\it coarse geometry} is to replace usual   homeomorphisms
and homotopies between topological spaces by more general {\it coarse
maps}  so that the Index Theory will classify  manifolds modulo the 
equivalence relation defined by such  maps.   The coarse geometry 
is present in Mostow's proof of the Rigidity Theorem, see 
[Mostow 1973]  \cite{MOS};   it has been studied in geometric group theory.
We refer the reader to    the books  by  [Roe 1996]   \cite{R1} and [Nowak \& Yu  2012]  \cite{NY}
for a detailed account.    
\begin{dfn}
If $X$ and $Y$ are metric spaces then a (generally discontinuous)  map $f: X\to Y$ 
is called coarse if:

\medskip
(i)  for each $R>0$ there exists $S>0$ such that $d(x,x')\le R$ implies $d(f(x), f(x'))\le S$;

\smallskip
(ii) for each bounded subset $B\subseteq Y$ the pre-image $f^{-1}(B)$ is bounded in $X$.
\end{dfn}
 \index{coarse map}
\begin{dfn}
Two coarse maps $f_0: X\to Y$ and $f_1: X\to Y$ are coarsely equivalent if there is a 
constant $K$ such that
\displaymath
d(f_0(x), f_1(x))\le K
\enddisplaymath
for all $x\in X$.  
\end{dfn}
\begin{dfn}
Two metric spaces $X$ and $Y$ are said to be coarsely equivalent
if there are maps from $X$ to $Y$ and from $Y$ to $X$ whose composition 
is coarsely equivalent to the identity map;   the coarse equivalence class 
of metrics on $X$ is called a coarse structure (or coarse geometry)  on $X$.   
\end{dfn}
\begin{rmk}
\textnormal{
It is useful to think of coarse  geometry as a ``blurry version'' of usual geometry;  
all the local data is washed out and only the large scale  features are preserved. 
Because the index of  Fredholm operator  on a Hilbert space ${\cal H}$ does not ``see''
the local geometry,  one can  think of the index as an abstract topological invariant
of coarse equivalence.       
}
\end{rmk}
\begin{rmk}
\textnormal{
Any complete riemannian manifold is a metric space;  thus it can be endowed 
with a coarse structure so that the Index Theory becomes a {\it functor} 
on such a structure. 
}
\end{rmk}
\begin{exm}
\textnormal{
The natural inclusion ${\Bbb Z}^n\to {\Bbb R}^n$ is a coarse equivalence.
}
\end{exm}
\begin{exm}
\textnormal{
The natural inclusion ${\Bbb R}^n\to {\Bbb R}^n\times [0,\infty)$ is a coarse equivalence.
}
\end{exm}
\begin{rmk}
\textnormal{
Since the tools of algebraic topology can be applied to the coarse structures on
manifolds,  one can reformulate all content of Sections 10.1-10.4 in terms of the
coarse geometry;   thus one arrives at the notion of a  coarse index,  coarse
assembly map,  coarse Baum-Connes Conjecture, {\it etc.}  We refer the interested 
reader to the monographs    by  [Roe 1996]   \cite{R1} and [Nowak \& Yu  2012]  \cite{NY}.  
}
\end{rmk}

 \index{coarse Baum-Connes Conjecture}

\vskip1cm\noindent
{\bf Guide to the literature.}
The Atiyah-Singer Theorem was announced in [Atiyah \& Singer 1963]  \cite{AtSi2}
and  proved in the series of papers [Atiyah \& Singer 1968, 1971]  \cite{AtSi1}.  
For the foundation of topological $K$-theory we refer the reader to 
monograph [Atiyah 1967]  \cite{A1}.   Atiyah's realization of $K$-homology
can be found in [Atiyah 1970]  \cite{Ati1}.     
The Brown-Douglas-Fillmore realization of $K$-homology appeared in
[Brown, Douglas \& Fillmore  1977]   \cite{BrDoFi1};  see also the monograph
by [Douglas 1980]  \cite{D1}.  
The Kasparov's  $KK$-theory can be found in [Kasparov 1980]  \cite{Kas1};
see the book by [Blackadar 1986]  \cite{B},  Chapter VIII  for a detailed account.    
 For the proof of special cases of  Novikov's  Conjecture using the $KK$-theory, see
 e.g.  [Kasparov 1984]  \cite{Kas2}.    
The published version of Baum-Connes Conjecture can be found in [Baum \& Connes 2000]
\cite{BaCo1}.  Rosenberg's Theorem on positive scalar curvature appears  in 
[Rosenberg 1983]  \cite{Ros1}.  The ideas and methods of coarse geometry are covered
in the books  by  [Roe 1996]   \cite{R1} and [Nowak \& Yu  2012]  \cite{NY}.    
An interesting link between topology and  operator algebras was 
studied  by [Hughes 2012]  \cite{Hug1}.





\chapter{Jones Polynomials}
The first non-trivial functor on the classical geometry (knots and
links)  with  values in the NCG  (the so-called {\it subfactors} of a finite-dimensional 
$W^*$-algebra) has been constructed in 1980's by 
V.~F.~R.~Jones.  The corresponding noncommutative invariant --
a {\it trace invariant} $V_L(t)$  --
revolutionized the low-dimensional topology.  The  $V_L(t)$
cannot be easily obtained otherwise despite some efforts 
to do so, see  e.g. the Khovanov homology;
this fact  demonstrates the power and beauty  of the NCG.  
Roughly speaking,  the trace invariant $V_L(t)$ comes from a functor
\displaymath
{\cal B}_n\to {\cal A}_n,
\enddisplaymath
 where ${\cal B}_n$ is the category of the $n$-string braid groups $B_n$ and 
 ${\cal A}_n$ a category of the $n$-dimensional $W^*$-algebras $A_n$.  Remarkably,
 the algebras $A_n$  have a trace behaving extraordinary nice with respect to the second
 Markov move of a braid: 
\displaymath
tr~(e_n x) =     {1\over   [{\cal M} : {\cal N}]} ~tr~(x), \quad\forall x\in A_n.
\enddisplaymath
(The lack of such a trace formula for other known representations of the 
braid groups was a major no-go for using braids for the topological classification
of  knots and links obtained by the closure of a braid.)  The trace invariant $V_L(t)$
is a normalization of the above trace formula.  

 \index{subfactor}
 \index{braid group}
 \index{braid}

\section{Subfactors}
We refer the reader to Chapter 8 for  an introduction and notation of the 
$W^*$-algebras.  Roughly speaking,  the theory of subfactors of the $W^*$-algebras
mimics the Galois theory of the finite field extensions. 
\begin{dfn}
By a subfactor of a factor ${\cal M}$ one understands the unital $W^*$-subalgebra of ${\cal M}$
which is a factor itself. 
\end{dfn}
\begin{thm}
{\bf (Galois Theory for type II$_1$ factors)}
If $G$ is a finite group of outer automorphisms of a type {\bf II}$_1$ factor ${\cal M}$
and let $H\subseteq G$ be its subgroup. Then  the formula 
\displaymath
H \leftrightarrow {\cal M}^H:=\{x\in {\cal M} ~|~ \alpha(x)=x, ~\forall\alpha\in H\} 
\enddisplaymath
gives a Galois correspondence between subgroups of $G$ and subfactors of ${\cal M}$
containing the fixed point subalgebra ${\cal M}^G\subseteq {\cal M}$.  
\end{thm}
Recall that to each subgroup $H\subseteq G$ one can assign an integer $[G : H]$
called {\it index}  of the subgroup;  the following definition extends the notion to
the subfactors.  
\begin{dfn}
Let ${\cal N}\subseteq {\cal M}$ be two factors of type {\bf II}$_1$;  then the index 
of ${\cal N}$ in ${\cal M}$ is defined by the formula
\displaymath
[{\cal M} : {\cal N}]={\dim L^2({\cal N}) \over  \dim L^2({\cal M})},
\enddisplaymath
where $L^2({\cal N})$ and $L^2({\cal M})$ are the Hilbert spaces associated to the
GNS construction for $W^*$-algebras ${\cal N}$ and ${\cal M}$,  respectively. 
\end{dfn}
Unlike the case of  groups,  the index of subfactors is not always an integer;  the following
theorem lists all possible values of the index for subfactors of type {\bf II}$_1$ factors.
 \index{Jones Theorem}
\begin{thm}
{\bf (V.~F.~R.~Jones)}
If ${\cal M}$ is  a  type {\bf II}$_1$ factor, then:

\medskip
(i)  if ${\cal N}$ is a subfactor with the same identity and $[{\cal M} : {\cal N}]<4$,
then 
\displaymath
[{\cal M} : {\cal N}]=4\cos^2\left({\pi\over n}\right)
\enddisplaymath
for some $n=3, 4, \dots$;

\smallskip
(ii) if ${\cal M}$ is the hyper-finite type {\bf II}$_1$ factor,  then the index of a
subfactor ${\cal N}$ takes the following values:
\displaymath
[{\cal M} : {\cal N}]=\{ 4\cos^2\left({\pi\over n}\right) ~:~ n\ge 3\} \cup [4,\infty]. 
\enddisplaymath
\end{thm}
\begin{rmk}
\textnormal{
The proof of  Jones'  Theorem exploits  the following {\it basic construction}.  
Take ${\cal M}$  and its GNS representation on the Hilbert space $L^2({\cal M})$;
suppose that ${\cal N}$ is a subfactor.  By uniqueness,   $tr~({\cal N})$ 
is the restriction of $tr~({\cal M})$ and thus $L^2({\cal N})$ is a sub-space of   
$L^2({\cal M})$. The projection 
\displaymath
e:   L^2({\cal M})\to L^2({\cal N})
\enddisplaymath
has the following properties:  (i)  ${\cal M}\cap \{e\}'={\cal N}$;
(ii)   $W^*$-algebra $\langle {\cal M}, e\rangle := ({\cal M}\cup \{e\})''$ 
is a factor;  (iii)  $[\langle {\cal M}, e\rangle : {\cal M}]=[{\cal M} : {\cal N}]$;
(iv)  $tr~(ae)={1\over  [{\cal M} : {\cal N}]} tr~(a)$. 
}
\end{rmk}
 \index{Jones projection}
\begin{dfn}
By the Jones projections one understands a sequence $\{e_i\}_{i=1}^{\infty}$ obtained as the 
result of iteration of the basic construction, i.e.
\displaymath
{\cal N}\subset  {\cal M}\subset  \langle {\cal M}, e_1\rangle \subset
\langle\langle {\cal M}, e_1\rangle, e_2\rangle\subset\dots 
\enddisplaymath
\end{dfn}
\begin{cor}
{\bf (V.~F.~R.~Jones)}
The Jones projections satisfy the following relations and a trace formula: 
\displaymath
\left\{
\begin{array}{cccc}
e_ie_{i\pm 1}e_i  &=& {1\over   [{\cal M} : {\cal N}]} ~e_i&\\ 
e_ie_j &=& e_je_i, & ~\hbox{if} ~|i-j|\ge 2\\ 
tr~(e_n x) &=&     {1\over   [{\cal M} : {\cal N}]} ~tr~(x), &~\forall x\in 
\langle {\cal M}, e_1,\dots, e_{n-1}\rangle.              
\end{array}
\right.
\enddisplaymath
\end{cor}

 \index{braid}
\section{Braids}
Motivated by a topological classification of knots and links in the 
three-dimensional space,  E.~Artin introduced the notion of a braid;
the braid is a simpler object than the knot or link and, moreover,  the
braids can be composed with each other so that one gets a 
finitely generated group.  
\begin{dfn}
By an $n$-string braid $b_n$ one understands two parallel copies of the plane ${\Bbb R}^2$ in 
${\Bbb R}^3$  with $n$ distinguished points taken together with  $n$ disjoint smooth paths (``strings'')
joining pairwise the distinguished points of the planes; the tangent vector to each string
is never parallel to the planes.  The braid $b_n$ can be given by a diagram by projecting 
its strings into a generic plane in ${\Bbb R}^3$ and indicating the over- and underpasses 
of the strings, see Fig. 11.1.   
\end{dfn}
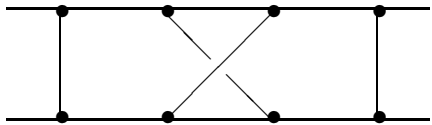
\begin{figure}[here]
\begin{picture}(500,110)(50,-5)

\put(165,82){\line(1,0){160}}
\put(165,40){\line(1,0){160}}

\put(305,80){\line(0,-1){40}}
\put(185,40){\line(0,1){40}}

\put(225,80){\line(1,-1){17}}
\put(248,57){\line(1,-1){18}}

\put(225,40){\line(1,1){40}}

\put(223,78){$\bullet$}
\put(263,78){$\bullet$}
\put(263,38){$\bullet$}
\put(223,38){$\bullet$}

\put(183,78){$\bullet$}
\put(303,78){$\bullet$}
\put(183,38){$\bullet$}
\put(303,38){$\bullet$}

\end{picture}

\caption{The diagram of a braid $b_4$.}
\end{figure}
\begin{rmk}
\textnormal{
The braids $b_n$ are endowed with a natural equivalence relation: 
two braids $b_n$ and $b_n'$ are equivalent if $b_n$ can be deformed 
into $b_n'$ without intersection of the strings and so that at each moment of the
deformation $b_n$ remains  a braid.   
}
\end{rmk}
 \index{braid closure}
\begin{dfn}
By a closure $\hat b_n$ of the braid $b_n$ one understands a link or knot in ${\Bbb R}^3$
obtained by gluing the endpoints of strings at the top of the braid with such at the bottom
of the braid.  
\end{dfn}
\begin{dfn}
By an $n$-string braid group $B_n$ one understands the set of all $n$-string braids $b_n$
endowed with a multiplication operation of the concatenation of  $b_n\in B_n$ and $b_n'\in B_n$, 
i.e the identification of the bottom of $b_n$ with the top of $b_n'$.  The group is non-commutative and the 
identity is given by the trivial braid,  see Fig. 11.2 
\end{dfn}
\begin{figure}[here]
\begin{picture}(500,110)(50,-5)

\put(165,82){\line(1,0){160}}
\put(165,40){\line(1,0){160}}

\put(305,80){\line(0,-1){40}}
\put(185,40){\line(0,1){40}}

\put(225,80){\line(1,-1){17}}
\put(248,57){\line(1,-1){18}}

\put(225,40){\line(1,1){40}}

\put(223,78){$\bullet$}
\put(263,78){$\bullet$}
\put(263,38){$\bullet$}
\put(223,38){$\bullet$}

\put(183,78){$\bullet$}
\put(303,78){$\bullet$}
\put(183,38){$\bullet$}
\put(303,38){$\bullet$}

\put(140,60){$b_4$}

\end{picture}

\begin{picture}(500,110)(50,-5)

\put(165,82){\line(1,0){160}}
\put(165,40){\line(1,0){160}}

\put(305,80){\line(0,-1){40}}
\put(185,40){\line(0,1){40}}

\put(225,80){\line(1,-1){40}}

\put(225,40){\line(1,1){17}}
\put(248,62){\line(1,1){17}}

\put(223,78){$\bullet$}
\put(263,78){$\bullet$}
\put(263,38){$\bullet$}
\put(223,38){$\bullet$}

\put(183,78){$\bullet$}
\put(303,78){$\bullet$}
\put(183,38){$\bullet$}
\put(303,38){$\bullet$}

\put(140,60){$b_4^{-1}$}

\put(245,110){$\times$}

\end{picture}

\begin{picture}(500,110)(50,-5)

\put(165,82){\line(1,0){160}}
\put(165,40){\line(1,0){160}}

\put(305,80){\line(0,-1){40}}
\put(185,40){\line(0,1){40}}

\put(265,40){\line(0,1){40}}
\put(225,40){\line(0,1){40}}

\put(223,78){$\bullet$}
\put(263,78){$\bullet$}
\put(263,38){$\bullet$}
\put(223,38){$\bullet$}

\put(183,78){$\bullet$}
\put(303,78){$\bullet$}
\put(183,38){$\bullet$}
\put(303,38){$\bullet$}

\put(140,60){$Id_4$}

\put(245,110){$||$}

\end{picture}

\caption{The concatenation  of a braids  $b_4$ and $b_4^{-1}$.}
\end{figure}
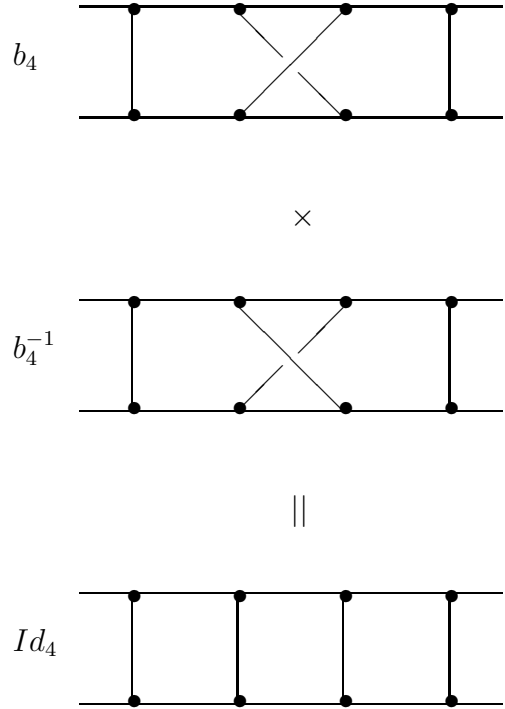
 \index{Artin's representation of braid group}
\begin{thm}
{\bf (E.~Artin)}
The braid group $B_n$ is isomorphic to a group on the standard generators 
$\sigma_1, \sigma_2,\dots,\sigma_{n-1}$ (like the one shown in Fig.11.1) satisfying
the following relations: 
\displaymath
\left\{
\begin{array}{cccc}
\sigma_i\sigma_{i+1}\sigma_i  &=&\sigma_{i+1}\sigma_i\sigma_{i+1} &\\ 
\sigma_i\sigma_j &=& \sigma_j\sigma_i, & ~\hbox{if} ~|i-j|\ge 2
\end{array}
\right.
\enddisplaymath
\end{thm}
 \index{Markov move}
\begin{dfn}
Let $b_n\in B_n$ be an $n$-string braid;   by a Markov move one understand the 
following two types of transformations of the braid $b_n$:

\medskip
(i) {\bf Type I}:  $b_n\mapsto a_nb_na_n^{-1}$ for a braid $a_n\in B_n$;

\smallskip
(ii) {\bf Type II}: $b_n\mapsto b_n\sigma_n^{\pm 1}\in B_{n+1}$, 
where $\sigma_n\in B_{n+1}$.  
\end{dfn}
 \index{Markov Theorem}
\begin{thm}
{\bf (A.~Markov)}
The closure of any two braids $b_n\in B_n$ and $b_m\in B_m$ give the same link 
$\hat b_n\cong \hat b_m$  in ${\Bbb R}^3$ if and only if  $b_n$ and $b_m$  can be connected by a 
sequence of the Markov moves of types I and II.  
\end{thm}

 \index{trace invariant}
 \index{Jones algebra}
\section{Trace invariant}
\begin{dfn}
By the Jones algebra $A_n$ we shall understand an $n$-dimensional $W^*$-algebra
generated by the identity and the Jones projections $e_1,e_2,\dots, e_n$;  the projections
are known to satisfy the  relations:
\displaymath
\left\{
\begin{array}{cccc}
e_ie_{i\pm 1}e_i  &=& {1\over   [{\cal M} : {\cal N}]} ~e_i&\\ 
e_ie_j &=& e_je_i, & ~\hbox{if} ~|i-j|\ge 2,
\end{array}
\right.
\enddisplaymath
and the trace formula:
\displaymath
tr~(e_n x) =     {1\over   [{\cal M} : {\cal N}]} ~tr~(x), \quad\forall x\in A_n.
\enddisplaymath
\end{dfn}
\begin{rmk}
\textnormal{
The reader can verify that  the relations for the Jones projections $e_i$ coincide with
such for the  generators $\sigma_i$ of the braid group  after  a minor adjustment 
of the notation:
\displaymath
\left\{
\begin{array}{ccc}
\sigma_i &\mapsto & \sqrt{t} [(t+1)e_i-1]\\ 
\left[  {\cal M} : {\cal N} \right] &=& 2+t+{1\over t}.
\end{array}
\right.
\enddisplaymath
 Thus one gets a family $\rho_t$ of  representations of the braid group $B_n$ into 
the Jones algebra $A_n$.   However,  to get  a topological invariant
of the closed braid $\hat b$ of $b\in B_n$ coming from the 
trace (a character) of the representation,  one needs to choose 
a representation whose trace is invariant under the first and the second 
Markov moves of the braid $b$.  There is no problem with the first move,
because  two similar matrices have the same trace
for any representation from the family  $\rho_t$.    For the second Markov move,
we luckily have the trace formula which (after obvious substitutions) takes the
form: 
\displaymath
\left\{
\begin{array}{ccc}
tr~(b\sigma_n) &=& -{1\over t+1} tr~(b)\\ 
tr~(b\sigma_n^{-1}) &=& -{t\over t+1} tr~(b).
\end{array}
\right.
\enddisplaymath
Thus in general $tr~(b\sigma_n^{\pm 1})\ne tr~(b)$,  but one can always
re-scale the trace to get the equality.  Indeed,  the second Markov move
takes the braid from $B_i$ and replaces it by a braid from $B_{i-1}$;  
there is  a finite number of such replacements because  the algorithm stops
for $B_1$.  Therefore  a finite number of re-scalings by the constants    
 $-{1\over t+1}$ and     $-{t\over t+1}$ will give a quantity invariant under
 the second Markov move;  the quantity is known as the {\it Jones polynomial}
 of the closed braid $\hat b$.
}
\end{rmk}
 \index{Jones Polynomial}
\begin{thm}
{\bf (V.~F.~R.~Jones)}
Let $b\in B_n$ be a braid and $\exp(b)$ be the sum of all  powers of generators $\sigma_i$
and $\sigma_i^{-1}$  in the word presentation of $b$;  let $L:=\hat b$ be the closure of $b$.
Then  the number 
\displaymath
V_{L}(t):=\left(-{t+1\over\sqrt{t}}\right)^{n-1}(\sqrt{t})^{\exp(b)}~tr~(b)
\enddisplaymath
is an isotopy  invariant of the link $L$.   
\end{thm}
\begin{exm}
\textnormal{
Let $b\in B_2$ be a braid whose closure is isotopic to the trefoil knot.  
Then 
\displaymath
V_{L}(t):=\left(-{t+1\over\sqrt{t}}\right)(\sqrt{t})^{3}\left({t^3-t^2-1\over t+1}\right)=-t^4+t^3+t.
\enddisplaymath
}
\end{exm}
\begin{exm}
\textnormal{
Let $b\in B_2$ be a braid whose closure is isotopic to a pair of linked circles $S^1\cup S^1$.  
Then 
\displaymath
V_{L}(t):=\left(-{t+1\over\sqrt{t}}\right)(\sqrt{t})^{2}\left({t^2+1\over t+1}\right)=-\sqrt{t}(t^2+1).
\enddisplaymath
}
\end{exm}

\vskip1cm\noindent
{\bf Guide to the literature.}
E.~Artin introduced the braid groups in the seminal paper [Artin 1925]  \cite{Art2},
see also more accessible [Artin 1947]  \cite{Art3}.  The equivalence classes of braids
were studied by [Markov 1935]  \cite{Mar1}.   
For the index of subfactors the reader is referred to [Jones 1983]  \cite{Jon1}.
The trace invariant was introduced in [Jones 1985]  \cite{Jon2}.  
For an extended exposition of subfactors and knots we refer the reader to the monograph 
[Jones 1991]  \cite{J1}.





 \index{non-commutative Pontryagin duality}

\chapter{Quantum Groups}
The quantum groups can be viewed as a special family  of non-commutative 
rings labeled by  a ``quantization parameter'' $q$ such that  $q=1$ corresponds
to a commutative ring.  The quantum groups are  closely related to the quantum mechanics 
and a non-commutative version of the Pontryagin duality for abelian groups.  
Let us consider the simplest example of such a group.  
\begin{exm}
\textnormal{
Let $A(k)$ be an associative algebra over the field $k$ and let 
\displaymath
M_q(2)=\left(\matrix{a & b\cr c & d}\right)
\enddisplaymath
 be the set of two-by-two martices with entries $a, b, c, d\in A(k)$  satisfying the 
 commutation relations
\displaymath
\left\{
\begin{array}{ccc}
ab &= &  q ~ba\\ 
ac &= &  q ~ca\\ 
bd &= &  q ~db\\ 
cd &= &  q ~dc\\
bc &= &  cb\\
ad-da &= & \left(q- {1\over q}\right) bc,  
\end{array}
\right.
\enddisplaymath
where $q\in k$ is a parameter.  The algebra $A(k)$ can be viewed as  the ``algebra of functions''
on $M_q(2)$ and we   denote it by   ${\cal F}(M_q(2))$.  
Notice that  if $q=1$,  then ${\cal F}(M_q(2))$ is commutative and, therefore,
the product of two matrices  gives a matrix with entries in  ${\cal F}(M_q(2))$;   if $q\ne 1$,  it is false 
because such a product gives a matrix whose  entries in general do not satisfy  the commutation
relations.  However,  if   ${\cal F}(M_q(2))$ is a {\it Hopf algebra},  i.e. a bi-algebra endowed 
with multiplication 
\displaymath
\mu:  {\cal F}(M_q(2))\otimes {\cal F}(M_q(2))\to {\cal F}(M_q(2)),
\enddisplaymath
and {\it co-multiplication} 
\displaymath
\Delta:  {\cal F}(M_q(2))\to  {\cal F}(M_q(2))\otimes {\cal F}(M_q(2)),
\enddisplaymath
then  $M_q(2)$  becomes the required semi-group.  To get the structure of a 
group, one also needs to invert the quantum determinant.  The resulting group
$M_q(2)$ is called a {\it quantum group}.      
}
\end{exm}

 \index{quantum group}

 \index{Manin's quantum plane}

\section{Manin's quantum plane}
\begin{dfn}
By  Manin's quantum plane $P_q$ one understands the Euclidean plane whose
coordinates $(x,y)$ satisfy the commutation relation
\displaymath
xy=q~yx
\enddisplaymath
Such a plane is shown in Fig. 12.1. 
\end{dfn}
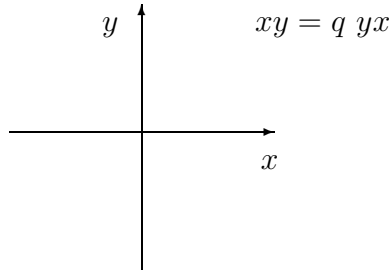
\begin{figure}[here]
\begin{picture}(500,110)(0,-5)

\put(150,62){\vector(1,0){100}}
\put(200,10){\vector(0,1){100}}
\put(245,48){$x$}
\put(185,100){$y$}
\put(243,100){$xy=q~yx$}

\end{picture}
\caption{Quantum plane $P_q$.}
\end{figure}
\begin{thm}
{\bf (Manin)}
The group of linear transformations of  quantum plane $P_q$ is isomorphic 
to the quantum group $M_q(2)$. 
\end{thm}
{\it Proof.}   The linear transformation of quantum plane $P_q$ can be written
in the form
\displaymath
\left(\matrix{x'\cr y'}\right)=
\left(\matrix{a & b\cr c & d}\right)\left(\matrix{x\cr y}\right). 
\enddisplaymath
Using the transformation,  one gets the following system of equations
\displaymath
\left\{
\begin{array}{ccc}
x'y' &= &  ac~x^2+bc~yx+ad~xy+bd~y^2\\ 
y'x' &= &  ca~x^2+da~yx+cb~xy+db~y^2.
\end{array}
\right.
\enddisplaymath
Because $x'y'=q~y'x'$,  we have the following identities
\displaymath
\left\{
\begin{array}{ccc}
ac &= & q~ca  \\ 
bc &= & q~da  \\
ad &= & q~cb  \\
bd &= & q~db.
\end{array}
\right.
\enddisplaymath
The first and the last equations are exactly the commutation relations
for  the quantum group $M_q(2)$;   it is easy to see that 
whenever $bc=cb$,  it follows from  $bc=q~da$ and  $ad=q~cb$, that
\displaymath
ad-da =  \left(q- {1\over q}\right) bc.  
\enddisplaymath
We leave it as exercise to the reader to prove remaining relations $bc=cb,  ab=q~ba$ 
and $cd=q~dc$. 
$\square$

\begin{rmk}
\textnormal{
The quantum group $M_q(2)$ can be viewed as a group of symmetries
of the quantum plane $P_q$.  
}
\end{rmk}

 \index{Hopf algebra}

\section{Hopf algebras}
\begin{dfn}
By  a Hopf algebra one understands a bi-algebra $A$ endowed with
a product
\displaymath
\mu:  A\otimes A\to A,
\enddisplaymath
and  co-product 
\displaymath
\Delta:  A\to  A\otimes A,
\enddisplaymath
such that $\mu$ and $\Delta$ are homomorphisms of (graded) algebras.
\end{dfn}
\begin{rmk}
\textnormal{
Roughly speaking, in case $A\cong {\cal F}(M_q(2))$ the product structure $\mu$ on the Hopf algebra $A$
ensures that $M_q(2)$ is a semi-group,  while the co-product structure $\Delta$ (together with an inverse of the
  determinant)  turns $M_q(2)$  into a  group.
}
\end{rmk}
 \index{co-product}
\begin{exm}
\textnormal{
Let $A$ be an algebra of regular functions on an affine algebraic group $G$;
let us define  a  product $\mu$ on $A$ coming from multiplication
\displaymath
 G\times G\to G,
\enddisplaymath
and a co-product $\Delta$ coming from the emebdding
\displaymath
\{e\}\to G,
\enddisplaymath
where $e$ is the unit of $G$.  Then the algebra $A$ is a Hopf algebra.    
}
\end{exm}
\begin{exm}
\textnormal{
Let $X$ be a connected topological group with multiplication
\displaymath
m: X\times X\to X
\enddisplaymath
and co-multiplication
\displaymath
\delta: X\to X\times X
\enddisplaymath
defined by the embedding $\iota: \{e\}\to X$.  Denote by $H^*(X)$ the cohomology
ring of $X$;  let 
\displaymath
\mu:=m^*: H^*(X)\otimes H^*(X)\to H^*(X)
\enddisplaymath
be a product induced by $m$ on the cohomology ring $H^*(X)$ and let 
\displaymath
\Delta:= \delta^*:  H^*(X)\to H^*(X)\otimes H^*(X)
\enddisplaymath
be the corresponding co-product.  The algebra $A\cong H^*(X)$ is a Hopf algebra.
}
\end{exm}
\begin{rmk}
\textnormal{
The algebra $H^*(X)$ has been introduced in algebraic topology by Heinz Hopf;  
hence the name. 
}
\end{rmk}

 \index{non-commutative Pontryagin duality}

\section{Operator algebras and quantum groups}
An analytic approach to quantum groups is based on an observation that the
locally compact {\it abelian} groups $G$ correspond  to commutative $C^*$-algebras
$C_0(G)$;  therefore a non-commutative version of {\it Pontryagin duality} 
must involve more general $C^*$-algebras. It was suggested by 
[Woronowicz 1987]  \cite{Wor1}  that such $C^*$-algebras come from an
embedding of the Hopf algebra into a $C^*$-algebra coherent with 
the co-product structure of the  Hopf algebra;  this $C^*$-algebra is known
as a {\it compact quantum group}.  To give some details,  recall that if $G$ is 
a locally compact group  then multiplication is a continuous map $G\times G$ to $G$
and therefore can be translated to a morphism
\displaymath
\Delta:  C_0(G)\to M(C_0(G\times G)),
\enddisplaymath
where $M(C_0(G\times G))$ is a {\it multiplier algebra} of  $C_0(G\times G)$,
i.e.  the maximal $C^*$-algebra containing $C_0(G\times G)$ as an essential ideal,
see e.g. [Blackadar  1986]  \cite{B},  Section 12;  the translation itself is given 
by the formula $\Delta(f)(p,q)=f(pq)$.   Define a morphism $\iota$
by a co-associativity formula $(\Delta ~\otimes ~\iota)\Delta=(\iota ~\otimes ~\Delta)\Delta$.   
 \index{compact quantum group}
\begin{dfn}
{\bf (Woronowicz)}
If $A$ is a unital $C^*$-algebra together with a unital
$\ast$-homomorphism $\Delta: A\to A\otimes A$ such that
 $(\Delta ~\otimes ~\iota)\Delta=(\iota ~\otimes ~\Delta)\Delta$
 and such that the spaces  $\Delta(A)(A\otimes 1)$ and $\Delta(A)(1\otimes A)$ are dense in
$A\otimes A$,   then the pair $(A,\Delta)$ is called a compact quantum 
group.
\end{dfn}
Let $\varphi$ be weight on a $C^*$-algebra $A$, see Section 9.1.1;  we denote by ${\cal M}{\varphi}^+$
the space of all integrable elements of $A$, i.e. 
${\cal M}_{\varphi}^+:=\{x\in A^+ ~|~ \varphi(x)<\infty\}$.   
The following definition generalizes  the compact quantum groups to the 
locally compact case. 
 \index{locally compact quantum group}
\begin{dfn}
{\bf (Kustermans \& Vaes)}
Suppose that $A$ is a unital $C^*$-algebra together with a unital
$\ast$-homomorphism $\Delta: A\to A\otimes A$ such that
 $(\Delta ~\otimes ~\iota)\Delta=(\iota ~\otimes ~\Delta)\Delta$
 and such that the spaces  $\Delta(A)(A\otimes 1)$ and $\Delta(A)(1\otimes A)$ are dense in
$A\otimes A$.   Moreover,  assume that

\medskip
(i) there exists a faithful KMS weight $\varphi$ (see Section 9.1.1)  on the compact quantum
group $(A,\Delta)$ such that $\varphi((\omega\otimes\iota)\Delta(x))=\varphi(x)\omega(1)$
for $\omega\in A_+$ and $x\in {\cal M}_{\varphi}^+$;

\smallskip
(ii) there exists a KMS weight $\psi$  on the compact quantum
group $(A,\Delta)$ such that $\psi((\iota\otimes\omega)\Delta(x))=\psi(x)\omega(1)$
for $\omega\in A_+$ and $x\in {\cal M}_{\psi}^+$. 

\medskip
Then the pair $(A,\Delta)$ is called a locally compact quantum group.
\end{dfn}

\vskip1cm\noindent
{\bf Guide to the literature.}
There exits an extensive body of literature on the quantum groups; 
the main references   are the  textbooks by [Kassel 1995]  \cite{KAS}
and [Timmermann 2008]  \cite{TIM}.
The Hopf algebras were introduced in  [Hopf  1941]   \cite{Hop1}.  
Compact quantum groups were defined by [Woronowicz 1987]  \cite{Wor1}
and locally compact quantum groups by  [Kustermans \& Vaes  2000]  \cite{KuVa2}.  
Locally compact quantum groups in  the context of $W^*$-algebras
were considered by [Kustermans \& Vaes  2003]  \cite{KuVa1}.





\chapter{Non-commutative Algebraic Geometry}
The subject grew from the works of  [Serre 1957] \cite{Ser1}
and [Sklyanin 1982] \cite{Skl1}.  Roughly speaking,  Sklyanin's construction
exploits a natural transformation {\bf det}  of the functor 
{\bf GL}$_n:$ {\bf CRng} $\to$ {\bf Grp},  see Chapter 2 for the notation;
let us outline the main ideas.  
Denote by {\bf CRng} $\cong k[x_1,\dots,x_m]$ a polynomial ring over the 
algebraically closed field $k$ and  $m$ variables $x_i$.  Consider the 
functor 
\displaymath
\hbox{{\bf GL}}_n: \hbox{{\bf CRng}} \to \hbox{{\bf Grp}}
\enddisplaymath
from $k[x_1,\dots,x_m]$ to the multiplicative group of $n\times n$
matrices with entries in $k[x_1,\dots,x_m]$.   Suppose that our 
matrix group can be given by a finite number of generators satisfying
a finite number of relations $w_1=w_2=\dots=w_r=I$.  
Then taking the determinants of the relations, one gets the following
system of polynomial equations
\displaymath
\left\{
\begin{array}{ccc}
det~(w_1)  &=& 1\\ 
det~(w_2) &=& 1\\ 
\vdots && \\
det~(w_r) &=& 1,                   
\end{array}
\right.
\enddisplaymath
which define an algebraic variety $V$ over the field $k$. 
 If $I\subset k[x_1,\dots,x_m]$ is an ideal generated by the above equations, 
 then the coordinate ring $R\cong k[x_1,\dots,x_m]/I$ of variety $V$  satisfies the natural transformation 
 diagram of Fig. 13.1;  the functor {\bf U}$_R: R\to R^*$ takes the commutative ring 
 $R$ to its group $R^*$  of units (invertible elements).  In this case {\bf det}$_R$ is a natural
 transformation (an isomorphism) between the functors {\bf GL}$_n$ and {\bf U}$_R$. 
\begin{figure}[here]
\begin{picture}(300,100)(-50,0)

\put(120,93){\vector(1,0){60}}
\put(110,80){\vector(1,-1){40}}
\put(160,40){\vector(1,1){40}}

\put(100,90){$R$}
\put(200,90){$R^*$}
\put(140,20){$GL_n(R)$}

\put(150,99){\hbox{{\bf U}}$_R$}
\put(110,50){\hbox{{\bf GL}}$_n$}
\put(190,50){\hbox{{\bf det}}$_R$}

\end{picture}
\caption{Sklyanin algebra $GL_n(R)$  as a natural transformation {\bf det}$_R$.}
\end{figure}
\begin{rmk}
\textnormal{
Sklyanin's original result deals with the functor {\bf GL}$_2$ and the group $GL_2(R)$
given by four generators and six quadratic relations;  the algebraic variety 
$V$ in this case is isomorphic to the product of two elliptic curves over the field
of complex numbers, i.e. $V={\cal E}({\Bbb C})\times {\cal E}({\Bbb C})$.  
}
\end{rmk}

 \index{Serre isomorphism}
 \index{Gelfand Duality}
 \index{quasi-coherent sheaf}
 \index{finitely generated graded module}
\section{Serre isomorphism}
If $X$ is a projective variety,  $\Coh~(X)$ a category of  the quasi-coherent
sheaves on $X$ and  $\Mod~(B)$ a category of 
the finitely generated graded modules  over $B$ factored by a torsion $\Tors$,
then the well known {\it Serre isomorphism}  says that
\displaymath
\Coh~(X)\cong \Mod~(B)~/~\Tors, 
\enddisplaymath
see  [Serre 1955]  \cite{Ser1}.  Because $\Mod~(B)$ is correctly defined for non-com\-muta\-tive rings
$B$ one can wonder if there are concrete examples of such rings for given variety $X$;
the answer is emphatically positive and plenty of such rings are known, 
see  ([Stafford \& van den Bergh 2001]  \cite{StaVdb1}, pp. 172-173)  for  details.
The  non-commutative ring $B$ satisfying the fundamental duality is called  
{\it  twisted  homogeneous coordinate ring} of   $X$;  an exact definition will be given in Section13.2.  
To convey  the basic idea,  consider the simplest example. 
If $X$ is  a compact Hausdorff  space and $C(X)$ the commutative algebra of  
continuous functions from $X$ to ${\Bbb C}$,  then topology of $X$ is determined   by
 algebra $C(X)$  (the Gelfand Duality);
in terms of  K-theory this  can be written as  
$K_0^{top}(X)\cong K_0^{alg}(C(X))$.   Taking  the two-by-two matrices 
with entries in $C(X)$,  one  gets  an  algebra 
$C(X)\otimes M_2({\Bbb C})$;  in view of stability of  K-theory
under  tensor products (e.g. [Blackadar 1986]  \cite{B}, Chapter 5),  it holds     
$K_0^{top}(X)\cong K_0^{alg}(C(X))\cong K_0^{alg}(C(X)\otimes M_2({\Bbb C}))$.  
In other words,  the topology of $X$ is defined by 
algebra $C(X)\otimes M_2({\Bbb C})$,  which is no longer a commutative algebra.  
In algebraic geometry,  one replaces $X$
by a  projective variety, $C(X)$ by its  coordinate ring,
$C(X)\otimes M_2({\Bbb C})$ by a twisted  homogeneous  
coordinate ring  of $X$ and  $K^{top}(X)$ by a category of the 
quasi-coherent sheaves on $X$.
The simplest concrete example of  $B$ is as follows. 
\begin{exm}
\textnormal{
Let $k$ be a field and $U_{\infty}(k)$ the algebra of polynomials
over $k$ in two non-commuting variables $x_1$ and $x_2$ satisfying  the quadratic relation
\displaymath
x_1x_2-x_2x_1-x_1^2=0.
\enddisplaymath
 If ${\Bbb P}^1(k)$ be the projective
line over $k$,   then $B=U_{\infty}(k)$ and $X={\Bbb P}^1(k)$
satisfy the fundamental duality $\Coh~(X)\cong \Mod~(B)~/~\Tors$.  
Notice, that $B$ is far from being a commutative ring.  
}
\end{exm}
\begin{exm}
\textnormal{
Let $k$ be a field and $U_{q}(k)$ the algebra of polynomials
over $k$ in two non-commuting variables $x_1$ and $x_2$ satisfying  the quadratic relation
\displaymath
x_1x_2=q x_2x_1,
\enddisplaymath
where $q\in k^*$ is a non-zero element of $k$.  
 If ${\Bbb P}^1(k)$ be the projective
line over $k$,   then $B=U_{q}(k)$ and $X={\Bbb P}^1(k)$
satisfy the fundamental duality $\Coh~(X)\cong \Mod~(B)~/~\Tors$
for all $q\in k^*$.   Again,  $B$ is an essentially non-commutative ring.  
}
\end{exm}
In general,  there exists a canonical non-commutative ring $B$,  attached to the
projective variety $X$ and an automorphism  $\alpha: X\to X$;
we refer the reader to [Stafford \& van den Bergh 2001]  \cite{StaVdb1},  pp. 180-182.
To give an idea,  let $X=Spec~(R)$ for a commutative graded
ring $R$.  One considers the ring $B:=R[t,t^{-1}; \alpha]$ of skew
Laurent polynomials defined by the commutation relation
\displaymath
b^{\alpha}t=tb,
\enddisplaymath
for all $b\in R$, where  $b^{\alpha}\in R$ is the image of $b$ under automorphism
$\alpha$;   then $B$ satisfies the isomorphism
$\Coh~(X)\cong \Mod~(B)~/~\Tors$.  The ring $B$ is non-commutative,  unless $\alpha$
is the  trivial automorphism of $X$.  
\begin{exm}
\textnormal{
The  ring  $B=U_{\infty}(k)$
 corresponds to the automorphism $\alpha(u)=u+1$ 
of the projective line ${\Bbb P}^1(k)$.   Indeed,  $u=x_2x_1^{-1}=x_1^{-1}x_2$
and,  therefore,  $\alpha$ maps $x_2$ to $x_1+x_2$;  if  one substitutes
$t=x_1,  b=x_2$ and $b^{\alpha}=x_1+x_2$  in equation $b^{\alpha}t=tb$,   then  one  gets the defining 
relation  $x_1x_2-x_2x_1-x_1^2=0$ for  the algebra    $U_{\infty}(k)$.   
}
\end{exm}
\begin{exm}
\textnormal{
The  ring  $B=U_{q}(k)$
 corresponds to the automorphism $\alpha(u)=qu$ 
of the projective line ${\Bbb P}^1(k)$.   Indeed,  $u=x_2x_1^{-1}=x_1^{-1}x_2$
and,  therefore,  $\alpha$ maps $x_2$ to $qx_2$;  if  one substitutes
$t=x_1,  b=x_2$ and $b^{\alpha}=qx_2$  in equation $b^{\alpha}t=tb$,   then  one  gets the defining 
relation  $x_1x_2=qx_2x_1$ for  the algebra    $U_{q}(k)$.  
}
\end{exm}

 \index{twisted homogeneous coordinate ring}

\section{Twisted homogeneous coordinate rings}
Let $X$ be a projective scheme over a field $k$, and let ${\cal L}$ 
be the invertible sheaf ${\cal O}_X(1)$ of linear forms on $X$.  Recall
that the homogeneous coordinate ring of $X$ is a graded $k$-algebra, 
which is isomorphic to the algebra
\displaymath
B(X, {\cal L})=\bigoplus_{n\ge 0} H^0(X, ~{\cal L}^{\otimes n}). 
\enddisplaymath
Denote by $\Coh$ the category of quasi-coherent sheaves on a scheme $X$
and by $\Mod$ the category of graded left modules over a graded ring $B$.  
If $M=\oplus M_n$ and $M_n=0$ for $n>>0$, then the graded module
$M$ is called {\it right bounded}.  The  direct limit  $M=\lim M_{\alpha}$
is called a {\it torsion}, if each $M_{\alpha}$ is a right bounded graded
module. Denote by $\Tors$ the full subcategory of $\Mod$ of the torsion
modules.  The following result is basic about the graded ring $B=B(X, {\cal L})$.   
\begin{thm}
{\bf (Serre)}
\quad $\Mod~(B) ~/~\Tors \cong \Coh~(X).$
\end{thm}
Let $\alpha$ be an automorphism of $X$.  The pullback of sheaf ${\cal L}$ 
along $\alpha$ will be denoted by ${\cal L}^{\alpha}$,  i.e. 
${\cal L}^{\alpha}(U):= {\cal L}(\alpha U)$ for every $U\subset X$. 
We shall set
\displaymath
B(X, {\cal L}, \alpha)=\bigoplus_{n\ge 0} H^0(X, ~{\cal L}\otimes {\cal L}^{\alpha}\otimes\dots
\otimes  {\cal L}^{\alpha^{ n}}). 
\enddisplaymath
  The multiplication of sections is defined by the rule
 \displaymath
 ab=a\otimes b^{\alpha^m},
 \enddisplaymath
 whenever $a\in B_m$ and $b\in B_n$. 
 \begin{dfn}
 Given a pair $(X,\alpha)$ consisting of a Noetherian scheme $X$ and 
 an automorphism $\alpha$ of $X$,  an invertible sheaf ${\cal L}$ on $X$
 is called {\it $\alpha$-ample}, if for every coherent sheaf ${\cal F}$ on $X$,
 the cohomology group 
 \displaymath
 H^q(X, ~{\cal L}\otimes {\cal L}^{\alpha}\otimes\dots
\otimes  {\cal L}^{\alpha^{ n-1}}\otimes {\cal F})
\enddisplaymath
vanishes for $q>0$ and $n>>0$.  
\end{dfn}
\begin{rmk}
\textnormal{
Notice,  that if $\alpha$ is trivial,  this definition is equivalent to the
usual definition of ample invertible sheaf [Serre 1955]  \cite{Ser1}.    
}
\end{rmk}
A  non-commutative generalization of the Serre theorem is as follows.
\begin{thm}
{\bf (Artin \& van den Bergh)}
Let $\alpha: X\to X$ be an automorphism of a projective scheme $X$
over $k$  and let ${\cal L}$ be a $\alpha$-ample invertible sheaf on $X$. If
$B(X, {\cal L}, \alpha)\cong \bigoplus_{n\ge 0} H^0(X, ~{\cal L}\otimes {\cal L}^{\alpha}\otimes\dots
\otimes  {\cal L}^{\alpha^{ n}})$,   then
\displaymath
\Mod~(B(X, {\cal L}, \alpha)) ~/~\Tors \cong \Coh~(X).  
\enddisplaymath
\end{thm}
 \index{Artin-van den Bergh Theorem}
\begin{rmk}
\textnormal{
The question of which invertible sheaves are $\alpha$-ample is fairly
subtle, and there is   no characterization of the automorphisms $\alpha$
for which such an invertible sheaf exists.  However, in many important
special cases this problem is solvable,  see   [Artin \& van den Bergh  1990]  \cite{ArtVdb1},  Corollary 1.6.
}
\end{rmk}
 \index{$\alpha$-ample invertible sheaf}
 \begin{dfn}
 For an automorphism  $\alpha: X\to X$  of a projective scheme $X$
and $\alpha$-ample invertible sheaf  ${\cal L}$ on $X$  the ring  
 \displaymath
B(X, {\cal L}, \alpha)\cong \bigoplus_{n\ge 0} H^0(X, ~{\cal L}\otimes {\cal L}^{\alpha}\otimes\dots
\otimes  {\cal L}^{\alpha^{ n}})
\enddisplaymath
is called  a twisted homogeneous coordinate ring of $X$.  
\end{dfn}

 \index{Sklyanin algebra}

\section{Sklyanin algebras}
\begin{dfn}
If  $k$ is a field of $char~k\ne 2$,  then by  a   Sklyanin algebra 
${\goth S}_{\alpha,\beta,\gamma}(k)$ one understands 
 a free $k$-algebra  on  four generators  $x_i$ and  six  quadratic relations
\displaymath
\left\{
\begin{array}{ccc}
x_1x_2-x_2x_1 &=& \alpha(x_3x_4+x_4x_3),\\
x_1x_2+x_2x_1 &=& x_3x_4-x_4x_3,\\
x_1x_3-x_3x_1 &=& \beta(x_4x_2+x_2x_4),\\
x_1x_3+x_3x_1 &=& x_4x_2-x_2x_4,\\
x_1x_4-x_4x_1 &=& \gamma(x_2x_3+x_3x_2),\\ 
x_1x_4+x_4x_1 &=& x_2x_3-x_3x_2,
\end{array}
\right.
\enddisplaymath
where $\alpha,\beta,\gamma\in k$ are such that
\displaymath
 \alpha+\beta+\gamma+\alpha\beta\gamma=0.
\enddisplaymath
\end{dfn}
Let $\alpha\not\in \{0;\pm 1\}$ and consider  a  non-singular  elliptic  curve ${\cal E}(k)\subset {\Bbb P}^3(k)$ 
given by the intersection of two quadrics  
\displaymath
\left\{
\begin{array}{ccc}
u^2+v^2+w^2+z^2 &=& 0,\\
{1-\alpha\over 1+\beta}v^2+{1+\alpha\over 1-\gamma}w^2+z^2 &=& 0,
\end{array}
\right.
\enddisplaymath
and  an automorphism $\sigma: {\cal E}(k)\to {\cal E}(k)$ acting on the points of ${\cal E}(k)$ according to 
the formula
\displaymath
\left\{
\begin{array}{ccc}
\sigma(u) &=& -2\alpha\beta\gamma vwz-u(-u^2+\beta\gamma v^2+\alpha\gamma w^2+\alpha\beta z^2,\\
\sigma(v) &=& 2\alpha uwz + v(u^2-\beta\gamma v^2+\alpha\gamma w^2+\alpha\beta z^2,\\
\sigma(w) &=& 2\beta vwz + w(u^2+\beta\gamma v^2 - \alpha\gamma w^2 + \alpha\beta z^2,\\
\sigma(z) &=&  2\gamma uvw + z(u^2 +\beta\gamma v^2 +\alpha\gamma w^2 - \alpha\beta z^2.
\end{array}
\right.
\enddisplaymath
 \index{Sklyanin Theorem}
\begin{thm}
{\bf (Sklyanin)}
The algebra ${\goth S}_{\alpha,\beta,\gamma}(k)$ is a twisted homogeneous coordinate ring 
of the elliptic curve ${\cal E}(k)$ and  automorphism $\sigma: {\cal E}(k)\to {\cal E}(k)$,  i.e.
\displaymath
\Mod~ ({\goth S}_{\alpha,\beta,\gamma}(k)~/~\Omega)\cong  \Coh~({\cal E}(k)),
\enddisplaymath
where $\Mod$ is a category of the quotient graded modules over 
the algebra ${\goth S}_{\alpha,\beta,\gamma}(k)$ modulo torsion, $\Coh$ a category 
of the quasi-coherent sheaves on ${\cal E}(k)$ and $\Omega\subset {\goth S}_{\alpha,\beta,\gamma}(k)$ 
a two-sided ideal generated by the  central  elements
\displaymath
\left\{
\begin{array}{ccc}
\Omega_1 &=& x_1^2+x_2^2+x_3^2+x_4^2,\\
\Omega_2 &=&  x_2^2+{1+\beta\over 1-\gamma}x_3^2+{1-\beta\over 1+\alpha}x_4^2.
\end{array}
\right.
\enddisplaymath
\end{thm}

\vskip1cm\noindent
{\bf Guide to the literature.}
The Serre isomorphism has been established in [Serre 1957] \cite{Ser1}.
The concrete example of non-commutative ring satisfying the Serre isomorphism
for complex elliptic curves was constructed in [Sklyanin 1982] \cite{Skl1};
for a subsequent development,  see   [Sklyanin 1982] \cite{Skl2}.  
For a generalization of the Sklyanin algebras we refer the reader to
the papers  by  [Feigin \& Odesskii 1989]  \cite{FeOd1}  and  [Feigin \& Odesskii  1993]  \cite{FeOd2}. 
A general approach to non-commutative algebraic geometry can be found in 
[Artin \& Schelter 1987]  \cite{ArtSch1}  and [Artin \& van den Bergh  1990]  \cite{ArtVdb1}.  
The topic inspired  geometers in [Artin, Tate \& van den Bergh 1990]  \cite{ArtTatVdb1}  
and grew further in [Smith \& Stafford  1992]  \cite{SmiSta1} and [Artin \& Zhang 1994]  \cite{ArtZha1}.  
For a survey of this already formidable area  of mathematics we refer the reader to 
[Stafford \& van den Bergh 2001]  \cite{StaVdb1} and [Odesski 2002]  \cite{Ode1}.





\chapter{Non-commutative Trends in Algebraic Geometry}
There exist several more or less independent approaches to   non-commutative algebraic geometry 
due to [Bondal \& Orlov 2001]  \cite{BonOrl1},   [Kapranov 1998] \cite{Kap1},  [Kontsevich 1997]  \cite{Kon1},
 [van~Oystaeyen 1975]   \cite{O1} and [van~Oystaeyen \&  Verschoren 1980]   \cite{OV};   
 we shall briefly review them below.  Let us especially  mention  a powerful  {\it reconstruction theory}  
 developed by A.~L.~Rosenberg  in  the monograph [Rosenberg 1995]  \cite{ROS}. 
An new interesting  trend appears in the paper by [Laudal 2000]  \cite{Lau1}.

 \index{derived category}
\section{Derived categories}
Recall that in algebraic geometry the category of finitely generated 
projective modules over the coordinate ring of variety $X$ is isomorphic 
to the category of coherent sheaves over $X$ (Serre's Theorem);  
in other words,  variety $X$ can be recovered from a ``derived'' 
category of modules over commutative rings. This observation  has been
formalized by A.~Bondal and D.~Orlov   and  can be adapted to  the case of 
non-commutative rings.  
\begin{dfn}
Let {\bf A}  be an abelian category and let {\bf Kom~(A)}  be the category of complexes 
over {\bf A},  i.e. a category whose objects are co-chain complexes of abelian groups
and arrows are  morphisms of the complexes.  By a derived category ${\cal D}${\bf (A)}
of abelian category {\bf A} one understands such a category that a functor 
\displaymath
Q: \hbox{{\bf Kom~(A)}}\longrightarrow {\cal D} \hbox{{\bf (A)}} 
\enddisplaymath
has the following properties:

\medskip
(i) $Q(f)$ is an isomorphism for any quasi-isomorphism $f$ in {\bf Kom ~(A)},
i.e. an arrow of {\bf A} which gives an isomorphism of  {\bf Kom ~(A)};

\smallskip
(ii)  any functor $F:$ {\bf Kom~(A)} $\to {\cal D}$  transforming quasi-isomorphism
into isomorphism can be uniquely factored through ${\cal D}${\bf (A)}. 
\end{dfn}
\begin{exm}
\textnormal{
If {\bf CRng} is the category of commutative rings and {\bf Mod} the category  of
modules over the rings,  then  
\displaymath
\hbox{{\bf Mod}} \cong {\cal D} \hbox{{\bf (CRng)}}. 
\enddisplaymath
}
\end{exm}
\begin{dfn}
If $X$ is the category of algebraic varieties, then by ${\cal D}_{coh} (X)$ one understands 
the derived category of coherent sheaves on $X$. 
\end{dfn}
The following result is known as the Reconstruction Theorem for smooth projective varieties;
this remarkable theorem has been proved in 1990's by Alexei Bondal and Dmitri Orlov. 
 \index{Bondal-Orlov Reconstruction Theorem}
\begin{thm}
{\bf (Bondal and Orlov)}
Let $X$ be a smooth irreducible projective variety with ample canonical and 
anti-canonical sheaf.   If for two projective varieties $X$ and $X'$ the 
categories 
\displaymath
{\cal D}_{coh} (X)\cong {\cal D}_{coh} (X')
 \enddisplaymath
are equivalent,  then $X\cong X'$ are isomorphic projective varieties.   
\end{thm}
\begin{rmk}
\textnormal{
It is known that certain projective abelian varieties and $K3$ surfaces without restriction 
on the sheaves can be non-isomorphic but have equivalent derived categories of coherent
sheaves;  thus the requirements on the sheaves cannot be dropped. 
}
\end{rmk}
\begin{rmk}
\textnormal{
The Reconstruction Theorem extends to the category of modules over the non-commutative rings;
therefore such a theorem can be regarded as a categorical foundation of the noncommutative algebraic 
geometry of M.~Artin, J.~Tate and M.~van den Bergh.  
}
\end{rmk}

 \index{non-commutative thickening}

\section{Non-commutative thickening}
Let $R$ be an associative algebra over ${\Bbb C}$ and 
$R_{ab}=R/[R,R]$  its abelianization by the commutator $[R,R]=\{xyx^{-1}y^{-1} ~| ~x,y\in R\}$.
Denote by  
\displaymath
X_{ab}= \hbox{{\bf Spec}} ~(R_{ab})
\enddisplaymath
the space of all prime ideals of $R_{ab}$ endowed with the Zariski topology.
\begin{rmk}
\textnormal{
The naive aim of non-commutative algebraic geometry consists in construction of 
an embedding
\displaymath
X_{ab}\hookrightarrow X\cong ~\hbox{{\bf Spec}} ~(R),
\enddisplaymath
where $X$ is some ``noncommutative space'' associated to the homomorphism 
of rings $R\to R_{ab}$. 
}
\end{rmk}
 \index{Kapranov Theorem}
\begin{dfn}
{\bf (Kapranov)}
By a non-commutative thickening of the space $X_{ab}$ one understands 
a sheaf of non-commutative rings ${\cal O}^{NC}$ (corresponding to the formal
small neighborhood  of $X_{ab}$ in $X$) which is the completion 
of the algebra $R$ by the iterated commutators of the form
\displaymath
\lim_{n\to\infty} [a_1, [a_2,\dots [a_{n-1},a_n]]],  \quad a_i\in R,
\enddisplaymath
such that the commutators are small in a completion of the Zariski topology. 
\end{dfn}
\begin{dfn}
{\bf (Kapranov)}
By a non-commutative scheme $X$ one understands the glued ringed spaces 
of the form $(X_{ab}, {\cal O}^{NC})$. 
\end{dfn}
\begin{thm}
{\bf (Uniqueness of non-commutative thickening)}
If $R$ is a finitely generated commutative algebra,  then its non-commutative
thickening $(X_{ab}, {\cal O}^{NC})$ is unique up to an isomorphism of $R$.
\end{thm}
\begin{rmk}
\textnormal{
If $X_{ab}=M$ is a smooth manifold,  then its non-commutative thickening 
 $(X_{ab}, {\cal O}^{NC})$ induces  to a new differential-geometric structure
 on $M$ which corresponds to  additional characteristic classes of $M$. 
}
\end{rmk}

 \index{deformation quantization}
 \index{Poisson bracket}
 \index{Poisson manifold}
\section{Deformation quantization of Poisson manifolds}
Roughly speaking,  deformation quantization generalizes Rieffel's construction 
of the $n$-dimensional noncommutative torus  from the algebra $C(T^n)$.  It was
shown by M.~Kontsevich that the Poisson bracket is all it takes  to obtain a family of non-commutative associative 
algebras  from the algebra of functions on a manifold;  manifolds with such a bracket  are called  {\it Poisson manifolds}. 
Below we give details of the construction.

Let $M$ be a manifold endowed with a smooth (or analytic, or algebraic) structure;   let ${\cal F}(M)$ 
be the (commutative) algebra of  function on $M$.  The idea is to keep the point-wise addition of functions
on $M$ and replace the point-wise multiplication by a non-commutative   (but associative) binary operation $\ast_{\hbar}$
depending on a deformation parameter $\hbar$ so that $\hbar=0$ corresponds to the usual point-wise multiplication  
of functions on $M$.  One can unfold the non-commutative multiplication in an infinite series in parameter $\hbar$
\displaymath
f\ast_{\hbar} g=fg+\{f,g\}\hbar +O(\hbar),
 \enddisplaymath
where the bracket
\displaymath
\{\bullet, \bullet\}:  {\cal F}(M)\otimes {\cal F}(M)\to {\cal F}(M)
 \enddisplaymath
is bilinear and without loss of generality one can assume that $\{f, g\}=-\{g,f\}$. 
It follows from the associativity of  operation $\ast_{\hbar}$ that
\displaymath
\left\{
\begin{array}{ccc}
\{f, gh\} &= &  \{fg\}h+\{f,h\}g\\ 
\{f, \{g,h\}\} &+& \{h, \{f, g\}\}+\{g, \{h,f\}\}.
\end{array}
\right.
\enddisplaymath
The first equation implies that $\{f,f\}=0$ and the second is the {\it Jacobi identity};
in other words,  the non-commutative multiplication $\ast_{\hbar}$  gives rise to 
 a {\it Poisson bracket}  $\{\bullet, \bullet\}$ on the manifold $M$.  
It is remarkable, that the converse is true.
 \index{Kontsevich Theorem}
\begin{thm}
{\bf (Kontsevich)}
If $M$ be a manifold endowed with the Poisson bracket $\{\bullet, \bullet\}$,
then there exists an associative non-commutative multiplication on ${\cal F}(M)$
defined by the formula
\displaymath
f\ast_{\hbar} g:= fg+\{f,g\}+\sum_{n\ge 2} B_n(f,g)\hbar^n,
 \enddisplaymath
where $B_n(f,g)$ are certain explicit operators defined by the Poisson bracket.     
\end{thm}
\begin{exm}
{\bf (Rieffel)}
\textnormal{
Let  $M\cong T^n$ and ${\cal F}(M)\cong C^{\infty}(T^n)$;  let us define  the Poisson bracket 
by the formula
\displaymath
\{f,g\}=\sum \theta_{ij} {\partial f\over\partial x_i}   {\partial g\over\partial x_j}, 
\enddisplaymath
where $\Theta=(\theta_{ij})$ is a real skew symmetric matrix.   Then the non-commutative multiplication $\ast_{{\hbar}=1}$  
defines a (smooth)  $n$-dimensional noncommutative torus  ${\cal A}_{\Theta}$. 
}
\end{exm}

 \index{Oystaeyen-Verschoren geometry}

\section{Algebraic geometry of non-commutative rings}
In 1980's Fred van Oystaeyen and Alain Verschoren initiated a vast program aimed to
develop a non-commutative version of algebraic geometry based on the notion of a 
spectrum of non-commutative ring $R$;  such spectra consist of the maximal (left and right)
ideals of $R$.  A difficulty here that non-commutative rings occurring in practice are simple
rings, i.e. has no ideals whatsoever;  the difficulty can be overcome, see the original works of
the above authors.  Below we briefly review some of the constructions.     
\begin{dfn}
For an (associative) algebra $A$ the maximal ideal spectrum {\bf Max}  $(A)$ of $A$
is the set of all maximal two-sided ideals $M\subset A$ equipped with the non-commutative
Zariski topology,  i.e.  a topology with the typical open set of {\bf Max} $(A)$ given by the formula
\displaymath
{\Bbb X}(I)=\{M\in \hbox{{\bf Max}} (A) ~|~ I\not\subset M\}.  
\enddisplaymath
\end{dfn}
The non-commutative sheaves ${\cal O}^{nc}_A$ associated to $A$ can be defined as
follows. 
\begin{dfn}
The  ${\cal O}^{nc}_A$ is defined by taking the sections over the typical open set 
${\Bbb X}(I)$ of a two-sided ideal $I\subset$ {\bf Max} $(A)$ according to the formula
\displaymath
\Gamma ({\Bbb X}(I), {\cal O}_A^{nc}):=
\{\delta\in\Sigma ~|~ \exists l\in {\Bbb N} ~:~ I^l\delta\subset A\},
\enddisplaymath
where $\Sigma$ is the central simple algebra of $A$. 
\end{dfn}
\begin{rmk}
\textnormal{
One can develop all essential features which  the ``non-com\-muta\-tive scheme'' 
({\bf Max} $(A), {\cal O}_A^{nc}$) must have;    the subject, however,  lies beyond
the scope of present notes and we refer the interested reader to the corresponding 
literature.  
}
\end{rmk}

 \index{Rosenberg Reconstruction Theory}

\vskip1cm\noindent
{\bf Guide to the literature.}
The derived categories were studied by [Bondal \& Orlov 2001]  \cite{BonOrl1},
see also [Bondal \& Orlov 2002]  \cite{BonOrl2}.   
The non-commutative thickening was introduced by  [Kapranov 1998] \cite{Kap1}.  
The deformation quantization of Poisson manifolds was studied by  [Kontsevich 1997]
\cite{Kon1}.  
For the basics of  algebraic geometry of  non-commutative rings we refer the reader
to the monographs by  [van~Oystaeyen 1975]   \cite{O1}   and  [van~Oystaeyen \&  Verschoren 1980]   \cite{OV}.  
For the {\it reconstruction theory}   of  A.~L.~Rosenberg,   see  the monograph [Rosenberg 1995]  \cite{ROS}. 
An new interesting  trend in non-commutative algebraic geometry can be found  in the paper by [Laudal 2003]  \cite{Lau1}.

 






\begin{theindex}

  \item $C^*$-algebra, 37
  \item $K$-rational elliptic curve, 161
  \item $K$-rational points, iii
  \item $L$-function of noncommutative torus, 170
  \item $W^*$-algebra, 219
  \item $\alpha$-ample invertible sheaf, 275
  \item $\ell$-adic cohomology, 191
  \item $p$-adic invariant, 79
  \item ${\Bbb Q}$-curve, 17, 164
  \item ${\Bbb Q}$-rank, 165

  \indexspace

  \item abelian semigroup, 44
  \item abelianized Handelman invariant, 85
  \item AF-algebra, 52
  \item Alexander polynomial, 22, 82
  \item almost-culminating period, 165
  \item Anosov automorphism, 19
  \item Anosov bundle, 89
  \item Anosov diffeomorphism, 92
  \item Anosov map, 69, 80
  \item arithmetic complexity, 17, 162
  \item arrows, 30
  \item Artin $L$-function, 189
  \item Artin reciprocity, 190
  \item Artin's representation of braid group, 260
  \item Artin-Mazur zeta function, 198
  \item Artin-van den Bergh Theorem, 275
  \item Atiyah's realization of K-homology, 243
  \item Atiyah-Singer Index Theorem, 240
  \item Atiyah-Singer Theorem for family of elliptic operators, 232
  \item Atkinson Theorem, 238

  \indexspace

  \item basic isomorphism, 114
  \item Baum-Connes Conjecture, 251
  \item Birch and Swinnerton-Dyer Conjecture, 169
  \item Bondal-Orlov Reconstruction Theorem, 278
  \item Bost-Connes system, 234
  \item Bost-Connes Theorem, 235
  \item Bott periodicity, 243
  \item braid, 257, 259
  \item braid closure, 260
  \item braid group, 257
  \item Bratteli diagram, 53
  \item Brown-Douglas-Fillmore Theorems, 245
  \item Brown-Douglas-Fillmore Theory, 244

  \indexspace

  \item category, 27
  \item Chebyshev polynomial, 203
  \item Chern character formula, 192, 243
  \item Chern character of elliptic operator, 239
  \item Chern-Connes character, 228, 230
  \item class field theory, 189
  \item co-product, 268
  \item coarse Baum-Connes Conjecture, 255
  \item coarse geometry, 253
  \item coarse map, 254
  \item commutative Desargues space, 218
  \item commutative diagram, 27
  \item compact quantum group, 269
  \item complex algebraic curve, 122
  \item complex modulus, 11, 24
  \item complex multiplication, 16, 24, 148
  \item complex projective variety, 129
  \item complex torus, 11, 108
  \item configurations, 213
  \item Connes geometry, 225
  \item Connes Invariant, 226
  \item Connes'  Index Theorem, 233
  \item Connes-Moscovici Theorem, 231
  \item continuous geometry, 219
  \item contravariant funcor, 32
  \item covariant functor, 32
  \item covariant representation, 41
  \item crossed product, 40
  \item culminating period, 165
  \item Cuntz algebra, 60
  \item Cuntz-Krieger algebra, 59, 87
  \item Cuntz-Krieger Crossed Product Theorem, 60
  \item Cuntz-Krieger invariant, 88
  \item Cuntz-Krieger Theorem, 60
  \item cyclic complex, 230
  \item cyclic homology, 228, 229

  \indexspace

  \item deformation quantization, 6, 279
  \item derived category, 277
  \item Desargues axiom, 217
  \item Desargues projective plane, 217
  \item Deuring Theorem, 180
  \item dimension function, 216
  \item dimension group, 54
  \item Dirichlet $L$-series, 189
  \item Dirichlet Unit Theorem, 199
  \item double commutant, 220
  \item Double Commutant Theorem, 220
  \item dyadic number, 59

  \indexspace

  \item Effros-Shen algebra, 54
  \item Eichler-Shimura theory, 170
  \item Elliott Theorem, 55
  \item elliptic curve, i, iii, 10, 107
  \item equivalent foliations, 72
  \item essential spectrum, 238

  \indexspace

  \item factor, 221
  \item faithful functor, 33
  \item Fermion algebra, 58
  \item finite field, 191
  \item finite geometry, 213
  \item finitely generated graded module, 272
  \item flow of weights, 227
  \item foliated space, 232
  \item Fredholm operator, 237
  \item Frobenius element, 189
  \item Frobenius endomorphism, 234
  \item Frobenius map, 148
  \item function $\pi(n)$, 180
  \item function $L({\cal A}_{RM}^{2n},s)$, 182
  \item functor, i, 31
  \item functor $K_0$, 45
  \item functor $K_0^+$, 44
  \item functor $K_1$, 46
  \item fundamental AF-algebra, 89
  \item fundamental phenomenon, 12
  \item fundamental theorem of projective geometry, 216

  \indexspace

  \item Gauss method, 25
  \item Gelfand Duality, 272
  \item Gelfand-Naimark Theorem, 40
  \item Gelfond-Schneider Theorem, 204
  \item geodesic spectrum, 157
  \item Geometrization Conjecture, 23
  \item Glimm Theorem, 59
  \item GNS-construction, 40
  \item golden mean, 56
  \item Gr\"ossencharacter, 159, 174, 178
  \item Gromov-Lawson Conjecture, 253
  \item Grothendieck map, 45

  \indexspace

  \item Handelman invariant, 21, 57, 68
  \item Hasse lemma, 175
  \item Hasse-Weil $L$-function, 170
  \item Hecke $C^*$-algebra, 234
  \item Hecke lemma, 159
  \item Hecke operator, 235
  \item higher-dimensional\linebreak  noncommutative torus, 47
  \item Hilbert module, 246
  \item Hirzebruch Signature, 250
  \item Hirzebruch-Riemann-Roch\linebreak   Formula, 240
  \item Hochschild homology, 228
  \item Hochschild-Kostant-Rosenberg Theorem, 228
  \item Hopf algebra, 267
  \item hyper-finite $W^*$-algebra, 221

  \indexspace

  \item imaginary quadratic number, 24
  \item inverse temperature, 226
  \item irrational rotation algebra, 10
  \item isogeny, 17, 148
  \item isomorphic categories, 33
  \item isotropy subgroup, 189

  \indexspace

  \item Jacobi elliptic curve, 11, 108
  \item Jacobi-Perron continued fraction, 68
  \item Jacobian of measured foliation, 70
  \item Jones algebra, 261
  \item Jones Polynomial, 262
  \item Jones projection, 259
  \item Jones Theorem, 258

  \indexspace

  \item K-homology, 237, 243
  \item K-theory, 42
  \item Kapranov Theorem, 279
  \item Kasparov module, 248
  \item Kasparov Theorem, 251
  \item Kasparov's KK-theory, 237, 246
  \item KK-groups, 247
  \item KMS condition, 226
  \item KMS state, 235
  \item Kontsevich Theorem, 280
  \item Kronecker-Weber Theorem, 16

  \indexspace

  \item Landstadt-Takai duality, 42
  \item Langalands Conjecture for noncommutative torus, 185
  \item Langlands program, iii, 148, 181
  \item localization formula, 171
  \item locally compact quantum group, 269

  \indexspace

  \item Manin's quantum plane, 266
  \item mapping class group, 139
  \item Markov category, 196
  \item Markov move, 260
  \item Markov Theorem, 261
  \item measured foliation, 70, 116
  \item Mischenko-Fomenko Theorem, 253
  \item Mordell-N\'eron Theorem, 162
  \item Morita equivalence, 9
  \item morphism of category, 30
  \item Muir symbols, 163

  \indexspace

  \item NCG, i
  \item non-commutative Pontryagin duality, 263, 268
  \item non-commutative reciprocity, 169
  \item non-commutative thickening, 278
  \item noncommutative torus, i, 5, 20, 50
  \item Novikov Conjecture, 249, 250
  \item Novikov Conjecture for hyperbolic\linebreak  groups, 230

  \indexspace

  \item objects, 30
  \item obstruction theory, 93
  \item Oystaeyen-Verschoren geometry, 281

  \indexspace

  \item Pappus axiom, 217
  \item Pappus projective space, 217
  \item partial isometry, 38
  \item Perron-Frobenius eigenvalue, 20
  \item Perron-Frobenius eigenvector, 20
  \item Pimsner-Voiculescu Theorem, 47
  \item Plante group, 94
  \item Poisson bracket, 6, 279
  \item Poisson manifold, 279
  \item Pontryagin characteristic classes, 250
  \item positive matrix, 20
  \item positive scalar curvature, 252
  \item projection, 38
  \item projective pseudo-lattice, 119
  \item projective space ${\Bbb P}^n(D)$, 213, 215
  \item pseudo-Anosov map, 65
  \item pseudo-lattice, 118, 173

  \indexspace

  \item quadratic irrationality, 10
  \item quadric surface, 10
  \item quantum group, 266
  \item quasi-coherent sheaf, 272

  \indexspace

  \item rank conjecture, 209
  \item rank of elliptic curve, 17, 162
  \item rational elliptic curve, 200
  \item rational identity, 213
  \item real multiplication, iii, 9, 21, 137, 148, 151
  \item Rieffel-Schwarz Theorem, 49
  \item Riemann Hypothesis, 234
  \item Riemann zeta function, 236
  \item robust torus bundle, 101
  \item Rosenberg Reconstruction Theory, 282

  \indexspace

  \item scaled unit, 114
  \item Serre $C^*$-algebra, 107, 129, 191
  \item Serre isomorphism, 272
  \item Seventh Hilbert Problem, 208
  \item shift automorphism, 56
  \item Shimura-Taniyama Conjecture, 181
  \item signature of pseudo-Anosov map, 104
  \item skew field, 213
  \item Sklyanin $\ast$-algebra, 112
  \item Sklyanin algebra, 13, 110, 275
  \item Sklyanin Theorem, 276
  \item spectrum of $C^*$-algebra, 61
  \item stable isomorphism, 9, 39
  \item stationary AF-algebra, 55, 57
  \item stationary dimension group, 57
  \item Strong Novikov Conjecture, 251
  \item subfactor, 257
  \item subshift of finite type, 85
  \item supernatural number, 59
  \item surface map, 65
  \item suspension, 46

  \indexspace

  \item Teichm\"uller space, 116
  \item tight hyperbolic matrix, 100
  \item Tomita-Takesaki Theorem, 226
  \item Tomita-Takesaki Theory, 225
  \item topological K-theory, 241
  \item toric AF-algebra, 107, 122
  \item torus bundle, 84
  \item trace invariant, 261
  \item trace of Frobenius endomorphism, 191
  \item transcendental number theory, 148, 203
  \item twisted homogeneous coordinate ring, 111, 274
  \item type {\bf III}$_0$,  {\bf III}$_{\lambda}$, {\bf III}$_1$-factors, 
		227
  \item type {\bf I}, {\bf II}$_1$, {\bf II}$_{\infty}$,  {\bf III}-factors, 
		221

  \indexspace

  \item UHF-algebra, 57
  \item unit of algebraic number field, 178
  \item unitary, 38

  \indexspace

  \item von Neumann geometry, 222

  \indexspace

  \item weak topology, 219
  \item Weierstrass $\wp$ function, 12, 109
  \item weight, 225
  \item Weil Conjectures, 191

\end{theindex}


\end{document}